
\ifx\shlhetal\undefinedcontrolsequence\let\shlhetal\relax\fi
\def\fmtname{AmS-TeX}

\def\fmtversion{2.1}
\catcode`\@=11
\ifx\amstexloaded@\relax\catcode`\@=\active
  \endinput\else\let\amstexloaded@\relax\fi
\newlinechar=`\^^J
\def\W@{\immediate\write\sixt@@n}
\def\CR@{\W@{^^J\fmtname - Version \fmtversion^^J%
COPYRIGHT 1985, 1990, 1991 - AMERICAN MATHEMATICAL SOCIETY^^J%
Use of this macro package is not restricted provided^^J%
each use is acknowledged upon publication.^^J}}
\CR@ \everyjob{\CR@}
\message{Loading definitions for}
\message{misc utility macros,}
\toksdef\toks@@=2
\long\def\rightappend@#1\to#2{\toks@{\\{#1}}\toks@@
 =\expandafter{#2}\xdef#2{\the\toks@@\the\toks@}\toks@{}\toks@@{}}
\def\alloclist@{}
\newif\ifalloc@
\def\showallocations{{\def\\{\immediate\write\m@ne}\alloclist@}\alloc@true}
\def\alloc@#1#2#3#4#5{\global\advance\count1#1by\@ne
 \ch@ck#1#4#2\allocationnumber=\count1#1
 \global#3#5=\allocationnumber
 \edef\next@{\string#5=\string#2\the\allocationnumber}%
 \expandafter\rightappend@\next@\to\alloclist@}
\newcount\count@@
\newcount\count@@@
\def\FN@{\futurelet\next}
\def\DN@{\def\next@}
\def\DNii@{\def\nextii@}
\def\RIfM@{\relax\ifmmode}
\def\RIfMIfI@{\relax\ifmmode\ifinner}
\def\setboxz@h{\setbox\z@\hbox}
\def\wdz@{\wd\z@}
\def\boxz@{\box\z@}
\def\setbox@ne{\setbox\@ne}
\def\wd@ne{\wd\@ne}
\def\iterate{\body\expandafter\iterate\else\fi}
\def\err@#1{\errmessage{AmS-TeX error: #1}}
\newhelp\defaulthelp@{Sorry, I already gave what help I could...^^J
Maybe you should try asking a human?^^J
An error might have occurred before I noticed any problems.^^J
``If all else fails, read the instructions.''}
\def\Err@{\errhelp\defaulthelp@\err@}
\def\eat@#1{}
\def\in@#1#2{\def\in@@##1#1##2##3\in@@{\ifx\in@##2\in@false\else\in@true\fi}%
 \in@@#2#1\in@\in@@}
\newif\ifin@
\def\space@.{\futurelet\space@\relax}
\space@. %
\newhelp\athelp@
{Only certain combinations beginning with @ make sense to me.^^J
Perhaps you wanted \string\@\space for a printed @?^^J
I've ignored the character or group after @.}
{\catcode`\~=\active 
 \lccode`\~=`\@ \lowercase{\gdef~{\FN@\at@}}}
\def\at@{\let\next@\at@@
 \ifcat\noexpand\next a\else\ifcat\noexpand\next0\else
 \ifcat\noexpand\next\relax\else
   \let\next\at@@@\fi\fi\fi
 \next@}
\def\at@@#1{\expandafter
 \ifx\csname\space @\string#1\endcsname\relax
  \expandafter\at@@@ \else
  \csname\space @\string#1\expandafter\endcsname\fi}
\def\at@@@#1{\errhelp\athelp@ \err@{\Invalid@@ @}}
\def\atdef@#1{\expandafter\def\csname\space @\string#1\endcsname}
\newhelp\defahelp@{If you typed \string\define\space cs instead of
\string\define\string\cs\space^^J
I've substituted an inaccessible control sequence so that your^^J
definition will be completed without mixing me up too badly.^^J
If you typed \string\define{\string\cs} the inaccessible control sequence^^J
was defined to be \string\cs, and the rest of your^^J
definition appears as input.}
\newhelp\defbhelp@{I've ignored your definition, because it might^^J
conflict with other uses that are important to me.}
\def\define{\FN@\define@}
\def\define@{\ifcat\noexpand\next\relax
 \expandafter\define@@\else\errhelp\defahelp@                               
 \err@{\string\define\space must be followed by a control
 sequence}\expandafter\def\expandafter\nextii@\fi}                          
\def\undefined@@@@@@@@@@{}
\def\preloaded@@@@@@@@@@{}
\def\next@@@@@@@@@@{}
\def\define@@#1{\ifx#1\relax\errhelp\defbhelp@                              
 \err@{\string#1\space is already defined}\DN@{\DNii@}\else
 \expandafter\ifx\csname\expandafter\eat@\string                            
 #1@@@@@@@@@@\endcsname\undefined@@@@@@@@@@\errhelp\defbhelp@
 \err@{\string#1\space can't be defined}\DN@{\DNii@}\else
 \expandafter\ifx\csname\expandafter\eat@\string#1\endcsname\relax          
 \global\let#1\undefined\DN@{\def#1}\else\errhelp\defbhelp@
 \err@{\string#1\space is already defined}\DN@{\DNii@}\fi
 \fi\fi\next@}

\def\predefine#1#2{\let#1#2}
\def\undefine#1{\let#1\undefined}
\message{page layout,}
\newdimen\captionwidth@
\captionwidth@\hsize
\advance\captionwidth@-1.5in
\def\pagewidth#1{\hsize#1\relax
 \captionwidth@\hsize\advance\captionwidth@-1.5in}
\def\pageheight#1{\vsize#1\relax}
\def\hcorrection#1{\advance\hoffset#1\relax}
\def\vcorrection#1{\advance\voffset#1\relax}
\message{accents/punctuation,}

\let\graveaccent\`
\let\acuteaccent\'
\let\tildeaccent\~
\let\hataccent\^
\let\underscore\_
\let\B\=
\let\D\.
\let\ic@\/
\def\/{\unskip\ic@}
\def\textfonti{\the\textfont\@ne}
\def\t#1#2{{\edef\next@{\the\font}\textfonti\accent"7F \next@#1#2}}
\def~{\unskip\nobreak\ \ignorespaces}
\def\.{.\spacefactor\@m}
\atdef@;{\leavevmode\null;}
\atdef@:{\leavevmode\null:}
\atdef@?{\leavevmode\null?}
\edef\@{\string @}
\def\relaxnext@{\let\next\relax}
\atdef@-{\relaxnext@\leavevmode
 \DN@{\ifx\next-\DN@-{\FN@\nextii@}\else
  \DN@{\leavevmode\hbox{-}}\fi\next@}%
 \DNii@{\ifx\next-\DN@-{\leavevmode\hbox{---}}\else
  \DN@{\leavevmode\hbox{--}}\fi\next@}%
 \FN@\next@}
\def\srdr@{\kern.16667em}
\def\drsr@{\kern.02778em}
\def\sldl@{\drsr@}
\def\dlsl@{\srdr@}
\atdef@"{\unskip\relaxnext@
 \DN@{\ifx\next\space@\DN@. {\FN@\nextii@}\else
  \DN@.{\FN@\nextii@}\fi\next@.}%
 \DNii@{\ifx\next`\DN@`{\FN@\nextiii@}\else
  \ifx\next\lq\DN@\lq{\FN@\nextiii@}\else
  \DN@####1{\FN@\nextiv@}\fi\fi\next@}%
 \def\nextiii@{\ifx\next`\DN@`{\sldl@``}\else\ifx\next\lq
  \DN@\lq{\sldl@``}\else\DN@{\dlsl@`}\fi\fi\next@}%
 \def\nextiv@{\ifx\next'\DN@'{\srdr@''}\else
  \ifx\next\rq\DN@\rq{\srdr@''}\else\DN@{\drsr@'}\fi\fi\next@}%
 \FN@\next@}

\def\textfontii{\the\textfont\tw@}
\def\lbrace@{\delimiter"4266308 }
\def\rbrace@{\delimiter"5267309 }
\def\{{\RIfM@\lbrace@\else{\textfontii f}\spacefactor\@m\fi}
\def\}{\RIfM@\rbrace@\else
 \let\@sf\empty\ifhmode\edef\@sf{\spacefactor\the\spacefactor}\fi
 {\textfontii g}\@sf\relax\fi}
\let\lbrace\{
\let\rbrace\}
\def\AmSTeX{{\textfontii A\kern-.1667em%
  \lower.5ex\hbox{M}\kern-.125emS}-\TeX}
\message{line and page breaks,}
\def\vmodeerr@#1{\Err@{\string#1\space not allowed between paragraphs}}
\def\mathmodeerr@#1{\Err@{\string#1\space not allowed in math mode}}
\def\linebreak{\RIfM@\mathmodeerr@\linebreak\else
 \ifhmode\unskip\unkern\break\else\vmodeerr@\linebreak\fi\fi}

\newskip\saveskip@
\def\allowlinebreak{\RIfM@\mathmodeerr@\allowlinebreak\else
 \ifhmode\saveskip@\lastskip\unskip
 \allowbreak\ifdim\saveskip@>\z@\hskip\saveskip@\fi
 \else\vmodeerr@\allowlinebreak\fi\fi}
\def\nolinebreak{\RIfM@\mathmodeerr@\nolinebreak\else
 \ifhmode\saveskip@\lastskip\unskip
 \nobreak\ifdim\saveskip@>\z@\hskip\saveskip@\fi
 \else\vmodeerr@\nolinebreak\fi\fi}
\def\newline{\relaxnext@
 \DN@{\RIfM@\expandafter\mathmodeerr@\expandafter\newline\else
  \ifhmode\ifx\next\par\else
  \expandafter\unskip\expandafter\null\expandafter\hfill\expandafter\break\fi
  \else
  \expandafter\vmodeerr@\expandafter\newline\fi\fi}%
 \FN@\next@}
\def\dmatherr@#1{\Err@{\string#1\space not allowed in display math mode}}
\def\nondmatherr@#1{\Err@{\string#1\space not allowed in non-display math
 mode}}
\def\onlydmatherr@#1{\Err@{\string#1\space allowed only in display math mode}}
\def\nonmatherr@#1{\Err@{\string#1\space allowed only in math mode}}
\def\mathbreak{\RIfMIfI@\break\else
 \dmatherr@\mathbreak\fi\else\nonmatherr@\mathbreak\fi}
\def\nomathbreak{\RIfMIfI@\nobreak\else
 \dmatherr@\nomathbreak\fi\else\nonmatherr@\nomathbreak\fi}
\def\allowmathbreak{\RIfMIfI@\allowbreak\else
 \dmatherr@\allowmathbreak\fi\else\nonmatherr@\allowmathbreak\fi}
\def\pagebreak{\RIfM@
 \ifinner\nondmatherr@\pagebreak\else\postdisplaypenalty-\@M\fi
 \else\ifvmode\removelastskip\break\else\vadjust{\break}\fi\fi}
\def\nopagebreak{\RIfM@
 \ifinner\nondmatherr@\nopagebreak\else\postdisplaypenalty\@M\fi
 \else\ifvmode\nobreak\else\vadjust{\nobreak}\fi\fi}
\def\nonvmodeerr@#1{\Err@{\string#1\space not allowed within a paragraph
 or in math}}
\def\vnonvmode@#1#2{\relaxnext@\DNii@{\ifx\next\par\DN@{#1}\else
 \DN@{#2}\fi\next@}%
 \ifvmode\DN@{#1}\else
 \DN@{\FN@\nextii@}\fi\next@}
\def\newpage{\vnonvmode@{\vfill\break}{\nonvmodeerr@\newpage}}
\def\smallpagebreak{\vnonvmode@\smallbreak{\nonvmodeerr@\smallpagebreak}}
\def\medpagebreak{\vnonvmode@\medbreak{\nonvmodeerr@\medpagebreak}}
\def\bigpagebreak{\vnonvmode@\bigbreak{\nonvmodeerr@\bigpagebreak}}
\def\NoBlackBoxes{\global\overfullrule\z@}
\def\BlackBoxes{\global\overfullrule5\p@}
\def\Invalid@#1{\def#1{\Err@{\Invalid@@\string#1}}}
\def\Invalid@@{Invalid use of }
\message{figures,}
\Invalid@\caption
\Invalid@\captionwidth
\newdimen\smallcaptionwidth@
\def\topspace{\mid@false\ins@}
\def\midspace{\mid@true\ins@}
\newif\ifmid@
\def\captionfont@{}
\def\ins@#1{\relaxnext@\allowbreak
 \smallcaptionwidth@\captionwidth@\gdef\thespace@{#1}%
 \DN@{\ifx\next\space@\DN@. {\FN@\nextii@}\else
  \DN@.{\FN@\nextii@}\fi\next@.}%
 \DNii@{\ifx\next\caption\DN@\caption{\FN@\nextiii@}%
  \else\let\next@\nextiv@\fi\next@}%
 \def\nextiv@{\vnonvmode@
  {\ifmid@\expandafter\midinsert\else\expandafter\topinsert\fi
   \vbox to\thespace@{}\endinsert}
  {\ifmid@\nonvmodeerr@\midspace\else\nonvmodeerr@\topspace\fi}}%
 \def\nextiii@{\ifx\next\captionwidth\expandafter\nextv@
  \else\expandafter\nextvi@\fi}%
 \def\nextv@\captionwidth##1##2{\smallcaptionwidth@##1\relax\nextvi@{##2}}%
 \def\nextvi@##1{\def\thecaption@{\captionfont@##1}%
  \DN@{\ifx\next\space@\DN@. {\FN@\nextvii@}\else
   \DN@.{\FN@\nextvii@}\fi\next@.}%
  \FN@\next@}%
 \def\nextvii@{\vnonvmode@
  {\ifmid@\expandafter\midinsert\else
  \expandafter\topinsert\fi\vbox to\thespace@{}\nobreak\smallskip
  \setboxz@h{\noindent\ignorespaces\thecaption@\unskip}%
  \ifdim\wdz@>\smallcaptionwidth@\centerline{\vbox{\hsize\smallcaptionwidth@
   \noindent\ignorespaces\thecaption@\unskip}}%
  \else\centerline{\boxz@}\fi\endinsert}
  {\ifmid@\nonvmodeerr@\midspace
  \else\nonvmodeerr@\topspace\fi}}%
 \FN@\next@}
\message{comments,}
\def\newcodes@{\catcode`\\12\catcode`\{12\catcode`\}12\catcode`\#12%
 \catcode`\%12\relax}
\def\oldcodes@{\catcode`\\0\catcode`\{1\catcode`\}2\catcode`\#6%
 \catcode`\%14\relax}
\def\comment{\newcodes@\endlinechar=10 \comment@}
{\lccode`\0=`\\
\lowercase{\gdef\comment@#1^^J{\comment@@#10endcomment\comment@@@}%
\gdef\comment@@#10endcomment{\FN@\comment@@@}%
\gdef\comment@@@#1\comment@@@{\ifx\next\comment@@@\let\next\comment@
 \else\def\next{\oldcodes@\endlinechar=`\^^M\relax}%
 \fi\next}}}
\def\pr@m@s{\ifx'\next\DN@##1{\prim@s}\else\let\next@\egroup\fi\next@}
\def\prime{{\null\prime@\null}}
\mathchardef\prime@="0230
\let\dsize\displaystyle

\let\ssize\scriptstyle

\message{math spacing,}
\def\,{\RIfM@\mskip\thinmuskip\relax\else\kern.16667em\fi}
\def\!{\RIfM@\mskip-\thinmuskip\relax\else\kern-.16667em\fi}
\let\thinspace\,
\let\negthinspace\!
\def\medspace{\RIfM@\mskip\medmuskip\relax\else\kern.222222em\fi}
\def\negmedspace{\RIfM@\mskip-\medmuskip\relax\else\kern-.222222em\fi}
\def\thickspace{\RIfM@\mskip\thickmuskip\relax\else\kern.27777em\fi}
\let\;\thickspace
\def\negthickspace{\RIfM@\mskip-\thickmuskip\relax\else
 \kern-.27777em\fi}
\atdef@,{\RIfM@\mskip.1\thinmuskip\else\leavevmode\null,\fi}
\atdef@!{\RIfM@\mskip-.1\thinmuskip\else\leavevmode\null!\fi}
\atdef@.{\RIfM@&&\else\leavevmode.\spacefactor3000 \fi}
\def\and{\DOTSB\;\mathchar"3026 \;}

\message{fractions,}
\def\frac#1#2{{#1\over#2}}

\newdimen\ex@
\ex@.2326ex
\Invalid@\thickness
\def\thickfrac{\relaxnext@
 \DN@{\ifx\next\thickness\let\next@\nextii@\else
 \DN@{\nextii@\thickness1}\fi\next@}%
 \DNii@\thickness##1##2##3{{##2\above##1\ex@##3}}%
 \FN@\next@}

\def\thickfracwithdelims#1#2{\relaxnext@\def\ldelim@{#1}\def\rdelim@{#2}%
 \DN@{\ifx\next\thickness\let\next@\nextii@\else
 \DN@{\nextii@\thickness1}\fi\next@}%
 \DNii@\thickness##1##2##3{{##2\abovewithdelims
 \ldelim@\rdelim@##1\ex@##3}}%
 \FN@\next@}
\def\binom#1#2{{#1\choose#2}}

\def\:{\nobreak\hskip.1111em\mathpunct{}\nonscript\mkern-\thinmuskip{:}\hskip
 .3333emplus.0555em\relax}
\def\snug{\unskip\kern-\mathsurround}
\message{smash commands,}
\def\topsmash{\top@true\bot@false\smash@}
\def\botsmash{\top@false\bot@true\smash@}
\newif\iftop@
\newif\ifbot@
\def\smash{\top@true\bot@true\smash@}
\def\smash@{\RIfM@\expandafter\mathpalette\expandafter\mathsm@sh\else
 \expandafter\makesm@sh\fi}
\def\finsm@sh{\iftop@\ht\z@\z@\fi\ifbot@\dp\z@\z@\fi\leavevmode\boxz@}
\message{large operator symbols,}
\def\LimitsOnSums{\global\let\slimits@\displaylimits}
\def\NoLimitsOnSums{\global\let\slimits@\nolimits}
\LimitsOnSums
\mathchardef\coprod@="1360       \def\coprod{\DOTSB\coprod@\slimits@}
\mathchardef\bigvee@="1357       \def\bigvee{\DOTSB\bigvee@\slimits@}
\mathchardef\bigwedge@="1356     \def\bigwedge{\DOTSB\bigwedge@\slimits@}
\mathchardef\biguplus@="1355     \def\biguplus{\DOTSB\biguplus@\slimits@}
\mathchardef\bigcap@="1354       \def\bigcap{\DOTSB\bigcap@\slimits@}
\mathchardef\bigcup@="1353       \def\bigcup{\DOTSB\bigcup@\slimits@}
\mathchardef\prod@="1351         \def\prod{\DOTSB\prod@\slimits@}
\mathchardef\sum@="1350          \def\sum{\DOTSB\sum@\slimits@}
\mathchardef\bigotimes@="134E    \def\bigotimes{\DOTSB\bigotimes@\slimits@}
\mathchardef\bigoplus@="134C     \def\bigoplus{\DOTSB\bigoplus@\slimits@}
\mathchardef\bigodot@="134A      \def\bigodot{\DOTSB\bigodot@\slimits@}
\mathchardef\bigsqcup@="1346     \def\bigsqcup{\DOTSB\bigsqcup@\slimits@}
\message{integrals,}
\def\LimitsOnInts{\global\let\ilimits@\displaylimits}
\def\NoLimitsOnInts{\global\let\ilimits@\nolimits}
\NoLimitsOnInts
\def\int{\DOTSI\intop\ilimits@}
\def\oint{\DOTSI\ointop\ilimits@}
\def\intic@{\mathchoice{\hskip.5em}{\hskip.4em}{\hskip.4em}{\hskip.4em}}
\def\negintic@{\mathchoice
 {\hskip-.5em}{\hskip-.4em}{\hskip-.4em}{\hskip-.4em}}
\def\intkern@{\mathchoice{\!\!\!}{\!\!}{\!\!}{\!\!}}
\def\intdots@{\mathchoice{\plaincdots@}
 {{\cdotp}\mkern1.5mu{\cdotp}\mkern1.5mu{\cdotp}}
 {{\cdotp}\mkern1mu{\cdotp}\mkern1mu{\cdotp}}
 {{\cdotp}\mkern1mu{\cdotp}\mkern1mu{\cdotp}}}
\newcount\intno@
\def\iint{\DOTSI\intno@\tw@\FN@\ints@}
\def\iiint{\DOTSI\intno@\thr@@\FN@\ints@}
\def\iiiint{\DOTSI\intno@4 \FN@\ints@}
\def\idotsint{\DOTSI\intno@\z@\FN@\ints@}
\def\ints@{\findlimits@\ints@@}
\newif\iflimtoken@
\newif\iflimits@
\def\findlimits@{\limtoken@true\ifx\next\limits\limits@true
 \else\ifx\next\nolimits\limits@false\else
 \limtoken@false\ifx\ilimits@\nolimits\limits@false\else
 \ifinner\limits@false\else\limits@true\fi\fi\fi\fi}
\def\multint@{\int\ifnum\intno@=\z@\intdots@                                
 \else\intkern@\fi                                                          
 \ifnum\intno@>\tw@\int\intkern@\fi                                         
 \ifnum\intno@>\thr@@\int\intkern@\fi                                       
 \int}                                                                      
\def\multintlimits@{\intop\ifnum\intno@=\z@\intdots@\else\intkern@\fi
 \ifnum\intno@>\tw@\intop\intkern@\fi
 \ifnum\intno@>\thr@@\intop\intkern@\fi\intop}
\def\ints@@{\iflimtoken@                                                    
 \def\ints@@@{\iflimits@\negintic@\mathop{\intic@\multintlimits@}\limits    
  \else\multint@\nolimits\fi                                                
  \eat@}                                                                    
 \else                                                                      
 \def\ints@@@{\iflimits@\negintic@
  \mathop{\intic@\multintlimits@}\limits\else
  \multint@\nolimits\fi}\fi\ints@@@}
\def\LimitsOnNames{\global\let\nlimits@\displaylimits}
\def\NoLimitsOnNames{\global\let\nlimits@\nolimits@}
\LimitsOnNames
\def\nolimits@{\relaxnext@
 \DN@{\ifx\next\limits\DN@\limits{\nolimits}\else
  \let\next@\nolimits\fi\next@}%
 \FN@\next@}
\message{operator names,}
\def\newmcodes@{\mathcode`\'"27\mathcode`\*"2A\mathcode`\."613A%
 \mathcode`\-"2D\mathcode`\/"2F\mathcode`\:"603A }
\def\operatorname#1{\mathop{\newmcodes@\kern\z@\fam\z@#1}\nolimits@}
\def\operatornamewithlimits#1{\mathop{\newmcodes@\kern\z@\fam\z@#1}\nlimits@}
\def\qopname@#1{\mathop{\fam\z@#1}\nolimits@}
\def\qopnamewl@#1{\mathop{\fam\z@#1}\nlimits@}
\def\arccos{\qopname@{arccos}}
\def\arcsin{\qopname@{arcsin}}
\def\arctan{\qopname@{arctan}}
\def\arg{\qopname@{arg}}
\def\cos{\qopname@{cos}}
\def\cosh{\qopname@{cosh}}
\def\cot{\qopname@{cot}}
\def\coth{\qopname@{coth}}
\def\csc{\qopname@{csc}}
\def\deg{\qopname@{deg}}
\def\det{\qopnamewl@{det}}
\def\dim{\qopname@{dim}}
\def\exp{\qopname@{exp}}
\def\gcd{\qopnamewl@{gcd}}
\def\hom{\qopname@{hom}}
\def\inf{\qopnamewl@{inf}}
\def\injlim{\qopnamewl@{inj\,lim}}
\def\ker{\qopname@{ker}}
\def\lg{\qopname@{lg}}
\def\lim{\qopnamewl@{lim}}
\def\liminf{\qopnamewl@{lim\,inf}}
\def\limsup{\qopnamewl@{lim\,sup}}
\def\ln{\qopname@{ln}}
\def\log{\qopname@{log}}
\def\max{\qopnamewl@{max}}
\def\min{\qopnamewl@{min}}
\def\Pr{\qopnamewl@{Pr}}
\def\projlim{\qopnamewl@{proj\,lim}}
\def\sec{\qopname@{sec}}
\def\sin{\qopname@{sin}}
\def\sinh{\qopname@{sinh}}
\def\sup{\qopnamewl@{sup}}
\def\tan{\qopname@{tan}}
\def\tanh{\qopname@{tanh}}
\def\varinjlim{\mathop{\vtop{\ialign{##\crcr
 \hfil\rm lim\hfil\crcr\noalign{\nointerlineskip}\rightarrowfill\crcr
 \noalign{\nointerlineskip\kern-\ex@}\crcr}}}}
\def\varprojlim{\mathop{\vtop{\ialign{##\crcr
 \hfil\rm lim\hfil\crcr\noalign{\nointerlineskip}\leftarrowfill\crcr
 \noalign{\nointerlineskip\kern-\ex@}\crcr}}}}
\def\varliminf{\mathop{\underline{\vrule height\z@ depth.2exwidth\z@
 \hbox{\rm lim}}}}

\newdimen\buffer@
\buffer@\fontdimen13 \tenex
\newdimen\buffer
\buffer\buffer@

\def\ResetBuffer{\fontdimen13 \tenex\buffer@\global\buffer\buffer@}
\def\shave#1{\mathop{\hbox{$\m@th\fontdimen13 \tenex\z@                     
 \displaystyle{#1}$}}\fontdimen13 \tenex\buffer}

\message{multilevel sub/superscripts,}
\Invalid@\\
\def\Let@{\relax\iffalse{\fi\let\\=\cr\iffalse}\fi}
\Invalid@\vspace
\def\vspace@{\def\vspace##1{\crcr\noalign{\vskip##1\relax}}}
\def\multilimits@{\bgroup\vspace@\Let@
 \baselineskip\fontdimen10 \scriptfont\tw@
 \advance\baselineskip\fontdimen12 \scriptfont\tw@
 \lineskip\thr@@\fontdimen8 \scriptfont\thr@@
 \lineskiplimit\lineskip
 \vbox\bgroup\ialign\bgroup\hfil$\m@th\scriptstyle{##}$\hfil\crcr}
\def\Sb{_\multilimits@}
\def\endSb{\crcr\egroup\egroup\egroup}
\def\Sp{^\multilimits@}

\def\spreadlines#1{\RIfMIfI@\onlydmatherr@\spreadlines\else
 \openup#1\relax\fi\else\onlydmatherr@\spreadlines\fi}
\def\Mathstrut@{\copy\Mathstrutbox@}
\newbox\Mathstrutbox@
\setbox\Mathstrutbox@\null
\setboxz@h{$\m@th($}
\ht\Mathstrutbox@\ht\z@
\dp\Mathstrutbox@\dp\z@
\message{matrices,}
\newdimen\spreadmlines@
\def\spreadmatrixlines#1{\RIfMIfI@
 \onlydmatherr@\spreadmatrixlines\else
 \spreadmlines@#1\relax\fi\else\onlydmatherr@\spreadmatrixlines\fi}
\def\matrix{\null\,\vcenter\bgroup\Let@\vspace@
 \normalbaselines\openup\spreadmlines@\ialign
 \bgroup\hfil$\m@th##$\hfil&&\quad\hfil$\m@th##$\hfil\crcr
 \Mathstrut@\crcr\noalign{\kern-\baselineskip}}
\def\endmatrix{\crcr\Mathstrut@\crcr\noalign{\kern-\baselineskip}\egroup
 \egroup\,}
\def\format{\crcr\egroup\iffalse{\fi\ifnum`}=0 \fi\format@}
\newtoks\hashtoks@
\hashtoks@{#}
\def\format@#1\\{\def\preamble@{#1}%
 \def\l{$\m@th\the\hashtoks@$\hfil}%
 \def\c{\hfil$\m@th\the\hashtoks@$\hfil}%
 \def\r{\hfil$\m@th\the\hashtoks@$}%
 \edef\preamble@@{\preamble@}\ifnum`{=0 \fi\iffalse}\fi
 \ialign\bgroup\span\preamble@@\crcr}
\def\smallmatrix{\null\,\vcenter\bgroup\vspace@\Let@
 \baselineskip9\ex@\lineskip\ex@
 \ialign\bgroup\hfil$\m@th\scriptstyle{##}$\hfil&&\thickspace\hfil
 $\m@th\scriptstyle{##}$\hfil\crcr}
\def\endsmallmatrix{\crcr\egroup\egroup\,}

\newmuskip\dotsspace@
\dotsspace@1.5mu
\def\strip@#1 {#1}
\def\spacehdots#1\for#2{\multispan{#2}\xleaders
 \hbox{$\m@th\mkern\strip@#1 \dotsspace@.\mkern\strip@#1 \dotsspace@$}\hfill}
\def\hdotsfor#1{\spacehdots\@ne\for{#1}}
\def\multispan@#1{\omit\mscount#1\unskip\loop\ifnum\mscount>\@ne\sp@n\repeat}
\def\spaceinnerhdots#1\for#2\after#3{\multispan@{\strip@#2 }#3\xleaders
 \hbox{$\m@th\mkern\strip@#1 \dotsspace@.\mkern\strip@#1 \dotsspace@$}\hfill}
\def\innerhdotsfor#1\after#2{\spaceinnerhdots\@ne\for#1\after{#2}}
\def\cases{\bgroup\spreadmlines@\jot\left\{\,\matrix\format\l&\quad\l\\}
\def\endcases{\endmatrix\right.\egroup}
\message{multiline displays,}
\newif\ifinany@
\newif\ifinalign@
\newif\ifingather@
\def\strut@{\copy\strutbox@}
\newbox\strutbox@
\setbox\strutbox@\hbox{\vrule height8\p@ depth3\p@ width\z@}
\def\topaligned{\null\,\vtop\aligned@}
\def\botaligned{\null\,\vbox\aligned@}
\def\aligned{\null\,\vcenter\aligned@}
\def\aligned@{\bgroup\vspace@\Let@
 \ifinany@\else\openup\jot\fi\ialign
 \bgroup\hfil\strut@$\m@th\displaystyle{##}$&
 $\m@th\displaystyle{{}##}$\hfil\crcr}
\def\endaligned{\crcr\egroup\egroup}

\def\alignedat#1{\null\,\vcenter\bgroup\doat@{#1}\vspace@\Let@
 \ifinany@\else\openup\jot\fi\ialign\bgroup\span\preamble@@\crcr}
\newcount\atcount@
\def\doat@#1{\toks@{\hfil\strut@$\m@th
 \displaystyle{\the\hashtoks@}$&$\m@th\displaystyle
 {{}\the\hashtoks@}$\hfil}
 \atcount@#1\relax\advance\atcount@\m@ne                                    
 \loop\ifnum\atcount@>\z@\toks@=\expandafter{\the\toks@&\hfil$\m@th
 \displaystyle{\the\hashtoks@}$&$\m@th
 \displaystyle{{}\the\hashtoks@}$\hfil}\advance
  \atcount@\m@ne\repeat                                                     
 \xdef\preamble@{\the\toks@}\xdef\preamble@@{\preamble@}}

\def\gathered{\null\,\vcenter\bgroup\vspace@\Let@
 \ifinany@\else\openup\jot\fi\ialign
 \bgroup\hfil\strut@$\m@th\displaystyle{##}$\hfil\crcr}
\def\endgathered{\crcr\egroup\egroup}
\newif\iftagsleft@
\def\TagsOnLeft{\global\tagsleft@true}
\def\TagsOnRight{\global\tagsleft@false}
\TagsOnLeft
\newif\ifmathtags@
\def\TagsAsMath{\global\mathtags@true}
\def\TagsAsText{\global\mathtags@false}
\TagsAsText
\def\tagform@#1{\hbox{\rm(\ignorespaces#1\unskip)}}
\def\thetag{\leavevmode\tagform@}
\def\tag#1$${\iftagsleft@\leqno\else\eqno\fi                                
 \maketag@#1\maketag@                                                       
 $$}                                                                        
\def\maketag@{\FN@\maketag@@}
\def\maketag@@{\ifx\next"\expandafter\maketag@@@\else\expandafter\maketag@@@@
 \fi}
\def\maketag@@@"#1"#2\maketag@{\hbox{\rm#1}}                                
\def\maketag@@@@#1\maketag@{\ifmathtags@\tagform@{$\m@th#1$}\else
 \tagform@{#1}\fi}
\interdisplaylinepenalty\@M
\def\allowdisplaybreaks{\RIfMIfI@
 \onlydmatherr@\allowdisplaybreaks\else
 \interdisplaylinepenalty\z@\fi\else\onlydmatherr@\allowdisplaybreaks\fi}
\Invalid@\allowdisplaybreak
\Invalid@\displaybreak
\Invalid@\intertext
\def\allowdisplaybreak@{\def\allowdisplaybreak{\crcr\noalign{\allowbreak}}}
\def\displaybreak@{\def\displaybreak{\crcr\noalign{\break}}}
\def\intertext@{\def\intertext##1{\crcr\noalign{%
 \penalty\postdisplaypenalty \vskip\belowdisplayskip
 \vbox{\normalbaselines\noindent##1}%
 \penalty\predisplaypenalty \vskip\abovedisplayskip}}}
\newskip\centering@
\centering@\z@ plus\@m\p@
\def\align{\relax\ifingather@\DN@{\csname align (in
  \string\gather)\endcsname}\else
 \ifmmode\ifinner\DN@{\onlydmatherr@\align}\else
  \let\next@\align@\fi
 \else\DN@{\onlydmatherr@\align}\fi\fi\next@}
\newhelp\andhelp@
{An extra & here is so disastrous that you should probably exit^^J
and fix things up.}
\newif\iftag@
\newcount\and@
\def\align@{\inalign@true\inany@true
 \vspace@\allowdisplaybreak@\displaybreak@\intertext@
 \def\tag{\global\tag@true\ifnum\and@=\z@\DN@{&&}\else
          \DN@{&}\fi\next@}%
 \iftagsleft@\DN@{\csname align \endcsname}\else
  \DN@{\csname align \space\endcsname}\fi\next@}
\def\Tag@{\iftag@\else\errhelp\andhelp@\err@{Extra & on this line}\fi}
\newdimen\lwidth@
\newdimen\rwidth@
\newdimen\maxlwidth@
\newdimen\maxrwidth@
\newdimen\totwidth@
\def\measure@#1\endalign{\lwidth@\z@\rwidth@\z@\maxlwidth@\z@\maxrwidth@\z@
 \global\and@\z@                                                            
 \setbox@ne\vbox                                                            
  {\everycr{\noalign{\global\tag@false\global\and@\z@}}\Let@                
  \halign{\setboxz@h{$\m@th\displaystyle{\@lign##}$}
   \global\lwidth@\wdz@                                                     
   \ifdim\lwidth@>\maxlwidth@\global\maxlwidth@\lwidth@\fi                  
   \global\advance\and@\@ne                                                 
   &\setboxz@h{$\m@th\displaystyle{{}\@lign##}$}\global\rwidth@\wdz@        
   \ifdim\rwidth@>\maxrwidth@\global\maxrwidth@\rwidth@\fi                  
   \global\advance\and@\@ne                                                
   &\Tag@
   \eat@{##}\crcr#1\crcr}}
 \totwidth@\maxlwidth@\advance\totwidth@\maxrwidth@}                       
\def\displ@y@{\global\dt@ptrue\openup\jot
 \everycr{\noalign{\global\tag@false\global\and@\z@\ifdt@p\global\dt@pfalse
 \vskip-\lineskiplimit\vskip\normallineskiplimit\else
 \penalty\interdisplaylinepenalty\fi}}}
\def\black@#1{\noalign{\ifdim#1>\displaywidth
 \dimen@\prevdepth\nointerlineskip                                          
 \vskip-\ht\strutbox@\vskip-\dp\strutbox@                                   
 \vbox{\noindent\hbox to#1{\strut@\hfill}}
 \prevdepth\dimen@                                                          
 \fi}}
\expandafter\def\csname align \space\endcsname#1\endalign
 {\measure@#1\endalign\global\and@\z@                                       
 \ifingather@\everycr{\noalign{\global\and@\z@}}\else\displ@y@\fi           
 \Let@\tabskip\centering@                                                   
 \halign to\displaywidth
  {\hfil\strut@\setboxz@h{$\m@th\displaystyle{\@lign##}$}
  \global\lwidth@\wdz@\boxz@\global\advance\and@\@ne                        
  \tabskip\z@skip                                                           
  &\setboxz@h{$\m@th\displaystyle{{}\@lign##}$}
  \global\rwidth@\wdz@\boxz@\hfill\global\advance\and@\@ne                  
  \tabskip\centering@                                                       
  &\setboxz@h{\@lign\strut@\maketag@##\maketag@}
  \dimen@\displaywidth\advance\dimen@-\totwidth@
  \divide\dimen@\tw@\advance\dimen@\maxrwidth@\advance\dimen@-\rwidth@     
  \ifdim\dimen@<\tw@\wdz@\llap{\vtop{\normalbaselines\null\boxz@}}
  \else\llap{\boxz@}\fi                                                    
  \tabskip\z@skip                                                          
  \crcr#1\crcr                                                             
  \black@\totwidth@}}                                                      
\newdimen\lineht@
\expandafter\def\csname align \endcsname#1\endalign{\measure@#1\endalign
 \global\and@\z@
 \ifdim\totwidth@>\displaywidth\let\displaywidth@\totwidth@\else
  \let\displaywidth@\displaywidth\fi                                        
 \ifingather@\everycr{\noalign{\global\and@\z@}}\else\displ@y@\fi
 \Let@\tabskip\centering@\halign to\displaywidth
  {\hfil\strut@\setboxz@h{$\m@th\displaystyle{\@lign##}$}%
  \global\lwidth@\wdz@\global\lineht@\ht\z@                                 
  \boxz@\global\advance\and@\@ne
  \tabskip\z@skip&\setboxz@h{$\m@th\displaystyle{{}\@lign##}$}%
  \global\rwidth@\wdz@\ifdim\ht\z@>\lineht@\global\lineht@\ht\z@\fi         
  \boxz@\hfil\global\advance\and@\@ne
  \tabskip\centering@&\kern-\displaywidth@                                  
  \setboxz@h{\@lign\strut@\maketag@##\maketag@}%
  \dimen@\displaywidth\advance\dimen@-\totwidth@
  \divide\dimen@\tw@\advance\dimen@\maxlwidth@\advance\dimen@-\lwidth@
  \ifdim\dimen@<\tw@\wdz@
   \rlap{\vbox{\normalbaselines\boxz@\vbox to\lineht@{}}}\else
   \rlap{\boxz@}\fi
  \tabskip\displaywidth@\crcr#1\crcr\black@\totwidth@}}
\expandafter\def\csname align (in \string\gather)\endcsname
  #1\endalign{\vcenter{\align@#1\endalign}}
\Invalid@\endalign
\newif\ifxat@
\def\alignat{\RIfMIfI@\DN@{\onlydmatherr@\alignat}\else
 \DN@{\csname alignat \endcsname}\fi\else
 \DN@{\onlydmatherr@\alignat}\fi\next@}
\newif\ifmeasuring@
\newbox\savealignat@
\expandafter\def\csname alignat \endcsname#1#2\endalignat                   
 {\inany@true\xat@false
 \def\tag{\global\tag@true\count@#1\relax\multiply\count@\tw@
  \xdef\tag@{}\loop\ifnum\count@>\and@\xdef\tag@{&\tag@}\advance\count@\m@ne
  \repeat\tag@}%
 \vspace@\allowdisplaybreak@\displaybreak@\intertext@
 \displ@y@\measuring@true                                                   
 \setbox\savealignat@\hbox{$\m@th\displaystyle\Let@
  \attag@{#1}
  \vbox{\halign{\span\preamble@@\crcr#2\crcr}}$}%
 \measuring@false                                                           
 \Let@\attag@{#1}
 \tabskip\centering@\halign to\displaywidth
  {\span\preamble@@\crcr#2\crcr                                             
  \black@{\wd\savealignat@}}}                                               
\Invalid@\endalignat
\def\xalignat{\RIfMIfI@
 \DN@{\onlydmatherr@\xalignat}\else
 \DN@{\csname xalignat \endcsname}\fi\else
 \DN@{\onlydmatherr@\xalignat}\fi\next@}
\expandafter\def\csname xalignat \endcsname#1#2\endxalignat
 {\inany@true\xat@true
 \def\tag{\global\tag@true\def\tag@{}\count@#1\relax\multiply\count@\tw@
  \loop\ifnum\count@>\and@\xdef\tag@{&\tag@}\advance\count@\m@ne\repeat\tag@}%
 \vspace@\allowdisplaybreak@\displaybreak@\intertext@
 \displ@y@\measuring@true\setbox\savealignat@\hbox{$\m@th\displaystyle\Let@
 \attag@{#1}\vbox{\halign{\span\preamble@@\crcr#2\crcr}}$}%
 \measuring@false\Let@
 \attag@{#1}\tabskip\centering@\halign to\displaywidth
 {\span\preamble@@\crcr#2\crcr\black@{\wd\savealignat@}}}
\def\attag@#1{\let\Maketag@\maketag@\let\TAG@\Tag@                          
 \let\Tag@=0\let\maketag@=0
 \ifmeasuring@\def\llap@##1{\setboxz@h{##1}\hbox to\tw@\wdz@{}}%
  \def\rlap@##1{\setboxz@h{##1}\hbox to\tw@\wdz@{}}\else
  \let\llap@\llap\let\rlap@\rlap\fi                                         
 \toks@{\hfil\strut@$\m@th\displaystyle{\@lign\the\hashtoks@}$\tabskip\z@skip
  \global\advance\and@\@ne&$\m@th\displaystyle{{}\@lign\the\hashtoks@}$\hfil
  \ifxat@\tabskip\centering@\fi\global\advance\and@\@ne}
 \iftagsleft@
  \toks@@{\tabskip\centering@&\Tag@\kern-\displaywidth
   \rlap@{\@lign\maketag@\the\hashtoks@\maketag@}%
   \global\advance\and@\@ne\tabskip\displaywidth}\else
  \toks@@{\tabskip\centering@&\Tag@\llap@{\@lign\maketag@
   \the\hashtoks@\maketag@}\global\advance\and@\@ne\tabskip\z@skip}\fi      
 \atcount@#1\relax\advance\atcount@\m@ne
 \loop\ifnum\atcount@>\z@
 \toks@=\expandafter{\the\toks@&\hfil$\m@th\displaystyle{\@lign
  \the\hashtoks@}$\global\advance\and@\@ne
  \tabskip\z@skip&$\m@th\displaystyle{{}\@lign\the\hashtoks@}$\hfil\ifxat@
  \tabskip\centering@\fi\global\advance\and@\@ne}\advance\atcount@\m@ne
 \repeat                                                                    
 \xdef\preamble@{\the\toks@\the\toks@@}
 \xdef\preamble@@{\preamble@}
 \let\maketag@\Maketag@\let\Tag@\TAG@}                                      
\Invalid@\endxalignat
\def\xxalignat{\RIfMIfI@
 \DN@{\onlydmatherr@\xxalignat}\else\DN@{\csname xxalignat
  \endcsname}\fi\else
 \DN@{\onlydmatherr@\xxalignat}\fi\next@}
\expandafter\def\csname xxalignat \endcsname#1#2\endxxalignat{\inany@true
 \vspace@\allowdisplaybreak@\displaybreak@\intertext@
 \displ@y\setbox\savealignat@\hbox{$\m@th\displaystyle\Let@
 \xxattag@{#1}\vbox{\halign{\span\preamble@@\crcr#2\crcr}}$}%
 \Let@\xxattag@{#1}\tabskip\z@skip\halign to\displaywidth
 {\span\preamble@@\crcr#2\crcr\black@{\wd\savealignat@}}}
\def\xxattag@#1{\toks@{\tabskip\z@skip\hfil\strut@
 $\m@th\displaystyle{\the\hashtoks@}$&%
 $\m@th\displaystyle{{}\the\hashtoks@}$\hfil\tabskip\centering@&}%
 \atcount@#1\relax\advance\atcount@\m@ne\loop\ifnum\atcount@>\z@
 \toks@=\expandafter{\the\toks@&\hfil$\m@th\displaystyle{\the\hashtoks@}$%
  \tabskip\z@skip&$\m@th\displaystyle{{}\the\hashtoks@}$\hfil
  \tabskip\centering@}\advance\atcount@\m@ne\repeat
 \xdef\preamble@{\the\toks@\tabskip\z@skip}\xdef\preamble@@{\preamble@}}
\Invalid@\endxxalignat
\newdimen\gwidth@
\newdimen\gmaxwidth@
\def\gmeasure@#1\endgather{\gwidth@\z@\gmaxwidth@\z@\setbox@ne\vbox{\Let@
 \halign{\setboxz@h{$\m@th\displaystyle{##}$}\global\gwidth@\wdz@
 \ifdim\gwidth@>\gmaxwidth@\global\gmaxwidth@\gwidth@\fi
 &\eat@{##}\crcr#1\crcr}}}
\def\gather{\RIfMIfI@\DN@{\onlydmatherr@\gather}\else
 \ingather@true\inany@true\def\tag{&}%
 \vspace@\allowdisplaybreak@\displaybreak@\intertext@
 \displ@y\Let@
 \iftagsleft@\DN@{\csname gather \endcsname}\else
  \DN@{\csname gather \space\endcsname}\fi\fi
 \else\DN@{\onlydmatherr@\gather}\fi\next@}
\expandafter\def\csname gather \space\endcsname#1\endgather
 {\gmeasure@#1\endgather\tabskip\centering@
 \halign to\displaywidth{\hfil\strut@\setboxz@h{$\m@th\displaystyle{##}$}%
 \global\gwidth@\wdz@\boxz@\hfil&
 \setboxz@h{\strut@{\maketag@##\maketag@}}%
 \dimen@\displaywidth\advance\dimen@-\gwidth@
 \ifdim\dimen@>\tw@\wdz@\llap{\boxz@}\else
 \llap{\vtop{\normalbaselines\null\boxz@}}\fi
 \tabskip\z@skip\crcr#1\crcr\black@\gmaxwidth@}}
\newdimen\glineht@
\expandafter\def\csname gather \endcsname#1\endgather{\gmeasure@#1\endgather
 \ifdim\gmaxwidth@>\displaywidth\let\gdisplaywidth@\gmaxwidth@\else
 \let\gdisplaywidth@\displaywidth\fi\tabskip\centering@\halign to\displaywidth
 {\hfil\strut@\setboxz@h{$\m@th\displaystyle{##}$}%
 \global\gwidth@\wdz@\global\glineht@\ht\z@\boxz@\hfil&\kern-\gdisplaywidth@
 \setboxz@h{\strut@{\maketag@##\maketag@}}%
 \dimen@\displaywidth\advance\dimen@-\gwidth@
 \ifdim\dimen@>\tw@\wdz@\rlap{\boxz@}\else
 \rlap{\vbox{\normalbaselines\boxz@\vbox to\glineht@{}}}\fi
 \tabskip\gdisplaywidth@\crcr#1\crcr\black@\gmaxwidth@}}
\newif\ifctagsplit@
\def\CenteredTagsOnSplits{\global\ctagsplit@true}
\def\TopOrBottomTagsOnSplits{\global\ctagsplit@false}
\TopOrBottomTagsOnSplits
\def\split{\relax\ifinany@\let\next@\insplit@\else
 \ifmmode\ifinner\def\next@{\onlydmatherr@\split}\else
 \let\next@\outsplit@\fi\else
 \def\next@{\onlydmatherr@\split}\fi\fi\next@}
\def\insplit@{\global\setbox\z@\vbox\bgroup\vspace@\Let@\ialign\bgroup
 \hfil\strut@$\m@th\displaystyle{##}$&$\m@th\displaystyle{{}##}$\hfill\crcr}
\def\endsplit{\crcr\egroup\egroup\iftagsleft@\expandafter\lendsplit@\else
 \expandafter\rendsplit@\fi}
\def\rendsplit@{\global\setbox9 \vbox
 {\unvcopy\z@\global\setbox8 \lastbox\unskip}
 \setbox@ne\hbox{\unhcopy8 \unskip\global\setbox\tw@\lastbox
 \unskip\global\setbox\thr@@\lastbox}
 \global\setbox7 \hbox{\unhbox\tw@\unskip}
 \ifinalign@\ifctagsplit@                                                   
  \gdef\split@{\hbox to\wd\thr@@{}&
   \vcenter{\vbox{\moveleft\wd\thr@@\boxz@}}}
 \else\gdef\split@{&\vbox{\moveleft\wd\thr@@\box9}\crcr
  \box\thr@@&\box7}\fi                                                      
 \else                                                                      
  \ifctagsplit@\gdef\split@{\vcenter{\boxz@}}\else
  \gdef\split@{\box9\crcr\hbox{\box\thr@@\box7}}\fi
 \fi
 \split@}                                                                   
\def\lendsplit@{\global\setbox9\vtop{\unvcopy\z@}
 \setbox@ne\vbox{\unvcopy\z@\global\setbox8\lastbox}
 \setbox@ne\hbox{\unhcopy8\unskip\setbox\tw@\lastbox
  \unskip\global\setbox\thr@@\lastbox}
 \ifinalign@\ifctagsplit@                                                   
  \gdef\split@{\hbox to\wd\thr@@{}&
  \vcenter{\vbox{\moveleft\wd\thr@@\box9}}}
  \else                                                                     
  \gdef\split@{\hbox to\wd\thr@@{}&\vbox{\moveleft\wd\thr@@\box9}}\fi
 \else
  \ifctagsplit@\gdef\split@{\vcenter{\box9}}\else
  \gdef\split@{\box9}\fi
 \fi\split@}
\def\outsplit@#1$${\align\insplit@#1\endalign$$}
\newdimen\multlinegap@
\multlinegap@1em
\newdimen\multlinetaggap@
\multlinetaggap@1em
\def\MultlineGap#1{\global\multlinegap@#1\relax}
\def\multlinegap#1{\RIfMIfI@\onlydmatherr@\multlinegap\else
 \multlinegap@#1\relax\fi\else\onlydmatherr@\multlinegap\fi}
\def\nomultlinegap{\multlinegap{\z@}}
\def\multline{\RIfMIfI@
 \DN@{\onlydmatherr@\multline}\else
 \DN@{\multline@}\fi\else
 \DN@{\onlydmatherr@\multline}\fi\next@}
\newif\iftagin@
\def\tagin@#1{\tagin@false\in@\tag{#1}\ifin@\tagin@true\fi}
\def\multline@#1$${\inany@true\vspace@\allowdisplaybreak@\displaybreak@
 \tagin@{#1}\iftagsleft@\DN@{\multline@l#1$$}\else
 \DN@{\multline@r#1$$}\fi\next@}
\newdimen\mwidth@
\def\rmmeasure@#1\endmultline{%
 \def\shoveleft##1{##1}\def\shoveright##1{##1}
 \setbox@ne\vbox{\Let@\halign{\setboxz@h
  {$\m@th\@lign\displaystyle{}##$}\global\mwidth@\wdz@
  \crcr#1\crcr}}}
\newdimen\mlineht@
\newif\ifzerocr@
\newif\ifonecr@
\def\lmmeasure@#1\endmultline{\global\zerocr@true\global\onecr@false
 \everycr{\noalign{\ifonecr@\global\onecr@false\fi
  \ifzerocr@\global\zerocr@false\global\onecr@true\fi}}
  \def\shoveleft##1{##1}\def\shoveright##1{##1}%
 \setbox@ne\vbox{\Let@\halign{\setboxz@h
  {$\m@th\@lign\displaystyle{}##$}\ifonecr@\global\mwidth@\wdz@
  \global\mlineht@\ht\z@\fi\crcr#1\crcr}}}
\newbox\mtagbox@
\newdimen\ltwidth@
\newdimen\rtwidth@
\def\multline@l#1$${\iftagin@\DN@{\lmultline@@#1$$}\else
 \DN@{\setbox\mtagbox@\null\ltwidth@\z@\rtwidth@\z@
  \lmultline@@@#1$$}\fi\next@}
\def\lmultline@@#1\endmultline\tag#2$${%
 \setbox\mtagbox@\hbox{\maketag@#2\maketag@}
 \lmmeasure@#1\endmultline\dimen@\mwidth@\advance\dimen@\wd\mtagbox@
 \advance\dimen@\multlinetaggap@                                            
 \ifdim\dimen@>\displaywidth\ltwidth@\z@\else\ltwidth@\wd\mtagbox@\fi       
 \lmultline@@@#1\endmultline$$}
\def\lmultline@@@{\displ@y
 \def\shoveright##1{##1\hfilneg\hskip\multlinegap@}%
 \def\shoveleft##1{\setboxz@h{$\m@th\displaystyle{}##1$}%
  \setbox@ne\hbox{$\m@th\displaystyle##1$}%
  \hfilneg
  \iftagin@
   \ifdim\ltwidth@>\z@\hskip\ltwidth@\hskip\multlinetaggap@\fi
  \else\hskip\multlinegap@\fi\hskip.5\wd@ne\hskip-.5\wdz@##1}
  \halign\bgroup\Let@\hbox to\displaywidth
   {\strut@$\m@th\displaystyle\hfil{}##\hfil$}\crcr
   \hfilneg                                                                 
   \iftagin@                                                                
    \ifdim\ltwidth@>\z@                                                     
     \box\mtagbox@\hskip\multlinetaggap@                                    
    \else
     \rlap{\vbox{\normalbaselines\hbox{\strut@\box\mtagbox@}%
     \vbox to\mlineht@{}}}\fi                                               
   \else\hskip\multlinegap@\fi}                                             
\def\multline@r#1$${\iftagin@\DN@{\rmultline@@#1$$}\else
 \DN@{\setbox\mtagbox@\null\ltwidth@\z@\rtwidth@\z@
  \rmultline@@@#1$$}\fi\next@}
\def\rmultline@@#1\endmultline\tag#2$${\ltwidth@\z@
 \setbox\mtagbox@\hbox{\maketag@#2\maketag@}%
 \rmmeasure@#1\endmultline\dimen@\mwidth@\advance\dimen@\wd\mtagbox@
 \advance\dimen@\multlinetaggap@
 \ifdim\dimen@>\displaywidth\rtwidth@\z@\else\rtwidth@\wd\mtagbox@\fi
 \rmultline@@@#1\endmultline$$}
\def\rmultline@@@{\displ@y
 \def\shoveright##1{##1\hfilneg\iftagin@\ifdim\rtwidth@>\z@
  \hskip\rtwidth@\hskip\multlinetaggap@\fi\else\hskip\multlinegap@\fi}%
 \def\shoveleft##1{\setboxz@h{$\m@th\displaystyle{}##1$}%
  \setbox@ne\hbox{$\m@th\displaystyle##1$}%
  \hfilneg\hskip\multlinegap@\hskip.5\wd@ne\hskip-.5\wdz@##1}%
 \halign\bgroup\Let@\hbox to\displaywidth
  {\strut@$\m@th\displaystyle\hfil{}##\hfil$}\crcr
 \hfilneg\hskip\multlinegap@}
\def\endmultline{\iftagsleft@\expandafter\lendmultline@\else
 \expandafter\rendmultline@\fi}
\def\lendmultline@{\hfilneg\hskip\multlinegap@\crcr\egroup}
\def\rendmultline@{\iftagin@                                                
 \ifdim\rtwidth@>\z@                                                        
  \hskip\multlinetaggap@\box\mtagbox@                                       
 \else\llap{\vtop{\normalbaselines\null\hbox{\strut@\box\mtagbox@}}}\fi     
 \else\hskip\multlinegap@\fi                                                
 \hfilneg\crcr\egroup}
\def\bmod{\mskip-\medmuskip\mkern5mu\mathbin{\fam\z@ mod}\penalty900
 \mkern5mu\mskip-\medmuskip}
\def\pmod#1{\allowbreak\ifinner\mkern8mu\else\mkern18mu\fi
 ({\fam\z@ mod}\,\,#1)}
\def\pod#1{\allowbreak\ifinner\mkern8mu\else\mkern18mu\fi(#1)}
\def\mod#1{\allowbreak\ifinner\mkern12mu\else\mkern18mu\fi{\fam\z@ mod}\,\,#1}
\message{continued fractions,}
\newcount\cfraccount@
\def\cfrac{\bgroup\bgroup\advance\cfraccount@\@ne\strut
 \iffalse{\fi\def\\{\over\displaystyle}\iffalse}\fi}
\def\lcfrac{\bgroup\bgroup\advance\cfraccount@\@ne\strut
 \iffalse{\fi\def\\{\hfill\over\displaystyle}\iffalse}\fi}
\def\rcfrac{\bgroup\bgroup\advance\cfraccount@\@ne\strut\hfill
 \iffalse{\fi\def\\{\over\displaystyle}\iffalse}\fi}
\def\gloop@#1\repeat{\gdef\body{#1}\iterate}
\def\endcfrac{\gloop@\ifnum\cfraccount@>\z@\global\advance\cfraccount@\m@ne
 \egroup\hskip-\nulldelimiterspace\egroup\repeat}
\message{compound symbols,}
\def\binrel@#1{\setboxz@h{\thinmuskip0mu
  \medmuskip\m@ne mu\thickmuskip\@ne mu$#1\m@th$}%
 \setbox@ne\hbox{\thinmuskip0mu\medmuskip\m@ne mu\thickmuskip
  \@ne mu${}#1{}\m@th$}%
 \setbox\tw@\hbox{\hskip\wd@ne\hskip-\wdz@}}
\def\overset#1\to#2{\binrel@{#2}\ifdim\wd\tw@<\z@
 \mathbin{\mathop{\kern\z@#2}\limits^{#1}}\else\ifdim\wd\tw@>\z@
 \mathrel{\mathop{\kern\z@#2}\limits^{#1}}\else
 {\mathop{\kern\z@#2}\limits^{#1}}{}\fi\fi}
\def\underset#1\to#2{\binrel@{#2}\ifdim\wd\tw@<\z@
 \mathbin{\mathop{\kern\z@#2}\limits_{#1}}\else\ifdim\wd\tw@>\z@
 \mathrel{\mathop{\kern\z@#2}\limits_{#1}}\else
 {\mathop{\kern\z@#2}\limits_{#1}}{}\fi\fi}
\def\oversetbrace#1\to#2{\overbrace{#2}^{#1}}
\def\undersetbrace#1\to#2{\underbrace{#2}_{#1}}
\def\sideset#1\and#2\to#3{%
 \setbox@ne\hbox{$\dsize{\vphantom{#3}}#1{#3}\m@th$}%
 \setbox\tw@\hbox{$\dsize{#3}#2\m@th$}%
 \hskip\wd@ne\hskip-\wd\tw@\mathop{\hskip\wd\tw@\hskip-\wd@ne
  {\vphantom{#3}}#1{#3}#2}}
\def\rightarrowfill@#1{\setboxz@h{$#1-\m@th$}\ht\z@\z@
  $#1\m@th\copy\z@\mkern-6mu\cleaders
  \hbox{$#1\mkern-2mu\box\z@\mkern-2mu$}\hfill
  \mkern-6mu\mathord\rightarrow$}
\def\leftarrowfill@#1{\setboxz@h{$#1-\m@th$}\ht\z@\z@
  $#1\m@th\mathord\leftarrow\mkern-6mu\cleaders
  \hbox{$#1\mkern-2mu\copy\z@\mkern-2mu$}\hfill
  \mkern-6mu\box\z@$}
\def\leftrightarrowfill@#1{\setboxz@h{$#1-\m@th$}\ht\z@\z@
  $#1\m@th\mathord\leftarrow\mkern-6mu\cleaders
  \hbox{$#1\mkern-2mu\box\z@\mkern-2mu$}\hfill
  \mkern-6mu\mathord\rightarrow$}
\def\overrightarrow{\mathpalette\overrightarrow@}
\def\overrightarrow@#1#2{\vbox{\ialign{##\crcr\rightarrowfill@#1\crcr
 \noalign{\kern-\ex@\nointerlineskip}$\m@th\hfil#1#2\hfil$\crcr}}}

\def\overleftarrow{\mathpalette\overleftarrow@}
\def\overleftarrow@#1#2{\vbox{\ialign{##\crcr\leftarrowfill@#1\crcr
 \noalign{\kern-\ex@\nointerlineskip}$\m@th\hfil#1#2\hfil$\crcr}}}
\def\overleftrightarrow{\mathpalette\overleftrightarrow@}
\def\overleftrightarrow@#1#2{\vbox{\ialign{##\crcr\leftrightarrowfill@#1\crcr
 \noalign{\kern-\ex@\nointerlineskip}$\m@th\hfil#1#2\hfil$\crcr}}}
\def\underrightarrow{\mathpalette\underrightarrow@}
\def\underrightarrow@#1#2{\vtop{\ialign{##\crcr$\m@th\hfil#1#2\hfil$\crcr
 \noalign{\nointerlineskip}\rightarrowfill@#1\crcr}}}

\def\underleftarrow{\mathpalette\underleftarrow@}
\def\underleftarrow@#1#2{\vtop{\ialign{##\crcr$\m@th\hfil#1#2\hfil$\crcr
 \noalign{\nointerlineskip}\leftarrowfill@#1\crcr}}}
\def\underleftrightarrow{\mathpalette\underleftrightarrow@}
\def\underleftrightarrow@#1#2{\vtop{\ialign{##\crcr$\m@th\hfil#1#2\hfil$\crcr
 \noalign{\nointerlineskip}\leftrightarrowfill@#1\crcr}}}
\message{various kinds of dots,}
\let\DOTSI\relax
\let\DOTSB\relax

\newif\ifmath@
{\uccode`7=`\\ \uccode`8=`m \uccode`9=`a \uccode`0=`t \uccode`!=`h
 \uppercase{\gdef\math@#1#2#3#4#5#6\math@{\global\math@false\ifx 7#1\ifx 8#2%
 \ifx 9#3\ifx 0#4\ifx !#5\xdef\meaning@{#6}\global\math@true\fi\fi\fi\fi\fi}}}
\newif\ifmathch@
{\uccode`7=`c \uccode`8=`h \uccode`9=`\"
 \uppercase{\gdef\mathch@#1#2#3#4#5#6\mathch@{\global\mathch@false
  \ifx 7#1\ifx 8#2\ifx 9#5\global\mathch@true\xdef\meaning@{9#6}\fi\fi\fi}}}
\newcount\classnum@
\def\getmathch@#1.#2\getmathch@{\classnum@#1 \divide\classnum@4096
 \ifcase\number\classnum@\or\or\gdef\thedots@{\dotsb@}\or
 \gdef\thedots@{\dotsb@}\fi}
\newif\ifmathbin@
{\uccode`4=`b \uccode`5=`i \uccode`6=`n
 \uppercase{\gdef\mathbin@#1#2#3{\relaxnext@
  \DNii@##1\mathbin@{\ifx\space@\next\global\mathbin@true\fi}%
 \global\mathbin@false\DN@##1\mathbin@{}%
 \ifx 4#1\ifx 5#2\ifx 6#3\DN@{\FN@\nextii@}\fi\fi\fi\next@}}}
\newif\ifmathrel@
{\uccode`4=`r \uccode`5=`e \uccode`6=`l
 \uppercase{\gdef\mathrel@#1#2#3{\relaxnext@
  \DNii@##1\mathrel@{\ifx\space@\next\global\mathrel@true\fi}%
 \global\mathrel@false\DN@##1\mathrel@{}%
 \ifx 4#1\ifx 5#2\ifx 6#3\DN@{\FN@\nextii@}\fi\fi\fi\next@}}}
\newif\ifmacro@
{\uccode`5=`m \uccode`6=`a \uccode`7=`c
 \uppercase{\gdef\macro@#1#2#3#4\macro@{\global\macro@false
  \ifx 5#1\ifx 6#2\ifx 7#3\global\macro@true
  \xdef\meaning@{\macro@@#4\macro@@}\fi\fi\fi}}}
\def\macro@@#1->#2\macro@@{#2}
\newif\ifDOTS@
\newcount\DOTSCASE@
{\uccode`6=`\\ \uccode`7=`D \uccode`8=`O \uccode`9=`T \uccode`0=`S
 \uppercase{\gdef\DOTS@#1#2#3#4#5{\global\DOTS@false\DN@##1\DOTS@{}%
  \ifx 6#1\ifx 7#2\ifx 8#3\ifx 9#4\ifx 0#5\let\next@\DOTS@@\fi\fi\fi\fi\fi
  \next@}}}
{\uccode`3=`B \uccode`4=`I \uccode`5=`X
 \uppercase{\gdef\DOTS@@#1{\relaxnext@
  \DNii@##1\DOTS@{\ifx\space@\next\global\DOTS@true\fi}%
  \DN@{\FN@\nextii@}%
  \ifx 3#1\global\DOTSCASE@\z@\else
  \ifx 4#1\global\DOTSCASE@\@ne\else
  \ifx 5#1\global\DOTSCASE@\tw@\else\DN@##1\DOTS@{}%
  \fi\fi\fi\next@}}}
\newif\ifnot@
{\uccode`5=`\\ \uccode`6=`n \uccode`7=`o \uccode`8=`t
 \uppercase{\gdef\not@#1#2#3#4{\relaxnext@
  \DNii@##1\not@{\ifx\space@\next\global\not@true\fi}%
 \global\not@false\DN@##1\not@{}%
 \ifx 5#1\ifx 6#2\ifx 7#3\ifx 8#4\DN@{\FN@\nextii@}\fi\fi\fi
 \fi\next@}}}
\newif\ifkeybin@
\def\keybin@{\keybin@true
 \ifx\next+\else\ifx\next=\else\ifx\next<\else\ifx\next>\else\ifx\next-\else
 \ifx\next*\else\ifx\next:\else\keybin@false\fi\fi\fi\fi\fi\fi\fi}
\def\dots{\RIfM@\expandafter\mdots@\else\expandafter\tdots@\fi}
\def\tdots@{\unskip\relaxnext@
 \DN@{$\m@th\mathinner{\ldotp\ldotp\ldotp}\,
   \ifx\next,\,$\else\ifx\next.\,$\else\ifx\next;\,$\else\ifx\next:\,$\else
   \ifx\next?\,$\else\ifx\next!\,$\else$ \fi\fi\fi\fi\fi\fi}%
 \ \FN@\next@}
\def\mdots@{\FN@\mdots@@}
\def\mdots@@{\gdef\thedots@{\dotso@}
 \ifx\next\boldkey\gdef\thedots@\boldkey{\boldkeydots@}\else                
 \ifx\next\boldsymbol\gdef\thedots@\boldsymbol{\boldsymboldots@}\else       
 \ifx,\next\gdef\thedots@{\dotsc}
 \else\ifx\not\next\gdef\thedots@{\dotsb@}
 \else\keybin@
 \ifkeybin@\gdef\thedots@{\dotsb@}
 \else\xdef\meaning@{\meaning\next..........}\xdef\meaning@@{\meaning@}
  \expandafter\math@\meaning@\math@
  \ifmath@
   \expandafter\mathch@\meaning@\mathch@
   \ifmathch@\expandafter\getmathch@\meaning@\getmathch@\fi                 
  \else\expandafter\macro@\meaning@@\macro@                                 
  \ifmacro@                                                                
   \expandafter\not@\meaning@\not@\ifnot@\gdef\thedots@{\dotsb@}
  \else\expandafter\DOTS@\meaning@\DOTS@
  \ifDOTS@
   \ifcase\number\DOTSCASE@\gdef\thedots@{\dotsb@}%
    \or\gdef\thedots@{\dotsi}\else\fi                                      
  \else\expandafter\math@\meaning@\math@                                   
  \ifmath@\expandafter\mathbin@\meaning@\mathbin@
  \ifmathbin@\gdef\thedots@{\dotsb@}
  \else\expandafter\mathrel@\meaning@\mathrel@
  \ifmathrel@\gdef\thedots@{\dotsb@}
  \fi\fi\fi\fi\fi\fi\fi\fi\fi\fi\fi\fi
 \thedots@}
\def\plainldots@{\mathinner{\ldotp\ldotp\ldotp}}
\def\plaincdots@{\mathinner{\cdotp\cdotp\cdotp}}
\def\dotsi{\!\plaincdots@}
\let\dotsb@\plaincdots@
\newif\ifextra@
\newif\ifrightdelim@
\def\rightdelim@{\global\rightdelim@true                                    
 \ifx\next)\else                                                            
 \ifx\next]\else
 \ifx\next\rbrack\else
 \ifx\next\}\else
 \ifx\next\rbrace\else
 \ifx\next\rangle\else
 \ifx\next\rceil\else
 \ifx\next\rfloor\else
 \ifx\next\rgroup\else
 \ifx\next\rmoustache\else
 \ifx\next\right\else
 \ifx\next\bigr\else
 \ifx\next\biggr\else
 \ifx\next\Bigr\else                                                        
 \ifx\next\Biggr\else\global\rightdelim@false
 \fi\fi\fi\fi\fi\fi\fi\fi\fi\fi\fi\fi\fi\fi\fi}
\def\extra@{%
 \global\extra@false\rightdelim@\ifrightdelim@\global\extra@true            
 \else\ifx\next$\global\extra@true                                          
 \else\xdef\meaning@{\meaning\next..........}
 \expandafter\macro@\meaning@\macro@\ifmacro@                               
 \expandafter\DOTS@\meaning@\DOTS@
 \ifDOTS@
 \ifnum\DOTSCASE@=\tw@\global\extra@true                                    
 \fi\fi\fi\fi\fi}
\newif\ifbold@
\def\dotso@{\relaxnext@
 \ifbold@
  \let\next\delayed@
  \DNii@{\extra@\plainldots@\ifextra@\,\fi}%
 \else
  \DNii@{\DN@{\extra@\plainldots@\ifextra@\,\fi}\FN@\next@}%
 \fi
 \nextii@}
\def\extrap@#1{%
 \ifx\next,\DN@{#1\,}\else
 \ifx\next;\DN@{#1\,}\else
 \ifx\next.\DN@{#1\,}\else\extra@
 \ifextra@\DN@{#1\,}\else
 \let\next@#1\fi\fi\fi\fi\next@}
\def\ldots{\DN@{\extrap@\plainldots@}%
 \FN@\next@}
\def\cdots{\DN@{\extrap@\plaincdots@}%
 \FN@\next@}

\def\dotsc{\relaxnext@
 \DN@{\ifx\next;\plainldots@\,\else
  \ifx\next.\plainldots@\,\else\extra@\plainldots@
  \ifextra@\,\fi\fi\fi}%
 \FN@\next@}
\def\cdot{\mathchar"2201 }
\def\longrightarrow{\DOTSB\relbar\joinrel\rightarrow}

\def\mapsto{\DOTSB\mapstochar\rightarrow}

\message{special superscripts,}
\def\dddot#1{{\mathop{#1}\limits^{\vbox to-1.4\ex@{\kern-\tw@\ex@
 \hbox{\rm...}\vss}}}}
\def\ddddot#1{{\mathop{#1}\limits^{\vbox to-1.4\ex@{\kern-\tw@\ex@
 \hbox{\rm....}\vss}}}}
\def\sphat{^{\mathchoice{}{}%
 {\,\,\botsmash{\hbox{\lower4\ex@\hbox{$\m@th\widehat{\null}$}}}}%
 {\,\botsmash{\hbox{\lower3\ex@\hbox{$\m@th\hat{\null}$}}}}}}

\def\spacute{^{\!\botsmash{\hbox{\lower\@ne ex\hbox{\'{}}}}}}
\def\spgrave{^{\mathchoice{}{}{}{\!}%
 \botsmash{\hbox{\lower\@ne ex\hbox{\`{}}}}}}
\def\spdot{^{\hbox{\raise\ex@\hbox{\rm.}}}}
\def\spddot{^{\hbox{\raise\ex@\hbox{\rm..}}}}
\def\spdddot{^{\hbox{\raise\ex@\hbox{\rm...}}}}
\def\spddddot{^{\hbox{\raise\ex@\hbox{\rm....}}}}
\def\spbreve{^{\!\botsmash{\hbox{\lower4\ex@\hbox{\u{}}}}}}

\message{\string\text,}
\def\textonlyfont@#1#2{\def#1{\RIfM@
 \Err@{Use \string#1\space only in text}\else#2\fi}}
\textonlyfont@\rm\tenrm
\textonlyfont@\it\tenit
\textonlyfont@\sl\tensl
\textonlyfont@\bf\tenbf
\def\oldnos#1{\RIfM@{\mathcode`\,="013B \fam\@ne#1}\else
 \leavevmode\hbox{$\m@th\mathcode`\,="013B \fam\@ne#1$}\fi}
\def\text{\RIfM@\expandafter\text@\else\expandafter\text@@\fi}
\def\text@@#1{\leavevmode\hbox{#1}}
\def\mathhexbox@#1#2#3{\text{$\m@th\mathchar"#1#2#3$}}
\def\dag{{\mathhexbox@279}}
\def\ddag{{\mathhexbox@27A}}
\def\S{{\mathhexbox@278}}
\def\P{{\mathhexbox@27B}}
\newif\iffirstchoice@
\firstchoice@true
\def\text@#1{\mathchoice
 {\hbox{\everymath{\displaystyle}\def\textfonti{\the\textfont\@ne}%
  \def\textfontii{\the\textfont\tw@}\textdef@@ T#1}}
 {\hbox{\firstchoice@false
  \everymath{\textstyle}\def\textfonti{\the\textfont\@ne}%
  \def\textfontii{\the\textfont\tw@}\textdef@@ T#1}}
 {\hbox{\firstchoice@false
  \everymath{\scriptstyle}\def\textfonti{\the\scriptfont\@ne}%
  \def\textfontii{\the\scriptfont\tw@}\textdef@@ S\rm#1}}
 {\hbox{\firstchoice@false
  \everymath{\scriptscriptstyle}\def\textfonti
  {\the\scriptscriptfont\@ne}%
  \def\textfontii{\the\scriptscriptfont\tw@}\textdef@@ s\rm#1}}}
\def\textdef@@#1{\textdef@#1\rm\textdef@#1\bf\textdef@#1\sl\textdef@#1\it}
\def\rmfam{0}
\def\textdef@#1#2{%
 \DN@{\csname\expandafter\eat@\string#2fam\endcsname}%
 \if S#1\edef#2{\the\scriptfont\next@\relax}%
 \else\if s#1\edef#2{\the\scriptscriptfont\next@\relax}%
 \else\edef#2{\the\textfont\next@\relax}\fi\fi}
\scriptfont\itfam\tenit \scriptscriptfont\itfam\tenit
\scriptfont\slfam\tensl \scriptscriptfont\slfam\tensl
\newif\iftopfolded@
\newif\ifbotfolded@
\def\topfoldedtext{\topfolded@true\botfolded@false\foldedtext@}
\def\botfoldedtext{\botfolded@true\topfolded@false\foldedtext@}
\def\foldedtext{\topfolded@false\botfolded@false\foldedtext@}
\Invalid@\foldedwidth
\def\foldedtext@{\relaxnext@
 \DN@{\ifx\next\foldedwidth\let\next@\nextii@\else
  \DN@{\nextii@\foldedwidth{.3\hsize}}\fi\next@}%
 \DNii@\foldedwidth##1##2{\setbox\z@\vbox
  {\normalbaselines\hsize##1\relax
  \tolerance1600 \noindent\ignorespaces##2}\ifbotfolded@\boxz@\else
  \iftopfolded@\vtop{\unvbox\z@}\else\vcenter{\boxz@}\fi\fi}%
 \FN@\next@}
\message{math font commands,}
\def\bold{\RIfM@\expandafter\bold@\else
 \expandafter\nonmatherr@\expandafter\bold\fi}
\def\bold@#1{{\bold@@{#1}}}
\def\bold@@#1{\fam\bffam\relax#1}
\def\slanted{\RIfM@\expandafter\slanted@\else
 \expandafter\nonmatherr@\expandafter\slanted\fi}
\def\slanted@#1{{\slanted@@{#1}}}
\def\slanted@@#1{\fam\slfam\relax#1}
\def\roman{\RIfM@\expandafter\roman@\else
 \expandafter\nonmatherr@\expandafter\roman\fi}
\def\roman@#1{{\roman@@{#1}}}
\def\roman@@#1{\fam\rmfam\relax#1}
\def\italic{\RIfM@\expandafter\italic@\else
 \expandafter\nonmatherr@\expandafter\italic\fi}
\def\italic@#1{{\italic@@{#1}}}
\def\italic@@#1{\fam\itfam\relax#1}
\def\Cal{\RIfM@\expandafter\Cal@\else
 \expandafter\nonmatherr@\expandafter\Cal\fi}
\def\Cal@#1{{\Cal@@{#1}}}
\def\Cal@@#1{\noaccents@\fam\tw@#1}
\mathchardef\Gamma="0000
\mathchardef\Delta="0001
\mathchardef\Theta="0002
\mathchardef\Lambda="0003
\mathchardef\Xi="0004
\mathchardef\Pi="0005
\mathchardef\Sigma="0006
\mathchardef\Upsilon="0007
\mathchardef\Phi="0008
\mathchardef\Psi="0009
\mathchardef\Omega="000A
\mathchardef\varGamma="0100
\mathchardef\varDelta="0101
\mathchardef\varTheta="0102
\mathchardef\varLambda="0103
\mathchardef\varXi="0104
\mathchardef\varPi="0105
\mathchardef\varSigma="0106
\mathchardef\varUpsilon="0107
\mathchardef\varPhi="0108
\mathchardef\varPsi="0109
\mathchardef\varOmega="010A
\let\alloc@@\alloc@
\def\hexnumber@#1{\ifcase#1 0\or 1\or 2\or 3\or 4\or 5\or 6\or 7\or 8\or
 9\or A\or B\or C\or D\or E\or F\fi}
\def\loadmsam{%
 \font@\tenmsa=msam10
 \font@\sevenmsa=msam7
 \font@\fivemsa=msam5
 \alloc@@8\fam\chardef\sixt@@n\msafam
 \textfont\msafam=\tenmsa
 \scriptfont\msafam=\sevenmsa
 \scriptscriptfont\msafam=\fivemsa
 \edef\next{\hexnumber@\msafam}%
 \mathchardef\dabar@"0\next39
 \edef\dashrightarrow{\mathrel{\dabar@\dabar@\mathchar"0\next4B}}%
 \edef\dashleftarrow{\mathrel{\mathchar"0\next4C\dabar@\dabar@}}%
 \let\dasharrow\dashrightarrow
 \edef\ulcorner{\delimiter"4\next70\next70 }%
 \edef\urcorner{\delimiter"5\next71\next71 }%
 \edef\llcorner{\delimiter"4\next78\next78 }%
 \edef\lrcorner{\delimiter"5\next79\next79 }%
 \edef\yen{{\noexpand\mathhexbox@\next55}}%
 \edef\checkmark{{\noexpand\mathhexbox@\next58}}%
 \edef\circledR{{\noexpand\mathhexbox@\next72}}%
 \edef\maltese{{\noexpand\mathhexbox@\next7A}}%
 \global\let\loadmsam\empty}%
\def\loadmsbm{%
 \font@\tenmsb=msbm10 \font@\sevenmsb=msbm7 \font@\fivemsb=msbm5
 \alloc@@8\fam\chardef\sixt@@n\msbfam
 \textfont\msbfam=\tenmsb
 \scriptfont\msbfam=\sevenmsb \scriptscriptfont\msbfam=\fivemsb
 \global\let\loadmsbm\empty
 }
\def\widehat#1{\ifx\undefined\msbfam \DN@{362}%
  \else \setboxz@h{$\m@th#1$}%
    \edef\next@{\ifdim\wdz@>\tw@ em%
        \hexnumber@\msbfam 5B%
      \else 362\fi}\fi
  \mathaccent"0\next@{#1}}
\def\widetilde#1{\ifx\undefined\msbfam \DN@{365}%
  \else \setboxz@h{$\m@th#1$}%
    \edef\next@{\ifdim\wdz@>\tw@ em%
        \hexnumber@\msbfam 5D%
      \else 365\fi}\fi
  \mathaccent"0\next@{#1}}
\message{\string\newsymbol,}
\def\newsymbol#1#2#3#4#5{\define#1{}%
  \count@#2\relax \advance\count@\m@ne 
 \ifcase\count@
   \ifx\undefined\msafam\loadmsam\fi \let\next@\msafam
 \or \ifx\undefined\msbfam\loadmsbm\fi \let\next@\msbfam
 \else  \Err@{\Invalid@@\string\newsymbol}\let\next@\tw@\fi
 \mathchardef#1="#3\hexnumber@\next@#4#5\space}

\def\Bbb{\RIfM@\expandafter\Bbb@\else
 \expandafter\nonmatherr@\expandafter\Bbb\fi}
\def\Bbb@#1{{\Bbb@@{#1}}}
\def\Bbb@@#1{\noaccents@\fam\msbfam\relax#1}
\message{bold Greek and bold symbols,}
\def\loadbold{%
 \font@\tencmmib=cmmib10 \font@\sevencmmib=cmmib7 \font@\fivecmmib=cmmib5
 \skewchar\tencmmib'177 \skewchar\sevencmmib'177 \skewchar\fivecmmib'177
 \alloc@@8\fam\chardef\sixt@@n\cmmibfam
 \textfont\cmmibfam\tencmmib
 \scriptfont\cmmibfam\sevencmmib \scriptscriptfont\cmmibfam\fivecmmib
 \font@\tencmbsy=cmbsy10 \font@\sevencmbsy=cmbsy7 \font@\fivecmbsy=cmbsy5
 \skewchar\tencmbsy'60 \skewchar\sevencmbsy'60 \skewchar\fivecmbsy'60
 \alloc@@8\fam\chardef\sixt@@n\cmbsyfam
 \textfont\cmbsyfam\tencmbsy
 \scriptfont\cmbsyfam\sevencmbsy \scriptscriptfont\cmbsyfam\fivecmbsy
 \let\loadbold\empty
}
\def\boldnotloaded#1{\Err@{\ifcase#1\or First\else Second\fi
       bold symbol font not loaded}}
\def\mathchari@#1#2#3{\ifx\undefined\cmmibfam
    \boldnotloaded@\@ne
  \else\mathchar"#1\hexnumber@\cmmibfam#2#3\space \fi}
\def\mathcharii@#1#2#3{\ifx\undefined\cmbsyfam
    \boldnotloaded\tw@
  \else \mathchar"#1\hexnumber@\cmbsyfam#2#3\space\fi}
\edef\bffam@{\hexnumber@\bffam}
\def\boldkey#1{\ifcat\noexpand#1A%
  \ifx\undefined\cmmibfam \boldnotloaded\@ne
  \else {\fam\cmmibfam#1}\fi
 \else
 \ifx#1!\mathchar"5\bffam@21 \else
 \ifx#1(\mathchar"4\bffam@28 \else\ifx#1)\mathchar"5\bffam@29 \else
 \ifx#1+\mathchar"2\bffam@2B \else\ifx#1:\mathchar"3\bffam@3A \else
 \ifx#1;\mathchar"6\bffam@3B \else\ifx#1=\mathchar"3\bffam@3D \else
 \ifx#1?\mathchar"5\bffam@3F \else\ifx#1[\mathchar"4\bffam@5B \else
 \ifx#1]\mathchar"5\bffam@5D \else
 \ifx#1,\mathchari@63B \else
 \ifx#1-\mathcharii@200 \else
 \ifx#1.\mathchari@03A \else
 \ifx#1/\mathchari@03D \else
 \ifx#1<\mathchari@33C \else
 \ifx#1>\mathchari@33E \else
 \ifx#1*\mathcharii@203 \else
 \ifx#1|\mathcharii@06A \else
 \ifx#10\bold0\else\ifx#11\bold1\else\ifx#12\bold2\else\ifx#13\bold3\else
 \ifx#14\bold4\else\ifx#15\bold5\else\ifx#16\bold6\else\ifx#17\bold7\else
 \ifx#18\bold8\else\ifx#19\bold9\else
  \Err@{\string\boldkey\space can't be used with #1}%
 \fi\fi\fi\fi\fi\fi\fi\fi\fi\fi\fi\fi\fi\fi\fi
 \fi\fi\fi\fi\fi\fi\fi\fi\fi\fi\fi\fi\fi\fi}
\def\boldsymbol#1{%
 \DN@{\Err@{You can't use \string\boldsymbol\space with \string#1}#1}%
 \ifcat\noexpand#1A%
   \let\next@\relax
   \ifx\undefined\cmmibfam \boldnotloaded\@ne
   \else {\fam\cmmibfam#1}\fi
 \else
  \xdef\meaning@{\meaning#1.........}%
  \expandafter\math@\meaning@\math@
  \ifmath@
   \expandafter\mathch@\meaning@\mathch@
   \ifmathch@
    \expandafter\boldsymbol@@\meaning@\boldsymbol@@
   \fi
  \else
   \expandafter\macro@\meaning@\macro@
   \expandafter\delim@\meaning@\delim@
   \ifdelim@
    \expandafter\delim@@\meaning@\delim@@
   \else
    \boldsymbol@{#1}%
   \fi
  \fi
 \fi
 \next@}
\def\mathhexboxii@#1#2{\ifx\undefined\cmbsyfam
    \boldnotloaded\tw@
  \else \mathhexbox@{\hexnumber@\cmbsyfam}{#1}{#2}\fi}
\def\boldsymbol@#1{\let\next@\relax\let\next#1%
 \ifx\next\cdot\mathcharii@201 \else
 \ifx\next\prime{{\null\mathcharii@030 \null}}\else
 \ifx\next\lbrack\mathchar"4\bffam@5B \else
 \ifx\next\rbrack\mathchar"5\bffam@5D \else
 \ifx\next\{\mathcharii@466 \else
 \ifx\next\lbrace\mathcharii@466 \else
 \ifx\next\}\mathcharii@567 \else
 \ifx\next\rbrace\mathcharii@567 \else
 \ifx\next\surd{{\mathcharii@170}}\else
 \ifx\next\S{{\mathhexboxii@78}}\else
 \ifx\next\P{{\mathhexboxii@7B}}\else
 \ifx\next\dag{{\mathhexboxii@79}}\else
 \ifx\next\ddag{{\mathhexboxii@7A}}\else
 \DN@{\Err@{You can't use \string\boldsymbol\space with \string#1}#1}%
 \fi\fi\fi\fi\fi\fi\fi\fi\fi\fi\fi\fi\fi}
\def\boldsymbol@@#1.#2\boldsymbol@@{\classnum@#1 \count@@@\classnum@        
 \divide\classnum@4096 \count@\classnum@                                    
 \multiply\count@4096 \advance\count@@@-\count@ \count@@\count@@@           
 \divide\count@@@\@cclvi \count@\count@@                                    
 \multiply\count@@@\@cclvi \advance\count@@-\count@@@                       
 \divide\count@@@\@cclvi                                                    
 \multiply\classnum@4096 \advance\classnum@\count@@                         
 \ifnum\count@@@=\z@                                                        
  \count@"\bffam@ \multiply\count@\@cclvi
  \advance\classnum@\count@
  \DN@{\mathchar\number\classnum@}%
 \else
  \ifnum\count@@@=\@ne                                                      
   \ifx\undefined\cmmibfam \DN@{\boldnotloaded\@ne}%
   \else \count@\cmmibfam \multiply\count@\@cclvi
     \advance\classnum@\count@
     \DN@{\mathchar\number\classnum@}\fi
  \else
   \ifnum\count@@@=\tw@                                                    
     \ifx\undefined\cmbsyfam
       \DN@{\boldnotloaded\tw@}%
     \else
       \count@\cmbsyfam \multiply\count@\@cclvi
       \advance\classnum@\count@
       \DN@{\mathchar\number\classnum@}%
     \fi
  \fi
 \fi
\fi}
\newif\ifdelim@
\newcount\delimcount@
{\uccode`6=`\\ \uccode`7=`d \uccode`8=`e \uccode`9=`l
 \uppercase{\gdef\delim@#1#2#3#4#5\delim@
  {\delim@false\ifx 6#1\ifx 7#2\ifx 8#3\ifx 9#4\delim@true
   \xdef\meaning@{#5}\fi\fi\fi\fi}}}
\def\delim@@#1"#2#3#4#5#6\delim@@{\if#32%
\let\next@\relax
 \ifx\undefined\cmbsyfam \boldnotloaded\@ne
 \else \mathcharii@#2#4#5\space \fi\fi}
\def\vert{\delimiter"026A30C }
\def\Vert{\delimiter"026B30D }
\let\|\Vert
\def\backslash{\delimiter"026E30F }
\def\boldkeydots@#1{\bold@true\let\next=#1\let\delayed@=#1\mdots@@
 \boldkey#1\bold@false}  
\def\boldsymboldots@#1{\bold@true\let\next#1\let\delayed@#1\mdots@@
 \boldsymbol#1\bold@false}
\message{Euler fonts,}

\def\frak{\mathfont@\frak}

\def\loadmathfont#1{%
   \expandafter\font@\csname ten#1\endcsname=#110
   \expandafter\font@\csname seven#1\endcsname=#17
   \expandafter\font@\csname five#1\endcsname=#15
   \edef\next{\noexpand\alloc@@8\fam\chardef\sixt@@n
     \expandafter\noexpand\csname#1fam\endcsname}%
   \next
   \textfont\csname#1fam\endcsname \csname ten#1\endcsname
   \scriptfont\csname#1fam\endcsname \csname seven#1\endcsname
   \scriptscriptfont\csname#1fam\endcsname \csname five#1\endcsname
   \expandafter\def\csname #1\expandafter\endcsname\expandafter{%
      \expandafter\mathfont@\csname#1\endcsname}%
 \expandafter\gdef\csname load#1\endcsname{}%
}
\def\mathfont@#1{\RIfM@\expandafter\mathfont@@\expandafter#1\else
  \expandafter\nonmatherr@\expandafter#1\fi}
\def\mathfont@@#1#2{{\mathfont@@@#1{#2}}}
\def\mathfont@@@#1#2{\noaccents@
   \fam\csname\expandafter\eat@\string#1fam\endcsname
   \relax#2}
\message{math accents,}
\def\accentclass@{7}
\def\noaccents@{\def\accentclass@{0}}
\def\makeacc@#1#2{\def#1{\mathaccent"\accentclass@#2 }}
\makeacc@\hat{05E}
\makeacc@\check{014}
\makeacc@\tilde{07E}
\makeacc@\acute{013}
\makeacc@\grave{012}
\makeacc@\dot{05F}
\makeacc@\ddot{07F}
\makeacc@\breve{015}
\makeacc@\bar{016}

\newcount\skewcharcount@
\newcount\familycount@
\def\theskewchar@{\familycount@\@ne
 \global\skewcharcount@\the\skewchar\textfont\@ne                           
 \ifnum\fam>\m@ne\ifnum\fam<16
  \global\familycount@\the\fam\relax
  \global\skewcharcount@\the\skewchar\textfont\the\fam\relax\fi\fi          
 \ifnum\skewcharcount@>\m@ne
  \ifnum\skewcharcount@<128
  \multiply\familycount@256
  \global\advance\skewcharcount@\familycount@
  \global\advance\skewcharcount@28672
  \mathchar\skewcharcount@\else
  \global\skewcharcount@\m@ne\fi\else
 \global\skewcharcount@\m@ne\fi}                                            
\newcount\pointcount@
\def\getpoints@#1.#2\getpoints@{\pointcount@#1 }
\newdimen\accentdimen@
\newcount\accentmu@
\def\dimentomu@{\multiply\accentdimen@ 100
 \expandafter\getpoints@\the\accentdimen@\getpoints@
 \multiply\pointcount@18
 \divide\pointcount@\@m
 \global\accentmu@\pointcount@}
\def\Makeacc@#1#2{\def#1{\RIfM@\DN@{\mathaccent@
 {"\accentclass@#2 }}\else\DN@{\nonmatherr@{#1}}\fi\next@}}
\def\unbracefonts@{\let\Cal@\Cal@@\let\roman@\roman@@\let\bold@\bold@@
 \let\slanted@\slanted@@}
\def\mathaccent@#1#2{\ifnum\fam=\m@ne\xdef\thefam@{1}\else
 \xdef\thefam@{\the\fam}\fi                                                 
 \accentdimen@\z@                                                           
 \setboxz@h{\unbracefonts@$\m@th\fam\thefam@\relax#2$}
 \ifdim\accentdimen@=\z@\DN@{\mathaccent#1{#2}}
  \setbox@ne\hbox{\unbracefonts@$\m@th\fam\thefam@\relax#2\theskewchar@$}
  \setbox\tw@\hbox{$\m@th\ifnum\skewcharcount@=\m@ne\else
   \mathchar\skewcharcount@\fi$}
  \global\accentdimen@\wd@ne\global\advance\accentdimen@-\wdz@
  \global\advance\accentdimen@-\wd\tw@                                     
  \global\multiply\accentdimen@\tw@
  \dimentomu@\global\advance\accentmu@\@ne                                 
 \else\DN@{{\mathaccent#1{#2\mkern\accentmu@ mu}%
    \mkern-\accentmu@ mu}{}}\fi                                             
 \next@}\Makeacc@\Hat{05E}
\Makeacc@\Check{014}
\Makeacc@\Tilde{07E}
\Makeacc@\Acute{013}
\Makeacc@\Grave{012}
\Makeacc@\Dot{05F}
\Makeacc@\Ddot{07F}
\Makeacc@\Breve{015}
\Makeacc@\Bar{016}
\def\Vec{\RIfM@\DN@{\mathaccent@{"017E }}\else
 \DN@{\nonmatherr@\Vec}\fi\next@}
\def\accentedsymbol#1#2{\csname newbox\expandafter\endcsname
  \csname\expandafter\eat@\string#1@box\endcsname
 \expandafter\setbox\csname\expandafter\eat@
  \string#1@box\endcsname\hbox{$\m@th#2$}\define
  #1{\copy\csname\expandafter\eat@\string#1@box\endcsname{}}}
\message{roots,}
\def\sqrt#1{\radical"270370 {#1}}
\let\underline@\underline
\let\overline@\overline
\def\underline#1{\underline@{#1}}
\def\overline#1{\overline@{#1}}
\Invalid@\leftroot
\Invalid@\uproot
\newcount\uproot@
\newcount\leftroot@
\def\root{\relaxnext@
  \DN@{\ifx\next\uproot\let\next@\nextii@\else
   \ifx\next\leftroot\let\next@\nextiii@\else
   \let\next@\plainroot@\fi\fi\next@}%
  \DNii@\uproot##1{\uproot@##1\relax\FN@\nextiv@}%
  \def\nextiv@{\ifx\next\space@\DN@. {\FN@\nextv@}\else
   \DN@.{\FN@\nextv@}\fi\next@.}%
  \def\nextv@{\ifx\next\leftroot\let\next@\nextvi@\else
   \let\next@\plainroot@\fi\next@}%
  \def\nextvi@\leftroot##1{\leftroot@##1\relax\plainroot@}%
   \def\nextiii@\leftroot##1{\leftroot@##1\relax\FN@\nextvii@}%
  \def\nextvii@{\ifx\next\space@
   \DN@. {\FN@\nextviii@}\else
   \DN@.{\FN@\nextviii@}\fi\next@.}%
  \def\nextviii@{\ifx\next\uproot\let\next@\nextix@\else
   \let\next@\plainroot@\fi\next@}%
  \def\nextix@\uproot##1{\uproot@##1\relax\plainroot@}%
  \bgroup\uproot@\z@\leftroot@\z@\FN@\next@}
\def\plainroot@#1\of#2{\setbox\rootbox\hbox{$\m@th\scriptscriptstyle{#1}$}%
 \mathchoice{\r@@t\displaystyle{#2}}{\r@@t\textstyle{#2}}
 {\r@@t\scriptstyle{#2}}{\r@@t\scriptscriptstyle{#2}}\egroup}
\def\r@@t#1#2{\setboxz@h{$\m@th#1\sqrt{#2}$}%
 \dimen@\ht\z@\advance\dimen@-\dp\z@
 \setbox@ne\hbox{$\m@th#1\mskip\uproot@ mu$}\advance\dimen@ 1.667\wd@ne
 \mkern-\leftroot@ mu\mkern5mu\raise.6\dimen@\copy\rootbox
 \mkern-10mu\mkern\leftroot@ mu\boxz@}
\def\boxed#1{\setboxz@h{$\m@th\displaystyle{#1}$}\dimen@.4\ex@
 \advance\dimen@3\ex@\advance\dimen@\dp\z@
 \hbox{\lower\dimen@\hbox{%
 \vbox{\hrule height.4\ex@
 \hbox{\vrule width.4\ex@\hskip3\ex@\vbox{\vskip3\ex@\boxz@\vskip3\ex@}%
 \hskip3\ex@\vrule width.4\ex@}\hrule height.4\ex@}%
 }}}
\message{commutative diagrams,}
\let\ampersand@\relax
\newdimen\minaw@
\minaw@11.11128\ex@
\newdimen\minCDaw@
\minCDaw@2.5pc
\def\minCDarrowwidth#1{\RIfMIfI@\onlydmatherr@\minCDarrowwidth
 \else\minCDaw@#1\relax\fi\else\onlydmatherr@\minCDarrowwidth\fi}
\newif\ifCD@
\def\CD{\bgroup\vspace@\relax\let\ampersand@&\iffalse}\fi
 \CD@true\vcenter\bgroup\Let@\tabskip\z@skip\baselineskip20\ex@
 \lineskip3\ex@\lineskiplimit3\ex@\halign\bgroup
 &\hfill$\m@th##$\hfill\crcr}
\def\endCD{\crcr\egroup\egroup\egroup}
\newdimen\bigaw@
\atdef@>#1>#2>{\ampersand@                                                  
 \setboxz@h{$\m@th\ssize\;{#1}\;\;$}
 \setbox@ne\hbox{$\m@th\ssize\;{#2}\;\;$}
 \setbox\tw@\hbox{$\m@th#2$}
 \ifCD@\global\bigaw@\minCDaw@\else\global\bigaw@\minaw@\fi                 
 \ifdim\wdz@>\bigaw@\global\bigaw@\wdz@\fi
 \ifdim\wd@ne>\bigaw@\global\bigaw@\wd@ne\fi                                
 \ifCD@\enskip\fi                                                           
 \ifdim\wd\tw@>\z@
  \mathrel{\mathop{\hbox to\bigaw@{\rightarrowfill@\displaystyle}}%
    \limits^{#1}_{#2}}
 \else\mathrel{\mathop{\hbox to\bigaw@{\rightarrowfill@\displaystyle}}%
    \limits^{#1}}\fi                                                        
 \ifCD@\enskip\fi                                                          
 \ampersand@}                                                              
\atdef@<#1<#2<{\ampersand@\setboxz@h{$\m@th\ssize\;\;{#1}\;$}%
 \setbox@ne\hbox{$\m@th\ssize\;\;{#2}\;$}\setbox\tw@\hbox{$\m@th#2$}%
 \ifCD@\global\bigaw@\minCDaw@\else\global\bigaw@\minaw@\fi
 \ifdim\wdz@>\bigaw@\global\bigaw@\wdz@\fi
 \ifdim\wd@ne>\bigaw@\global\bigaw@\wd@ne\fi
 \ifCD@\enskip\fi
 \ifdim\wd\tw@>\z@
  \mathrel{\mathop{\hbox to\bigaw@{\leftarrowfill@\displaystyle}}%
       \limits^{#1}_{#2}}\else
  \mathrel{\mathop{\hbox to\bigaw@{\leftarrowfill@\displaystyle}}%
       \limits^{#1}}\fi
 \ifCD@\enskip\fi\ampersand@}
\begingroup
 \catcode`\~=\active \lccode`\~=`\@
 \lowercase{%
  \global\atdef@)#1)#2){~>#1>#2>}
  \global\atdef@(#1(#2({~<#1<#2<}}
\endgroup
\atdef@ A#1A#2A{\llap{$\m@th\vcenter{\hbox
 {$\ssize#1$}}$}\Big\uparrow\rlap{$\m@th\vcenter{\hbox{$\ssize#2$}}$}&&}
\atdef@ V#1V#2V{\llap{$\m@th\vcenter{\hbox
 {$\ssize#1$}}$}\Big\downarrow\rlap{$\m@th\vcenter{\hbox{$\ssize#2$}}$}&&}
\atdef@={&\enskip\mathrel
 {\vbox{\hrule width\minCDaw@\vskip3\ex@\hrule width
 \minCDaw@}}\enskip&}
\atdef@|{\Big\Vert&&}
\atdef@\vert{\Big\Vert&&}
\def\pretend#1\haswidth#2{\setboxz@h{$\m@th\scriptstyle{#2}$}\hbox
 to\wdz@{\hfill$\m@th\scriptstyle{#1}$\hfill}}
\message{poor man's bold,}
\def\pmb{\RIfM@\expandafter\mathpalette\expandafter\pmb@\else
 \expandafter\pmb@@\fi}
\def\pmb@@#1{\leavevmode\setboxz@h{#1}%
   \dimen@-\wdz@
   \kern-.5\ex@\copy\z@
   \kern\dimen@\kern.25\ex@\raise.4\ex@\copy\z@
   \kern\dimen@\kern.25\ex@\box\z@
}
\def\binrel@@#1{\ifdim\wd2<\z@\mathbin{#1}\else\ifdim\wd\tw@>\z@
 \mathrel{#1}\else{#1}\fi\fi}
\newdimen\pmbraise@
\def\pmb@#1#2{\setbox\thr@@\hbox{$\m@th#1{#2}$}%
 \setbox4\hbox{$\m@th#1\mkern.5mu$}\pmbraise@\wd4\relax
 \binrel@{#2}%
 \dimen@-\wd\thr@@
   \binrel@@{%
   \mkern-.8mu\copy\thr@@
   \kern\dimen@\mkern.4mu\raise\pmbraise@\copy\thr@@
   \kern\dimen@\mkern.4mu\box\thr@@
}}
\def\documentstyle#1{\W@{}\input #1.sty\relax}
\message{syntax check,}
\font\dummyft@=dummy
\fontdimen1 \dummyft@=\z@
\fontdimen2 \dummyft@=\z@
\fontdimen3 \dummyft@=\z@
\fontdimen4 \dummyft@=\z@
\fontdimen5 \dummyft@=\z@
\fontdimen6 \dummyft@=\z@
\fontdimen7 \dummyft@=\z@
\fontdimen8 \dummyft@=\z@
\fontdimen9 \dummyft@=\z@
\fontdimen10 \dummyft@=\z@
\fontdimen11 \dummyft@=\z@
\fontdimen12 \dummyft@=\z@
\fontdimen13 \dummyft@=\z@
\fontdimen14 \dummyft@=\z@
\fontdimen15 \dummyft@=\z@
\fontdimen16 \dummyft@=\z@
\fontdimen17 \dummyft@=\z@
\fontdimen18 \dummyft@=\z@
\fontdimen19 \dummyft@=\z@
\fontdimen20 \dummyft@=\z@
\fontdimen21 \dummyft@=\z@
\fontdimen22 \dummyft@=\z@
\def\fontlist@{\\{\tenrm}\\{\sevenrm}\\{\fiverm}\\{\teni}\\{\seveni}%
 \\{\fivei}\\{\tensy}\\{\sevensy}\\{\fivesy}\\{\tenex}\\{\tenbf}\\{\sevenbf}%
 \\{\fivebf}\\{\tensl}\\{\tenit}}
\def\font@#1=#2 {\rightappend@#1\to\fontlist@\font#1=#2 }
\def\dodummy@{{\def\\##1{\global\let##1\dummyft@}\fontlist@}}
\def\nopages@{\output{\setbox\z@\box\@cclv \deadcycles\z@}%
 \alloc@5\toks\toksdef\@cclvi\output}
\let\galleys\nopages@
\newif\ifsyntax@
\newcount\countxviii@
\def\syntax{\syntax@true\dodummy@\countxviii@\count18
 \loop\ifnum\countxviii@>\m@ne\textfont\countxviii@=\dummyft@
 \scriptfont\countxviii@=\dummyft@\scriptscriptfont\countxviii@=\dummyft@
 \advance\countxviii@\m@ne\repeat                                           
 \dummyft@\tracinglostchars\z@\nopages@\frenchspacing\hbadness\@M}
\def\first@#1#2\end{#1}
\def\printoptions{\W@{Do you want S(yntax check),
  G(alleys) or P(ages)?}%
 \message{Type S, G or P, followed by <return>: }%
 \begingroup 
 \endlinechar\m@ne 
 \read\m@ne to\ans@
 \edef\ans@{\uppercase{\def\noexpand\ans@{%
   \expandafter\first@\ans@ P\end}}}%
 \expandafter\endgroup\ans@
 \if\ans@ P
 \else \if\ans@ S\syntax
 \else \if\ans@ G\galleys
 \else\message{? Unknown option: \ans@; using the `pages' option.}%
 \fi\fi\fi}
\def\alloc@#1#2#3#4#5{\global\advance\count1#1by\@ne
 \ch@ck#1#4#2\allocationnumber=\count1#1
 \global#3#5=\allocationnumber
 \ifalloc@\wlog{\string#5=\string#2\the\allocationnumber}\fi}
\def\document{\def\alloclist@{}\def\fontlist@{}}
\let\enddocument\bye

\let\proclaim\undefined
\let\footnote\undefined
\let\=\undefined
\let\>\undefined

\catcode`\@=\active
\message{... finished}

\expandafter\ifx\csname mathdefs.tex\endcsname\relax
  \expandafter\gdef\csname mathdefs.tex\endcsname{}
\else \message{Hey!  Apparently you were trying to
  \string\input{mathdefs.tex} twice.   This does not make sense.} 
\errmessage{Please edit your file (probably \jobname.tex) and remove
any duplicate ``\string\input'' lines}\endinput\fi




\catcode`\X=12\catcode`\@=11

\def\n@wcount{\alloc@0\count\countdef\insc@unt}
\def\n@wwrite{\alloc@7\write\chardef\sixt@@n}
\def\n@wread{\alloc@6\read\chardef\sixt@@n}
\def\r@s@t{\relax}\def\v@idline{\par}\def\@mputate#1/{#1}
\def\l@c@l#1X{\firstpart.#1}\def\gl@b@l#1X{#1}\def\t@d@l#1X{{}}

\def\crossrefs#1{\ifx\all#1\let\tr@ce=\all\else\def\tr@ce{#1,}\fi
   \n@wwrite\cit@tionsout\openout\cit@tionsout=\jobname.cit 
   \write\cit@tionsout{\tr@ce}\expandafter\setfl@gs\tr@ce,}
\def\setfl@gs#1,{\def\@{#1}\ifx\@\empty\let\next=\relax
   \else\let\next=\setfl@gs\expandafter\xdef
   \csname#1tr@cetrue\endcsname{}\fi\next}
\def\m@ketag#1#2{\expandafter\n@wcount\csname#2tagno\endcsname
     \csname#2tagno\endcsname=0\let\tail=\all\xdef\all{\tail#2,}
   \ifx#1\l@c@l\let\tail=\r@s@t\xdef\r@s@t{\csname#2tagno\endcsname=0\tail}\fi
   \expandafter\gdef\csname#2cite\endcsname##1{\expandafter
     \ifx\csname#2tag##1\endcsname\relax?\else\csname#2tag##1\endcsname\fi
     \expandafter\ifx\csname#2tr@cetrue\endcsname\relax\else
     \write\cit@tionsout{#2tag ##1 cited on page \folio.}\fi}
   \expandafter\gdef\csname#2page\endcsname##1{\expandafter
     \ifx\csname#2page##1\endcsname\relax?\else\csname#2page##1\endcsname\fi
     \expandafter\ifx\csname#2tr@cetrue\endcsname\relax\else
     \write\cit@tionsout{#2tag ##1 cited on page \folio.}\fi}
   \expandafter\gdef\csname#2tag\endcsname##1{\expandafter
      \ifx\csname#2check##1\endcsname\relax
      \expandafter\xdef\csname#2check##1\endcsname{}%
      \else\immediate\write16{Warning: #2tag ##1 used more than once.}\fi
      \multit@g{#1}{#2}##1/X%
      \write\t@gsout{#2tag ##1 assigned number \csname#2tag##1\endcsname\space
      on page \number\count0.}%
   \csname#2tag##1\endcsname}}

\def\multit@g#1#2#3/#4X{\def\t@mp{#4}\ifx\t@mp\empty%
      \global\advance\csname#2tagno\endcsname by 1 
      \expandafter\xdef\csname#2tag#3\endcsname
      {#1\number\csname#2tagno\endcsnameX}%
   \else\expandafter\ifx\csname#2last#3\endcsname\relax
      \expandafter\n@wcount\csname#2last#3\endcsname
      \global\advance\csname#2tagno\endcsname by 1 
      \expandafter\xdef\csname#2tag#3\endcsname
      {#1\number\csname#2tagno\endcsnameX}
      \write\t@gsout{#2tag #3 assigned number \csname#2tag#3\endcsname\space
      on page \number\count0.}\fi
   \global\advance\csname#2last#3\endcsname by 1
   \def\t@mp{\expandafter\xdef\csname#2tag#3/}%
   \expandafter\t@mp\@mputate#4\endcsname
   {\csname#2tag#3\endcsname\lastpart{\csname#2last#3\endcsname}}\fi}
\def\t@gs#1{\def\all{}\m@ketag#1e\m@ketag#1s\m@ketag\t@d@l p
\let\realscite\scite
\let\realstag\stag
   \m@ketag\gl@b@l r \n@wread\t@gsin
   \openin\t@gsin=\jobname.tgs \re@der \closein\t@gsin
   \n@wwrite\t@gsout\openout\t@gsout=\jobname.tgs }
\outer\def\localtags{\t@gs\l@c@l}
\outer\def\globaltags{\t@gs\gl@b@l}
\outer\def\newlocaltag#1{\m@ketag\l@c@l{#1}}
\outer\def\newglobaltag#1{\m@ketag\gl@b@l{#1}}

\newif\ifpr@ 
\def\m@kecs #1tag #2 assigned number #3 on page #4.%
   {\expandafter\gdef\csname#1tag#2\endcsname{#3}
   \expandafter\gdef\csname#1page#2\endcsname{#4}
   \ifpr@\expandafter\xdef\csname#1check#2\endcsname{}\fi}
\def\re@der{\ifeof\t@gsin\let\next=\relax\else
   \read\t@gsin to\t@gline\ifx\t@gline\v@idline\else
   \expandafter\m@kecs \t@gline\fi\let \next=\re@der\fi\next}
\def\pretags#1{\pr@true\pret@gs#1,,}
\def\pret@gs#1,{\def\@{#1}\ifx\@\empty\let\n@xtfile=\relax
   \else\let\n@xtfile=\pret@gs \openin\t@gsin=#1.tgs \message{#1} \re@der 
   \closein\t@gsin\fi \n@xtfile}

\newcount\sectno\sectno=0\newcount\subsectno\subsectno=0
\newif\ifultr@local \def\ultralocal{\ultr@localtrue}
\def\firstpart{\number\sectno}
\def\lastpart#1{\ifcase#1 \or a\or b\or c\or d\or e\or f\or g\or h\or 
   i\or k\or l\or m\or n\or o\or p\or q\or r\or s\or t\or u\or v\or w\or 
   x\or y\or z \fi}

\def\resetall{\global\advance\sectno by 1\subsectno=0
   \gdef\firstpart{\number\sectno}\r@s@t}
\def\resetsub{\global\advance\subsectno by 1
   \gdef\firstpart{\number\sectno.\number\subsectno}\r@s@t}
\def\newsection#1\par{\resetall\vskip0pt plus.3\vsize\penalty-250
   \vskip0pt plus-.3\vsize\bigskip\bigskip
   \message{#1}\leftline{\bf#1}\nobreak\bigskip}
\def\subsection#1\par{\ifultr@local\resetsub\fi
   \vskip0pt plus.2\vsize\penalty-250\vskip0pt plus-.2\vsize
   \bigskip\smallskip\message{#1}\leftline{\bf#1}\nobreak\medskip}


\newdimen\marginshift

\newdimen\margindelta
\newdimen\marginmax
\newdimen\marginmin

\def\margininit{       
\marginmax=3 true cm                  
				      
\margindelta=0.1 true cm              
\marginmin=0.1true cm                 
\marginshift=\marginmin
}    

\def\t@gsjj#1,{\def\@{#1}\ifx\@\empty\let\next=\relax\else\let\next=\t@gsjj
   \def\@@{p}\ifx\@\@@\else
   \expandafter\gdef\csname#1cite\endcsname##1{\citejj{##1}}
   \expandafter\gdef\csname#1page\endcsname##1{?}
   \expandafter\gdef\csname#1tag\endcsname##1{\tagjj{##1}}\fi\fi\next}
\newif\ifshowstuffinmargin
\showstuffinmarginfalse
\def\jjtags{\ifx\shlhetal\relax 
  \else
\ifx\shlhetal\undefinedcontrolseq
\else
\showstuffinmargintrue
\ifx\all\relax\else\expandafter\t@gsjj\all,\fi\fi \fi
}

\def\tagjj#1{\realstag{#1}\mginpar{\zeigen{#1}}}
\def\citejj#1{\rechnen{#1}\mginpar{\zeigen{#1}}}     

\def\rechnen#1{\expandafter\ifx\csname stag#1\endcsname\relax ??\else
                           \csname stag#1\endcsname\fi}

\newdimen\theight

\def\marginfont{\sevenrm}

\def\trymarginbox#1{\setbox0=\hbox{\marginfont\hskip\marginshift #1}%
		\global\marginshift\wd0 
		\global\advance\marginshift\margindelta}

\def \mginpar#1{%
\ifvmode\setbox0\hbox to \hsize{\hfill\rlap{\marginfont\quad#1}}%
\ht0 0cm
\dp0 0cm
\box0\vskip-\baselineskip
\else 
             \vadjust{\trymarginbox{#1}%
		\ifdim\marginshift>\marginmax \global\marginshift\marginmin
			\trymarginbox{#1}%
                \fi
             \theight=\ht0
             \advance\theight by \dp0    \advance\theight by \lineskip
             \kern -\theight \vbox to \theight{\rightline{\rlap{\box0}}%
\vss}}\fi}


\def\t@gsoff#1,{\def\@{#1}\ifx\@\empty\let\next=\relax\else\let\next=\t@gsoff
   \def\@@{p}\ifx\@\@@\else
   \expandafter\gdef\csname#1cite\endcsname##1{\zeigen{##1}}
   \expandafter\gdef\csname#1page\endcsname##1{?}
   \expandafter\gdef\csname#1tag\endcsname##1{\zeigen{##1}}\fi\fi\next}
\def\verbatimtags{\showstuffinmarginfalse
\ifx\all\relax\else\expandafter\t@gsoff\all,\fi}
\def\zeigen#1{\hbox{$\langle$}#1\hbox{$\rangle$}}

\def\margincite#1{\ifshowstuffinmargin\mginpar{\zeigen{#1}}\fi}

\def\margintag#1{\ifshowstuffinmargin\mginpar{\zeigen{#1}}\fi}

\def\marginplain#1{\ifshowstuffinmargin\mginpar{{#1}}\fi}
\def\marginbf#1{\marginplain{{\bf \ \ #1}}}

\def\(#1){\edef\dot@g{\ifmmode\ifinner(\hbox{\noexpand\etag{#1}})
   \else\noexpand\eqno(\hbox{\noexpand\etag{#1}})\fi
   \else(\noexpand\ecite{#1})\fi}\dot@g}

\newif\ifbr@ck
\def\eat#1{}
\def\[#1]{\br@cktrue[\br@cket#1'X]}
\def\br@cket#1'#2X{\def\temp{#2}\ifx\temp\empty\let\next\eat
   \else\let\next\br@cket\fi
   \ifbr@ck\br@ckfalse\br@ck@t#1,X\else\br@cktrue#1\fi\next#2X}
\def\br@ck@t#1,#2X{\def\temp{#2}\ifx\temp\empty\let\neext\eat
   \else\let\neext\br@ck@t\def\temp{,}\fi
   \def\teemp{#1}\ifx\teemp\empty\else\rcite{#1}\fi\temp\neext#2X}
\def\resetbr@cket{\gdef\[##1]{[\rtag{##1}]}}
\def\references{\resetbr@cket\newsection References\par}

\newtoks\symb@ls\newtoks\s@mb@ls\newtoks\p@gelist\n@wcount\ftn@mber
    \ftn@mber=1\newif\ifftn@mbers\ftn@mbersfalse\newif\ifbyp@ge\byp@gefalse
\def\defm@rk{\ifftn@mbers\n@mberm@rk\else\symb@lm@rk\fi}
\def\n@mberm@rk{\xdef\m@rk{{\the\ftn@mber}}%
    \global\advance\ftn@mber by 1 }
\def\rot@te#1{\let\temp=#1\global#1=\expandafter\r@t@te\the\temp,X}
\def\r@t@te#1,#2X{{#2#1}\xdef\m@rk{{#1}}}
\def\b@@st#1{{$^{#1}$}}\def\str@p#1{#1}
\def\symb@lm@rk{\ifbyp@ge\rot@te\p@gelist\ifnum\expandafter\str@p\m@rk=1 
    \s@mb@ls=\symb@ls\fi\write\f@nsout{\number\count0}\fi \rot@te\s@mb@ls}
\def\byp@ge{\byp@getrue\n@wwrite\f@nsin\openin\f@nsin=\jobname.fns 
    \n@wcount\currentp@ge\currentp@ge=0\p@gelist={0}
    \re@dfns\closein\f@nsin\rot@te\p@gelist
    \n@wread\f@nsout\openout\f@nsout=\jobname.fns }
\def\m@kelist#1X#2{{#1,#2}}
\def\re@dfns{\ifeof\f@nsin\let\next=\relax\else\read\f@nsin to \f@nline
    \ifx\f@nline\v@idline\else\let\t@mplist=\p@gelist
    \ifnum\currentp@ge=\f@nline
    \global\p@gelist=\expandafter\m@kelist\the\t@mplistX0
    \else\currentp@ge=\f@nline
    \global\p@gelist=\expandafter\m@kelist\the\t@mplistX1\fi\fi
    \let\next=\re@dfns\fi\next}
\def\symbols#1{\symb@ls={#1}\s@mb@ls=\symb@ls} 
\def\bigsymbol{\textstyle}
\symbols{\bigsymbol\ast,\dagger,\ddagger,\sharp,\flat,\natural,\star}
\def\ftnumbers{\ftn@mberstrue} \def\ftsymbols{\ftn@mbersfalse}
\def\paginal{\byp@ge} \def\resetftnumbers{\ftn@mber=1}
\def\ftnote#1{\defm@rk\expandafter\expandafter\expandafter\footnote
    \expandafter\b@@st\m@rk{#1}}

\long\def\jump#1\endjump{}
\def\ssum{\mathop{\lower .1em\hbox{$\textstyle\Sigma$}}\nolimits}

\def\qed{\nobreak\kern 1em \vrule height .5em width .5em depth 0em}
\def\newneq{\hbox{\rlap{\hbox to 1\wd9{\hss$=$\hss}}\raise .1em 
   \hbox to 1\wd9{\hss$\scriptscriptstyle/$\hss}}}
\def\subsetne{\setbox9 = \hbox{$\subset$}\mathrel{\hbox{\rlap
   {\lower .4em \newneq}\raise .13em \hbox{$\subset$}}}}
\def\supsetne{\setbox9 = \hbox{$\subset$}\mathrel{\hbox{\rlap
   {\lower .4em \newneq}\raise .13em \hbox{$\supset$}}}}

\def\vbar{\mathchoice{\vrule height6.3ptdepth-.5ptwidth.8pt\kern-.8pt}
   {\vrule height6.3ptdepth-.5ptwidth.8pt\kern-.8pt}
   {\vrule height4.1ptdepth-.35ptwidth.6pt\kern-.6pt}
   {\vrule height3.1ptdepth-.25ptwidth.5pt\kern-.5pt}}
\def\f@dge{\mathchoice{}{}{\mkern.5mu}{\mkern.8mu}}
\def\b@c#1#2{{\rm \mkern#2mu\vbar\mkern-#2mu#1}}
\def\b@b#1{{\rm I\mkern-3.5mu #1}}
\def\b@a#1#2{{\rm #1\mkern-#2mu\f@dge #1}}
\def\bb#1{{\count4=`#1 \advance\count4by-64 \ifcase\count4\or\b@a A{11.5}\or
   \b@b B\or\b@c C{5}\or\b@b D\or\b@b E\or\b@b F \or\b@c G{5}\or\b@b H\or
   \b@b I\or\b@c J{3}\or\b@b K\or\b@b L \or\b@b M\or\b@b N\or\b@c O{5} \or
   \b@b P\or\b@c Q{5}\or\b@b R\or\b@a S{8}\or\b@a T{10.5}\or\b@c U{5}\or
   \b@a V{12}\or\b@a W{16.5}\or\b@a X{11}\or\b@a Y{11.7}\or\b@a Z{7.5}\fi}}

\catcode`\X=11 \catcode`\@=12




\let\thischap\jobname

\def\partof#1{\csname returnthe#1part\endcsname}
\def\chapof#1{\csname returnthe#1chap\endcsname}

\def\setchapter#1,#2,#3.{%
  \expandafter\def\csname returnthe#1part\endcsname{#2}%
  \expandafter\def\csname returnthe#1chap\endcsname{#3}%
}

\setchapter 300a,A,I.
\setchapter 300b,A,II.
\setchapter 300c,A,III.
\setchapter 300d,A,IV.
\setchapter 300e,A,V.
\setchapter 300f,A,VI.
\setchapter 300g,A,VII.
\setchapter   88,B,I.
\setchapter  600,B,II.
\setchapter  705,B,III.

\def\cprefix#1{
\edef\theotherpart{\partof{#1}}\edef\theotherchap{\chapof{#1}}%
\ifx\theotherpart\thispart
   \ifx\theotherchap\thischap 
    \else 
     \theotherchap%
    \fi
   \else 
     \theotherpart.\theotherchap\fi}

\def\sectioncite[#1]#2{%
     \cprefix{#2}#1}

\def\chaptercite#1{Chapter \cprefix{#1}}

\edef\thispart{\partof{\thischap}}
\edef\thischap{\chapof{\thischap}}


\def\spuriousreset{}


\expandafter\ifx\csname citeadd.tex\endcsname\relax
\expandafter\gdef\csname citeadd.tex\endcsname{}
\else \message{Hey!  Apparently you were trying to
\string\input{citeadd.tex} twice.   This does not make sense.} 
\errmessage{Please edit your file (probably \jobname.tex) and remove
any duplicate ``\string\input'' lines}\endinput\fi

\def\sciteu{\sciteerror{undefined}}

\def\sciteerror#1#2{{\mathortextbf{\scite{#2}}}\complainaboutcitation{#1}{#2}}
\def\mathortextbf#1{\hbox{\bf #1}}
\def\complainaboutcitation#1#2{%
\vadjust{\line{\llap{---$\!\!>$ }\qquad scite$\{$#2$\}$ #1\hfil}}}

\sectno=-1   
\localtags
\jjtags
\NoBlackBoxes
\newbox\noforkbox \newdimen\forklinewidth
\forklinewidth=0.3pt   
\setbox0\hbox{$\textstyle\bigcup$}
\setbox1\hbox to \wd0{\hfil\vrule width \forklinewidth depth \dp0
                        height \ht0 \hfil}
\wd1=0 cm
\setbox\noforkbox\hbox{\box1\box0\relax}
\def\unionstick{\mathop{\copy\noforkbox}\limits}
\def\nonfork#1#2_#3{#1\unionstick_{\textstyle #3}#2}
\def\nonforkin#1#2_#3^#4{#1\unionstick_{\textstyle #3}^{\textstyle #4}#2}     
%
\setbox0\hbox{$\textstyle\bigcup$}
\setbox1\hbox to \wd0{\hfil{\sl /\/}\hfil}
\setbox2\hbox to \wd0{\hfil\vrule height \ht0 depth \dp0 width
                                \forklinewidth\hfil}
\wd1=0cm
\wd2=0cm
\newbox\doesforkbox
\setbox\doesforkbox\hbox{\box1\box0\relax}
\def\nunionstick{\mathop{\copy\doesforkbox}\limits}

\def\fork#1#2_#3{#1\nunionstick_{\textstyle #3}#2}
\def\forkin#1#2_#3^#4{#1\nunionstick_{\textstyle #3}^{\textstyle #4}#2}     
\define\mr{\medskip\roster}
\define\sn{\smallskip\noindent}
\define\mn{\medskip\noindent}
\define\bn{\bigskip\noindent}
\define\ub{\underbar}
\define\wilog{\text{without loss of generality}}
\define\ermn{\endroster\medskip\noindent}
\define\dbca{\dsize\bigcap}
\define\dbcu{\dsize\bigcup}
\define \nl{\newline}
\magnification=\magstep 1
\documentstyle{amsppt}

{    
\catcode`@11

\ifx\alicetwothousandloaded@\relax
  \endinput\else\global\let\alicetwothousandloaded@\relax\fi

\gdef\subjclass{\let\savedef@\subjclass
 \def\subjclass##1\endsubjclass{\let\subjclass\savedef@
   \toks@{\def\usualspace{{\rm\enspace}}\eightpoint}%
   \toks@@{##1\unskip.}%
   \edef\thesubjclass@{\the\toks@
     \frills@{{\noexpand\rm2000 {\noexpand\it Mathematics Subject
       Classification}.\noexpand\enspace}}%
     \the\toks@@}}%
  \nofrillscheck\subjclass}
} 


\expandafter\ifx\csname alice2jlem.tex\endcsname\relax
  \expandafter\xdef\csname alice2jlem.tex\endcsname{\the\catcode`@}
\else \message{Hey!  Apparently you were trying to
\string\input{alice2jlem.tex}  twice.   This does not make sense.}
\errmessage{Please edit your file (probably \jobname.tex) and remove
any duplicate ``\string\input'' lines}\endinput\fi

\expandafter\ifx\csname bib4plain.tex\endcsname\relax
  \expandafter\gdef\csname bib4plain.tex\endcsname{}
\else \message{Hey!  Apparently you were trying to \string\input
  bib4plain.tex twice.   This does not make sense.}
\errmessage{Please edit your file (probably \jobname.tex) and remove
any duplicate ``\string\input'' lines}\endinput\fi

\def\renewcommand{\newcommand}	       
\edef\cite{\the\catcode`@}%
\catcode`@ = 11
\let\@oldatcatcode = \cite
\chardef\@letter = 11
\chardef\@other = 12
%
%
%
%
\def\@innerdef#1#2{\edef#1{\expandafter\noexpand\csname #2\endcsname}}%
%
%
\@innerdef\@innernewcount{newcount}%
\@innerdef\@innernewdimen{newdimen}%
\@innerdef\@innernewif{newif}%
\@innerdef\@innernewwrite{newwrite}%
%
%
%
\def\@gobble#1{}%
%
%
%
\ifx\inputlineno\@undefined
   \let\@linenumber = \empty 
\else
   \def\@linenumber{\the\inputlineno:\space}%
\fi
%
%
%
\def\@futurenonspacelet#1{\def\cs{#1}%
   \afterassignment\@stepone\let\@nexttoken=
}%
\begingroup 
\def\\{\global\let\@stoken= }%
\\ 
\endgroup
\def\@stepone{\expandafter\futurelet\cs\@steptwo}%
\def\@steptwo{\expandafter\ifx\cs\@stoken\let\@@next=\@stepthree
   \else\let\@@next=\@nexttoken\fi \@@next}%
\def\@stepthree{\afterassignment\@stepone\let\@@next= }%
%
%
%
\def\@getoptionalarg#1{%
   \let\@optionaltemp = #1%
   \let\@optionalnext = \relax
   \@futurenonspacelet\@optionalnext\@bracketcheck
}%
%
%
\def\@bracketcheck{%
   \ifx [\@optionalnext
      \expandafter\@@getoptionalarg
   \else
      \let\@optionalarg = \empty
      \expandafter\@optionaltemp
   \fi
}%
\def\@@getoptionalarg[#1]{%
   \def\@optionalarg{#1}%
   \@optionaltemp
}%
%
%
%
\def\@nnil{\@nil}%
\def\@fornoop#1\@@#2#3{}%
\def\@for#1:=#2\do#3{%
   \edef\@fortmp{#2}%
   \ifx\@fortmp\empty \else
      \expandafter\@forloop#2,\@nil,\@nil\@@#1{#3}%
   \fi
}%
\def\@forloop#1,#2,#3\@@#4#5{\def#4{#1}\ifx #4\@nnil \else
       #5\def#4{#2}\ifx #4\@nnil \else#5\@iforloop #3\@@#4{#5}\fi\fi
}%
\def\@iforloop#1,#2\@@#3#4{\def#3{#1}\ifx #3\@nnil
       \let\@nextwhile=\@fornoop \else
      #4\relax\let\@nextwhile=\@iforloop\fi\@nextwhile#2\@@#3{#4}%
}%
%
%
%
\@innernewif\if@fileexists
\def\@testfileexistence{\@getoptionalarg\@finishtestfileexistence}%
\def\@finishtestfileexistence#1{%
   \begingroup
      \def\extension{#1}%
      \immediate\openin0 =
         \ifx\@optionalarg\empty\jobname\else\@optionalarg\fi
         \ifx\extension\empty \else .#1\fi
         \space
      \ifeof 0
         \global\@fileexistsfalse
      \else
         \global\@fileexiststrue
      \fi
      \immediate\closein0
   \endgroup
}%
%
%
%
%
\def\bibliographystyle#1{%
   \@readauxfile
   \@writeaux{\string\bibstyle{#1}}%
}%
\let\bibstyle = \@gobble
%
%
\let\bblfilebasename = \jobname
\def\bibliography#1{%
   \@readauxfile
   \@writeaux{\string\bibdata{#1}}%
   \@testfileexistence[\bblfilebasename]{bbl}%
   \if@fileexists
      \nobreak
      \@readbblfile
   \fi
}%
\let\bibdata = \@gobble
%
%
\def\nocite#1{%
   \@readauxfile
   \@writeaux{\string\citation{#1}}%
}%
\@innernewif\if@notfirstcitation
%
%
\def\cite{\@getoptionalarg\@cite}%
%
%
\def\@cite#1{%
   \let\@citenotetext = \@optionalarg
   \printcitestart
   \nocite{#1}%
   \@notfirstcitationfalse
   \@for \@citation :=#1\do
   {%
      \expandafter\@onecitation\@citation\@@
   }%
   \ifx\empty\@citenotetext\else
      \printcitenote{\@citenotetext}%
   \fi
   \printcitefinish
}%
\newif\ifweareinprivate
\weareinprivatetrue
\ifx\shlhetal\undefinedcontrolseq\weareinprivatefalse\fi
\ifx\shlhetal\relax\weareinprivatefalse\fi
\def\@onecitation#1\@@{%
   \if@notfirstcitation
      \printbetweencitations
   \fi
   \expandafter \ifx \csname\@citelabel{#1}\endcsname \relax
      \if@citewarning
         \message{\@linenumber Undefined citation `#1'.}%
      \fi
     \ifweareinprivate
      \expandafter\gdef\csname\@citelabel{#1}\endcsname{%
\strut 
\vadjust{\vskip-\dp\strutbox
\vbox to 0pt{\vss\parindent0cm \leftskip=\hsize 
\advance\leftskip3mm
\advance\hsize 4cm\strut\openup-4pt 
\rightskip 0cm plus 1cm minus 0.5cm ?  #1 ?\strut}}
         {\tt
            \escapechar = -1
            \nobreak\hskip0pt\pfeilsw
            \expandafter\string\csname#1\endcsname
             \pfeilso
            \nobreak\hskip0pt
         }%
      }%
     \else  
      \expandafter\gdef\csname\@citelabel{#1}\endcsname{%
            {\tt\expandafter\string\csname#1\endcsname}
      }%
     \fi  
   \fi
   \csname\@citelabel{#1}\endcsname
   \@notfirstcitationtrue
}%
%
%
\def\@citelabel#1{b@#1}%
%
%
\def\@citedef#1#2{\expandafter\gdef\csname\@citelabel{#1}\endcsname{#2}}%
%
%
%
\def\@readbblfile{%
   \ifx\@itemnum\@undefined
      \@innernewcount\@itemnum
   \fi
   \begingroup
      \def\begin##1##2{%
         \setbox0 = \hbox{\biblabelcontents{##2}}%
         \biblabelwidth = \wd0
      }%
      \def\end##1{}
      %
      %
      \@itemnum = 0
      \def\bibitem{\@getoptionalarg\@bibitem}%
      \def\@bibitem{%
         \ifx\@optionalarg\empty
            \expandafter\@numberedbibitem
         \else
            \expandafter\@alphabibitem
         \fi
      }%
      \def\@alphabibitem##1{%
         \expandafter \xdef\csname\@citelabel{##1}\endcsname {\@optionalarg}%
         \ifx\biblabelprecontents\@undefined
            \let\biblabelprecontents = \relax
         \fi
         \ifx\biblabelpostcontents\@undefined
            \let\biblabelpostcontents = \hss
         \fi
         \@finishbibitem{##1}%
      }%
      \def\@numberedbibitem##1{%
         \advance\@itemnum by 1
         \expandafter \xdef\csname\@citelabel{##1}\endcsname{\number\@itemnum}%
         \ifx\biblabelprecontents\@undefined
            \let\biblabelprecontents = \hss
         \fi
         \ifx\biblabelpostcontents\@undefined
            \let\biblabelpostcontents = \relax
         \fi
         \@finishbibitem{##1}%
      }%
      \def\@finishbibitem##1{%
         \biblabelprint{\csname\@citelabel{##1}\endcsname}%
         \@writeaux{\string\@citedef{##1}{\csname\@citelabel{##1}\endcsname}}%
         \ignorespaces
      }%
      %
      %
      \let\em = \bblem
      \let\newblock = \bblnewblock
      \let\sc = \bblsc
      \frenchspacing
      \clubpenalty = 4000 \widowpenalty = 4000
      \tolerance = 10000 \hfuzz = .5pt
      \everypar = {\hangindent = \biblabelwidth
                      \advance\hangindent by \biblabelextraspace}%
      \bblrm
      \parskip = 1.5ex plus .5ex minus .5ex
      \biblabelextraspace = .5em
      \bblhook
      \input \bblfilebasename.bbl
   \endgroup
}%
%
%
\@innernewdimen\biblabelwidth
\@innernewdimen\biblabelextraspace
%
%
%
\def\biblabelprint#1{%
   \noindent
   \hbox to \biblabelwidth{%
      \biblabelprecontents
      \biblabelcontents{#1}%
      \biblabelpostcontents
   }%
   \kern\biblabelextraspace
}%
%
%
%
\def\biblabelcontents#1{{\bblrm [#1]}}%
%
%
\def\bblrm{\rm}%
%
%
\def\bblem{\it}%
%
%
\def\bblsc{\ifx\@scfont\@undefined
              \font\@scfont = cmcsc10
           \fi
           \@scfont
}%
%
%
\def\bblnewblock{\hskip .11em plus .33em minus .07em }%
%
%
\let\bblhook = \empty
%
%
%
\def\printcitestart{[}
\def\printcitefinish{]}
\def\printbetweencitations{, }
\def\printcitenote#1{, #1}
%
%
%
\let\citation = \@gobble
%
%
%
\@innernewcount\@numparams
%
%
\def\newcommand#1{%
   \def\@commandname{#1}%
   \@getoptionalarg\@continuenewcommand
}%
%
%
\def\@continuenewcommand{%
   \@numparams = \ifx\@optionalarg\empty 0\else\@optionalarg \fi \relax
   \@newcommand
}%
%
%
\def\@newcommand#1{%
   \def\@startdef{\expandafter\edef\@commandname}%
   \ifnum\@numparams=0
      \let\@paramdef = \empty
   \else
      \ifnum\@numparams>9
         \errmessage{\the\@numparams\space is too many parameters}%
      \else
         \ifnum\@numparams<0
            \errmessage{\the\@numparams\space is too few parameters}%
         \else
            \edef\@paramdef{%
               \ifcase\@numparams
                  \empty  No arguments.
               \or ####1%
               \or ####1####2%
               \or ####1####2####3%
               \or ####1####2####3####4%
               \or ####1####2####3####4####5%
               \or ####1####2####3####4####5####6%
               \or ####1####2####3####4####5####6####7%
               \or ####1####2####3####4####5####6####7####8%
               \or ####1####2####3####4####5####6####7####8####9%
               \fi
            }%
         \fi
      \fi
   \fi
   \expandafter\@startdef\@paramdef{#1}%
}%
%
%
%
%
\def\@readauxfile{%
   \if@auxfiledone \else 
      \global\@auxfiledonetrue
      \@testfileexistence{aux}%
      \if@fileexists
         \begingroup
            \endlinechar = -1
            \catcode`@ = 11
            \input \jobname.aux
         \endgroup
      \else
         \message{\@undefinedmessage}%
         \global\@citewarningfalse
      \fi
      \immediate\openout\@auxfile = \jobname.aux
   \fi
}%
%
%
\newif\if@auxfiledone
\ifx\noauxfile\@undefined \else \@auxfiledonetrue\fi
%
%
%
%
\@innernewwrite\@auxfile
\def\@writeaux#1{\ifx\noauxfile\@undefined \write\@auxfile{#1}\fi}%
%
%
%
\ifx\@undefinedmessage\@undefined
   \def\@undefinedmessage{No .aux file; I won't give you warnings about
                          undefined citations.}%
\fi
%
%
\@innernewif\if@citewarning
\ifx\noauxfile\@undefined \@citewarningtrue\fi
%
%
%
\catcode`@ = \@oldatcatcode

\def\pfeilso{\leavevmode
            \vrule width 1pt height9pt depth 0pt\relax
           \vrule width 1pt height8.7pt depth 0pt\relax
           \vrule width 1pt height8.3pt depth 0pt\relax
           \vrule width 1pt height8.0pt depth 0pt\relax
           \vrule width 1pt height7.7pt depth 0pt\relax
            \vrule width 1pt height7.3pt depth 0pt\relax
            \vrule width 1pt height7.0pt depth 0pt\relax
            \vrule width 1pt height6.7pt depth 0pt\relax
            \vrule width 1pt height6.3pt depth 0pt\relax
            \vrule width 1pt height6.0pt depth 0pt\relax
            \vrule width 1pt height5.7pt depth 0pt\relax
            \vrule width 1pt height5.3pt depth 0pt\relax
            \vrule width 1pt height5.0pt depth 0pt\relax
            \vrule width 1pt height4.7pt depth 0pt\relax
            \vrule width 1pt height4.3pt depth 0pt\relax
            \vrule width 1pt height4.0pt depth 0pt\relax
            \vrule width 1pt height3.7pt depth 0pt\relax
            \vrule width 1pt height3.3pt depth 0pt\relax
            \vrule width 1pt height3.0pt depth 0pt\relax
            \vrule width 1pt height2.7pt depth 0pt\relax
            \vrule width 1pt height2.3pt depth 0pt\relax
            \vrule width 1pt height2.0pt depth 0pt\relax
            \vrule width 1pt height1.7pt depth 0pt\relax
            \vrule width 1pt height1.3pt depth 0pt\relax
            \vrule width 1pt height1.0pt depth 0pt\relax
            \vrule width 1pt height0.7pt depth 0pt\relax
            \vrule width 1pt height0.3pt depth 0pt\relax}

\def\pfeilsw{ \leavevmode 
            \vrule width 1pt height0.3pt depth 0pt\relax
            \vrule width 1pt height0.7pt depth 0pt\relax
            \vrule width 1pt height1.0pt depth 0pt\relax
            \vrule width 1pt height1.3pt depth 0pt\relax
            \vrule width 1pt height1.7pt depth 0pt\relax
            \vrule width 1pt height2.0pt depth 0pt\relax
            \vrule width 1pt height2.3pt depth 0pt\relax
            \vrule width 1pt height2.7pt depth 0pt\relax
            \vrule width 1pt height3.0pt depth 0pt\relax
            \vrule width 1pt height3.3pt depth 0pt\relax
            \vrule width 1pt height3.7pt depth 0pt\relax
            \vrule width 1pt height4.0pt depth 0pt\relax
            \vrule width 1pt height4.3pt depth 0pt\relax
            \vrule width 1pt height4.7pt depth 0pt\relax
            \vrule width 1pt height5.0pt depth 0pt\relax
            \vrule width 1pt height5.3pt depth 0pt\relax
            \vrule width 1pt height5.7pt depth 0pt\relax
            \vrule width 1pt height6.0pt depth 0pt\relax
            \vrule width 1pt height6.3pt depth 0pt\relax
            \vrule width 1pt height6.7pt depth 0pt\relax
            \vrule width 1pt height7.0pt depth 0pt\relax
            \vrule width 1pt height7.3pt depth 0pt\relax
            \vrule width 1pt height7.7pt depth 0pt\relax
            \vrule width 1pt height8.0pt depth 0pt\relax
            \vrule width 1pt height8.3pt depth 0pt\relax
            \vrule width 1pt height8.7pt depth 0pt\relax
            \vrule width 1pt height9pt depth 0pt\relax
      }


\def\widestnumber#1#2{}

\def\citewarning#1{\ifx\shlhetal\relax 
    \else
    \par{#1}\par
    \fi
}

\def\rm{\fam0 \tenrm}

\def\fakesubhead#1\endsubhead{\bigskip\noindent{\bf#1}\par}



%
%
%

%

\font\textrsfs=rsfs10
\font\scriptrsfs=rsfs7
\font\scriptscriptrsfs=rsfs5

\newfam\rsfsfam
\textfont\rsfsfam=\textrsfs
\scriptfont\rsfsfam=\scriptrsfs
\scriptscriptfont\rsfsfam=\scriptscriptrsfs

\edef\oldcatcodeofat{\the\catcode`\@}
\catcode`\@11

\def\Cal@@#1{\noaccents@ \fam \rsfsfam #1}

\catcode`\@\oldcatcodeofat


\expandafter\ifx \csname margininit\endcsname \relax\else\margininit\fi

\long\def\red#1\endred{}
\long\def\green#1\endgreen{}
\long\def\blue#1\endblue{}

\def\endred{ \unmatched endred! }
\def\endgreen{ \unmatched endgreen! }
\def\endblue{ \unmatched endblue! }

\ifx\latexcolors\undefinedcs\def\latexcolors{}\fi

\def\emptycs{}
\def\evaluatelatexcolors{%
        \ifx\latexcolors\emptycs\else
        \expandafter\xxevaluate\latexcolors\xxfertig\evaluatelatexcolors\fi}
\def\xxevaluate#1,#2\xxfertig{\setupthiscolor{#1}%
        \def\latexcolors{#2}}

\font\smallfont=cmsl7
\def\rutgerscolor{\ifmmode\else\endgraf\fi\smallfont
\advance\leftskip0.5cm\relax}
\def\setupthiscolor#1{\edef\tmptmpcs{\noexpand\bgroup\noexpand\rutgerscolor
\noexpand\def\noexpand\currentcolor{#1}%
\noexpand}%
\expandafter\let\csname#1\endcsname\tmptmpcs
\def\tmptmpcs{\checkColorUnmatched{#1}\popthecolor}
\expandafter\let\csname end#1\endcsname\tmptmpcs}

\def\checkColorUnmatched#1{\def\expectcolor{#1}%
    \ifx\expectcolor\currentcolor   
    \else \edef\failhere{\noexpand\tryingToClose '\currentcolor' with end\expectcolor}\failhere\fi}

\def\currentcolor{???}

\def\popthecolor{\ifmmode\else\endgraf\fi\egroup}

\expandafter\def\csname#1\endcsname{}

\evaluatelatexcolors

 \let\outerhead\head
 \def\head{\innerhead}
 \let\innerhead\outerhead

 \let\outersubhead\subhead
 \def\subhead{\innersubhead}
 \let\innersubhead\outersubhead

 \let\outersubsubhead\subsubhead
 \def\subsubhead{\innersubsubhead}
 \let\innersubsubhead\outersubsubhead

 \def\proclaim{\innerproclaim}
 \let\innerproclaim\outerproclaim

 %
 %
 %
 %

\def\demo#1{\medskip\noindent{\it #1.\/}}
\def\enddemo{\smallskip}

\def\remark#1{\medskip\noindent{\it #1.\/}}
\def\endremark{\smallskip}

\pageheight{8.5truein} 
\topmatter
\title{Toward classification theory of good $\lambda$-frames and
abstract elementary classes} \endtitle
\rightheadtext{Classification theory of frames and classes}
\author {Saharon Shelah \thanks {\null\newline I would like to thank 
Alice Leonhardt for the beautiful typing \null\newline
This research was supported by The Israel Science Foundation. 
Publication 705. \null\newline
} \endthanks} \endauthor 

\affil{The Hebrew University of Jerusalem\\
Einstein Institute of Mathematics\\
Edmond J. Safra Campus, Givat Ram \\
Jerusalem 91904, Israel
 \medskip
 Department of Mathematics \\
 Hill Center-Busch Campus \\
Rutgers, The State University of New Jersey \\ 
110 Frelinghuysen Road \\
Piscataway, NJ 08854-8019  USA} \endaffil
\endtopmatter
\document  
 
\pretags{300a,300b,300c,300d,300e,300f,300g,88,600}

\newpage

\head {\S0 Introduction} \endhead  \resetall \sectno=0
 \spuriousreset
\bn
Now $\lambda$-good frame is for us a
parallel of the class of models of a superstable theory.  Our main
line is to start with $\lambda$-good$^+$
frame ${\frak s}$, categorical in $\lambda,n$-successful for $n$ large 
enough and try to have parallel
of stability theory for ${\frak K}_{{\frak s}(+ \ell)}$ for $\ell < n$
not too large.  Characteristically from time to time we have to increase $n$
relative to $\ell$
to get our desirable properties; we do not critically mind the exact $n$,
so you can think of an $\omega$-successful ${\frak s}$.  Usually each claim
or definition is for a fixed ${\frak s}$, assumed to be successful
enough.  So using assumptions on $\lambda^{+2}$ rather than
$\lambda^{++}$ is not so crucial now.  \nl
But a postriori we are interested in the model theory of such classes
${\frak K}_{\frak s}$ per-se, and see as a test for this theory,
that in the $\omega$-successful case we can understand also the model in
higher cardinals, e.g., prove that ${\frak K}^{\frak s}_\mu \ne
\emptyset$ for every $\mu \ge \lambda$.  Recall there are reasonable
$\lambda$-frames which are not $n$-excellent but still we can say alot
on models in ${\frak K}_{{\frak s}(+ \ell)}$ for $\ell < n$. \nl
Moving from $\lambda$ to $\lambda^+$ we would have preferred 
not to restrict ourselves to saturated models but at present we do not
know it.  However, in the $\omega$-excellent case we can understand
the class of $\lambda^{+ \omega}$-saturated models in ${\frak K}_{\frak s}$,
i.e., ${\frak K}^{{\frak s}(+ \omega)}$.  This fits well the thesis
that it is reasonable to first analyze the quite saturated case. \nl
Why are we interested in ${\frak K}_{\frak s}$ (\scite{705-0.2A}(1))? we
can ``blow it up" by \sectioncite[\S1]{600} but for good frames this is not
so. 
\nl
Concerning the framework note that the unidimensional (or just non
multi-dimensional) case is easier.  In the characteristic unidimensional 
case, each $p \in {\Cal S}^{\text{bs}}(M)$ is minimal and any $p,q \in
{\Cal S}^{\text{bs}}(M)$ are not orthogonal.  In the characteristic non
multi-dimensional case for any $M \in K_\lambda,{\Cal S}^{\text{bs}}(M)$
contains up to nonorthogonality every $p \in {\Cal S}^{\text{bs}}(N),
M \le_{\frak K} N \in K_\lambda$. \nl
Generally the unidimensional case is easiest and 
is enough to continue \cite{Sh:576}, and to deal with categoricity.
\nl
A drawback in \sectioncite[\S5]{600} is that we need to
assume that the normal ideal WDmId$(\lambda^+)$ is not
$\lambda^{++}$-saturated.  This will be improved in \cite{Sh:F603}; it
is easier to do it when we have the theory developed here.
\bn
\margintag{705-0.A}\ub{\stag{705-0.A} Notation}:  Let ${\frak s}$ denote a good frame
(usually) or just a pre-frame, that is
\definition{\stag{705-0.2A} Definition}  1) We say ${\frak s}$ is a
pre-frame (or pre $\lambda$-frame) if
${\frak s} = ({\frak K}_{\frak s},{\Cal S}^{\text{bs}}_{\frak s},
\nonfork{}{}_{\frak s})$ with ${\frak K}_{\frak s} = 
{\frak K}({\frak s})$ a $\lambda_{\frak s}$-a.e.c., 
${\Cal S}^{\text{bs}}[{\frak s}] =
{\Cal S}^{bs}_{\frak s},\nonfork{}{}_{}[{\frak s}] = 
\nonfork{}{}_{\frak s},\le_{\frak s} = \le_{{\frak K}_{\frak s}}$ and
satisfying axioms (A), (C), (D)(a), (E)(a)(b).  We say ${\frak s}$ is
a frame if it satisfies axioms
(A), (B), (C), (D)(a),(b), (E)(a),(b)
from \sectioncite[\S2]{600}.  Recall that ${\frak s}$ is a good frame if it
satisfies all the axioms there.
\nl
2) For a frame ${\frak s}$ let 
${\frak K}^{\frak s} = {\frak K}[{\frak s}]$ 
be the a.e.c. derived from ${\frak K}_{\frak s}$ and ${\frak K}^{\frak
s}_\mu = ({\frak K}[{\frak s}])_\mu$ so ${\frak K}_{\frak s} 
= {\frak K}^{\frak s}_{\lambda({\frak s})}$.  Recall that 
if ${\frak K}$ is a $\lambda$-a.c.e.,
\ub{then} the a.e.c.-derived from it, ${\frak K}^{\text{up}}$ 
is the unique a.e.c. 
${\frak K}'$ with $\tau({\frak K}') = \tau({\frak K}),
\text{LS}({\frak K}') = \lambda,
{\frak K}'_\lambda = {\frak K}_{\frak s}$, see \sectioncite[\S1]{600}. \nl
3) For a frame ${\frak s}$ let $\le^{\frak s} = \le^{\frak s}_{\frak
K}$ be $\le_{{\frak K}^{\frak s}}$.
\enddefinition
\bn
\margintag{705-0.2B}\ub{\stag{705-0.2B} Convention}:  For simplicity 
we assume ${\frak K}$ is such that if $\bar a \in {}^{\omega >}M,
M \in K$ then $\bar a$ can be considered an element of $M$.
\bigskip

\definition{\stag{705-0.X} Definition}  Let ${\frak s}$ be a good
frame and $\mu \ge \lambda_{\frak s}$. \nl
1) Let ${\frak s} \langle \mu \rangle
= ({\frak K}^{\frak s}_\mu,{\Cal S}^{\text{bs}}_{{\frak s},\mu},
\nonfork{}{}_{{\frak s},\mu})$ with ${\Cal S}^{\text{bs}}_{{\frak
s},\mu},\nonfork{}{}_{{\frak s},\mu}$
as defined in \sectioncite[\S2]{600}; also ${\Cal S}^{\text{bs}}_{{\frak s},<
\mu},{\Cal S}^{\text{bs}}_{{\frak s},< \infty},\nonfork{}{}_{{\frak s},<
\mu}$ and $\nonfork{}{}_{{\frak s},< \infty}$ are from there. \nl
2) Let ${\frak s}(\mu) =: ({\frak K}_{{\frak s}(\mu)},{\Cal
S}^{\text{bs}}_{{\frak s},\mu},\nonfork{}{}_{{\frak s}(\mu)})$ where
$K_{{\frak s}(\mu)} =: \{M \in K^{\frak s}_\mu:M$ is superlimit in
${\frak K}^{\frak s}_\mu\},\le_{K_{{\frak s}(\mu)}} = \le_{\frak
K}^{\frak s} \restriction K_{{\frak s}(\mu)}$,
and of course ${\Cal S}^{\text{bs}}_{{\frak s}(\mu)} 
= {\Cal S}^{\text{bs}}_{{\frak s},\mu} \restriction K_{{\frak
s}(\mu)},\nonfork{}{}_{{\frak s},\mu} 
= \nonfork{}{}_{} \restriction K^{{\frak s},s}_\mu$.
\nl
3) Let ${\frak s}(\mu) = ({\frak K}_{{\frak s}[\mu]},
\nonfork{}{}_{{\frak s}(mu)})$ where $K_{{\frak s}[\mu]} = \{M \in
K^{\frak s}_\mu:\text{ if } N \le_{{\frak K}[{\frak s}]} M,N \in
K_{\frak s}$ and $p \in {\Cal S}^{\text{bs}}_{\frak s}(N)$ then we can
find $\bold I \subseteq M \backslash N$ such that, if $N \le_{\frak s}
N_1 \le_{{\frak K}[{\frak s}]} M$ then $N_1 \cap \bold I$ is
independent in $(N,N_1)$. \nl
4) If ${\frak s}$ is $\omega$-successful let ${\frak s}^{+ \omega} =
{\frak s}(+ \omega)$ is ${\frak s}[\omega]$.
\enddefinition
\bigskip

\remark{\stag{705-0.xA} Remark}  1) In \S12, ${\frak K}_{{\frak s}[\mu]}$
under strong assumptions is proved to be good (and categorical in
$\mu$).
\sn
Note, if ${\frak K}_\mu$ have the JEP \ub{then} we can show that
any two superlimit models $M \in K^{\frak s}_\mu$ are isomorphic but a
priori not in general.  However, here there is a canonical way to
proceed: if ${\frak K}$ has superlimit model $M_0$ in $\lambda$ we can
try to choose by induction on $\ell < \omega$ a superlimit model
$M_\ell$ in $K_{\lambda^{+\ell}}$ such that $M_{\ell +1}$ is a model
in ${\frak K}_{[M_\ell]}$ which has cardinality $\lambda^{+\ell +1}$
and is saturated in ${\frak K}_{[M_\ell]}$.  It is not clear a priori
if $M_{\ell +1}$ exists but at least it is unique (we can continue in
$\lambda^{+\alpha}$ for ordinals with more case). 
\nl
3) Let $<^+_{\frak s}$ be the following two place relation on
$K_{\frak s}:M <^+_{\frak s}$ iff $M \le_{\frak s} N$ and $N$ is
${\frak K}_{\frak s}$-universal over $M$. \nl
4) Note that ${\frak K}_{{\frak s}[\mu]}$ includes ${\frak K}_{{\frak
s}[\mu]}$ and ${\frak K}_{{\frak s}[\mu]}$ but it is not clear how to
compare ${\frak K}_{{\frak s}[\mu]},{\frak K}_{{\frak s}[\mu]}$.
\endremark
\bigskip

\proclaim{\stag{705-0.xB} Claim}  If ${\frak s}$ is $\omega$-successful
then ${\frak s}^{+ \omega}$ is a good $\lambda^{+ \omega}$-frame
categorical in $(\lambda^{+ \omega})$.
\endproclaim
\newpage

\head {\S1 Good$^+$ frames; basics} \endhead  \resetall \sectno=1
 \spuriousreset
\bn
In \marginbf{!!}{\cprefix{600}.\scite{600-rg.7}}(4) there was what may look like a minor
drawback: moving from $\lambda$ to $\lambda^+$
the derived class not
only have fewer models of cardinality $\ge \lambda^+$ but also the notion
of being a submodel changes; this is fine there, and surely unavoidable in
some circumstances.  More specifically, for proving the main theorem, 
it was enough to move from ${\frak s}$ to a 
good frame ${\frak t}$ satisfying $\lambda_{\frak t} = \lambda^+_{\frak
s},\lambda_{\frak s} < \mu < \lambda^{+ \omega} \Rightarrow
I(\lambda,K^{\frak t}) \le I(\lambda,K^{\frak s})$ and forget ${\frak s}$.
But for us now this is undesirable (as arriving to $\lambda^{+ \omega}$ we
have forgotten everything)  and toward this we consider a (quite mild)
strengthening of $\lambda^+$-good.
\bigskip

\demo{\stag{705-stg.0} Hypothesis}  ${\frak s} = ({\frak K}_\lambda,
{\Cal S}^{\text{bs}},\nonfork{}{}_{\lambda})$ is a good $\lambda$- frame.
\enddemo
\bigskip

\definition{\stag{705-stg.0A} Definition}  1) We say 
${\frak s}$ is successful \ub{if} the conclusions 
under ``no nonstructure assumptions in $\lambda^{++}$" hold, that is:
\mr
\item "{$(*)(a)$}"  $K^{3,\text{uq}}_\lambda$ is dense;
i.e. for every $M \in K_\lambda$ and $p \in {\Cal S}^{\text{bs}}(M)$ 
there is
\nl
$(M,N,a) \in K^{3,\text{uq}}_\lambda$ such that tp$(a,M,N) = p$ 
\sn
\item "{$(b)$}"   if $\langle N_i:i \le \delta \rangle$ is
$\le^*_{\lambda^+}[{\frak s}]$-increasing continuous 
in $K^{\text{nice}}_{\lambda^+}[{\frak s}]$ and \nl
$i < \delta \Rightarrow N_i \le^*_{\lambda^+}
N \in K^{\text{nice}}_{\lambda^+}$ \ub{then} 
$N_\delta \le^*_{\lambda^+} N$ (see \marginbf{!!}{\cprefix{600}.\scite{600-rg.7}}(1)).
\ermn
2) We say (the $\lambda$-good frame) 
${\frak s}$ is weakly successful \ub{if}
clause (a) of $(*)$ above holds.  
\enddefinition
\bn
Usually at least ``${\frak s}$ is weakly successful" is used, but
sometimes less suffice (and is helpful though not crucial).
\definition{\stag{705-1.2A} Definition}  1) Let $K^{3,\text{puq}}_{\frak
s}$ be 
the class of $(M,N,a) \in K^{3,\text{bs}}_{\frak s}$ such that:
if $(M,N,a) \le_{\text{bs}} (M_1,N_1,a)$ and $b \in \bold I_{M,N_1}$
then tp$(b,N,N_1)$ does not fork over $M$.  We may say $(M,N,a)$ has
pseudo uniqueness. \nl
2) We say that ${\frak s}$ has weakly pseudo uniqueness if
$K^{3,\text{puq}}_{\frak s}$ is dense in $(K^{3,\text{bs}}_{\frak
s},\le_{\text{bs}})$.
\nl
3) We say that ${\frak s}$ has semi-pseudo uniqueness if $\alpha <
\lambda^+,\langle M_i:i \le \alpha \rangle$ is $\le_{\frak
s}$-increasing continuous, $(M_i,M_{i+1},a_i) \in
K^{3,\text{bs}}_{\frak s}$ then for some $\le_{\frak s}$-increasing
continuous $\langle N_i:i \le \alpha \rangle$ we have
$(M_i,M_{i+1},a_i) \le_{\text{bs}} (N_i,N_{i+1},a_i) \in
K^{3,\text{puq}}_{\frak s}$.
\enddefinition
\bigskip

\remark{Remark}  For successful ${\frak s}$ we can define
a successor, ${\frak s}^+ = {\frak s}(+)$, 
a $\lambda^+$-good frame (see
\scite{705-stg.3A} below), but not with the most desirable
$\le_{{\frak K}_{{\frak s}(+)}}$, for rectifying this 
we consider below good$^+$
frames.  Together with locality of types for models in
${\frak K}^{{\frak s}(+)}_{\lambda^+_{\frak s}}$ this seems 
to be in the
right direction.  Less crucial, still worthwhile, is that ${\frak s}'
= {\frak s}^+$ has a strong property we call saturative, which can be
used in several cases as an alternative of another version which 
says ${\frak s}'$ has the form ${\frak s}^+$ with ${\frak s}$ a 
successful good$^+$ frame.
\endremark
\bigskip

\definition{\stag{705-stg.1} Definition}  1) We say that ${\frak s} = 
({\frak K}_\lambda,{\Cal S}^{\text{bs}},\nonfork{}{}_{\lambda}) =
(K_\lambda,\le_{{\frak K}_\lambda},{\Cal S}^{\text{bs}},
\nonfork{}{}_{\frak s})$ is a good$^+ \lambda$-frame or 
good$^1 \lambda$-frame \ub{if}:
\mr
\item "{$(a)$}"  ${\frak s}$ is a good $\lambda$-frame
\sn
\item "{$(b)$}"  the following is impossible
{\roster
\itemitem{ $(*)$ }   $\langle M_i:i < \lambda^+ \rangle$ is
$\le_{\frak s}$-increasing continuous in $K_\lambda$ and $\langle N_i:i <
\lambda^+ \rangle$ is $\le_{\frak s}$-increasing continuous in $K_\lambda$
and $M_i \le_{\frak s} N_i,p^* \in {\Cal S}^{\text{bs}}_{\frak s}(M_0)$ and for each $i <
\lambda^+$ we have: \nl
$a_{i+1} \in M_{i+2} \backslash M_{i+1}$, tp$_{\frak s}(a_{i+1},M_{i+1},M_{i+2})$ is
a nonforking extension of $p^*$ but tp$(a_{i+1},N_0,N_{i+1})$ is
not. \nl
We then say $\langle M_i,N_i,a_i:i < \lambda^+ \rangle$ is a
counterexample (ignoring $a_i$ being defined only for successor $i$;
we could demand this for every $i$ by monotonicity of nonforking).
\endroster}
\ermn
2) We say a good $\lambda$-frame ${\frak s}$ is saturative \ub{if}:
\mr
\item "{$(a)$}"  every $M \in {\frak K}_{\frak s}$ is
$(\lambda,*)$-brimmed (for ${\frak K}_{\frak s}$), equivalently: 
is superlimit
\sn
\item "{$(b)$}"  if $M_0 \le_{\frak s} M_1 \le_{\frak s} M_2$ and
$M_1$ is $(\lambda,*)$-brimmed over $M_0$ \ub{then} $M_2$ is
$(\lambda,*)$-brimmed over $M_0$.
\ermn
The ``${\frak s}$ is saturative" is a relative of
``non-multidimensional", but be careful, see \scite{705-stg.2}(3) below. \nl
Well, do we lose much by adopting the good$^+$ version?  
First, are the old cases
covered?  Yes, by the following claim
\enddefinition
\bigskip

\proclaim{\stag{705-stg.2} Claim}  1) In \sectioncite[\S3]{600} all the cases where we
prove ``good $\lambda$-frame" we actually get ``good$^+ \lambda$-frames. \nl
2) In fact those frames are also saturative. \nl
3) If $T$ is a complete superstable first order theory stable in
$\lambda$ and $\kappa \le \lambda$ or $\kappa = \aleph_\varepsilon$
(in an abuse of notation) or
$\kappa = 0,\lambda \ge |T|$ and $[\kappa > 0 \Rightarrow T \text{
stable in } \lambda]$ and 
${\frak s} = {\frak s}^\kappa_{T,\lambda}$ (so 
$K_{\frak s} = (\{M:M \models T,\|M\| \ge \lambda$ 
and $M$ is $\kappa$-brimmed$\},\prec)$, that is
${\frak s}$ is defined in \marginbf{!!}{\cprefix{600}.\scite{600-Ex.0}}, \ub{then}
\mr
\widestnumber\item{$(iii)$}
\item "{$(i)$}"  ${\frak s}$ is a good$^+ \lambda$-frame
\sn
\item "{$(ii)$}"  assume $0 < \kappa < \lambda$, \ub{then}: ${\frak s}$ is
saturative iff $T$ is non-multidimensional (see \scite{705-stg.5}(5))
\sn
\item "{$(iii)$}"  if ${\frak s}' = {\frak s}^\kappa_{T,\lambda}[M]$
where $M \in {\frak K}_\lambda$ is the superlimit model (i.e., the
saturated one), \ub{then} ${\frak s}'$ is saturative
\sn
\item "{$(iv)$}"  if $\kappa = \lambda,{\frak s}$ is saturative.
\endroster
\endproclaim
\bigskip

\demo{Proof}  1), 2).
\enddemo
\bn
\ub{Case 1}:  \ub{Concerning}  Claim \marginbf{!!}{\cprefix{600}.\scite{600-Ex.4}}.

So $2^\lambda < 2^{\lambda^+} < 2^{\lambda^{++}},{\frak K}$ is an abstract
elementary class categorical in $\lambda,\lambda^+$ and 
$1 \le I(\lambda^{++},K) <
2^{\lambda^{++}}$, with LS$({\frak K}) \le \lambda$,  
WDmId$(\lambda^+)$ is not $\lambda^+$-brimmed (or
just a model theoretic consequence).  We let 
$\lambda_{\frak s} = \lambda^+,{\frak
K}_{\frak s} = {\frak K}_{\lambda^+}$ and for $M \in K_{\frak s}$ we  
let ${\Cal S}^{\text{bs}}(M) = \{p \in {\Cal S}(M):p$ is not algebraic and
for some $M_0 \le_{\frak K} M,M_0 \in K_\lambda$ and $p \restriction
M_0$ is minimal$\}$ and
$\nonfork{}{}_{}(M_0,M_1,a,M_3)$ iff $M_0 \le_{\frak K} M_1
\le_{\frak K} M_2$ are in $K_{\frak s}$ and $a \in M_3 \backslash M_1$
and for some $M'_0 \le_{\frak K} M_0$ from $K_{\frak s}$
the type tp$(a,M_0,M_3) \in {\Cal S}^{\text{bs}}(M)$ is minimal.

So let $\langle (M_i,N_i,a_i):i < \lambda^+_{\frak s} 
\rangle$ be as in $(*)$ of 
clause (b) of Definition \scite{705-stg.1}.  So for success 
$i < \lambda^+_{\frak s}$, tp$_{\frak s}(a_i,M_0,N_{i+1}) 
\in {\Cal S}^{bs}_{\frak s}(M_0)$ so is minimal while
tp$_{\frak s}(a_i,N_0,N_{i+1})$ is not its nonforking extension, hence
necessarily $i < \lambda^+_{\frak s} \Rightarrow a_{i+1} 
\in N_0 \backslash M_0$, so $\langle a_{i+1}:i < \lambda^+_{\frak s}
\rangle$ is a sequence with no repetitions of members from 
$N_0 \backslash M_0$ while $N_0 \in K_{\lambda_{\frak s}}$ so $\|N_0\|
= \lambda_{\frak s}$, contradiction.  Also saturativity should be clear.
\bn
\ub{Case 2}:  Claim \marginbf{!!}{\cprefix{600}.\scite{600-Ex.1}}; actually from \cite{88r}.

So let $\langle (M_i,N_i,a_i):i < \omega_1 \rangle$ be as in clause (b) of
\scite{705-stg.1}.  So there is a finite $A_0 \subseteq M_0$ such that
tp$(a_{i+1},M_{i+1},M_{i+2})$ is definable over $A_0$ (see
\marginbf{!!}{\cprefix{88}.\scite{88-5.11}}), but tp$(a_{i+1},N_0,N_{i+2})$ does not have the same
definition hence it splits over $A_0$.  By \marginbf{!!}{\cprefix{88}.\scite{88-5.11}} for a
club of $E$ of $\omega_1$ we have
\mr
\item "{$\boxtimes$}"  $\delta \in E \and \delta < \alpha < \omega_1 \and
\bar a \in N_\delta \Rightarrow \text{ tp}(\bar a,M_\alpha,N_\alpha)$ is
definable over some finite $B_{\bar a} \subseteq M_\delta$.
\ermn
We get a contradiction to the symmetry (i.e., 
(E)(f) proved in the proof
of \marginbf{!!}{\cprefix{600}.\scite{600-Ex.1t}} that is by \marginbf{!!}{\cprefix{88}.\scite{88-5.20}}).  Also
saturatively should be clear.
\bn
\ub{Case 3}:  Claim \marginbf{!!}{\cprefix{600}.\scite{600-Ex.1A}}; actually from \cite{Sh:48}.

Similar to case 2.  \nl
3) First we prove clause (b) of \scite{705-stg.1}(1), i.e. we prove
${\frak s}$ is a good$^+$ frame; so assume toward contradiction that
$\langle M_i:i < \lambda^+ \rangle,\langle N_i:i < \lambda^+
\rangle,p^*$ and $a_{i+1}$ for $i < \lambda^+$ are as in $(*)$ of
\scite{705-stg.1}(2) clause (b).  Let $M = \cup\{M_i:i < \lambda^+\}$ and
$N = \cup\{N_i:i < \lambda^+\}$. Now for every finite sequence $\bar
c$ from $N_0$, there is $i_c < \lambda^+$ such that tp$(\bar c,M,N)$
does not fork over $M_{i_{\bar c}}$, and let $i^* = \sup\{i_{\bar
c}:\bar c \in {}^{\omega >}(N_0)\}$ so $i^* < \lambda^+$ and easily $i
\in [i^*,\lambda^+) \and \bar c \in {}^{\omega >}(N_0) \Rightarrow \text{ tp}(\bar
c,M_{i+2},N)$ does not fork over $M_{i+1}$, hence by symmetry and
local character (\cite[III,\S0]{Sh:c}) we have tp$(\bar a_{i+1},M_{i+1}
\cup N_0,N)$ does not fork over $M_i$ hence (transitivity) over $M_0$,
contradiction.  So clause (i) holds.

As for saturativeness in clauses $(iii),(iv)$ in \scite{705-stg.1}(2), to show
that the model $(M_2,c)_{c \in M_0}$ is saturated it is enough to show
that for every $A \subseteq M_2,|A| < \aleph_0$
and regular $p \in {\Cal S}^1(A \cup M_0)$, we have dim$(p,M_2) =
\lambda$.  
Why this holds?  If $p \perp M_0$ then we can find a regular $q \in
{\Cal S}(M_0),q \perp p$ and as $M_1$ is $(\lambda,*)$-brimmed over 
$M_0$, dim$(q,M_2) 
\ge \text{ dim}(q,M_1) = \lambda$ and easily dim$(p,M_2) = \text{
dim}(q,M_2)$.  If $p \perp M_0$ see \cite{Sh:225}.  The proof of
clause (ii) is easy too (for the case $\kappa =0$ when for some
$M,Th(M,c)_{c \in M}$ is categorical in $\lambda^+$, see \cite{Sh:c}
and properties as in \cite{ShHM:158} and the analysis of Laskowski of
model of $T$ in $\lambda = |T|$ when $T$ is categorical in $\lambda^+$). \nl
${{}}$  \hfill$\square_{\scite{705-stg.2}}$\margincite{705-stg.2}
\bn
Also in the main result of \chaptercite{600} we get good$^+$
\proclaim{\stag{705-stg.2A} Claim}  If ${\frak s}$ is a successful
$\lambda$-good frame, \ub{then} ${\frak s}(+)$ is a
$\lambda^+$-good$^+$ frame and is saturative.
\endproclaim
\bigskip

\demo{Proof}  Like Case 1 of the proof of Claim \scite{705-stg.2}(1).
\enddemo
\bigskip

\proclaim{\stag{705-stg.3} Goodness Plus Claim}  Assume that ${\frak s} =
({\frak K}_\lambda,{\Cal S}^{\text{bs}},\nonfork{}{}_{\lambda})$ is a
good$^+ \lambda$-frame.  \ub{Then}:
\mr
\item "{$(a)$}"  if $M^*_1 \le_{\frak K} 
M^*_2$ are from $K^{\text{nice}}
_{\lambda^+}$ and $M^*_1 \nleq^*_{\lambda^+} M^*_2$ \ub{then}
$(**)_{M^*_1,M^*_2}$ from \marginbf{!!}{\cprefix{600}.\scite{600-sg.10}} holds
\sn
\item "{$(b)$}"  if $\boxtimes$ from \marginbf{!!}{\cprefix{600}.\scite{600-rg.7}} (which holds
if $I(\lambda^{++},K) < 2^{\lambda^{++}}$) \ub{then} 
{\roster
\itemitem{ $(\alpha)$ }  $\le^*_{\lambda^+},\le_{\frak K}$ agree on
$K^{\text{nice}}_{\lambda^+}$
\sn
\itemitem{ $(\beta)$ }  $(K^{\text{nice}}_{\lambda^+},\le_{{\frak K}_\lambda}
{\Cal S}^{\text{bs}}_{\lambda^+},\nonfork{}{}_{\lambda^+})$ 
(as defined there, called ${\frak s}^+$ below) is a
$\lambda^+$-good$^+$ frame.
\endroster}
\ermn
Recall in \chaptercite{600} we get a weak version of $(\alpha)$ of $(b)$:
$\le^*_\lambda,\le^\otimes_{\frak K}$ agree on
$K^{\text{nice}}_{\lambda^+}$.
\endproclaim
\bn
Before proving \scite{705-stg.3} we see a conclusion
\definition{\stag{705-stg.3A} Definition}  1) For a good $\lambda$-frame
${\frak s} = ({\frak K},{\Cal S}^{\text{bs}},\nonfork{}{}_{\lambda})$ define
${\frak s}^+ = {\frak s}(+)$, a $\lambda^+$-frame, as follows (so $\lambda
[{\frak s}^+] = \lambda^+$):
\mr
\item "{$(a)$}"  $K_{\lambda^+}[{\frak s}^+] =$ the class of
$\lambda^+$-saturated models from $K[{\frak s}]$
\sn
\item "{$(b)$}"  $\le_{{\frak K}[{\frak s}^+]} = \le^*_{\lambda^+} 
\restriction K_{\lambda^+}[{\frak s}^+]$
\sn
\item "{$(c)$}"  ${\Cal S}^{\text{bs}}[{\frak s}^+] = {\Cal S}^{\text{bs}}
_{{\frak s}(+)} = \{\text{tp}_{{\frak s}(+)}(a,M_1,M_2):
M_1 \le_{{\frak s}(+)} M_2$ are
from $K_{\lambda^+}[{\frak s}^+],a \in M_2 \backslash M_1$ and there is
$N_1 \le_{\frak s} M_1,N_1 \in K_\lambda$, such that $N_1 \le_{\frak
K} N \le_{\frak K} M_1 \and N \in K_{\frak s} 
\Rightarrow \text{ tp}_{\frak s}(a,N,M_2) \in {\Cal S}^{\text{bs}}_{\frak s}
(N)$ does not fork over $N_1$ (in ${\frak s}$'s sense)$\}$
\sn
\item "{$(d)$}"  $\nonfork{}{}_{{\frak s}(+)} = \{(M_0,M_1,a,M_2):M_0
\le_{{\frak s}(+)} M_1 \le_{{\frak s}(+)} M_2$ are
of cardinality $\lambda^+$, \nl

\hskip60pt $a \in M_2 \backslash M_1$ and
tp$_{{\frak s}(+)}(a,M_1,M_2) \in {\Cal S}^{\text{bs}}
_{{\frak s}(+)}(M_1)$  \nl

\hskip60pt has a witness $N_1 \le_{\frak K} M_0\}$.
\ermn
2) If $a,M_1,M_2,N_1$ are as in clause (c) then 
we call $N_1$ or tp$_{\frak s}(a,N_1,M_2)$, a witness for
tp$_{{\frak s}(+)}(a,M,N)$; we may abuse our notation and say that
tp$_{{\frak s}(+)}(a,M_1,M_2)$ 
does not fork over $N_1$.  Similarly for stationarization 
(= nonforking extension).
\enddefinition
\bigskip

\demo{\stag{705-stg.3B} Conclusion}  1) If ${\frak s}$ is a 
good$^+ \lambda$-frame and is 
successful (see Definition \scite{705-stg.0A}), \ub{then}
$\le_{{\frak s}(+)} = \le_{\frak K} \restriction K_{{\frak s}(+)}$ 
and ${\frak s}^+$ is a $\lambda^+$-good$^+$ frame. \nl
2) For $M_1,M_2 \in K^{\text{nice}}_{\lambda^+}[{\frak s}]$ 
we have $M_1 \le^+_{\lambda^+} M_2$ for ${\frak s}$ (see 
Definition \marginbf{!!}{\cprefix{600}.\scite{600-ne.1}}(3)) 
iff $M_2$ is $(\lambda^+,*)$-brimmed over $M_1$ for ${\frak s}^+$.
\enddemo
\bigskip

\demo{Proof}  1) By \sectioncite[\S7]{600} particularly the proof of
\marginbf{!!}{\cprefix{600}.\scite{600-rg.7}} and \scite{705-stg.3}. \nl
2) By the proof of \marginbf{!!}{\cprefix{600}.\scite{600-ne.4}}(2).   \hfill$\square_{\scite{705-stg.3B}}$\margincite{705-stg.3B}
\enddemo
\bn
We shall use this conclusion freely.
\demo{Proof of \scite{705-stg.3}}
\sn
\ub{Clause $(a)$}:

Assume that $(**)_{M^*_1,M^*_2}$ fails, 
\ub{then} by the assumptions of clause (a), from the
clauses of $(**)_{M^*_1,M^*_2}$ only clause (iv) there may fail.  So we 
can find $N^*_1 \le_{\frak s} N^*_2$ from 
$K_\lambda,N^*_\ell \le_{\frak K} M^*_\ell$ and $p \in {\Cal S}^{\text{bs}}_{\frak s}
(N^*_2)$ which does not fork over $N^*_1$ such that no $a \in M^*_1$ realizes $p$
in $M^*_2$.
Let $\langle M^\ell_\alpha:\alpha < \lambda^+ \rangle$ be a 
$\le_{\frak s}$-representation of $M^*_\ell$ for $\ell = 1,2$.  
Without loss of generality $N^*_1 \le_{\frak s} M^2_0$ and $M^2_\alpha
\cap M^*_1 = M^1_\alpha$ and as $M^*_1 \in K^{\text{nice}}_{\lambda^+}$ also
$M^1_{\alpha +1}$ is $(\lambda,*)$-brimmed over $M^1_\alpha$
for $\alpha < \lambda^+$.  For each $\alpha < \lambda^+$ the type $p \restriction
N^*_1 \in {\Cal S}^{bs}_{\frak s}(N^*_1)$ 
has a nonforking extension $p_\alpha \in {\Cal S}^{\text{bs}}_{\frak s}(M^1_\alpha)$.  
As $M^1_{\alpha +1}$ is $(\lambda,*)$-brimmed over $M^1_\alpha$
clearly for every $\alpha < \lambda^+$ there is $a_\alpha \in M^1_{\alpha +1}
\backslash M^1_\alpha$ realizing $p_\alpha$.

Let $M_\alpha =: M^1_\alpha,N_\alpha =: M^2_\alpha,p^* = p_0$; note that
$N^*_1 \le_{\frak s} M_0,N^*_2 \le_{\frak s} N_0$, so all the demands in
$(*)$ of clause (b) of Definition \scite{705-stg.1}(1) holds, in
particular $a_{i+1} \in M_{i+2} \backslash M_{i+1}$, tp$_{\frak
s}(a_{i+1},M_{i+1},M_{i+2}) = p_{i+1} \in {\Cal S}^{\text{bs}}_{\frak
s}(M_{i+1})$ is a nonforking extension of $p \restriction N^*_1$ hence
of $p^* =: p_0$ but if tp$(a_{i+1},N_0,N_{i+1})$ does not fork over
$M_0$ then it is also a non forking extension of $p$ (recall $N^*_2
\le_{\frak s} N_0$).  So this contradicts
``${\frak s}$ is $\lambda$-good$^+$", i.e., clause (b) of Definition
\scite{705-stg.1}, in other words, some $a_i$ realizes $p$ so actually
clause $(iv)$ of $(**)$ of \marginbf{!!}{\cprefix{600}.\scite{600-sg.10}} holds.
\mn
\ub{Clause $(b)$}:  \ub{Subclause $(\alpha)$}.

Straight by clause (a) and \marginbf{!!}{\cprefix{600}.\scite{600-rg.7}}.
\mn
\ub{Clause $(b)$}:  \ub{Subclause $(\beta)$}.

Recalling \marginbf{!!}{\cprefix{600}.\scite{600-rg.7}}, the only new point is the $+$ of
the good$^+$ (for ${\frak s}^+$) and the forking extension of $p$ in
${\Cal S}^{\text{bs}}_{{\frak s}(t)}(M_{i+1})$. \nl 
So assume $\langle (M_i,N_i,a_i):i < \lambda^{++} \rangle$ be a 
counterexample to the ``${\frak s}^+$ is a good$^+ \,\lambda^+$-frame".  
So in particular $M_i,N_i \in K_{{\frak s}(+)}$ and 
$p \in {\Cal S}^{\text{bs}}_{{\frak s}(t)}(M_0)$ and
$p_i = \text{ tp}(a_i,M_{i+1},M_{i+2})$ the type for 
${\frak K}^{\text{nice}}_{\lambda^+}$ of course, which by clause (a)
is $(K^{\text{nice}}_{\lambda^+},\le^*_{\frak K} \restriction
K^{\text{nice}}_{\lambda^+})$.  As 
$p =: p_i \restriction M_0 \in {\Cal S}^{\text{bs}}_{{\frak s}(+)}(M_0)$ 
there are $M' \le_{\frak K} M_0,M' \in K_\lambda$ and 
$q \in {\Cal S}^{\text{bs}}(M')$
which witness $p \in {\Cal S}^{\text{bs}}_{{\frak s}(+)}(M_0)$ 
(see \scite{705-stg.3A}).  Let $\langle N'_\varepsilon:
\varepsilon < \lambda^+ \rangle$ 
be a sequence which $\le_{\frak K}$-represents $N_0$.
For each $i < \lambda^{++}$, as $p'_i = 
\text{ tp}_{\frak s}(a_{i+1},N_0,N_{i+2})$ is not
a non-forking extension of $p$ necessarily there is $\varepsilon =
\varepsilon_i < \lambda^+$ such that $M' \le_{\frak K} N'_\varepsilon$ and
$p'_i \restriction N'_\varepsilon = \text{ tp}
(a_{i+1},N'_\varepsilon,
N_{i+2})$ is not a nonforking extension of $q$.  So
for some $\varepsilon < \lambda^+$ the set $B_\varepsilon = \{i <
\lambda^{++}:\varepsilon_i = \varepsilon\}$ is unbounded in $\lambda^{++}$.
We now choose by induction on $\zeta < \lambda^+$ a triple $(i_\zeta,
M^*_\zeta,N^*_\zeta)$ such that: 
\mr
\item "{$(a)$}"  $i_\zeta < \lambda^{++}$ is increasing continuous
\sn
\item "{$(b)$}"  $\zeta = \xi +1 \Rightarrow i_\zeta \in B_\varepsilon$
\sn
\item "{$(c)$}"  $M^*_\zeta \le_{\frak K} M_{i_\zeta}$
\sn
\item "{$(d)$}"  $M^*_\zeta \in K_{\frak s}$ is $\le_{\frak s}$-increasing
continuous
\sn
\item "{$(e)$}"  $N^*_\zeta \le_{\frak K} N_{i_\zeta}$
\sn
\item "{$(f)$}"  $N^*_\zeta \in K_{\frak s}$ is $\le_{\frak s}$-increasing
continuous
\sn
\item "{$(g)$}"  $M^*_\zeta = N^*_\zeta \cap M_{i_\zeta}$
\sn
\item "{$(h)$}"  $\zeta = \xi +1 \Rightarrow a_{i_\xi} \in M^*_{\zeta
+1}$
\sn
\item "{$(i)$}"  $M^*_0 = M'$ and $N^*_0 = N'_\varepsilon$ and $i_0 =
\text{ Min}(B_\varepsilon)$.
\ermn
There is no problem to do this and letting $a^*_\zeta = a_{i_\zeta}$ clearly
$\langle (M^*_\zeta,N^*_\zeta,a^*_\zeta):\zeta < \lambda^+ \rangle$
contradict ``${\frak s}$ is $\lambda$-good$^+$".
\hfill$\square_{\scite{705-stg.3}}$\margincite{705-stg.3}
\enddemo
\bn
The following claim sums up the ``localness" of the basic types for
${\frak s}(+)$, i.e., how to translate their properties to ones in
${\frak s}$.
\proclaim{\stag{705-stg.11} Claim}  [${\frak s}$ is a successful
good $\lambda$-frame].

Assume $\langle M_\alpha:\alpha < \lambda^+ \rangle$ is a
$\le_{\frak s}$-representation of $M \in K_{{\frak s}(+)}$. \nl
1) If $p_1,p_2 \in {\Cal S}^{\text{bs}}_{{\frak s}(+)}(M)$ 
\ub{then} $p_1 = p_2 \Leftrightarrow
\dsize \bigwedge_{\alpha < \lambda^+} p_1 \restriction M_\alpha = p_2
\restriction M_\alpha \Leftrightarrow (\exists^{\lambda^+} \alpha)(p_1
\restriction M_2 = p_2 \restriction M_\alpha) \Leftrightarrow (\exists
\beta < \lambda^+)[p_1 \restriction M_\beta = p_2 \restriction M_\beta \and
(\forall \alpha)(\beta \le \alpha < \lambda^+ \rightarrow p 
\restriction M_\alpha \in {\Cal S}^{\text{bs}}_{\frak s}(M_\alpha)$ does not
fork over $M_\beta)]$. \nl
2) If $S \subseteq \lambda^+ = \sup(S),
\alpha_* \le \,{\text{\rm min\/}}(S)$
and for $\alpha \in S,p_\alpha \in {\Cal S}^{\text{bs}}_{\frak s}(M_\alpha)$ does
not fork over $M_{\alpha_*}$ and $p_\alpha \restriction M_{\alpha_*} = p_*$,
\ub{then}
\mr
\item "{$(a)$}"  there is $p \in {\Cal S}^{\text{bs}}_{{\frak s}(+)}(M)$ satisfying
$\alpha \in S \Rightarrow p \restriction M_\alpha = p_\alpha$
\sn
\item "{$(b)$}"   there is no $p' \in {\Cal S}^{\text{bs}}_{{\frak s}(+)}(M)
\backslash \{p\}$ satisfying this, i.e., $p$ is unique
\ermn
3) If $p = { \text{\rm tp\/}}_{{\frak s}(+)}(a,M,N) \in {\Cal S}_{{\frak
s}(+)}(M)$ so $M \le_{{\frak s}(+)} N$ and $a \in N$, \ub{then} \nl
$p \in {\Cal S}^{\text{bs}}_{{\frak s}(+)}(M) \Leftrightarrow$ [for every $\alpha <
\lambda^+$ large enough, ${\text{\rm tp\/}}_{\frak s}(a,M_\alpha,N) \in {\Cal
S}^{\text{bs}}_{\frak s}(M_\alpha)] \Leftrightarrow [\text{for
stationarily may } \alpha < \lambda^+,{\text{\rm tp\/}}_{\frak
s}(a,M_\alpha,N) \in {\Cal S}^{\text{bs}}_{\frak s}(M_\alpha)]$. \nl
4) Assume $M_1 \le_{{\frak s}(+)} M_2 \le_{{\frak s}(+)} M_3$ and $\langle
M^\ell_\alpha:\alpha < \lambda^+ \rangle$ is a 
$\le_{\frak s}$-representation of $M_\ell$ and 
assume $a \in M_3$.  \ub{Then} ${\text{\rm tp\/}}_{{\frak s}(+)}
(a,M_2,M_3) \in {\Cal S}^{\text{bs}}_{{\frak s}(+)}(M_2)$ does not fork
over $M_1$ (for ${\frak s}^+$) \ub{iff} for some club $E$ of
$\lambda^+$, for every $\alpha < \beta$ from $E,{\text{\rm tp\/}}_{\frak s}
(a,M^2_\beta,M^3_\beta) \in {\Cal S}^{\text{bs}}_{\frak s}
(M^2_\beta)$ does not fork over $M^1_\alpha$ (for ${\frak s}$) 
\ub{iff} for some
stationary subset $S \subseteq \lambda^+$ for every $\delta$
from $S,{\text{\rm tp\/}}_{\frak s}(a,M^2_\delta,M^3_\delta) \in
{\Cal S}^{\text{bs}}_{\frak s}(M^2_\delta)$ does not fork over $M^1_i$. 
\endproclaim
\bigskip

\demo{Proof}  1) Among the four statements which we have to prove
equivalent, the first implies the second, trivially, the second
implies the third trivially, the third implies the second easily
and it implies the fourth by our assumption ``$p_1,p_2 
\in {\Cal S}^{\text{bs}}_{{\frak s}(+)}(M)$" and the definition of ${\Cal
S}^{\text{bs}}_{{\frak s}(+)}(M)$.  To finish we shall prove that the fourth
implies the first, so assume toward contradiction that this fails.

Let $M^0 = M,M^0_\alpha = M_\alpha$, let $M^\ell,a_\ell$ for $\ell = 1,2$
be such that $M^0 \le_{{\frak s}(+)} M^\ell,a_\ell \in M^\ell$ and
tp$(a_\ell,M^0,M^\ell) = p_\ell$ and let $\langle M^\ell_\alpha:\alpha
< \lambda^+ \rangle$ be a $\le_{\frak K}$-representation of $M_\ell$
for $\ell =1,2$.  Without loss of generality $\alpha < \beta \Rightarrow$
NF$_{\frak s}(M^0_\alpha,M^\ell_\alpha,M^0_\beta,M^\ell_\beta),M^\ell_{\alpha +1}$
is $(\lambda,*)$-brimmed over $M^0_{\alpha +1} \cup M^\ell_\alpha$ and
$a_\ell \in M^\ell_0$.  Clearly tp$_{\frak s}(a_1,M^0_0,M^1_0) =
\text{ tp}_{\frak s}(a_2,M^0_0,M^2_0)$ and  
we can build an isomorphism $f$ from $M^1_i$
onto $M_2$ over $M_0$ mapping $a_1$ to $a_2$ by choosing $f \restriction
M^1_\alpha:M^1_\alpha \underset{\text{onto}} {}\to \longrightarrow M^2_\alpha$
by induction on $\alpha < \lambda^+$.  \nl
2) Also easy. \nl
3), 4)  By the definition of 
${\frak s}^+$ and properties of NF.   \hfill$\square_{\scite{705-stg.11}}$\margincite{705-stg.11}
\enddemo
\bn
We may like in \scite{705-stg.11} to replace basic types by any types
(later this is needed and more is done):
\proclaim{\stag{705-stg.31} Claim}  [${\frak s}$ is a successful
good $\lambda$-frame].  ($\lambda^+$-locality)

Assume $\langle M_\alpha:\alpha < \lambda^+ \rangle$ is a $\le_{\frak
K}$-representation of $M \in K_{{\frak s}(+)}$. \nl
1) For any $p_1,p_2 \in {\Cal S}_{{\frak s}(+)}(M)$, \ub{then} $p_1 =
p_2 \Leftrightarrow (\forall \alpha)(p_1 \restriction M_\alpha = p_2
\restriction M_\alpha) \Leftrightarrow (\exists^{\lambda^+}
\alpha)(p_1 \restriction M_\alpha = p_2 \restriction M_\alpha)$.
\endproclaim
\bigskip

\demo{Proof}  See \chaptercite{600} on raising automorphisms: i.e., the
proof of \marginbf{!!}{\cprefix{600}.\scite{600-ne.4}}(2).  \nl
${{}}$  \hfill$\square_{\scite{705-stg.31}}$\margincite{705-stg.31}
\enddemo
\bn
Recall
\proclaim{\stag{705-stg.32} Claim}  If ${\text{\rm NF\/}}_{\frak s}(M_0,M_1,M_2,M_3)$ and
$(M_0,M_2,a) \in K^{3,\text{bs}}_{\frak s}$ \ub{then} {\rm tp}$(a,M_1,M_3) \in {\Cal
S}^{\text{bs}}_{\frak s}(M_1)$ does not fork over $M_0$.
\endproclaim
\bigskip

\demo{Proof}  See \marginbf{!!}{\cprefix{600}.\scite{600-nf.19}}.
\enddemo
\bigskip

\definition{\stag{705-e.A} Definition}  1) We define by induction on $n$:
\mr
\item "{$(a)$}"  ${\frak s}$ is $n$-successful
\sn
\item "{$(b)$}"  ${\frak s}^{+m} = {\frak s}(+m)$ for $m \le n$.
\ermn
\ub{For $n=0$}:  We say ${\frak s}$ is $0$-successful if it is
$\lambda({\frak s})$-good.

Let ${\frak s}^{+0} = {\frak s}$. 
\sn
\ub{For $n=1$}:  We say ${\frak s}$ is 1-successful \ub{if} it is
$\lambda_{\frak s}$-good and successful; let ${\frak s}^{+1} = {\frak s}^+$.
\sn
\ub{For $n=m+1$}:  We say ${\frak s}$ is $n$-successful if ${\frak s}^{+m}$
is $1$-successful.

We let ${\frak s}^{+(n+1)} = ({\frak s}^{+m})^+$. \nl
2) We say ${\frak s}$ is $(n + \frac 12)$-successful or say is weakly
$(n+1)$-successful \ub{if} it is 
$n$-successful and ${\frak s}^{+n}$ satisfies clause (a) of \scite{705-stg.0A}.
\nl
3) We say ${\frak s}$ is $\omega$-successful if it is $n$-successful for
every $n$. \nl
4) If ${\frak s}^{+n}$ is well defined let ${\frak B}_n = {\frak
B}^{\frak s}_n = {\frak B}^{{\frak s}(+n)}_n = 
{\frak B}({\frak s}^{+n})$ be a superlimit model in
${\frak K}[{\frak s}^{+n}]$; it is defined only up to isomorphism.
\enddefinition
\bigskip

\proclaim{\stag{705-e.B} Claim}  Assume ${\frak s}$ is an $n$-successful
good frame. \nl
1) {\rm NF}$_{{\frak s}[+n]} = { \text{\rm NF\/}}[{\frak s}^{+n}]$ is well defined \ub{if}
${\frak s}$ is $(n + \frac 12)$-successful. \nl
2) There is ${\frak B}^{\frak s}_n \in K_{\lambda^{+n}}$, that is a
${\frak K}[{\frak s}^{+n}]$-superlimit is well defined. \nl
3) ${\frak K}_{{\frak s}(+n)} = {\frak K}_{\frak s}[{\frak B}
^{\frak s}_n]$. \nl
4) If $k +m=n$ \ub{then} ${\frak s}^{+k}$ is $m$-successful good frame and
$(({\frak s}^{+k})^{+m}) = {\frak s}^{+m}$. \nl
5) ${\frak s}^k$ is $m$-successful iff ${\frak s}$ is
$(k+n)$-successful; and if this holds then $({\frak s}^{+k})^{+m} =
{\frak s}^{+(k+m)}$. 
\endproclaim
\bigskip

\remark{Remark}  If ${\frak s}$ is also a good$^+ \lambda$-frame
then ${\frak B}^{\frak s}_n$ is almost superlimit also in ${\frak
K}^{\frak s}_{\lambda^{+n}}$: the universality is not clear.
\endremark
\bigskip

\demo{Proof}  Easy, by induction on $n$ (and for (5) on $k+m$).
\hfill$\square_{\scite{705-e.B}}$\margincite{705-e.B}
\enddemo
\bigskip

\demo{\stag{705-stg.4} Conclusion}  In the main lemma
\marginbf{!!}{\cprefix{600}.\scite{600-fc.1}}, 
if we strengthen the assumption to ``${\frak s} =
({\frak K}_\lambda,{\Cal S}^{\text{bs}},\nonfork{}{}_{\lambda})$ is a
$\lambda$-good$^+$ frame", \ub{then} we can strengthen the conclusion to
\mr
\item "{$(\alpha)$}"  ${\frak s}_\ell = ({\frak K}^\ell_{{\frak s}_\ell},
{\Cal S}^{\text{bs}}_{{\frak s}_\ell},\nonfork{}{}_{{\frak s}_\ell})$
is $\lambda^{+ \ell}$-good
\sn
\item "{$(\beta)$}"  ${\frak s}_\ell = {\frak s}(\lambda^{+ \ell})$,
hence for some $M^*_\ell \in K_{\lambda^{+ \ell}}$ we have:
{\roster
\itemitem{ $(i)$ }  $M^*_\ell$ is superlimit (in 
${\frak K}^{\frak s}_{\lambda^{+ \ell}})$
\sn
\itemitem{ $(ii)$ }  ${\frak K}_{{\frak s}(\ell)}$ is 
${\frak K}^{[M^*_\ell]}$,
so $\le_{{\frak K}^\ell_{\lambda^{+ \ell}}} =
\le_{\frak K} \restriction K^\ell_{\lambda^{+ \ell}}$ (new
``$\le_{\frak K}$" is the new point)
\sn
\itemitem{ $(iii)$ }  ${\Cal S}^{\text{bs}}_{{\frak s}_\ell},
\nonfork{}{}_{{\frak s}_\ell}$ are defined as in \sectioncite[\S1]{600} but
restricted to ${\frak K}_{{\frak s}(\ell)}$ of course.
\endroster}
\endroster
\enddemo
\bigskip

\demo{Proof}  Should be clear (or combine \marginbf{!!}{\cprefix{600}.\scite{600-fc.1}} with
\scite{705-stg.3}).  \hfill$\square_{\scite{705-stg.4}}$\margincite{705-stg.4}
\enddemo
\bigskip

\proclaim{\stag{705-1.15} Claim}  Assume ${\frak s}$ is a weakly successful
good $\lambda$-frame.  Let $\delta < \lambda^+$, be a limit ordinal
and $\langle M_i:i \le \delta +1 \rangle$ be $\le_{\frak s}$-increasing
continuous. \nl
1) If $b \in M_{\delta +1}$ satisfies {\rm tp}$_{\frak s}(b,M_i,M_{\delta
+1}) \in {\Cal S}^{\text{bs}}_{\frak s}(M_i)$ for arbitrarily large $i <
\delta$, \ub{then} {\rm tp}$_{\frak s}(b,M_\delta,M_{\delta +1}) \in {\Cal
S}^{\text{bs}}_{\frak s}(M_\delta)$ hence does not fork over $M_i$ for every
$i < \delta$ large enough. 
\endproclaim
\bigskip

\demo{Proof}  1) Let $\langle N_i:i \le \delta \rangle$ be as in Claim
\scite{705-1.15K} below, so in particular $N_\delta$ is $(\lambda,*)$-brimmed over
$M_\delta$ hence $\le_{\frak s}$-universal over $N_\delta$, so \wilog
\, $M_{\delta +1} \le_{\frak s} N_\delta$.  So for some $i < \delta$
we have $b \in N_i$, so \wilog \, tp$_{\frak s}(b,M_i,M_{\delta +1})
\in {\Cal S}^{\text{bs}}_{\frak s}(M_i)$, now as NF$_{\frak s}
(M_i,M_\delta,N_i,N_\delta)$ by \scite{705-stg.32} we have tp$_{\frak
s}(b,M_\delta,N_\delta) = \text{ tp}_{\frak s}(b,M_\delta,M_{\delta
+1})$ is a nonforking extension of tp$_{\frak s}(b,M_i,N_i)$, and so
we are done.   \hfill$\square_{\scite{705-1.15}}$\margincite{705-1.15}
\enddemo
\bigskip

\proclaim{\stag{705-1.15K} Claim}  1) If $\langle M_i:i \le \delta \rangle$
is $\le_{\frak s}$-increasing continuous, then we can find $\langle
N_i:i \le \delta \rangle$ such that $M_i \le N_i,i < \delta
\Rightarrow { \text{\rm NF\/}}_{\frak s}(N_i,M_i,N_{i+1},M_{i+1})$ and
$N_{i+1}$ is universal over $N_i \cup M_{i+1}$. \nl
2) If $\langle M_i,N_i:i \le \delta \rangle$ are as in part (1), then
$N_\delta$ is $i \le j \le \delta \Rightarrow { \text{\rm NF\/}}_{\frak
s}(M_i,N_i,M_j,N_j)$ and $N_\delta$ is $(\lambda,*)$-brimmed over $M_\delta$.
\endproclaim
\bigskip

\demo{Proof}  By \marginbf{!!}{\cprefix{600}.\scite{600-nf.17}}.
\enddemo
\bigskip

\proclaim{\stag{705-1.15M} Claim}  If $(M_0,N_0,a) \le_{\text{bs}}
(M_j,N_1,a) \in K^{3,\text{bs}}_{\frak s}$ and $(M_0,N_0,a) \in
K^{3,\text{uq}}_{\frak s}$ then {\rm NF}$_{\frak s}(M_0,N_0,M_1,N_1)$.
\endproclaim
\bigskip

\demo{Proof}  By \marginbf{!!}{\cprefix{600}.\scite{600-nf.20}}.
\enddemo
\bigskip
 
\proclaim{\stag{705-1.17} Claim}   [${\frak s}$ is a weakly successful
good frame.]  Assume 
$M_0 \le_{\frak s} M_1 \le_{\frak s} M_2,a \in M_2$
and {\rm tp}$(a,M_1,M_2)$ does not fork over $M_0$ and $b \in M_1$,
{\rm tp}$(b,M_0,M_1) \in {\Cal S}^{\text{bs}}(M_0)$.  
\ub{Then} there are
$M^*_1,M^*_2$ such that $M_2 \le_{\frak s} M^*_2,M_0 \le_{\frak s}
M^*_1 \le_{\frak s} M^*_2,(M_0,M^*_1,a) \in 
K^{3,\text{uq}}_{\frak s}$ and
{\rm tp}$(b,M^*_1,M^*_2)$ does not fork over $M^*_0$.
\endproclaim
\bigskip

\demo{Proof}  By NF calculus (and the symmetry axiom).
\hfill$\square_{\scite{705-1.17}}$\margincite{705-1.17}
\enddemo
\bigskip

\proclaim{\stag{705-1.18} Claim}  Assume ${\frak s}$ is successful good$^+$
frame.  \ub{Then} ${\frak s}^+ = {\frak s}(\lambda^+_{\frak s})$ 
where on ${\frak s}(\lambda^+_{\frak s})$ see Definition \scite{705-0.X}(1).
\endproclaim
\bigskip

\demo{Proof}  Easy.
\enddemo
\newpage

\head {\S2  Unidimensionality and nonsplitting} \endhead  \resetall \sectno=2
 \spuriousreset
\bn
We may wonder how to define ``unidimensional" and whether: 
${\frak s}$ is categorical in $\lambda$ and is 
unidimensional and ${\frak s}$ is $n$-successful (see \scite{705-e.A}),
\ub{then} $K^{\frak s}_{\lambda^{+n}}$ is categorical.
By \scite{705-stg.10A} below the answer is yes.
\bn
We may consider a more restricted framework closed to categoricity.
Note that ``saturative" is closed to non multi-dimensional.
\bigskip

\demo{\stag{705-2.0} Hypothesis}  ${\frak s}$ is a $\lambda$-good frame.
\enddemo
\bigskip

\definition{\stag{705-stg.5} Definition}  1) We say ${\frak s}$ is
semi$^{\text{bs}}$
unidimensional \ub{when} for any model $M \in K_{\frak s}$, if 
$M <_{\frak s} N_k \in K_{\frak s}$ for $k=1,2$ 
\ub{then} some $p \in 
{\Cal S}^{\text{bs}}_{\frak s}(M)$ is realized in $N_1$
and in $N_2$. 
Let ``${\frak s}$ is semi$^{\text{na}}$ unidimensional" be 
defined similarly but we allow $p
\in {\Cal S}^{\text{na}}_{\frak s}(M)$ 
and ``${\frak s}$ is semi$^1$ unidimensional" be called
``${\frak s}$ is semi unidimensional".  Instead of na,bs we may write
$0,1$ respectively. \nl
2) We say ${\frak s}$ is almost unidimensional \ub{if} for any 
model $M \in K_{\frak s}$, there is an unavoidable $p \in {\Cal S}
^{\text{bs}}_{\frak s}(M)$ (see below). \nl
3) For $M \in K_{\frak s}$ we say $p \in {\Cal S}_{\frak s}(M)$ 
is $({\frak s})$-unavoidable, \ub{if} for every 
$N,M <_{\frak s} N \in K_{\frak s}$, some $a \in N$ realizes $p$.
\nl
4) We say ${\frak s}$ is explicitly unidimensional \ub{if} every $p \in 
{\Cal S}^{\text{bs}}_{\frak s}(M)$ where $M \in K_\lambda$, is unavoidable. \nl
5) We call ${\frak s}$ non-multi-dimensional \ub{if} for every $M_0 \in
K_\lambda$ whenever $M_0 <_{\frak s} M_1 <_{\frak s} M_2$, there is
$p \in {\Cal S}^{\text{bs}}_{\frak s}(M_1)$ which does not fork over $M_0$ and 
is realized in $M_2$. \nl
6) We say 
${\frak s}$ is weakly unidimensional \ub{if} for every $M <_{\frak s}
M_\ell$ for $\ell = 1,2$, there is $c \in M_2 \backslash M$ such that
tp$_{\frak s}(c,M,M_2)$ belongs to ${\Cal S}^{\text{bs}}_{\frak s}(M)$ and has 
more than one extension in ${\Cal S}^{\text{all}}_{\frak s}(M_1)$.
\enddefinition
\bn
On meaning in the first order case see \scite{705-stg.6}(5) below. \nl
We naturally first look at the natural implications.  Note that being
``semi$^1$/almost explicitly/weakly unidimensional" is influenced by the choice of the basic
types (compare \scite{705-stg.6}(1) with \scite{705-stg.6}(2)) as well as
non-multi-dimensional but not so semi$^0$-unidimensional.

\proclaim{\stag{705-stg.5A} Claim}  1) If 
${\frak s}$ is explicitly unidimensional, \ub{then}
${\frak s}$ is almost unidimensional. \nl
2) If ${\frak s}$ is almost unidimensional, \ub{then} is
semi$^{\,\ell}$-unidimensional for $\ell=0,1$;
if ${\frak s}$ is semi$^{\,1}$-unidimensional \ub{then} ${\frak s}$ is 
semi$^{\,0}$-unidimensional.
\nl
3) If ${\frak s}$ is semi$^{\,0}$-unidimensional, \ub{then} $K^{\frak s}$ is
categorical in $\lambda^+_{\frak s}$. \nl
4) If ${\frak s}$ is weakly unidimensional, \ub{then}
$K^{\frak s}$ is categorical in $\lambda^+_{\frak s}$. \nl
5) If ${\frak s}$ is semi-unidimensional, \ub{then} ${\frak s}$ is
weakly unidimensional.
\endproclaim
\bigskip

\remark{Remark}  1) Concerning 
non multi-dimensionality, does it follow from weak unidimensionality
under reasonable assumptions?  Yes, see \scite{705-c.3}.
\nl 
2) See more in \scite{705-stg.10}.
\endremark
\bigskip

\demo{Proof}  1) By Ax(D)(c) of $\lambda$-good frames (density of
basic types). \nl
2)  Check the definitions.
\nl
3)  Let $M_0,M_1 \in K^{\frak s}_{\lambda^+}$ and we shall prove that they
are isomorphic.  Let $\langle M^\ell_\alpha:\alpha < \lambda^+ \rangle$
be $<_{\frak s}$-representation of $M_\ell$ such that $\alpha <
\lambda^+ \Rightarrow M^\ell_\alpha \ne M^\ell_{\alpha +1}$ 
for $\ell = 0,1$.  Let
$\langle a^\ell_i:i < \lambda^+ \rangle$ list the elements of $M_\ell$.  We
choose by induction on $\varepsilon < \lambda^+$ a tuple $(N_\varepsilon,
\alpha^1_\varepsilon,f^1_\varepsilon,\alpha^2_\varepsilon,f^2_\varepsilon)$
such that:
\mr
\item "{$(a)$}"  $N_\varepsilon \in 
K_{\frak s}$ is $\le_{\frak s}$-increasing continuous
\sn
\item "{$(b)_\ell$}"  $\alpha^\ell_\varepsilon < \lambda^+$ is increasing
continuous
\sn
\item "{$(c)_\ell$}"  $f^\ell_\varepsilon$ is a $\le_{\frak s}$-embedding of
$M^\ell_{\alpha^\ell_\varepsilon}$ into $N_\varepsilon$
\sn
\item "{$(d)_\ell$}"  $f^\ell_\varepsilon$ is increasing continuous
\sn
\item "{$(e)_\ell$}"  if $\varepsilon = 4 \zeta + \ell,\ell \in \{0,1\}$ and 
$j_\varepsilon = \text{ Min}
\{i:f^\ell_{4 \zeta}( \text{tp}
(a^\ell_i,M^\ell_{\alpha_\varepsilon},M_\ell))$ is
realized by some $d \in N_{4 \zeta} \backslash \text{
Rang}(f^\ell_\varepsilon)\}$
is well defined \ub{then} $a^\ell_{j_\varepsilon} \in \text{
Dom}(f^\ell_{\varepsilon +1})$ and $f^\ell_{\varepsilon +1}
(a^\ell_{j_\varepsilon}) \in N_{4 \zeta}$ and 
$f^{1- \ell}_{\varepsilon +1} = f^{1- \ell}_\varepsilon\}$
\sn
\item "{$(f)_\ell$}"  if $\varepsilon 
= 4 \zeta +2 + \ell,\ell \in \{0,1\}$ \ub{then}
$\alpha^\ell_{\varepsilon +1} > \alpha^\ell_\varepsilon$.
\ermn
If we succeed, then for some club $E$ of $\lambda^+$ we have
$a^\ell_\alpha \in M^\ell_\delta \Leftrightarrow \alpha < \delta$ and
$\delta \in E \Rightarrow  N_\delta \cap
\dbcu_{\varepsilon < \lambda^+} \text{ Rang}(f^\ell_\varepsilon)
 = \text{ Rang}(f^\ell_\delta)$.  If $\delta \in E$ and $\ell \in
\{0,1\}$ but Rang$(f^\ell_\delta) \ne N_\delta$ then by the assumption
(``${\frak s}$ is semi$^0$-unidimensional") see 
Definition \scite{705-stg.5}(1), for
some $c \in M^\ell_{\alpha^\ell_\delta +1} \backslash 
M^\ell_{\alpha^\ell_\delta}$ and $d \in N_\delta
\backslash \text{ Rang}(f^\ell_\delta)$ we have tp$(d,\text{Rang}(f^\ell
_\delta),N_\delta) = f^\ell_\delta(\text{tp}(c,M^\ell_{\alpha_\delta},
M_\ell))$ so by clause $(e)_\ell$ (as $f^\ell_{\delta + \ell} = 
f^\ell_\delta$) we have Rang$(f^\ell_{\delta + 2}) \cap N_\delta \backslash
\text{ Rang}(f^\ell_\delta) \ne \emptyset$ contradiction.  So $\delta  \in E
\wedge \ell \in \{0,1\} \Rightarrow \text{ Rang}(f^\ell_\delta) = N_\delta$
hence $f_\ell =: \dbcu_{\delta \in E} f^\ell_\delta$ is an isomorphism from
$M_\ell$ onto $N =: \dbcu_{\delta \in E} N_\delta$, so $M_1 \cong N \cong
M_2$ and we are done.

So we have just to carry the induction, which is straight as $K_{\frak s}$ is
a $\lambda_{\frak s}$-a.e.c. with amalgamation and the hypothesis
\footnote{actually the ``semi$^0$-unidimensional" can be weakened - in
Definition \scite{705-stg.5}(1) we may ask $N_1 \in K_{{\frak s},\lambda^+}$}. \nl
4) The proof is similar to that of part (3) but we replace clause
$(e)_\ell$ by
\mr
\item "{$(e)^*_\ell$}"  if $\varepsilon = 4 \zeta + \ell$ and for some
$c \in N_{4 \zeta} \backslash \text{ Rang}(f^\ell_\varepsilon)$ we have
tp$(c,f^\ell_\varepsilon(M^\ell_{\alpha^\ell_\varepsilon}),N_\varepsilon)
\in {\Cal S}^{\text{bs}}(f^\ell_\varepsilon
(M^\ell_{\alpha^\ell_\varepsilon}))$ and 
$(f^\ell_\varepsilon)^{-1}(\text{tp}(c,f^\ell_\varepsilon
(M^\ell_{\alpha^\ell_\varepsilon}),N_\varepsilon))$ has at least two 
extensions in ${\Cal S}^{\text{all}}(M^\ell_\beta)$ for some $\beta \in
(\alpha^\ell_\varepsilon,\lambda^+)$ \ub{then} 
tp$(c,f^\ell_{\varepsilon +1}(M^\ell_{\alpha^\ell_{\varepsilon +1}}),
N_{\varepsilon +1})$ is not the nonforking extension of tp$(c,
f^\ell_\varepsilon(M^\ell_{\alpha^\ell_\varepsilon}),N_\varepsilon)$ in
${\Cal S}^{\text{bs}}(f^\ell_{\varepsilon +1}(M^\ell_{\alpha^\ell
_{\varepsilon +1}}))$.
\ermn
Again there is no problem to carry the definition.  So it is enough to
prove that $\cup\{\text{Rang}(f^\ell_\varepsilon):\varepsilon < \lambda^+\}
= N_{\lambda^+}$, and for this it suffices to prove that $S_\ell = \{\delta <
\lambda^+:N_\delta \ne \text{ Rang}(f^\ell_\delta)\}$ is not stationary.
For every $\delta \in S_\ell$ by the assumption ``${\frak s}$ is weakly
unidimensional" we know that the assumption of $(e)^*_\ell$ holds
hence there is $c = c^\ell_\delta$ as there.
By Fodor lemma for some $c_\ell$ the set $S'_\ell = \{\delta \in
S_\ell:c^\ell_\delta = c_\ell\}$ is stationary.  Choose $\delta^* \in
S'_\ell,\delta^* = \sup(\delta^* \cap S'_\ell)$ and use ``${\frak s}$
is good".
\nl
5) Easy.   \hfill$\square_{\scite{705-stg.5A}}$\margincite{705-stg.5A}
\enddemo
\bn
A conclusion is (see Definition \scite{705-0.X}(r)): \nl
\margintag{705-stg.5B}\ub{\stag{705-stg.5B} Conclusion}:  [${\frak s}$ is successful good$^+ \,
\lambda$-frame.] \nl
1) If ${\frak s}$ is semi$^0$-unidimensional \ub{then} ${\frak s}\langle 
\lambda^+ \rangle = {\frak s}^+$ so ${\frak K}^{\frak s}_{\lambda^+} =
{\frak K}_{{\frak s}(+)}$.
\nl
2) Similarly if $K^{\frak s}$ is categorical in $\lambda^+_{\frak s}$. \nl
3) If ${\frak s}$ is non-multi-dimensional, \ub{then} so is ${\frak s}^+$.
\bigskip

\demo{Proof}  1) Clearly $K_{{\frak s}(+)} \subseteq K^{\frak
s}_{\lambda^+}$ and by \scite{705-stg.3B} we know $\le_{K^{\frak s}}
\restriction K_{{\frak s}(+)} = \le_{{\frak s}(+)}$.  But by 
\scite{705-stg.5A}(3), $K_{{\frak s}(+)}$ is categorical in
$\lambda_{{\frak s}(+)} = \lambda^+$.  Hence $K_{{\frak s}(+)} =
K^{\frak s}_{\lambda^+}$ and even ${\frak
K}_{{\frak s}(+)} = {\frak K}^{\frak s}_{\lambda^+}$ and check
similarly above $\nonfork{}{}_{}$ and ${\Cal S}^{\text{bs}}$. \nl
2) So again $K_{{\frak s}(+)} = K^{\frak s}_{\lambda^+}$ and just
check. \nl
3) By \scite{705-stg.11}.    \hfill$\square_{\scite{705-stg.5B}}$\margincite{705-stg.5B} 
\enddemo
\bigskip

\remark{Remark}  But we may need ``non-multi-dimensional" (defined below),
does it follow?

We also note that the cases we have dealt with categoricity
hypothesis, give not just good frames but even unidimensional ones.
\endremark
\bigskip

\proclaim{\stag{705-stg.6} Claim}  1) In \marginbf{!!}{\cprefix{600}.\scite{600-Ex.4}} 
(= \scite{705-stg.2}(1) Case 1 above) we can add:  
the ${\frak s}$ obtained there is explicitly unidimensional. \nl
2) In \marginbf{!!}{\cprefix{600}.\scite{600-Ex.1}}, if ${\frak K}$ is categorical in $\aleph_1$
(see above \scite{705-stg.2}(1), Case 2)
\ub{then} we get almost unidimensionality for the ${\frak s}$ obtained there. \nl
3) In \marginbf{!!}{\cprefix{600}.\scite{600-Ex.1A}}, if $\psi$ is categorical in $\aleph_1$
(see above \scite{705-stg.2}(1), Case 3) \ub{then} we get almost
unidimensionality for the ${\frak s}$ obtained there and is
unidimensional. \nl
4) In \marginbf{!!}{\cprefix{600}.\scite{600-rg.7}}(3), if ${\frak s} = ({\frak K}_\lambda,
{\Cal S}^{\text{bs}},\nonfork{}{}_{\lambda})$ is almost unidimensional (see above
\scite{705-stg.2}(3)), \ub{then} we get that also the $\lambda^+$-good frame
system ${\frak s}^+$ obtained there is almost unidimensional.  \nl
5) If $T$ is complete superstable first order and ${\frak s} = {\frak
s}^\kappa_{T,\lambda}$ (see \scite{705-stg.2}(3)) and $\lambda \ge |T| +
\kappa^+ (\kappa > 0 \Rightarrow T$ stable in $\lambda$) then:
\mr
\item "{$(i)$}"  ${\frak s}$ is saturative iff $T$ is
non-multidimensional (see \cite{Sh:c}; this is (ii) of
\scite{705-stg.2}(3))
\sn
\item "{$(ii)$}"  ${\frak s}$ is categorical in $\lambda^+$ \ub{iff}
$T$ is unidimensional \ub{iff} ${\frak s}$ is almost unidimensional.
\endroster 
\endproclaim
\bigskip

\demo{Proof}  Easy.  (in (2) we have to use a minimal type see
\marginbf{!!}{\cprefix{88}.\scite{88-xx}}, we get more on (5) follows from the claims below).
\enddemo
\bn
Now we consider those properties and how they are related in ${\frak
s}$ and ${\frak s}^+$.
\proclaim{\stag{705-stg.6A} Claim}  [$K_{\frak s}$ categorical in 
$\lambda_{\frak s}$]. 
If $K^{\frak s}_{\lambda^+}$ is categorical in $\lambda^+$ 
\ub{then} \,${\frak s}$ is weakly unidimensional. 
\endproclaim
\bigskip

\demo{Proof}  Assume toward contradiction that 
${\frak s}$ is not weakly unidimensional,
hence we can find $M_0 <_{\frak s} M_\ell$ for $\ell = 1,2$ such that:
if $c \in M_2 \backslash M_0$, tp$(c,M_0,M_2) \in {\Cal S}^{\text{bs}}(M_0)$
then it has a unique extension in ${\Cal S}(M_1)$. By Axiom (D)(d) of
$\lambda$-good frames (existence)
we can choose $c \in M_2 \backslash M_0$ such that $p = \text{ tp}(c,M_0,M_2)
\in {\Cal S}^{\text{bs}}(M_0)$.  Now we choose by induction on $\alpha <
\lambda^+$ a model $N_\alpha \in K_{\frak s},\le_{\frak s}$-increasing
continuous, $N_\alpha \ne N_{\alpha +1},N_0 = M_0$ and $p$ has a unique
extension in ${\Cal S}(N_\alpha)$, call it $p_\beta$ and by Axiom (E)(g)
(extension) we know that 
$p_\alpha \in {\Cal S}^{\text{bs}}(N_\alpha)$ does not
fork over $N_0$.  For $\alpha =0$ this is trivial, for $\alpha = \beta+1$ by
\scite{705-stg.9} below (noting that every $M \in K_{\frak s}$ is
isomorphic to $M_0$ and is $(\lambda,*)$-brimmed as $K_{\frak s}$ is categorical in $\lambda$)
there is an isomorphism $f_\beta$ from $N_0 = M_0$ onto 
$N_\beta$ such that $f_\beta(p_0) = p_\beta$, so we can find $N_\alpha =
N_{\beta+1}$ and isomorphism $g_\beta$ from $M_1$ onto $N_\alpha$ extending
$f_\beta$.  Hence $f(p) = p_\beta$ and so $p_\beta$ has a unique extension
in ${\Cal S}(g_\beta(M_1)) = {\Cal S}_{\frak s}(N_\alpha)$ as required.  For $\beta$
limit use Axiom (E)(h), continuity.

Now $N =: \dbcu_{\alpha < \lambda^+} N_\alpha \in K^{\frak s}_{\lambda^+}$
(recall $N_\alpha \ne N_{\alpha +1}$), and $p_\beta$ is not realized in
$N_\beta$ for $\beta < \lambda^+$ hence $p = p_0$ is not realized in $N$,
so $N \in K^{\frak s}_{\lambda^+}$ is not saturated contradicting
categoricity in $\lambda^+$.  \hfill$\square_{\scite{705-stg.6A}}$\margincite{705-stg.6A}
\enddemo
\bigskip

\proclaim{\stag{705-stg.10} Claim}  Assume 
${\frak s}$ is a successful $\lambda$-good frame. \nl
1) If ${\frak s}$ is semi$^x$-unidimensional where $x \in \{{\text{\rm
na,bs\/}}\}$ \ub{then} ${\frak s}^+$ is semi$^\ell$ dimensional. \nl
2) If ${\frak s}$ is weakly unidimensional, \ub{then} ${\frak s}^+$ is
weakly unidimensional. \nl
3) If ${\frak s}$ is almost unidimensional, \ub{then} ${\frak s}^+$ is almost
unidimensional. \nl
4) If ${\frak s}$ is explicitly unidimensional, \ub{then} 
${\frak s}^+$ is explicitly unidimensional.
\endproclaim
\bigskip

\demo{Proof}  1) First let $x = \text{ na}$.  
Assume toward contradiction, that ${\frak s}^+$ is not 
semi$^\ell$-unidimensional, so we can find $M_0 <_{{\frak s}(+)} M_\ell$ 
for $\ell = 1,2$ such that $c_1 \in M_1 \backslash
M_0 \and c_2 \in M_2 \backslash M_0 \Rightarrow \text{ tp}(c_1,M_0,M_1) \ne
\text{ tp}(c_2,M_0,M_2)$.  Let $\langle M^\ell_\alpha:\alpha < \lambda^+
\rangle$ be a $\le_{\frak s}$-representation of $M_\ell$ for $\ell <
3$, hence for some club $E$ of $\lambda^+$ the following holds:
\mr
\item "{$(*)$}"   for $\ell \in \{1,2\}$ and
$\alpha < \beta$ in $E$, we have NF$_{\frak s}(M^0_\alpha,M^\ell_\alpha,
M^0_\beta,M^\ell_\beta)$ (hence $M^\ell_\alpha \cap M_0 = M^0_\alpha$) and
$M^0_\alpha \ne M^\ell_\alpha$ and $c_1 \in M^1_\alpha \backslash M^0_\alpha
\and c_2 \in M^2_\alpha \backslash M^0_\alpha \and (\exists \gamma <
\lambda^+)(\text{tp}(c_1,M^0_\gamma,M^1_\gamma) \ne \text{tp}(c_2,M^0_\gamma,
M^2_\beta) \Rightarrow \text{ tp}(c_1,M^0_\alpha,M^1_\alpha) \ne
\text{ tp}(c_2,M^0_\alpha,M^2_\alpha)$.
\ermn
Let $\delta = \text{ Min}(E)$ and apply the assumption so there are
$c_1 \in M^1_\delta \backslash M^0_\delta,c_2 \in M^2_\delta \backslash
M^0_\delta$ satisfying tp$(c_1,M^0_\delta,M^1_\delta) = \text{ tp}(c_2,
M^0_\delta,M^2_\delta)$.  By the choice of $E$ we have $\beta < \lambda^+
\Rightarrow \text{ tp}(c_1,M^0_\beta,M_1) = \text{ tp}(c_2,M^0_\beta,M_2)$
and use \scite{705-stg.31}(1).  \nl
If $x = \text{ bs}$ the proof if similar using basic types (and
\scite{705-stg.11}(3)) or use \scite{705-stg.5A}(2), second clause. \nl
2)  So toward contradition assume $M_0 <_{{\frak s}(+)} M_\ell$ for
$\ell = 1,2$ are such that: if $c \in M_1 \backslash M_0$ and $p =
\text{ tp}_{{\frak s}(+)}(c,M_0,M_1) \in 
{\Cal S}^{\text{bs}}_{{\frak s}(+)}(M_0)$,
then $p$ has a unique extension in ${\Cal S}^{\text{bs}}_{{\frak s}
(+)}(M_2)$.
Let $\langle M^\ell_\alpha:\alpha < \lambda^+ \rangle$ be a
$\le_{\frak s}$-representation of $M_\ell$ and $E$ a thin 
enough club of $\lambda^+$.

For each $\delta \in E$ 
we have $M^0_\delta <_{\frak s} M^\ell_\delta$ for
$\ell =1,2$ but ${\frak s}$ is weakly unidimensional hence for some
$c_\delta \in M^1_\delta \backslash M^0_\delta$ the type $p_\delta =
\text{ tp}_{\frak s}(c,M^0_\delta,M^1_\delta)$ 
belongs to ${\Cal S}^{\text{bs}}
_{\frak s}(M^0_\delta)$ and has more than one 
extension in ${\Cal S}(M^2_\delta)$.
Clearly there is $\alpha_\delta < \delta$ such that $p_\delta$ does not
fork for ${\frak s}$ over $M^0_{\alpha_\delta}$ and trivially $p_\delta
\restriction M^0_{\alpha_\delta} \in {\Cal S}^{\text{bs}}_{\frak s}
(M^0_{\alpha_\delta})$
which has cardinality $\le 
\lambda$.  By Fodor lemma for some stationary
$S \subseteq \lambda^+$ we have $\delta \in S \Rightarrow c_\delta = c_*
\and \alpha_\delta = \alpha_* \and p_\delta \restriction M^0_{\alpha_*} =
p_*$.  By \scite{705-stg.31} the type
$q =: \text{ tp}_{{\frak s}(+)}(c_*,M_0,M_1)$ belongs to
${\Cal S}^{\text{bs}}_{{\frak s}(+)}(M_0)$  hence by 
the choice of $M_0,M_1,M_2$ (``toward contradiction") $q$ has at least two distinct extensions
$q_1,q_2 \in {\Cal S}_{{\frak s}(+)}(M_2)$.  Now by the choice of $S$
and as for each $\delta \in S$ there is $q'_\delta \in {\Cal
S}^{\text{bs}}_{\frak s}(M^2_\delta)$ which is a non-forking extension
of $p_*$ clearly
\mr
\item "{$(**)$}"  if $\delta \in S \and \ell \in \{1,2\}$ then $q_\ell
\restriction M^2_\delta$ belongs to ${\Cal S}^{\text{bs}}_{\frak s}(M^2_\delta)$
and does not fork over $M^0_\delta$ and extends $p_* = p_\delta \restriction
M^0_{\alpha_*}$ hence does not fork over $M^0_{\alpha_*}$.
\ermn
By \scite{705-stg.11} we get a contradiction. \nl
3) As for $M_0,\langle M^0_\alpha:\alpha < \lambda^+ \rangle$ as above
there is $p \in {\Cal S}^{\text{bs}}_{{\frak s}(+)}(M_0$ such that
$\{\delta < \lambda^+:p \restriction M^0_\alpha \in {\Cal
S}^{\text{bs}}_{\frak s}(M^0_\alpha)$ is unavoidable$\}$ is
stationary. \nl
4) Straightforward.  \hfill$\square_{\scite{705-stg.10}}$\margincite{705-stg.10}
\enddemo
\bn
Together we can ``close the circle" to continuing ``up" we shall get
more (see more in \scite{705-2e.3}). \nl
\margintag{705-stg.10A}\ub{\stag{705-stg.10A} Conclusion}:  Assume ${\frak s}$ is a
$\lambda$-good frame categorical in $\lambda$. \nl
1) Then ${\frak s}$ is weakly-unidimensional \ub{iff} 
$K^{\frak s}_{\lambda^+}$ is categorical in $\lambda^+$. \nl
2) If ${\frak s}$ is successful we can add: iff ${\frak s}^+$ is
weakly unidimensional and ${\frak K}^{\frak s}_{\lambda^+} = {\frak
K}_{s(+)}$.
\bigskip

\demo{Proof}  1) The second condition implies the first by
\scite{705-stg.6A}, the first condition implies the second by
\scite{705-stg.5A}(4). \nl
2) The first implies the third as by \scite{705-stg.10}(2) we have ${\frak
s}^+$ is weakly unidimensional and ${\frak K}
^{\frak s}_{\lambda^+} = {\frak K}_{s(+)}$ by \scite{705-0.X}(1) 
and the third condition implies the second 
as ${\frak K}_{{\frak s}(+)}$ is categorical in $\lambda^+$ by the
definition of ${\frak K}_{{\frak s}(+)}$.  \hfill$\square_{\scite{705-stg.10A}}$\margincite{705-stg.10A}
\enddemo
\bn
Earlier (say in \cite{Sh:576}) minimal type were central, so let us
mention them:
\definition{\stag{705-stg.7} Definition}  1) We say 
${\frak s}$ is (a $\lambda$-good
frame) of minimals when the following holds: if
$p \in {\Cal S}^{\text{bs}}(M_0),M_0 \le_{\frak s} M_1
\le_{\frak s} N_1,M_0 \le_{\frak s} N_0 \le_{\frak s} N_1,a \in M_1,
p = \text{ tp}(a,M_0,M_1)$ and $a \notin N_0$ \ub{then} tp$(a,N_0,N_1)$ is a
nonforking extension of $p$. \nl
So the triple $(M,N,a)$ is called ${\frak s}$-minimal. \nl
2) For an $\lambda$-a.e.c. ${\frak K},p \in {\Cal S}(M)$ is
minimal \ub{if} for every $N,M \le_{\frak K} N \in K_\lambda,p$ has one and
only one extension in ${\Cal S}(N)$.
\enddefinition
\bigskip

\proclaim{\stag{705-stg.8} Claim}  1) In \marginbf{!!}{\cprefix{600}.\scite{600-Ex.4}} 
(= above in \scite{705-stg.2}(1))
we can add:  ${\frak s}$ is a $\lambda$-good frame of minimals. \nl
2) Similarly in \scite{705-stg.6}(2),(3). \nl
3) In \scite{705-stg.3B}, if ${\frak s}$ is a frame of minimals \ub{then} so is
${\frak s}^+$. \nl
4) If $(M_0,M_1,a)$ is ${\frak s}$-minimal 
(i.e., as in Definition \scite{705-stg.7}(1)) \ub{then}:
\mr
\widestnumber\item{$(iii)$}
\item "{$(i)$}"  $p = { \text{\rm tp\/}}_{{\frak K}_{\frak s}}
(a,M_0,M_1)$ is minimal, 
\sn
\item "{$(ii)$}"  if $M_0 \le_{\frak s} M_1$ and $q \in {\Cal
S}_{\frak s}(M_1)$
extends $p$ is not algebraic \ub{then} $q$ does not fork over $M_0$;
hence, in particular, $\in {\Cal S}^{\text{bs}}_{\frak s}(M_0)$ 
\sn
\item "{$(iii)$}"  if $M_0 \le_{\frak s} M_1$ and 
$q = \text{ {\rm tp\/}}(b,M_0,M_1)$ satisfies clauses (i) and (ii) \ub{then} 
$q$ is minimal.
\endroster
\endproclaim
\bn
\centerline{$* \qquad * \qquad *$}
\bn
We could have mentioned in \chaptercite{600}:
\proclaim{\stag{705-stg.9} Claim} [${\frak s}$ is a $\lambda$-good frame]. \nl
Assume $M_1 \le_{\frak s} M_2$ are superlimit in $K_{\frak s}$
and $p_i \in {\Cal S}^{\text{bs}}(M_2)$ does not fork 
over $M_1$ for $i < \alpha < \lambda_{\frak s}$.  \ub{Then} there 
is an isomorphism $f$ from $M_1$
onto $M_2$ such that $i < \alpha \Rightarrow f(p_i \restriction M_1) = p_i$.
\endproclaim
\bigskip

\demo{Proof}  First assume that $M_2$ is $(\lambda,*)$-brimmed over $M_1$.
Clearly we can find a regular cardinal $\theta$ such that $\alpha < \theta
\le \lambda$.  Now we can find a sequence $\langle N_\beta:\beta < \theta
\rangle$ which is $\le_{\frak s}$-increasing continuous, $N_{\beta +1}$ being
$(\lambda,*)$-brimmed over $N_\beta$ (of course, we are using
\sectioncite[\S4]{600}).  Clearly $\dbcu_{\beta < \theta} N_\beta \in
K_{\frak s}$ is $(\lambda,*)$-brimmed over $N_0$, so without loss of generality is equal to $M_1$.

So for each $i < \alpha$ for some $\beta(i) < \theta$ the type $p_i
\restriction M_1$ which $\in {\Cal S}^{\text{bs}}(M_1)$ 
does not fork over $N_{\beta(i)}$, so $\beta = \sup\{\beta
(i):i < \alpha\} < \theta$, hence by transitivity and monotonicity of
nonforking $i < \alpha \Rightarrow p_i$ does not fork
over $N_\beta$.  Clearly also $M_2$ is $(\lambda,*)$-brimmed over $N_\beta$
and by the choice of $\langle N_\gamma:\gamma < \theta \rangle$ also 
$M_1$ is $(\lambda,*)$-brimmed over $N_\beta$
hence there is an isomorphic $f$ from $M_2$ onto $M_1$ over $N_\beta$.  Now
for $i < \alpha$ the types $p_i \restriction M_1$ and $f(p_i)$ are members of
${\Cal S}^{\text{bs}}(M_1)$ which does not fork over $N_\beta$ and has the
same restriction to $N_\beta$ hence are equal.  So $f^{-1}$ is as required.

Second without the assumption ``$M_2$ is $(\lambda,*)$-brimmed over $M_1$"
we can find $M_3 \in K_{\frak s}$ which is $(\lambda,*)$-brimmed over
$M_2$ hence also over $M_1$ and let $q_i \in {\Cal S}^{\text{bs}}(M_3)$ be a
nonforking extension of $p_i$.

Applying what we have already proved to the pair 
$(M_1,M_3)$ there is an isomorphism
$f_1$ from $M_1$ onto $M_3$ mapping $p_i \restriction M_1 = q_i
\restriction M_1$ to $q_i$ for
$i < \alpha$.  Applying what we have already proved to the pair
$(M_2,M_3)$, there
is an isomorphism $f_2$ from $M_2$ onto $M_3$ mapping $p_i$ to $q_i$ for
$i < \alpha$. Now $f^{-1}_2 \circ f_1$ is as required. 
\hfill$\square_{\scite{705-stg.9}}$\margincite{705-stg.9}
\enddemo
\bn
Recalling Definition \scite{705-0.X}(3) note:
\proclaim{\stag{705-stg.33} Claim}  1) Assume $\bar M = \langle M_i:i \le
\delta +1 \rangle$ is $<_{\frak s}$-increasing continuous.  If
$\Gamma \subseteq {\Cal S}^{\text{bs}}(M_{\delta +1}),|\Gamma| < 
\,{\text{\rm cf\/}}(\delta)$ 
and $M_\delta \le_{\frak s} N$ and $p \in \Gamma \Rightarrow
p$ does not fork over $M_\delta$.  \ub{Then} for every large enough $i
< \delta$ there is an isomorphism $f$ from $N$ onto $M_{\delta +1}$ over
$M_i$ such that $p \in \Gamma \Rightarrow f(p) = p \restriction
M_\delta$ provided that $i < \delta \Rightarrow M_\delta,N$ are
$(\lambda,*)$-brimmed over $M_i$. \nl
2) Instead $|\Gamma| < {\text{\rm cf\/}}(\delta)$ it is enough to
demand: $|\Gamma| < \lambda_{\frak s}$ and $p \in \Gamma \Rightarrow
p$ does not fork over $M_i$.
\endproclaim
\bigskip

\demo{Proof}  1) For $p \in \Gamma$, choose $i(p) < \delta$ such that $p
\in \Gamma \Rightarrow p \restriction M_\delta$ does not fork over
$M_{i(p)}$ and let $i(*) = \sup\{i(p):p \in \Gamma\}$ it is $< \delta$
or $|\Gamma| < \text{ cf}(\delta)$. \nl
2) Similar.  \hfill$\square_{\scite{705-stg.33}}$\margincite{705-stg.33}
\enddemo
\bn
\centerline{$* \qquad * \qquad *$}
\bigskip

\definition{\stag{705-gr.1} Definition}  Let ${\frak K}$ be a
$\lambda$-a.e.c. (so ${\frak K} = {\frak K}_\lambda$)(normally with
amalgamation (in $\lambda$)). \nl
1) We say that $p \in {\Cal S}_{\frak K}(M_1) \, \alpha$-splits or $(\alpha,{\frak K})$-split over 
$A \subseteq M_1$ \ub{if} there are 
$\bar a_1,\bar a_2 \in {}^\alpha(M_1)$ such that:
\mr
\item "{$(\alpha)$}"  $\bar a_1,\bar a_2$ realize the same type over $A$
inside $M_1$ that is,
{\roster
\itemitem{ $(*)$ }  for some $M_2,f$ we have:

$$
M_1 \le_{\frak K} M_2
$$

$$
f \text{ is an automorphism of } M_2 \text{ over } A \text{ mapping }
\bar a_1 \text{ to } \bar a_2
$$
\endroster}
\item "{$(\beta)$}"  if $M_1 \le_{\frak K} M_2$ and $c \in M_2$ realizes $p$
inside $M_2$ \ub{then} $\bar a_1,\bar a_2$ do not realize the same
type over $A \cup \{c\}$ inside $M_2$, that is 
 for no $M_3,f$ do we have $M_2 \le_{\frak K} M_3$
and $f$ is an automorphism of $M_3$ over $M$ mapping $\bar a_1
\char 94 \langle c \rangle$ to $\bar a_2 \char 94 \langle c \rangle$.
\ermn
3) We may write $\bar a$ instead of $A = \text{ Rang}(\bar a)$ and
$M_0$ instead of $A = |M_0|$.  If we omit $\alpha$ (and write split or
${\frak K}$-split) we mean ``for some $\alpha$". \nl
4) We say ${\frak K}$ has $\chi$-nonsplitting \ub{if} for 
every $M \in K_\lambda$ and $p \in {\Cal S}_{{\frak
K}_\lambda}(M)$  
there is $A \subseteq M,|A| \le \chi$ such
that $p$ does not split over $A$ (in ${\frak K}$). \nl
5) We say ${\frak s}$ has $\chi$-nonsplitting \ub{if} ${\frak K}_{\frak s}$
has basically $\chi$-nonsplitting which means that 
this holds for $p \in {\Cal S}^{\text{bs}}_{\frak s}(M)$. \nl
6) In part (1), (2), (3) though not (4) writing ${\frak s}$ instead of
${\frak K}$ means ${\frak K}_{\frak s}$.
\enddefinition
\bigskip

\proclaim{\stag{705-gr.5} Claim}  1) If {\rm NF}$_{\frak
s}(M_0,M_1,M_2,M_3)$ and $\bar c \subseteq M_2$ \ub{then} 
{\rm tp}$_{\frak s}(\bar c,M_1,M_3)$ does not split over $M_0$. \nl
2) Similarly for $\bar c \in {}^\alpha(M_2)$.
\endproclaim
\bigskip

\demo{Proof}  Straightforward (by uniqueness of NF).
\enddemo
\bn
We could have noted earlier:
\proclaim{\stag{705-gr.6} Claim}  Assume
\mr
\item "{$(a)$}"  $\delta < \lambda^+_{\frak s}$ is a limit ordinal
\sn
\item "{$(b)$}"  $\langle M_\alpha:\alpha \le \delta \rangle$ is
$\le_{\frak s}$-increasing continuous
\sn
\item "{$(c)$}"  $M_{\alpha +1}$ is $(\lambda,*)$-brimmed over
$M_\alpha$
\sn
\item "{$(d)$}"  $p \in {\Cal S}_{\frak s}(M_\delta)$.
\ermn
\ub{Then} for some $i < \delta$ the type $p$ 
does not $\lambda$-split over $M_i$ for ${\frak K}_{\frak s}$.
\endproclaim
\bigskip

\demo{Proof}  We can find a $<_{\frak s}$-increasing continuous
sequence $\langle N_\alpha:\alpha \le \delta \rangle$ such that
$M_\alpha \le_{\frak s} N_\alpha$ and $N_{\alpha +1}$ is
$(\lambda,*)$-brimmed over $M_{\alpha +1} \cup N_\alpha$ and
NF$_{\frak s}(M_\alpha,N_\alpha,M_{\alpha +1},N_{\alpha +1})$.  We
know (\sectioncite[\S6]{600}) that $N_\delta$ is $(\lambda,*)$-brimmed over
$M_\delta$, hence some $c \in N_\delta$ realizes $p$, so for some $i <
\delta,c \in N_i$ and this $\alpha$ is as required.
\hfill$\square_{\scite{705-gr.6}}$\margincite{705-gr.6}
\enddemo
\bn
We define rank as in \cite{Sh:394}.
\definition{\stag{705-10b.1} Definition} 
rk = rk$_{\frak s}$, e.g. is defined as follows:

rk$_{\frak s}(p)$ is defined if $p \in {\Cal S}_{\frak s}(M)$ for some
$M \in K_{\frak s}$

it is an ordinal or $\infty$

rk$_{\frak s}(p) \ge \alpha$ iff for every $\beta < \alpha$ we can find
$(M_1,p_1)$ such that
\mr
\item "{{}}"  $M \le_{\frak s} M_1,p_1 \in {\Cal S}_{\frak s}(M_1)$ is an
extension of $p$ which splits over $M$ and rk$_{\frak s}(p_1) \ge \beta$.
\ermn
Lastly, rk$_{\frak s}(p) = \alpha$ iff rk$_{\frak s}(p) \ge \alpha$ and
rk$_{\frak s}(p) \ngeq \alpha + 1$.
\enddefinition
\bn
Basic properties of rk$_{\frak s}$ are
\proclaim{\stag{705-10b.4} Claim}  Assume ${\frak s}$ is weakly successful
good.
If $M \in K_{\frak s}$ and $p \in {\Cal S}_{\frak s}(M)$, 
\ub{then} {\rm rk}$_{\frak s}(p) < \infty$.
\endproclaim
\bigskip

\remark{Remark}  So this applies in \scite{705-6a.1} to ${\frak s}^*$ but
rk$_{{\frak s}(+)} = \text{ rk}_{{\frak s}(*)}$ so it applies to
${\frak s}^+$, too.
\endremark
\bigskip

\demo{Proof}  Assume rk$_{\frak s}(p) = \infty$ we can choose by
induction on $n$ a triple $(M_n,N_n,a),M_n \le_{\frak s} N_n,a \in
N_n$, rk$_{\frak s}(\text{tp}(a,M_n,N_n)) = \infty$ and $M_n
\le_{\frak s} M_{n+1},N_n \le_{\frak s} N_{n+1}$ and
tp$(a,M_{n+1},N_{n+1})$ does $\lambda$-split over $M_n$ (in the induction
step we use amalgamation and having $\le 2^{\lambda_{\frak s}}$
possible isomorphism types for $(M_{n+1},N_{n+1},a)$ over $M_n$).
Clearly we can find $\langle N^+_n:n < \omega \rangle$ such that
$M_n \le_{\frak s} N^+_n$ and NF$_{\frak s}
(M_n,M_{n+1},N^+_n,N^+_n),N^+_{n+1}$ is $(\lambda,*)$-brimmed over
$M_{n+1} \cup N_n$.  By \scite{705-1.15K} we know that 
$N^+_\omega = \cup \{N^+_n:n <
\omega\}$ is $(\lambda,*)$-brimmed over $M_\omega = \cup\{M_n:n <
\omega\}$, hence we can embed $N_\omega = \cup \{N_n:n < \omega\}$ into
$N^+_\omega$ over $M_n$ so \wilog \, $n < \omega \Rightarrow N_n
\le_{\frak s} N^+_\omega$.  So for some $n < \omega$ we have $a \in
N^+_n$, and by long transitivity for NF we have NF$_{\frak
s}(M,N^+_n,M_\omega,N^+_\omega)$.  We get easy contradiction to
tp$(a,M_{n+1},N^+_{n+1}) = \text{ tp}_{\frak s}(a,M_{n+1},N^+_\omega)
= \text{ tp}(a,M_{n+1},N_{n+1})$ does $\lambda$-split over $M_n$.
\hfill$\square_{\scite{705-10b.4}}$\margincite{705-10b.4}
\enddemo
\bigskip

\remark{\stag{705-10b.5} Remark}  An important point is that for any $\langle M_i:i
\le \delta \rangle$ which is $\le_{\frak s}$-increasing continuous and
$p_i \in {\Cal S}_{\frak s}(M_i)$ for $i < \delta$ such that $i < j \Rightarrow
p_i = p_j \restriction M_i$ in general there is no $p \in {\Cal
S}_{\frak s}(\cup\{M_i:i < \delta\})$ such that $i < \delta \Rightarrow p_i = p
\restriction M_i$, but for $\delta = \omega$ there is.
\endremark
\bigskip

\proclaim{\stag{705-gr.12} Claim}  1)  {\rm rk}$_{\frak s}(p)$ is a well
defined ordinal ($< \infty$) \ub{if}
$p \in {\Cal S}_{\frak s}(M),M \in K_{\frak s}$. \nl
2) If $M <_{\frak s} N$ and 
$p \in {\Cal S}_{\frak s}(N)$ splits over $M$ and
{\rm rk}$_{\frak s}(p) < \infty$, \ub{then} 
{\rm rk}$_{\frak s}(p) < { \text{\rm rk\/}}_{\frak s}
(p \restriction M)$. \nl
3) If {\rm NF}$_{\frak s}(M_0,M_1,M_2,M_3)$ and
$a \in M_4$, \ub{then} {\rm rk}$_{\frak s}
({\text{\rm tp\/}}_{\frak s}(a,M_0,M_3)) =
{ \text{\rm rk\/}}_{\frak s}({\text{\rm tp\/}}_{\frak s}(a,M_1,M_3))$. 
\nl
4) If $M \le_{\frak s} N$ and $p \in
{\Cal S}^{\text{bs}}(N)$ does not fork over 
$M$ \ub{then} $p$ does not split over
$M$ and {\rm rk}$_{\frak s}(p) = 
{ \text{\rm rk\/}}_{\frak s}(p \restriction M)$. 
\endproclaim
\bigskip

\demo{Proof}  1) Immediate by \scite{705-gr.6}. \nl
2) - 4) Easy.  \hfill$\square_{\scite{705-gr.12}}$\margincite{705-gr.12}
\enddemo
\bn
We may like to translate ranks between ${\frak s}$ and ${\frak s}^+$.
\proclaim{\stag{705-10b.3} Claim}  [${\frak s}$ is a successful good$^+$-frame]

Assume $N_1 <_{{\frak K}[{\frak s}]} M_1,
N_1 \in K_{\frak s},M_1 \in K_{{\frak s}(+)},p \in {\Cal S}_{{\frak s}(+)}(M_1)$. \nl
1) If $p$ does not $\lambda$-split over $N_1$ for ${\frak s}$, 
\ub{then} {\rm rk}$_{{\frak s}(+)}(p) = { \text{\rm rk\/}}_{\frak s}(p \restriction N_1)$. \nl
2) Also the inverse holds. \nl
3) If $p \in {\Cal S}^{\text{bs}}_{{\frak s}(+)}(M_1)$ and $N_1$
witnesses it then $p$ does not $\lambda$-split over $N_1$ and moreover
does not split over $N_1$. \nl
4) If $p \in {\Cal S}_{{\frak s}(+)}(M_1)$ \ub{then} for some $N_0
<_{{\frak K}[{\frak s}]} M_1$ of cardinality $\lambda,p$ does not
split over $N_0$ and even does not split over $N_0$; we call such $N_0$
a witness for $p$. \nl
5) If $p \in {\Cal S}^{\text{bs}}_{{\frak s}(+)}(M_1)$ \ub{then} $N_1$
is a witness for $p \in {\Cal S}^{\text{bs}}_{{\frak s}(+)}(M_1)$ iff
$p$ does not $\lambda$-split over $N_1$.
\endproclaim
\bigskip

\remark{\stag{705-10b.3A} Remark}  No real harm in assuming ``${\frak s}$
is type full" (see Definition \scite{705-stg.13}).
\endremark
\bn
\margintag{705-10b.3B}\ub{\stag{705-10b.3B} Conclusion}  If $M_0 <_{{\frak s}(+)} M_1$ and $p \in
{\Cal S}^{\text{bs}}_{{\frak s}(+)} (M_1)$ 
does ${\frak s}(*)$-fork over $M_0$
\ub{then} rk$_{{\frak s}(+)}(p) < \text{ rk}_{{\frak s}(+)}
(p \restriction M_0)$.
\bigskip

\demo{Proof of \scite{705-10b.3}}  3),4),5) should be clear. \nl
1) We prove by induction $\alpha$ that
\mr
\item "{$\circledast_\alpha$}"  for any such $(N_1,M_1,p)$ we have \nl
rk$_{{\frak s}(+)}(p) \ge \alpha \Leftrightarrow \text{ rk}_{\frak s}
(p \restriction N_1) \ge \alpha$
\ermn
This clearly suffices. \nl
For $\alpha = 0$ and $\alpha$ limit there are no problems.  So assume
$\alpha = \beta + 1$.  First assume rk$_{{\frak s}(+)}
(p) \ge \alpha$ hence
by the definition of rk$_{{\frak s}(+)}$ we can find $p,M_2$ such that
$M_1 \le_{{\frak s}(+)} M_2,q \in 
{\Cal S}_{{\frak s}(+)}(M_2),q \restriction
M_1 = p$ and $q$ does $\lambda^+$-split over $M_1$ and rk$_{{\frak
s}(+)}(q) \ge \beta$.  Hence $q$ is not witnessed by $N_1$ 
hence $q$ does $\lambda$-split over
$N_1$ hence for some $N_2 \in K_{\frak s}$ we have 
$N_1 \le_{\frak s} N_2 \le_{{\frak K}[{\frak s}]} M_2$ and
$q \restriction N_2$ does $\lambda$-split over $N_1$ 
for ${\frak K}_{\frak s}$ and
\wilog \, $N_2$ is a witness for $q$.  So by the induction hypothesis
rk$_{{\frak s}(+)}(q) \ge \beta \Leftrightarrow \text{ rk}_{\frak s}(q
\restriction N_2) \ge \beta$.

But by the choice of $q$, rk$_{{\frak s}(+)}(q) \ge \beta$ hence
rk$_{\frak s}(q \restriction N_2) \ge \beta$.  By the definition of
rk$_{\frak s}$, as $q \restriction N_2$ does $\lambda$-split over $N_1$ for
${\frak s}$, we get rk$_{\frak s}
(p \restriction N_1) > \text{ rk}_{\frak s}
(q \restriction N_2) \ge \beta$ so rk$_{\frak s}
(p \restriction N_1) \ge \beta +1 = \alpha$ as required.

Second assume rk$_{\frak s}(p \restriction N_1) \ge \alpha$ so we can find
$N_2,N_3,a$ such that $N_1 \le_{\frak s} N_2 \le_{\frak s} N_3,a \in N_3,
q = \text{ tp}(a,N_2,N_3)$ is a $\lambda$-splitting (for ${\frak s}$) extension of $p^- = 
p \restriction N_1$.  We use NF amalgamation to lift this to
$M_2,p$. \nl
2) It is enough to prove rk$_{{\frak s}(+)}(p) < \text{ rk}_{\frak
s}(p \restriction N_1)$ assuming the $p$ does $\lambda$-split over
$N_1$.  Now we can find $N_2 \in K_{\frak s}$ such that $N_1
\le_{\frak s} N_2 \le {{\frak K}[{\frak s}]}-M_1$ and $p$ does not
$\lambda$-split over $N_2$ but $p \restriction N_2$ does
$\lambda$-split over $N_1$.  So by part (1) we have rk$_{{\frak
s}(+)}(p) = \text{ rk}_{\frak s}(p \restriction N_2)$, and by the
definition of rk$_{\frak s}$ we know that rk$_{\frak s}(p \restriction
N_2) \le \text{ rk}_{\frak s}(p \restriction N_1)$.  
Together we are done. \nl
3),4) Left to the reader.
\hfill$\square_{\scite{705-10b.3}}$\margincite{705-10b.3}
\enddemo
\newpage

\head {\S3 primes triples} \endhead  \resetall \sectno=3
 \spuriousreset
\bigskip

\demo{\stag{705-c.0} Hypothesis}   ${\frak s}$ is a good $\lambda$-frame.
\enddemo
\bigskip

\definition{\stag{705-c.1} Definition}  1) Assume ${\frak s} = ({\frak K},
\nonfork{}{}_{},{\Cal S}^{\text{bs}})$ is a good $\lambda$-frame.  Let
$K^{3,\text{pr}}_\lambda = K^{3,\text{pr}}_{\frak s}$ 
be the family (pr stands for prime) of triples
$(M,N,a) \in K^{3,\text{bs}}_\lambda = K^{3,\text{bs}}_{\frak s}$ 
such that: if 
$(M,N',a') \in 
K^{3,\text{bs}}_\lambda$ and tp$(a,M,N) = \text{ tp}(a',M,N')$
\ub{then} there is a $\le_{\frak s}$-embedding $f:N \rightarrow N'$ 
over $M$ satisfying $f(a) = a'$. 
So such triples are called prime.
\nl
2) We say that
${\frak s} = ({\frak K},\nonfork{}{}_{},{\Cal S}^{\text{bs}})$ is
$\lambda$-good$^2$ or that ${\frak s}$ has primes \ub{if} 
${\frak s}$ is $\lambda$-good and
\mr
\item "{$(a)$}"  if 
$M \in K_\lambda,p \in {\Cal S}^{\text{bs}}(M)$ \ub{then}
for some $N,a$ we have $(M,N,a) \in K^{3,\text{pr}}_\lambda$ 
and $p = \text{ tp}(a,M,N)$.
\ermn
3) $(M,N,a)$ is model-minimal \ub{if} it belongs to 
$K^{3,\text{bs}}_\lambda$ and there is no $N'$ such that
$M <_{\frak s} N' <_{\frak s} N$ and $a \in N'$ (this notion is
close to ``tp$(a,M,N)$ is of depth zero,
$N$ prime over $M \cup \{a\}$" in the context of \cite{Sh:c}).
\nl
4) We say ${\frak s}$ has [model]-minimality \ub{if} for every $M \in
K_{\frak s},p \in {\Cal S}^{\text{bs}}(M)$ there is $(M,N,a) \in
K^{3,\text{bs}}_\lambda$ in which $a$ realizes $p$ and 
$(M,N,a)$ is [model]-minimal (see Definition \scite{705-stg.7}). 
\enddefinition
\bigskip

\definition{\stag{705-c.1A} Definition}  1) We say $\langle M_i,a_j:i \le \alpha,
j < \alpha \rangle$ is a pr-decomposition of $N$ over $M$
or of $(M,N)$ if: $M_i$ is
$\le_{\frak s}$-increasing continuous, $(M_i,M_{i+1},a_i) \in 
K^{3,\text{pr}}_{\frak s},M_0 = M$ and $M_\alpha = N$; we may allow $N \in K^{\frak
s}_{\lambda^+}$ but $i < \alpha \Rightarrow M_i \in K_{\frak s}$. 
If we demand just $M_\alpha \le_{\frak s} N$ we say ``inside $N$"
instead of ``of $N$".  If we also allow $M \le_{\frak s} M_0,
M_\alpha \le_{\frak s} N$ we
say in $(M,N)$. Instead ``over $M$" we can say $M$-based.
We call $\alpha$ the length of the decomposition.
\nl
2) Similarly for uq ($K^{3,\text{uq}}_\lambda$ is from
\marginbf{!!}{\cprefix{600}.\scite{600-nu.1A}}) and we define uq-decomposition.  
We may write just decomposition (or 
${\frak s}$-decomposition) instead pr-decomposition.
\enddefinition
\bigskip

\proclaim{\stag{705-c.2} Claim}  1) If $(M,N,a) 
\in K^{3,\text{pr}}_\lambda$ and
$M \cup \{a\} \subseteq N' \le_{\frak s} N$ \ub{then} $(M,N,a) \in
K^{3,\text{pr}}_\lambda$. \nl
2) Similarly for $K^{3,\text{uq}}_\lambda$. \nl 
3) If $(M,N_1,a_1) \in K^{3,\text{bs}}_\lambda$ is model-minimal and
$(M,N_2,a_2) \in  K^{3,\text{pr}}_\lambda$ and 
$p = { \text{\rm tp\/}}(a_1,M,N_1)
= { \text{\rm tp\/}}(a_2,M,N_2)$ \ub{then} there is an isomorphism from
$N_1$ onto $N_2$ over $M$, mapping $a_1$ to $a_2$ (so both triples are 
model-minimal and prime and so if $(M,N',a')$ is prime or is 
model minimal with {\rm tp}$(a',M,N') = { \text{\rm tp\/}}
(a_\ell,M,N_\ell)$ then for $\ell =1,2$
there is an isomorphism $f_\ell$ from $N'$ onto $N_\ell$ mapping $a'$
to $a_\ell$ and being the identity on $M$). \nl
4) Assume ${\frak s}$ is weakly successful. 
If $M_0 \le_{\frak s} M_\ell \le_{\frak s} M_3$ for $\ell =1,2$ and
$(M_0,M_1,a)$ belong to $K^{3,\text{uq}}_{\frak s}$ and 
{\rm tp}$(a,M_2,M_3)$ does
not fork over $M_0$ (e.g. {\rm tp}$(a,M_0,M_1)$ has unique extension in
${\Cal S}(M_2)$) \ub{then} NF$_{\frak s}(M_0,M_1,M_2,M_3)$. \nl
5) If $M_0 \le_{\frak s} M_1 \le_{\frak s} M_2,a_\ell \in M_{\ell +1}$
and {\rm tp}$_{\frak s}(a_\ell,M_\ell,M_{\ell +1}) \in {\Cal
S}^{\text{bs}}
(M_\ell)$ does not fork over $M_0$ for $\ell = 0,1$ \ub{then}
$(M_0,M_2,a_0) \notin K^{3,\text{uq}}_\lambda$.
\endproclaim
\bigskip

\demo{Proof}  Easy (e.g. (3) is \scite{705-1.15M} and (4) is by the definition of 
$K^{3,\text{uq}}_{\frak s}$ and the 
existence of NF$_{\frak s}$-amalgamation). \nl
${{}}$  \hfill$\square_{\scite{705-c.2}}$\margincite{705-c.2}
\enddemo
\bigskip

\proclaim{\stag{705-c.4} Claim}  1) Assume that ${\frak s}$ has primes; if
$(M,N,a) \in K^{3,\text{uq}}_\lambda$ \ub{then} for some $N',M \cup 
\{a\} \subseteq N' \le_{\frak K} N$ and $(M,N',a) \in 
K^{3,\text{pr}}_\lambda \cap K^{3,\text{uq}}_\lambda$. \nl
2) If $(M,N,a) \in K^{3,\text{pr}}_\lambda$ and $K^{3,\text{uq}}_\lambda$
is dense (e.g. if ${\frak s}$ is weakly successful) \ub{then} 
$(M,N,a) \in K^{3,\text{uq}}_\lambda$. 
\endproclaim
\bigskip

\demo{Proof}  Immediate: part (1) by the definition and monotonicity
of $K^{3,\text{uq}}_\lambda$, part (2) by the
proof of part (1).  \hfill$\square_{\scite{705-c.4}}$\margincite{705-c.4}
\enddemo
\bigskip

\proclaim{\stag{705-c.3} Claim}  If 
${\frak s}$ is non-multi-dimensional weakly successful and has primes
\ub{then} ${\frak s}$ has model-minimality and all $(M,N,a) \in 
K^{3,\text{pr}}_{\frak s}$ are model minimal.
\endproclaim
\bigskip

\demo{Proof}  Let 
$M \in K_{\frak s}$ and $p \in {\Cal S}^{\text{bs}}_{\frak s}(M)$.  
We know (by Definition \scite{705-c.1}(2)(a)) 
that there is $(M,N_2,a) \in K^{3,\text{bs}}_\lambda$ which is prime 
and $p = \text{ tp}(a,M,N_2)$.  If $(M,N_2,a)$ is model-minimal 
we are done, otherwise there is $N_1$ satisfying 
$M \cup \{a\} \subseteq N_1 <_{\frak s} N_2$.
As $(M,N_2,a)$ is prime there is an $\le_{\frak s}$-embedding 
$f$ of $N_2$ into $N_1$ over
$M \cup \{a\}$, let $N_0 = f(N_2)$ hence $M \cup \{a\} \subseteq N_0 
\le_{\frak s} N_1$.  So $M \cup \{a\} \subseteq N_0 <_{\frak s} N_2$, 
so by non-multi-dimensionality
there is $b \in N_2 \backslash N_0$ such that tp$(b,N_0,N_2) \in
{\Cal S}^{\text{bs}}_{\frak s}(N_0)$ 
does not fork over $M$ hence by \scite{705-c.2}(5) we have $(M_a,N_2,a)
\notin K^{3,\text{uq}}_\lambda$ (the $M,N_0,N_2,a,b$ here correspond
to $M_0,M_1,M_2,a_0,a_1$.  This easily 
contradicts ``$(M,N_2,a) \in K^{3,\text{pr}}_{\frak s} 
\subseteq K^{3,\text{uq}}_{\frak s}$" which holds by
\scite{705-c.4}(2).  \hfill$\square_{\scite{705-c.3}}$\margincite{705-c.3}
\enddemo
\bigskip

\proclaim{\stag{705-zm.1} Claim}  [${\frak s}$ is a good, weakly
successful $\lambda$-frame].
\nl
1) Assume $M_0 \le_{\frak s} M_\ell \le_{\frak s} M_3,a_\ell \in
M_\ell$ for $\ell =1,2$ and
$(M_0,M_\ell,a_\ell) \in K^{3,\text{uq}}_{\frak s}$ for $\ell = 1,2$.

\ub{Then} {\rm tp}$(a_2,M_1,M_3)$
does not fork over $M_0$ \ub{iff} {\rm tp}
$(a_1,M_2,M_3)$ does not fork over $M_0$. \nl
2) Assume $M_0 \le_{\frak s} M_\ell \le_{\frak s} M_3$ and $a_\ell \in
M_\ell$ for $\ell = 1,2$ and $(M_0,M_1,a_1) \in 
K^{3,\text{uq}}_{\frak s}$
and {\rm tp}$(a_2,M_0,M_2) \in {\Cal S}^{\text{bs}}(M_0)$.
If {\rm tp}$_{\frak s}(a_1,M_2,M_3)$ does not fork over $M_0$ then
{\rm tp}$_{\frak s}(a_2,M_1,M_3)$ does not fork over $M_0$.
\endproclaim
\bigskip

\demo{Proof}  1) By the symmetry in the claim it is 
enough to prove the if part, so assume that
tp$(a_1,M_2,M_3)$ does not fork over $M_0$ (see \sectioncite[\S5]{600} on
$\le_{\text{bs}}$).  As $(M_0,M_1,a_1) \in K^{3,\text{uq}}_{\frak s}$ by
\scite{705-c.2}(4) it follows that
NF$_{\frak s}(M_0,M_1,M_2,M_3)$, hence by symmetry of 
NF$_{\frak s}$ (see \marginbf{!!}{\cprefix{600}.\scite{600-nf.14}}) we have 
NF$_{\frak s}(M_0,M_2,M_1,M_2)$ which implies that 
tp$(a_2,M_1,M_3)$ does not fork over $M_0$ by \scite{705-stg.32}.  \nl
2) The proof is included in proof of part (1).  
\hfill$\square_{\scite{705-zm.1}}$\margincite{705-zm.1}
\enddemo
\bigskip

\proclaim{\stag{705-c.6} Claim}  Assume ${\frak s}$ has primes.
\mr
\item  If $M \le_{\frak s} N$ \ub{then} there is a 
decomposition of $N$ over 
$M$ (see Definition \scite{705-c.1A}(1),(2)).  Moreover, if $(M,N,a) \in
K^{3,\text{bs}}_{\frak s}$ \ub{then} without loss of generality $a_0=a$.
\sn
\item  If $M \le_{{\frak K}[{\frak s}]} N, M \in K_{\frak s},N \in
K_{\lambda^+}^{\frak s}$, \ub{then} there is a decomposition of $N$
over $M$ 
\sn
\item   If $N \in K^{\frak s}_{\lambda^+}$ the length of the
decomposition is $\lambda^+$
\sn
\item  In part (1) there is a decomposition of $N$ over $M$ of length
$\le \lambda$
\sn
\item  In part (1) if $N$ is $(\lambda,*)$-brimmed over $M$,
\ub{then} there is a decomposition of $N$ over $M$ of length exactly $\lambda$.
\endroster
\endproclaim
\bigskip

\demo{Proof}  1) By the 
definition of $\lambda$-good frame there is $a \in N$ such
that tp$_{\frak s}(a,M,N) \in {\Cal S}^{\text{bs}}_{\frak s}(M)$, so
it is enough to prove the second sentence, so \wilog \, $a$ is well defined.  
Choose $a_i,M_i$ by induction $i < \lambda^+$.  
Arriving to $i$, if $i=0,M_i = M,
a_i = a$. If $i$ is limit let $M_i = \cup \{M_j:j<i\}$: 
if $M_i = N$ we are done, if not then 
for some $a_i \in N \backslash M_i$ we have tp$_{\frak s}(a_i,M_i,N) \in
{\Cal S}^{\text{bs}}_{\frak s}(M_i)$.  
If $i = j+1$ then $M_j,a_j$ are well defined and we
know (as ${\frak s}$ has primes) that there is $M_{j+1}
\le_{\frak s} N$ such that $(M_j,M_{j+1},a_j) \in K^{3,\text{pr}}_{\frak s}$;
again if $M_i=N$ we are done and otherwise we can choose $a_i \in N
\backslash M_i$ such that tp$(a_i,M_i,N) \in {\Cal S}^{\text{bs}}(M_i)$.  So
by cardinality consideration at some point we are stuck, i.e., $M_i = N$.
\nl
2)  Let $\langle b_\varepsilon:\varepsilon < \lambda^+ \rangle$
list the elements of $N$.  Repeating the proof of part (1), now in
choosing $a_i$ when $i>0$ we 
can choose any $a \in \bold I_i = \{a \in N \backslash
M_i:\text{tp}_{\frak s}(a,M_i,N) \in {\Cal S}^{\text{bs}}_{\frak s}
(M_i)\}$ so we can demand that
$a_i = b_{\varepsilon_1} \and b_{\varepsilon_2} \in \bold I_i \Rightarrow
\varepsilon_1 \le \varepsilon_2$.  It suffice to show 
that $M_{\lambda^+} = \cup\{M_i:i < \lambda^+\}$ is not equal to $N$,
obviously $M_{\lambda^+} \le_{{\frak K}[{\frak s}]} N$.
Otherwise we can find $a \in N \backslash
M_{\lambda^+}$ such that $S = \{i:\text{tp}(a,M_i,N) \in {\Cal
S}^{\text{bs}}(M_i)\}$ is stationary and \wilog \, $i \in S \Rightarrow
\text{ tp}(a,M_i,N)$ does not fork over $M_{i(*)},i(*) = \text{
Min}(S)$ so \wilog \, $S = [i(*),\lambda^+)$.
The type tp$_{\frak s}(a,M_{\lambda^+},N)$ does not fork over 
$M_{i(*)}$, so $i \in [1 + i(*),
\lambda^+) \Rightarrow a \in \bold I_i$, so if $a =
b_{\varepsilon(*)}$ then $i \in [1+i(*),\lambda^+) \Rightarrow a_i \in
\{b_\varepsilon:\varepsilon < \varepsilon(*)\}$, so we have a 1-to-1 function
from $[1+i(*),\lambda)$ into $[0,\varepsilon(*))$, contradiction. \nl
3)-5) Left to the reader.  \hfill$\square_{\scite{705-c.6}}$\margincite{705-c.6}
\enddemo
\bigskip

\proclaim{\stag{705-6.xob1} Claim}  1) [${\frak s}$ is a (good)
weakly successful $\lambda$-frame with primes]. \nl
If ${\frak C} \in K^{\frak s}_{\lambda^+}$ is $\lambda^+$-saturated
(over $\lambda$ of course), $M \in K_{\frak s},M \le_{{\frak K}[{\frak
s}]} {\frak C}$ and $a_1,a_2 \in {\frak C}$ satisfy
{\rm tp}$(a_\ell,M,{\frak C}) \in {\Cal S}^{\text{bs}}(M)$ 
for $\ell=1,2$,  \ub{then} the following are equivalent:
\mr
\item "{$(a)$}"  there are $M_1,M_2$ from $K_{\frak s}$ such 
that {\rm NF}$_{\frak s}(M,M_1,M_2,{\frak C})$ and $a_1 \in 
M_1,a_2 \in M_2$
(the meaning of {\rm NF} above is 
for some $M_3 \le_{{\frak K}[{\frak s}]} {\frak C}$ from
$K_{\frak s}$ we have {\rm NF}$_{\frak s}(M,M_1,M_2,M_3))$
\sn
\item "{$(b)_\ell$}"  there is $M_\ell \le_{{\frak K}
[{\frak s}]} {\frak C}$
from $K_{\frak s}$ satisfying 
$M \le_{\frak s} M_\ell \le_{{\frak K}[{\frak s}]}
{\frak C}$ such that 
$a_\ell \in M_\ell$ and {\rm tp}$(a_{3 - \ell},M_\ell,
{\frak C})$ does not fork over $M$
\sn
\item "{$(c)_\ell$}"  if $(M,M_\ell,a_\ell) \in 
K^{3,\text{uq}}_\lambda$ and
$M_\ell \le_{{\frak K}[{\frak s}]} {\frak C}$ \ub{then} 
{\rm tp}$(a_{3 - \ell},M_\ell,{\frak C})$ does not fork over $M$
\sn
\item "{$(d)_\ell$}"  if $(M,M_\ell,a_\ell) \in 
K^{3,\text{pr}}_\lambda$ and
$M_\ell \le_{{\frak K}[{\frak s}]} {\frak C}$ \ub{then} 
{\rm tp}$(a_{3 - \ell},M_\ell,{\frak C})$ does not fork over $M$.
\ermn
2) [${\frak s}$ is a (good) weakly successful $\lambda$-frame.]  Above
$(a) \Leftrightarrow (b)_\ell \Leftrightarrow (c)_\ell \Rightarrow (d)_\ell$.
\endproclaim
\bigskip

\demo{Proof}  1) 
\mn
\ub{$(a) \Rightarrow 
(b)_\ell$} by \scite{705-stg.32} (and the symmetry of NF).
\mn
\ub{$(b)_\ell \Rightarrow (a) + (c)_{3 -\ell}$}.  To prove 
$(c)_{3 - \ell}$ assume $(M,M_{3 - \ell},a_{3 - \ell}) \in
K^{3,\text{uq}}_\lambda$.  As we assume $(b)_\ell$ for some 
$M_\ell \le_{\frak K} {\frak C}$ in $K_\lambda$, we have 
tp$(a_{3 - \ell},M_\ell,{\frak C})$ does not fork over $M$ and $a_\ell
\in M_\ell$, so as $M \le_{\frak s} M_\ell$ and
tp$(a_\ell,M,M_\ell) \in {\Cal S}^{\text{bs}}(M)$ clearly
$(M,M_\ell,a_\ell) \in K^{3,\text{bs}}$.  By \scite{705-c.2}(4)
we have NF$_{\frak s}(M,M_\ell,M_{3 - \ell},
{\frak C})$ hence by \scite{705-stg.32} the desired conclusion of
$(c)_{3 - \ell}$ holds.  This proves also clause (a) if we note, as
${\frak s}$ is weakly successful for some $(M,N,b) \in
K^{3,\text{uq}}$, tp$(b,M,N) = \text{ tp}(a_{3 -\ell},M,{\frak C})$,
so as ${\frak C}$ is $\lambda^+$-saturated?? \wilog \, $a_{3 -\ell} =
b,N \le_{{\frak K}[{\frak s}]} {\frak C}$.
\mn
\ub{$(c)_{3-\ell} \Rightarrow (d)_{3-\ell}$}: to 
prove $(d)_{3-\ell}$ assume $(M,M_\ell,a_\ell) \in
K^{3,\text{pr}}_\lambda$ and $M_\ell \le_{{\frak K}
[{\frak s}]} {\frak C}$; now ``${\frak s}$ is weakly successful" and
\scite{705-c.4}(2) implies $(M,M_\ell,a_\ell) \in
K^{3,\text{uq}}_\lambda$, and we can apply clause $(c)_{3-\ell}$ to get the
desired conclusion of $(d)_{3 -\ell}$.
\mn
\ub{$(d)_{3 - \ell} \Rightarrow (b)_{3 -\ell}$}: as ${\frak s}$ has
primes there is $M_{3-\ell}$ such that $(M,M_{3 - \ell},a_{3 - \ell}) 
\in K^{3,\text{pr}}_\lambda$ and use \scite{705-c.4}(2).
\nl
Clearly those implications are enough. \nl
2) The proof is included in the proof of part (1) except
\ub{$(c)_\ell \Rightarrow (b)_\ell$} like $(d)_{3 -\ell} \Rightarrow
(b)_{3 - \ell}$ using ``weakly successful".
\hfill$\square_{\scite{705-6.xob1}}$\margincite{705-6.xob1}
\enddemo
\bigskip

\proclaim{\stag{705-e.1} Claim}  Assume 
${\frak s}$ is good$^+$ and $n$-successful and $n>0$. \nl
1) ${\frak s}^{+n}$ is a $\lambda^{+n}$-good$^+$ frame. \nl
2) ${\Cal S}^{\text{bs}}_{{\frak s}[+n]} 
= {\Cal S}^{\text{bs}}_{{\frak s}<\lambda^{+n}>}
\restriction K^{{\frak s}^{+n}}_{\lambda^{+n}}$ 
(see \scite{705-0.X}, the ${\Cal S}^{\text{bs}}_{\frak s}$ is from
\sectioncite[\S2]{600} and is $\{p \in {\Cal S}^{\text{bs}}
(M,{\frak s}):M \in K^{{\frak s}^{+n}}_{\lambda^{+n}}\})$.
\nl
3) If $M \in K^{{\frak s}^{+n}}_{\lambda^{+n}}$ and $p \in 
{\Cal S}^{\text{bs}}_{{\frak s}[+n]}(M)$ \ub{then} for some
$(M,N,a) \in K^{3,\text{pr}}_{\lambda^{+n}}[{\frak s}^{+n}]$ 
we have {\rm tp}$(a,M,N) = p$. \nl
4) If $(M,N,a) \in K^{3,\text{pr}}_{\lambda^{+n}}
[{\frak s}^{+n}]$ \ub{then} $(M,N,a) \in
K^{3,\text{uq}}_{\lambda^{+n}}[{\frak s}^{+n}]$. \nl
5) If $M \le_{\frak K} N^1 \le_{\frak K} N^2$ are in 
$K^{{\frak s}^{+n}}_{\lambda^{+n}}$ and $a \in N^1$ \ub{then}
$(M,N^2,a) \in K^{3,\text{uq}}_{\lambda^{+n}}[{\frak s}^{+n}] 
\Rightarrow (M,N^1,a) \in K^{3,\text{uq}}_{\lambda^{+n}}
[{\frak s}^{+n}]$ and 
$(M,N^2,a) \in K^{3,\text{pr}}_{\lambda^{+n}}[{\frak s}^{+n}]
\Rightarrow (M,N^1,a) \in K^{3,\text{pr}}_{\lambda^{+n}}
[{\frak s}^{+n}] \Rightarrow (M,N^1,a) \in 
K^{3,\text{uq}}_{\lambda^{+n}}[{\frak s}^{+n}]$. \nl
6) Assume $n=m+1$, and $(M_0,M_1,a) \in K^{3,\text{bs}}_{\lambda^{+n}}
[{\frak s}^{+n}]$ and $\bar M_\ell = \langle M_{\ell,\alpha}:
\alpha < \lambda^{+n} \rangle$ a $\le_{\frak K}$-representation of 
$M_\ell$ for $\ell = 1,2$.
\ub{Then}:
\mr
\item "{$(*)$}"  $(M_0,M_1,a) 
\in K^{3,\text{pr}}_{\lambda^{+n}}[{\frak s}^{+n}]$ 
\ub{iff} for some club $E$ of $\lambda^+$ we have: 
for $\alpha < \beta$ in $E,(M_{0,\alpha},M_{1,\alpha},
a) \in K^{3,\text{uq}}_\lambda[{\frak s}^{+m}]$ 
and $(M_{0,\alpha},M_{1,\alpha},a) <^{*,{\frak s}^{+m}}_{\text{bs}} 
(M_{0,\beta,M_{1,\beta}},a)$, see \scite{705-d.0A}. 
\endroster
\endproclaim
\bigskip

\demo{Proof}  Straight; all by induction on $n$; part (3), (6) by
\scite{705-d.8} + \scite{705-d.1} below, part (4) by \scite{705-c.4}, part (5) 
by \scite{705-c.2}(1),\scite{705-c.2}(2).   \hfill$\square_{\scite{705-e.1}}$\margincite{705-e.1}
\enddemo
\bigskip

\remark{\stag{705-e.2} Remark}  1) If we assume ${\frak s}_0$ is
unidimensional (see \S2),
life is easier: $(M,N,a) \in K^{3,\text{pr}}_{\lambda^{+n}}$ 
implies model-minimality, see \S3, \scite{705-c.3}.  
On categoricity see \scite{705-2e.3}.
\nl
2) For \scite{705-e.1}(6), note that $M_1$ is saturated (in 
${\frak K}^{{\frak s}^{+(n-1)}}$ above $\lambda^{+n-1}$) if $(M,M_2,b)
\in K^{3,\text{bs}}_{\lambda^{+n}}[{\frak s}^{+n}]$ and tp$_{{\frak
s}[+n]}(b,M_0,M_2) = \text{ tp}_{{\frak s}[+n]}(a,M_0,M_1)$ then we
can choose an $\le_{K[{\frak s}]}$-embedding $f_\alpha$ of
$M_{1,\alpha}$ into $M_{2,\alpha}$, increasing continuous with
$\alpha$, mapping $a$ to $b$.  For $\alpha =0$ use the saturation, for
$\alpha = \beta +1$ use $(M_{1,\alpha},M_{2,\alpha},a) \in
K^{3,\text{uq}}_{{\frak s}[n-1]} +$ saturation. 
\endremark
\bigskip

\proclaim{\stag{705-2e.3} Claim} 
If ${\frak s}$ is an $n$-successful 
good$^+ \lambda$-frame and 
weakly unidimensional and categorical in $\lambda$, \ub{then}
\mr
\item "{$(i)$}"  ${\frak K}^{{\frak s}(+n)} = {\frak K}^{\frak s}
_{\ge \lambda^{+n}}$
\sn
\item "{$(ii)$}"  ${\frak s}(+n)$ is weakly unidimensional, and
categorical in $\lambda^{+n}$.
\endroster
\endproclaim
\bigskip

\demo{Proof}  By \scite{705-e.1} and \scite{705-stg.10}.
\hfill$\square_{\scite{705-2e.3}}$\margincite{705-2e.3} 
\enddemo
\newpage

\head {\S4 Prime existence} \endhead  \resetall \sectno=4
 \spuriousreset
\bigskip

We give some easy properties of primes for ${\frak s}^+$.  A major
point is \scite{705-d.8}: existence of primes.  We also note how various
properties reflect from $K_{{\frak s}^+}$ to $K_{\frak s}$.
\bigskip

\demo{\stag{705-d.0} Hypothesis}  1) ${\frak s} = ({\frak K}_{\frak s},
\nonfork{}{}_{},{\Cal S}^{\text{bs}})$ is a successful 
good$^+ \lambda$-frame, ${\frak K} = {\frak K}[{\frak s}]$ as usual.
\enddemo
\bn
Recall
\definition{\stag{705-d.0A} Definition}  1) We let $\le_{\text{bs}} =
\le^{\frak s}_{\text{bs}}$ be the following relation (really quasi order)
on $K^{3,\text{bs}}_\lambda:(M,N,a) \le_{\text{bs}}
(M',N',a)$ if both are in $K^{3,\text{bs}}_{\frak s},
M \le_{\frak s} M',N \le_{\frak s} N'$ and tp$(a,M',N')$ does not fork over $M$. \nl
2)  $\le^*_{\text{bs}} = \le^{*,{\frak s}}_{\text{bs}}$ 
is the following quasi order
on $K^{3,\text{bs}}_\lambda:(M,N,a) 
\le^*_{\text{bs}} (M',N',a)$ if (they are in 
$K^{3,\text{bs}}_{\frak s}$ and) 
$(M,N,a) \le_{\text{bs}} (M',N',a)$ and if they are not equal
then $M',N'$ is universal over $M,N$ respectively (and $<^*_{\text{bs}}$ 
has the obvious meaning).
\enddefinition
\bigskip

\proclaim{\stag{705-d.1} Claim}  Assume $M_0 \in K_{{\frak s}(+)}$ and
$p \in {\Cal S}^{\text{bs}}_{{\frak s}(+)}(M_0)$.  \ub{Then} we can find
$a,\bar M_0,M_1,\bar M_1$ such that:
\mr
\widestnumber\item{$(viii)$}
\item "{$(i)$}"  $M_0 \le_{{\frak K}[{\frak s}]} 
M_1 \in K^{\frak s}_{\lambda^+}$
\sn
\item "{$(ii)$}"  $M_1 \in K^{\frak s}_{\lambda^+}$ is 
saturated, for ${\frak s}$, equivalently $M_1 \in K_{{\frak s}(+)}$
\sn
\item "{$(iii)$}"   $a \in M_1$ and $p = 
{ \text{\rm tp\/}}_{{\frak s}(+)}(a,M_0,M_1)$
\sn
\item "{$(iv)$}"  $\bar M_\ell = 
\langle M_{\ell,\alpha}:\alpha < \lambda^+
\rangle$ is a $\le_{{\frak K}[{\frak s}]}$-representation of $M_\ell$
for  $\ell=0,1$
\sn
\item "{$(v)$}"  $a \in M_{1,0}$
\sn
\item "{$(vi)$}"  $(M_{0,\alpha},M_{1,\alpha},a) \in 
K^{3,\text{uq}}_\lambda$ for every $\alpha < \lambda^+$
\sn
\item "{$(vii)$}"  $M_{\ell,i+1}$ is $(\lambda,*)$-brimmed over
$M_{\ell,i}$ for $i < \lambda^+,\ell < 2$
\sn
\item "{$(viii)$}"  $(M_{0,\alpha},M_{1,\alpha},a)$ is $<^{\frak
s}_{\text{bs}}$-increasing 
\endroster
\endproclaim
\bigskip

\definition{\stag{705-d.1A} Definition}  If we say 
$(M_0,M_1,a)$ is canonically
${\frak s}^+$-prime \ub{if} there are $\bar M^0,\bar M^1$ are as in
claim \scite{705-d.1} above (see \scite{705-d.8} below, formally this depends
on ${\frak s}$, but our ${\frak s}$ is constant).
\enddefinition
\bigskip

\demo{Proof}  Let $M_{0,0} \le_{\frak K} M_0,M_{0,0} \in K_\lambda$ be
such that $M_{0,0}$ is a witness for $p$ and $\langle b_\alpha:\alpha
< \lambda^+ \rangle$ list $|M_0|$

We choose by induction on $\alpha < \lambda^+$, a pair 
$(M_{0,\alpha},M_{1,\alpha})$ such that:
\mr
\item "{$(a)$}"  $(M_{0,\alpha},M_{1,\alpha},a) 
\in K^{3,\text{bs}}_\lambda$
\sn
\item "{$(b)$}"  $(M_{0,\beta},M_{1,\beta},a) \le_{\text{bs}} (M_{0,\alpha},
M_{1,\alpha},a)$ for $\beta < \alpha$
\sn
\item "{$(c)$}"  if $\alpha$ is a limit ordinal then $M_{\ell,\alpha} =
\dbcu_{\beta < \alpha} M_{\ell,\beta}$ for $\ell = 0,1$
\sn
\item "{$(d)$}"  for every even $\alpha$, if 
$(M_{0,\alpha},M_{1,\alpha},a) \notin
K^{3,\text{uq}}_\lambda$ then \nl
$\neg \text{ NF}_{\frak s}(M_{0,\alpha},
M_{1,\alpha},M_{0,\alpha +1},M_{1,\alpha +1})$
\sn
\item "{$(e)$}"  for odd $\alpha,M_{\ell,\alpha +1}$ is brimmed over
$M_{\ell,\alpha}$ for ${\frak s}$, for $\ell =1,2$.
\ermn
There is no problem to carry the definition (concerning clause (d), it
follows by \scite{705-d.2} below).
\enddemo
\bn
Before we continue note
\proclaim{\stag{705-d.2} Claim}  1) If $(M,N,a) \in K^{3,\text{bs}}_\lambda$
\ub{then} $(M,N,a) \notin K^{3,\text{uq}}_\lambda$ 
\ub{iff} for some $(M',N',a) \in K^{3,\text{bs}}_\lambda$ we 
have $(M,N,a) \le_{\text{bs}} (M',N',a)$ and 
$\neg {\text{\rm NF\/}}_{\frak s}(M,N,M',N')$. \nl
2) If $(M_\ell,N_\ell,a) <^*_{\text{bs}} (M_{\ell+1},N_{\ell +1},a)$ for
$\ell=0,1$ \ub{then} $(M_2,N_2,a)$ is universal over
$(M_0,N_0,a)$ for $\le_{\text{bs}}$.
\endproclaim
\bigskip

\demo{Proof}  By the definition of 
$K^{3,\text{uq}}_\lambda$ and the uniqueness and existence of
NF-amalgamation by \sectioncite[\S6]{600} (i.e., any 
two NF ones are compatible), 
the conclusion follows.   \hfill$\square_{\scite{705-d.2}}$\margincite{705-d.2}
\enddemo
\bn
\ub{Continuation of the proof of \scite{705-d.1}}

By clause (e), necessarily $M'_\ell =: \dbcu_{\alpha < \lambda^+} 
M_{\ell,\alpha} \in K^{\frak s}_{\lambda^+}$ are saturated for 
$\ell =0,1$.  Also $M'_0 \le_{{\frak K}[{\frak s}]} M'_1$ and by
clause (b) we have tp$_{{\frak s}(+)}(a,M'_0,M'_1) \in 
{\Cal S}^{\text{bs}}_{{\frak s}(+)}(M'_0)$ 
is a stationarization of $p \restriction M_{0,0}$, i.e., is witnessed
by it.  So without loss of generality 
$M'_0 = M_0$ and tp$(a,M'_0,M'_1) = p$ and by \scite{705-stg.3}, 
clause (b) we
have $M_0 \le^*_{\lambda^+} M_1$, so by its definition (see
\marginbf{!!}{\cprefix{600}.\scite{600-ne.1}}) for some club $E$ of $\lambda^+$ we have
$\alpha \in E \and \alpha < \beta \in E \Rightarrow 
\text{ NF}_{\frak s}(M_{0,\alpha},M_{1,\alpha},M_{0,\beta},
M_{1,\beta})$, hence by
monotonicity of NF$_{\frak s}$ and clause (d) of the construction we
have $\alpha \in E \Rightarrow
(M_{0,\alpha},M_{1,\alpha},a) \in K^{3,\text{uq}}_\lambda$.  By
renaming we get the conclusion.  \hfill$\square_{\scite{705-d.1}}$\margincite{705-d.1}
\bigskip
 
\proclaim{\stag{705-d.2A} Claim}
1) Assume $\beta < \lambda^+,\langle M_i:i \le \beta \rangle$ is
$\le_{\frak s}$-increasing continuous and $(M_i,M_{i+1},a_i) \in
K^{3,\text{bs}}_{\frak s}$ 
for $i < \beta$ and $M_0 \le_{\frak s} M^+$.  \ub{Then} we can
find $\langle N_i:i \le \beta \rangle$ such that:
\mr
\widestnumber\item{$(iii)$}
\item "{$(i)$}"  $M_i \le_{\frak s} N_i$
\sn
\item "{$(ii)$}"  $N_i$ is $\le_{\frak s}$-increasing continuous
\sn
\item "{$(iii)$}"  {\rm tp}$_{\frak s}
(a_i,N_i,N_{i+1})$ does not fork over $M_i$
\sn
\item "{$(iv)$}"  $(N_i,N_{i+1},a_i) \in K^{3,\text{uq}}_{\frak s}$
\sn
\item "{$(v)$}"  $M^+$ can be $\le_{\frak s}$-embedded into $N_0$ over
$M_0$
\sn
\item "{$(vi)$}"  $N_i$ is $(\lambda,*)$-brimmed over $M_i$ for $i
\le \beta$
\ermn
2) Assume further {\rm NF}$_{\frak s}(M_0,M^+,M_\beta,M^*)$, e.g. $M^+ =
M_0,M_\beta \le_{\frak s} M^*$, \ub{then} we can
replace (v) by
\mr
\item "{$(v)^+$}"  $M^+ \le_{\frak s} N_0$ and $M^* \le_{\frak s}
N_\beta$.
\endroster
\endproclaim
\bigskip

\demo{Proof}  1)  We try to choose by induction on $\zeta < \lambda^+$ a sequence
$\bar M^\zeta = \langle M^\zeta_i:i \le \beta \rangle$ such that
\mr
\item "{$(a)$}"  $\bar M^\zeta$ is $\le_{\frak s}$-increasing continuous
\sn
\item "{$(b)$}"  $\bar M^0 = \langle M_i:i \le \beta \rangle$
\sn
\item "{$(c)$}"  for each $i < \beta$ the sequence $\langle M^\varepsilon_i:
\varepsilon \le \zeta \rangle$ is $\le_{\frak s}$-increasing continuous
\sn
\item "{$(d)$}"  tp$_{\frak s}(a_i,M^\zeta_i,M^\zeta_{i+1})$ 
belongs to ${\Cal S}^{\text{bs}}_{\frak s}(M^\zeta_i)$ and 
does not fork over $M^0_i$
\sn
\item "{$(e)$}"  if $\zeta = 1$ then $M^+$ can be $\le_{\frak s}$-embedded 
into $M^\zeta_0$ over $M_0 = M^0_0$
\sn
\item "{$(f)$}"   if $\zeta = \varepsilon +1$ and $\varepsilon$ limit, \ub{then}
for some $i < \beta$  we have $\neg 
\text{ NF}_{\frak s}(M^\varepsilon_i,
M^\varepsilon_{i+1},M^\zeta_i,M^\zeta_{i+1})$
\sn
\item "{$(g)$}"  if $\zeta = \varepsilon +2$ \ub{then} $M^\zeta_i$ is
$(\lambda,*)$-brimmed over $M^{\varepsilon +1}_i$.
\ermn
There is no problem to define for $\zeta =0,\zeta=1$ and $\zeta$ limit.
For $\zeta = \varepsilon +1,\varepsilon$ not limit straightforward as in
\chaptercite{600}, so assume $\varepsilon$ is limit, 
if we cannot proceed then $\varepsilon$ is a
limit ordinal and we are done: let $\langle N_i:i \le \beta \rangle = \langle
M^\zeta_i:i \le \beta \rangle$.  If we succeed to carry the induction
for all $\zeta < \lambda^+$.
Now we have ${\frak s}$ is successful (we use the second demand (b) of
\sciteu{705-stg.0B}(3)), for each $i < \beta$ for some club $E_i$ of
$\lambda^+$, for every $\varepsilon < \zeta$ from $E$ we have
NF$_{\frak
s}(M^\varepsilon_i,M^\varepsilon_{i+1},M^\zeta_i,M^\zeta_{i+1})$.  Let
$\varepsilon < \zeta$ be successive members of $E = \cap\{E_i:i <
\beta\}$, so by monotonicity of non-forking we have $i < \beta
\Rightarrow \text{ NF}_{\frak
s}(M^\varepsilon_i,N^\varepsilon_{i+1},M^{\varepsilon
+1}_i,M^{\varepsilon +1}_{i+1})$, contradiction to the construction.
\nl
2) Similarly with the following changes: we let $M^0_{\beta +1} =
M^0_\beta$, in the demand above we have $\bar M^\zeta = \langle
M^\zeta_i:i \le \beta + 1 \rangle$ and
\mr
\item "{$(e)'$}"  if $\zeta = 1$ then $M^+ \le_{\frak s} M^\zeta_0,M^*
\le_{\frak s} M^\zeta_{\beta +1}$
\sn
\item "{$(h)$}"  if $\zeta = \varepsilon +2$ and $M^\varepsilon_\beta
\ne M^\varepsilon_{\beta +1}$ then for some $b \in M^\varepsilon_{\beta
+1} \backslash M^\varepsilon_\beta$ the type tp$(b,M^{\varepsilon
+1}_\beta,M^{\varepsilon + 1}_{\beta +1}) \in {\Cal S}_{\frak s}
(M^{\varepsilon +1}_\beta)$ forks $M^\varepsilon_\beta$).
\ermn
Alternatively for every $\zeta \ge 1$, NF$_{\frak
s}(M^0_0,M^\zeta_0,M^0_\beta,M^\zeta_\beta)$ hence NF$_{\frak
s}(M^0_0,M^1_0,M^0_\beta,M^\zeta_\beta)$ over $M^1_0 \cup B^0_\beta$
(using uniqueness of NF).

As $\cup\{M^\zeta_\beta:\zeta < \lambda^+\}$ is saturated we can embed
$M^*$ into this union hence into $M^{\zeta(1)}_\beta$ for some
$\zeta(1)$. In the end of the proof of (1) demand $\zeta > \zeta(1)$.
${{}}$  \hfill$\square_{\scite{705-d.2A}}$\margincite{705-d.2A}
\enddemo
\bigskip

\proclaim{\stag{705-d.2B} Claim}  Assume 
$(*)_{\bar M}$ holds (see below) and
$p \in {\Cal S}^{\text{bs}}_{\frak s}(M_0)$ 
\ub{then} we can find $\bar N$
and $a$ such that $(*)_{\bar N,\bar M}$ holds, where:
\mr
\item "{$(*)_{\bar N,\bar M}$}"  $\ell g(\bar N) 
= \ell g(\bar M),M_i \le_{\frak s} N_i$ are from 
$K_{\frak s},\bar N$ is $<_{\frak s}$-increasing, 
$a \in N_0$, {\rm tp}$(a,M_i,N_i)$ is a nonforking extension of $p$ and
$(M_i,N_i,a) \in K^{3,\text{uq}}_{\frak s}$ 
for every $i < \ell g(\bar M)$
and $N_{i+1}$ is $(\lambda,*)$-brimmed over 
$N_i$ and $N_0$ is $(\lambda,*)$-brimmed
\sn
\item "{$(*)_{\bar M}$}"  $\bar M = \langle M_i:i < \alpha \rangle$ is
$\le_{\frak s}$-increasing continuous, $\alpha \le \lambda^+,M_0$ is
$(\lambda,*)$-brimmed and $M_{i+1}$ is $(\lambda,*)$-brimmed over $M_i$
for $i$ such that $i+1 < \alpha$ (hence $[i < j < \alpha 
\Rightarrow M_j$ is 
$(\lambda,*)$-brimmed over $M_i])$ and if $\alpha$ is limit, 
then $[i < \alpha \Rightarrow \cup\{M_j:j < \alpha\}$ is 
$(\lambda,*)$-brimmed over $M_i])$.
\endroster
\endproclaim
\bigskip

\demo{Proof}  By the proof of \scite{705-d.1} (let $\langle
\alpha_\zeta:\zeta \in [1,\delta]\rangle$ be increasing continuous
with $\alpha_\zeta \in E,\alpha_1 > \text{ Min}(E)$, let $\alpha_0 =
0$, let $N_\zeta = M_{1,\alpha_\zeta}$, for $\zeta > 0$ and lastly
find $N_0 \le_{\frak s} M_{1,\text{Min}(E)}$).  \hfill$\square_{\scite{705-d.2B}}$\margincite{705-d.2B}
\enddemo
\bn
\margintag{705-d.7}\ub{\stag{705-d.7} Conclusion}:  1) $K^{3,\text{uq}}_{\frak s}$ 
is $\le^*_{\text{bs}}$-dense in
$K^{3,\text{bs}}_{\frak s}$, also $<^*_{\text{bs}} 
\subseteq \le_{\text{bs}}$ (both by their definitions). \nl
2) If $(M_1,N_1,a) <^*_{\text{bs}} (M_2,N_2,a) 
\le_{\text{bs}} (M_3,N_3,a)$ \ub{then}
$(M_1,N_1,a) <^*_{\text{bs}} (M_3,N_3,a)$. \nl
3) If $(M_1,N_1,a) \le_{\text{bs}} (M_2,N_2,a) 
<^*_{\text{bs}} (M_3,N_3,a)$ \ub{then}
$(M_1,N_1,a) <^*_{\text{bs}} (M_3,N_3,a)$ and 
$(M_1,N_1,a) \le_{\text{bs}} (M_3,N_3,a)$. \nl
4) If $\delta < \lambda^+$ is a limit ordinal, $(M_i,N_i,a) 
\in K^{3,\text{bs}}_{\frak s}$ for 
$i < \delta$ is $<^*_{\text{bs}}$-increasing \ub{then}
$(M_i,N_i,a) <^*_{\text{bs}} 
(\dbcu_{j < \delta} M_j,\dbcu_{j < \delta} N_j,a) \in 
K^{3,\text{bs}}_\lambda$.
\bigskip

\proclaim{\stag{705-d.8} Claim}  (Prime Existence) 
1) If $M \in K_{{\frak s}(+)},p \in
{\Cal S}^{\text{bs}}_{{\frak s}(+)}(M)$, \ub{then} there are $N \in 
K_{{\frak s}(+)}$ and an element
$a$ satisfying $(M,N,a) \in K^{3,\text{pr}}_{{\frak s}(+)}$.
This means that if $M \le_{\frak K} M' \in K_{{\frak s}(+)},
a' \in N'$ realizes $p$ \ub{then} there is a 
$\le_{\frak K}$-embedding $f$ of $N$
into $M$ 
such that $f \restriction M = { \text{\rm id\/}}_M,f(a) = a'$.  \nl
2) In fact
if $(M,N,a)$ is like $(M_0,M_1,a)$ of \scite{705-d.1} 
then this holds, i.e., $(M,N,a)$ is
canonical ${\frak s}^+$-primes are primes for ${\frak s}^+$. 
\endproclaim
\bigskip

\demo{Proof of \scite{705-d.8}}    Let 
$M_0 = M$ and let $\bar M_0,M_1,\bar M_1,a$ be as in
\scite{705-d.1} and let $N = M_1$ and we shall prove that $N,a$ are as required.
So let $M \le_{\frak K} M' \in K_{{\frak s}(+)},a' \in M'$,
tp$_{{\frak s}(+)}(a',M,M') = p$.  We define by induction on $\alpha$ a 
$\le_{\frak K}$-embedding $f_\alpha$ of $M_{1,\alpha}$ into $M'$, such
that $f_\alpha(a) = a',f_\alpha$ is increasing continuous and $f_\alpha
\restriction M_{0,\alpha} \equiv \text{ id}_{M_0,\alpha}$.  For $\alpha =0$,
as $M'$ is saturated in ${\frak K}^{\frak s}_{\lambda^+}$ (over
$\lambda$)  and $a'$ 
realizes in $M'$ the type $p \restriction
M_{0,0}$ this should be clear. For $\alpha$ limit take the unions.  For
$\alpha = \beta +1$, there is a model $N_\alpha \le_{\frak K} M'$ from
$K_\lambda$ which includes $f_\alpha(M_{1,\beta}) \cup M_{0,\beta +1}$ and
is $(\lambda,*)$-brimmed for ${\frak s}$ over this set, there is such
$N$ as $M'$ is saturated over $\lambda$.  So as $(M_{0,\beta},M_{1,\beta},a)
\in K^{3,\text{uq}}_\lambda$ and tp$_{\frak s}(a,M_{0,\beta +1},M')$ does not
fork over $M_{0,\beta}$, we have
NF$_{\frak s}(M_{0,\beta},f_\beta(M_{1,\beta}),M_{0,\beta +1},N_\alpha)$ hence
by the definition of $K^{3,\text{uq}}_{\frak s}$ 
we can extend $f_\beta \cup \text{ id}_{M_0,\alpha}$
to a $\le_{\frak K}$-embedding 
$f_\alpha$ of $M_{1,\alpha}$ into $N_\alpha$.

So having carried the induction, $f = \cup\{f_\alpha:\alpha < \lambda^+\}$
is a $\le_{\frak K}$-embedding of $M_1 = \dbcu_{\alpha < \lambda^+}
M_{1,\alpha}$ into $M'$ over $M = M_0$ mapping $a$ to $a'$, so we are done.
\hfill$\square_{\scite{705-d.8}}$\margincite{705-d.8}
\enddemo
\bigskip

\proclaim{\stag{705-d.5} Claim}  $(K^{3,\text{uq}}_{\frak s},
\le^{\frak s}_{\text{bs}})$ is $\lambda^+$-strategically closed.
Moreover, if $\delta < \lambda^+$ is 
a limit ordinal, $\langle (N_{0,i},N_{1,i},a):i < \delta \rangle$ is
$\le_{\text{bs}}$-increasing 
continuous in $K^{3,\text{uq}}_\lambda$ and $N_{0,i+1},
N_{1,i+1}$ is $(\lambda,*)$-brimmed over $N_{0,i},N_{1,i}$ respectively
for each $i < \delta$ (equivalently, 
the sequence is $\le^*_{\text{bs}}$-increasing
continuous), \ub{then} $(\dbcu_{i < \delta} N_{0,i} \, 
\dbcu_{i < \delta} N_{1,i},a) \in K^{3,\text{uq}}_\lambda$.
\endproclaim
\bn
Recalling
\definition{\stag{705-d.8K} Definition}  A partial order $I$ is
$\delta$-strategically closed if in the following game the COM player
has a winning strategy.  A play last $\delta$-moves, for $\alpha <
\delta$ the INC player chooses $s_\alpha \in I$ such that $\beta <
\alpha \Rightarrow t_\beta \le_I s_\alpha$ and then the player COM
chooses $t_{\alpha,s}$ such that $s_\alpha \le_I t_\alpha$.  The
player INC wins the play if or some $\alpha < \delta$ he has no legal
move; otherwise, the player COM wins the play.
\enddefinition
\bigskip

\demo{Proof}  It suffices to prove the second sentence
by \scite{705-d.2B} or \scite{705-d.7}. 
Let $\langle M_{0,i},M_{1,i},a):i < \lambda^+ \rangle$ and $E$
be as constructed in \scite{705-d.1} be such that
$i = 0 \Rightarrow (M_{0,i},M_{1,i},a)
= (N_{0,i},N_{1,i},a)$.

We now by induction on $i$ choose $f_i,\alpha_i$ such that:
\mr
\item "{$(a)$}"  $f_i$ is an $\le_{\frak s}$-embedding of $N_{1,i}$ into
$M_{1,\alpha_i}$
\sn
\item "{$(b)$}"  $f_i$ increasing continuous in $i$
\sn
\item "{$(c)$}"  $f_i$ maps $N_{0,i}$ into $M_{0,\alpha_i}$
\sn
\item "{$(d)$}"  $\alpha_i$ is increasing continuous, $i > 0 \Rightarrow
\alpha_i \in E$ 
\sn
\item "{$(e)$}"  $\alpha_0=0$ and $f_0$ is the identity
\sn
\item "{$(f)$}"  $f_i(N_{1,i}) \cap \dbcu_{\gamma < \lambda^+}
M_{0,\gamma} = f_i(N_{0,i})$
\sn
\item "{$(g)$}"  $N_{0,\alpha_i} \le_{\frak K} f_i(M_{0,i+1})$.
\ermn
Note that clause (f) follows automatically as $(N_{0,i},N_{1,i},a) \in
K^{3,\text{uq}}_\lambda$ and non-forking amalgamation is disjoint.  
Also for limit $i$, clearly $f_i(N_{0,i}) = 
M_{0,\alpha_i}$ (check) and $\{a\} \cup M_{0,\alpha_i} \subseteq f_i
(N_{1,i}) \le_{\frak s} M_{1,\alpha_i}$. 

Lastly $f_i$ maps $(\dbcu_{i < \delta} N_{0,i},\dbcu_{i < \delta} \,
N_{1,\delta},a)$ isomorphically into $(M_{0,\alpha_\delta},
M_{1,\alpha_\delta},a)$ and maps $\dbcu_{i < \delta} N_{0,i}$ onto
$M_{0,\alpha_\delta}$ (see clause (g), i.e., the previous paragraph).  
The latter belongs to 
$K^{3,\text{uq}}_\lambda$, so 
(by \scite{705-c.2}(2) monotonicity for $K^{3,uq}_{\frak s}$)
we are done.  \hfill$\square_{\scite{705-d.5}}$\margincite{705-d.5}
\enddemo
\bn
In \scite{705-zm.2}(2) below we show how relevant situations in ${\frak
s}^+$ reflect to ${\frak s}$.
\proclaim{\stag{705-zm.2} Claim}  1) If $(M_0,M_1,a) \in 
K^{3,\text{pr}}_{{\frak s}(+)}$ and $\bar M^\ell = \langle
M^\ell_\alpha:\alpha < \lambda \rangle$ is a 
$\le_{\frak s}$-representation of $M_\ell$ for
$\ell = 1,2$ \ub{then} for some club $E$ of $\lambda^+$ we have
\mr
\item "{$(*)$}"  if $\alpha \in E$ then $(M^0_\alpha,M^1_\alpha,a) \in
K^{3,\text{uq}}_{\frak s}$.
\ermn
2) Assume $M_0 \le_{{\frak s}(+)} 
M_\ell \le_{{\frak s}(+)} M_3$ and $a_\ell
\in M_\ell$ and $p_\ell = { \text{\rm tp\/}}_{{\frak s}(+)}
(a_\ell,M_0,M_\ell)$ for $\ell = 1,2$.  \ub{Then} for a club 
$E$ of $\lambda^+$ for every $\delta \in E$ we have
\mr
\widestnumber\item{$(iii)$}
\item "{$(i)$}"  if $p_\ell \in {\Cal S}^{\text{bs}}_{{\frak s}(+)}
(M_0)$ then $p_{\ell,\delta} = { \text{\rm tp\/}}_{\frak s}
(a_\ell,M_{0,\delta},M_{\ell,\delta}) \in {\Cal S}^{\text{bs}}
(M_{0,\delta})$ and $M_{\ell,\delta}$ (and $p_{\ell,\delta}$)
are witnesses for $p_\ell$
\sn
\item "{$(ii)$}"  if $\ell \in \{1,2\}$ and
{\rm tp}$_{{\frak s}(+)}(a_\ell,M_{3 - \ell},M_3)$ is an
${\frak s}^+$-nonforking extension of $p_\ell$ \ub{then} 
{\rm tp}$_{\frak s}(a_\ell,M_{3 -\ell,\delta},
M_{3,\delta})$ is an ${\frak s}$-nonforking extension
of {\rm tp}$_{\frak s}(a_\ell,M_{0,\delta},M_{\ell,\delta})$ (hence of
{\rm tp}$_{\frak s}(a_\ell,M_{0,\text{min}(E)},M_{\ell,\text{min}(E)})$
\sn
\item "{$(iii)$}"  if $\alpha \in \delta \cap E$ then
$M^\ell_{\alpha +1}$ is $(\lambda,*)$-brimmed over $M^\ell_\alpha$
for $\ell < 4$
\sn
\item "{$(iv)$}"  if $(M_0,M_\ell,a_\ell) 
\in K^{3,\text{uq}}_{{\frak s}+1}$ then $(M^0_\delta,
M^\ell_\delta,a_\ell) \in K^{3,uq}_{\frak s}$.
\endroster
\endproclaim
\bigskip

\demo{Proof}  1) We can find $M_2,a_2$ such that $(M_0,M_2,a_2) \in 
K^{3,\text{pr}}_{{\frak s}(+)}$ and this triple is canonically ${\frak
s}^+$-prime, i.e., is as in \scite{705-d.1}
(with $M_0,M_2,a_2$ here standing for $M_0,M_1,a$ there, of course) and
tp$_{{\frak s}(+)}(a_2,M_0,M_2) = \text{ tp}_{{\frak s}(+)}(a,M_0,M_1)$.  
There is a $\le_{{\frak s}(+)}$-embedding of $M_1$ into $M_2$ over $M_0$ 
mapping $a$ to $a_2$ so \wilog \, $a_2 = a \and M_1 \le_{\frak s} M_2$.  
Let $\bar M^\ell$ be a $\le_{\frak s}$-representation of $M_\ell$ for 
$\ell < 3$ and $E$ a thin enough club of $\lambda^+$.  As 
$M_\ell \le_{{\frak s}(+)} M_{\ell +1}$ for any 
$\alpha < \beta$ from $E$ we have NF$_{\frak s}
(M^\ell_\alpha,M^{\ell +1}_\alpha,M^\ell_\beta,M^{\ell +1}_\beta)$ 
for $\ell = 0,1$ and by the choice of $M_2$ for any $\alpha < \beta$ from
$E$ we have $(M^0_\alpha,M^2_\alpha,a) \in 
K^{3,\text{uq}}_{\frak s}$.  By monotonicity of 
$K^{3,\text{uq}}_{\frak s}$, i.e., \scite{705-c.2}(2) 
as $M^0_\alpha \cup \{a\} \subseteq M^1_\alpha
\le_{{\frak K}[{\frak s}]} M^2_\alpha$ we get $(M^0_\alpha,M^1_\alpha,a)
\in K^{3,\text{uq}}_{\frak s}$ as required. \nl
2) Straightforward (for (iii) recall that $M_\ell \in K^{\frak
s}_{\lambda^+}$ is saturated (above $\lambda$) for ${\frak s}$).   \nl
${{}}$  \hfill$\square_{\scite{705-zm.2}}$\margincite{705-zm.2}
\enddemo
\bigskip

\proclaim{\stag{705-d.9} Claim}  1) $(M,N,a) \in 
K^{3,\text{pr}}_{{\frak s}(+)}$,
\ub{then} it is canonically ${\frak s}^+$-prime. \nl
2) Uniqueness: if $(M,N_\ell,a_\ell) \in 
K^{3,\text{pr}}_{{\frak s}(+)}$ and
{\rm tp}$_{{\frak s}(+)}(a_1,M,N_1) = 
{ \text{\rm tp\/}}_{{\frak s}(+)}(a_1,M,N_2)$
\ub{then} there is an isomorphism $f$
from $N_1$ onto $N_2$ over $M$ satisfying $f(a_1)= a_2$.
\endproclaim
\bigskip

\demo{Proof}  1) By \scite{705-c.2}(2) and \scite{705-zm.2}. \nl
2) By part (1), we know that $(M,N_\ell,a_\ell)$ is canonically
${\frak s}^+$-prime.  Now we build the isomorphism by 
hence and forth as in the proof of \scite{705-d.8}.
\hfill$\square_{\scite{705-d.0}}$\margincite{705-d.0}
\enddemo
\bn
It is good to know that also NF$_{{\frak s}(+)}$ reflect down (when we
have it).
\proclaim{\stag{705-d.4} Claim}  Assume that also ${\frak s}^+$ is weakly
successful so {\rm NF}$_{{\frak s}(+)}$ is well defined.  
If $M_\ell \in K_{{\frak s}(+)}$ for $\ell = 0,1,2,3$ and
{\rm NF}$_{{\frak s}(+)}(M_0,M_1,M_2,M_3)$ and $\bar M_\ell = \langle
M_{\ell,\alpha}:\alpha < \lambda^+ \rangle $ does $\le_{\frak
s}$-represent $M_\ell$ for $\ell < 4$, \ub{then} for a club of $\delta
< \lambda^+$ we have 
{\rm NF}$_{\frak s}
(M_{0,\delta},M_{1,\delta},M_{2,\delta},M_{3,\delta})$.
\endproclaim
\bigskip

\demo{Proof}  Without loss of generality $M_\ell$ is
$(\lambda^+,*)$-brimmed over $M_0$ for ${\frak s}^+$ for $\ell =
1,2$ and $M_3$ is $(\lambda^+,*)$-brimmed over $M_1 \cup M_2$ (by
density of $(\lambda^+,*)$-brimmed recalling ${\frak s}^+$ is a
weakly successful $\lambda^+$-good frame, and the existence property
of NF$_{{\frak s}(+)}$ and the monotonicity of NF$_{{\frak s}(+)}$ and
NF$_{\frak s}$).

Let $\langle N_{\ell,\alpha},
a^\ell_\alpha:\alpha < \lambda^+ \rangle$ be as
in \scite{705-c.6} applied to ${\frak s}(+)$ 
for $(M_0,M_\ell)$ for $\ell = 1,2$ (the length being
$\lambda^+$ is somewhat more transparent and is allowed as $M_\ell$ is
$(\lambda^+,*)$-brimmed over $M_0$ for ${\frak s}^+$ by
\scite{705-c.6}(6)).  As $M_3$ is
$(\lambda^+,*)$-brimmed over $M_1 \cup M_2$, \wilog \, we have
$\langle N_{3,\alpha}:\alpha \le \lambda^+ \rangle$ which
is $\le_{{\frak s}(+)}$-increasing continuous, $N_{3,0} = M_2,N_{3,\lambda^+}
= M_3,N_{3,\alpha} \cap M_1 = N_{1,\alpha}$ and
$(N_{1,\alpha},N_{1,\alpha +1},a^1_\alpha) \le^{{\frak s}(+)}_{\text{bs}}
(N_{3,\alpha},N_{3,\alpha +1},a^1_\alpha)$; see \sectioncite[\S6]{600}.
For each $\alpha < \lambda^+,\ell \le 2,3$ let $\langle
N_{\ell,\alpha,i}:i < \lambda^+ \rangle \le_{\frak s}$-represent
$N_{\ell,\alpha}$. Let $E_{1,\alpha}$ be a club of $\alpha$ such that
$i \in E_{1,\alpha}$ implies $(N_{1,\alpha,i},N_{1,\alpha
+1,i},a^1_\alpha) \in K^{3,\text{uq}}_\lambda$ and
tp$_{\frak s}(a_\alpha,N_{1,\alpha,i},N_{1,\alpha +1,i})$ does not fork over
$N_{1,\text{Min}(E_{1,\alpha})}$ hence
tp$_{\frak s}(a^1_\alpha,N_{3,\alpha,i},N_{3,\alpha +1,i})$ does not fork over
$N_{1,\text{Min}(E_{1,\alpha +1})}$; note that $E_{1,\alpha}$ exists by
\scite{705-zm.2}.  Let $E = \{\delta < \lambda^+:\delta 
\text{ limit},\delta \in \dbca_{\alpha < \delta} E_{1,\alpha}$ and
$M_{\ell,\delta} \cap M_m = M_{\ell,\delta}$ for $\ell < m \le
3,(\ell,m) \ne (1,2),M_{1,\delta} = N_{1,\delta},M_{3,\delta} =
\dbcu_{\alpha < \delta} N_{3,\alpha,\delta} = \dbcu_{\alpha < \delta}
N_{3,\alpha,\alpha},M_{0,\delta} = N_{1,0,\delta}$ and $M_{2,\delta} =
N_{3,0,\delta}\}$.  Now check.  \hfill$\square_{\scite{705-d.4}}$\margincite{705-d.4}
\enddemo
\bn
\centerline {$* \qquad * \qquad *$}
\bn
\ub{Discussion}:  By the 
above in the cases we construct good $\lambda$-frames 
${\frak s}$ in \marginbf{!!}{\cprefix{600}.\scite{600-Ex.4}}, ${\frak s}$ is essentially ${\frak
t}^+,{\frak t}$ is almost good and we need that ${\frak s}$ is
successful to get ${\frak s}^+$ has primes; but ${\frak t}$ is close
enough to being good and successful so that ${\frak s}$ itself has
primes.  In the other two cases of \scite{705-stg.2}(1) there are primes for
different reason: $\aleph_0$ is easier.  As for the case we use
\cite{Sh:576} this is not clear.
\bigskip

\proclaim{\stag{705-d.13} Claim}  1) If ${\frak s}$ is as in
\marginbf{!!}{\cprefix{600}.\scite{600-Ex.4}}, \ub{then} ${\frak s}$ has primes. \nl
2) If ${\frak s}$ is as in \marginbf{!!}{\cprefix{600}.\scite{600-Ex.1}},\marginbf{!!}{\cprefix{600}.\scite{600-Ex.1a}}, i.e. Cases 2,3 of
\scite{705-stg.2}(1), \ub{then} ${\frak s}$ has primes.
\endproclaim
\bigskip

\demo{Proof}  1) Like the proof of \scite{705-d.8}. \nl
2) Using stability in $\aleph_0$ we can construct primes directly.
\hfill$\square_{\scite{705-d.13}}$\margincite{705-d.13}
\enddemo
\newpage

\head {\S5 Independence} \endhead  \resetall   \sectno=5
 \spuriousreset 
\bigskip

Here we make a real step forward: independence (of set of elements
realizing basic types) can be defined and proved to be as required.
In an earlier version we have used existence of primes but eventually
eliminate it.
\bigskip

\demo{\stag{705-5.0} Hypothesis}   ${\frak s}$ is a good weakly
successful $\lambda$-frame.
\enddemo
\bigskip

\remark{Remark}  E.g. ${\frak s} = {\frak t}^+,{\frak t}$ 
is $\lambda^+$-good and 1.5-successful, $\lambda = \lambda_{\frak s}$.
\endremark
\bigskip

\definition{\stag{705-5.1} Definition}  Let $M \le_{\frak s} N$ 
(hence from $K_{\frak s} = K^{\frak s}_\lambda$).
\nl
1) Let $\bold I_{M,N} = \{a \in N:\text{tp}_{\frak s}(a,M,N) 
\in {\Cal S}^{\text{bs}}_{\frak s}(M)\}$. \nl
2) We say that $\bold J$ is 
independent in $(M,A,N)$ if $(*)$ below holds;
when: $A = N'$ we may write $N'$ instead of $A$; 
if $N$ is understood from the context we may
write ``over $(M,A)$"; if $A=M$ and we may omit it and 
then we say ``in $(M,N)$" or
``for $(M,N)$"; where:
\mr
\item "{$(*)$}"  $\bold J \subseteq \bold I_{M,N},M \le_{\frak s} N,M
\subseteq A \subseteq N$ and we can find a witness
$\langle M_i,a_j:i \le \alpha,j < \alpha \rangle$ 
and $N^+$ which means:
{\roster
\itemitem{ $(a)$ }  $\langle M_i:i \le \alpha \rangle$ is
$\le_{\frak K}$-increasing continuous \footnote{note that omitting
``continuous" makes no difference}
\sn
\itemitem{ $(b)$ }  $M \cup A \subseteq M_i \le_{\frak s} N^+$, (usually
$M \subseteq A$) and $N \le_{\frak s} N^+$
\sn
\itemitem{ $(c)$ }  $a_i \in M_{i+1} \backslash M_i$, (if we forget to
mention $M_\alpha$ we may stipulate $M_\alpha = N^+$)
\sn
\itemitem{ $(d)$ }  tp$(a_i,M_i,M_{i+1})$ does not fork over $M$, 
\sn
\itemitem{ $(e)$ }  $\bold J = \{a_i:i < \alpha\}$.
\endroster}
\endroster
\enddefinition
\bn
The name independent indicates we expect various properties, like
finite character, so we start to prove them.
\proclaim{\stag{705-5.1R} Claim}  Assume that {\rm NF}$_{\frak s}
(M_0,M_1,M_2,M_3)$.  \ub{Then} $\bold J$ is independent in $(M_0,M_2)$
iff $\bold J$ is independent in $(M_0,M_1,M_3)$.
\endproclaim
\bigskip

\demo{Proof}  The if implication is trivial.  The only if implication
is easy by chasing arrows and using NF$_{\frak s}$-uniqueness (and
existence and \scite{705-stg.32}).
\hfill$\square_{\scite{705-5.12}}$\margincite{705-5.12}
\enddemo
\bigskip

\proclaim{\stag{705-5.2} Theorem}  Assume $M \le_{\frak s} N$ so are from 
$K_{\frak s}$ and $\bold J \subseteq \bold I_{M,N}$. \nl
1) The following are equivalent:
\mr
\item "{$(*)_0$}"  every finite $\bold J' \subseteq \bold J$ is independent
in $(M,N)$
\sn
\item "{$(*)_1$}"  $\bold J$ is independent in $(M,N)$
\sn
\item "{$(*)_2$}"  like $(*)$ of Definition \scite{705-5.2} adding
{\roster
\itemitem{ $(f)$ }  $M_{i+1}$ is $(\lambda,*)$-brimmed over $M_i
\cup \{a_i\}$ for $i < \alpha$
\endroster}
\item "{$(*)_3$}"  for every ordinal $\beta,|\beta| 
= |\bold J|$ and list
$\langle a_i:i < \beta \rangle$ with no repetitions of $\bold J$ 
\ub{there are} $M_i$ (for $i \le \beta$) such that $\langle M_j,a_i:
j \le \beta,i < \beta \rangle$ and $N^+$ which satisfy:
{\roster
\itemitem{ $(a)$ }  $M_i$ is 
$\le_{\frak s}$-increasing continuous, $M_0=M$
\sn
\itemitem{ $(b)$ }  $M \le_{\frak s} M_i \le_{\frak s} N^+$ and $N
\le_{\frak s} N^+$
\sn
\itemitem{ $(c)$ }  $a_i \in M_{i+1} \backslash M_i$
\sn
\itemitem{ $(d)$ }  {\rm tp}$_{\frak s}
(a_i,M_i,M_{i+1})$ does not fork over $M$
\sn
\itemitem{ $(e)$ }  $\bold J = \{a_i:i < \beta\}$
\sn
\itemitem{ $(f)_3$ }  $M_{i+1}$ is $(\lambda,*)$-brimmed
\footnote{omitting $\{a_i\}$ give an equivalent condition}
over $M_i \cup \{a_i\}$
\endroster}
\item "{$(*)_4$}"  like $(*)_3$ replacing $(f)_3$ by
{\roster
\itemitem{ $(f)_4$ }  $(M_i,M_{i+1},a_i) \in 
K^{3,\text{uq}}_{\frak s}$ 
\endroster}
\item "{$(*)_5$}"  like $(*)_4$ adding
{\roster 
\itemitem{ $(g)$ }  if there is $M'_{i+1} \le_{\frak s} M_{i+1}$ such
that $(M_i,M'_{i+1},a_i) \in K^{3,\text{pr}}_{\frak s}$ then
$(M_i,M_{i+1},a_i) \in K^{3,\text{pr}}_{\frak s}$
\endroster}
\item "{$(*)_6$}"  like $(*)_5$ adding
{\roster
\itemitem{ $(h)$ }  if $(M_i,M_{i+1},a_i) \in 
K^{3,\text{pr}}_{\frak s}$ for each $i < \beta$ (holds e.g. 
if ${\frak s}$ has primes) then $M_\alpha \le_{\frak s} N$.
\endroster}
\ermn
(so in $(*)_6$ 
$\langle (M_i,a_i:i < \beta \rangle$ is a witness for ``$\bold J$
independent for $(M,N)$ and it is an $M$-based uq-decomposition inside
$(M,N)$ see Definition \scite{705-c.1A}).
\nl
2) If $M \cup \bold J \subseteq N^- \le_{\frak s} N$ and 
$N \le_{\frak s} N^+ \in K_{\frak s}$, \ub{then}: $\bold J$ is independent 
in $(M,N^+)$ \ub{iff} $\bold J$ is
independent in $(M,N)$ \ub{iff} $\bold J$ is independent in $(M,N^-)$. \nl
3) If $\bold J$ is independent in $(M,N)$ and
$a_i \in \bold J$ for $i < \beta$ are with no repetitions, $M_0 =M$
and $\langle M_i:i \le \beta \rangle,\langle a_i:i < \beta \rangle$ are
as in $(*)_4$ clauses (a)-(e),(f)$_4$ from 
part (1) and $M_\beta \le_{\frak s} N$, 
\ub{then} $\bold J \backslash \{a_i:i < \beta\}$
is independent in $(M,M_\beta,N)$. \nl
4) If $M^- \le_{\frak s} M,M^- \in K_{\frak s},\bold J$ is independent in
$(M,N)$ and $\bold J' \subseteq \bold J$ and $[a \in \bold J'
\Rightarrow { \text{\rm tp\/}}_{\frak s}(a,M,N)$ does not fork 
over $M^-]$, \ub{then} 
$\bold J'$ is independent in $(M^-,N)$; moreover if $M^- \subseteq A
\subseteq M$ then $\bold J'$ is independent in $(M^-,A,N)$.
\endproclaim
\bigskip

\demo{Proof}  First note that part (2) is immediate by the
amalgamation property, and also part (4) is straightforward so it is
enough to prove part (1) + (3).  Clearly
\mr
\item "{$\boxtimes_0$}"  $(*)_4 \Rightarrow (*)_3$ \nl
[Why?  We choose $(f_i,M'_i)$ by induction on $i \le \alpha$ such that
$M'_i$ is $\le_{\frak s}$-increasing continuous, $f_i$ is a
$\le_{\frak s}$-embedding of $M_i$ into $M'_i,f_i$ is increasing
continuous, $f_0 = \text{ id}_{M_0}$, tp$(f_{i+1}(a_i),M'_i,M'_{i+1})$
does not fork over $M,M'_{i+1}$ is $(\lambda,*)$-brimmed over $M'_i
\cup \{f_{i+1}(a_i)\}$.  By amalgamation without loss of generality 
there are $(g,N')$
such that $M'_\alpha \le N',g$ is a $\le_{\frak s}$-embedding of $N^+$
into $N'$ extending $f_\alpha$.  Renaming $g = \text{ id}_{N^+}$ so
clearly we are done.]
\sn
\item "{$\boxtimes_1$}"   $(*)_6 \Rightarrow (*)_5 \Rightarrow
(*)_4 \Rightarrow (*)_3 \Rightarrow (*)_2
\Rightarrow (*)_1 \Rightarrow (*)_0$ \nl
[Why?  The implication $(*)_4 \Rightarrow (*)_3$ by $\boxtimes_0$,
$(*)_5 \Rightarrow (*)_4$ as $K^{3,\text{pr}} \subseteq 
K^{3,\text{uq}}$ by \scite{705-c.4}(2) (as we assume \scite{705-5.0}, i.e.,
${\frak s}$ is weakly successful).  For the
others, just read them.]
\sn
\item "{$\boxtimes_2$}"  if $\langle M'_i:i \le \alpha'
\rangle,\langle a_i:i < \alpha' \rangle$ and
$N^+$ are as in $(*)_2$  witnessing
``$\bold J$ is independent in $(M,N)$" then we can find $M''_i
\le_{\frak s} M'_i$ for
$i \le \alpha'$ such that $\langle M''_i:i \le \alpha' \rangle,N^+$ 
is as required in $(*)_4$ clauses (a)-(e),(f)$_4$ of part (1) \nl
[why?  choose $M''_i \le_{\frak s} M'_i$ by induction on $i$ such that
$M'_i$ is $\le_{\frak s}$-increasing continuous, 
$(M''_i,M''_{i+1},a_i) \in K^{3,\text{uq}}_{\frak s}$ and 
$i \le \alpha \Rightarrow M''_i 
\le_{\frak K} M'_i$, using the hypothesis 
``${\frak s}$ is weakly successful" and ``$M'_{i+1}$ is
$(\lambda,*)$-brimmed over $M'_i$"].
\ermn
Hence
\mr
\item "{$\boxtimes_3$}"  $(*)_3 \Rightarrow (*)_4$ \nl
and similarly (recalling $K^{3,\text{pr}}_{\frak s} \subseteq
K^{3,\text{uq}}_{\frak s}$)
\sn
\item "{$\boxtimes_4$}"  $(*)_4 \Rightarrow (*)_5$.
\ermn
Also if $\langle M_i:i \le \alpha \rangle,N^+,\langle a_i:i < \alpha 
\rangle$ are as in $(*)_1$, i.e.,
satisfy clauses (a)-(e) of $(*)_3$ 
with $\alpha$ instead of $\beta$ then we
can choose $(M^+_i,f_i)$ by induction on $i \le \alpha$ 
such that $M^+_i$ is
$\le_{\frak s}$-increasing continuous, $f_i$ is a $\le_{\frak
s}$-embedding of $M_i$ into $M^+_i,f_i$ is increasing continuous,
NF$(f_i(M_i),M^+_i,f_{i+1}(M_i),M^+_{i+1})$ and $M^+_i$ is
$(\lambda,*)$-brimmed over $M^+_i \cup \{a_i\}$.  Without loss of
generality $f_i = \text{ id}_{M_i}$ for $i \le \alpha$.  By renaming
\wilog \, $f_i = \text{ id}_{M_i}$ and by
amalgamation \wilog \, $M^+_i \le_{\frak s} N^+$.  So
\mr
\item "{$\boxtimes_5$}"  $(*)_1 \Rightarrow (*)_2$.
\ermn
We prove part (1) + (3) by induction on $|\bold J|$.
\enddemo
\bn
\ub{Case 1}:  $|\bold J| \le 1$.

Trivial.
\bn
\ub{Case 2}:  $n = |\bold J|$ finite $> 1$.

As $\bold J$ is finite (and monotonicity of independence) clearly
\mr
\item "{$\otimes_1$}"  $(*)_0 \Leftrightarrow (*)_1$.
\ermn
We first show
\mr
\item "{$\otimes_2$}"  $(*)_2 \Rightarrow (*)_3$.
\ermn
As the permutations exchanging $m,m+1$ generate all permutations of
$\{0,\dotsc,n-1\}$, above it is enough to show
\mr
\item "{$\otimes_3$}"  if $\langle M_k,a_\ell:k \le n,\ell < n
\rangle$ and $N^+$ are as in $(*)_2$ and $m < n-1$ and $a'_\ell$ is $a_\ell$
if $\ell < n \and \ell \ne m \and \ell \ne m+1$, is $a_{m+1}$ if $\ell = m$
and is $a_m$ if $\ell = m+1$, \ub{then} for some $M'_\ell$ for $\ell <
n$ we have $\langle M'_k,a'_\ell:k \le n,\ell < n \rangle$ and $N^+$
are as in $(*)_3$.
\ermn
Why does $\otimes_3$ hold?   
Let $M'_\ell$ be $M_\ell$ if $\ell \le m \vee
\ell \ge m+2$.  
\nl
As ${\frak s}$ is a good frame and
tp$(a_{m+1},M_{m+1},M_{m+2}) \in {\Cal S}^{\text{bs}}_{\frak
s}(M_{m+1})$ does not fork over $M$ and $M \le_{\frak s} M_m
\le_{\frak s} M_{m+1}$ clearly tp$_{\frak s}
(a_{m+1},M_{m+1},M_{m+2}) \in {\Cal S}^{\text{bs}}_{\frak s}
(M_{m+1})$ does not fork over $M_m$ and similarly tp$(a_m,M_m,M_{m+1})
\in {\Cal S}^{\text{bs}}(M_m)$.  Hence there are
by symmetry $M',M''$ such that $M_{m+2} \le_{\frak s} M'',M_m
\le_{\frak s} M' \le_{\frak s} M'',a_{m+1} \in M'$ 
and tp$_{\frak s}(a_m,M',M'')$ does not fork over $M_m$ and \wilog \,
$M'$ is $(\lambda,*)$-saturated?? over $M_m \cup \{a_{m+1}\}$.
As $M_{m+2}$ is $(\lambda,*)$-brimmed 
over $M_{m+1} \cup \{a_{m+1}\}$ which include
$M_m \cup \{a_m\}$, \wilog \, $M'' \le_{\frak s} M_{m+2}$ and, moreover,
\wilog \, $M_{m+2}$ is $(\lambda,*$)-brimmed over $M''$ hence over 
$M' \cup \{a_m\} = M' \cup \{a'_{m+1}\}$
(equivalently over $M'$); just think on the definitions.  Let
$M'_{m+1} = M'$ so $\otimes_1$ holds. \nl
So (in the present case) we have
\mr
\item "{$\otimes_4$}"  $(*)_0 \Leftrightarrow (*)_1 \Leftrightarrow (*)_2
\Leftrightarrow (*)_3 \Leftrightarrow (*)_4 \Leftrightarrow (*)_5$.
\ermn
Next
\mr
\item "{$\otimes_5$}"  part (3) holds.
\ermn  
Why?  By the induction hypothesis, it is enough to deal with the
case $\beta = 1$, i.e., to prove that $\bold J
\backslash \{a_0\}$ is independent in $(M_1,N)$ assuming $(M,M_1,a_0)
\in K^{3,\text{uq}}_{\frak s},a_0 \in \bold J$, also without loss of
generality $\bold J \backslash \{a_0\} \ne \emptyset$ (otherwise 
the conclusion is trivial) hence $n \ge 2$.  Choose $b_0,\dotsc,
b_{n-2}$ such that they list $\bold J \backslash \{a_0\}$ and 
we can let $b_{n-1} = a_0$ and possibly increasing $N$ 
let $\langle M'_\ell:\ell \le n \rangle$ be such that
$\langle b_\ell:\ell < n \rangle,\langle M'_\ell:\ell \le n \rangle$ are as
in $(*)_4$ clauses (a)-(e),(f)$_4$, they 
exist as we have already proved most of part (1) in $\otimes_4$ 
above in the present 
case.  So tp$_{\frak s}(a_0,M'_{n-1},M'_n) = 
\text{ tp}_{\frak s}(b_{n-1},M'_{n-1},N)$ does not 
fork over $M'_0 = M$ so as
$(M,M_1,a_0) \in K^{3,\text{uq}}_\lambda$ we can deduce
that NF$_{\frak s}(M,M'_{n-1},M_1,N)$ by \cite[\S6]{600},\marginbf{!!}{\cprefix{600}.\scite{600-nf.20t}}.  
Hence (or see \scite{705-5.3}(2)) easily by the NF calculus
for some $N^+,N \le_{\frak K} N^+ \in K_{\frak s}$ 
we can find $M_2,\dotsc,M_n$ such that 
$M_1 \le_{\frak s} M_2
\le_{\frak s} \ldots \le_{\frak s} M_n 
\le_{\frak K} N^+,\ell \in \{1,2,\dotsc,n-1\} \Rightarrow 
\text{  NF}_{\frak s}(M'_\ell,
M_{\ell+1},M'_{\ell +1},M_{\ell +2})$.  By \scite{705-stg.32} tp$_{\frak
s}(b_\ell,M_{\ell +1},M_{\ell +2})$ does not fork over $M'_\ell$
hence by transitivity of nonforking for 
$\ell < n-1$ the type tp$_{\frak s}(b_\ell,M_{\ell
+1},M_{\ell+2})$ does not fork over $M$.
So $\langle M_{1 + \ell}:\ell \le n-2 \rangle$ 
witness that $\langle b_0,\dotsc,b_{n-2} \rangle$ is independent
in $(M_1,N^+)$.  So by part (2), i.e., 
for $n-1$, clearly $\langle b_0,\dotsc,
b_{n-2} \rangle$ is independent in $(M_1,N)$, so by part (4) 
we have shown part (3), i.e., $\otimes_3$.
\nl
To complete the proof in the present case we need
\mr
\item "{$\otimes_6$}"  $(*)_5 \Rightarrow (*)_6$.
\ermn
We do more: we prove this in the general case provided that part (3) has been
proved. \nl
So let $\langle a_i:i < \beta \rangle$ list $\bold J$ with no
repetitions and let $N^+,\langle M_i:i \le \beta \rangle$
be as in $(*)_5$.  The only nontrivial case is when $i < \beta
\Rightarrow (M_i,M_{i+1},a_\beta) \in 
K^{3,\text{pr}}_{\frak s}$.  We now
choose by induction on $i \le \beta$ a $\le_{\frak s}$-embedding $f_i$ of
$M_i$ into $N$ such that
$f_0 = \text{ id}_{M_0}$ and $f_{i+1}(a_i) = a_i$.
\nl
Now $f_0$ is defined, if $f_i$ is defined, then by part (3)
which we have already proved for this case we know that 
$\bold J \backslash \{a_j:j < i\}$ is independent in $(M,f(M_i),N^+)$
so tp$_{\frak s}(a_i,f(M_i),N^+)$ does not fork over $M$
hence $f_i(\text{tp}(a_i,M_i,M_{i+1})) = 
\text{ tp}_{\frak s}(a_i,f(M_i),N^+)$.  Hence
$f_{i +1}$ exists as $(M_i,M_{i +1},a_\ell) \in
K^{3,\text{pr}}_{\frak s}$ and the definition of 
$K^{3,\text{pr}}_{\frak s}$.
Lastly $\langle f_\ell(M_\ell):\ell \le n \rangle$ witnesses
that $(*)_6$ holds.
\bn 
\ub{Case 3}:  $|\bold J| = \mu \ge \aleph_0$.

Now we first prove what we now call (3)$^-$, a weaker variant of part (3), 
which is: replacing in the conclusion ``independent" by
``every finite subset is independent", this will be subsequently used to
prove the other parts, so by (1) we shall get 
part (3) itself;  we prove (3)$^-$
by induction on the ordinal $\beta$ (for all
possibilities) and for a fixed $\beta$ by induction on $|\bold J
\backslash \{a_i:i < \beta\}|$ which \wilog \, is finite.

First, for $\beta = 0$ it is trivial. \nl
Second, assume $\beta = \gamma +1$
and let $b_0,\dotsc,b_{n-1} \in \bold J \backslash \{a_i:i < \beta\}$ be
pairwise distinct, and we should prove that $\{b_0,\dotsc,b_{n-1}\}$ is
independent in $(M,M_\beta,N)$.  
Now by the induction hypothesis on $\beta$ applied
to $\langle M_j,a_i:j \le \gamma,i < \gamma \rangle$ and $\{b_0,\dotsc,
b_{n-1},a_\gamma\}$ we deduce that $\{b_0,\dotsc,b_{n-1},a_\gamma\}$ is
independent in $(M,M_\gamma,N)$.  Now the desired conclusion follows 
from the case with $\bold J$ finite.

Lastly, assume $\beta$ is a limit ordinal and let $n < \omega,b_0,\dotsc,
b_{n-1} \in \bold J' =: 
\bold J \backslash \{a_i:i < \beta\}$ be pairwise distinct.
We should prove that 
$\{b_0,\dotsc,b_{n-1}\}$ is independent in $(M,M_\beta,N)$.
By the induction hypothesis on $\beta$ we have $\varepsilon < \beta$ implies
tp$(b_0,M_\varepsilon,N)$ does not fork
over $M$, hence tp$(b_0,M_\beta,N)$ does not fork over $M$.  By Claim
\scite{705-5.3}(2) below we can find $\langle M'_\varepsilon:\varepsilon \le 
\beta \rangle,\le_{\frak s}$-increasing continuous and $N^+$ satisfying 
$N \le_{\frak s} N^+ \in K_{\frak s}$ and NF$_{\frak s}
(M_\varepsilon,M_\zeta,M'_{1+\varepsilon},M'_{1+\zeta})$ 
for $\varepsilon
< \zeta \le \beta,M'_\beta \le_{\frak s} N^+,M_0 \le_{\frak s} M'_0$ 
and tp$_{\frak s}(a_\varepsilon,M'_{1 + \varepsilon},
M'_{1 + \varepsilon +1})$ does not fork over $M$  
such that letting $a'_0 = b'_0,
a'_{1 + \varepsilon} = a_\varepsilon,\bar M' = \langle M'_\varepsilon,
a'_\zeta:\varepsilon \le 1 + \beta,\zeta < \beta \rangle$
we have $(M'_\varepsilon,M'_{\varepsilon +1},a'_\varepsilon) \in
K^{3,\text{uq}}_{\frak s}$ and $M_\varepsilon \le_{\frak s} 
M'_{1 + \varepsilon}$ for $\varepsilon < \beta$; 
note that $1 + \beta = \beta$ as $\beta$ is a limit ordinal.  \nl
Now
\mr
\item "{$\otimes_7$}"  $\{M'_\varepsilon,a_\zeta:\varepsilon \le
\beta,\zeta < \beta \rangle$ is as in $(*)_4$ for $(M'_0,N^+)$.
\ermn
[Why?  As NF$_s(M_0,M'_0,N,N^+)$ by \scite{705-5.1R}.] \nl
Hence by the induction hypothesis on $n,\{b_1,\dotsc,b_{n-1}\}$ 
is independent in $(M'_0,M'_\beta,N^+)$ so by part (4) also in
$(M,M'_{1 + \beta},N^+)$ and recall tp$(b_0,M_\beta,N^+)$
does not fork over $M$ while $b_0 \in M'_\beta,M \le_{\frak s} M_\beta
\le_{\frak s} M'_{1 + \beta}$, so 
easily $\{b_0,\dotsc,b_{n-1}\}$ is independent in
$(M,M_\beta,N^+)$, i.e. by \scite{705-5.3A}(1) below hence by the induction 
hypothesis (using parts (2) + (4) for
the cases of finite $\bold J$) the set
$\{b_0,\dotsc,b_{n-1}\}$ is independent in $(M,M_\beta,N)$ as required
so we have proved (3)$^-$
\mn
Next we prove $(*)_0 \Rightarrow (*)_4$ in part (1).

For proving $(*)_4$ let $\langle a_i:i < \beta \rangle$
be a given list of $\bold J$ and we will find $N^+,\langle M_i:i \le
\beta \rangle$ as required; we do it by induction on $\beta$, i.e.
we prove $(*)_{4,\beta}$.  
We now choose by induction on $i$ a pair of models 
$M_i \le_{\frak s} N_i$ such that
\mr
\item "{$\boxdot$}"   $N_i$ is $\le_{\frak s}$-increasing
continuous, $N_0 = N,M_i$ is 
$\le_{\frak s}$-increasing continuous, $M_0 = M$ and
$i = j+1$ implies $(M_j,M_i,a_j) \in 
K^{3,\text{uq}}_\lambda$ and every finite 
subset of $\bold J \backslash \{a_j:j < i\}$ is independent in
$(M,M_i,N_i)$ and $N_0 = N$ and $N_i$ is $\le_{\frak s}$-increasing
continuous.
\endroster
\bn
\ub{Subcase a}:  For $i=0$ there is no problem.
\bn
\ub{Subcase b}:    For $i$ limit let $M_i = \dbcu_{j < i}M_j$,
the least trivial part is the clause in $\boxdot$ on 
independence.  As for $i=\beta$ this clause is trivial, in fact
$(*)_{4,\beta}$ is already proved.  We can
assume $i < \beta$ and let $\{a_{i_0},\dotsc,a_{i_{n-1}}\} \subseteq
\bold J \backslash \{a_j:j < i\}$ be with no repetitions.  Now if 
$i < \mu$ then by renaming without loss of generality
max$[\{i_\ell:\ell < n\} \cup \{i\}] < \mu$ so we can use our induction
hypothesis on $\mu$.  So we can assume $\beta > \mu \vee (i = \beta =
\mu)$ and the case $i = \beta = \mu$ is trivial and can be forgotten, 
so as we are inducting
on all listings of $\bold J$ of a given length we have actually proved 
$(*)_0 \Rightarrow (*)_{4,\beta}$ for $\beta = \mu$, 
so $\bold J$ in independent in $(M,N)$ (as witnessed 
by some list of length
$\mu$!).  So we can apply (3)$^-$ and get that $\{a_{i_0},\dotsc,
a_{i_{n-1}}\}$ is independent in $(M,M_i,N_i)$ as required.
\bn
\ub{Subcase c}: For $i=j+1$, as in ${\frak s}$ we know that 
$K^{3,\text{uq}}_{\frak s}$ is dense (as ${\frak s}$ is weakly
successful and good) we can find $N_i,M_i$ as 
required and the induction assumption on $i$ 
holds as we have proved the claims for finite $\bold J$.  \nl
So we have finished the induction on $i$, thus
proving $(*)_0 \Rightarrow (*)_4$ for $\bold J$ of cardinality $\le
\mu$.  Hence we have part (1) for $\mu$ by $\boxtimes_1 +
\boxtimes_4 + \boxtimes_5$ 
from the beginning, as $(*)_0 \Rightarrow (*)_4$ was just
proved and $(*)_5 \Rightarrow (*)_6$ was proved inside the proof of
the finite case.  Hence we have proved part (1).
Now part (3) for $\mu$ follows from (1) + (3)$^-$.  So we have
finished the induction step for $\mu$ also in the infinite case (case
3) so have finished the proof.  \hfill$\square_{\scite{705-5.2}}$\margincite{705-5.2}
\bn
Still to finish the proof of \scite{705-5.2} we have to show
\scite{705-5.3}(2), \scite{705-5.3A}(1) below.
\proclaim{\stag{705-5.3} Claim}  1) [Assume ${\frak s}$ has primes.]
If $\langle N_i,a_j:i \le \alpha,j < \alpha
\rangle$ is a $M$-based pr-decomposition for ${\frak s}$ inside $N$ (so
$N_\alpha \le_{\frak s} N$, see Definition \scite{705-c.1A}), 
$b \in N$, {\rm tp}$_{\frak s}(b,N_\alpha,N)$ is a nonforking
extension of $p \in {\Cal S}^{\text{bs}}_{\frak s}(N_0)$, 
\ub{then} we can find
$\langle N'_i,a'_j:i \le 1 + \alpha,j < 1 + \alpha \rangle$ an 
$M$-based pr-decomposition for ${\frak s}$ inside $N^+$, 
such that $N_i \le N'_{1+i},N_0 = N'_0,
b = a'_0,a_i = a'_{1+i}$, 
{\rm tp}$_{\frak s}(a_i,N'_{1+i},N'_{1+i})$ does
not fork over $N_i,N'_\alpha \le_{\frak s} N^+,N \le_{\frak K}
N^+$. \nl
2) Similarly for {\rm uq}-decomposition only $N_0 \le_{\frak s} N'_0$
(instead equality).
\endproclaim
\bigskip

\demo{Proof}  1) Chasing arrows: first ignore $b = a'_0$, demand just 
tp$_{\frak s}(a'_0,N_0,N'_1) = \text{tp}(b,N_0,N)$ 
and ignore $N \le_{\frak s} N^+$.  After proving this we
can use equality of types.

In details, we choose by induction on $i \le \alpha$ a pair $(N^*_i,f_i)$ and
$b^*,a^*_i$ (if $i < \alpha$) such that:
\mr
\item "{$(a)$}"  $f_i$ is a $\le_{\frak s}$-embedding of $N_i$ into $N^*_i$
\sn
\item "{$(b)$}"  $N^*_i$ is $\le_{\frak s}$-increasing continuous
\sn
\item "{$(c)$}"  $N^*_0$ satisfies $(N_0,N^*_0,b^*) \in
K^{3,\text{pr}}_{\frak s}$ and {\rm tp}$_{\frak s}
(b^*,N_0,N^*_0) = \text{ tp}_{\frak s}
(b,N_0,N)$
\sn
\item "{$(d)$}"  $f_0 = \text{ id}_{N_0}$
\sn
\item "{$(e)$}"  if 
$i=j+1$ then $(N^*_j,N^*_i,a^*_j) \in K^{3,\text{pr}}_{\frak s}$
and $f_i(a_j) = a^*_j$
\sn
\item "{$(f)$}"  tp$_{\frak s}(b^*,f_i(N_i),N^*_i)$ 
does not fork over $N_0$.
\ermn
For $i=0$ just use ``${\frak s}$ has primes".

For $i=j+1$ first choose $p_j \in 
{\Cal S}^{\text{bs}}_{\frak s}(N^*_j)$, a
nonforking extension of $p^-_j = 
f_j(\text{tp}_{\frak s}(a_j,N_j,N_{j+1}))$ and
second choose $N^*_i,a^*_j$ such that $(N^*_j,N^*_i,a^*_j) 
\in K^{3,\text{pr}}_{\frak s}$ and 
tp$_{\frak s}(a^*_j,N^*_j,N^*_i) = p_j$ (using ``${\frak s}$ 
has primes"), lastly choose $f_i \supseteq f_j$ mapping $a_j$ to
$a^*_j$ using the assumption $(N_j,N_i,a_j) \in 
K^{3,\text{pr}}_{\frak s}$ and $a^*_j$'s realizing $p_j$).

For $i$ limit take union. \nl
Having finished the induction 
\wilog \, each $f_i$ is the identity on $N_i$
so $j < \alpha \Rightarrow a^*_j = a_j$.  
So $N_\alpha \le_{\frak s} N^*_\alpha$ and $N_\alpha
\le_{\frak s} N$ and tp$(b^*,N_\alpha,N^*_\alpha)$ 
does not fork over $N_0$
(by clause (f)) and extend tp$_{\frak s}(b,N_0,N)$; also
tp$(b,N_\alpha,N)$ satisfies this so as $K_{\frak s}$ has amalgamation \wilog \,
$b = b^*$ and for some $N^+ \in K_{\frak s}$ we have $N \le_{\frak s} N^+
\and  N^*_\alpha \le_{\frak s} N^+$. \nl
Letting $N'_0 = N',N'_{1+i} = N^*_i$ we are done. \nl
2) As in \sectioncite[\S6]{600}.  \hfill$\square_{\scite{705-5.3}}$\margincite{705-5.3}
\enddemo
\bn
Some trivial properties are:
\proclaim{\stag{705-5.3A} Claim}  1) If $\langle M_i:i \le \alpha \rangle$ is
$\le_{\frak s}$-increasing continuous, and $\bold J_i$ is independent in
$(M_0,M_i,M_{i+1})$ \ub{then} $\cup\{\bold J_i:i < \alpha\}$ is independent in
$(M_0,M_\alpha)$. \nl
2) If {\rm NF}$_{\frak s}(M_0,M_1,M_2,M_3)$
and $\bold J$ is independent in $(M_0,M_1)$ \ub{then} $\bold J$ is
independent in $(M_2,M_3)$. 
If {\rm NF}$_{\frak s}(M_0,M_1,M_2,M_3),
M_0 \le_{\frak s} M^-_1 \le_{\frak s} M_1,
\bold J_1$ independent in $(M^-_1,M_1)$ and $\bold J_2$
independent in $(M_0,M_2)$ \ub{then} $\bold J_1 \cup \bold J_2$ is
independent in $(M^-_1,M_3)$.
\nl
3) [Monotonicity]  If 
$\bold J$ is independent in $(M,A,N)$ and $\bold I \subseteq
\bold J$, \ub{then} $\bold I$ is independent in $(M,A,N)$. \nl
4) If $\bold J$ is independent in $(M_1,M_3,N)$ and $M_0 \le_{\frak s}
M_\ell \le_{\frak s} M_3$ for $\ell = 1,2$ and $c \in \bold J
\Rightarrow { \text{\rm tp\/}}_{\frak s}(c,M_1,N)$ 
does not fork over $M_0$,
\ub{then} $\bold J$ is independent in $(M_0,M_2,N)$. \nl
5) [${\frak s}$ has primes].
Assume that {\rm NF}$_{\frak s}(M_0,M_1,M_2,M_3)$ and $\langle
M_{0,i},a_j:i \le \alpha,j < \alpha \rangle$ 
is a decomposition of $M_2$ over $M_0$.  \ub{Then} we
can find $M^+_3$ satisfying $M_3 \le_{\frak s} M^+_3$ and $\langle
M_{1,i}:i \le \alpha \rangle$ such that $M_{1,\alpha} \le_{\frak s}
M^+_3,\langle M_{1,\alpha},a_i:i < \alpha \rangle$ is a decomposition
of $M_{1,\alpha}$ over $M_1$ and {\rm tp}$(a_i,M_{1,i},M^+_3)$ 
does not fork over $M_{0,i}$. \nl
6) Similarly for {\rm uq}-decompositions except that $M_1 \le_{\frak s}
M_{1,0}$ (not necessarily equal). \nl
7) The set $\{a\}$ is 
${\frak s}$-independent in $(M,N)$ \ub{iff} $(M,N,a) \in
K^{3,\text{bs}}_{\frak s}$. 
\endproclaim
\bigskip

\demo{Proof}  1) Should be clear (e.g., \wilog \, $M_{i+1}$ is
$(\lambda,*)$-brimmed over $M_i \cup \bold J_i$). \nl
2) The first phrase is \scite{705-5.1}(2).  
The second phrase by chasing arrows.
\nl
3) Trivial. \nl
4) Easy by the nonforking calculus. \nl
5) As in the proof of the previous claim there is a sequence $\langle
M_{1,i}:i \le \alpha \rangle$ and a $\le_{\frak s}$-embedding
$f_\alpha$ of $M_2 = M_{0,\alpha}$ into $M_{1,\alpha}$ such that
$M_{1,0} = M_1, f \restriction M_0 = \text{ id}_{M_0}$ and
$f_i(M_{0,i}) \le_{\frak s} M_{1,i}$ and tp$_{\frak s}
(f(a_i),M_{1,i},M_{1,i+1})$ does not fork over $f(M_{0,i})$ and
$(M_{1,i},M_{1,i+1},f_{i+1}(a_i)) \in K^{3,\text{pr}}_{\frak s}$ 
(e.g., just
choose $M_{1,i},f_i = f \restriction M_{0,i}$ by induction on $i \le
\alpha$).    As
$(M_{0,i},M_{0,i+1},a_i) \in K^{3,\text{pr}}_{\frak s}$ also
$(f(M_{0,i}),f(M_{0,i}),a_i)$ belongs to $K^{3,\text{pr}}_{\frak s}$ 
hence to $K^{3,\text{uq}}_{\frak s}$ hence NF$_{\frak
s}(f(M_{0,i}),f(M_{0,i+1}),M_{1,i},M_{1,i+1})$ hence by long
transitivity NF$(f(M_{0,0}),f(M_{0,\alpha}),M_{1,0},M_{1,i})$.  By the
uniqueness of NF, \wilog \, $f$ the identity and for some $M^+_3$ we
have $M_{1,\alpha} \le_{\frak s} M^+_3$ and $M_3 \le_{\frak s}
M^+_3$. \nl
6) As in \sectioncite[\S6]{600}. \nl
7) Trivial by the definitions.  \hfill$\square_{\scite{705-5.3A}}$\margincite{705-5.3A}
\enddemo
\bn
\centerline {$* \qquad * \qquad *$}
\bigskip
 
\definition{\stag{705-5.4} Definition}  1) We say $N$ is prime over $M \cup 
\bold J$ (for ${\frak s}$), or $(N,M,\bold J) \in 
K^{3,\text{qr}}_{\frak s}$ if:
\mr
\item "{$(a)$}"  $M \le_{\frak s} N$ in $K_\lambda$
\sn
\item "{$(b)$}"  $\bold J \subseteq \bold I_{M,N}$
\sn
\item "{$(c)$}"  if $M \le_{\frak s} N',\bold J' \subseteq \bold
I_{M,N'}$ and $\bold J'$ is independent in $(M,N')$ and
$h$ is a one to one mapping from $\bold J$ onto $\bold J'$ such that
tp$(a,M,N) = \text{ tp}(h(a),M,N')$ for every $a \in \bold J$, \ub{then}
there is a $\le_{\frak s}$-embedding of $N$ into $N'$ over $M$ extending
$h$.
\ermn
2) Let $(M,N,\bold J) \in K^{3,\text{bs}}_{\frak s}$ 
means that $M \le_{\frak s} N$ and $\bold J$ is independent 
in $(M,N)$, (see \scite{705-5.1}).
\enddefinition
\bn
Some basic properties are
\proclaim{\stag{705-5.5} Claim}  1) [Assume ${\frak s}$ has primes.]
If $M \le_{\frak s} N$ (in $K_\lambda$) and
$\bold J \subseteq \bold I_{M,N}$ is independent in $(M,N)$, \ub{then} there
is $N' \le_{\frak s} N$ which is prime over $M \cup \bold J$. \nl
2) If $\bold J$ is independent in $(M,N)$ and 
$\langle M_i,a_j:i \le \alpha,j < \alpha \rangle$ is 
an $M$-based pr-decomposition for 
$\bold J$ of $(M_0,N)$ (see Definition \scite{705-5.4}(2)),
\ub{then} $M_\alpha$ is prime over $M_0 \cup \bold J$. \nl
3) $(M,N,\{a\}) \in K^{3,\text{qr}}_{\frak s}$ iff $(M,N,a) \in
K^{3,\text{pr}}_{\frak s}$. \nl
4) If $(M,N,\bold J) \in K^{3,\text{qr}}_{\frak s}$ and $M \cup \bold J
\subseteq N^- \le_{\frak s} N$ \ub{then} $(M,N^-,\bold J) \in 
K^{3,\text{qr}}_{\frak s}$.
\endproclaim
\bigskip

\demo{Proof}  1) By \scite{705-5.2}(1), $(*)_1 \Leftrightarrow (*)_6$,
letting $\langle a_i:i < \alpha \rangle$ list $\bold J$ we can find $M_i
\le_{\frak s} N$ for $i < \alpha$ such that $\langle M_i,a_j:i \le
\alpha,j < \alpha \rangle$ as in $(*)_6$ of \scite{705-5.2}.  Now we can
use part (2).   
\nl
2) Let $M_0 \le_{\frak s} N^*$ and a one-to-one function 
$h:\bold J \rightarrow \bold J' \subseteq N^*$ satisfying $c \in \bold
J \Rightarrow \text{ tp}_{\frak s}(h(c),M_0,N^*) = \text{ tp}_{\frak
s}(c,M_0,M_\alpha)$ and $\bold J'$ independent in $(M_0,N^*)$ be
given.  Let $h(a_j) = c_j$.
We now define by induction on 
$i \le \alpha$ a $\le_{\frak s}$-embedding $f_i$
of $M_i$ into $N^*$, increasing continuous with $i$ 
and mapping $a_j$ to $c_j$.
For $i=0$ this is given, for $i$ limit take union.  For $i=j+1$, we know
that tp$(c_j,f_j(M_j),N^*)$ does not fork over $N_0 = f_0(M_0)$ by
\scite{705-5.2}(3) (because $(M_j,M_{j+1},a_j) \in 
K^{3,\text{uq}}_{\frak s}$ by
\scite{705-c.4}(2)) and so as tp$(a_j,M_j,N_1)$ does not fork over $M_0$ and
$f_0[\text{tp}(a_j,M_0,M_1)] = \text{ tp}(c_j,N_0,N_1)$ clearly 
$f_j[\text{tp}(a_j,M_j,M_\alpha)] = \text{ tp}(c_j,f_j(M_j),N^*)$.  But
$(M_j,M_{j+1},a_j) \in K^{3,\text{pr}}_\lambda$ so we can find $f_i \supseteq f_j$
as required.  So $f_\alpha$ is as required in Definition
\scite{705-5.4}. \nl
3) By the definition. \nl
4) Easy, like in \scite{705-c.2}(1).   \hfill$\square_{\scite{705-5.5}}$\margincite{705-5.5}
\enddemo
\bn
\centerline {$* \qquad * \qquad *$}
\bigskip

\proclaim{\stag{705-4.8} Claim}  1) If $(M_0,M_1,\bold J) \in
K^{3,\text{qr}}_{\frak s}$ 
hence $\bold J$ is independent in $(M_0,M_1)$ and $N_0
\le_{\frak s} N_1,f_0$ is an isomorphism from $M_0$ onto $N_0$ and 
$\{c_a:a \in \bold J\}$ is an independent set in $(N_0,N_1)$ 
satisfying {\rm tp}$(c_a,N_0,N_1)
= f_0[{\text{\rm tp\/}}(a,M_0,M_1)]$ for $a \in \bold J$ and of course
$\langle c_a:a \in \bold J \rangle$ is with no repetitions 
\ub{then} there is a $\le_{\frak s}$-embedding $f$ of $M_1$ 
into $N_1$ extending $f_0$ and mapping each $b \in
\bold J$ to $c_b$. \nl
2) [Assume ${\frak s}$ has primes].   
Assume $(M_0,M_1,\bold J) \in K^{3,\text{qr}}_{\frak s}$ and 
$M_0 \le_{\frak s} M_2 \le_{\frak s} M_3,M_0 \le_{\frak s} M_1 
\le_{\frak s} M_3$ and $\bold J$ is
independent in $(M_0,M_2,M_3)$.  \ub{Then} 
{\rm NF}$_{\frak s}(M_0,M_1,M_2,M_3)$.
\endproclaim
\bigskip

\demo{Proof}  1) This just rephrases Definition \scite{705-5.4}. \nl
2) We are allowed to increase $M_1,M_3$, i.e. if $M_3 \le_{\frak s}
M'_3,M'_1 
\le_{\frak s} M'_3,M_1 \le_{\frak s} M'_1$ and $(M_0,M'_1,\bold J) \in 
K^{3,\text{qr}}_{\frak s}$ 
then we can replace $M_1,M_3$ by $M'_1,M'_3$.  
So by \scite{705-5.5}(2) amalgamation in ${\frak K}_{\frak s}$ and 
the definition of primes \wilog \, we can find 
$\langle M^0_i,a_j:i \le \alpha,i <
\alpha \rangle$ such that $M^0_0 = M_0,M^0_\alpha = M_1,M^0_i$ is
$\le_{\frak s}$-increasing continuous, $(M^0_i,M^0_{i+1},a_i) \in
K^{3,\text{pr}}_{\frak s} \subseteq K^{3,\text{uq}}_{\frak s}$ 
for $i < \alpha$ and $\langle a_i:i < \alpha \rangle$
list $\bold J$ with no repetitions.  We now choose by induction on
$i \le \alpha,M^2_i,M^3_i$ such that:
\mr
\item "{$(\alpha)$}"  $M^2_i$ is $\le_{\frak s}$-increasing continuous
\sn
\item "{$(\beta)$}"  $M^3_i$ is $\le_{\frak s}$-increasing continuous
\sn
\item "{$(\gamma)$}"  $M^2_0 = M_2,M^3_0 = M_3$
\sn
\item "{$(\delta)$}"  $M^0_i \le_{\frak s} M^2_i \le_{\frak s} M^3_i$
\sn
\item "{$(\varepsilon)$}"  $(M^2_i,M^2_{i+1},a_i) \in 
K^{3,\text{pr}}_{\frak s}$
\sn
\item "{$(\zeta)$}"  tp$_{\frak s}
(a_i,M^2_i,M^3_i)$ does not fork over $M^0_i$
\sn
\item "{$(\eta)$}"  $M^0_i \le_{\frak s} M^3_i$.
\ermn
Why is this enough?  For each $i$ we have $(M^0_i,M^0_{i+1},a_i) \in
K^{3,\text{uq}}_{\frak s}$ 
(by the choice of $M^0_i,M^0_{i+1},a_i$ we have $(M^0_i,M^0_{i+1},a_i) 
\in K^{3,\text{pr}}_{\frak s}$ and use claim \scite{705-c.4}(2)).  By
this, \cite[\S6]{600},\marginbf{!!}{\cprefix{600}.\scite{600-nf.20t}} and clauses 
$(\delta), (\varepsilon), (\zeta)$ and the definition of
$K^{3,\text{uq}}_{\frak s}$ we have \nl
NF$_{\frak s}(M^0_i,M^0_{i+1},M^2_i,M^2_{i+1})$.
By the symmetry property of NF$_{\frak s}$ we have \nl
NF$_{\frak s}(M^0_i,M^2_i,M^0_{i+1},M^2_{i+1})$.
As this holds for every $i < \alpha$ and clauses 
$(\alpha) + (\beta)$ by the long transitivity 
property of NF$_{\frak s}$ (see \cite[\S6]{600},\marginbf{!!}{\cprefix{600}.\scite{600-nf.16t}}(1))
we get NF$_{\frak s}(M^0_0,M^2_0,M^0_\alpha,M^2_\alpha)$, which means
NF$_{\frak s}(M_0,M_2,M_1,M^2_\alpha)$.  Now 
by monotonicity we
can replace $M^2_\alpha$ first by $M^3_\alpha$ then by $M_3$ so we got
NF$_{\frak s}(M_0,M_2,M_1,M_3)$ as required. \nl
Why is it possible to carry 
the induction?  Having arrive to $i$ the type
tp$_{\frak s}(a_i,M^2_i,M^3_i)$ does not fork 
over $M$ by \scite{705-5.2}(3) and then
we can find a $M^2_{i+1} \le_{\frak s} 
M^3_{i+1}$ such that $M^3_i \le_{\frak s} 
M^3_{i+1}$ and $(M^2_i,M^2_{i+1},a_i) \in K^{3,\text{pr}}_{\frak s}$.
Now by the definition of prime, there is a $\le_{\frak s}$-embedding
$f_i$ of $M^0_{i+1}$ into $M^3_{i+1}$ over 
$M^0_i$ satisfying $f_i(a_i) = a_i$.
Chasing arrows \wilog \, $M^0_{i+1} \le_{\frak s} M^2_i$.
\hfill$\square_{\scite{705-4.8}}$\margincite{705-4.8}
\enddemo
\bigskip

\proclaim{\stag{705-4.11} Claim}   Assume $\langle M_i:i \le \delta
+1 \rangle$ is $\le_{\frak s}$-increasing continuous and $\bold J \subseteq
\bold I_{M_\delta,M_{\delta +1}}$. \nl
1) If $|\bold J| < { \text{\rm cf\/}}(\delta)$ 
and $\bold J$ is independent in
$(M_\delta,M_{\delta +1})$ \ub{then} for every $i < \delta$ large enough,
$\bold J$ is independent in $(M_i,M_\delta,M_{\delta +1})$. \nl
2) If $\bold J \subseteq \bold I_{M_\delta,M_{\delta +1}}$ is
independent in $(M_i,M_{\delta +1})$ for every $i < \delta$, 
\ub{then} $\bold J$ is independent in $(M_\delta,M_{\delta +1})$. 
\endproclaim
\bn
\ub{Discussion}:  At this point, if $\langle M_i:i \le \alpha \rangle$
is $\le_{\frak s}$-increasing continuous $M_\alpha \le N,a \in N$,
tp$(a,M_\alpha,N) \in {\Cal S}^{\text{bs}}(M_\alpha)$ does not fork over
$M_0$ we do not know if there is $\langle N_i:i \le \alpha \rangle$
which is $\le_{\frak s}$-increasing continuous $a \in N_0,M_i
\le_{\frak s} N_i$ and $(M_i,N_i,a) \in 
K^{3,\text{pr}}_{\frak s}$.  So we go
around this.  This claim is used in \scite{705-5.12}.
\bigskip

\demo{Proof}  1) For each $c \in \bold J$ for some $i_c \in \delta$,
tp$_{\frak s}(c,M_\delta,M_{\delta +1})$ does not fork over $M_{i_c}$
let $i(*) = \sup\{i_c:c \in \bold J\}$ and use \scite{705-5.2}(4). \nl
2) By \scite{705-5.2}(1) it suffices to deal with finite $\bold J$, say
$\bold J = \{b_\ell:\ell < n\}$ with no repetitions. \nl
By the NF$_{\frak s}$-calculus there is a $\le_{\frak
s}$-increasing continuous sequence $\langle M^+_i:i \le \delta + 1
\rangle$ such that NF$_{\frak s}(M_i,M^+_i,M_j,M^+_j)$ for any $i < j
\le \delta +1$ and $M^+_{i+1}$ is $(\lambda,*)$-brimmed over
$M_{i+1} \cup M^+_i$ for $i \le \delta$ hence $M^+_\delta$ is
$(\lambda,*$)-brimmed over $M_\delta$ see \scite{705-1.15K}.  
Hence there is a $\le_{\frak s}$-embedding 
$h$ of $M_{\delta +1}$ into $M^+_\delta$ over
$M_\delta$, so \wilog \, $M_{\delta +1} \le_{\frak s} M^+_\delta$.  As
$\bold J$ is finite and $M^+_\delta$ is the union of the $\le_{\frak
s}$-increasing sequence $\langle M^+_i:i < \delta \rangle$ clearly for
some $i < \delta$ we have $\bold J \subseteq M^+_i$ hence ``$\bold J$
is independent in $(M_i,M^+_i)$.  But
NF$(M_i,M^+_i,M_\delta,M^+_\delta)$ hence by Claim \scite{705-5.1R} we
deduce ``$\bold J$ is independent in $(M_i,M_\delta,M^+_\delta)$ hence
in $(M_i,M_\delta,M_{\delta +1})$ as required.
\hfill$\square_{\scite{705-4.11}}$\margincite{705-4.11}
\enddemo
\bigskip

\proclaim{\stag{705-4.11A} Claim}  Assume 
${\frak s} = {\frak t}^+,{\frak t}$ a
successful $\lambda_{\frak t}$-good$^+$ frame, $\lambda = 
\lambda^+_{\frak t}$. \nl
Assume further $M_\ell \in K_{\frak s}$ and 
$\langle M^\ell_\alpha:\alpha < \lambda \rangle$ is
a $\le_{{\frak K}[{\frak t}]}$-representation of $M_\ell$ for $\ell =1,2$.

If $M_0 \le_{\frak s} M_1$ and $\bold J \subseteq \bold I_{M_0,M_1}$
\ub{then}: $\bold J$ is independent in $(M_0,M_1)$ \ub{iff} 
for a club of $\delta < \lambda$ the set 
$\bold J \cap M_{1,\delta}$ is independent in $(M_{0,\delta},
M_{1,\delta})$ (for ${\frak t}$) \ub{iff} for stationarily many $\delta 
< \lambda,\bold J \cap M_{1,\delta}$ is independent in $(M_{0,\delta},
M_{1,\delta})$ (for ${\frak t}$).
\endproclaim
\bigskip

\demo{Proof}   By \scite{705-5.2}(1), applied to ${\frak s}$ and to
${\frak t}$ \wilog \, $\bold J$ is finite.

Using $(*)$ of Definition \scite{705-5.1}, the first clause implies 
the second clause (by the definition and \scite{705-stg.11}) 
and trivially second implies third.  So it suffices to
assume the failure of the first and show the failure of the third.  
Let $\bold J = \{a_\ell:\ell <
n\}$ without repetitions. We can try to choose by induction on $\ell$ a 
model $M'_\ell \le_{\frak s} M_1$ 
such that $M'_0 = M_0,(M'_\ell,M'_{\ell +1},
a_\ell) \in K^{3,\text{pr}}_{\frak s}$, moreover is as 
constructed in \scite{705-d.8}(1) + \scite{705-d.1}
and tp$_{\frak s}(a_\ell,M'_\ell,M_1) \in 
{\Cal S}^{\text{bs}}_{\frak s}(M'_\ell)$ does not
fork over $M_0$.  We cannot succeed so for some $m <n$ we have $M'_0,\dotsc,
M'_m$ as above but tp$(a_m,M_m,N)$ forks over $M_0$.  Rename $M_1$ as
$M'_{m+1}$ and let $\langle M'_{\ell,\alpha}:\alpha \le \lambda
\rangle$ be a $\le_{\frak t}$-representation of $M'_\ell$ for $\ell
\le m+1$ and $M'_{0,\alpha} = M^0_\alpha,M'_{m+1} = M^1_\alpha$.
Now by \scite{705-d.9}(1) for some club $E$ of
$\lambda$, if $\delta$ is from $E$ and $\ell < m$ then $(M'_{\ell,\delta},
M'_{\ell +1,\delta},a_\ell) \in K^{3,\text{uq}}_{\frak s}$ 
and $a_m \in 
M'_{m+1,\delta}$ and tp$_{\frak s}(a_m,M_{m,\delta},M_{m+1,\delta})$ forks over 
$M_{0,\delta}$ while $M_{m,\delta} \le_{\frak t} M_{m+1,\delta}$.
By \scite{705-5.2}(3) for ${\frak t}$
we get $\{a_\ell:\ell \le m\} \subseteq \bold I_{M_{0,\delta},M_{m,\delta}}$
is not independent.  So we have gotten the
failure of the third clause.  \hfill$\square_{\scite{705-4.11A}}$\margincite{705-4.11A}
\enddemo
\bigskip

We can deal with dimension as in \cite[Ch.III]{Sh:c}.
\definition{\stag{705-4.25Y} Definition}  1) Assume that $M \le_{\frak s}
N$ and $p \in {\Cal S}^{\text{bs}}(M)$, then we let 

$$
\align
\text{dim}(p,N) = \text{ Min}\{|bold J|:&\bold J \text{ satisfies} \\
  &(i) \quad \bold J \text{ is a subset of } \{c \in
  N:\text{tp}(c,M,N) \text{ is equal to } p\}, \\
  &(ii) \quad \text{ the triple } (M,N,\bold J) \text{ belongs to }
K^{3,\text{bs}}_{\frak s} \text{ and } \\ 
  &(iii) \quad \bold J \text{ is maximal under those restrictions}\}.
\endalign
$$
\mn
We shall say more on dim after we understand regular types but for now ??
\enddefinition
\bigskip

\proclaim{\stag{705-4.25} Claim}  Assume $M \in K_{\frak s}$ 
and $\bold J$ is independent in $(M,N^*)$.
\nl
1) If {\rm tp}$(a,M,N^*) \in {\Cal S}^{\text{bs}}(M)$, 
\ub{then} for some finite $\bold J' \subseteq \bold J$ the set
$(\bold J \backslash \bold J') \cup \{a\}$ is independent in
$(M,N^*)$.
\nl
2) If $a \in N^*$ \ub{then} for some finite $\bold J' \subseteq \bold J$ and
$M'$ we have: $M \cup \{a\} \subseteq M' \le_{\frak s} N$ and $\bold J
\backslash \bold J'$ is independent over $(M,M')$ in $N^*$. \nl
\endproclaim
\bigskip

\demo{Proof}  1) Let 
$\bold J = \{a_i:i < \alpha\}$, we prove the statement by
induction on 
$\alpha$. For $\alpha = 0,\alpha$ successor this is trivial.  For
$\alpha$ limit $< \lambda^+$ by the definition there is 
a $\le_{\frak s}$-increasing continuous sequence $\langle M_i:i \le
\alpha \rangle$ such that $M_0 = M,M_i \le_{\frak s} N^+,N^*
\le_{\frak s} N^+$ and tp$(a_i,M_i,M_{i+1}) \in {\Cal
S}^{\text{bs}}_{\frak s}(M_i)$ does not fork over $M_i$.  
As in the proof of \scite{705-4.11}(2) \wilog \, $M_\alpha = N^*$.

Now for some $\beta < \alpha,a \in M_\beta$ and by the induction hypothesis
on $\alpha$ for some finite $u \subseteq \beta$ we have 
$\{a_i:i \in \beta \backslash u\} \cup \{a\}$ is independent in 
$(M,M_\beta)$.  Clearly by the Definition \scite{705-5.1} the set
$\{a_i:i \in \alpha \backslash \beta\}$ 
is independent in $(M,M_\alpha,M_\beta)$.
By \scite{705-5.3A}(1) and
the last two sentences $(\{a_i:i \in \beta \backslash u\} \cup
\{a\}) \cup (\{a_i:i \in \alpha \backslash \beta\})$ is 
independent in $(M,M_\alpha)$
hence in $(M,N^+)$ hence in $(M,N^*)$ by \scite{705-5.2}(2).  
But the set is $\{a_i:i \in \alpha \backslash u\} \cup \{a\}$ so
we are done.
\nl
2) Similar.  \hfill$\square_{\scite{705-4.25}}$\margincite{705-4.25}
\enddemo
\bigskip

\proclaim{\stag{705-4.25A} Claim}  Assume that $M \le_{\frak s} N$ and $p
\in {\Cal S}^{\text{bs}}(M)$.  Then  any two sets $\bold J$ satisfying
the demands (i) + (ii) from \scite{705-4.9}(2) have the same cardinality
or are both finite.
\endproclaim
\bigskip

\demo{Proof}  By \scite{705-4.25}.
\enddemo
\bigskip

\definition{\stag{705-4.9} Definition}  Let 
$(M,N,\bold J) \in K^{3,\text{vq}}_{\frak s}$
means:
\mr
\item "{$(a)$}"  $M \le_{\frak s} N$
\sn
\item "{$(b)$}"  $\bold J$ is independent in $(M,N)$
\sn
\item "{$(c)$}"  if $N \le_{\frak s} M_3,M \le_{\frak s} M_2 \le_{\frak s}
M_3$ and $\bold J$ is independent in $(M_0,M_2,M_3)$ \ub{then}
NF$_{\frak s}(M,M_1,N,M_3)$.
\endroster
\enddefinition
\bn
So by \scite{705-4.8}(2)
\proclaim{\stag{705-4.10} Claim}  1) If 
$(M,N,\bold J) \in K^{3,\text{qr}}_{\frak s}$ \ub{then}
$(M,N,\bold J) \in K^{3,\text{vq}}_{\frak s}$. \nl
2) If $(M_i,N_i,\bold J_i) \in K^{3,\text{qr}}_{\frak s}$ for $i <
\delta,\delta$ a limit ordinal $< \lambda^+_{\frak s},M_i$ is
$\le_{\frak s}$-increasing, $N_i$ is $\le_{\frak s}$-increasing,
$\bold J_i$ is $\subseteq$-increasing \ub{then}
\mr
\item "{$(a)$}"  $(M_\delta,N_\delta,\bold J_\delta) \in
K^{3,\text{qr}}_{\frak s}$ when we let $M_\delta = \cup\{M_i:i <
\delta\},N_\delta =$ \nl
$\cup\{N_i:i < \delta\},\bold J_\delta = \cup\{\bold J_i:i < \delta\}$
\sn
\item "{$(b)$}"  for $k<j \le \delta$,{\rm NF}$_{\frak
s}(M_i,N_i,M_j,N_j)$
\sn  
\item "{$(c)$}"  $\alpha \le j \le \delta$ is a limit ordinal then
NF$_{\frak s}(\dbcu_{i < \alpha} M_i,\dbcu_{i < \delta}
M_i,M_j,N_j)$. 
\endroster
\endproclaim
\bigskip

\demo{Proof}  1),2) Clauses (b),(c) hold by \scite{705-4.8}(2).
\enddemo
\bigskip

\proclaim{\stag{705-4.11K} Claim}  [${\frak s} = {\frak t}^+,{\frak t}$ 
a successful good $\lambda_t$-frame, so $\lambda = 
\lambda^+_{\frak t}$ and ${\frak t}$ satisfies Hypothesis \scite{705-d.0}].
\nl
Assume
\mr
\item "{$(a)$}"  $\langle M^\ell_\alpha:\alpha < \lambda_{\frak s}
\rangle$ is a $\le_{\frak s}$-representation $M_\ell \in K_{\frak s}$ 
for $\ell = 1,2$
\sn
\item "{$(b)$}"  $M_1 \le_{\frak s} M_2$
\sn
\item "{$(c)$}"  $\bold J \subseteq \bold I_{M_1,M_2}$ is independent
in $(M_1,M_2)$ and let $\bold J_\alpha = \bold J \cap M^2_\alpha$.
\ermn
\ub{Then} $(M_1,M_2,\bold J) \in K^{3,\text{qr}}_{\frak s}$ \ub{iff} 
for a club of $\delta < \lambda_{\frak s}$ the triple
$(M^1_\alpha,M^2_\alpha,\bold J_\alpha) \in 
K^{3,\text{vq}}_{\frak t}$
\ub{iff} for stationarily many $\delta < \lambda_{\frak s}$ the triple
$(M^1_\alpha,M^2_\alpha,\bold J_\alpha)$ belongs to 
$K^{3,\text{vq}}_{\frak t}$.
\endproclaim
\bigskip

\demo{Proof}  Like \scite{705-d.8} we can prove ``second implies first".
The second implies the third trivially, the third implies the second
by \scite{705-4.10}(2), clause (a).  
If the third fails, without loss of generality the failure is for
every $\alpha < \lambda^+_{\frak s}$ in particular, $\bold J_\alpha
\ne \emptyset$ and \wilog \, $M^\ell_{\alpha +1}$ is
$(\lambda,*)$-brimmed over $M^\ell_\alpha$ for $\alpha <
\lambda^+_{\frak s},\ell=1,2$.

By the definition of ``$(M^1_0,M^2_0,\bold J_0) \in
K^{3,\text{vq}}_{\frak t}"$ we can find $N_1,N_2$ such that:
\mr
\item "{$(*)_1$}"  $(a) \quad M^1_0 \le_{\frak t} N_1 \le_{\frak t}
N_2$
\sn
\item "{${{}}$}"  $(b) \quad M^2_0 \le_{\frak t} N_2$
\sn
\item "{${{}}$}"  $(c) \quad \bold J_0$ is independent in
$(M^1_0,N_1,N_2)$
\sn
\item "{${{}}$}"  $(d) \quad \neg \text{NF}_t(M^1_0,M^2_0,N_1,N_2)$.
\ermn
Without loss of generality
\mr
\item "{$(*)_2$}"  $N_1$ is $(\lambda,*)$-brimmed over $M^1_0$ in
${\frak t}$.
\ermn
Hence we can find $\langle N^\ell_i:i \le \lambda \rangle$ for
$\ell=1,2$ and $\langle a_i:i < \lambda \rangle$ such that
\mr
\item "{$(*)_3$}"  $(a) \quad \langle N^\ell_i:i \le \lambda \rangle$
is $\le_{\frak t}$-increasing continuous
\sn
\item "{${{}}$}"  $(b) \quad N^\ell_0 = M^\ell_0,N^\ell_\lambda =
N_\ell$ for $\ell=1,2$
\sn
\item "{${{}}$}"  $(c) \quad N^1_i \le_{\frak s} N^2_i$
\sn
\item "{${{}}$}"  $(d) \quad (N^1_i,N^2_i,a_i) \in
K^{3,\text{uq}}_{\frak t}$ \nl
(why? see \marginbf{!!}{\cprefix{600}.\scite{600-nf.0at}}).
\ermn
Now by induction on $i \le \lambda$ we choose $(M_{1,i},M_{2,i},f_i)$
and if $i$ is a successor also $b_{i-1}$ such that
\mr
\item "{$(*)_4$}"  $(a) \quad \langle M_{\ell,i}:i \le \lambda
\rangle$ is $\le_{\frak s}$-increasing continuous for $\ell=1,2$
\sn
\item "{${{}}$}"  $(b) \quad M_{1,i} \le_{\frak s} M_{2,i}$
\sn
\item "{${{}}$}"  $(c) \quad M_{\ell,0} = M_\ell$ or $\ell=1,2$
\sn
\item "{${{}}$}"  $(d) \quad (M_{ell,i},M_{\ell,i+1},b_i) \in
K^{3,\text{pr}}_{\frak s}$ for $\ell=1$
\sn
\item "{${{}}$}"  $(e) \quad$ tp$_{\frak s}(b_i,M_{2i},M_{2,i+1})$
does not fork over $M_{1,i}$
\sn
\item "{${{}}$}"  $(f) \quad f_i$ is a $\le_{{\frak K}[{\frak
t}]}$-embedding of $N^2_i$ into $M_{2,i}$
\sn
\item "{${{}}$}"  $(g) \quad f_i$ maps $N^1_i$ into $M_{1,i}$ and maps
$a_j$ to $b_j$ for $j<i$
\sn
\item "{${{}}$}"  $(h) \quad$ if $j<i$ then $f_i(\text{tp}_{\frak
t}(a_j,N^1_j,N^1_{j+1}))$ is a witness for tp$_{\frak
s}(b_j,M_{1,j},m_{1,j+1})$.
\ermn
For $i=0$ this is trivial, for $i$ limit take union.  For $i=j+1,q_j =
f_j(\text{tp}_{\frak t},N^1_j,N^1_i)) \in {\Cal S}^{\text{bs}}_{\frak
t}(f_j(^1_j))$ is well defined and so there is $p_j \in {\Cal
S}^{\text{bs}}_{\frak s}(M_{1,i})$ such that $q_j$ witnesses it.  Now
use the existence of primes (?) for ${\frak s}$ to choose
$b_i,M_{1,i}$ and as there we can choose $f_i$.  The rest should be clear.
\enddemo
\bigskip

\proclaim{\stag{705-5.30} Claim}  1) If $M_0 \le_{\frak s} M_1 \le_{\frak
s} N$ and $\bold J$ is independent in $(M_0,N)$ and $(M_0,M_1,M_1 \cap
\bold J) \in K^{3,\text{vq}}_{\frak s}$, 
\ub{then} $\bold J \backslash (M_1
\cap \bold J)$ is independent in $(M_0,M_1,N)$. \nl
2) [${\frak s}$ is successful.]
If $\bold J$ is independent in $(M_0,N_0)$, \ub{then} we can find
$(M_1,N_1)$ such that:
\mr
\item "{$(a)$}"  $M_0 \le_{\frak s} M_1 \le_{\frak s} N_1,M_0
\le_{\frak s} N_0 \le_{\frak s} N_1$
\sn
\item "{$(b)$}"  $M_1$ is $(\lambda,*)$-brimmed over $M_0$
\sn
\item "{$(c)$}"  $N_1$ is $(\lambda,*)$-brimmed over $N_0$
\sn
\item "{$(d)$}"  $(M_1,N_1,\bold J) \in K^{3,\text{vq}}_{\frak s}$
\sn
\item "{$(e)$}"  {\rm tp}
$(c,M_1,N_1)$ does not fork over $M_0$ for every $c
\in \bold J$ so $\bold J$ is independent in $(M_0,M_1,N_1)$.
\ermn
[Similarly for towers?]
\endproclaim
\bigskip

\demo{Proof}  1) Let $\bold J_1 = M_1 \cap \bold J$ and let it be
$\{a_i:i < \alpha\}$.  As $\bold J_1$ is independent in $(M_0,M_1)$
(by monotonicity \scite{705-5.3A}(3)) by Theorem \scite{705-5.2} we can find a
$\le_{\frak s}$-increasing continuous sequence $\langle M_{1,i}:i \le
\alpha \rangle$ such that $M_0 \le_{\frak s} M_{1,0},M_1 \le_{\frak
s} M_{1,\alpha}$ and $(M_{1,i},M_{1,\alpha},a_i) \in 
K^{3,\text{uq}}_{\frak s}$ 
and tp$(a_i,M_{1,i},M_{1,i+1})$ does not fork over $M_0$.

By the existence of NF-amalgamation \wilog \, for some $N^+$ we have
NF$_{\frak s}(M_1,M_{1,\alpha},N,N^+)$.

Now $\langle M_{1,i}:i \le \alpha \rangle$ witness then $\bold J_1 =
\{a_i:i < \alpha\}$ is independent in $(M_0,M_{1,0},M_{1,\alpha})$ so
as $(M_0,M_1,\bold J_1) \in K^{3,\text{vq}}_{\frak s}$ 
by the assumption of the claim definition we
have NF$_{\frak s}(M_0,M_{1,0},M_1,M_{1,\alpha})$.  By the
transitivity of NF we get NF$_{\frak s}(M_0,M_{1,0},N,N^+)$.  So as
$\bold J$ is independent in $(M_0,N)$ by Claim \scite{705-5.1R} 
we get that
$\bold J$ is independent in $(M_0,M_{1,\alpha},N^+)$. \nl
By \scite{705-5.2}(3) and the properties of $\langle M_{1,i},a_j:i \le
\alpha,j < \alpha \rangle$ we can deduce that $\bold J \backslash
\{a_j:j < \alpha\} = \bold J \backslash \bold J_1 = \bold J \backslash
(M_1 \cap \bold J)$ is independent in $(M_0,M_{1,\alpha},N^+)$ 
hence by
monotonicity in $(M_0,N)$ as required.
\nl
2) We try to choose by induction on $\zeta < \lambda^+_{\frak s}$, a
pair $(M'_\zeta,N'_\zeta)$ such that:
\mr
\item "{$(*)_1$}"  $(M'_0,N'_0) = (M_0,N_0)$
\sn 
\item "{$(*)_2$}"  $M'_\zeta$ is 
$\le_{\frak s}$-increasing continuous,
\sn
\item "{$(*)_3$}"  $N_\zeta$ is $\le_{\frak s}$-increasing continuous
\sn
\item "{$(*)_4$}"  $M'_\zeta 
\le_{\frak s} N'_\zeta,\bold J$ is independent in
$(M_0,M'_\zeta,N'_\zeta)$
\sn
\item "{$(*)_5$}"  $\neg \text{ NF}_{\frak s}(M'_\zeta,N'_\zeta,
M'_{\zeta +1},N'_{\zeta +1})$ for $\zeta$ even
\sn
\item "{$(*)_6$}"  $M'_{\zeta+1}$ is saturated over $M'_\zeta$ and
$N'_{\zeta +1}$ is saturated over $N'_\zeta$ if $\zeta$ is odd.
\ermn
There are no problems in successor stages and for limit stages use
\scite{705-4.11}(2). \nl
We necessarily (by \S1) get stuck for 
some $\zeta$ and $(M'_\zeta,N'_\zeta)$
can serve as $(M_1,N_1)$.  \nl
${{}}$  \hfill$\square_{\scite{705-5.30}}$\margincite{705-5.30}
\enddemo
\newpage

\head {\S6 Orthogonality} \endhead  \resetall  \sectno=6
 \spuriousreset 
\bigskip

Note that presently the case ``orthogonality = weak orthogonality" is
the main one for us.  In the latter part of the section ``${\frak s}$
has primes" is usually used and we shall later weaken this, but this
is not a serious flaw here.
\demo{\stag{705-5a.0} Hypothesis}  1) ${\frak s}$ is a $\lambda$-good$^+$
frame, weakly successful. 
\enddemo
\bigskip

\definition{\stag{705-5a.1} Definition}   1) For $p,q \in 
{\Cal S}^{\text{bs}}_{\frak s}(M)$ we say that they are weakly orthogonal,
$p \underset{\text{wk}} {}\to \bot q$ when: if $(M,N,b) \in 
K^{3,\text{uq}}_\lambda$ and tp$(b,M,N) = q$ \ub{then} $p$ 
has a unique extension in 
${\Cal S}_{\frak s}(N)$ equivalently, every extension of $p$ in 
${\Cal S}_{\frak s}(N)$ does not fork over $M$; note: the order of
$p,q$ is seemingly important.  (In the first order case the symmetry is
essentially by the definition and here it will be proved). \nl
2) For $p,q \in {\Cal S}^{\text{bs}}_{\frak s}(M)$ 
we say that they are strongly
orthogonal, $p \underset{\text{st}} {}\to \bot q$ or $p \bot q$
\ub{if} $p_1,q_1$ are weakly orthogonal
whenever $M \le_{\frak s} M_1$ and 
$p_1,q_1 \in {\Cal S}^{\text{bs}}(M_1)$ is
a nonforking extension of $p,q$ respectively. \nl
3) If $p \in {\Cal S}^{\text{bs}}(M),
q \in {\Cal S}^{\text{bs}}(N),M \le_{\frak s} N$ 
orthogonality of $p$ and $q$ means $p',q$ are orthogonal
where $p' \in {\Cal S}^{\text{bs}}_{\frak s}(N)$ is 
the unique nonforking extension of $p$.  Similarly, if
$p \in {\Cal S}^{\text{bs}}_{\frak s}(N),
q \in {\Cal S}^{\text{bs}}_{\frak s}(M)$.
\enddefinition
\bn
Naturally we now show that the definition is equivalent to some 
variants (e.g. for some such pair $(N,a)$ 
rather than all such $(N,a)$).
\proclaim{\stag{705-5a.2} Claim}  Assume that $p,q \in 
{\Cal S}^{\text{bs}}_{\frak s}(M)$ and
$(M,N,b) \in K^{3,\text{uq}}_{\frak s}$ and $q = { \text{\rm
tp\/}}(b,M,N)$.  \ub{Then} $p \underset{\text{wk}} {}\to
\bot q$ iff $p$ has a unique extension in ${\Cal S}_{\frak s}(N)$.
\endproclaim
\bigskip

\demo{Proof}  The implication $\Rightarrow$ holds by the ``every" in
the definition.  So assume $p$ has a unique extension in ${\Cal
S}^{\text{bs}}(N)$ and we 
shall prove $p \underset{\text{wk}} {}\to \bot q$. \nl
So assume $(M,N_2,b_2) \in K^{3,\text{uq}}_{\frak s}$ 
and tp$_{\frak s}(b_2,M,N_2) = q$
and let $p_2 \in {\Cal S}_{\frak s}(N_2)$ extend $p$.  So there is
$N^+ \in K_{\frak s}$ and $a \in N^+$ such that $N_2 \le_{\frak s}
N^+$ and $a$ realizes $p_2$ in $N^+$.  As ${\frak s}$ is weakly
successful, possibly replacing $N^+$ by a $\le_{\frak s}$-extension,
there is $N_1 \le_{\frak s} N^+$ such that $(M,N_1,a) \in
K^{3,\text{uq}}_{\frak s}$ and \wilog \, $b_2 = b$
(as $N,N^+$ are $\le_{\frak s}$-extensions of $M$ and $b \in N,b_2\in
N_2 \le_{\frak s} N^+$ realizes the same type so we can amalgamate)
also \wilog \, $N \le_{\frak s} N^+$.

Now as $a$ realizes $p_2$ in $N^+$ it also realizes $p_2 \restriction
M$ which is $p$ so tp$_{\frak s}(a,N,N^+)$ is an extension of $p$ in
${\Cal S}_{\frak s}(N)$ hence by our present assumption it does not
fork over $M$.  We can conclude by Claim \scite{705-1.15M} 
that NF$_{\frak s}(M,N,N_1,N^+)$,
but $b \in N$, tp$(b,M,N) \in {\Cal S}^{\text{bs}}_{\frak s}(M)$ 
hence by \scite{705-stg.32} tp$(b,N_1,N^+)$ does not fork over
$M$.  But we have $(N,N_2,b) \in K^{3,\text{uq}}_{\frak s}$ 
hence NF$_{\frak s}(M,N_2,N_1,N^+)$ from which (as $a \in N_1$, 
tp$(a,M,N^+) \in {\Cal S}^{\text{bs}}(M))$ we 
deduce tp$(a,N_2,N^+)$ does not fork over $M$, but this
last type is $p_2$ so we are done.  \hfill$\square_{\scite{705-5a.2}}$\margincite{705-5a.2}
\enddemo
\bigskip

\proclaim{\stag{705-5a.3Y} Claim}  1) Assume $(M,N,b) 
\in K^{3,\text{pr}}_{\frak s}$
or just $(M,N,b) \in K^{3,\text{uq}}_{\frak s}$ and $q
= { \text{\rm tp\/}}(b,M,N)$ and 
$p \in {\Cal S}^{\text{bs}}_{\frak s}(M)$. \nl
1) If $p$ is realized in $N$, \ub{then} $p,q$ are not
orthogonal and even not weakly orthogonal. \nl
2) Let $p,q \in {\Cal S}^{\text{bs}}_{\frak s}(M)$.  Then
$p \underset{\text{wk}} {}\to \pm q$ \ub{iff} for 
some $a,N'$ we have $N \le_{\frak s} N',a \in N'$
realizes $p$ and $\{a,b\}$ is not independent in $(M,N')$.
\endproclaim
\bigskip

\demo{Proof}  1) The type $p$ has at least 
two extensions in ${\Cal S}_{\frak s}(N)$: one algebraic, is
tp$_{\frak s}(a,N,N)$ where $a \in N$ realizes $p$ 
and the second is in ${\Cal S}^{\text{bs}}_{\frak s}(N)$, 
hence nonalgebraic, in fact a nonforking extension of $p$.  
So by \scite{705-5a.2} we have 
$p \underset{\text{wk}} {}\to \pm q$. \nl
2) Let $(M,N,a) \in K^{3,\text{uq}}_{\frak s},q = \text{ tp}(a,M,N)$. If $p \underset{\text{wk}} {}\to \pm q$ 
then by \scite{705-5a.2} there is $p_1 \in {\Cal S}_{\frak
s}(N)$ extending $p$ forking over $M$, and let $a,N'$ be such that $N
\le_{\frak s} N'$, and $p_1 = \text{ tp}(a,N,N')$; now by
\scite{705-5.2}(3) we get $\{a,b\}$ is not independent
in $(M,N')$; so first phrase implies the second.  If $p \bot q$ then
$\{a,b\}$ is independent over $M$ inside $N'$ by the Definition \scite{705-5.1}.
 \hfill$\square_{\scite{705-5a.3Y}}$\margincite{705-5a.3Y} 
\enddemo
\bigskip

\definition{\stag{705-5a.4} Definition}  Fixing ${\frak C} \in K^{\frak
s}$, if $p_\ell \in {\Cal S}^{\text{bs}}_{\frak s}
(M_\ell)$ and $M_\ell \le_{{\frak K}[{\frak s}]} {\frak C}$ 
for $\ell = 1,2$, \ub{then} let
$p_1 \| p_2$, in words $p_1,p_2$ are parallel inside ${\frak C}$, 
mean that for some $M,p$ we have $M_1
\cup M_2 \subseteq M <_{{\frak K}
[{\frak s}]} {\frak C}$ and $M_2 \in K_{\frak s},
p \in {\Cal S}^{\text{bs}}(M)$ 
does not fork over $M_\ell$ and extend $p_\ell$
for $\ell =1,2$. 
\enddefinition
\bigskip

\remark{Remark}  If $M_\ell \le_{{\frak K}[{\frak s}]} {\frak C}$ and
$p_\ell \in {\Cal S}^{\text{bs}}_{\frak s}
(M_\ell)$ for $\ell = 1,2$ we can
define when $p_1 \bot p_2$ and prove the natural properties.
Similarly for $p \bot M$ defined in \scite{705-5.7} below.
\endremark
\bn
Obvious properties of parallelism are
\proclaim{\stag{705-5a.5Y} Claim}  1) Parallelism 
inside ${\frak C} \in K^{\frak s}$ 
is an equivalence relation. \nl
2) If $M \in K_{\frak s}$, \ub{then} 
on ${\Cal S}^{\text{bs}}_{\frak s}(M)$, parallelism is equality. \nl
\endproclaim
\bigskip

\demo{Proof}  Easy.
\enddemo
\bigskip

\proclaim{\stag{705-5.6} Claim}  1) If 
$p,q \in {\Cal S}^{\text{bs}}_{\frak s}(M)$ and
$f$ as an isomorphism from $M$ onto $N$ \ub{then} 
$p \underset{\text{wk}} {}\to \bot q 
\Leftrightarrow f(p) \underset{\text{wk}} {}\to \bot f(q)$.  
Similarly for $\bot$. \nl
2) If $p,q \in {\Cal S}^{\text{bs}}_{\frak s}(M)$ 
\ub{then} $p \underset{\text{wk}} {}\to \bot 
q \Leftrightarrow q \underset{\text{wk}} {}\to \bot p$.  
Similarly for $\bot$.
\nl
3) Assume that $M,N \in K_{\frak s}$ are brimmed (e.g. $K_{\frak s}$
categorical).  If $M \le_{\frak s} N$, and $p,q \in 
{\Cal S}^{\text{bs}}_{\frak s}(N)$ do not
fork over $M$, \ub{then} $p \underset{\text{wk}} {}\to \bot q 
\Leftrightarrow (p \restriction M) \bot_{\text{wk}} 
(q \restriction M)$. \nl
4) Assume $M,N \in K_{\frak s}$ are brimmed.
If $p_1,p_2 \in {\Cal S}^{\text{bs}}_{\frak s}(M)$ and $q_1,q_2 \in 
{\Cal S}^{\text{bs}}_{\frak s}(N)$ and $M \le_{{\frak K}[{\frak s}]}
{\frak C},N \le_{{\frak K}[{\frak s}]} {\frak C}$
and $p_1\|q_1,p_2\|q_2$ inside ${\frak C}$
\ub{then} $p_1 \underset{\text{wk}} {}\to \bot 
p_2 \Leftrightarrow q_1 \underset{\text{wk}} {}\to \bot q_2$. \nl
Hence $p_1 \bot p_2 \Leftrightarrow q_1 \bot q_2 \Leftrightarrow p_1
\underset{\text{wk}} {}\to \bot p_2$. \nl
5) If $\langle M_i,a_j:i \le \alpha,j < \alpha \rangle$ is an
$M_0$-based {\rm pr}-decomposition or just {\rm uq}-decomposition of 
$(M_0,M_\alpha)$ and $p \in {\Cal S}^{\text{bs}}_{\frak s}(M_0)$ 
is weakly orthogonal to {\rm tp}$_{\frak s}(a_j,M_j,M_{j+1})$ for
every $j < \alpha$ \ub{then} $p$ has a unique 
extension in ${\Cal S}_{\frak s}(M_\alpha)$. \nl
6) Assume $M_0 \le_{\frak s} M_1$ and $p,q \in 
{\Cal S}^{\text{bs}}(M_1)$
does not fork over $M_0$.  If $p \underset{\text{wk}} {}\to \bot q$ 
\ub{then} $(p \restriction M_0) \underset{\text{wk}} {}\to \bot 
(q \restriction M_0)$.
\nl
7) If $M <_{\frak s} N,N$ is universal over $M,p,q \in {\Cal
S}^{\text{bs}}(N)$ does 
not fork over $M$ and $p \underset{\text{wk}} {}\to \bot q$
\ub{then} $p \bot q$ (hence $(p \restriction M) \bot q$). \nl
8) If $\langle M_i:i \le \delta \rangle$ is 
$\le_{\frak s}$-increasing
continuous; $\delta < \lambda^+_{\frak s}$ and $p,q \in {\Cal
S}^{\text{bs}}
(M_\delta)$ \ub{then} $p \underset{\text{wk}} {}\to \bot q$ iff for
every large enough $i < \delta,(p \restriction M_i) 
\underset{\text{wk}} {}\to \bot (p \restriction M_j)$.
\endproclaim
\bigskip

\demo{Proof}  1) Immediate. \nl
2) By \scite{705-zm.1} (alternatively follows 
from \scite{705-5a.3Y}(2) + \scite{705-5.2}). \nl
3) By \scite{705-stg.9} there is an isomorphism $f$ from $M$ onto $N$
such that $f(p \restriction M) = p,f(g \restriction M)=q$ so the
results holds by part (1). 
\nl
4) The first phrase follows by (3) and the definitions (using a third
model $(\le_{{\frak K}[{\frak s}]},{\frak C})$ extending $M$ and
$N$).  The second phrase follows by the first and Definition \scite{705-5a.1}(2).
\nl
5) Let $M_\alpha \le_{\frak s} N$ and $c \in N$ realizes $p$.  We
prove by induction on $\beta \le \alpha$ then tp$(c,M_\beta,N) \in
{\Cal S}^{\text{bs}}_{\frak s}(M_\beta)$ does not fork over $M_0$. \nl
6) Assume toward contradiction 
$(p \restriction M_0) \underset{\text{wk}}
{}\to \pm (q \restriction M_0)$.  So there are 
$M_0 \le_{\frak s} N_1 \le_{\frak s} N_2,b \in N_1$
realizes $q \restriction M_0,(M_0,N_1,b) \in K^{3,\text{uq}}_{\frak
s}$ and $a \in N_2$ realizes $p \restriction M_0$ but
tp$(a,N_1,N_2)$ forks over $N_0$.  By our knowledge on NF \wilog \,
for some $N_3$ we have NF$_s(M_0,M_1,N_2,N_3)$ hence $a,b$ realizes
$p,q$ in $N_3$ respectively (see \scite{705-stg.32}).  
By \scite{705-5.3A}(4), if $\{a,b\}$ is independent in $(M_1,N_3)$ then
$\{a,b\}$ is independent in $(M_0,N_3)$.
So $N^+_1,N_3,a,b$ witness $p \underset{\text{wk}} {}\to \bot q$ by
\scite{705-5a.3Y}(2). \nl
7) Follows by part (6). \nl
8) Easy.  \hfill$\square_{\scite{705-5.6}}$\margincite{705-5.6}
\enddemo
\bigskip

\proclaim{\stag{705-5.6A} Claim}  1) If $p,q \in 
{\Cal S}^{\text{bs}}_{\frak s}(M_1)$ does not fork over $M_0$ where
$M_0 \le_{\frak s} M_1$ \ub{then} $(p \bot q) \Leftrightarrow (p
\restriction M_0) \bot (q \restriction M_0)$. \nl
2) If $p \bot q$ \ub{then} $p \underset{\text{wk}} {}\to \bot q$. \nl
3) Assume that $\langle M_\alpha:\alpha \le \delta
\rangle$ is $\le_{\frak s}$-increasing continuous, $\delta <
\lambda^+$ limit ordinal and $p,q \in {\Cal S}^{\text{bs}}
(M_\delta)$.  \ub{Then}
$p \underset{\text{wk}} {}\to \bot q$ 
iff for every $\alpha < \delta$ large
enough $(p \restriction M_\alpha) \underset{\text{wk}} {}\to \bot (p
\restriction M_\alpha)$.  Similarly $\bot$.
\nl
4) If $p_1,p_2 \in {\Cal S}^{\text{bs}}_{\frak s}(M)$ 
and $q_1,q_2 \in {\Cal S}^{\text{bs}}_{\frak s}(N)$ and 
$M \le_{{\frak K}[{\frak s}]} {\frak C},N
\le_{{\frak K}[{\frak s}]} {\frak C}$ and $p_1 \| q_1,p_2 \| q_2$
inside ${\frak C}$ then $p_1 \bot p_2 \Leftrightarrow q_1 \bot
q_2$. \nl
5) If $K_{\frak s}$ is categorical \ub{then} 
$\bot,\underset{\text{wk}} {}\to \bot$ are equal. \nl
6) If $M \in K_{\frak s}$ is brimmed and $p,q \in {\Cal
S}^{\text{bs}}_{\frak s}(M)$ \ub{then} 
$p \underset{\text{wk}} {}\to \bot q \Leftrightarrow p \bot q$. 
\endproclaim
\bigskip

\demo{Proof}  1) The implication $\Leftarrow$ is by the definition.
For the other direction assume $p \bot q$ and $M_0 \le_{\frak s} M_2$
and $p_2,q_2 \in {\Cal S}^{\text{bs}}(M_2)$ are 
nonforking extensions of $p \restriction M_0,q
\restriction M_0$ respectively.  Without loss of 
generality for some $M_3$ we have 
$M_2 \le_{\frak s} M_3,M_1 \le_{\frak s} M_3$ and let $p_3,q_3 \in
{\Cal S}^{\text{bs}}(M_3)$ be 
nonforking extensions of $p_2,q_2$ respectively
hence of $p \restriction M_0,q \restriction M_0$ respectively.  As $p
\bot q$ we have $p_3 \underset{\text{wk}} {}\to \bot q_3$ and by
\scite{705-5.6}(6) also $p_3 \restriction M_2) \underset{\text{wk}} {}\to
\bot (q_3 \restriction M_2)$ which means
$p_2 \underset{\text{wk}} {}\to \bot q_2$, as
required. \nl
2) Read the definitions. \nl
3) By \scite{705-5a.2} and \scite{705-4.11}. For the last phrase use the proof
of part (1).  \nl
4) By part (1).  \nl
5) By (6). \nl
6) First assume $\neg(p \underset{\text{wk}} {}\to \bot q)$ 
then by the definitions $\neg(p \bot q)$,
the counterexample is $M$ itself. \nl
Second, assume $p \underset{\text{wk}} {}\to \bot q$ and 
let $N_1$ be such
that $M \le_{\frak s} N_1$ and $p_1,q_1 \in 
{\Cal S}^{\text{bs}}(N_1)$ be
nonforking extensions of $p,q$ respectively; we shall prove $p_1
\underset{\text{wk}} {}\to \bot q_1$, 
this suffices for $p \bot q$ hence finish the proof.  Now there are
$N_2,p_2,q_2$ such that $N_1 \le_s N_2,N_2$ is $(\lambda,*)$-brimmed
and $p_2,q_2 \in {\Cal S}^{\text{bs}}(N_2)$ are nonforking extensions of
$p,q$ respectively hence $p_2 \restriction N_1 = p_1,q_2 \restriction
N_1 = q_1$.  By \scite{705-5.6}(4) we have 
$(p \underset{\text{wk}} {}\to \bot q)
\equiv (p_2 \underset{\text{wk}} {}\to \bot q_2)$ so by our present
assumption $p_2 \underset{\text{wk}} {}\to \bot q_2$ 
hence by \scite{705-5.6}(6) we have 
$p_1 \underset{\text{wk}} {}\to \bot q_1$ so we are done.
\hfill$\square_{\scite{705-5.6A}}$\margincite{705-5.6A}
\enddemo
\bigskip

\definition{\stag{705-5.7} Definition}  Assuming 
$M \le_{\frak s} N$ and $p \in
{\Cal S}^{\text{bs}}_{\frak s}(N)$, we 
let $p \bot M$ ($p$ orthogonal to $M$) mean
that: for any $q$, if 
$q \in {\Cal S}^{\text{bs}}_{\frak s}(N)$ does not fork over 
$M$ \ub{then} $p \bot q$ (but see \scite{705-5.8}(1) below).  Similarly
for $\underset{\text{wk}} {}\to \bot$.
\enddefinition
\bigskip

\proclaim{\stag{705-5.8} Claim}  0) Automorphism of any ${\frak C} \in
K^{\frak s}_{\ge \lambda}$ preserves
$p\|q,p \underset{x} {}\to \bot q,p \underset{x} {}\to \bot M$
for $x \in \{${\rm wk,st}$\}$. \nl
1) If $M \le_{\frak s} N_\ell \, (\le_{{\frak K}[{\frak s}]} {\frak C}),
p_\ell \in {\Cal S}^{\text{bs}}_{\frak s}(N_\ell)$ for $\ell = 1,2$ and
$p_1\|p_2$ \ub{then} $p_1 \bot M \Leftrightarrow p_2 \bot M$ (so we can write
$p \bot N$ if for some $p' \in {\Cal S}^{\text{bs}}_{\frak s}(N')$ 
parallel to $p,M \le_{\frak s} N' 
\and p' \underset{x} {}\to \bot N$). \nl
2) If $\langle M_\alpha:\alpha \le \delta \rangle$ is 
$\le_{\frak s}$-increasing continuous, $p \in 
{\Cal S}^{\text{bs}}_{\frak s}(N)$ where $N \le_{{\frak K}[{\frak s}]}
{\frak C},M_\delta \le {\frak C}$ (so $N \in K_{\frak s}$) and
$\alpha < \delta \Rightarrow p \bot M_\alpha$ \ub{then} 
$p \bot M_\delta$. \nl
3) If $M \le_{\frak s} N_\ell \le_{\frak s} N$ and 
$p_\ell \in {\Cal S}^{\text{bs}}_{\frak s}
(N_\ell)$ for $\ell=1,2$ and {\rm NF}$_{\frak s}(M,N_1,N_2,N)$ 
and $p_2 \bot M$ \ub{then}
$p_2 \bot p_1$ (hence $p_2 \bot N_1$). \nl
4) If $p \in {\Cal S}^{\text{bs}}(M_3),
M_0 \le_{\frak s} M_\ell \le_{\frak s}
M_3$ for $\ell = 1,2,p$ does not fork over $M_2$ and
$p \underset{x} {}\to \bot M_1$ \ub{then} $p \restriction M_2
\underset{x} {}\to \bot M_0$ when $x \in \{{\text{\rm st,wk\/}}\}$.
\endproclaim
\bigskip

\demo{Proof}  0)  Trivial.  \nl
1)  Just note that if $q \in {\Cal S}^{\text{bs}}_{\frak s}(M),q_\ell
\in {\Cal S}^{\text{bs}}_{\frak s}(N_\ell)$ is a nonforking extension 
of $q$ for
$\ell = 1,2$ then by \scite{705-5.6A}(4) because by $q_1\|q_2$ we have 
$p_1 \bot q_1 \Leftrightarrow p_2 \bot q_2$ and we are done.
\nl
2) Easy by the local character (i.e., Axiom (E)(c) of good frames) and
\scite{705-5.6A}(3). \nl
3) By part (1) without loss of generality
\mr
\item "{$(*)$}"  $N_2$ is 
$(\lambda,*)$-brimmed over $M$ and $M$ is brimmed.  
\ermn
Also we can find $\langle M_n,N_{2,n}:n < \omega \rangle$ such that:
NF$_{\frak s}(M_n,N_{2,n},M_{n+1},
N_{2,n+1}),M_{n+1}$ is $(\lambda,*)$-brimmed over
$M_n,N_{2,n+1}$ is $(\lambda,*)$-brimmed over 
$M_{n+1} \cup N_{2,n}$ (by NF calculus).  So 
(\marginbf{!!}{\cprefix{600}.\scite{600-nf.17}}) $\dbcu_{n < \omega} N_{2,n}$ is
$(\lambda,*)$-brimmed over $\dbcu_{n < \omega} M_n$ so by
$(*)$ above without loss of generality 
$\dbcu_{n < \omega} N_{2,n} = N_2$ and $\dbcu_{n < \omega} M_n = M$.  
So for some $k < \omega$ the type $p_2$ 
does not fork over $N_{2,k}$.  By the NF$_{\frak s}$ calculus we have
NF$_{\frak s}(M_k,N_{2,k},N_1,N)$.  Recall ${\frak C}$ is
$(\lambda,*)$-brimmed over $N$, so we can find an automorphism $f$ of 
${\frak C}$ such that $f \restriction N_{2,k} = 
\text{ id}_{N_{2,k}},f(N_1) \subseteq M_{k+1} \subseteq M$.
Let $p'_1 = f(p_1) \in {\Cal S}^{\text{bs}}(f_1(N_1))$ and let $p''_1 \in
{\Cal S}^{\text{bs}}(M)$ be a nonforking extension of $p'_1$ as $f_1(N_1)
\le_{\frak s} M$.  Now $p''_1 \underset x {}\to \bot p_2$ as 
$p_2 \bot M$, hence $p'_1
\bot (p_2 \restriction M_k)$ by \scite{705-5.6A}(4).  By part (0) we have
$p_1 \bot (p_2 \restriction M_k)$ and 
lastly $p_1 \bot p_2$ by \scite{705-5.6A}(4). \nl
4) Easy.  \hfill$\square_{\scite{705-5.8}}$\margincite{705-5.8}
\enddemo
\bn
Naturally, we would like to reduce orthogonality for ${\frak s} =
{\frak t}^+$, to orthogonality for ${\frak t}$.
\proclaim{\stag{705-6.xobO} Claim}  Assume ${\frak s} 
= {\frak t}^+,{\frak t}$
a successful $\lambda_{\frak t}$-good$^+$ frame, so $\lambda =
\lambda_{\frak s} = \lambda^+_{\frak t}$.
\nl
Assume further $M_\ell \in K_{\frak s}$ and 
$\langle M^\ell_\alpha:\alpha < \lambda \rangle$ is
a $\le_{\frak t}$-representation of $M$ and for simplicity each
$M^\ell_\alpha$ is brimmed (for ${\frak t}$) and
$M_\ell \le_{\frak s} {\frak C}_{\frak s}$.  \nl
0) For ${\frak s}$ we have $\bot = \underset{\text{wk}} {}\to \bot =
\underset{\text{st}} {}\to \bot$. \nl
1) If $p_1,p_2 \in {\Cal S}^{\text{bs}}_{\frak s}(M_0)$ \ub{then}:
\footnote{here we use ``$K^{3,\text{uq}}_\lambda$ is 
dense" in definition of $\bot$} \nl
$p_1 \bot_{\frak s} p_2$ 
\ub{iff} for unboundedly many $\alpha < \lambda$ we have
$(p_1 \restriction M^0_\alpha) 
{\underset{\text{wk}} {}\to \bot_{\frak t}}
(p_2 \restriction M^0_\alpha)$ \ub{iff} for every large enough $\alpha
< \lambda$, we have $(p_1 \restriction M^0_\alpha) \bot_{\frak t} 
(p_2 \restriction M^0_\alpha)$. \nl
2) If $M_0 \le_{\frak s} M_1 \le_{\frak s} M_2$ 
and $a \in M_2 \backslash M_1$, \ub{then}: 
{\rm tp}$_{\frak s}(a,M_1,M_2) \in 
{\Cal S}^{\text{bs}}_{\frak s}(M_1)$ and
is orthogonal (for ${\frak s}$) to 
$M_0$ \ub{iff} for a club of ordinals $\delta < \lambda$
we have {\rm tp}$_{\frak t}
(a,M^1_\delta,M^2_\delta) \in {\Cal S}^{\text{bs}}_{\frak t}
(M^1_\delta)$ and is orthogonal (for ${\frak t}$) to 
$M^0_\delta$ \ub{iff} for a stationary
set of ordinals $\delta < \lambda$ we 
have {\rm tp}$_{{\frak t}_0}(a,M^1_\delta,
M^2_\delta) \in {\Cal S}^{\text{bs}}_{\frak t}(M^1_\delta)$ and is
$\bot_{\frak t} M^0_\delta$. \nl
3) In part (2), ``for all but boundedly many $\delta < \lambda$", ``for
unboundedly many $\delta < \lambda$" can replace ``club of $\delta <
\lambda$", ``stationarily many $\delta < \lambda$" respectively.
\endproclaim
\bigskip

\demo{Proof}  0) By \scite{705-5.6A}(5) and the definition of ${\frak
t}^+$.  \nl
1) Without loss of generality for every $\alpha < \lambda,
\ell < 2$ we have $M^0_\alpha$ is ${\frak t}$-brimmed and
$p_\ell \restriction M^0_\alpha \in 
{\Cal S}^{bs}_{\frak t}(M^0_\alpha)$ does not fork 
over $M^0_0$ (and so $p_\ell \restriction 
M^0_\alpha$ is a witness for $p_\ell$) 
hence by \scite{705-5.6A}(4) for every
$\alpha < \lambda$ we have $(p_1 \restriction M^0_\alpha)
\bot_{\frak t} (p_2 \restriction M^0_\alpha) 
\Leftrightarrow 
(p_1 \restriction M^0_0) \bot_{\frak t} (p_2 \restriction M^0_0)$ and
by \scite{705-5.6A}(6) + transitivity of equivalence 
$(p_1 \restriction M^0_\alpha) 
{\underset{\text{wk}} {}\to \bot_{\frak t}} (p_2 
\restriction M^\alpha_\beta) \Leftrightarrow (p_1
\restriction M^0_\beta) \bot_{\frak t} (p_2 \restriction M^0_\beta)$ for
$\alpha,\beta < \lambda$.
\enddemo
\bn
\ub{Case 1}:  Assume that 
$(p_1 \restriction M^0_0) \bot_{\frak t} (p_2 \restriction M^0_0)$.

As ${\frak s}$ has primes (by \scite{705-d.8}), we can assume that 
$M_0 \le_{\frak s} M_1 \le_{\frak s} M_2,a_1 \in M_1 \backslash
M_0,p_1 = \text{ tp}_{\frak s}(a_1,M_0,M_1)$ 
and $(M_0,M_1,a_1) \in K^{3,\text{pr}}_{\frak s}$
and $a_2 \in M_2,p_2 = \text{ tp}_{\frak s}
(a_2,M_0,M_2)$ and it suffices to prove that
tp$_{\frak s}(a_2,M_1,M_2)$ is a nonforking extension of $p_2$. 
By \scite{705-5a.2} without loss of generality 
$(M_0,M_1,a)$ is as in \scite{705-d.8}(1),
\scite{705-d.1} i.e., is canonically prime.  So for a club $E$ of
$\lambda$, for every $\delta \in E$ we have:
\mr
\item "{$(*)$}"  $(M^0_\delta,M^1_\delta,a_1) 
\in K^{3,\text{uq}}_{\frak t}$ and $M^1_\delta 
\le_{\frak t} M^2_\delta,a_2 \in M^2_\delta$,
tp$_{\frak t}(a_2,M^0_\delta,M^2_\delta)$ is a nonforking
extension of $p_2 \restriction M^0_0$ and 
tp$_{\frak t}(a_1,M^0_\delta,M^1_\delta)$
is a nonforking extension of $p_1 \restriction M^0_0$.
\ermn
As we are assuming $(p_1 \restriction M^0_0) 
\bot_{\frak t} (p_2 \restriction M^0_0)$ hence $(p_1 \restriction
M^0_\delta) \bot_{\frak t} (p_2 \restriction M^0_\delta)$ so we get by 
$(*)$ and the definition of orthogonality that 
tp$_{\frak t}(a_2,M^1_\delta,M^2_\delta)$ is a 
nonforking extension of tp$_{\frak t}(a_2,M^0_\delta,M^2_\delta)$
hence it does not fork over $M^0_0$.
As this holds for every $\delta \in E$ clearly 
$M^0_0$ witness that tp$_{\frak s}(a_2,M_1,M_2)$
does not fork over $M_0$ as required in this case.
\bn
\ub{Case 2}:  Assume that
$(p_1 \restriction M^0_0) \pm_{\frak t} (p_2 \restriction M^0_0)$.
\nl
We shall prove that $p_1 \pm p_2$.  
Let $M_0 <_{\frak s} M_1,a_1 \in M_1,
(M_0,M_1,a_1) \in K^{3,\text{pr}}_{\frak s}$ (recall that ${\frak s}$ has
primes being ${\frak t}^+$) and 
$p_1 = \text{ tp}_{\frak s}(a_1,M_0,M_1)$
and so as ${\frak t}$ is successful, \wilog \, $\alpha < \beta
< \lambda^+ \Rightarrow {\text{\rm NF\/}}_{\frak s}(M^0_\alpha,
M^1_\alpha,M^0_\beta,
M_{1,\beta})$.  By easy manipulation there is $q_2 \in {\Cal S}_{\frak t}
(M^1_0)$ extending $p_2 \restriction M^0_0$ which fork over $M^0_0$.  
We can choose $M^2_0,a_2$ such that $M^1_0 \le_{\frak t} M^2_0,a_2 \in M^2_0$
and $q_2 = \text{ tp}_{\frak t}(a_2,M^1_0,M^2_0)$.  
Now we can choose inductively $f_\alpha,M^2_\alpha$ such 
that $M^2_\alpha$ is $\le_{\frak t}$-increasing
continuous, $f_\alpha$ is a $\le_{\frak t}$-embedding of $M^1_\alpha$ into
$M^2_\alpha$, increasing continuous with $\alpha,f_0 = \text{ id}_{M_{2,0}}$
and $\alpha = \beta +1 \Rightarrow \text{ NF}_{\frak t}(M^1_\beta,
M^2_\beta,M^1_\alpha,M^2_\alpha)$.  No problem to do it and at the end
\wilog \, $\dbcu_{\alpha < \lambda_{\frak s}}
f_\alpha = \text{ id}_{M_1}$ and let $M_2 =
\cup \{M^2_\alpha:\alpha < \lambda\}$.  Easily $M_1,a_1,M_2,a_2$ exemplifies
$p_1 \pm p_2$; that is by \scite{705-stg.32} for every $\alpha$, 
tp$_{\frak t}(a_2,M^0_\alpha,M^2_\alpha)$ is a nonforking
extension of tp$_{\frak t}(a_2,M^0_0,M^2_0) = p_2 \restriction M^0_0$,
hence tp$_{\frak s}(a_2,M_0,M_2) = p_2$.  We conclude tp$_{\frak
s}(a_2,M_1,M_2)$ extends $p_2$ but is not its nonforking extension in
${\Cal S}_{\frak s}(M_1)$ as required for proving 
$p_1 \pm_{\frak s} p_2$. \nl
2) Without loss of generality $a \in M^2_0$ and for 
$\ell < 2,\alpha < \beta < \lambda$ we have \nl
NF$_{\frak t}(M^\ell_\alpha,M^{\ell +1}_\alpha,M^\ell_\beta,
M^{\ell +1}_\beta)$ hence $M^0_0$ is a witness for $p$.  
Clearly tp$_{\frak s}(a,M_1,M_2) 
\bot_{\frak s} M_0$ iff for every 
$q \in {\Cal S}^{\text{bs}}_{\frak s}(M_0)$ 
we have tp$_{\frak s}(a,M_1,M_2) \bot_{\frak s} q$
iff for each $\alpha < \lambda$ for every $q \in 
{\Cal S}^{\text{bs}}_{\frak s}(M_0)$ which does not 
fork over $M^0_\alpha$ we have tp$_{\frak s}(a,M_1,M_2)
\bot_{\frak s} q$ \ub{iff} for
each $\alpha < \lambda$ for every $q \in {\Cal S}^{\text{bs}}_{\frak t}
(M^0_\alpha)$, for every $\beta \in [\alpha,\lambda)$, the types $p
\restriction M^1_\beta$ is ${\frak t}$-orthogonal to the nonforking
extension of $q$ in ${\Cal S}^{\text{bs}}_{\frak t}
(M^0_\beta)$ \ub{iff} for each $\gamma < \lambda$ we 
have $(p \restriction {\Cal S}^{\text{bs}}_{\frak t}
(M^1_\gamma)) \bot_{\frak t} M^0_\gamma$.  Thus we finish. \nl
3) By monotonicity and \scite{705-1.15}.  \hfill$\square_{\scite{705-6.xobO}}$\margincite{705-6.xobO}
\bigskip

\proclaim{\stag{705-6.4A} Claim}  1) Assume $(M,N,a) \in 
K^{3,\text{uq}}_{\frak s}$. \nl
If $M \cup \{a\} \subseteq N'<_{\frak s} N,
b \in N \backslash N'$ and $q = 
{ \text{\rm tp\/}}_{\frak s}(b,N',N) 
\in {\Cal S}^{\text{bs}}_{\frak s}(N')$ 
\ub{then} $q$ is weakly orthogonal to $M$. \nl
2) [${\frak s}$ has primes.]  Assume $(M,N,a) 
\in K^{3,\text{bs}}_{\frak s}$.
We can find $\langle M_i,a_i:i < \alpha \rangle$ 
for some $\alpha < \lambda^+$, which is a pr-decomposition of $N$ over
$M$ with $a_0 = a$, i.e., such that:
\mr
\item "{$(a)$}"  $a_0 = a$
\sn
\item "{$(b)$}"  $M_0 = M$.
\sn
\item "{$(c)$}"  $M_i \le_{\frak s} N$ is $\le_{\frak s}$-increasing
continuous
\sn
\item "{$(d)$}"  {\rm tp}$_{\frak s}
(a_i,M_i,N) \in {\Cal S}^{\text{bs}}_{\frak s}(M_i)$
\sn
\item "{$(e)$}"  stipulating $M_\alpha = N$ we have: \nl
$(M_i,M_{i+1},a_i) \in K^{3,\text{pr}}_{\frak s}$ for $i < \alpha$.
\ermn
3) [${\frak s}$ has primes].  In part (2) 
if also $(M,N,a) \in K^{3,\text{uq}}_{\frak s}$ 
and $\langle M_i,a_i:i < \alpha 
\rangle$ is as there \ub{then} we can add
\mr 
\item "{$(f)$}"  if 
$i > 0$ then {\rm tp}$_{\frak s}(a_i,M_i,N)$ is weakly orthogonal to $M$.
\ermn
4) If $(M,N,\bold J) \in K^{3,\text{uq}}_{\frak s}$ and $M \cup \bold J
\subseteq N' \le_{\frak s} N$ and $b \in N \backslash N'$ and $q =
{ \text{\rm tp\/}}_{\frak s}(b,N',N) \in {\Cal S}^{\text{bs}}(N')$ 
\ub{then} $q$ is weakly orthogonal to $M$.
\endproclaim
\bigskip

\demo{Proof}  (1)  If $q \underset{\text{wk}} {}\to \pm M$ 
then for some $c,N^+_1,r$ we have $N \le_{\frak s} N^+,
c \in N^+_1,r =: \text{ tp}_{\frak s}(c,N',N^+_1) 
\in {\Cal S}^{\text{bs}}_{\frak s}(N')$ does not fork over
$M$ but $\{b,c\}$ is not independent in $(N',N^+)$ (or $b=c$).
Possibly $\le_{\frak s}$-increasing $N^+$ as tp$(c,N',N^+)$ does not
fork over $M \le_{\frak s} N'$, clearly there is $M'$ such that $M
\cup \{c\} \subseteq M'$ and NF$_{\frak s}(M,M',N',N^+)$.  As $a \in
N'$ and tp$(a,N') \in {\Cal S}^{\text{bs}}(M)$ this implies that
tp$(a,M',N^+) \in {\Cal S}^{\text{bs}}(M')$ does not fork over $M$.  As
$(M,N,a) \in K^{3,\text{uq}}_{\frak s}$ it follows that NF$_{\frak
s}(M,N,M',N^+)$, and this implies that $\{b,c\}$ is independent
in $(N',N^+)$, by \scite{705-5.3A}(2), contradicting the choice of $c$.
\nl
2) This is \scite{705-c.6}(1).
\nl
3) Follows by (part (2) and) part (1). \nl
4) Like part (1).  \hfill$\square_{\scite{705-6.4A}}$\margincite{705-6.4A}
\enddemo
\bigskip

\proclaim{\stag{705-5.11} Claim}   1) [${\frak s}$ has primes].
Assume $(M,N,a) \in K^{3,\text{uq}}_{\frak s}$.
\ub{Then} we can find $\langle M_i,a_i:i < \alpha \rangle$ such that:
\mr
\item "{$(a)-(d)$}"  as in \scite{705-6.4A}(2)
\sn
\item "{$(f)$}"  as in \scite{705-6.4A}(3)
\sn
\item "{$(g)$}"  $\alpha \le \lambda$.
\ermn
2)  If in addition ${\frak s} = {\frak t}^+,{\frak t}$ is a
successful $\lambda_{\frak t}$-good$^+$ frame (so $\lambda = 
\lambda^+_{\frak t})$ \ub{then} we can add
\mr
\item "{$(h)$}"  for each $i < \alpha$, 
for any $<_{\frak t}$-representations
$\langle M^i_\varepsilon:\varepsilon < \lambda^+_{\frak t} \rangle,
\langle M^{i+1}_\varepsilon:\varepsilon < \lambda^+ \rangle$ of 
$M^i,M^{i+1}$ respectively, for a club of ordinals
$\delta < \lambda^+_{\frak t}$ we have $(M^i_\delta,
M^{i+1}_\delta,a_i)\in K^{3,\text{uq}}_{\frak t}$.
\endroster
\endproclaim
\bigskip

\demo{Proof}  1)  Exactly as in \scite{705-c.6}(5), i.e., in the proof of
\scite{705-c.6}(1) use a bookkeeping in order to get clause (g). \nl
2) By \scite{705-d.9}(1).
\hfill$\square_{\scite{705-5.11}}$\margincite{705-5.11}
\enddemo
\bigskip

\proclaim{\stag{705-5.12} Claim}  [${\frak s}$ with primes.]
Assume
\mr
\item "{$(a)$}"  {\rm NF}$_{\frak s}(M_0,M^+_0,M_1,M_3)$
\sn
\item "{$(b)$}"  $\bold J$ is independent in $(M_1,M_3)$
\sn
\item "{$(c)$}"  {\rm tp}$_{\frak s}(c,M_1,M_3)$ is orthogonal to $M_0$ for every $c \in
\bold J$.
\ermn
\ub{Then} we can find $M^+_1,M^+_3$ such that:
\mr
\item "{$(\alpha)$}"  $M_3 \le_{\frak s} M^+_3$
\sn
\item "{$(\beta)$}"  $M_1 \cup M^+_0 \subseteq M^+_1 \le_{\frak s} M^+_3$
\sn
\item "{$(\gamma)$}"  {\rm tp}$_{\frak s}
(c,M^+_1,M^+_3)$ does not fork over $M_1$ for $c \in \bold J$
\sn
\item "{$(\delta)$}"  $\bold J$ is independent \footnote{so $(\gamma)
+ (\delta)$ says that $\bold J$ is independent in $(M_1,M^+_1,M^+_3)$}
 in $(M^+_1,M^+_3)$.
\endroster
\endproclaim
\bn
Before we prove \scite{705-5.12} note that:
\demo{\stag{705-5.13} Conclusion}  [${\frak s}$ with primes.]
If to the assumptions of \scite{705-5.12} we add
\mr
\item "{$(d)$}"  $(M_1,M_2,\bold J) \in K^{3,\text{qr}}_{\frak s}$ or
just $(M_1,M_2,\bold J) \in K^{3,\text{vq}}_{\frak s}$ and $M_2
\le_{\frak s} M_3$,
\ermn
\ub{then} we can add to the conclusion (in fact follows from it):
\mr
\item "{$(\varepsilon)$}"  NF$_{\frak s}(M_1,M_2,M^+_1,M^+_3)$.
\endroster
\enddemo
\bigskip

\demo{Proof}  By \scite{705-5.12} and \scite{705-4.8}(2) or by the definition.
\enddemo
\bigskip

\demo{Proof of \scite{705-5.12}}  We can find $\langle M^0_i,a_i:i <
\alpha \rangle$ 
which is a decomposition of $M^+_0$ over $M_0$ and
stipulate $M^0_\alpha = M^+_0$.  We can
now choose by induction on $i,(M^1_i,f_i)$ such that 
$M^1_i \in K_{\frak s}$
is $\le_{\frak s}$-increasing continuous, $M^1_0 = M_1,
f_0 = \text{ id}_{M_0},f_i$ is an 
$\le_{\frak s}$-embedding of $M^0_i$ into $M^1_i$, increasing
continuous with $i$ and $(M^1_i,M^1_{i+1},f_{i+1}(a_i)) \in 
K^{3,\text{pr}}_{\frak s}$ 
and tp$_{\frak s}(f_{i+1}(a_i),M^1_i,M^1_{i+1})$ does not
fork over $f_i(M^0_i)$.  There is no problem to do this, (as in stage
$i=j+1$ first choose $p_i = f_i(\text{tp}_{\frak
s}(a_i,M^0_i,M^0_{i+1}))$ and then $M^1_{i+1}$ such that some $b_i \in
M_{i+1}$ realizes $p_i$ and as $(M^0_i,M^0_{i+1},a_i) \in
K^{3,\text{pr}}_{\frak s}$ we can choose a $\le_{\frak
K}$-embedding of $M^0_{i+1}$ into $M^1_{i+1}$ extending $f_i$ and
mapping $a_i$ to $b_i$).  As 
$K^{3,\text{pr}}_{\frak s} \subseteq K^{3,\text{uq}}_{\frak s}$ 
and the definition of $K^{3,\text{uq}}_{\frak s}$ 
easily NF$_{\frak s}(f_i(M^0_i)$, 
$f_{i+1}(M^0_{i+1}),M^1_i,M^1_{i+1})$ 
hence by NF$_{\frak s}$-symmetry NF$_{\frak s}(f_i(M^0_i),
M^1_i,f_{i+1}(M^0_{i+1}),M^1_{i+1})$ for every $i$ hence
by long transitivity 
NF$_{\frak s}(f_0(M_0),M^1_0,f_\alpha(M^0_\alpha),M^1_\alpha)$, and recalling
$f_0 = \text{ id}_{M_0},M^0_0 = M_0,M^1_0 = M_1,M^0_\alpha = M^+_0$
this means NF$_{\frak s}(M_0,M_1,f_\alpha(M^+_0),M^1_\alpha)$.  
But also we assume
NF$_{\frak s}(M_0,M_1,M^+_0,M_3)$,  hence by NF$_{\frak s}$ uniqueness
without loss of generality for some $M^+_3,M_3
\le_{\frak s} M^+_3,f_i = \text{ id}_{M^0_i}$ and $M^1_\alpha
\le_{\frak s} M^+_3$.  

For $i < \alpha$ for each $c \in \bold J$, note that
tp$_{\frak s}(c,M_1,M_3)$ is orthogonal to $M_0$ (by a hypothesis).
We prove by induction on $i \le \alpha$ that $\bold J$ is independent
over $(M_1,M^1_i)$ inside $M^+_3$ and for every $c \in \bold J$,
tp$_{\frak s}(c,M^1_i,M^+_3)$ (does not fork over $M^1_0 = M_1$ and) is
orthogonal to $M^0_i$.  For $i=0$ this is trivial for does not fork
given for orthogonal.  For $i$ limit easy.  For $i+1$, as
tp$_{\frak s}(a_i,M^1_i,M^+_3)$ does not
fork over $M^0_i$, it is orthogonal to tp$_{\frak s}(c,M^1_i,M^+_3)$ 
for $c \in \bold J$ hence $\bold J \cup \{a_i\}$ is independent over
$M^1_i$.  As $(M^1_i,M^1_{i+1},a_i) \in K^{3,\text{pr}}_{\frak s}$ we
get tp$_{\frak s}(c,M^1_{i+1},M^+_3)$ does not fork over $M^1_i$ hence
over $M_1$.   Let $M^+_1$ be chosen as $M^1_\alpha$. \nl
By \scite{705-5.13B} below we are done.  \hfill$\square_{\scite{705-5.12}}$\margincite{705-5.12}
\enddemo
\bigskip

\remark{\stag{705-5.13A} Remark}  We can phrase the proof of $(\alpha) +
(\beta)$ as a subclaim.
\endremark

\proclaim{\stag{705-5.13B} Claim}  [${\frak s}$ has primes.]  Assume
\mr
\item "{$(a)$}"  $M_0 \le_{\frak s} M_\ell \le_{\frak s} M_3$ for
$\ell = 1,2$
\sn
\item "{$(b)$}"  $\langle M_{0,i},a_i:i < \alpha \rangle$ is a
decomposition of $M_2$ over $M_0$
\sn
\item "{$(c)$}"  $\bold J$ is independent in $(M_0,M_1)$
\sn
\item "{$(d)$}"  {\rm tp}$(c,M_0,M_1) 
\bot { \text{\rm tp\/}}(a_i,M_{0,i},M_2)$ for
$i < \alpha$ and $c \in \bold J$.
\ermn
\ub{Then} $\bold J$ is independent in $(M_2,M_3)$ moreovr in
$(M_0,M_2,M_3)$.
\endproclaim
\bigskip

\demo{Proof}  We prove this by induction on $|\bold J|$.

By \scite{705-5.2} \wilog \, $\bold J$ is finite and let $\bold J$ be
$\{b_\ell:\ell < n\}$ where $n = |\bold J|$.  If $n=0$ this is trivial
and let $\langle M_{\ell,0}:\ell \le n \rangle$ be such that $M_{0,0}
= M_0,M_{\ell,0} \le_{\frak s} M_1$ and $(M_{\ell,0},M_{\ell
+1},b_\ell) \in K^{3,\text{pr}}_{\frak s}$, so
tp$(b_\ell,M_{0,\ell},M_{0,\ell+1})$ does not fork over $M_0$ hence is
orthogonal to tp$(a_i,M_{0,i},M_{0,i+1})$ for $i < \alpha$.

We choose by induction on $\ell \le n$, a sequence $\bar M_\ell =
\langle M_{\ell,i}:i \le \alpha \rangle$ and $N^*_\ell$ such that:

$N^*_0 = M_3$

$\bar M_\ell$ is $\le_{\frak s}$-increasing continuous

$M_{\ell,0}$ as above

$M_{0,i}$ as above hence $M_{0,\alpha} = M_2$

$(M_{\ell,i},M_{\ell,i+1},a_i) \in K^{3,pr}_{\frak s}$

tp$(a_i,M_{\ell,i},M_{\ell,i+1})$ does not fork over $M_{0,i}$

For $\ell=0$ this is clear.  For $\ell+1$,
tp$(b_\ell,M_{\ell,0},M_{\ell+1,0}) = \text{
tp}(b_\ell,M_{\ell,0},N^*_\ell)$ is orthogonal to each
tp$(a_i,M_{0,i},M_{0,i+1}) = \text{ tp}(a_i,M_{0,i},N^*_\ell)$ hence
to tp$(a_i,M_{\ell,i},N^*_\ell)$ so by \scite{705-5.6}(5) we deduce that
tp$(b_\ell,M_{\ell,\alpha},N^*_\ell)$ does not fork over $M_{\ell,0}$
hence NF$_{\frak s}(M_{0,\ell},M_{0,\ell +1},M_{\ell,\alpha},
N^*_\ell)$ so by \scite{705-5.3A}(5) we can choose
$\bar M_{\ell +1},N^*_{\ell +1}$.  \hfill$\square_{\scite{705-5.13B}}$\margincite{705-5.13B}
\enddemo
\bn
Below the restriction $\gamma \le \omega$ may seem quite undesirable
but it will be used as a stepping stone for better things.  Note that
in the proof of \scite{705-5.14}(1), clause $(c)$ in the 
induction hypothesis
on $\langle M^n_i:i \le \alpha \rangle$, primeness, is not proved to
hold for $\langle M^\omega_i:i \le \alpha \rangle$, though enough is
proved to finish the proof, this is why
the proof does not naturally work for
$\gamma > \omega$.
\proclaim{\stag{705-5.14} Claim}  [${\frak s}$ has primes.]
\nl
1) Assume
\mr
\item "{$(a)$}"  $\langle M_\beta:\beta \le \gamma \rangle$ is
$\le_{\frak s}$-increasing continuous with $\gamma \le \omega$
\sn
\item "{$(b)$}"  $M_0 \le_{\frak s} M^+_0 \le_{\frak s} M^*_3$ and
$M_\gamma \le_{\frak s} M^*_3$
\sn
\item "{$(c)$}"  {\rm NF}$_{\frak s}(M_0,M_1,M^+_0,M^*_3)$
\sn
\item "{$(d)$}"  if $0 < \beta < \gamma$ then $(M_\beta,M_{\beta +1},
\bold J_\beta) \in K^{3,\text{qr}}_{\frak s}$ 
(so $\bold J_\beta$ is independent
in $(M_\beta,M_{\beta +1})$
\sn
\item "{$(e)$}"  for every 
$\beta \in (0,\gamma)$ and $a \in \bold J_\beta$
the type {\rm tp}$_{\frak s}
(a,M_\beta,M_{\beta +1})$ is orthogonal to $M_0$.
\ermn
\ub{Then} {\rm NF}$_{\frak s}(M_0,M_\gamma,M^+_0,M^*_3)$. \nl
2) If $\langle M_\beta:
\beta \le \gamma \rangle,\langle \bold J_\beta:0 <
\beta < \gamma \rangle$ satisfy clauses (a), (d), 
(e) above and $(M_0,M_1,
\bold J) \in K^{3,\text{vq}}_{\frak s}$ 
\ub{then} $(M_0,M_\gamma,\bold J) \in K^{3,\text{vq}}_{\frak s}$.
\endproclaim
\bigskip

\demo{Proof}  1) We choose $\langle M^0_i,a_j:i \le \alpha,j < \alpha 
\rangle$, a decomposition of $M^+_0$ over $M_0$ (as in the proof of
\scite{705-5.12}).  Now by induction on $n \le \gamma,n < \omega$ we choose
$N^3_n,\bar M^n = \langle M^n_i:i \le \alpha \rangle$ such that:
\mr
\item "{$(a)$}"  $N^3_0 = M^*_3$ and $\bar M^0 = \langle M^0_i:i \le
\alpha \rangle$
\sn
\item "{$(b)$}"  $N^3_n \le_{\frak s} N^3_{n+1}$
\sn
\item "{$(c)$}"  $\langle 
M^n_i,a_j:i \le \alpha,j < \alpha \rangle$ is a
decomposition inside $N^3_n$ over $M_n$ i.e., $M^n_0 = M_n$
\sn
\item "{$(d)$}"  $M^n_i \le_{\frak s} M^{n+1}_i$
\sn
\item "{$(e)$}"  tp$_{\frak s}(a_j,M^n_i,N^3_n)$ does not fork over $M^0_i$. 
\ermn
For $n=0$ this is done.  The step from $n$ to $n+1$ is by the first
paragraph of the proof of \scite{705-5.12} or \scite{705-5.3A}(5), \ub{but}
for this we need to know that NF$_{\frak
s}(M_n,M_{n+1},M^n_\alpha,N^*_n)$.  First if $n=0$ this holds by clause
(c) of the assumption as $M^0_\alpha = M^+_0$.  Second if 
$n>0$ then holds by clause $(\varepsilon)$ of \scite{705-5.13} 
using \scite{705-5.8}(1).  
If $\gamma < \omega$ we are done.  So assume $\gamma =
\omega$, and let for $i < \alpha,M^\omega_i =: 
\dbcu_{n < \omega} M^n_i$.  Now for each
$i < \alpha$, and $n < \omega$ we have (see the proof of \scite{705-5.12}),
NF$_{\frak s}(M^n_i,M^n_{i+1},M^{n+1}_i,M^{n+1}_{i+1})$,  hence by long
transitivity of NF$_{\frak s}$ 
(see \marginbf{!!}{\cprefix{600}.\scite{600-nf.16}}) we have NF$_{\frak s}(M^0_i,
M^0_{i+1},M^\omega_i,M^\omega_{i+1})$.  By symmetry we get NF$_{\frak s}
(M^0_i,M^\omega_i,M^0_{i+1},M^\omega_{i+1})$ for $i < \omega$.  As
$\langle M^0_i:i \le \alpha \rangle,\langle M^\omega_i:i \le \alpha \rangle$
are $\le_{\frak s}$-increasing continuous, by long transitivity of
NF$_{\frak s}$ we get
NF$_{\frak s}(M^0_0,M^\omega_0,M^0_\alpha,M^\omega_\alpha)$ which means
NF$_{\frak s}(M_0,M_\omega,M^+_0,M^\omega_\alpha)$ so by using monotonicity
twice we get \nl
NF$_{\frak s}(M_0,M_\omega,M^+_0,M^*_3)$ as required. \nl
2) By definition \scite{705-4.9}, Claim \scite{705-4.8}(2) and the first part
of the claim.  \hfill$\square_{\scite{705-5.14}}$\margincite{705-5.14}
\enddemo
\bn
We could have noted earlier
\proclaim{\stag{705-5.22} Claim}  Assume $p_i = { \text{\rm tp\/}}
(a_i,M,N) \in {\Cal S}^{\text{bs}}_{\frak s}(M)$ for 
$i < \alpha$ are pairwise orthogonal.  \ub{Then} 
$\{a_i:i < \alpha\}$ is independent in $(M,N)$.
\endproclaim
\bigskip

\demo{Proof}  By \scite{705-5.2} and renaming it is enough to deal with finite
$\alpha$, (not really used). \nl
We now choose a pair $(M_\ell,N_\ell)$ by induction on $\ell \le
\alpha$ such that
\mr
\widestnumber\item{$(*)(iii)$}
\item "{$\circledast(i)$}"  $M_\ell \le_{\frak s} N_\ell$
\sn
\item "{$(ii)$}"  $M_0 = M,N_0 = N$
\sn
\item "{$(iii)$}"  if $m < \ell$ then $M_m \le_{\frak s} M_\ell$ and
$N_m \le_{\frak s} N_\ell$
\sn
\item "{$(iv)$}"   if $\ell = m+1$ then $(M_m,M_\ell,a_m) \in 
K^{3,\text{uq}}_{\frak s}$
\sn
\item "{$(v)$}"   tp$(a_m,M_m,M_{m+1})$ does not fork over $M_0$.
\ermn
For $\ell=0$ this is trivial.  For $\ell = m+1$, first we prove by induction
on $k \le m$ that $p^k_m = \text{ tp}(a_m,M_k,N_m)$ is the nonforking
extension of $p_m$ in ${\Cal S}(M_k)$, now for $k=0$ this is trivial
by the choice of $p_m$ and for $k+1 \le m$ we use the assumption $p_k
\bot p_m$ and $(M_k,M_{k+1},a_k) \in K^{3,\text{uq}}_{\frak s}$ noting
that by the induction hypothesis on $k,\text{tp}(a_m,M_k,N_m)$ is a non-forking
extension of $p_n$ and by clause (v) for $k,\text{tp}_{\frak
s}(a_k,M_k,M_n)$ is a non-forking extension of $p_k$.
\nl
Second, as ${\frak s}$ is weakly successful there are $b_m,M^*_\ell$ such
that $(M_m,M^*_\ell,b_m) \in K^{3,\text{uq}}_{\frak s}$ and tp$_{\frak
s}(b_m,M_m,M^*_\ell) = p_m$.  By the definition of types and as
$K_{\frak s}$ has amalgamation by renaming there is $N_\ell$ such that
$M_\ell \le_{\frak s} N_\ell,N_m \le_{\frak s} N_\ell$ and $b_m =
a_m$.  So we can define $(M_\ell,N_\ell)$ for $\ell \le n$ as in
$\circledast$.   By the definition of independents we are done.
\nl
${{}}$  \hfill$\square_{\scite{705-5.22}}$\margincite{705-5.22}
\enddemo
\bigskip

\proclaim{\stag{705-5.27} Claim}  If $p_i \in 
{\Cal S}^{\text{bs}}(M)$ for $i < \alpha$
are pairwise orthogonal and $q \pm p_i$ for $i < \alpha$ \ub{then}
$\alpha < \omega$.
\endproclaim
\bigskip

\demo{Proof}  By \scite{705-4.25} and \scite{705-5.22}.  That is assume $\alpha
\ge \omega$, let $q = \text{ tp}_{\frak s}(b,M,N_0)$ and we can find
$N_n(n < \omega) \le_{\frak K}$-increasing and $a_{n+1} \in N_{n+1}$
reazling $p_n$ such that $\{b,a_n\}$ is not independent.  By
\scite{705-5.22}, $\{a_n:n < \omega\}$ is independent in $N_\omega =
\cup\{N_n:n <\omega\}$ and so by \scite{705-4.25}, we get a contradiction.
\hfill$\square_{\scite{705-5.27}}$\margincite{705-5.27}
\enddemo
\bigskip

\proclaim{\stag{705-5.28Y} Claim}  [${\frak s}$ has primes].
Assume that $(M_0,M_1,\bold J) \in K^{3,\text{vq}}_{\frak s}$ 
and $(M_1,M_2,a) \in K^{3,\text{pr}}_{\frak s}$ or just
$\in K^{3,\text{uq}}_{\frak s}$ and {\rm tp}
$(a,M_1,M_2)$ is orthogonal to $M_0$. \nl
\ub{Then} $(M_0,M_2,\bold J) \in K^{3,\text{vq}}_{\frak s}$.
\endproclaim
\bigskip

\demo{Proof}  Assume $M_0 \le_{\frak s} N_0 \le_{\frak s} N_2,M_2
\le_{\frak s} N_2$ such that $\bold J$ is independent in
$(M_0,N_0,N_2)$ and we should prove that NF$_{\frak
s}(M_0,N_0,M_2,N_2)$, this suffices.

We can find a decomposition $\langle M_{0,i},a_i:i < \alpha \rangle$
of $(M_0,N_0)$.  By \scite{705-5.13B} we can find $N^+_2$ and an
$\le_{\frak s}$-increasing continuous $\langle M_{1,i}:i \le \alpha
\rangle$ such that $N_2 \le_{\frak s} N^+_2,M_{0,i} \le_{\frak s}
M_{1,i}$ stipulating $M_{0,\alpha} = N_0$ and $M_{1,0} = M_1$,
tp$(a_i,M_{1,i},N^+)$ does not fork over $M_{0,i}$ and
$(M_{1,i},M_{1,i+1},a_i) \in K^{3,\text{pr}}_{\frak s}$.  Now we prove by
induction on $i \le \alpha$ that tp$(b,N_{1,i},N^+_2)$ does not fork
over $M_1 = M_{1,0}$.  As usual NF$(M_0,M_1,M_{0,i},M_{1,i})$ for $i
\le \alpha$.  For $i=0$ trivial for $i$ limit by Axiom (E)(h) and for
$i=j+1$ just note that by \scite{705-5.8} tp$(b,M_{1,0},N^+_2) = \text{
tp}(b,M_1,M_2)$ is orthogonal to $M_{0,i}$ hence
tp$(a_i,M_{0,i},M_{0,i+1})$ hence to tp$(a_i,M_{1,i},M_{1,i+1})$.  So
tp$(b,M_{1,\alpha},N^+_2)$ does not fork over $M_1 = M_{0,\alpha}$ and
$(M_1,M_2,a) \in K^{3,\text{uq}}_{\frak s}$ so NF$_{\frak
s}(M_1,M_{1,\alpha},M_2,N^+_2)$.  By transitivity of NF we have
NF$_{\frak s}(M_0,N_0,M_2,N^+_3)$ hence NF$_{\frak
s}(M_0,N_0,M_2,N_3)$ is as required. \nl
${{}}$    \hfill$\square_{\scite{705-5.28Y}}$\margincite{705-5.28Y}
\enddemo
\newpage

\head {\S7 Understanding $K^{3,\text{vq}}_{\frak s}$} \endhead  \resetall  \sectno=7
 \spuriousreset 
\bn
We would like to show that $K^{3,\text{vq}}_{\frak s} = 
K^{3,\text{qr}}_{\frak s}$ and $K^{3,\text{pr}}_{\frak s} = 
K^{3,\text{uq}}_{\frak s}$ and more
remembering that every $M \in K_{\frak s}$ is saturated. \nl
The hypothesis below holds if ${\frak t}$ is successful, ${\frak s} =
{\frak t}$ \footnote{where do we use successful rather than weakly
successful?  E.g. in \scite{705-4.23}(3).  This can be somewhat weakened:
replacing a club of $\lambda^+_{\frak s}$ a member of a normal filter
on $\lambda^+_{\frak s}$} is successful.
\demo{\stag{705-5.14X} Hypothesis}  
\mr
\item "{$(a)$}"  ${\frak s}$ is a $\lambda$-good frame
\sn
\item "{$(b)$}"   ${\frak s}$ is successful
\sn
\item "{$(c)$}"  ${\frak s}$ has primes
\sn
\item "{$(d)$}"  $\bot = \underset{\text{wk}} {}\to \bot$.
\endroster
\enddemo
\bn
In the definition below note that our aim is to analyze $(M,N,\bold
J_0)$ so $\bold J_0$ has a special role.

\definition{\stag{705-4.22} Definition}  1) ${\Cal W}_\alpha = \{(N,\bar M,
\bar{\bold J}):\bar M = \langle M_i:i < \alpha \rangle$ is 
$\le_{\frak s}$-increasing continuous, $M_i \le_{\frak s} N$ and
$\bar{\bold J} = \langle \bold J_i:i < \alpha \rangle$ and $\bold J_i$
is independent in $(M_i,M_{i+1})$ stipulating $M_\alpha =N$ 
and we let

$$
{\Cal W} = \dbcu_{\alpha < \lambda^+} {\Cal W}_\alpha
$$
\mn
2) $\le_{\Cal W} = \le_{{\Cal W}[{\frak s}]}$ 
is the following two place relation on ${\Cal W}$:

$$
\align
(N^1,\bar M^1,\bar{\bold J}^1) \le_{\Cal W} &(N^2,\bar M^2,\bar{\bold J}^2)
\text{ \ub{iff} } (a) + (b) \text{ where} \\
  &(a) \quad N^1 \le_{\frak s} N^2,\ell g(\bar M^1) \le
\ell g(\bar M^2),i < \ell g(\bar M^1) \Rightarrow M^1_i \le_{\frak s} M^2_i 
  \and \bold J^1_i \subseteq \bold J^2_i \text{ and} \\
  &(b) \quad a \in \bold J^1_i \Rightarrow \text{ tp}_{\frak s}
(a,M^2_i,M^2_{i+1}) \text{ does not fork over } M^1_i
\endalign
$$
\mn
3) $\le^{\text{fx}}_{\Cal W}$ is defined like $\le_{\Cal W}$ but also
$\bar{\bold J}^1 = 
\bar{\bold J}^2$ (so $\ell g(\bar M^1) = \ell g(\bar M^2)$ in particular).
\enddefinition
\bigskip

\proclaim{\stag{705-4.23} Claim}  1) $\le_{\Cal W}$ is a partial order. \nl
2) If $\delta < \lambda^+_{\frak s}$ is a limit ordinal and $\langle (N^\alpha,
\bar M^\alpha,\bar{\bold J}^\alpha):\alpha < \delta \rangle$ is
$\le_{\Cal W}$-increasing, \ub{then} this sequence has a 
$\le_{\Cal W}$-{\rm lub}
$(N,\bar M,\bold J)$, with $\ell g(\bar M) = \sup\{\ell g(\bar M
^\alpha):\alpha < \delta\},N = \cup \{N^\alpha:\alpha < \delta\},M_i =
\cup\{M^\alpha_i:\alpha \text{ satisfies that } 
i < \ell g(\bar M^\alpha)$ and
$\alpha < \delta\},\bold J^*_i = \cup
\{\bold J^\alpha_i:\alpha \text{ satisfies that }i < \ell g
(\bar M^\alpha)$ and $\alpha < \delta\}$. \nl
3) If $(N^1,\bar M^1,\bold{\bar J}^1) \in {\Cal W}_\alpha$ \ub{then} for some
$(N^2,\bar M^2,\bar{\bold J}^2)$ we have
\mr
\item "{$(\alpha)$}"  $(N^1,\bar M^1,\bar{\bold J}^1) 
\le^{\text{fx}}_{\Cal W}
(N^2,\bar M^2,\bar{\bold J}^2)$ 
\sn
\item "{$(\beta)$}"    $(M^2_i,M^2_{i+1},\bold J^2_i) \in 
K^{3,\text{vq}}_{\frak s}$ for each $i < \ell g(\bar M)$
\sn
\item "{$(\gamma)$}"  $N^2 = \cup\{M^2_i:i < \ell g(\bar M^2)\}$.
\endroster
\endproclaim
\bigskip

\demo{Proof}  Straight: part (1) is trivial, part (2) holds by
\scite{705-4.11}(2), and part (3) is proved repeating in the proof of
\scite{705-5.30}(2) but using part (2) here.  
\hfill$\square_{\scite{705-4.23}}$\margincite{705-4.23}
\enddemo
\bn
We are interested in ``nice" such sequences; we define several
variants.
\definition{\stag{705-5.15} Definition}  1) $K^{\text{or}}_{\frak s} 
= \{(N,\bar M,
\bar{\bold J}) \in {\Cal W}_\omega$: if $a \in \bold J_{n+1}$ \ub{then}
tp$_{\frak s}(a,M_{n+1},M_{n+2})$ is 
orthogonal to $M_0\}$, if we omit $N$ we mean
$N = \cup\{M_n:n < \omega\}$. \nl
2) $K^{\text{ar}}_{\frak s} 
= \{(N,\bar M,\bar{\bold J}) \in {\Cal W}_\omega$: if 
$a \in \bold J_{n+1}$ then tp$_{\frak s}(a,M_{n+1},M_{n+2})$ 
is orthogonal to $M_n\}$. 
\nl 
3) $K^{\text{br}}_{\frak s} 
= \{(N,\bar M,\bar{\bold J}) \in K^{\text{or}}_{\frak s}$: if
$b \in \bold J_{n+1}$ then for some $m=m(b) \le n$ we have
tp$(b,M_{n+1},M_{n+2})$ does not fork over $M_{m+1}$ and is 
orthogonal to $M_m\}$. \nl
4) We say that $(N,\bar M,\bar{\bold J})$ is 
$K^{\text{or}}_{\frak s}$-full \ub{if} it belongs
to $K^{\text{or}}_{\frak s}$ and
\footnote{Note that on $\bold J_0$ there are no demands}
$p \underset{\text{wk}} {}\to \bot M_0 
\and p \in {\Cal S}^{\text{bs}}(M_{n+1})
\Rightarrow \lambda_{\frak s} 
= |\{c \in \bold J_{n+1}:p = \text{ tp}_{\frak s}(c,
M_{n+1},M_{n+2})\}|$ and $N = \cup\{M_n:n < \omega\}$.  
\nl
5) We say that $(N,\bar M,\bold J)$ is $K^{\text{ar}}_{\frak s}$-full
if it
$\in K^{\text{ar}}_{\frak s}$ 
and $p \in {\Cal S}^{\text{bs}}(M_{n+1}) \and 
p \underset{\text{wk}} {}\to \bot M_n
\Rightarrow \lambda_{\frak s} 
= |\{c \in \bold J_{n+1}:p = \text{ tp}_{\frak s}(c,
M_{n+1},M_{n+2})\}|$ and $N = \cup\{M_n:n < \omega\}$. \nl
6) We say that $(N,\bar M,\bold J)$ is $K^{\text{br}}_{\frak s}$-full 
if it $\in K^{\text{br}}_{\frak s}$ and 
$p \in {\Cal S}^{\text{bs}}(M_{n+1}) \and p$ does not
fork over $M_{m+1} \and p \underset{\text{wk}} {}\to \bot 
M_m \Rightarrow \lambda_{\frak s} = 
|\{c \in \bold J_{n+1}:c$ realizes $p\}|$ and 
$N = \cup\{M_n:n < \omega\}$.
\nl
7) $\le_{\text{or}} = \le^{\frak s}_{\text{or}}$ 
is the following two place relation over $K^{\text{or}}_{\frak s}$:

$$
\align
(N^1,\bar M^1,\bar{\bold J}^1) \le_{\text{or}} 
(N^2,\bar M^2,\bar{\bold J}^2) \text{ \ub{iff} } &(N^1,
\bar M^1,\bar{\bold J}^1) \le_{\Cal W} (N^2,\bar M^2,\bar{\bold J}^2) \text{ and}\\
  &\bold J^1_0 = \bold J^2_0.
\endalign
$$
\mn
8) We say that $(N,\bar M,\bar{\bold J}) \in {\Cal W}$ is
prime if $(M_n,M_{n+1},\bold J_n) \in K^{3,\text{qr}}_{\frak s}$ 
for $n < \ell g(\bar M)$.
\enddefinition
\bigskip

\definition{\stag{705-5.15A} Definition}  We say ${\frak s}$ has enough regulars
when: for any $\bar M = 
\langle M_\alpha:\alpha \le \delta +1 \rangle$ which is
$\le_{\frak s}$-increasing continuous, $M_{\delta +1} \ne M_\delta$, 
\ub{there are} $\alpha < \delta$ and $c \in M_{\delta +1}$ such that
tp$_{\frak s}(c,M_\delta,M_{\delta +1})
\in {\Cal S}^{\text{bs}}_{\frak s}
(M_\delta)$ does not fork over $M_\alpha$ 
and $\alpha =0$
or $\alpha = \beta +1 \and p \bot M_\beta$ for some $\beta$.
\enddefinition
\bigskip

\remark{Remark}  If we are dealing with ${\frak s}^\kappa_{T,\lambda}$,
see \scite{705-stg.2}(3), then using regular types this property holds.
\endremark
\bigskip

\proclaim{\stag{705-5.16} Claim}  1) $K^{\text{ar}}_{\frak s} \subseteq 
K^{\text{br}}_{\frak s} \subseteq K^{\text{or}}_{\frak s}$. \nl
2)  $\le_{\text{or}}$ 
is a partial order on $K^{\text{or}}_{\frak s}$; for an
$\le_{\text{or}}$-increasing sequence of length $< \lambda^+_{\frak s}$, it
has a $\le_{\text{or}}$-lub which is a $\le_{\Cal W}$-{\rm lub}.
\nl
3) If $(N^\alpha,\bar M^\alpha,\bold J^\alpha) \in 
K^{\text{or}}_{\frak s}$ for
$\alpha < \delta < \lambda^+$ is $\le_{\text{or}}$-increasing, 
\ub{then} its
$\le_{\Cal W}$-{\rm lub} (see \scite{705-4.22}) is its 
$\le_{\text{or}}$-{\rm lub} (so it belongs to 
$K^{\text{or}}_{\frak s}$). \nl
4) In part (2) if $(N^\alpha,
\bar M^\alpha,\bar{\bold J}^\alpha) \in K^{\text{ar}}
_{\frak s}$ for $\alpha < \delta$ is $\le_{\text{or}}$-increasing 
\ub{then} the $\le_{\text{or}}$-{\rm lub} belongs to
$K^{\text{ar}}_{\frak s}$. \nl
5) In part (2) if $(N^\alpha,\bar M^\alpha,\bar{\bold J}^\alpha) \in
K^{\text{br}}_{\frak s}$ for $\alpha < \delta$ is 
$\le_{\text{or}}$-increasing, 
\ub{then} the $\le_{\text{or}}$-{\rm lub} 
(of this sequence) belongs to $K^{\text{br}}_{\frak s}$. \nl
6) If $(N^1,\bar M^1,\bar{\bold J}^1) \in 
K^{\text{or}}_{\frak s}$ \ub{then} there
is a $K^{\text{or}}_{\frak s}$-full 
$(N^2,\bar M^2,\bold J^2)$ such that $(N^1,
\bar M^1,\bar{\bold J}^1) \le_{\text{or}} 
(N^2,\bar M^2,\bar{\bold J}^2)$. \nl
7) If $(N^1,\bar M^1,\bar{\bold J}^1) 
\in K^{\text{ar}}_{\frak s}$ \ub{then} there
is a $K^{\text{ar}}_{\frak s}$-full $(N^2,
\bar M^2,\bar{\bold J}^2)$ such that $(N^1,
\bar M^1,\bar{\bold J}^1) \le_{\text{or}} 
(N^2,\bar M^2,\bar{\bold J}^2)$. \nl
8) If $(N^1,\bar M^1,\bar{\bold J}^1) \in K^{\text{br}}_{\frak s}$ 
\ub{then} there is a $K^{\text{br}}_{\frak s}$-full 
$(N^2,\bar M^2,\bar{\bold J}^2)$ such that $(N^1,
\bar M^1,\bar{\bold J}^1) \le_{\text{or}} 
(N^2,\bar M^2,\bar{\bold J}^2)$.  \nl
9) Like parts (3), (4), (5) for $K^{\text{or}}_{\frak s}$-full,
$K^{\text{ar}}_{\frak s}$-full, $K^{\text{br}}_{\frak s}$-full triples.
\endproclaim
\bigskip

\demo{Proof}  Straight (for (9) use the local character of non-forking.
\enddemo
\bigskip

\proclaim{\stag{705-5.17} Claim}  1) Assume that $(M,N,\bold J) \in
K^{3,\text{vq}}_{\frak s}$ or 
at least $(*)_{(M,N,\bold J)}$ below.  \ub{Then} we can find
$(\bar M,\bar{\bold J})$ such that $(**)_{(M,N,a),\bar M,\bar{\bold J}}$
below holds, where
\mr
\item "{$(*)_{(M,N,\bold J)}$}"  $(M,N,\bold J) 
\in K^{3,\text{bs}}_{\frak s}$ and
for no $N',b$ do we have $\bold J \subseteq N',M \le_{\frak s} N'
\le_{\frak s} N$, \nl
$b \in N \backslash N'$ and {\rm tp}$(b,N',N) \underset{wk} {}\to \pm M$
\sn
\item "{$(**)_{(M,N,\bold J),\bar M,\bar{\bold J}}$}" $\qquad (a) \quad
(M,N,\bold J) \in K^{3,\text{bs}}_{\frak s}$
\sn
\item "{${{}}$}" $\qquad (b) \quad
\bar M = \langle M_n:n < \omega \rangle,M_n \le_{\frak s} M_{n+1}$
\sn
\item "{${{}}$}"  $\qquad (c) \quad 
M_0=M$ and $\cup\{M_n:n < \omega\} = N$
\sn
\item "{${{}}$}"  $\qquad (d) \quad
(\bar M,\bar{\bold J}) \in K^{\text{or}}_{\frak s}$
\sn
\item "{${{}}$}"  $\qquad (e) \quad
(M_n,M_{n+1},\bold J_n) \in K^{3,\text{qr}}_{\frak s}$
\sn
\item "{${{}}$}"  $\qquad (f) \quad \bold J_0 = \bold J$
\sn
\item "{${{}}$}"  $\qquad (g) \quad$ if $n < \omega$ and $b \in \bold
J_{n+1}$ then {\rm tp}$(b,M_n,M_{n+1}) \underset{\text{wk}} {}\to \bot M_0$
(follows by (clause (d)).
\ermn
2) If $M \le_{\frak s} N$ we can find $\bar M,\bar{\bold J}$ such that
(a)-(e),(g) above holds.
\endproclaim
\bigskip

\demo{Proof}  1) By \scite{705-6.4A}(4), we 
know that $(*)_{(M,N,\bold J)}$ holds in
both cases.

We shall 
choose $M_n,\bold J_n$ by induction on $n$ satisfying the relevant
clauses in $(**)$.  Let $M_0 = M$, let $\bold J_0 = \bold J$ and let $M_1
\le_{\frak s} N$ be such that $(M_0,M_1,\bold J) \in 
K^{3,\text{qr}}_{\frak s}$, exists by \scite{705-5.5}(1).  If
$M_n \le_{\frak s} N$ is 
well defined, $n \ge 1$ let $\bold J_n$ be a maximal
subset of $\bold I_{M_n,N}$ independent in $(M_n,N)$ such that
$b \in \bold J_n \Rightarrow \text{ tp}_{\frak s}(b,M_n,N) \bot M_0$.

Lastly, let $M_{n+1} \le_{\frak s} N$ 
be such that $(M_n,M_{n+1},\bold J_n)
\in K^{3,\text{qr}}_{\frak s}$, exists by \scite{705-5.5}(1).  
To finish we need to prove
that $M_\omega =: \dbcu_{n < \omega} M_n$ is equal to $N$.  
Clearly $M_\omega \le_{\frak s} N$,
if $M_\omega \ne N$ then for some $b \in N \backslash M_\omega$ we have
tp$_{\frak s}(b,M_\omega,N) \in {\Cal S}^{\text{bs}}(M_\omega)$, by
$(*)_{(M,N,\bold J)}$ clearly tp$(b,M_\omega,N) \bot M_0$ and 
clearly for some $n < \omega$, tp$_{\frak s}(b,M_\omega,
N)$ does not fork over $M_n$ (and necessarily $n \ge 1$), and
similarly we have tp$(b,M_n,N) \underset{\text{wk}} {}\to \bot M_0$
so $b$ contradicts
the choice of $\bold J_n$ (as maximal such that ...).  So we are done. \nl
2) Should be clear.  \hfill$\square_{\scite{705-5.17}}$\margincite{705-5.17}
\enddemo
\bigskip

\proclaim{\stag{705-5.18} Claim}  [${\frak s}$ has enough regulars]. \nl 
1)  In \scite{705-5.17} we can get
\medskip

$(**)^+_{(M,N,\bold J),\bar M,\bar{\bold J}} \quad$  (a)-(f)  
as in \scite{705-5.17}
\smallskip

$\qquad \qquad \qquad (g)^+ \,\,
(\bar M,\bar{\bold J}) \in K^{\text{ar}}_{\frak s}$ (i.e. we
strengthen clause (d)).
\mn
2) In \scite{705-5.19} below we can add
\medskip

$\qquad \qquad \qquad (B)^+ \,\,$ like $(B)$ adding $(N,\bar
M,\bar{\bold J}) \in K^{\text{ar}}_{\frak s}$.
\endproclaim
\bigskip

\demo{Proof}  1) Similar to \scite{705-5.17} using the definition of
``${\frak s}$ has enough regulars". \nl
2) In \scite{705-5.19} note that $(C) \Rightarrow (B)^+$ by
\scite{705-5.18}(1) and $(B)^+ \Rightarrow (B)$ trivially.  \hfill$\square_{\scite{705-5.18}}$\margincite{705-5.18}
\enddemo
\bn
Now we arrive to ``understanding $K^{3,\text{vq}}_{\frak s}$".
\demo{\stag{705-5.19} Conclusion}  For every triple 
$(M,N,\bold J)$, the following are equivalent:
\mr
\item "{$(A)$}"  $(M,N,\bold J) \in K^{3,\text{vq}}_{\frak s}$
\sn
\item "{$(B)$}"  We can find $\bar M,\bar{\bold J}$ such that $(N,\bar M,
\bar{\bold J}) \in K^{\text{or}}_{\frak s},\bold J_0 = \bold J,M_0 = M,
N = \dbcu_{n < \omega} M_n$ and $(M_n,M_{n+1},\bold J_n) \in 
K^{3,\text{qr}}_{\frak s}$ (we then say $(N,\bar M,\bar{\bold J})$ 
is prime), note that
necessarily tp$_{\frak s}(b,M_{n+1},M_{n+2})$ is orthogonal to $M_0$
for every $n < \omega,b \in \bold J_{n+1}$
\sn
\item "{$(C)$}" $\qquad (a) \quad M \le_{\frak s} N$
\sn
\item "{${{}}$}"  $\qquad (b) \quad \bold J$ is independent in $(M,N)$
\sn
\item "{${{}}$}"  $\qquad (c) \quad$ if 
$M \cup \bold J \subseteq N' \le_{\frak s}
N$ and $b \in N \backslash N'$ and 
tp$_{\frak s}(b,N',N) \in {\Cal S}^{\text{bs}}_{\frak s}(N')$
\nl

\hskip40pt  \ub{then} tp$_{\frak s}(b,N',N) \bot M$
\sn
\item "{$(D)$}"  $\qquad (a),(b) \quad$ as above
\sn
\item "{${{}}$}"  $\qquad (c) \quad$ if $N \le_{\frak s} N^+,b \in N^+
\backslash M \backslash \bold J$ and $\bold J 
\cup \{b\}$ is independent in $(M,N^+)$ \nl

$\quad \quad \,\,$ \ub{then} 
tp$_{\frak s}(b,N,N^+) \in {\Cal S}^{\text{bs}}_{\frak s}(N)$ does not
fork over $M$
\sn
\item "{$(E)$}"  there is a {\rm uq}-decomposition $\langle M_i,a_j:i \le
\alpha,j < \alpha \rangle$ of $(M_0,N)$ such that $(M,M_0,\bold J)
\in K^{3,\text{vq}}_{\frak s},M_\alpha = N$ and each tp$_{\frak
s}(a_j,M_j,M_{j+1})$ is orthogonal to $M$.
\endroster
\enddemo
\bn
\relax From this we shall deduce (after the proof of \scite{705-5.19}):
\demo{\stag{705-5.20} Conclusion}  [${\frak s} = {\frak t}^+,{\frak t}$ is a
good$^+$ and successful frame.]
\nl
1)  $K^{3,\text{uq}}_{\frak s} = K^{3,\text{pr}}_{\frak s}$; so 
together with \scite{705-d.9}(2) we get uniqueness. \nl
2) $(M,N,\bold J) \in K^{3,\text{qr}}_{\frak s}$ iff $(M,N,\bold J) \in
K^{3,\text{vq}}_{\frak s}$; so together with \scite{705-4.11}(k) we get
uniqueness.
\enddemo
\bigskip

\proclaim{\stag{705-5.20A} Claim}   If 
$\delta < \lambda^+_{\frak s}$ and $\langle M_i:i \le \delta
\rangle$ is $\le_{\frak s}$-increasing continuous, and $\langle \bold
J_i:i \le \delta \rangle$ is increasing continuous and $(M,M_i,\bold
J_i)$ belongs to $K^{3,\text{vq}}_{\frak s}$ for $i < \delta$, \ub{then}
$(M,M_\delta,\bold J_\delta)$ belongs to $K^{3,\text{vq}}_{\frak s}$.
\endproclaim
\bn
\ub{Question}:  Can we assume less than ${\frak s} = {\frak t}^+$?
\bigskip

\demo{Proof of \scite{705-5.19}}  The following implications clearly suffice.
\mn
\ub{$(A) \Rightarrow (E)$}:  Let $\alpha=0,M_0 = N$.
\mn
\ub{$(E) \Rightarrow (D)$}:  Clauses (a), (b) are obvious, so let us
turn to (c).
Assume $b,N^+$ are as in clause (c) of (D), so by Claim \scite{705-5.30}(1)
we know that tp$_{\frak s}(b,M_0,N^+)$ does not fork over $M$; now we 
prove by induction on
$i \le \alpha$ that tp$(b,M_i,N^+)$ does not fork over $M$, for $i=0$ see
above, for $i$ limit use Axiom (E)(h), 
for $i$ successor by the definition of orthogonality.  For $i =
\alpha,M_\alpha = N$.  So we are done. 
\mn
\ub{$(C) \Rightarrow (B)$}: by \scite{705-5.17}
\mn
\ub{$(B) \Rightarrow (A)$}: by \scite{705-5.14}(2)
\mn
\ub{$(A) \Rightarrow (D)$}:  Clauses (a), (b) are obvious.  For clause (c),
as tp$_{\frak s}(b,M,N^+) \in {\Cal S}^{\text{bs}}_{\frak s}(M)$,
there is $M' \le_{\frak s} N^+$ 
such that $(M,M',b) \in K^{3,\text{pr}}_{\frak s}$ (recalling
${\frak s}$ has primes).  By \scite{705-5.2} we know $\bold J$ is
independent over $(M,M',N^+)$ hence by Definition \scite{705-4.9}, we have
NF$_{\frak s}(M,M',N,N^+)$ hence \scite{705-stg.32}) 
we get that tp$_{\frak s}(b,N,N^+)
\in {\Cal S}^{bs}(N)$ does not fork over $M$ as required.  (Actually
not used). 
\mn
\ub{$(D) \Rightarrow (C)$}:

Again the problem is to 
prove clause (c) of (C) so toward contradiction assume
that $M \cup \bold J \subseteq N' <_{\frak s} N$ and $b 
\in N \backslash N'$
and $p = \text{ tp}_{\frak s}(b',N',N) \in {\Cal S}^{\text{bs}}(N')$ is 
not orthogonal to $M$ hence by Hypothesis \scite{705-5.14}(d)
is not weakly orthogonal to $M$.  So for some
$q \in {\Cal S}^{\text{bs}}_{\frak s}(M)$ 
we have $p \pm q$, and let $q_1 \in {\Cal S}^{\text{bs}}_{\frak s}
(N')$ be a nonforking extension of $q$.  We can find $N_2$ such that
$N' \cup \{b\} \subseteq N_2 \le_{\frak s} N$ and $(N',N_2,b) \in
K^{3,\text{pr}}_{\frak s}$.  
So (see \scite{705-5a.2}) $q_1$ has some extension $q_2 \in
{\Cal S}^{\text{bs}}_{\frak s}(N_2)$ which is not a 
nonforking extension of $q$, and so we can
find $N_4$ and $c$ such that $N \le_{\frak s} N_4$ 
and $q_2 = \text{ tp}_{\frak s}(c,N_2,N_4)$.  
Now as $c$ realizes $q_1$ clearly tp$_{\frak s}(c,N',N)$ does not fork
over $M$ hence $\bold J \cup \{c\}$ is independent in $(M,N_4)$; but 
as $c$
realizes $q_2$ clearly tp$_{\frak s}(c,N_2,N_4)$ does not fork over
$M$, hene as $N_2 \le_{\frak s} N \le_{\frak s} N$
also tp$_{\frak s}(c,N,N_4)$ forks over
$M$.  So we have gotten a contradiction to clause (c) of (D).
\hfill$\square_{\scite{705-5.19}}$\margincite{705-5.19}
\enddemo
\bigskip

\demo{Proof of \scite{705-5.20}}  1) This is a special case of (2). \nl
2) The ``only" if implication we already proved in
\scite{705-c.4}(2), more exactly \scite{705-4.10}.  For the other direction assume 
$(M,N,\bold J) \in K^{3,\text{vq}}_{\frak s}$ and by
\scite{705-5.17}(1) applied to
$(M,N,\bold J)$ we get $(\bar M,\bar{\bold J})$ satisfying
$(**)_{(M,N,\bold J),\bar M,\bar{\bold J}}$ of \scite{705-5.17}(1) hence
it is as in clause (B) of \scite{705-5.19} in particular
$\bold J_0 = \bold J$.  Let $\langle M^\beta_i:i < \lambda_{\frak
s}\rangle$ be
$\le_{\frak t}$-representations of $M_\beta$ for $\beta \le \omega$.  Now
by \scite{705-4.11K} there is a club $E$ of 
$\lambda_{\frak s} = \lambda^+_{\frak t}$ such 
that for $\delta \in E,(\langle M^\beta_\delta:\beta 
\le \omega \rangle,\langle \bold J_n \cap M^{n+1}_\delta:
n < \omega \rangle)$ are as in \scite{705-5.19} clause (B) for ${\frak t}$,
hence clause (A) so 
$(M^0_\delta,M^\omega_\delta,\bold J \cap M^\omega_\delta) 
\in K^{3,\text{vq}}_{\frak t}$,  hence (see
\scite{705-d.8} as $M_0,M_\omega \in K_{\frak s}$, more exactly by
\scite{705-zm.2}(1) if $\bold J$ is a singleton, 
by \scite{705-4.11K} in general) the triple
$(M_0,M_\omega,\bold J)$ belongs to $K^{3,\text{qr}}_{\frak s}$, as
required.   \hfill$\square_{\scite{705-5.20}}$\margincite{705-5.20}
\enddemo
\bigskip

\demo{Proof of \scite{705-5.20A}}  It is 
enough to check clause $(D)$ of \scite{705-5.19}, now clause (a)
is trivial, clause (b) holds by \scite{705-5.2}(3).  For proving clause (c)
we assume $M_\delta \le_{\frak s} N^+,c \in N^+ \backslash 
\bold J \backslash M$ and 
$\bold J \cup \{b\}$ is independent in $(M,N^+)$, and we should
prove that ``tp$_{\frak s}(c,M_\delta,N)$ 
belongs to ${\Cal S}^{\text{bs}}_{\frak s}(M_\delta)$ and does not
fork over $M$".  Now clearly for each $i < \delta$ the set 
$\bold J_i \cup \{b\}$
is independent in $(M,N^+)$ by monotonicity of independence.  Hence by
\scite{705-5.19} $(A) \Rightarrow (D)$ as we are assuming $(M,M_i,\bold
J_i)  \in K^{3,\text{vq}}_{\frak s}$ we can conclude that
tp$_{\frak s}(c,M_i,N) \in {\Cal S}^{\text{bs}}_{\frak s}(M_i)$ does not
fork over $M$; so this holds for every $i < \delta$.
Now tp$_{\frak s}(c,M_\delta,N^+)$ does not fork over $M$ by
Axiom (E)(h).   \nl
${{}}$   \hfill$\square_{\scite{705-5.20A}}$\margincite{705-5.20A}
\enddemo
\bn
\centerline {$* \qquad * \qquad *$}
\bn
We now try to show that there is a parallel to universal homogeneous or
saturated among $\{(M,N,a) \in
K^{3,\text{uq}}_{\frak s}:\text{tp}_{\frak s}(a,M,N)=p\}$
\definition{\stag{705-5.23} Definition}  1) If $M \in K_{\frak s},p \in
{\Cal S}^{\text{bs}}_{\frak s}(M)$ let 
$K^{3,\text{uq}}_{{\frak s},p} = \{(M^*,N^*,a^*) \in
K^{3,\text{uq}}_{\frak s}:M = M^*$ 
and $p = \text{ tp}_{\frak s}(a,M,N)\}$, we identify
$(M,N,a)$ with $(M,N,\{a\})$ and recall $M^* = \text{ Dom}(p^*)$. \nl
2) We say $(M,N,\bold J) \in K^{3,\text{vq}}_{\frak s}$ is full or is
$K^{3,\text{vq}}_{\frak s}$-full \ub{if} there
is a full $(N,\bar M,\bold{\bar J}) \in 
K^{3,\text{or}}_{\frak s}$ satisfying $M_0 = M,\cup\{M_n:n < \omega\} = 
N_0$.  If $\bold J = \{a\}$ we may write $a$ instead of $\bold J$ and say
$(M,N,\bold J)$ is $K^{3,\text{uq}}_{\frak s}$-full and 
if $p = \text{ tp}(a,M,N)$ we say $(M,N,\bold J)$ is 
$K^{3,\text{uq}}_{{\frak s},p}$-full.
\enddefinition
\bigskip

\proclaim{\stag{705-5.24} Claim}  1) If 
$(M,N,a) \in K^{3,\text{uq}}_{{\frak s},p^*}$ 
and $(N^1,\bar M^1,\bar{\bold J})$
is $K^{\text{or}}_{\frak s}$-full, $\bold J_0 = \{a^*\}$, \nl
$(M^1_0,M^1_1,a^*) \in K^{3,\text{uq}}_{{\frak s},p^*}$ and 
$M = M^1_0$ \ub{then} 
there is a $\le_{\frak s}$-embedding
$f$ of $N$ into $N^1 = \cup\{M^1_n:n < \omega\}$ over $M$ mapping $a$ to
$a^*$. \nl
2) If 
$(N^\ell,\bar M^\ell,\bar{\bold J}^\ell)$ is a 
$K^{\text{or}}_{\frak s}$-full and prime (see \scite{705-5.19}, clause (B)), 
$\bold J^\ell_0 = \{a_\ell\}$, {\rm tp}$_{\frak s}
(a_\ell,M^\ell_0,M^\ell_1)=p$
for $\ell=1,2$ (so $M^\ell_0 = \,{\text{\rm Dom\/}}(p)$ 
does not depend on $\ell$)
\ub{then} there is an isomorphism $f$ from $N^1$ onto $N^2$ over 
${\text{\rm Dom\/}}(p)$
which maps $a_1$ to $a_2$. \nl
3) Similar to (2) with $\bold J^\ell_0 = \bold J,M^\ell_0 = M_0$
and {\rm tp}$_{\frak s}(c,M^1_0,M^1_1) = { \text{\rm tp\/}}{\frak s}
(c,M^2_0,M^2_1)$ for $c \in \bold J$  for $\ell = 1,2$. \nl
4) If $(M,N,\bold J) \in K^{3,\text{vq}}_{\frak s}$ and
$(N^1,\bar M^1,\bar{\bold J}^1)$ is $K^{\text{or}}_{\frak s}$-full and 
${\Cal J} = \bold J^1_0$ and
$c \in \bold J \Rightarrow { \text{\rm tp\/}}_{\frak s}(c,M,N) =
{ \text{\rm tp\/}}_{\frak s}
(c,M^1_0,N^1)$ \ub{then} there is an embedding of
$N$ into $\cup\{M^1_n:n < \omega\}$ which is the identity on $M \cup
\bold J$. 
\endproclaim
\bigskip

\remark{\stag{705-5.24A} Remark}  In \scite{705-5.24} we can allow
stronger demands on $f$.  In part (1) if $M \cup \{a\} \subseteq N'
\le_{\frak s} N,f'$ a $\le_{\frak K}$-embedding of $N'$ into $M^1_n,f'
\subseteq \text{ id}_M,f'(a) = a^*$, \ub{then} we can require $f' \subseteq
f$.  \nl
In part (2) of \scite{705-5.24}, 
if $M^\ell_0 \cup \{a_\ell\} \subseteq M'_\ell \le M^\ell_n$
for $\ell = 1,2,f'$ an isomorphism from $M'_1$ onto $M'_2$ extending
id$_{\text{Dom}(p)} \cup \{\langle a_2,a_2 \rangle\}$ then we can require
$f' \subseteq f$.
\endremark
\bigskip

\demo{Proof}  1)  By \scite{705-5.19} $(A) \Rightarrow (B)$, we can find
$(N,\bar M,\bar{\bold J}) \in K^{\text{or}}_{\frak s}$ with $M_0 = M,
\bold J_0 = \{a\}$,
as in \scite{705-5.19}, clause (B).  Now we choose 
by induction on $n < \omega$
a $\le_{\frak s}$-embedding $f_n$ of $M_n$ into $M^1_n$ increasing with
$n,f_0 = \text{ id}_{M_0},f_1(a) = a^*$.  For $n=0$ this is trivial, for
$n=1$ note that tp$_{\frak s}(a,M,M_1) = \text{ tp}_{\frak s}(a^*,M,M^*_1)$ 
and recall the
definition of $(M,M_1,a) \in K^{3,\text{pr}}_{\frak s}$.  
For $n = m+1 > 1$, by the definition 
of ``\ub{full}" in \scite{705-5.15}(4) we can find a one-to-one mapping
$h_m$ from $\bold J_m$ into $\bold J^*_m$ such that
\mr
\item "{$(i)$}"  $b \in \bold J_m \Rightarrow \text{ tp}_{\frak s}(h_m(b),
M^*_m,M^*_{m+1})$ does not fork over Rang$(f_m)$
\sn
\item "{$(ii)$}"  $b \in \bold J_m \Rightarrow f_m(\text{tp}_{\frak s}
(b,M_m,M_{m+1})) = 
\text{ tp}_{\frak s}(h_m(b),\text{ Rang}(f_m),M^*_{m+1})$.
\ermn
Then 
choose $f_n \supseteq f_m,f_m(c) = h_m(c)$ for $c \in \bold J_m$ by the
definition of $(M_m,M_{m+1},\bold J_m) \in K^{3,\text{qr}}_{\frak s}$.
\nl
2) We choose by induction on $n$ a tuple $(N^1_n,N^2_n,f_n,\bold
I^1_n)$ such that (with $\ell \in \{1,2\}$):
\mr
\item "{$(a)$}"   $N^\ell_n \le_{\frak s} M^\ell_n$
\sn
\item "{$(b)$}"   $f_n$ is an isomorphism from $N^1_n$ onto $N^2_n$
\sn
\item "{$(c)$}"   $N^\ell_n \le_{\frak s} N^\ell_{n+1}$ and $f_n
\subseteq f_{n+1}$
\sn
\item "{$(d)$}"   $N^\ell_n = M^\ell_0$ and $f_0$ is the identity
\sn
\item "{$(e)$}"   $(N^1_0,N^1_1,\bold J) \in K^{3,\text{qr}}_{\frak s}$ and
$f_1(a_1) = a_2$
\sn
\item "{$(f)$}"   if $n = \ell \text{ mod } 2$ (where $\ell \in
\{1,2\}$) then
{\roster
\itemitem{ $(\alpha)$ }  $\bold I_n$ is a maximal subset of $\{b \in
\bold I_{N^\ell_n,M^\ell_n}:\text{ tp}(b,N^\ell_n,M^\ell_n) \bot
M_0\}$
\sn 
\itemitem{ $(\beta)$ }  $(N^\ell_n,N^\ell_{n+1},\bold J_n) \in
K^{3,\text{qr}}_{\frak s}$
\sn
\itemitem{ $(\gamma)_1$ }  if $\ell=1,f_{n+1} \restriction \bold I_n$
is a one-to-one mapping from $\bold I_n$ into $\bold J^2_{n+1}$
\sn
\itemitem{ $(\gamma)_2$ }     if 
$\ell=2,f^{-1}_{n+1} \restriction \bold I_n$
is a one-to-one mapping from $\bold I_n$ into $\bold J^1_{n+1}$.
\endroster}
\ermn
There is no problem to carry the induction and $f = \dbcu_{n < \omega}
f_n$ is as required as in the proof of \scite{705-5.17}. \nl
3), 4)  Similarly.  \hfill$\square_{\scite{705-5.24}}$\margincite{705-5.24}
\enddemo
\bigskip

\proclaim{\stag{705-5.25} Claim}  1) If 
$(M,N_\ell,a) \in K^{3,\text{uq}}_{{\frak s},p}$ 
is full for $\ell = 1,2$ \ub{then}
$N_1,N_2$ are isomorphic over $M \cup \{a\}$. \nl
2) Similarly for $K^{3,\text{vq}}_{\frak s}$. \nl
3) If $(M,N,\bold J) \in K^{3,\text{bs}}_{\frak s}$ \ub{then} for some $N'
\le_{\frak s} N''$ we have $J \subseteq N'$ and
$N \le_{\frak s} N''$ and $(M,N',\bold J)$
is $K^{3,\text{vq}}_{\frak s}$-full. 
\endproclaim
\bigskip

\demo{Proof}  By the proof of \scite{705-5.24}(2) (note that under somewhat
stronger assumption in \scite{705-5.20} we get uniqueness even without
assuming fullness).  \hfill$\square_{\scite{705-5.25}}$\margincite{705-5.25}
\enddemo
\bigskip

\proclaim{\stag{705-5.28} Claim}  1) Assume $(M_i,N_i,\bold J_i) \in
K^{3,\text{vq}}_{\frak s}$ 
for $i < \delta$ where $\delta < \lambda^+_{\frak s}$ and
$i < j < \delta \and c \in \bold J_i \Rightarrow \,{\text{\rm tp\/}}_{\frak
s}(c,M_j,N_j)$ does not fork over $M_i$.  Assume further that $\langle
M_i:i < \delta \rangle$ is $\le_{\frak s}$-increasing continuous,
$\langle N_i:i < \delta \rangle$ is $\subseteq$-increasing
continuous and $\langle \bold J_i:i < \delta \rangle$ is 
$\subseteq_{\frak s}$-increasing.  
Let $M_\delta = \cup\{M_i:i < \delta\},N_\delta =
\cup\{N_i:i < \delta\}$ and $\bold J_\delta = \cup\{\bold J_i:i <
\delta\}$.  \ub{Then} $(M_\delta,N_\delta,\bold J_\delta) \in
K^{3,\text{vq}}_{\frak s}$.
\endproclaim
\bigskip

\demo{Proof}  Note that $\bold J_\delta$ is independent in
$(M_\delta,N_\delta)$. \nl
[Why?  For each $i<j < \delta$ as $\bold J_j$ is independent in
$(M_j,N_j)$, but $\bold J_i \subseteq \bold J_j$ hence also $\bold
J_i$ is independent in $(M_j,N_j)$.  However, $c \in \bold J_i
\Rightarrow \text{ tp}_{\frak s}(c,M_j,N_j)$ does not fork over $M_i$,
hence $\bold J_i$ is independnet in $(M_i,M_j,N_j)$, so as $N_j
\le_{\frak s} N_\delta$, clearly $\bold J_i$ is independent in
$(M_i,M_j,N_\delta)$.  As fixing $i < \delta$ this holds for every $j
\in (i,\delta)$ and as $\langle M_j:j \in (i,\delta) \rangle$ is
increasing by \scite{705-4.11} as get that $\bold J_i$ is independent in
$(M_i,\dbcu_{j \in (i,\delta)} M_j,N_\delta)$ which means that it is
independent in $(M_i,M_\delta,N_\delta)$.  So $\bold J_i$ is
independent in $(M_\delta,N_\delta)$.  As $\bold J_i$ is increasing
with $i < \delta$ by \scite{705-5.2}(3), $\bold J_\delta = \cup\{\bold
J_i:i < \delta\}$ is independent in $(M_\delta,N_\delta)$ as
required.]

We shall use Claim \scite{705-5.19}, our desired conclusion is clause (A) for
$(M_\delta,N_\delta,\bold J_\delta)$ so it is enough to check clause
(D).  So let $M_\delta \le_{\frak s} N^+_\delta,b \in M^+ \backslash \bold
J_\delta \backslash M_\delta$ and 
assume that $\bold J_\delta \cup \{b\}$ is independent
in $(M_\delta,N^+_\delta)$.  
So tp$(b,M_\delta,N^+_\delta) \in {\Cal S}^{bs}_{\frak s}$ hence for some
$i(*) < \delta$ the type tp$(b,M_\delta,N^+_\delta)$ does not fork over
$M_{i(*)}$.  It is enough to prove that for every $i \in
[i(*),\delta)$, the type tp$(b,N_i,N^+_\delta) \in 
{\Cal S}^{\text{bs}}_{\frak s}(N_i)$ does not fork 
over $M_i$ (as then the nonforking extension $q
\in {\Cal S}^{\text{bs}}_{\frak s}(N_\delta)$ of tp$(b,M_\delta,N^+_\delta)$
satisfies $i \in [i(*),\delta) \Rightarrow q \restriction N_i = 
\text{ tp}(b,N_i,N^+_\delta))$.  
As $\bold J_i \cup \{b\} \subseteq \bold J_\delta \cup \{\delta\}$ and
$c \in \bold J_i \cup \{b\} \Rightarrow \text{
tp}(c,M_\delta,N^+_\delta)$ does not fork over $N_i$ it follows that
$\bold J_i \cup \{b\}$ is independent in $(M_i,M_\delta,N_i)$.
As $(M_i,N_i,\bold J_i) \in K^{3,\text{vq}}_{\frak s}$ 
and \scite{705-5.19} clearly tp$(b,N_i,N^+_\delta)$ does not fork over
$M_i$, and as said earlier this suffices.  
\nl
${{}}$  \hfill$\square_{\scite{705-5.28}}$\margincite{705-5.28}
\enddemo
\bigskip

\proclaim{\stag{705-7n.16} Claim}  Assume 
$\langle M_i:i \le \alpha \rangle$ is 
$\le_{\frak s}$-increasing continuous and $(M_i,M_{i+1},\bold J_i) \in
K^{3,\text{vq}}_{\frak s}$ and 
$\bold J_i$ is independent in $(M_0,M_i,M_{i+1})$.  \ub{Then}
$(M_0,M_\alpha,\cup \{\bold J_i:i < \alpha\}) \in 
K^{3,\text{vq}}_{\frak s}$.
\endproclaim
\bigskip

\demo{Proof}  We 
prove the statement by induction on $\alpha$.  Let $\bold J =
\cup\{\bold J_i:i < \alpha\}$.  First note that $\bold J$ is
independent in $(M_0,M_\alpha)$ by \scite{705-4.10}[??].  To prove that
$(M_0,M_\alpha,\bold J) \in K^{3,\text{vq}}_{\frak s}$ by \scite{705-5.20A} it
suffices to prove clause (D) there, so assume $M_\alpha \le_{\frak s}
N$ and $b \in N \backslash \bold J \backslash M_0$ and $\bold J \cup
\{b\}$ is independent in $(M_0,N)$.  If $\alpha = 0$ this is trivial.
For $\alpha$ limit it is enough to show for every $i < \alpha$ that
tp$(b,M_i,N)$.  But each $i < \alpha$, clearly $\{\bold J_j:
j < i\} \cup \{b\}$ is independent in $(M_0,N^+)$ hence by the
induction hypothesis and \scite{705-5.19} we know that tp$(b,M_i,N)$ does
not fork over $M_i$; as this holds for every $i < \alpha$, we can
deduce that tp$(b,M_\alpha,N)$ does not fork over $M_0$ as required.
\nl
So we are left with the case $\alpha = \beta +1$.  So $\bold J \cup
\{b\}$ is independent in $(M_0,N)$ and clearly $(\bold J \cup \{b\})
\cap M_\beta = \cup \{\bold J_i:i < \beta\}$, so by part
\scite{705-5.30}(1) we know that
$\bold J \cup \{b\} \backslash M_\beta = \bold J_\beta \cup \{b\}$ is
independent in $(M_0,M_\beta,N)$.  As $(M_\beta,M_\alpha,\bold
J_\beta) \in K^{3,\text{vq}}_{\bold s}$ by \scite{705-5.30}(1), 
tp$(b,M_\alpha,N)$ does not
fork over $M_\beta$ hence over $M_0$.  \hfill$\square_{\scite{705-5.30}}$\margincite{705-5.30}  
\enddemo
\newpage

\head {\S8 Tries to decompose and independence of sequences of models} \endhead  \resetall \sectno=8
 \spuriousreset
\bigskip

We try to find smooth or otherwise good decompositions; at present only
\scite{705-6.1} works.  We shall get really what we want after using
fullness + regular types.
This assumption having ${\Cal S}^{\text{na}}$ equal to 
${\Cal S}^{\text{bs}}$, i.e., fullness is ``soft", see \S9.
\demo{\stag{705-6.0} Hypothesis}
\mr
\item "{$(a)$}"    ${\frak s}$ is a successful good$^+$ frame, 
\sn
\item "{$(b)$}"  ${\frak s}$ has primes,
\sn
\item "{$(c)$}"  $K^{3,\text{uq}}_{\frak s} = K^{3,\text{pr}}_{\frak s}$, 
moreover $K^{3,\text{vq}}_{\frak s} = K^{\text{qr}}_{\frak s}$ and
\sn
\item "{$(d)$}"  $\bot = \underset{\text{wk}} {}\to \bot$.
\ermn
The last hypothesis is reasonable by the last section and also by
\scite{705-6.2}? below, but in some
examples it holds without going to a 
successor, so we use it as an hypothesis.
\enddemo
\bigskip

\definition{\stag{705-5.11A} Definition}   We say that $\langle M_i,a_j:
i \le \alpha,j < \alpha \rangle$ is a smooth decomposition inside $N$
over $M$ if:
\mr
\item "{$(a)$}"  it is a decomposition inside $N$ over $M$ (see
Definition \scite{705-c.1A}), which mean
\nl

$\langle M_i:i \le \alpha \rangle$ is $\le_{\frak s}$-increasing
continuous \nl

$M_\alpha \le_{\frak s} N$ \nl

$M = M_0$ \nl

tp$_{\frak s}(a_i,M_i,M_{i+1}) \in {\Cal S}^{\text{bs}}_{\frak s}
(M_i)$ \nl

$(M_i,M_{i+1},a_i) \in K^{3,\text{pr}}_{\frak s}$
\sn
\item "{$(b)$}"   for every $i < \beta$ there is
nonlimit $j \le i$ such that tp$_{\frak s}
(a_i,M_i,M_{i+1})$ does not fork over $M_j$
and is orthogonal to $M_{j-1}$ if $j > 0$.
\endroster
\enddefinition
\bigskip

\proclaim{\stag{705-6.1} Claim}  [No need of \scite{705-6.0}(b).]
If $M \le_{\frak s} N$ \ub{then}
we can find $\bar M = \langle M_i,a_j:i \le \alpha,j < \alpha \rangle$ 
such that
\mr
\item "{$\boxtimes_{M,N,\bar M},\bar a(a)$}"  $M_0 = M$ 
and $N \le_{\frak K} M_\alpha$
\sn
\item "{$(b)$}"  $M_i$ is $\le_{\frak K}$-increasing continuous
\sn
\item "{$(c)$}"  $(M_i,M_{i+1},a_i) \in K^{3,\text{uq}}_{\frak s}$ and 
$\in K^{3,\text{pr}}_{\frak s}$ if ${\frak s}$ has primes
\sn
\item "{$(d)$}"  for each 
$j < \alpha$ either {\rm tp}$(a_j,M_j,M_{j+1})$ does not
fork over $M_0$ \ub{or} it is weakly orthogonal to $M_0$.
\endroster
\endproclaim
\bigskip

\demo{Proof}  We try to choose by induction on $i$ a pair $(M_i,N_i)$ and
if $i = j+1$ also $a_j$ such that:  $M_i \le_{\frak s} N_i,M_i$ is
$\le_{\frak s}$-increasing continuous, $N_i$ is
$\le_{\frak s}$-increasing continuous, and the $M_i,a_j$ satisfy the relevant
cases of clauses (b), (c), (d) and $M_0 = M,N_0 = N$ 
and $i = j+1 \Rightarrow \neg$NF$_{\frak s}(M_j,N_j,M_i,N_i)$.  
We cannot succeed (as ${\frak s}$ is good$^+$ and successful, see \S1)
and we can
define for $i=0$ and $i$ limit.  Hence for some $i$ we have $(M_i,N_i)$
but cannot choose $M_{i+1},N_{i+1},a_i$.  If $M_i \ne N_i$ there is
$b_i \in N_i$ such that tp$(b,M_i,N_i) \in {\Cal S}^{\text{bs}}(M_i)$,
and one of the following cases occurs.
\enddemo
\bn
\ub{Case 1}:  tp$(b_i,M_i,N_i)$ is orthogonal to $M_0$.  \nl
Then we let
$N_{i+1} = N_i,a_i = b_i$ and we can find 
$M_{i+1} \le_{\frak s} N_{i+1}$ such that $N_i \le_{\frak s} N_{i+1}$ and
$(M_i,M_{i+1},a_i) \in K^{3,\text{uq}}_{\frak s}$ and if ${\frak s}$ has
primes even $\in K^{3,\text{pr}}_{\frak s}$.  
All the induction demands hold and
$\neg$NF$(M_i,N_i,M_{i+1},N_i)$ as NF implies disjointness (by the
definition of NF$_{\frak s}$ see \chaptercite{600}).
\bn
\ub{Case 2}:  tp$(b_i,M_i,N_i)$ is not orthogonal to $M_0$.

So there is $p_i \in {\Cal S}^{\text{bs}}(M_i)$ which does not fork
over $M_0$ such that $p_i$,
tp$_{\frak s}(b_i,M_i,N_i)$ are not weakly orthogonal.
Hence we can find $N_{i+1} \in
K_{\frak s}$ such that $N_i \le_{\frak s} N_{i+1}$ and some $a_i \in
N_{i+1}$ realizing $p_i$ in $N_{i+1}$ and tp$_{\frak s}
(a_i,N_i,N_{i+1})$ is not the
nonforking extension of $p_i$ in ${\Cal S}^{\text{bs}}(N_i)$.  As we
can increase $N_{i+1}$ \wilog \, there is
$M_{i+1} \le_{\frak s} N_{i+1}$ such that $(M_i,M_{i+1},a_i)$ is in
$K^{3,\text{pr}}_{\frak s}$ if 
possible but always in $K^{3,\text{uq}}_{\frak s}$, 
so clearly $\neg$NF$_{\frak s}(M_i,N_i,M_{i+1},N_{i+1})$.  
So all the demands hold.  \nl
So if 
$M_i \ne N_i$ then we can continue the induction, contradiction, hence
$M_i = N_i$ and so $\alpha = i,\langle M_j:j \le \alpha \rangle,
\langle a_j:j < \alpha \rangle$ are as required.
\hfill$\square_{\scite{705-6.1}}$\margincite{705-6.1}
\bigskip

\proclaim{\stag{705-6.1A} Claim}   If $\langle M_i:i \le \alpha \rangle$
is $\le_{\frak s}$-increasing continuous $p \in 
{\Cal S}^{\text{bs}}(M_\alpha)$ does
not fork over $M_0$, \ub{then} we can find an $\le_{\frak s}$-increasing continuous sequence $\langle N_i:i \le \alpha \rangle$ and $a$
such that $M_i \le_{\frak s} N_i,a \in N_0$, a 
{\rm tp}$_{\frak s}(a,M_\alpha,N_\alpha)=p$ and 
$(M_i,N_i,a)$ is ${\frak s}$-prime for $i \le \alpha$. 
\endproclaim
\bigskip

\demo{Proof}  We choose by induction on $i \le \alpha$ a pair
$(N_i,f_i)$ a $a$ such that $N_i$ is $\le_{\frak s}$-increasing, $f_i$
is an $\le_{\frak s}$-embedding of $M_i$ into $N_i,f_i$ is increasing
continuous, $f_0 = \text{ id}_{M_0}$, tp$_{\frak s}(a,M_0,N_0) = p
\restriction M_0,(f_i(M_i),N_i,a) \in K^{3,\text{pr}}_{\frak s}$, 
tp$_{\frak s}(a,f_i(M_i),N_i)$ does not fork over $f_0(M_0) = M_0$.  
For $i=0$ use existence of primes.

For $i$ limit use \scite{705-5.28} and the hypothesis (c) of \scite{705-6.0}, 
and for $i = j+1$
use the definition of prime chasing arrows.  In the end, renaming
\wilog \, $f_i = \text{ id}_{M_i}$ for $i \le \alpha$. 
\hfill$\square_{\scite{705-6.1A}}$\margincite{705-6.1A} 
\enddemo
\bigskip

\proclaim{\stag{705-6.2} Claim}  [No use of \scite{705-6.0}(b),(c),(d).]
If $\boxtimes_{M,N,\bar M,\bar a}$ from Claim \scite{705-6.1} holds 
and $\bold J= \{a_i:{\text{\rm tp\/}}(a_i,M_i,M_{i+1})$ does not
fork over $M = M_0\}$ \ub{then} 
$(M,N,\bold J)$ belongs to $K^{3,\text{vq}}_{\frak s}$.
\endproclaim
\bigskip

\demo{Proof}  We shall use Claim \scite{705-5.19}, 
now the desired conclusion is
clause (A) there, so it suffices 
to prove clause (D) there.  Now subclauses
(a), (b) are obvious so let us prove subclause 
(c).  For this we prove by
induction on $\beta \le \alpha = \ell g(\bar a)$ that letting
$\bold J_\beta = \{a_i:i < \beta \text{ and tp}(a_i,M_i,M_{i+1})$ does not
fork over $M_0\}$, we have:
\mr
\item "{$\circledast_\beta$}"  if $M_\beta \le_{\frak s} N^+,n <
\omega$ and $b_\ell \in N^+$ for $\ell < n$ and 
$\bold J_\beta \cup \{b_0,\dotsc,b_{n-1}\}$ is independent in $(M_0,N^+)$ (and
$b_\ell \notin \bold J_\beta,\ell \ne k \Rightarrow b_\ell \ne b_k$ 
of course), \ub{then} $\{b_0,\dotsc,b_{n-1}\}$ is independent in
$(M,M_\beta,N^+)$.
\ermn
For $\beta=0$ this is trivial.  For $\beta$ limit this holds by
\scite{705-4.11}(2) as $(M_0,M_\gamma,\{b_\ell:\ell < n\})$ is independent
by the induction hypothesis for each $\gamma < \beta$.  
Lastly, let $\beta = \gamma +1$;
then by \scite{705-5.19}, as we have proved $\circledast_\gamma$, we have
$(M_0,M_\gamma,\bold J_\gamma) \in K^{3,\text{uq}}_{\frak s}$.
First assume $a_\gamma \in \bold J_\beta$.  So
\mr
\item "{$(*)$}"  we are given $n,b_0,\dotsc,b_{n-1}$ we let
$b_n=a_\gamma$ and we are assuming $\bold J_\beta \cup
\{b_0,\dotsc,b_{n-1}\}$ is independent. \nl
Hence $\bold J_\gamma \cup \{b_0,\dotsc,b_n\}$ is independent.
\ermn
Hence apply $(*)_\gamma$ for $n+1,b_0,\dotsc,b_{n-1},b_n$, so
$\{b_0,\dotsc,b_n\}$ is independent in $(M,M_\gamma,N^+)$
as $(M_\gamma,M_\beta,a_{\gamma +1}) \in K^{3,\text{pr}}_{\frak s}$
by \scite{705-5.2} we deduce that $\{b_0,\dotsc,b_{n-1}\}$ 
is independent in $(M_0,M_\beta,N^+)$, which
gives the desired conclusion.
\nl
Second, assume $a_\gamma \notin \bold J_\beta$ and
$n,b_0,\dotsc,b_{n-1}$ are given as in $(*)$.  By the induction
hypothesis $\{b_0,\dotsc,b_{n-1}\}$ is independent in
$(M,M_\gamma,N^+)$.  But 
tp$(a_\gamma,M_\gamma,M_\beta) = \text{ tp}(a_\gamma,M_\gamma,N^+)$ is
orthogonal to $M_0$, tp$(b_\ell,M_\gamma,N^+)$ does not fork over $M,M_0$ and
$(M_\gamma,M_\beta,a_\gamma) \in K^{3,pr}_{\frak s}$ so by
\scite{705-5.13} necessarily tp$(b,M_\beta,N^+)$ does not 
fork over $M_0$, as required.

Having carried the induction we got $\circledast_\alpha$ which for
$n=1$ is the statement $(D)$ of \scite{705-5.19} hence gives
the desired conclusion.  \hfill$\square_{\scite{705-6.2}}$\margincite{705-6.2}
\enddemo
\bn
Now we can show that any type in a sense is below a nonforking
combination of basic ones.  (Compare with \scite{705-6.3}).
\proclaim{\stag{705-6.3} Claim}  If $p \in {\Cal S}_{\frak s}(M)$ \ub{then}
for some $n,M_\ell(\ell \le n),a_k(k < n)$ and $b$ we have:
\mr
\item "{$(a)$}"  $M_0 = M$
\sn
\item "{$(b)$}"  $(M_\ell,M_{\ell +1},a_\ell) \in 
K^{3,\text{pr}}_{\frak s}$
\sn
\item "{$(c)$}"  {\rm tp}$_{\frak s}
(a_\ell,M_\ell,M_{\ell +1})$ does not fork over $M_0$, so
\sn
\item "{$(d)$}"  $\{a_\ell:\ell < n\}$ is independent in $(M_0,M_n)$ and
$(M_0,M_n,\{a_\ell:\ell < n\}) \in K^{3,\text{qr}}_{\frak s}$
\sn
\item "{$(e)$}"  $b \in M_n$ and $p = { \text{\rm tp\/}}(b,M,M_n)$.
\endroster
\endproclaim
\bigskip

\demo{Proof}  We can find $N_0,b$ such that $M \le_{\frak s} N_0$ and
$b \in N_0$ and tp$_{\frak s}(b,M,N_0) = p$.  By \scite{705-6.1}, \wilog, \, possibly
increasing $N$, we have $\boxtimes_{M,N,\bar M,\bar a}$ for some $\bar M,\bar a$
as there.  By \scite{705-6.2} we have $(M,N,\bold J) \in K^{3,vq}_{\frak s}$ for
some $\bold J$.  Among all such pairs $(N,\bold J)$ choose one $(N^*,
\bold J^*)$ with the cardinality of $\bold J^*$ being minimal.  
Let $\bold J^* = \{a_i:i < \theta\}$; by hypothesis \scite{705-6.0}(c) we
know $(M,N^*,\bold J^*) \in K^{3,\text{qr}}_{\frak s}$ so we 
can find an $M$-based pr-decomposition $\langle M_i,b_j:i \le
\theta,j < \theta \rangle$ over $M$, i.e. $M_0 = M$ such that 
tp$_{\frak s}(b_j,M,M_{j+1}) = \text{ tp}_{\frak s}(a_i,M,N^*)$ 
and tp$_{\frak s}(b_j,M_j,M_{j+1})$ does not
fork over $M$, of course.  So by the section hypothesis \scite{705-6.0} there is a
$\le_{\frak s}$-embedding of $N^*$ into $M_\theta$ over $M$ so \wilog \,
$N^* \le_{\frak s} M_\theta$.  Now if $\theta < \aleph_0$ we have gotten the
desired conclusion, otherwise $b \in N^* \subseteq M_\theta =
\dbcu_{i < \theta} M_i$ so for some $\beta < \theta$ we have $b \in N_\beta$
and has clearly $(M_0,M_\beta,\{a_i:i < \beta\}) \in 
K^{3,\text{qr}}_{\frak s}$ by \scite{705-6.0} so we
have gotten a contradiction to the choice of $(N^*,\bold J^*)$.
\hfill$\square_{\scite{705-6.3}}$\margincite{705-6.3}
\enddemo
\bigskip

\definition{\stag{705-6p.5} Definition}  1) We say that $\langle M_i:i <
\alpha \rangle$ is ${\frak s}$-independent over $M$ inside $N$ if $M
\le_{\frak s} M_i \le_{\frak s} N$ and we can find a $\le_{\frak
s}$-increasing sequence $\langle N_i:i \le \alpha \rangle$ such that
$N_0 = M,N \le_{\frak s} N_\alpha$ and NF$_{\frak
s}(M,N_i,M_i,N_{i+1})$.  We call $\langle N_i:i \le \alpha \rangle$ a
witness. \nl
2) For $\alpha = \lambda^+$ we define similarly. 
\enddefinition
\bigskip

\proclaim{\stag{705-6p.6} Weak Uniqueness Claim}  Assume
\mr
\item "{$(a)$}"  for $\ell = 1,2,\langle M^\ell_i:i < \alpha \rangle$
is ${\frak s}$-independent over $M_\ell$ inside $N_\ell$ as witnessed
by $\langle N^\ell_i:i \le \alpha \rangle$
\sn
\item "{$(b)$}"  $f$ is an isomorphism from $M_1$ onto $M_2$
\sn
\item "{$(c)$}"  for $i < \alpha,f_i$ is an isomorphism from $M^1_i$
onto $M^2_i$ extending $f$.
\ermn
\ub{Then} there is $N_3$ such that $N_2 \le_{\frak s} N_3$ and a
$\le_{\frak s}$-embedding $f^*$ of $N_1$ into $N_3$ extending every $f_i$.
\endproclaim
\bigskip

\demo{Proof}  Without loss of generality $M_1 = M_2$ call it $M$ and
$f_0$ is the identity on $M$, so $N^\ell_0 = M$.

We choose by induction on $i \le \alpha$ the tuple
$(N^3_i,g^1_i,g^2_i)$ such that
\mr
\item "{$(\alpha)$}"  $N^3_0 = M,g^\ell_0 = \text{ id}_M$
\sn
\item "{$(\beta)$}"  $N^3_i$ is $\le_{\frak s}$-increasing continuous
\sn
\item "{$(\gamma)$}"  $g^\ell_i$ is a $\le_{\frak s}$-embedding of
$N^\ell_i$ into $N^3_i$ for $\ell =1,2$
\sn
\item "{$(\delta)$}"  $g^\ell_i$ is increasing continuous for $i$
\sn
\item "{$(\varepsilon)$}"  $(\forall a \in M^1_i)(g^2_{i+1}(f_i(a)) =
g^1_{i+1}(a))$.
\ermn
For $i=0,i$ limit this is obvious.  For $i = j+1$ use the uniqueness
of NF$_{\frak s}$-amalgamation.  Having carried the induction, by
renaming we get the conclusion. \nl
${{}}$    \hfill$\square_{\scite{705-6p.6}}$\margincite{705-6p.6}
\enddemo
\bigskip

\demo{Proof}  Straightforward.
\enddemo
\bigskip

\proclaim{\stag{705-6p.7A} Claim}  1) Assume that $M \le_{\frak s} M_i
\le_{\frak s} N$ for $i < \alpha,\alpha < \lambda^+$.  \ub{Then} we can find
$N^+,\langle M^+_i:i < \alpha \rangle$ such that $\langle N^+_i:i \le
\alpha \rangle$
\mr
\item "{$(a)$}"  $M = N^+_0,N \le_{\frak s} N^+ = N^+_\alpha$
\sn
\item "{$(b)$}"  $\langle N^+_i:i \le \alpha \rangle$ is $\le_{\frak
s}$-increasing continuous
\sn
\item "{$(c)$}"  $M_i \le_{\frak s} M^+_i \le N_{i+1}$
\sn
\item "{$(d)$}"  {\rm NF}$_{\frak s}(M_i,N^+_i,M^+_i,N^+_{i+1})$
\sn
\item "{$(e)$}"  $M^+_i$ is $(\lambda,*)$-brimmed over $M_i$ hence
over $N$.
\ermn
2) In part (1) it follows that
\mr
\item "{$(f)$}"  there is $\bold J_i$ such that $(M,M^+_i,\bold J_i)
\in K^{3,\text{vq}}_{\frak s}$
\sn
\item "{$(g)$}"  if $u \subseteq \alpha$ and
$\{M_i:i \in u\}$ is independent over $M$ inside $N$ then
$\{M^+_i:i \in u\}$ is independent over $N$ inside $N^+$.
\ermn
3) If $M \le_{\frak s} M'_i$ for $i < \alpha$ then we can find $N$ and
$\bar M = \langle M_i:i < \alpha \rangle$ such that $\bar M$ is
independent inside $(M,N)$ and $M_i,M'_i$ are isomorphic over $M$ for
$i < \alpha$. 
\endproclaim
\bigskip

\demo{Proof}  Easy.
\enddemo
\bigskip

\proclaim{\stag{705-6p.7} Claim}  Assume that $\langle M_i:i < \alpha
\rangle$ is ${\frak s}$-independent 
over $M$, inside $N$ with $\bar N = \langle N_i:i \le \alpha \rangle$,
a witness. \nl
1) If $M \le_{\frak s} M'_i \le_{\frak s} M_i$ for $i < \alpha$, 
\ub{then} $\langle M'_i:i <
\alpha \rangle$ is ${\frak s}$-independent over $M$ inside $N$. \nl
2) In part (1), $\bar N$ 
is also a witness for $\langle M'_i:i < \alpha \rangle$
being ${\frak s}$-independent over $M$ inside $N$. \nl
3) If $M_i \le_{\frak s} M^+_i$ \ub{then} we can find $N^+,\langle N^+_i:i
\le \alpha \rangle,\langle f_i:i < \alpha \rangle$ such that:
\mr
\item "{$(a)$}"  $N_i \le_{\frak s} N^+_i$ and $N \le_{\frak s} N^+$
\sn
\item "{$(b)$}"  $\langle N^+_i:i \le \alpha \rangle$ is $\le_{\frak
s}$-increasing continuous
\sn
\item "{$(c)$}"   $f_i$ is a $\le_{\frak s}$-embedding of $M^+_i$ into
$N_{i+1}$ over $M_i$
\sn
\item "{$(d)$}"  $\langle N^+_i:i \le \alpha \rangle$ witness $\langle
f_i(M^+_i):i < \alpha \rangle$ is ${\frak s}$-independent over $M$
inside $N^+_1$.
\ermn
4) There are $\langle M^*_i:i < \alpha \rangle,\langle N^+_i:i
\le \alpha \rangle,N^+,\langle \bold J_i:i < \alpha \rangle$ such
that:
\mr
\item "{$(a)$}"  $M_i \le_{\frak s} M^*_i,N \le_{\frak s} N^+$ and
$N_i \le N^+_i$ for $i < \alpha$
\sn
\item "{$(b)$}"  $\langle N^+_i:i \le \alpha \rangle$ witness that
$\langle M^*_i:i < \alpha \rangle$ is ${\frak s}$-independent over $M$
inside $N^+$
\sn
\item "{$(c)$}"  $(M,M^*_i,\bold J_i) \in K^{3,\text{vq}}_{\frak s}$.
\ermn
5) If $N \le_{\frak s} N^+$ then $\langle M_i:i < \alpha \rangle$ is
${\frak s}$-independent over $M$ inside $N^+$.
\endproclaim
\bigskip

\demo{Proof}  1), 2) Straightforward. \nl
3) By \scite{705-6p.7A}(2) and then use part (1). \nl
4) By \scite{705-6p.7A}(2). \nl
5) Trivial.
\enddemo
\bigskip

\proclaim{\stag{705-6p.8} Claim}  Assume $M \le_{\frak s} M_i \le_{\frak s}
N$ for $i < \alpha$. \nl
1) For any $\bar M' = \langle M'_i:i < \alpha' \rangle$, a
permutation of $\bar M = \langle M_i:i < \alpha \rangle$ (that is for some
one to one function $\pi$ from $\alpha$ onto $\alpha',M_i =
M'_{\pi(i)}$) we have: $\bar M$ is ${\frak s}$-independent over $M$
inside $N$ \ub{iff} $\bar M'$ is ${\frak s}$-independent over $M$
inside $N$. \nl
2) $\bar M = \langle M_i:i < \alpha \rangle$ 
is 
${\frak s}$-independent over $M$ inside $N$ \ub{iff} every finite subsequence
$\bar M'$ of $\bar M$ is ${\frak s}$-independent over $M$ inside
$N$. \nl
3) Assume $(M,M_i,\bold J_i) \in K^{3,\text{vq}}_{\frak s}$.  
Then: $\langle
M_i:i < \alpha \rangle$ is ${\frak s}$-independent over $M$ inside
$N$ \ub{iff} $\cup\{\bold J_i:i < \alpha\}$ is independent in
$(M,N)$ and, of course, the $\bold J_i$ are pairwise disjoint.
\endproclaim
\bigskip

\demo{Proof}  1) By the symmetry assume $\langle M_i:i < \alpha
\rangle$ is ${\frak s}$-independent over $M$ inside $N$.
By \scite{705-6p.7A}(2) + \scite{705-6p.7}(1) \wilog, \, 
for each $i < \alpha$ for some $\bold J_i$
we have $(M,M_i,\bold J_i) \in K^{3,\text{vq}}_{\frak s}$.  Now using (3),
part (1) is translated to parts of \scite{705-5.2}. \nl
2) Similarly (note \scite{705-6p.7A}(2), clause (g)). \nl
3) First assume that $\langle M_i:
i < \alpha \rangle$ is ${\frak s}$-independent 
over $M$ inside $N$, let $\langle N_i:i \le \alpha
\rangle$ witness this; of course, $i \ne j \Rightarrow \bold J_i \cap
\bold J_j = 0,i < j \Rightarrow M_i \cap M_j = M$ (by properties of
NF$_{\frak s}$).  For $\beta = \alpha$ we get that $\cup\{\bold J_i:i
< \alpha\}$ is independent in $(M,N)$ as required.

We prove by induction on $\beta \le \alpha$ that $\cup\{\bold J_i:i <
\beta\}$ is independent in $(M,N_\beta) = (N_0,N_\beta)$; of course,
we can increase $N_\beta$ (see \scite{705-5.2}(2)).  For $\beta = 0$ this
is trivial, for $\beta$ limit use \scite{705-5.2}(1), for $\beta =
\gamma+1$, by \scite{705-5.3A}(2) we have $\bold J_\gamma$ is independent
in $(M,N_\gamma,N_\beta)$ and so by \scite{705-5.3A}(1) we have
$(\cup\{\bold J_i:i < \gamma\}) \cup \bold J_\gamma = \cup\{\bold
J_i:i < \beta\}$ is independent in $(M,N_\beta)$.

Second assume that the $\bold J_i$-s are pairwise disjoint and
$\cup\{\bold J_i:i < \alpha\}$ is independent in $(M,N)$.  Let $\bold
J_{< \beta} = \cup\{\bold J_i:i < \beta\}$, so $\langle \bold J_{<
\beta}:\beta \le \alpha \rangle$ is $\subseteq$-increasing continuous.

We now choose by induction on $\beta \le \alpha$, the tuple
\footnote{the sequence $\langle M^*_\beta:\beta \le \alpha \rangle$
will witness that $\langle M_i:i < \alpha \rangle$ is independent in $(M,N)$}
$(M^*_\beta,N^*_\beta,\bold J^*_\beta)$ such that:
\mr
\item "{$(a)$}"  $M^*_\beta$ is $\le_{\frak s}$-increasing continuous
\sn
\item "{$(b)$}"  $N^*_\beta$ is $\le_{\frak s}$-increasing continuous
\sn
\item "{$(c)$}"  $M^*_0 = M,N^*_0 = N$
\sn
\item "{$(d)$}"  $M^*_\beta \le_{\frak s} N^*_\beta$
\sn
\item "{$(e)$}"  $M_i \le_{\frak s} M^*_\beta$ for $i < \beta$
\sn
\item "{$(f)$}"  $\bold J^*_\beta \subseteq N^*_\beta \backslash M
\backslash \bold J_{< \beta}$
\sn
\item "{$(g)$}"  $\bold J^*_\beta$ is $\subseteq$-increasing continuous
\sn
\item "{$(h)$}"  $(M,M^*_\beta,\bold J^*_\beta \cup \bold J_{<
\beta})$ belongs to $K^{3,\text{vq}}_{\frak s}$
\sn
\item "{$(i)$}"   $\bold J^*_\beta \cup \bold J_{< \alpha}$ is
independent in $(M,N^*_\beta)$
\sn
\item "{$(j)$}"   NF$_{\frak s}(M,M^*_\gamma,M_\gamma,N^*_\gamma$.
\ermn
Note that by clauses (a),(e),(j) this is enough to prove
$(M^*_\beta:\beta \le \alpha)$ witness that $\langle M_i:i < \alpha
\rangle$ is independent in $(M,N)$, (see Definition \scite{705-6p.5}) as
required.
\nl
For $\beta=0$ let $M^*_\beta = M,N^*_\beta = N$ and $\bold J^*_\beta =
\emptyset$; easy to check.

For $\beta$ a limit ordinal let $M^*_\beta = \cup\{M^*_\gamma:\gamma <
\beta\},N^*_\beta = \cup\{N^*_\gamma:\gamma < \beta\}$ and $\bold
J^*_\beta = \cup\{\bold J^*_\gamma:\gamma < \beta\}$ the least obvious
point is clause (h) which holds by \scite{705-5.28} 
and clause (i) which holds by \scite{705-5.2}.

Lastly, for $\beta = \gamma +1$ as $(M,M^*_\gamma,\bold J^*_\gamma \cup
\bold J_{< \gamma}) \in K^{3,\text{vq}}_{\frak s}$ by clause (h), and $\bold
J^*_\gamma \cup \bold J_{< \gamma} \subseteq \bold J^*_\gamma 
\cup \bold J_{< \beta}$ 
is independent in $(M,N^*_\gamma)$ by clause (i), by \scite{705-5.30}(1)
we deduce that $\bold J_\gamma = \bold J^*_\gamma \cup \bold J_{<
\beta} \backslash \bold J^*_\gamma \cup \bold J_{< \gamma}$ 
is independent in $(M,M^*_\gamma,N^*_\gamma)$.
So as $(M,M_\gamma,\bold J_\gamma) \in K^{3,\text{vq}}_{\frak s}$, by the
definition of $K^{3,vq}_{\frak s}$ (see Definition \scite{705-4.9}) we get
NF$_{\frak s}(M,M^*_\gamma,M_\gamma,N^*_\gamma)$ hence $\bold
J_\gamma$ is independent in $(M,M^*_\gamma,N^*_\gamma)$ by
\scite{705-5.3A}(2).  By \scite{705-5.30}(2) we can find
$(M'_\gamma,N'_\gamma)$ such that: $M^*_\gamma \le_{\frak s} M'_\gamma
\le_{\frak s} N'_\gamma,N^*_\gamma \le_{\frak s} N'_\gamma,M'_\gamma$
is $(\lambda,*)$-brimmed over $M^*_\gamma,N'_\gamma$ is
$(\lambda,*)$-saturated over $N^*_\gamma$ and
$(M'_\gamma,N'_\gamma,\bold J_\gamma) \in K^{3,\text{vq}}_{\frak s}$ and
$\bold J_\gamma$ is independent in $(M,M'_\gamma,N'_\gamma)$.  By
Definition \scite{705-4.9} (of $K^{3,\text{vq}}_{\bold s}$) this implies NF$_{\frak
s}(M,M_\gamma,M'_\gamma,N'_\gamma)$.  There are also
$(f_\beta,N''_\gamma)$ such that $N'_\gamma \le_{\frak s}
N''_\gamma,f_\beta$ is a $\le_{\frak s}$-embedding of $M'_\gamma$ into
$N''_\gamma$ over $M^*_\gamma$ and NF$_{\frak
s}(M^*_\gamma,N^*_\gamma,f_\beta(M'_\gamma),N''_\gamma)$, simply by
the existence of NF$_{\frak s}$-amalagmation.

As we also have NF$_{\frak s}
(M,M_\gamma,M^*_\gamma,N^*_\gamma)$ (see above), by
transitivity for NF$_{\frak s}$ we have NF$_{\frak
s}(M,M_\gamma,f_\beta(M'_\gamma),N''_\gamma)$.  As we have NF$_{\frak
s}(M,M_\gamma,M'_\gamma,N'_\gamma)$ by the uniqueness of NF$_{\frak
s}$-amalgamation, possibly increasing $N''_\beta$, we can extend
$f_\beta$ to $f'_\beta$, a $\le_{\frak s}$-embedding of $N'_\gamma$
into $N''_\gamma$ such that id$_{M_\gamma} \subseteq f'_\beta$.  
Let $N^*_\beta = N''_\beta,M^*_\beta =
f'_\beta(N'_\gamma)$, note that $f'_\beta(M'_\gamma)$ is
$(\lambda,*)$-brimmed over $f'_\beta(M^*_\gamma) = M^*_\gamma$ and
$\bold J^*_\gamma \cup \bold J_{< \gamma}$ is independent in
$(M,f'_\beta(M'_\gamma))$ and $\bold J^*_\gamma \cup \bold J_{< \gamma}
\subseteq f'_\beta(M^*_\gamma) = M^*_\gamma$, hence we can find 
$\bold J'_\gamma \subseteq
M^*_\beta \backslash (\bold J^*_\gamma \cup \bold J_{< \gamma})
\backslash M$ such that: $\bold J'_\gamma \cup 
\bold J^*_\gamma \cup \bold J_{< \gamma}$ is independent in 
$(M,f'_\beta(M'_\gamma))$ and $(M,f'_\beta(M'_\gamma),
\bold J'_\gamma \cup J^*_\gamma \cup \bold J_{< \gamma}) 
\in K^{3,\text{vq}}_{\frak s}$ by \sciteu{705-xxX}.  As $M \le_{\frak s}
f'_\beta(M'_\gamma) \le_{\frak s} f'_\beta(N'_\gamma) = M^*_\beta$ and
$(M,f_\beta(M'_\gamma),\bold J'_\gamma \cup \bold J^*_\gamma \cup
\bold J_{< \gamma}) \in K^{3,\text{vq}}_{\frak s}$ and
$(f'_\beta(M'_\gamma),f'_\beta(N'_\gamma),\bold J_\gamma) \in
K^{3,\text{vq}}_{\frak s}$ we get by \scite{705-5.30}(3) 
that $(M,f'_\beta(N'_\gamma),\bold J'_\gamma \cup
\bold J^*_\gamma \cup \bold J_{< \gamma} \cup \bold J_\gamma) =
(M,M^*_\beta,(\bold J'_\gamma \cup \bold J^*_\gamma) \cup \bold J_{<
\beta}) \in K^{3,\text{vq}}_{\frak s}$.

Let $\bold J^*_\beta = \bold J'_\gamma \cup \bold J^*_\gamma$, so we have
almost finished proving the induction step, we still need: $\bold J^*_\beta$
disjoint to $\bold J^*_{< \alpha}$ and $\bold J^*_\beta \cup \bold
J_{< \alpha}$ is independent in $(M,N^*_\beta)$; for this we know that
$\bold J^*_\beta \cup \bold J_{< \gamma}$ is independent in
$(M,f'_\beta(M'_\gamma))$ and $M \le_{\frak s} f'_\beta(M'_\gamma)
<_{\frak s} N^*_\beta$ and $(\bold J^*_\beta \cap \bold J_{< \alpha})
\backslash f'_\beta(M'_\gamma) = (\bold J_{< \alpha} \backslash \bold
J_{< \gamma})$ hence it suffices to prove that $\bold J_{< \alpha}
\backslash \bold J_{< \gamma}$ is independent in
$(M,f_\beta(M'_\gamma),N^*_\beta)$.  But NF$_{\frak
s}(M^*_\gamma,f_\beta(M'_\gamma),N^*_\gamma,N''_\gamma)$ and $\bold
J_{< \alpha} \backslash \bold J_{< \gamma} \subseteq N^*_\gamma$ is
independent in $(M,M^*_\gamma,N^*_\gamma)$ as stated above, so by
\scite{705-5.3A}(2) we are done. \nl
${{}}$    \hfill$\square_{\scite{705-6p.8}}$\margincite{705-6p.8}     
\enddemo
\bigskip

\demo{\stag{705-6p.9} Conclusion}  1) If 
$\langle M_i:i < \alpha \rangle$ is ${\frak s}$-independent over $M$
inside $N,f_i$ is an isomorphism from $M_0$ onto $M_i$ over $M$ for
$i < \alpha$ and $\pi$ is a permutation of $\alpha$ and $N^+$ is
$(\lambda,*)$-brimmed over $N$, \ub{then} for some automorphism $f$ of
$N^+$ over $M$ we have $\pi(i) = j \Rightarrow f_j \circ f^{-1}_i
\subseteq f$. \nl
2) Assume that $\langle M^\ell_i:i < \alpha \rangle$ is 
${\frak s}$-indiscernible over $M_\ell$ inside $N_\ell$ for $\ell=1,2$
and $f_i \supseteq f$ is an isomorphism from $M^1_i$ onto $M^2_i$ for
$i < \delta,f$ an isomorphism from $M_1$ onto $M_2$ and $N_\ell$ is
$(\lambda,*$)- brimmed over $\cup\{M^\ell_i:i < \alpha\}$.  \ub{Then}
there is an isomorphism from $N_1$ onto $N_2$ extending $\cup\{f_i:i <
\alpha\}$. 
\enddemo
\bigskip

\demo{Proof}  1) By (2). \nl
2) By \scite{705-6p.7}(3) and uniqueness of the
$(\lambda,*)$-saturated model over a model in $K_{\frak s}$.    
\hfill$\square_{\scite{705-6p.9}}$\margincite{705-6p.9}
\enddemo  
\bigskip

\definition{\stag{705-6p.9P} Definition}  1) We say that $\langle f_i:i <
\alpha \rangle$ is $(< \theta)$-indiscernible over $M$ if: for some
sequence $\bar{\bold a} = \langle a_\varepsilon:\varepsilon < \zeta
\rangle$ (possibly infinite), Dom$(f_i) = \{a_\varepsilon:\varepsilon
< \zeta\}$ for every $i < \alpha$ and for every partial one to one
function $\pi$ such that Dom$(\pi) \cup \text{ Rang}(\pi) \subseteq
\alpha$ and $|\text{Dom}(\pi)| < \theta$ there are $N^+,g$ such that
$N \le_{\frak s} N^+,g$ is an automorphism for $N^+$ over $M$ and for
every $i \in \text{ Dom}(\pi)$ we have: $g$ maps $f_i(\bar{\bold a})$
to $f_{\pi(i)}(\bar{\bold a})$. 
\enddefinition
\bigskip

\proclaim{\stag{705-6p.10} Claim}  Assume $\langle M_i:i < \alpha \rangle$
is ${\frak s}$-independent over $M$ inside $N,f_i \, (i < \alpha)$ is
an isomorphism from $M_0$ onto $M_i$ over $M$,
and $p \in {\Cal S}^{\text{bs}}(M_0)$.  The
following are equivalent:
\mr
\item "{$(A)$}"  $p \bot M$
\sn
\item "{$(B)$}"  $p \bot f_1(p)$
\sn
\item "{$(C)$}"  for some $i < j < \alpha$ we have $f_i(p) \bot f_j(p)$
\endroster
\endproclaim
\bigskip

\demo{Proof}  By \scite{705-6p.7A}(3) [NO!: ad proof] it is 
easy to show that \wilog \,
$M$ is $(\lambda,*)$-saturated and by \sciteu{705-xxX} \wilog \, $\alpha$
is infinite. \nl
\ub{$(B) \Leftrightarrow (C)$}: by the indiscernibility 
(\scite{705-6p.9}(1)).
\mn
\ub{$\neg(C) \Rightarrow \neg(A)$}:  By the indiscernibility we have
$i < j < \alpha \Rightarrow f_i(p) \pm f_j(p)$. \nl
First we can find $\langle M^*_n:
n < \omega \rangle,\le_{\frak s}$-increasing, $M^*_{n+1}$ is
$(\lambda,*)$-saturated over $M^*_n$ such that 
$\cup\{M^*_n:n < \omega\} = M$.  Second, we can also find
$\langle M^*_{0,n}:n < \omega \rangle \, \le_{\frak s}$-increasing
such that NF$_{\frak s}(M^*_n,M^*_{0,n},M^*_{n+1},M^*_{0,n+1})$ 
and $M_0 \le_{\frak s} \cup
M^*_{0,n}$ by \sciteu{705-yyY}.  
Third, \wilog \, $\cup\{M^*_{0,n}:n < \omega\} = M_0$.  
Fourth,
for some $n,p$ does not fork over $M^*_{0,n}$ so without loss of generality
$n=0$.  Hence  we
can consider $M^*_0,\langle f_i(M^*_{0,0}):i < \alpha \rangle,
\langle f'_i = f_i \restriction M^*_{0,0}:i < \alpha \rangle,\langle
p'_i = f_i(p \restriction M^*_{0,0}):i < \alpha \rangle$ and can choose
$f'_\alpha,f'_\alpha(M^*_{0,0}) \le_{\frak s} M_0$ such that $\langle f'_i
\restriction M^*_{0,0}:i \le \alpha \rangle$ is indiscernible over
$M^*_0$.  By the indiscernibility clearly 
$p'_0 \pm p'_\alpha$ but $p'_0\|p_0$ and there is $q \in
{\Cal S}^{\text{bs}}(M), q\|p'_\alpha$, so we are done.  
[This is similar to \scite{705-5.8}].
\mn
\ub{$(C) \Rightarrow (A)$}:

Assume $q \in {\Cal S}^{bs}(M),q \pm p$, so $q \pm p_i$ for $i < \alpha$
where $p_i = f_i(p)$.  Let $N^+,b$ be such that
$N \le_{\frak s} N^+,b \in N^+$ and tp$(b,N,N^+)$ is
a nonforking extension of $q$.  So by \scite{705-5a.3Y}(2), possibly increasing
$N^+$, for each $i < \alpha$ there is $a_i \in N^+$ such that tp$(a_i,N,N^+)$
is a nonforking extension of $p_i$ and $\{b,a_i\}$ is not independent in
$(N,N^+)$.  But by \scite{705-5.22}(1) $\{a_i:i < \alpha\}$ is independent
in $(N,N^+)$, contradicting \scite{705-4.25}(1) as $\alpha$ is infinite.
\hfill$\square_{\scite{705-6p.10}}$\margincite{705-6p.10}
\enddemo
\bn
The following will be used in \S9[?].
\proclaim{\stag{705-6p.for9} Claim}  Assume:
\mr
\item "{$(a)$}"  $\langle N^\ell_\alpha:\alpha \le \delta \rangle$ is
$\le_{\frak s}$-increasing continuous for $\ell=1,2$ and
$\lambda|\delta$
\sn
\item "{$(b)$}"  $M^\ell_i \le_{\frak s} N^\ell_0$ for $i < i^*$
\sn
\item "{$(c)$}"  if $N^\ell_0 \le_{\frak s} M \le_{\frak s}
N^\ell_\delta,c \in N^\ell_\delta \backslash M$ and $p 
= { \text{\rm tp\/}}_{\frak s}(c,M,N^\ell_\lambda) 
\in {\Cal S}^{\text{bs}}_{\frak s}(M)$ 
\ub{then} $p \pm M^\ell_i$ for some $i < i^*$
\sn
\item "{$(d)$}"  if $p \in {\Cal S}^{\text{bs}}(N^\ell_\alpha),\alpha <
\delta,\ell \in \{1,2\}$ and $p \pm M^\ell_i$ for some $i < i^*$ and $p$
is regular, \ub{then} for $\lambda$ ordinals $\beta \in
(\alpha,\delta)$ there is $c \in M^\ell_{\beta +1}$ such that
{\rm tp}$_{\frak s}(c,N^\ell_\beta,N^\ell_{\beta +1})$ is a nonforking
extension of $p$
\sn
\item "{$(e)$}"  $f_0$ is an isomorphism from $N^1_0$ onto $N^2_0$
mapping $M^1_i$ onto $M^2_i$.
\ermn
\ub{Then} there is an isomorphism from $N^1_\delta$ onto $N^2_\delta$
extending $f_0$.
\endproclaim
\bigskip

\demo{Proof}  Hence and forth, as usual.
\enddemo
\bigskip

\proclaim{\stag{705-6p.fr10} Claim}  1) If $\langle M_i:i < \alpha \rangle$
is independent over $M$ and $a_i \in M_i$, {\rm tp}$_{\frak
s}(a_i,M,M_i) \in {\Cal S}^{\text{bs}}_{\frak s}(M)$ \ub{then}
$\{a_i:i < \alpha\}$ is independent over $M$. \nl
2) If above $(M,M_i,\bold J_i) \in K^{3,\text{bs}}_{\frak s},i <
\alpha \Rightarrow M_i \le_{\frak s} N$ \ub{then} $\cup\{\bold J_i:i <
\alpha\}$ is independent in $(M,N)$.
\endproclaim
\newpage

\head {\S9 Between cardinals, Nonsplitting and getting fullness} \endhead  \resetall \sectno=9
 \spuriousreset
\bigskip

Our aim is to get full $\lambda$-good$^+$-systems.
Fullness seems naturally desirable (being closed to superstability)
and help in proving the existence of enough regular types.
\definition{\stag{705-6a.0} Hypothesis}  $1) {\frak s}$ is a successful
$\lambda$-good$^+$ system.  
\enddefinition
\bigskip

\definition{\stag{705-9f.3} Definition}  1) For $\lambda < \mu$ we define
${\frak s}$ is a $[\lambda,\mu)$-frame as in \chaptercite{600}, except
that ${\frak K}_{\frak s}$ is a $[\lambda,\mu)$-a.e.c., i.e., ${\frak
K}_{\frak s} = K^{{\frak s}'}_\lambda \restriction \cup \{K^{{\frak
s}'}_\kappa:\kappa \in [\lambda,\mu)\}$.  Let $\lambda_{\frak s},
\mu_{\frak s}$ be $\lambda,\mu$, respectively. \nl
2) We define ``a good $[\lambda,\mu)$-frame", $K^{3,\text{bs}}_{\frak
s},K^{3,\text{pr}}_{\frak s},K^{3,\text{uq}}_{\frak s},
K^{3,\text{qr}}_{\frak s},K^{3,\text{vq}}_{\frak s}$ similarly. \nl
3) For a $\lambda$-frame ${\frak s}$ and $\mu > \lambda$ we define
${\frak s}[\lambda,\mu)$ naturally; so ${\frak K}_{{\frak
s}[\lambda,\mu)} = K^{\frak s} \restriction \cup 
\{K^{\frak s}_\kappa:\kappa \in [\lambda,\mu)\}$.
\enddefinition
\bigskip

\proclaim{\stag{705-9f.4} Claim}  The claims on good $\lambda$-frames hold
for $[\lambda,\mu)$-frames with some obvious changes.
\endproclaim
\bn
\centerline {$* \qquad * \qquad *$}
\bn
Not central but we may note
\definition{\stag{705-stg.13} Definition}  1) We say that ${\frak s}$ is
type-full \ub{if} ${\Cal S}^{\text{bs}}_{\frak s} = 
{\Cal S}^{\text{na}}_{\frak s}$. \nl
2) Let $\nonforkin{a}{b}_{M}^{N}$ means that 
tp$(a,M,N)$, tp$(b,M,N) \in {\Cal S}^{\text{bs}}_{\frak s}(N)$ and
for some $M_\ell (\ell < 3)$
we have $M_0 = M,N \le_{\frak s} M_2,M_0 \le_{\frak s} M_1 \le_{\frak s}
M_2,a \in M_1$ and tp$(b,M_1,M_2)$ does not fork over $M_0$. \nl
3) We may allow not to distinguish types of elements and of finite
tuples, 
so we use ${\Cal S}_{\frak s}(M) = \cup\{{\Cal S}^m_{\frak s}(M):m <
\omega\}$ we call such ${\frak s}$ of $(< \omega)$-type, then
let $ab$ be the concatanation (this makes no real difference). \nl
4)  We say that ${\frak s}$ is type-closed \ub{if} ${\frak s}$ is of
$(< \omega)$-types and 
tp$(a_\ell,M,N) \in {\Cal S}^{\text{bs}}(M)$ for $\ell=1,2$ and
$\nonforkin{a_1}{a_2}_{M}^{N}$ implies tp$(a_1a_2,M,N) \in 
{\Cal S}^{\text{bs}}_{\frak s}(M)$.  [See Claim x].
\enddefinition
\bigskip

\definition{\stag{705-stg.14} Definition}  For a $\lambda$-good frame ${\frak s}$
(not necessarily satisfying Hypothesis \scite{705-5.0}??)
we define a frame ${\frak t} = {\frak s}^{tc}$:

$$
\lambda_{\frak t} =  
\lambda_{\frak s},{\frak K}_{\frak t} = {\frak K}_{\frak s}
$$

$$
\align
{\Cal S}^{\text{bs}}_{\frak t} = 
\bigl\{ \text{tp}(a_0 \ldots a_{n-1},M,N):&M
\le_{\frak s} N,a_\ell \in N \backslash M \text{ for } \ell < n \\
  &\text{and } \nonforkin{}{}_{M}^{N} \{a_\ell:\ell < n\}, 
\text{ i.e. we can find } M_\ell \text{ satisfying} \\ 
 &M = M_0 \le_{\frak s} M_1 \le_{\frak s} \ldots \le_{\frak s} M_n = N
\text{ and} \\
  &a_\ell \in M_{\ell +1},\text{ tp}(a_\ell,M_\ell,M_{\ell +1}) \in
{\Cal S}^{\text{bs}}_{\frak s}(M_\ell) \\
  &\text{does not fork over } M_0 \text{ for } \ell < n \bigr\}
\endalign
$$

$$
\align
\nonfork{}{}_{} = \bigl\{(M_0,M_1,\bar a,M_3):
&\text{ for some } n \text{ and } 
\langle M^*_\ell:\ell \le n \rangle \text{ we have} \\
  &a = a_0 \ldots a_{n-1},M_3 \le_{\frak s} M^*_n,M_1 = M^*_0
\le_{\frak s} M^*_1 \le_{\frak s} \ldots \le_{\frak s} M^*_n \\
  &\text{ and tp}(a_\ell,M^*_\ell,M^*_{\ell +1}) \in
{\Cal S}^{\text{bs}}_{\frak s}(M_\ell) \\
  &\text{does not fork over } M_0 \bigr\}.
\endalign
$$
\enddefinition
\bigskip

\proclaim{\stag{705-stg.15} Claim}  1) If ${\frak s}$ is a good 
$\lambda$-frame as in Hypothesis \scite{705-5.0}[?], \ub{then} 
${\frak s}^{\text{tc}}$ is a good $\lambda$-frame as in Hypothesis
\scite{705-5.0} and deal with $(< \omega)$-type. \nl
2) If ${\frak s}$ is a good$^+ \lambda$-frame, 
\ub{then} ${\frak s}^{\text{tc}}$ is a $\lambda$-good$^+$ frame. \nl
3) Primes, $K^{3,\text{vq}}_{\frak s}$ dense, etc., are lifted [FILL!].
\endproclaim
\bigskip

\demo{Proof}  Saharon:  We may wonder (see answer later).
\enddemo
\bn
\centerline{ $* \qquad * \qquad *$}
\bigskip

\definition{\stag{705-gr.2} Definition}  1) We say ${\frak K}^1$ is
situated above ${\frak K}^0$ \ub{if} for 
some $\lambda_0 \le \lambda_1$ we have:
${\frak K}^1$ is a $\lambda_1$-a.e.c., with amalgamation,
${\frak K}^0$ is an $\lambda_0$-a.e.c. with amalgamation 
and letting ${\frak K}^{0,\text{up}}$ be the lifting (from
\sectioncite[\S1]{600}) we have $K^1_{\lambda_1} \subseteq
K^{0,\text{up}},
\le_{{\frak K}^1_\lambda} = \le_{{\frak K}^{0,\text{up}}} \restriction
K^1_\lambda$.
We may write $\lambda_\ell = \lambda({\frak K}^\ell)$.  We may treat
tp$_{{\frak K}^1}(a,M,N)$ as 
tp$_{{\frak K}^0}(a,M,N)$ when no confusion arises, e.g. $p
\restriction N, N \in K_{\frak t}$.  \nl
2) We say ${\frak K}^1_\lambda$ 
is weakly ${\frak K}^0$-local (or weakly local above ${\frak K}^0$) 
\ub{if}:
\mr
\item "{$(a)$}"  ${\frak K}^1_\lambda$ is situated above 
${\frak K}^0$
\sn
\item "{$(b)$}"  if $M \in {\frak K}^1_{\lambda_1}$ and $p,q \in {\Cal
S}_{{\frak K}^1_\lambda}(M)$, \ub{then}
{\roster
\itemitem{ $\boxtimes$ }  $p = q \Leftrightarrow (\forall N)[N
\le_{{\frak K}^{0,\text{up}}} M 
\and N \in K_{\lambda({\frak K}^0)} \rightarrow p
\restriction N = q \restriction N (\in {\Cal S}_{{\frak K}^0}(N))]$
which means if $p = \text{ tp}_{{\frak K}^1_\lambda}(a,M,M_1)$ and $q
= \text{ tp}_{{\frak K}^1}(b,M,M_2)$ \ub{then}: $p=q$ iff for every $N
\in {\frak K}^0$ satisfying $N \le_{{\frak K}^{0,\text{up}}} M$ we have 
tp$_{{\frak K}^{0,\text{up}}}
(a,N,M_1) = \text{ tp}_{{\frak K}^{0,\text{up}}}(b,N,M_2)$
recalling that this means that $N \le_{{\frak K}^{0,\text{up}}} N_\ell \le
M_\ell,a \in N_1,b \in N_2 \Rightarrow \text{ tp}_{{\frak K}^0}
(a,N,N_1) = \text{ tp}_{{\frak K}^0}(b,N,N_1)$.
\endroster}
\ermn
2A) ${\frak K}^1_\lambda$ is basically weak ${\frak K}^0$-local if
(a),(b) above but in (b), $p,q \in {\Cal S}^{\text{bs}}(M)$. \nl
3) We say ${\frak K}^1$ is ${\frak K}^0$-local (or local above
${\frak K}^0$) \ub{if}
\mr
\item "{$(\alpha)$}"  ${\frak K}^1$ is situated above ${\frak K}^0$
\sn  
\item "{$(\beta)$}"  ${\frak K}^1$ is weakly ${\frak K}^0$-local
\sn
\item "{$(\gamma)$}"   if $M \in {\frak K}^1$ and $p = \text{
tp}_{{\frak K}^1}(a,M,M_1) \in {\Cal S}_{{\frak
K}^1}(M)$, \ub{then} for some $A \subseteq M,|A| \le 
\lambda({\frak K}^0)$ the type $p$ does not split over $A$ which means
that: if $A \subseteq N \in {\frak K}^0,N_\ell 
\le_{{\frak K}^{0,\text{up}}} M$ for $\ell=1,2$ and $f$ is an
isomorphism from $N_1$ onto $N_2$ over $A$,
\ub{then} tp$_{{\frak K}^0}(a,N_2,M_1) = f(\text{tp}_{{\frak
K}^0}(a,N_1,M_1)$ (by Definition \scite{705-gr.1} this is equivalent to:
if $A \subseteq N \in {\frak K}^0,N \le_{{\frak K}^{0,\text{up}}} M$
then $(\text{tp}_{{\frak K}^{0,\text{uq}}}(a,N,M_1)$ does not split over $A$.
\ermn
4) We say ${\frak s}$ is situated above ${\frak t}$ \ub{if}:
\mr
\item "{$(a)$}"  ${\frak s},{\frak t}$ are good frames
\sn
\item "{$(b)$}"  $\lambda_{\frak t} \le \lambda_{\frak s}$
\sn
\item "{$(c)$}"  ${\frak K}_{\frak s}$ is situated 
above ${\frak K}_{\frak t}$
\sn
\item "{$(d)$}"  if $M \in {\frak K}_{\frak s}$ and 
$p \in {\Cal S}(M)$, \ub{then}:
$p \in {\Cal S}^{\text{bs}}(M)$ \ub{iff} for 
every large enough $N \le_{{\frak K}[{\frak t}]} M$ of cardinality
$\lambda_{\frak t}$ we have 
$p \restriction N \in {\Cal S}^{\text{bs}}_{\frak s}(N)$
\sn
\item "{$(e)$}"  similarly ``$p \in {\Cal S}^{\text{bs}}_{\frak s}(M_1)$
does not fork over $M_0 \le_{\frak s} M_1$" \ub{iff} there are
$N_0 \le_{{\frak K}_{\frak t}} N_1$ 
(so $N_\ell \in {\frak K}_{\frak t}$) such that 
$N_\ell \le_{{\frak K}[{\frak t}]} M_\ell$:
{\roster
\itemitem{ $(*)$ }  if $N_1 \le_{{\frak K}_{\frak t}} N^+_1 \le_{{\frak
K}[{\frak t}]} M_1$ then $p \restriction N^+_1 \in {\Cal
S}^{\text{bs}}_{\frak t}(N^+_1)$ does not fork over $N_0$.
\endroster}  
\ermn
5) We say ${\frak s}$ is weakly ${\frak t}$-local (or local above 
${\frak t}$)
\mr
\item "{$(\alpha)$}"  ${\frak s}$ is situated above ${\frak t}$
\sn
\item "{$(\beta)$}"  ${\frak K}_{\frak s}$ is weakly ${\frak K}
^{\frak t}$-local, i.e., clause (b) of (2) (which deals with not
necessarily basic types); but even for basic types this (b) of (2) is
not implies by ``${\frak s}$ situated above ${\frak t}$".
\ermn
5A) We say ${\frak s}$ is ${\frak t}$-local (or local above ${\frak t}$) if
\mr
\item "{$(\alpha)$}"  ${\frak s}$ is situated above ${\frak t}$
\sn
\item "{$(\beta)$}"  ${\frak K}_{\frak s}$ is ${\frak K}^{\frak
t}$-local (see (3)).
\ermn
6) In (5) we add basically if in (b) of (c) we demand this. \nl
6A) In (5A) we add the adjective ``basically" if we weaken clauses
$(\beta), (\gamma)$ of part (3) to $(\beta)^{\text{bs}} + 
(\gamma)^{\text{bs}}$ respectively; that is ${\frak s}$ is basically 
${\frak t}$-local \ub{if}:
\mr
\item "{$(\alpha)$}"  ${\frak s}$ is situated
\sn
\item "{$(\beta)^{\text{bs}}$}"  if $M \in {\frak K}_{\frak s}$ and 
$p,q \in {\Cal S}^{\text{bs}}_{\frak s}(M)$, 
\ub{then} $p = q \Leftrightarrow (\forall N)[N \le_{{\frak K}[{\frak
t}]} M \and N \in {\frak K}_{\frak t} \rightarrow
p \restriction N = q \restriction N (\in {\Cal S}_{\frak t}(N)]$
\sn
\item "{$(\gamma)^{\text{bs}}$}"  if $M \in {\frak K}_{\frak s},p \in {\Cal
S}^{\text{bs}}_{\frak s}(M)$, 
\ub{then} for some $A \subseteq M$ we have $|A|
\le \lambda_{\frak t}$ and $p$ does not split over $A$.
\ermn
7) We say ${\frak s}$ is a strongly ${\frak t}$-local \ub{if}:
\mr
\item "{$(\alpha)$}"  ${\frak s}$ is basically ${\frak t}$-local
\sn
\item "{$(\beta)$}"  ${\frak s},{\frak t}$ are weakly successful
\sn
\item "{$(\gamma)$}"  if $M \le_{\frak s} N$, \ub{then} for almost
every $X \in [N]^{\le \lambda[{\frak t}]}$ the type tp$_{K[{\frak
t}]}(X,M,N)$ does not split over $X \cap N$, see Definition xxx
\sn
\item "{$(\delta)$}"  for $M_\ell \in K_{\frak s}$ for $\ell < 4$, we
have NF$_s(M_0,M_1,M_2,M_3)$ \ub{iff} for almost all $Y \in [M_0 \cup
M_1 \cup M_2 \cup M_3]^{\lambda_{\frak t}}$, NF$_{\frak t}(M_0
\restriction Y,M_1 \restriction Y,M_2 \restriction Y,
M_3 \restriction Y)$
\sn
\item "{$(\delta)'$}"  NF have 
\footnote{is it reasonable to demand this instead for
$K^{3,\text{uq}}_{\frak s}$?  Then for NF it follows.} 
reflection from ${\frak s}$ to ${\frak
t}$ and lifting from ${\frak t}$ to ${\frak s}$ and x.x.
\ermn
8) We say ${\frak s}$ is super ${\frak t}$-local \ub{if}:
\mr
\item "{$(\alpha)$}"  ${\frak s},{\frak t}$ are weakly successful
frames
\sn 
\item "{$(\beta)$}"  ${\frak s},{\frak t}$ have primes 
\footnote{To demand every $(M_0,M_1,a) \in K^{3,\text{uq}}_{\frak s}$
is reflected incomparable.  To demand for a dense family is enough and weak.}
and $\bot = \underset{\text{wk}} {}\to \bot$ (for both)
\sn
\item "{$(\gamma)$}"   ${\frak s}$ is strongly ${\frak t}$-local
\sn
\item "{$(\delta)$}"  if $(M_0,M_1,a) \in K^{3,\text{pr}}_{\frak s}$ 
then for
almost all $Y \in [M_1]^{\lambda({\frak t})}$ we have \nl
$(M_0 \restriction Y,M_1 \restriction Y,a) \in K^{3,\text{pr}}_{\frak t}$.
\ermn
8A) We say that ${\frak s}$ is super-${\frak t}$-local \ub{if}:
\mr
\item "{$(\alpha)$}"   ${\frak s},{\frak t}$ are weakly successful
frames
\sn
\item "{$(\beta)$}"   ${\frak s}$ is strongly ${\frak t}$-local
\sn
\item "{$(\gamma)$}"   if $M_0 \in {\frak K}_{\frak s},p \in {\Cal
S}^{\text{bs}}(M)$ then there are $M_1,a$ such that
{\roster
\itemitem{ $(i)$ }   $(M_0,M_1,a) \in K^{3,\text{uq}}_{\frak s}$
\sn
\itemitem{ $(ii)$ }   tp$(a,M_0,M_1) = p$
\sn
\itemitem{ $(iii)$ }  for almost all $Y \in [M_1]^{\lambda({\frak
t})}$ we have $(M_0 \restriction Y,M_1 \restriction Y,a) \in
K^{3,\text{uq}}_{\frak t}$.
\endroster}
\ermn
9) In any of the above we add ``fully" if ${\frak K}^1_{\lambda_1} =
{\frak K}^{0,\text{up}}_{\lambda_1}$ (for a.e.c.) or ${\frak K}_{\frak s} =
(K[{\frak t}])_{\lambda({\frak s})}$ (for frames).
\enddefinition
\bigskip

\proclaim{\stag{705-gr.3} Claim}   Assume ${\frak s} = {\frak t}^+$ is
weakly successful and ${\frak t}$ is $\lambda_{\frak t}$-good$^+$ and
successful. \nl
1) If $M \in {\frak K}_{\frak s}$ and $p \in {\Cal S}_{\frak s}(M)$, \ub{then}
$p$ does not split over some $A \subseteq M$ of cardinality
$\lambda_{\frak t}$ [even in a stronger sense: for $N_\ell \le_{{\frak
K}[{\frak t}]} M,N_\ell$ in any cardinality]. \nl
2) ${\frak s}$ is strongly ${\frak t}$-local. \nl
3) [??] If $M \in {\frak K}_s,A \subseteq M,|A| < \lambda_{\frak s}$ and $p
\in {\Cal S}_{\frak s}(M)$, \ub{then}

$$
\align
p \text{ does not split over } A &\text{ iff } (\forall N \in {\frak
K}_{\frak t}) \\
  &[N \le_{{\frak K}^{\frak t}} M \and A \subseteq N \rightarrow p
\restriction N \text{ does not } \lambda_{\frak t} \text{-split over } A]
\endalign
$$
\mn
4) ${\frak s}$ is super ${\frak t}$-local.
\endproclaim
\bigskip

\demo{Proof}  Essentially from \marginbf{!!}{\cprefix{600}.\scite{600-stg.30t}}.
\enddemo
\bn
Some of the obvious properties are
\proclaim{\stag{705-gr.4} Claim}  1) If ${\frak s}_2$ is ${\frak
s}_1$-local and ${\frak s}_1$ is ${\frak s}_0$-local, \ub{then}
${\frak s}_2$ is ${\frak s}_0$-local. \nl
2) If ${\frak s}_2$ is strongly ${\frak s}_1$-local, ${\frak s}_1$ has
$\chi$-nonsplitting, \ub{then} ${\frak s}_2$ has
$\chi$-nonsplitting. \nl
3) If ${\frak s}_2$ is strongly ${\frak s}_1$-local and ${\frak s}_1$
is strongly ${\frak s}_0$-local, \ub{then} ${\frak s}_2$ is strongly
${\frak s}_0$-local. 
\nl
4) Similar results for the ``basic" version. \nl
5) Similar results for ``${\frak K}^\ell$ a
$\lambda_\ell$-a.e.c.". \nl
6) Similar results of ``${\frak s}$ is super ${\frak t}$-local".
\endproclaim
\bigskip

\demo{Proof}  Easy.
\enddemo
\bn
We now investigate the reflection of the following properties and
their negation: independent, orthogonal $K^{3,\text{qr}},
K^{3,\text{vq}}$.
\proclaim{\stag{705-gr.4A} Claim}  1) Assume 
${\frak K}^2$ is situated above ${\frak K}^1$.
If $M_1 \le_{{\frak K}^2} M_2$ \ub{then} for almost all $Y \in
[M_2]^{\lambda[{\frak t}]}$ we have $M_1 \restriction Y \le_{\frak t}
M_2$; similarly for $M_1 \nleq_{{\frak K}^2} M_2 \and \{M_1,M_2\} \in
{\frak K}^2$. \nl
2) [${\frak s}$ is situated above ${\frak t}$]
\mr
\item "{$(a)$}"  $p \in {\Cal S}^{\text{bs}}_{\frak s}(M_1)$ does not
fork over $M_0 \le_{\frak s} M_1$ \ub{iff} for almost every $Y \in
[M_1]^{\lambda({\frak t})}$ we have: $p \restriction Y \in {\Cal
S}^{\text{bs}}_{\frak t}(M_1 \restriction Y)$
\sn
\item "{$(b)$}"  if $M_0 
\le_{\frak s} M_1 \le_{\frak s} M_2$ and $\bold J \subseteq
\bold I_{M_1,M_2}$ in independent in $(M_0,M_1,M_2)$ for ${\frak s}$,
\ub{then} for almost all $Y \in [M_2]^{\lambda({\frak t})}$ we have:
$\bold J \cap Y \subseteq \bold I_{M_1 \restriction Y,M_2 \restriction
Y}$ is independent in $(M_0 \restriction Y,M_1 \restriction Y,M_2
\restriction Y)$ for ${\frak t}$,
\sn
\item "{$(c)$}"  if $p,q \in {\Cal S}^{\text{bs}}_{\frak s}(M)$, 
\ub{then} for almost all $Y \in [M]^{\lambda({\frak t})}$ 
we have $p \underset{\text{wk}} {}\to \pm_{\frak s}
q \Rightarrow p \restriction (M \restriction Y) 
\underset{\text{wk}_{\frak t}} {}\to \pm 
q \restriction (M \restriction Y)$
\sn
\item "{$(d)$}"  if $M_\ell \le_{\frak s} M,p_\ell \in {\Cal
S}^{\text{bs}}_s(M_\ell)$ then $p_1 \| p_2$ \ub{iff} for almost every
$Y \in [M]^{\lambda({\frak t})},M_\ell \restriction Y \le_{\frak t} M$
and $p_1 \restriction (M_1 \restriction Y)\|p_2 \restriction (M_2
\restriction Y)$.
\endroster
\endproclaim
\bigskip

\proclaim{\stag{705-gr.4Q} Claim}  Assume ${\frak s}$ is super ${\frak
t}$-local. \nl
1) If $M_0 \le_{\frak s} M_1 \le_{\frak s} M_2$ and $\bold J \subseteq
\bold I^{\frak s}_{M_1,M_2}$ \ub{then} $\bold J$ is independent in
$(M_0,M_1,M_2)$ \ub{iff} for almost every $Y \in [M_2]^{\lambda({\frak
t})}$ the set $\bold J \cap Y$ is independent in $(M_0 \restriction
Y,M_1 \restriction Y,M_2 \restriction Y)$. \nl
2) If $p,q \in {\Cal S}^{\text{bs}}_s(M)$ \ub{then} $p
{\underset{\text{wk}}_{\frak s} {}\to \pm q}$ iff for almost every $Y
\in [M]^{\lambda({\frak t})},p \restriction (M \restriction Y)
{\underset{\text{wk}}_{\frak t} {}\to \pm q} \restriction Y)$. \nl
3) If $(M_0,M_1,\bold J) \in K^{3,\text{qr}}_{\frak s}$ 
\ub{then} for almost all
$Y \in [M_0]^{\lambda({\frak t})}$ we have $(M_0 \restriction Y,M_1
\restriction Y,\bold J \cap Y) \in K^{3,\text{\text{qr}}}_{\frak t}$, 
similarly for $\notin K^{3,\text{qr}}_{\frak s}$. \nl
4) Like (3) for $K^{3,\text{vq}}_{\frak s}$. \nl
5) If $M_\ell \in K_{\frak s}$ for $\ell < 4$ then: 
{\rm NF}$_{\frak s}
(M_0,M_1,M_2,M_3)$ \ub{iff} for almost all $Y \in [\dbcu_{\ell < 4}
M_\ell]^{\lambda){\frak t})}$ we have {\rm NF}$_{\frak t}(M_0 \restriction
Y,M_1 \restriction Y,M_2 \restriction Y,M_3 \restriction Y)$;
similarly for the negation. \nl
6) If $M_0 \le_{{\frak K}[{\frak t}]} M_2 \in K_{\frak s},M_0
\le_{\frak t} M_1,M_1 \cap M_2 = M_0$ \ub{then} we can find $M_3$ such
that: $M_2 \le_{\frak s} M_3,M_1 \le_{{\frak K}[{\frak t}]} M_3$ and
for almost all $Y \in [M_3]^{\lambda({\frak t})}$ we have {\rm NF}
$_{\frak t}(M_0 \restriction Y,M_1 \restriction Y,M_2 \restriction Y,M_3
\restriction Y)$. We express this as {\rm NF}$_{{\frak t},{\frak s}}
(M_0,M_1,M_2,M_3)$.  Similar we define ``$p \in {\Cal S}
^{\text{bs}}_{\frak s}(M_2)$ does not fork over $M_0$ for $({\frak t},
{\frak s})$. \nl
7) If $M_1 \le_{\frak s} M_2$ and $\bar{\bold a} \in
{}^{\lambda({\frak t})}(M_2)$, \ub{then} for some $A \in
[M_1]^{\lambda({\frak t})}$ we have {\rm tp}$_{\frak s}(\bar{\bold
a},M_1,M_2)$ does not $\lambda_{\frak s}$-split over $A$. \nl
8) The parallel of \scite{705-6p.10} (fill).
\endproclaim
\bigskip

\demo{Proof}  Fill (particularly (7)).

For some club ${\Cal C}_\ell \subseteq [M_\ell]^{\lambda_{\frak t}}$ 
we have $Y \in {\Cal C}_\ell \Rightarrow
M_\ell \restriction Y \le_{{\frak K}[{\frak t}]} M_\ell$.

Now \wilog \, $Y \in {\Cal C}_2 \Rightarrow Y \cap M_1 \in {\Cal C}_1$
so now for every $Y \in {\Cal C}_2$ we have $M_1 \restriction Y = M_1
\restriction (Y \cap M_1) \le_{{\frak K}[{\frak t}]} M_1 \le_{{\frak
K}[{\frak t}]} M_2$ and $M_2 \restriction Y \le_{{\frak K}[{\frak t}]}
M_2$ and $M_1 \restriction Y \subseteq M_2 \restriction Y$ hence by
axiom V of a.e.c. we have $M_1 \restriction Y \le_{[{\frak K}[{\frak
t}]} M_2 \restriction Y$ which means that $M_1 \restriction Y
\le_{\frak t} M_2 \restriction Y$.  \hfill$\square_{\scite{705-gr.4Q}}$\margincite{705-gr.4Q}
\enddemo
\bigskip

\proclaim{\stag{705-gr.4B} Claim}  1) We say $\bar{\frak s} = \langle
{\frak s}_\theta:\lambda_0 \le \theta \le \lambda_1 \rangle$ is a
super local sequence \ub{if}:
\mr
\item "{$(a)$}"  each ${\frak s}_\theta$ is a weakly successful good
$\theta$-frame
\sn
\item "{$(b)$}"  if $\lambda_0 \le \theta_0 < \theta_1 <
\lambda$ \ub{then} ${\frak s}_{\theta_1}$ is fully 
super ${\frak s}_{\theta_1}$-local.
\ermn
2) For such $\bar{\frak s}$ let: ${\frak K}_{\bar{\frak s}} = ({\frak
K}[\bar s_{\lambda_0}])_{[\lambda_0,\lambda^+_1)}$ and we write
$\le_{\frak s}$ for $\le_{K_{\bar{\frak s}}}$.  Let 
{\rm NF}$_{\bar{\frak s}}
(M_0,M_1,M_2,M_3)$ \ub{iff} $M_0 \le_{\bar{\frak s}} M_\ell
\le_{\bar{\frak s}} M_3$ for $\ell =1,2$ and for almost every $Y \in
[M_3]^{\lambda_0}$ we have {\rm NF}$_{{\frak s}_{\lambda_0}}(M_0
\restriction Y,M_1 \restriction Y,
M_2 \restriction Y,M_3 \restriction
Y)$.
\endproclaim
\bigskip

\proclaim{\stag{705-gr.4C} Claim}  1) If $\bar{\frak s} = \langle {\frak
s}_\theta:\theta \in [\lambda_0,\lambda_1) \rangle$ is a super local
sequence and $\lambda_1$ is a limit cardinal \ub{then} for one and only one
${\frak s}_{\lambda_1}$ the sequence $\bar{\frak s}' = \langle {\frak
s}_\theta:\theta \in [\lambda_0,\lambda^+_1) \rangle$ is a super local
sequence. \nl
2) If $\bar{\frak s} = \langle {\frak s}_\theta:\theta \in
[\lambda_0,\lambda^+_1) \rangle$ is a super local sequence and ${\frak
s}_{\lambda_1}$ is successful, and we define ${\frak s}_{\lambda^+_1}
= {\frak s}^+_{\lambda_1}$ \ub{then} the sequence ${\frak s}' = \langle
{\frak s}_\theta:\theta \in [\lambda_0,\lambda^{++}_1) \rangle$ is a
super local sequence (but not necessarily full!).
\endproclaim
\bigskip

\demo{Proof}  FILL!  \hfill$\square_{\scite{705-gr.4C}}$\margincite{705-gr.4C}
\enddemo
\bn
Now we return to trying to deal with all types in ${\Cal S}(M)$, i.e. fullness.
\definition{\stag{705-6a.1} Definition}  [${\frak s}$ is a $\lambda$-frame] \nl
1) For $M \in K_{\frak s}$ let ${\Cal S}^{\text{na}}_{\frak s}(M) 
= \{\text{tp}(b,M,N):M \le_{\frak s} N$ and
$b \in N \backslash M\}$. \nl
2) The $\lambda$-frame ${\frak s}$ is type-full \ub{if} 
$M \in K_{\frak s} \Rightarrow 
{\Cal S}^{\text{bs}}_{\frak s}(M) = {\Cal S}^{\text{na}}_{\frak s}(M)$. 
\nl
3) [${\frak s}$ is a weakly successful good $\lambda$-frame] \nl
Let ${\frak s}^{\text{nf}} = {\frak s}(\text{nf})$ 
be the following frame (see below)
\mr
\item "{$(a)$}"  ${\frak K}_{{\frak s}(\text{nf})} = {\frak K}_{\frak s}$
\sn
\item "{$(b)$}"  ${\Cal S}^{\text{bs}}_{{\frak s}(\text{nf})}(M) 
= {\Cal S}^{\text{na}}_{\frak s}(M)$
\sn
\item "{$(c)$}"  $\nonfork{}{}_{{\frak s}(\text{nf})}(M_0,M_1,a,M_3)$ holds 
\ub{iff}
$M_0 \le_{\frak s} M_1 \le_{\frak s} M_3$ and 
$a \in M_3 \backslash M_1$ and
there are $M'_3,M_2$ such that $M_0 \le_{\frak s} M_2 \le_{\frak s} M'_3,
M_3 \le_{\frak s} 
M'_3$ and $a \in M_2$ and NF$_{\frak s}(M_0,M_1,M_2,M'_3)$.
\ermn
4) [${\frak s}$ as in (3)] \nl
Let ${\frak s}^{+\text{nsp}}$ be the $\lambda^+_{\frak s}$-frame
which we also denote by ${\frak s}(*)$ or ${\frak s}(+\text{nsp})$, 
defined by
\mr
\item "{$(a)$}"  ${\frak K}_{{\frak s}(*)} = {\frak K}_{{\frak s}(+)}$
\sn
\item "{$(b)$}"  ${\Cal S}^{\text{bs}}_{{\frak s}(*)}(M) 
= {\Cal S}^{\text{nsp}}_{{\frak s}(+)}
(M) =: \{p \in {\Cal S}^{\text{na}}_{{\frak s}(*)}(M):
\text{ for some } M_0
\le_{{\frak K}[{\frak s}]} M$ from $K^{\frak s}_\lambda,p$ does not 
$\lambda_{\frak s}$-split over $M_0$
(see Definition \scite{705-gr.1}(1) in our case as $M$ is $K^{\frak x}$-saturated
over $\lambda$, this means that every automorphism 
$g$ of $M$ over $M_1$ maps $p$ to itself), if
${\frak s}$ is successful this is ${\Cal S}^{\text{na}}_{s(*)}(M)$
\sn
\item "{$(c)$}"  $\nonfork{}{}_{{\frak s}(*)}(M_0,M_1,a,M_3)$ \ub{if}
$M_0 \le_{\frak s} M_1 \le_{\frak s} M_3,a \in M_3 \backslash M_1$ and
tp$_{{\frak s}(*)}(a,M_1,M_3)$ does not $\lambda_{\frak s}$-split over some
$N_0 \le_{{\frak K}[{\frak s}]} M_0,N_0 \in K_{\frak s}$.
\endroster
\enddefinition
\bigskip

\proclaim{\stag{705-9.11A} Claim}  Assume ${\frak s}$ is a successful
good$^+$-frame.  If $M \in K_{{\frak s}(+)}$ and $p,q \in 
{\Cal S}_{\frak s}(M)$
and $N \le_{{\frak K}[{\frak s}]} M \and N \in K_{\frak s} \Rightarrow
(p \restriction N) = (q \restriction N)$ \ub{then} $p = q$.
\endproclaim
\bigskip

\demo{Proof}  Similar to \scite{705-gr.6}.
\enddemo
\bigskip

\proclaim{\stag{705-6a.2} Lemma}   Assume that ${\frak s}$ is a successful
good$^+$ frame.  \ub{Then} the frame ${\frak s}(+{\text{\rm nsp\/}}) = 
{\frak s}(*)$ is a good$^+ \, \lambda^+_{\frak s}$-frame and ${\Cal
S}^{\text{bs}}_{{\frak s}(*)}(M) = {\Cal S}_{{\frak s}(*)}(M),{\frak s}(*)$
has primes and weak orthogonality and is equivalent to orthogonality
for it and ${\frak K}_{{\frak s}(*)}$ is categorical.
\endproclaim
\bigskip

\demo{Proof}  We have to check the axioms there.
\sn
\ub{Axioms}:  (A),(B),(C).

As ${\frak K}_{{\frak s}(*)} = 
{\frak K}_{{\frak s}(+)}$ this follows from \scite{705-stg.3B}(1).
\sn
\ub{Axiom}:  (D),(a),(b) by the Definition of 
${\Cal S}^{\text{bs}}_{{\frak s}(*)}$.
\sn
\ub{Axiom (D)(c)}:  Any 
$a \in N \backslash M$ is O.K. by \cite[\S6]{600},\marginbf{!!}{\cprefix{600}.\scite{600-nf.20t}}
(2) (use
representation, remember that every model $M \in K_{{\frak s}(*)}$ is
a saturated in ${\frak K}_{\frak s}$ above $\lambda_{\frak s}$ of
cardinality $\lambda^+_{\frak s} = \lambda_{{\frak s}(*)}$).  Hence we
get also $\bar{\Cal S}^{\text{bs}}_{{\frak s}(*)}(M) 
= {\Cal S}^{\text{na}}_{{\frak s}(*)}(M)$. 
\sn
\ub{Axiom (D)(d)}:  This holds as ${\frak K}_{{\frak s}(+)}$ is stable in 
$\lambda^+_{\frak s}$ by \scite{705-stg.3B}(2) + \cite[\S4]{600},\marginbf{!!}{\cprefix{600}.\scite{600-4a.1t}} but 
${\frak K}_{{\frak s}(*)} = {\frak K}_{{\frak s}(+)}$, 
alternately use \cite[\S6]{600},\marginbf{!!}{\cprefix{600}.\scite{600-nf.20t}}(1)(?).   
\sn
\ub{Axiom (E)(a)}:  By the definitions.
\sn
\ub{Axiom (E)(b)}:  [monotonicity].

So Assume $M_0 \le_{{\frak s}(*)} M'_0 \le_{{\frak s}(*)} M'_1
\le_{{\frak s}(*)} M_1 \le_{{\frak s}(*)} M_3 \le_{{\frak s}(*)}
M'_3$ and $\nonfork{}{}_{{\frak s}(*)} (M_0,M_1,a,M_3)$ so it is witnessed
by some $N_0 \le_{K[{\frak s}]} M_0$ with $N_0 \in K_{\frak s}$.

Now the same $N_0$ witness $\nonfork{}{}_{{\frak s}(*)} (M'_0,M'_1,a,M_3)$.
The other statement $(\nonfork{}{}_{{\frak s}(*)} (M_0,M_1,a,M_3)
\Leftrightarrow \nonfork{}{}_{{\frak s}(*)}(M_0,M_1,a,M'_3))$ 
is immediate by
tp$_{{\frak s}(*)}(a,M'_1,M'_3) = \text{tp}_{{\frak s}(*)}(a,M'_1,M_3)$.
\sn
\ub{Axiom (E)(c)}:  (local character).

So assume that $\langle M_i:i \le \delta +1 \rangle$ is 
${\frak s}(*)$-increasing continuous, $\delta < (\lambda_{{\frak
s}(*)})^+ = \lambda^{++},c \in
M_{\delta +1} \backslash M_\delta$, and assume toward contradiction
that tp$_{{\frak s}(*)}(c,M_\delta,M_{\delta +1}) \in 
{\Cal S}^{\text{bs}}_{{\frak s}(*)}(M_\delta)$ is a counterexample. 
Without loss of generality $\delta = \text{ cf}(\delta)$, so $\delta
\le \lambda^+$.   Let
$\bar M^i = \langle M^i_\alpha:\alpha < \lambda^+ \rangle$ be a
$\le_{\frak s}$-representation 
of $M_i,E$ a thin enough club of $\lambda^+$,
so e.g. 

$$
\alpha \in E \Rightarrow c \in M^{\delta +1}_\alpha
$$
\mn
and

$$
\alpha \in E \and \alpha < \beta \in E \and i \le \delta
\Rightarrow \text{ tp}_{\frak s}(c,M^i_\beta,M^{\delta +1}_\alpha)
\text{ does not } \lambda\text{-split over } M^i_\alpha
$$
\mn
and

$$
\alpha < \beta \in E \and i < \delta \Rightarrow \text{ tp}_{\frak s}
(c,M^\delta_\beta,M^{\delta +1}_\beta) \text{ does } 
\lambda\text{-split over } M^i_\alpha.
$$
\mn
Choose $\varepsilon_i \in E$ for $i \le \delta$, increasing continuous, so
$\langle M^i_{\varepsilon_i}:
i \le \delta \rangle$ is $<_{\frak s}$-increasing
continuous, each $M^i_{\varepsilon_i}$ is $(\lambda,*)$-brimmed for
${\frak s}$ and $i < j \le \delta \Rightarrow M^j_{\varepsilon_j}$ is
$(\lambda,*)$-brimmed over $M^i_{\varepsilon_i}$ for ${\frak s}$.
If $\delta < \lambda^+$, by Subclaim \scite{705-gr.6}, 
for some $i < \delta$, tp$_{\frak s}(c,
M^\delta_{\varepsilon_\delta},M^{\delta +1}_{\varepsilon_{\delta +1}})$ does
not $\lambda$-split over $M^i_{\varepsilon_i}$ for $K_{\frak s}$,
contradiction to the choice of $E$ above (and obvious monotonicity of
nonsplitting).  If $\delta = \lambda^+$, use Fodor's lemma.
\bn
\ub{Axiom (E)(d)}: [transitivity]

Assume
\mr
\item "{$(\alpha)$}"  $M_1 \le_{{\frak s}(*)} M_2 \le_{{\frak s}(*)}
M_3 \le_{{\frak s}(*)} M_4$
\sn
\item "{$(\beta)$}"  $a \in M_4 \backslash M_3$
\sn
\item "{$(\gamma)$}"  tp$_{{\frak s}(*)}(a,M_2,M_4),{\frak s}(*)$-does not
fork over $M_1$ and
\sn
\item "{$(\delta)$}"  tp$_{{\frak s}(*)}(a,M_3,M_4),{\frak s}(*)$-does not
fork over $M_2$.
\ermn
Let $\bar M^\ell = \langle M^\ell_\zeta:
\zeta < \lambda^+ \rangle$, for $\ell = 1,2,3,4$ 
be a $\le_{\frak s}$-representation of $M_\ell$ such that  
$a \in M^4_0$ and \wilog \, $\alpha < \beta < \lambda^+ \and 1 \le
\ell < m \le 4 \Rightarrow \text{ NF}_{\frak s}
(M^\ell_\alpha,M^m_\alpha,M^\ell_\beta,M^m_\beta)$ and for $\ell
=1,2$
\mr
\item "{$\boxtimes_\ell$}"  $M^\ell_0$ witness tp$_{{\frak
s}(*)}(a,M_{\ell +1},M_4),{\frak s}(*)$-does not fork over $M_\ell$.  
\ermn
Let $\bar{\bold a}_\ell$ list $M^\ell_0$ so
$\bar{\bold a}_\ell \in {}^\lambda(M_\ell)$.  Now assume $\bar{\bold b},
\bar{\bold c} \in {}^\lambda(M_3)$ are such that
\mr
\item "{$(\varepsilon)$}"  tp$_{{\frak K}[{\frak s}]}
(\bar{\bold b},M^0_1,M_4) =
\text{ tp}_{{\frak K}[{\frak s}]}(\bar{\bold c},M^0_1,M_4)$.
\ermn
As $M_2$ is ${\frak K}[{\frak s}]$-brimmed above $\lambda_{\frak s}$
we can find $\bar{\bold b}' \in {}^\lambda(M_2)$
such that
\mr
\item "{$(\zeta)$}"  tp$_{{\frak K}[{\frak s}]}(\bar{\bold b}',M^0_2,M_4) = 
\text{ tp}(\bar{\bold b},M^0_2,M_4)$ 
\ermn
similarly we can find $\bar{\bold c}' \in {}^\lambda(M_2)$ such that
\mr
\item "{$(\eta)$}"  tp$_{{\frak K}[{\frak s}]}(\bar{\bold c}',M^0_2,M_4) = 
\text{ tp}(\bar{\bold c},M^0_2,M_4)$.
\ermn
Chasing equalities $(\varepsilon) + (\zeta) + (\eta)$, as $M^0_1
\subseteq M^0_2$,  clearly
tp$_{{\frak K}[{\frak s}]}(\bar{\bold b}',M^0_1,M_4) = \text{ tp}
_{{\frak K}[{\frak s}]}(\bar{\bold c}',M^0_1,M_4)$, hence by clause
$(\gamma)$ more exactly by $\boxtimes_1$ we have
\mr
\item "{$(\theta)$}"  tp$_{{\frak K}[{\frak s}]}(\langle a \rangle \char 94
\bar{\bold b}',M^0_1,M_4) = \text{ tp}_{{\frak K}[{\frak s}]}
(\langle a \rangle \char 94 \bar{\bold c}',M^0_1,M_4)$.
\ermn
By clause $(\delta)$, i.e., by $\boxtimes_2$ and the statement $(\eta)$
we have
\mr
\item "{$(\iota)$}"  tp$_{{\frak K}[{\frak s}]}(\langle a \rangle \char 94
\bar{\bold c},M^0_2,M_4) = \text{ tp}_{{\frak K}[{\frak s}]}(\langle a \rangle
\char 94 \bar{\bold c}',M^0_2,M_4)$
\ermn
and similarly by $\boxtimes_2$ and $(\zeta)$ 
\mr
\item "{$(\kappa)$}"  tp$_{{\frak K}[{\frak s}]}(\langle a \rangle \char 94
\bar{\bold b},M^0_2,M_4) = \text{ tp}_{{\frak K}[{\frak s}]}(\langle a \rangle
\char 94 \bar{\bold b}',M^0_2,M_4)$.
\ermn 
By chasing the equalities $(\theta) + (\iota) + (\kappa)$ we get
tp$_{{\frak K}[{\frak s}]}(\langle a \rangle \char 94 \bar{\bold b}, 
M^0_0,M_4) = \text{ tp}_{{\frak K}[{\frak s}]}(\langle a \rangle
\char 94 \bar{\bold c},M^0_0,M_4)$ as required.
\sn
\ub{Axiom (E)(e)}:  [Unique nonforking extension].

So let $M_0 <_{{\frak s}(*)} M_1$ and $p,q \in 
{\Cal S}^{\text{bs}}_{{\frak s}(*)}(M_1)$ does 
not fork over $M_0$ and $p \restriction M_0 = q \restriction M_0$.
Let $M_1 <_{{\frak s}(*)} M_2$ and $a_1,a_2 \in M_2$ be such that
tp$_{{\frak s}(*)}(a_1,M_1,M_2) = p$ and 
tp$_{{\frak s}(*)}(a_2,M_1,M_2) = q$.
Let $\langle M_{\ell,\zeta}:\zeta < \lambda^+ \rangle$
be a $\le_{\frak s}$-representation of $M_\ell$ for $\ell =0,1,2$ with
$a_1,a_2 \in M_{2,0}$.  
By the assumption and the definition of ${\frak s}(*)$ for a club $E$ of
$\lambda^+$ we have, for $\zeta \in E,p \restriction
M_{0,\zeta} = q \restriction M_{0,\zeta} \in {\Cal S}_{\frak
s}(M_{0,\zeta})$ call it $r_\zeta$ and $p \restriction
M_{1,\zeta},q \restriction M_{1,\zeta}$ belong to ${\Cal S}_{\frak
s}(M_{1,\zeta})$ and does not $\lambda$-split over $M_{0,0}$ and they
extend $r_\zeta$; also for $\zeta < \xi$ in $E$
and $\ell < 2$ we have NF$_{\frak s}(M_{\ell,\zeta},M_{\ell
+1,\zeta},M_{\ell,\xi},M_{\ell +1,\xi})$ and for $\ell <
\xi,M_{\ell,\xi}$ is $(\lambda,*)$-brimmed over $M_{\ell,\zeta}$ for
${\frak s}$.  Now for $\zeta < \xi$ from $E$, let $g$ be an
isomorphism from $M_{1,\xi}$ onto $M_{2,\xi}$ over $M_{1,\zeta}$,
(clearly exists).  Let $\bar{\bold a}$ list $M_{1,\xi}$, and
$\bar{\bold b} = g(\bar{\bold a})$, so as for $\ell =1,2$, tp$_{\frak
s}(\bar a_\ell,M_{1,\xi},M_{2,\xi})$ does not $\lambda$-split over
$M_{0,0} \le_{\frak s} M_{1,\zeta}$ (by the hypothesis) there is 
an extension $g_\ell$ of
$g$ to an automorphism of $M_2$ such that $g_\ell(a_\ell) = a_\ell$,
so clearly $p \restriction M_{1,\xi} = q \restriction M_{1,\xi}$.  
By \scite{705-gr.3}(2) we get $p=q$ (recall 
that ${\frak K}_{{\frak s}(+)} = {\frak K}_{{\frak s}(*)}$).
\sn
\ub{Axiom (E)(f)}:  [Symmetry].

The proof relies on (E)(g) proved below.

We use freely ${\frak K}_{{\frak s}(*)} = {\frak K}_{{\frak s}(+)}$. \nl
By the symmetry (in the axiom) assume clause (b) there and we shall prove
clause (a).  So there are $M_4,M_5$ such that $M_3 \le_{{\frak s}(*)}
M_5$ and $M_0 \cup \{a_2\} \subseteq M_4 \le_{{\frak s}(*)} M_5$ and
tp$_{{\frak s}(*)}(a_1,M_4,M'_5)$ does not fork over $M_0$ for
${\frak s}(*)$.  By (E)(g) we can find $M_6,M_7 \in K_{{\frak s}(*)}$ such
that $M_4 \le_{{\frak s}(*)} M_6 \le_{{\frak s}(*)} M_7,M_5 
\le_{{\frak s}(*)} M_7$ and $M_6$ is $(\lambda^+_{\frak s},*)$-brimmed 
over $M_5$
for ${\frak s}^+$ (equivalently for ${\frak s}(*)$) and 
tp$_{{\frak K}_{{\frak s}(*)}}(a_1,M_6,M_7)$ does not fork over $M_4$ for
${\frak s}(*)$.  

By transitivity (i.e. clause (E)(d)) we know that
tp$_{s(*)}(a_1,M_6,M_7)$ does not fork over $M_0$ for ${\frak s}(*)$.  By
\scite{705-6.3} applied to ${\frak s}^+$ we can find $M_2,\bold J_2$ such that
$M_0 \le_{{\frak s}(*)} M_2 \le_{{\frak s}(*)} M_6,\bold J_2 \subseteq
\bold I^{s(*)}_{M_0,M_2}$ is finite and $(M_0,M_2,\bold J_2) \in
K^{3,\text{qr}}_{{\frak s}(+)}$ and $a_2 \in M_2$.  
Also by \scite{705-6.3} for ${\frak s}^+$ \wilog \,
we can find $M^*_1,a^*_1$ such that $M_0 \le_{{\frak s}(*)} M^*_1
\le_{{\frak s}(*)} M_7$ and $a^*_1 \in M^*_1$ such that
tp$_{{\frak s}(+)}(a^*_1,M_0,M^*_1) = \text{ tp}_{{\frak s}(+)}(a_1,M_0,
M_3)$ (equivalently, for tp$_{{\frak s}(*)}$)
and NF$_{{\frak s}(+)}(M_0,M^*_1,M_6,M_7)$ and $M^*_1$ is
$(\lambda^+,*)$-brimmed over $M_0$ for ${\frak s}^+$.  We also can find $M_1$
such that $M_0 \cup \{a^*_1\} \subseteq M_1 \le_{s(*)} M^*_1,\bold J_2
\subseteq \bold I^{{\frak s}(*)}_{M_0,M^*_1}$ is finite and $(M_0,M_1,
\bold J_1) \in K^{3,\text{qr}}_{{\frak s}(+)}$.
\sn
Now
\mr
\item "{$\circledast_1$}"  $\bold J_1 \cup \bold J_2$ is independent in
$(M_0,M_7)$ for ${\frak s}^+$ \nl
[why?  see \scite{705-6p.8}(3)]
\sn
\item "{$\circledast_2$}"  tp$_{{\frak s}(*)}(a^*_1,M_6,M_7)$ does not fork
over $M_0$ for ${\frak s}(*)$
\sn
\item "{$\circledast_3$}"  tp$_{{\frak s}(*)}(a^*_1,M_6,M_7) = 
\text{ tp}_{{\frak s}(*)}(a_1,M_6,M_7)$ \nl
[why?  as both are ${\frak s}(*)$-nonforking extensions of
tp$_{{\frak s}(*)}(a_1,M_0,M_3)$.]
\ermn
We are done by subclaim \sciteu{705-x.X} below.
\bn
\ub{Axiom(E)(g)}:  [extension existence]

So let $M_0 \le_{{\frak s}(*)} M_1$ and $p \in 
{\Cal S}^{\text{bs}}_{{\frak s}(*)}
(M_0)$ so for some $N_0 \in {\frak K}_{\frak s},N_0 
\le_{{\frak K}[{\frak s}]} M_0$ and $N_0$ is a 
witness for $p$.  So $M_0,M_1$ are
saturated models in $\lambda^+$ for ${\frak K}^{\frak s}$ hence there
is an isomorphism $f$ from $M_0$ onto $M_1$ over $N_0$ and $f(p) \in
{\Cal S}^{\text{bs}}_{{\frak s}(*)}(M_1)$ is witnessed by $N_0$ 
and extend $p \restriction N_0$ hence $p$.  [Saharon: say more]
\mn
\ub{Axiom (E)(h)}:  By claim \scite{705-6a.4}(1) below.
\mn
\ub{Axiom (E)(i)}:  By \cite[\S2]{600},\marginbf{!!}{\cprefix{600}.\scite{600-1.16t}} it follows.
\mn
Lastly \nl
\ub{${\frak s}(*)$ is good$^+$}:

So assume $\bar M^\ell = \langle M^\ell_\alpha:\alpha < \lambda^{++} \rangle$
is $\le_{{\frak s}(*)}$-increasing continuous, $M^0_\alpha \le_{\frak s}
M^1_\alpha,a_\alpha \in M^0_{\alpha +1}$, tp$(a_{\alpha +1},M^0_{\alpha +1},
M^0_{\alpha +2}) \in {\Cal S}^{\text{bs}}_{{\frak s}(*)}(M^0_{\alpha +1})$ is an
${\frak s}(*)$-nonforking extension of 
$p^* \in {\Cal S}^{\text{bs}}_{{\frak s}(*)}(M^0
_\delta)$ but tp$_{{\frak s}(*)}(a_{\alpha +1},M^1_0,M^1_{\alpha +2})$ \,
does ${\frak s}(*)$-fork over $M^0_0$ and we shall get a
contradiction.  

As $p^* \in {\Cal S}^{\text{bs}}_{{\frak s}(*)}(M^0_0)$ clearly for some
$N^* \in K_{\frak s}$ we have $N^* \le_{{\frak K}[{\frak s}]} M^0_0$
and $p^*$ does not $\lambda$-split over $N^*$ hence (by
\scite{705-6a.4}(2) below) also tp$_{{\frak s}(*)}(a_{\alpha+2},
M^0_{\alpha +1},M^0_{\alpha +2})$ does not $\lambda$-split over
$N^*$.  let $\langle N_\varepsilon:\varepsilon < \lambda^+ \rangle$ be
a $\le_{\frak s}$-representation of $M^1_0$, and \wilog \, $N^*
\le_{\frak s} N_0$.

Now for each $\alpha < \lambda^{++}$ the type tp$_{{\frak
s}(*)}(a_{\alpha +1},M^1_0,M^1_{\alpha +2})$ does ${\frak s}(*)$-fork
over $M^1_0$ hence it does $\lambda$-split over $N^*$, but clearly for
some $\zeta_\alpha < \lambda^+$ it does not $\lambda$-split over
$N_{\zeta(\alpha)}$.  So for some $\zeta^* < \lambda^+$ the set $W =
\{\alpha < \lambda^{++}:\zeta_\alpha = \zeta^*\}$ is unbounded in
$\lambda^{++}$.  Now choose by induction on $\varepsilon < \lambda$ a
triple $(\alpha_\varepsilon,M_{0,\varepsilon},M_{1,\varepsilon})$ such
that:
\mr
\item "{$(a)$}"  $\alpha_\varepsilon \in S$ is increasing
\sn
\item "{$(b)$}"  $M_{0,\varepsilon} \le_{{\frak K}[{\frak s}]}
M^0_{\alpha_\varepsilon}$ is $\le_{\frak s}$-increasing continuous
\sn
\item "{$(c)$}"  $M_{1,\varepsilon} \le_{{\frak K}[{\frak s}]}
M^1_{\alpha_\varepsilon}$ is $\le_{\frak s}$-increasing continuous
\sn
\item "{$(d)$}"  $M_{0,\varepsilon} \le_{\frak s} M_{1,\varepsilon}$
\sn
\item "{$(e)$}"  $a_{\varepsilon(\alpha)} \in M_{0,\varepsilon +1}$
\sn
\item "{$(f)$}"  $N^* \subseteq M_{0,\varepsilon},N_{\zeta^*}
\subseteq M_{1,\varepsilon}$.
\ermn
There is no problem to carry the definition and $\langle
(M_{0,\varepsilon},M_{1,\varepsilon};a_\varepsilon):\varepsilon <
\lambda^+ \rangle$ provide a counterexample to ``${\frak s}$ is good$^+$".
\enddemo
\hfill$\square_{\scite{705-6a.2}}$\margincite{705-6a.2}
\bigskip

\proclaim{\stag{705-6a.4} Claim}  1) Assume the pre-$\lambda$-frame ${\frak r}$
(see \scite{705-0.2A}) satisfies axiom (E)(c),(d) (of good frames of 
\marginbf{!!}{\cprefix{600}.\scite{600-1.1t}}) and 
${\Cal S}^{\text{bs}}_{\frak r} 
= {\Cal S}^{\text{na}}_{{\frak K}_{\frak r}}$.

\ub{Then} it satisfies (E)(h), too.  [have appeared?] \nl
2) In clause (c) of Definition \scite{705-6a.1}(4), an equivalent
condition is
\mr
\item "{$(*)$}"  if $N_0 \le_{{\frak K}[{\frak s}]} M_0,N_0 \in
K_{\frak s}$ and tp$_{{\frak s}(*)}(a,M_0,M_3)$ does not
$\lambda$-split over $N_0$ then also tp$_{{\frak s}(*)}(a,M_1,M_3)$
does not $\lambda$-split over it.
\endroster
\endproclaim
\bigskip

\demo{Proof}  So assume $\langle M_i:i \le \delta \rangle$ be
$\le_{\frak r}$-increasing continuous, $p \in {\Cal S}^{\text{bs}}_{\frak r}
(M_\delta)$ and 
$p \restriction M_i$ belongs to ${\Cal S}^{\text{bs}}_{\frak r}(M_i)$, does
not ${\frak r}$-fork over $M_0$ for $i < \delta$.  
As $p \restriction M_i$ is not realized in $M_i$ for $i < \delta$ this
holds for $i = \delta$, too.
Let $(M_\delta,M_{\delta +1},a) \in K^{3,\text{bs}}_{\frak r}$ 
be such that $p =
\text{tp}_{\frak r}(a,M_\delta,M_{\delta +1})$; now by an assumption
$p$ belongs to ${\Cal S}^{\text{bs}}_{\frak r}(M_\delta)$ hence
by (E)(c) for some $i < \delta$ the type $p = \text{ tp}_{\frak r}(a,
M_\delta,M_{\delta +1})$, does not ${\frak r}$-fork over $M_i$.  But we
assume that $p \restriction M_i$ does not ${\frak r}$-fork over $M_0$.  So
by Axiom (E)(d) together we get $p$ does not ${\frak r}$-fork over $M_0$.
\hfill$\square_{\scite{705-6a.4}}$\margincite{705-6a.4}
\enddemo
\bn
Of course
\proclaim{\stag{705-6a.5} Claim}  If 
${\frak s}$ is a successful type full $\lambda$-good$^+$-frame, 
\ub{then} ${\frak s}^+$ is a full $\lambda^+$-good frame and
${\frak s}^* = {\frak s}^+$.
\endproclaim
\bigskip

\demo{Proof}  Easy.
\enddemo
\bigskip

\proclaim{\stag{705-6a.6} Claim}  [${\frak s}$ is successful and type-full 
$\lambda$-good frame].

If $M_0 \le_{\frak s} M_1 \le_{\frak s} M_2$ and $(M_\ell,M_{\ell +1},
a_\ell) \in K^{3,\text{pr}}_\lambda$ 
for $\ell = 1,2$, \ub{then} $(M_0,M_1,a_0a_1) 
\in K^{3,\text{pr}}_\lambda$.
\endproclaim
\bigskip

\demo{Proof}  Use \scite{705-5.5}(2). \nl
[But $\oplus$-closed suffices.]
\enddemo
\bn
\margintag{705-gr.8}\ub{\stag{705-gr.8} Question}:  Add on ${\frak s}^{nf}$ for ${\frak s}$
saturative.  Saharon!
\bn
\centerline{$* \qquad * \qquad *$}
\bigskip

\proclaim{\stag{705-10b.6} Claim}  [Here?] Assume ${\frak s}$ is super ${\frak
t}$-local, both full.

If $M_0 \le_{K[{\frak t}]} M_1,M_0 \in K_{\frak t},M_1 \in {\frak
K}_{\frak s},p_1 \in {\Cal S}(M_1)$ does not fork over $M_0,p_0 = p
\restriction M_0 < {\Cal S}_{\frak t}(M_0)$, \ub{then} 
{\rm rk}$_{\frak s}(p_1) = { \text{\rm rk\/}}_{\frak t}(p_0)$.
\endproclaim
\newpage

\head {\S10 Regular types} \endhead  \resetall \sectno=10
 \spuriousreset
\bigskip

\demo{\stag{705-7.0} Hypothesis}  ${\frak s}$ is a $\lambda$-good weakly
successful frame with primes such that ${\frak s}$ is type-full.
\enddemo
\bigskip

\definition{\stag{705-7.1} Definition}  1) We say that
$p \in {\Cal S}^{\text{bs}}_{\frak s}(M)$ is
regular \ub{if} there are $M_0,M_1,a,M_2$ such that:
\mr
\item "{$(a)$}"  $M_\ell$ is $(\lambda,*)$-brimmed over $M_0$ for
$\ell =1,2$
\sn
\item "{$(b)$}"  $M_0,M \le_{\frak s} M_1 \le_{\frak s} M_2,a \in M_2$
\sn
\item "{$(c)$}"  $p' = \text{ tp}(a,M_1,M_2)$ is parallel to $p$
\sn
\item "{$(d)$}"  $p'$ does not fork over $M_0$
\sn
\item "{$(e)$}"  if $c \in M_2 \backslash M_1$ realizes $p' \restriction M_0$
\ub{then} $c$ realizes $p'$.
\ermn
2) We say that $p \in {\Cal S}^{\text{bs}}_{\frak s}
(M)$ is regular$^+$ \ub{if} there
are $M_1,M_2,a$ such that clauses (a)-(d) above holds and (see \S2)
\mr
\item "{$(e)'$}"  if 
$c \in M_2 \backslash M_1$, \ub{then} rk(tp$(c,M,M_2)) \ge
\text{ rk}(p)$.
\ermn
3) We add ``directly", if
\mr
\item "{$(g)$}"  $M_1 = M$.
\endroster
\enddefinition
\bigskip

\remark{Remark}  Note that regular $\ne$ regular$^+$, e.g. $T$ is the
first order theory of $M = (\omega \times \omega \cup
\omega,P^M,Q^M,F^M)$ when $P^M = \omega \times \omega,Q^M =
\omega,F^M((n,m)) = n,F^M(n) = n$.

By and for our purposes every regular type is ``equivalent to a
regular$^+$ type so those suffice.
\endremark
\bigskip

\proclaim{\stag{705-7.2} Claim}  1)
\mr
\item "{$(a)$}"   If $p_1\|p_2$ \ub{then}
$p_1$ is regular \ub{iff} $p_2$ is regular
\sn
\item "{$(b)$}"  if $M$ is $(\lambda,*)$-brimmed (trivially holds if ${\frak
K}_{\frak s}$ is categorical) and $p \in {\Cal S}^{\text{bs}}_{\frak s}(M)$,
\ub{then} $p$ is regular iff it is directly regular.
\ermn
2) If $p \in {\Cal S}^{\text{bs}}(M)$ is regular, $M$ is
$(\lambda,*)$-brimmed over $M_0,p$ does not fork over $M_0$ and
$(M,M_2,a) \in K^{3,\text{pr}}_{\frak s}$, 
{\rm tp}$(a,M,M_2) = p$ and we let $M_1 = M$. \ub{Then} $p$ is regular
\ub{iff} clause (e) of \scite{705-7.1} holds.  So trivially 
(a), (b), (c), (d)) of \scite{705-7.1} hold, i.e., holds for $M_0,M_1,M_2,a$.\nl
3) The parallel of parts (1),(2) holds for regular$^+$. \nl
4) If $(M,N,a) \in K^{3,\text{pr}}_{\frak s}$ 
and $p = { \text{\rm tp\/}}_{\frak s}(a,M,N)$ is 
regular$^+,M$  is $(\lambda,*)$-brimmed, \ub{then}
$c \in N \backslash M \Rightarrow { \text{\rm rk\/}}
({\text{\rm tp\/}}(c,M,N)) \ge { \text{\rm rk\/}}(p)$. \nl
5) If $p$ is regular$^+$ \ub{then} $p$ is regular. \nl
\endproclaim
\bigskip

\demo{Proof}  1a)  So assume that $M' \le_{\frak s} M$
and $M'' \le_{\frak s} M$ and $p' \in {\Cal S}^{\text{bs}}(M'),p'' \in
{\Cal S}^{\text{bs}}(M'')$ are parallel, that is some $p \in {\Cal
S}^{\text{bs}}(M)$ 
does not fork over $M'$ and over $M''$ and $p \restriction
M' = p',p \restriction M'' = p''$ and we should prove that $p'$ is
regular iff $p''$ is regular.  By the symmetry it suffices to show
that $p'$ is regular iff $p$ is regular.  Now the ``if" direction is
trivial (the same witnesses $M_0,M_1,M_2,a$ works).  For the ``only
if" direction, let $(M'_0,M'_1,M'_2,a)$ witness $p'$ is regular and \wilog
\, $M'_2$ is $(\lambda,*)$-brimmed over $M'_1$.

As ${\frak K}_{\frak s}$ has amalgamation and $M' \le_{\frak s} M,M'
\le_{\frak s} M'_2$ \wilog \, for some $M_1$ we have $M'_1 \le_{\frak
s} M_1$ and $M \le_{\frak s} M_1$ and \wilog \, $M_1$ is
$(\lambda,*)$-brimmed over $M'_1 \cup M$.  There is an
isomorphism $f$ from $M'_1$ onto $M_1$ over $M'_0$ as both are
$(\lambda,*)$-brimmed over it, and we can find $f^+,M_2,a^*$ such
that $M_1 \le_{\frak s} M_2,f^* \supseteq f,f^*$ an isomorphism from $M'_2$
onto $M_2$ and $f^*(a) = a^*$.

Now $f^*(M'_0),M_1 = f^*(M'_1),M_2 = f^*(M'_2)$ and $a^*$ witnesses
the regularity of $p$. 
\nl
1b)  The if direction is obvious (same witnesses).

For the other direction assume that $M_0,M_1,M_2,a$ witness that $p
\in {\Cal S}^{\text{bs}}(M)$ is regular.  There is $M'_0$ such that $M_0
\le_{\frak s} M'_0 \le_{\frak s} M_1$ such that $M'_0$ is
$(\lambda,*)$-brimmed over $M_0$ and $M_1$ is $(\lambda,*)$-brimmed
over $M'_0$.  Clearly there is an isomorphism $f$ from $M'_0$ onto $M$
and it can be extended to an isomorphism $f^+$ from $M_1$ onto $M_1$.
Without loss of generality $f^+(\text{tp}(a,M_1,M_2)) = \text{
tp}(a,M_1,M_2)$ (as both types does not fork over $M_0$ and
$M'_0,M_1$ are isomorphic over $M_0$ hence
there is $g \in \text{ Aut}(M'_0)$ over $M_0$ such that
$g(\text{tp}(a,M'_0,M_2)) = f^{-1}(p)$.  Hence replacing $f$ by $f
\circ g$ we are ``done".  Using $f^+$ we can find $f^*,M'_2,a'$ such that
$f^* \supseteq f^+$ is an isomorphism from $M_2$ onto $M'_2,f^*(a) =
a'$ and $(f^*(M_0),M_1,M'_2,a')$ is a witness to $p$ directly regular. \nl
2), 3), 4)  Similar. \nl
5) Because if $M_1$ is $(\lambda,*)$-brimmed over $M_0$ and $p \in
{\Cal S}^{\text{bs}}(M_1)$ 
and $q \in {\Cal S}(M_1),q \ne p,q \restriction M_0
= p \restriction M_0$ then rk$(p) = \text{ rk}(p \restriction M_0) >
\text{ rk}(q)$, see \scite{705-gr.12}(4) + \scite{705-10b.3B}.
\hfill$\square_{\scite{705-7.2}}$\margincite{705-7.2}
\enddemo
\bigskip

\proclaim{\stag{705-7.3} Claim}  1) If $M <_{\frak s} N$ and $M$ is
$(\lambda,*)$-brimmed, \ub{then} for some $c \in
N \backslash M$ the type {\rm tp}$(c,M,N)$ 
is regular$^+$ (hence regular). \nl
2) If $M_0 <_{\frak s} 
M <_{\frak s} N,M$ is $(\lambda,*)$-brimmed and $p \in
{\Cal S}^{\text{bs}}_{\frak s}(M_0)$ is realized by some member of $N
\backslash M$ and $M$ is $(\lambda,*)$-brimmed over $M_0$, 
\ub{then} for some $c_1 \in N \backslash M$ realizing $p$ we have 
{\rm tp}$_{\frak s}(c_1,M,N)$ is regular.
\endproclaim
\bigskip

\demo{Proof}  1) Choose $c \in N \backslash M$ such that rk$_{\frak
s}(\text{tp}_{\frak s}(c,M,N))$ is minimal.  Then choose 
$(\lambda,*)$-brimmed
$M_0 <_{\frak s} M$ such that $M$ is $(\lambda,*)$-brimmed over
$M_0$ and tp$(c,M,N)$ does not fork over $M_0$. \nl
2) Choose $c \in N \backslash M$ realizing $p$ with rk$_{\frak
s}(\text{tp}_{\frak s}(c,M,N))$ minimal. 
Let $a_1 \in M^1_0$ realize $p \restriction M_1$, let $f$ be an
isomorphism from $M_1$ onto $M$ (exists as there is $M^- <_{\frak s}
M$ such that $M$ is $(\lambda,*)$-brimmed over $M^-$ and $p$ does not
fork over $M^-$.  \hfill$\square_{\scite{705-7.3}}$\margincite{705-7.3}
\enddemo
\bigskip

\proclaim{\stag{705-7.4} Claim}  [${\frak s} = {\frak t}^+,{\frak t}$ is
$\lambda$-good$^+$ successful with primes \ub{or} just ${\frak s}$ is super
${\frak t}$-local ${\frak t}$ with primes [Saharon]].

Assume $(M,N,a) \in K^{3,\text{bs}}_{\frak s},p 
= { \text{\rm tp\/}}_{\frak s}(a,M,N)$ and $q \in
{\Cal S}^{\text{bs}}_{\frak s}(M)$ and $M_0 \in K_{\frak t},M_0 
\le_{{\frak K}[{\frak t}]} M$ (so $M$ is $(\lambda,*)$-brimmed
[if ${\frak s} = {\frak t}^+$]). \nl
1) If $p$ does not fork over $M_0$ \ub{then} 
\mr
\item "{$(a)$}"  $p$ is regular (for ${\frak s}$)
iff $p \restriction M_0$ is regular (for ${\frak t}$)
\sn
\item "{$(b)$}"  Similarly for regular$^+$
\sn 
\item "{$(c)$}"  if $(M,N,a) \in K^{3,\text{pr}}_{\frak s}$ and $c \in N
\backslash M$ realizes $p \restriction M_0$ \ub{then} $c$ realizes $p$.
\ermn
2) There is a regular$^+$ type $p_1 \in {\Cal S}^{\text{bs}}_{\frak s}
(M)$ not orthogonal to
$p$, and realized in $N$ such that {\rm rk}$_{\frak s}(p_1) \le 
{ \text{\rm rk\/}}_{\frak s}(p)$ and {\rm rk}$_{\frak s}(r) 
< { \text{\rm rk\/}}_{\frak s}(p_1)
\Rightarrow r \bot p$ for $r \in {\Cal S}_{\frak s}(M)$ 
or just $r \in {\Cal S}_{\frak s}(M'),M \le_{\frak s} M'$.
\ub{Question}:  Is ${\frak s}$ not ${\frak t}$? \nl
3) If $p,q$ are regular$^+$ not orthogonal, \ub{then} $q$ is realized in
$N$ and rk$_{\frak s}(q) = { \text{\rm rk\/}}_{\frak s}(p)$. \nl
4) If $M^* \le_{\frak s} M,p \restriction M^* = 
q \restriction M^*,p \ne q,
p$ does not fork over $M^*$ and $p$ is regular \ub{then} $p \bot q$. \nl
5) If $p$ is regular$^+$ and {\rm rk}$_{\frak s}(q) < 
{ \text{\rm rk\/}}_{\frak s}(p)$ \ub{then} $p \bot q$. \nl
6) Let $p_1 \in {\Cal S}^{\text{bs}}_{\frak s}(M)$ be not 
orthogonal to $p$ with minimal rank.
\ub{Then}
\mr
\item "{$(\alpha)$}"  $p_1$ is realized in $N$ and is regular$^+$
\sn
\item "{$(\beta)$}"  if $p$ is regular and
$(M,N^1,a^1) \in K^{3,\text{pr}}_{{\frak s},p_1}$ \ub{then} $p$
is realized in $N^1$.
\ermn
7) If $a_1 \in N \backslash M$ and $p \bot q$ and $(M,N,a) \in 
K^{3,\text{pr}}_{\frak s}$ \ub{then} {\rm tp}$_{\frak s}(a_1,M,N) 
\bot q$ and $M' <_{\frak s} M \and p \bot M' 
\Rightarrow { \text{\rm tp\/}}_{\frak s}(a_1,M,N) \bot M'$.
\nl
8) If $p,q$ are regular not orthogonal \ub{then} $q$ is realized in $N$.
\endproclaim
\bigskip

\demo{Proof}  For the case ${\frak s} = {\frak t}^+$ let 
$\langle M_\alpha:\alpha < \lambda^+ \rangle,\langle
N_\alpha:\alpha < \lambda^+ \rangle$ 
be $\le_{{\frak K}[{\frak t}]}$-representation of
$M,N$ respectively.  Without loss of generality $(M,N,a) 
\in K^{3,\text{pr}}_{\frak s}$ and $\alpha < \lambda^+ 
\Rightarrow (M_\alpha,N_\alpha,a) \in
K^{3,\text{uq}}_{\frak t}$ and $\alpha < \beta \Rightarrow \text{ NF}
_{\frak t}(M_\alpha,N_\alpha,M_\beta,N_\beta)$ hence $(M_\alpha,N_\alpha,
a) \in K^{3,\text{pr}}_{\frak t}$ [used?] and 
$\alpha < \beta \Rightarrow N_\beta,M_\beta$
is $(\lambda,*)-{\frak t}$-brimmed over $N_\alpha,M_\alpha$
respectively. \nl
1) Easy.  (see \scite{705-10b.3} which 
deals with rk, but the proof works). \nl
2) Choose $p_1 \in {\Cal S}^{\text{bs}}_{\frak s}(M)$ 
realized by some $c_1 \in N
\backslash M$ with rk$_{\frak s}
(p_1)$ minimal and let $N_1 \le_{\frak s} N$ be
such that $(M,N_1,c_1) \in K^{3,\text{pr}}_{\frak s}$.   
Now $p_1$ is regular$^+$
by \scite{705-7.2}(2),(3) it is realized in $N$ as exemplified by $c_1$.  Also
{\rm rk}$(p_1) \le { \text{\rm rk\/}}(p)$ by the minimality of 
{\rm rk}$(p_1)$.

Lastly, assume $r \in {\Cal S}^{\text{bs}}(M'),r \pm p,M \le_{\frak s} M'$;
\wilog \, $M = M'$ (see \scite{705-6p.10} and its proof).  
Without loss of generality
$p,p_1,r$ do not fork over $M_0$, 
so as $N$ is $\lambda^+_{\frak t}$-saturated (for
$K^{\frak t}$) there is $c_2 \in N$ realizing $r \restriction M_0$ such
that $\{a,c_2\}$ is not independent over $M_0$ inside $N$.
Now this implies $c_2 \notin M$ hence
rk$_{\frak s}(\text{tp}_{\frak s}(c_2,M,N)) \ge \text{ rk}_{\frak
s}(p_1)$ so necessarily using \scite{705-10b.3}? rk$_{\frak s}(r)= 
\text{ rk}_{\frak t}(r \restriction M_0)= 
\text{ rk}_{\frak t}(\text{tp}_{\frak t}(c_2,M_0,N)) 
\ge \text{ rk}_{\frak s}(\text{tp}_{\frak s}(c_2,M,N) \ge \text{ rk}_{\frak s}(p_1)$ as 
required. \nl
3)  First assume rk$(q) \le \text{ rk}(p)$.

Without loss of generality 
$p,q$ does not fork over $M_0$, so $p \restriction M_0,q
\restriction M_0$ are regular$^+$ not orthogonal (by part (1) and by
\scite{705-6.xobO} respectively).  As $a \in N$
realizes $p \restriction M_0$ and $N$ is $\lambda^+_{\frak t}$-saturated
there is $c \in N$ realizing $q \restriction M_0$ such that 
$\{a,c\}$ is not independent over $M_0$ inside $N$ for ${\frak t}$.  Hence 
$c \notin M$, hence rk$_{\frak s}(\text{tp}_{\frak s}(c,M,N)) \le 
\text{ rk}_{\frak t}(\text{tp}_{\frak t}(c,M_0,N)) = \text{ rk}_{\frak
s}(q) \le \text{ rk}_{\frak s}(p)$ so rk$_{\frak s}(\text{tp}_{\frak
s}(c,M,N)) = \text{rk}(p)$.

As $p$ is regular$^+$ and $(M,N,a) \in K^{3,\text{pr}}_{{\frak s},p}$ 
this implies 
that rk$_{\frak s}(\text{tp}_{\frak s}(c,M,N))$ is equal to rk$_{\frak
s}(p)$, hence to rk$_{\frak s}(q)$, so necessarily
tp$_{\frak s}(c,M,N) = q$ and rk$_{\frak s}(p) = \text{ rk}_{\frak
s}(q)$. 
We are left with the case rk$_{\frak s}(p) < \text{ rk}_{\frak s}(q)$.
But then interchanging $p$ and $q$ (and replacing $N,a$ by
others) we get a contradiction. \nl
4) Without loss of generality $M^*_0 =: M^* \cap M_0 <_{\frak t} M_1$
and $p$ does not fork over it.

Without loss of generality and $p$ and $q$ does not fork over (are
witnessed by) $M_0$ and $p \restriction M_1$ and $p \bot q
\Leftrightarrow (p \restriction M_0 \bot q \restriction M_0)$.  Assume toward
contradition that $p \pm q$ hence for some $c \in N$ realizing
$q \restriction M_0,\{a,c\}$ is not independent over $M_0$ inside $N$, 
hence $c \in N \backslash M$.  So
choose $c' \in N \backslash M$ realizing $q \restriction M_0$ with 
rk$_{\frak s}(\text{tp}_{\frak s}(c',M,N))$ minimal, and choose
$\alpha < \lambda^+_{\frak t}$ such that $M_\alpha \le_{\frak t} N'
\le_{\frak t} N_\alpha$ and tp$(c',M,N)$ does not fork over $M_\alpha$ 
and $c' \in N'_\alpha$ and $(M_\alpha,N',a)$ belongs 
to $K^{3,\text{uq}}_{\frak t}$
hence to $K^{3,\text{pr}}_{\frak t}$, 
and $N'$ is $(\lambda,*)-{\frak t}$-saturated over $M_0$. \nl
So $M_0,M_\alpha,N',a,p \restriction M_\alpha,c'$ 
contradict ``$p \restriction M_1$
is ${\frak t}$-regular (check Definition \scite{705-7.1}(1)). \nl
5) Proof similar to (3). \nl 
6) Clause $(\alpha)$: as in the proof of part (2).
\nl
\ub{Proof of Clause $(\beta)$}:

By clause $(\alpha)$ we know that $p_1$ is realized in $N$, so \wilog
\, $N^1 \le_{\frak s} N$ hence $a^1 \in N$.  Without loss of
generality both $p$ and $p_1$ does not fork over $M_0$ and so as in
earlier cases there is $c \in N^1$ realizing $p \restriction M_0$ such
that $\{c,a^1\}$ is not independent in $N$ over $M_0$.  This implies
$c \in N^1 \backslash M$ hence by clause (c) of part (1) we know that
$c$ realizes $p$, as required. \nl
7)  Easy. \nl
8)  Let $p_1 \in {\Cal S}^{\text{bs}}(M)$ 
be not orthogonal to $p$ of minimal
rank, and let $q_1 \in {\Cal S}^{\text{bs}}_{\frak s}(M)$ 
be not orthogonal to $q$ of minimal rank.  
By part (6), clause $(\alpha)$ $p_1,q_1$ are regular$^+$.  By part (6)
clause $(\alpha),p_1$ is realized in $N$ say by $a_1$.  Now $p_1 \pm q_1$
(apply twice part (6)$(\beta)$+(7) to get 
$(p_1 \bot q_1 \Rightarrow p \bot q))$, hence by
part (3) some $b_1 \in N$ realizes $q_1$ hence by 
part (6) clause $(\beta)$
(applied to $q,q_1$) some $b \in N$ realizes $q$, as required.  \nl
${{}}$  \hfill$\square_{\scite{705-7.4}}$\margincite{705-7.4} 
\enddemo
\bigskip

\demo{\stag{705-7.5} Conclusion}  1) Non-orthogonality among 
regular types is an equivalence relation.   \nl
2) For regular $p,q \in {\Cal S}^{\text{bs}}(M)$ and $r \in {\Cal S}(M)$ we
have $p \pm q,q \pm r \Rightarrow p \pm r$. \nl
3) For $p,q,r \in {\Cal S}^{\text{bs}}(M),q$ regular, $p \pm q,q \pm r$ we
have $p \pm r$. \nl
4) For nonorthogonal $p,q \in {\Cal S}^{\text{bs}}(N)$
and $M \le_{\frak s} N$, we have $p \bot M \Leftrightarrow q \bot
M$). \nl
[Saharon: with categoricity use $({\frak s},{\frak t}),{\frak s}$
saturated over ${\frak t}$.
\enddemo
\bigskip

\proclaim{\stag{705-7.6} Claim}  [${\frak s}$ categorical 
in $\lambda_{\frak s}]$. \nl
1) If $M \le_{\frak s} N$ and $p \in {\Cal S}^{\text{bs}}_{\frak s}(N)$ is 
not orthogonal to $M$ \ub{then} there is $q \in 
{\Cal S}^{\text{bs}}_{\frak s}(N)$ 
not orthogonal to $p$,
conjugate to $p$ (i.e., $f(p) = q$ for 
some $f \in \text{ Aut}(M))$ and $q$ does not fork over $M$. \nl
2) If $\langle M_i:i \le \delta +1 \rangle$ is $\le_{\frak s}$-increasing
continuous and $M_\delta \ne M_{\delta +1}$, \ub{then} for some $c \in
M_{\delta +1} \backslash M_\delta$ and nonlimit $i < \delta$, we have
{\rm tp}$(c,M_\delta,M_{\delta +1})$ 
does not fork over $M_i$, and is orthogonal
to $M_{i-1}$ if $i > 0$ so ${\frak s}$ has enough regulars (see
Definition x.x). \nl
3) If in part (2), $q \in {\Cal S}^{\text{bs}}_{\frak s}(M_\delta)$ 
is regular realized by
some member of $M_{\delta +1}$, \ub{then} 
we can demand {\rm tp}$(c,M_i,M_{\delta +1})$ is conjugate to $q$. \nl
4) In part (1), if for some $M_0$ is $(\lambda,*)$-brimmed over $M_0$,
\ub{then} we can get $q$ conjugate to $p$ over $M_0$.
\endproclaim
\bigskip

\demo{Proof}  1) Let $r \in {\Cal S}^{\text{bs}}(M)$ 
be not orthogonal to $p$.
Let $\langle M_\alpha:\alpha \le \omega \rangle,\langle
N_\alpha:\alpha \le \omega \rangle$ be as in the proof of
\scite{705-6p.10}, i.e., $M = M_\omega = \cup \{M_n:n < \omega\},N =
N_\omega = \cup\{N_n:n < \omega\}$, NF$_{\frak s}
(M_n,N_n,M_{n+1},M_{n+1},N_{n+1})$ and $M_{n+1},N_{n+1}$ is
$(\lambda,*)$-brimmed over $M_n,N_n$ respectively;
\wilog \, $p$ does not fork over $N_0$ and $r$ does not fork
over $M_0$.  We can find $\langle f_i:i < \lambda^+ \rangle$ such that
$f_{1+i}$ is a $\le_{\frak s}$-embedding of $N_0$ into $M$ over $M_0,
f_0 = \text{ id}_{N_0}$, such that $\langle f_i(N_0):i < \omega
\rangle$ is independent over $M_0$ (see \scite{705-6p.10}) and
clearly $f_i(p \restriction N_0) \pm q \restriction M_0$, hence $f_i(p
\restriction N_0) \pm M_0$.  By \scite{705-6p.10} clearly $p \restriction N_0
\pm f_1(p \restriction N_0)$ and let $q \in 
{\Cal S}^{\text{bs}}(N)$ be a nonforking
extension of $f_1(p \restriction N_0)$. \nl
2) By \scite{705-7.3}(1) for some $d \in M_{\delta +1} \backslash M_\delta$, the
type tp$(d,M_\delta,M_{\delta +1})$ is regular, and apply part (3). \nl
3) Let $j = \text{ Min}\{i \le \delta:q \pm M_i\}$, as $q \pm M_\delta$
clearly $j$ is well 
defined.  By \scite{705-5.8}(2), $j$ is a nonlimit ordinal and
by part (1) there is $r \in {\Cal S}^{\text{bs}}(M_\delta)$ not forking over
$M_j$ not orthogonal to $q$ and
conjugate to $q$ hence $r$ is regular and by \scite{705-7.5}(4)a is orthogonal to
$M_{j_1}$ for $j_1 < j$ but not orthogonal to $p$.

By \scite{705-7.4}(8) some $c \in M_{\delta +1} \backslash M_\delta$ realizes
$r$. \nl
4) Easy.    \hfill$\square_{\scite{705-7.6}}$\margincite{705-7.6}
\enddemo
\bigskip

\proclaim{\stag{705-7.6R} Claim}  If {\rm NF}$_{\frak s}(M_0,M_1,M_2,M_3)$
and $p \in {\Cal S}^{\text{bs}}_{\frak s}(M_3)$ is regular and $p \pm M_1,p \pm
M_2$ then $p \pm M_0$.
\endproclaim
\bigskip

\demo{Proof}  FILL (used in \scite{705-12b.3}).
\enddemo
\bigskip

\proclaim{\stag{705-7.6A} Claim}  1) Assume
\mr
\item "{$(a)$}"  ${\frak s}$ is super ${\frak t}$-local
\sn
\item "{$(b)$}"  ${\frak t}$ is weakly successful with primes,
categorical in $\lambda_{\frak t}$.
\ermn
\ub{Then} the conclusion of \scite{705-7.6} holds. \nl
2) If ${\frak s}$ satisfies \scite{705-7.6}(1) then it satisfies
\scite{705-7.6}(2), (3), (4).
\endproclaim
\bigskip

\demo{Proof}  1) By part (2) it suffices to prove \scite{705-7.6}(1), it
holds by the same proof using \scite{705-gr.4A}(8) instead of
\scite{705-6p.10}. \nl
2) Same proof.  \hfill$\square_{\scite{705-7.6A}}$\margincite{705-7.6A}
\enddemo
\bigskip

\definition{\stag{705-7.7} Definition}  We call $(\bar M,\bar{\bold J}) \in
{\Cal W}$ regular if $c \in \bold J_i \Rightarrow \text{ tp}(c,M_i,M_{i+1})$
is regular; we say ``regular except $\bold J$" if $c \in \bold J$ are
excluded.
\enddefinition
\bigskip

\proclaim{\stag{705-7.8} Claim}  [${\frak s}$ categorical in
$\lambda_{\frak s}$ or just the conclusions of \scite{705-7.6}.]
\nl
1) Assume $M \le_{\frak s} N$ and $\bold J
\subseteq \bold I_{M,N}$ is independent in $(M,N)$.  \ub{Then} we can find
a prime $(\bar M,\bar{\bold J}) \in K^{3,\text{ar}}_{\frak s}$ 
with $\bold J
\subseteq \bold J_0,M_0 = M,N = \cup\{M_n:n < \omega\}$ and $(\bar M,
\bar{\bold J})$ is regular except (possibly) $\bold J$. \nl
2) If $(M,N,a) \in K^{3,\text{uq}}_{\frak s}$ 
\ub{then} we can find prime regular
$(\bar M,\bar{\bold J}) \in K^{3,\text{ar}}_{\frak s}$ 
with $\bold J_0 = \{a\},
M_0 = M,N = \cup\{M_n:n < \omega\}$. \nl
3) If $N_0 \le_{\frak s} N_1 \le_{\frak s} N_2,c \in
N_2 \backslash N_1$,{\rm tp}$(c,N_1,N_2) \pm N_0$ \ub{then} for some 
$b \in N_2 \backslash N_1$ the type {\rm tp}$(b,N_1,N_2)$ 
does not fork over $N_0$ and is regular.
\endproclaim
\bigskip

\demo{Proof}  1) Like the proof of \scite{705-c.6} using \scite{705-7.6}(2).
\nl
2) Follows. \nl
3) Apply part (1) with $(N_1,N_2,\emptyset)$ here standing for
$(M,N,\bold J)$ there and get
$(\bar M,\bar{\bold J})$ as there.  If for some $c \in \bold J_0$,
tp$(c,M_0,M_1) \pm N_0$, 
then by \scite{705-7.6}(3) we get the desired conclusion.
Otherwise, we get contradiction by claim \scite{705-7.9} below.
\hfill$\square_{\scite{705-7.8}}$\margincite{705-7.8}
\enddemo
\bigskip

\proclaim{\stag{705-7.9} Claim}  If 
$(\bar M,\bar{\bold J}) \in K^{3,\text{ar}}_{\frak s}$
is prime, $N <_{\frak s} M_0$ and $c \in \bold J_0 \Rightarrow 
{ \text{\rm tp\/}}(c,M_0,M_1) \bot N$ \ub{then} $M_0 <^{\frak s}_N 
\dbcu_n M_n$ (see Definition \scite{705-5.4}(2)), i.e., every
$q \in {\Cal S}^{\text{bs}}(M_0)$ 
not forking over $N$ has a unique extension in 
${\Cal S}^{\text{bs}}(\dbcu_{n < \omega} M_n)$.
\endproclaim
\bigskip

\demo{Proof}  Easy (put in \S5?).
\enddemo
\bigskip

\definition{\stag{705-10b.10} Definition}  1) For $M \le_{\frak s} N$ let
$\bold I^{\text{reg}}_{M,N} = \{c \in N:\text{ tp}(c,M,N)$ is
regular$\}$. \nl
2) For $M <_{\frak s} N$, we define on $\bold I^{\text{reg}}_{M,N}$ a
dependence relation called the $(M,N)$-dependence relation by:
\mr
\item "{$(a)$}"  $\bold J \subseteq \bold I^{\text{reg}}_{M,N}$ is
$(M,N)$-independent if it is independent
\sn
\item "{$(b)$}"  $c \in \bold J^{\text{reg}}_{M,N}$ is
$(M,N)$-dependent on $\bold J \subseteq \bold J^{\text{reg}}_{M,N}$
\ub{if} there is an independent $\bold J' \subseteq \bold J$ such that
$c \in \bold J'$ or $\bold J \cup \{c\}$ is not independent.
\ermn We omit $(M,N)$ if clear.
\enddefinition
\bigskip

\remark{Remark}  We can  use only regular$^+$ types; somewhat simplify.
\endremark
\bigskip

\proclaim{\stag{705-10b.11} Claim}  Assume $M \le_{\frak s} N$. \nl
1) The relations in \scite{705-10b.10} and their negations are preserved
\ub{if} we replace $N$ by a $\le_{\frak s}$-extension. \nl
2) If $\bold J_1,\bold J_2 \subseteq \bold I^{\text{reg}}_{M,N}$ are
$(M,N)$-independent, every $b \in \bold J_2$ does $(M,N)$-depend on $\bold
J_1$ and $c \in \bold I^{\text{reg}}_{M,N}$ depend on $\bold J_2$,
\ub{then} $c \in (M,N)$-depend on $\bold J_1$. \nl
3) The $(M,N)$-dependence relation satisfies the axioms of dependence
relation, such that dimension is well defined. \nl
4) $\bold J \subseteq \bold I^{\text{reg}}_{M,N}$ is a maximal
$(M,N)$-independent subset of $\bold I^{\text{reg}}_{M,N}$ \ub{iff}
$(M,N,\bold J) \in K^{3,\text{uq}}_{\frak s}$ \ub{iff} $(M,N,\bold J) \in
K^{3,\text{qr}}_{\frak s}$. \nl
5) If ${\Cal P} \subseteq \{p \in {\Cal S}^{\text{bs}}(M):p$ 
regular$\}$ is a
maximal set of pairwise orthogonal types and $\bold J \subseteq \bold
I^{\text{reg}}_{M,N}$, \ub{Then} we can find $\bold J',f$ such that:
\mr
\item "{$(a)$}"  $\bold J' 
\subseteq \bold I^{\text{reg}}_{M,N}$ is $(M,N)$-independent
\sn
\item "{$(b)$}"  $h$ is a function from $\bold J$ onto $\bold J'$ such
that $h(c),(M,N)$-depend on $\{c\}$
\sn
\item "{$(c)$}"  $c \in \bold J' \Rightarrow { \text{\rm tp\/}}_{\frak
s}(c,M,N) \in {\Cal P}$.
\endroster
\endproclaim
\bigskip

\demo{Proof}  Straight. [Details???]
\enddemo
\bigskip

\proclaim{\stag{705-10b.12} Claim}  1) Assume $\bold J_i \subseteq \bold
I_{M,N}$ for $i < i^*$ and $i \ne j \and a \in \bold J_i \and b \in
\bold J_j \Rightarrow { \text{\rm tp\/}}_{\frak s}(a,M,N) 
\perp { \text{\rm tp\/}}(b,M,N)$. \nl
\ub{Then}
\mr
\item "{$(\alpha)$}"  $i \ne j \Rightarrow \bold J_i \cap \bold J_j =
\emptyset$
\sn
\item "{$(\beta)$}"  $\cup \{\bold J_i:i < i^*\}$ is independent in
$(M,N)$ iff for each $i,\bold J_i$ is independent in $(M,N)$.
\ermn
2) Assume $\bold J \subseteq \bold I^{\text{reg}}_{M,N}$ and ${\Cal
E}$ is the following equivalence relation on $\bold I_{M,N}:a Eb
\Leftrightarrow { \text{\rm tp\/}}(a,M,N) 
\pm { \text{\rm tp\/}}(b,M,N)$.  Then $\bold J$ is independent 
in $(M,N)$ \ub{iff} $\bold J/(a/{\Cal E})$ is independent
in $(M,N)$.
\endproclaim
\bigskip

\demo{Proof}  Easy.
\enddemo
\bn
\remark{\stag{705-10b.13} Remark}:  1) Say on weight and simple; so \scite{705-6.1} and
$(M,N,c) \in K^{3,\text{pr}}$; here? \nl
2) Where [${\frak s}$ saturative] if $M_0 \le_{\frak s} M_1 \le_{\frak
s} M_2$, tp$(c,M_1,M_2) \perp M_0$ and $(M_0,M_1,\bold J) \in
K^{3,\text{vq}}_{\frak s}$ 
then $(M_0,M_2,\bold J) \in K^{3,\text{vq}}_{\frak s}$?
\sn
For more on groups see \cite{Sh:F569}.
\endremark
\bigskip

\proclaim{\stag{705-10p.16} Claim}  Assume
\mr
\item "{$(a)$}"  $M \le_{\frak s} N$ are superlimit
\sn
\item "{$(b)$}"  if $p \in {\Cal S}^{\text{bs}}(M)$ is regular then
for some regular $q \in {\Cal S}^{\text{br}}(M)$ we have {\rm
dim}$(q,N) = \lambda_{\frak s}$ and $p \pm q$, (see Definition
\scite{705-4.25Y}, \sciteu{705-5.1A} .
\ermn
\ub{Then} $N$ is $(\lambda,*$)-brimmed over $M$. \nl
[Used in \scite{705-12f.4H}.  (is 6.orth.x)]
\endproclaim
\bigskip

\demo{Proof}  FILL!
\enddemo
\bigskip

\proclaim{\stag{705-10p.17} Claim}  If NF$_{\frak s}(M_0,M_1,M_2,M_3),p \in
{\Cal S}^{\text{bs}}(M_3)$ is regular and $p \pm M_1,p \pm M_2$ then $p
\pm M_0$.
\endproclaim
\newpage

\head {\S11 DOP} \endhead  \resetall \sectno=11
 \spuriousreset
\bn
Note that this is meaningful for non-excellent frame. \nl
\ub{Question}:  Change the framework to super local $\bar{\frak s}$?
or ${\frak s}$ super ${\frak t}$-local? \nl
On weight: see \scite{705-8.2A}, \scite{705-8.2B}, \scite{705-8.2C}.  Main gap for
describing $M_a \in K_{\frak s}$ where $<_{\frak s}$-extend a fixed $N$.
\bigskip

\demo{\stag{705-8.0} Hypothesis}  ${\frak s}$ is $\lambda$-good$^+$ frame which
is successful, with prime models (see below), 
such that $K^{3,\text{uq}}_\lambda = K^{3,\text{pr}}_\lambda$.  [check] \nl
Let ${\frak C} \in K^{\frak s}_{\lambda^+}$ be saturated over $\lambda$.
\enddemo
\bigskip

\definition{\stag{705-8.1} Definition}  1) We say ${\frak s}$ has DOP \ub{if}: 
we can find $M_\ell$ (for $\ell < 4$) and $a_\ell$ (for $\ell = 1,2$) 
and $q$ which exemplifies it, which means
\mr
\item "{$(a)$}"  NF$_{\frak s}(M_0,M_1,M_2,M_3)$
\sn
\item "{$(b)$}"  $(M_0,M_\ell,a_\ell) 
\in K^{3,\text{uq}}_\lambda$ for $\ell =1,2$
\sn
\item "{$(c)$}"  $(M_1,M_3,a_2) \in K^{3,\text{uq}}_\lambda$
\sn
\item "{$(d)$}"  $q \in {\Cal S}^{\text{bs}}(M_3)$ is orthogonal to $M_1$
and to $M_2$.
\ermn
2) Above we also say that $(p_1,p_2)$ has the DOP if there are $M_\ell \,
(\ell < 4),a_\ell \, (\ell = 1,2)$ exemplifying it which means
exemplifying DOP and {\rm tp}$_{\frak s}(a_\ell,M_0,M_\ell)\|p_\ell$;
we say $(p_1,p_2)$ has the explicit DOP if tp$_{\frak
s}(a_\ell,M_0,M_i) = p_\ell$ for $\ell=1,2$. \nl
3) We say ${\frak s}$ has NDOP \ub{if} it fails to have DOP.
\enddefinition
\bigskip

\proclaim{\stag{705-8.2} Claim}  1) If 
$M_\ell \, (\ell < 3),a_\ell(\ell = 1,2)$
satisfies clauses (a), (b), (c) of Definition \scite{705-8.1} 
and $(p_1,p_2) = ({\text{\rm tp\/}}(a_1,M_0,M_1)$, 
{\rm tp}$(a_2,M_0,M_2))$ has the explicit DOP,
\ub{then} replacing $M_3$ by some $M'_3,M_1 \cup M_2 \subseteq M'_3
\le_{\frak s} M_3$ clauses (c), (d) hold.
\endproclaim
\bigskip

\demo{Proof}  1) Let $M'_\ell \, (\ell < 4),a'_\ell \, (\ell = 1,2),q'$
exemplifies $(p_1,p_2)$ has explicit DOP, so
$M'_0 = M$, tp$(a'_\ell,M'_0,M'_\ell)
= p_\ell$. By uniqueness of primes we are done. \nl
[Saharon: need!] \nl
[Or first find $M_3$ such that 
$(M_1,M_3,a_2) \in K^{3,\text{uq}}_\lambda$ then as $(M_0,M_2,a_2) 
\in K^{3,\text{pr}}_\lambda$ 
there is $\le_{\frak K}$-embedding $f$ of $M_2$ into $M_3$ over $M_0 \cup
\{a_2\}$, let $g \in \text{AUT}({\frak C})$ extend $f \cup 
\text{ id}_{M_1}$
(exists as NF$(M,M_1,M_3,{\frak C})$ and NF$(M,M_1,f(M_2),{\frak C})$), so
$g^{-1}(M_3)$ is as required.]  \hfill$\square_{\scite{705-8.2}}$\margincite{705-8.2}
\enddemo
\bigskip

\proclaim{\stag{705-8.1A} Claim}  1) [${\frak s}$ has {\rm NDOP} see Definition
\scite{705-8.1}.]  Assume
\mr
\item "{$(a)$}"  {\rm NF}$_{\frak s}(M_0,M_1,M_2,M_3)$
\sn
\item  "{$(b)$}"  $(M_0,M_\ell,a_\ell) \in 
K^{3,\text{uq}}_{\frak s}$ for $\ell=1,2$
\sn
\item "{$(c)$}"  $(M_1,M_3,a_2) \in K^{3,\text{uq}}_{\frak s}$.
\ermn
\ub{Then}
\mr
\item "{$(d)$}"  $(M_2,M_3,a_1) \in K^{3,\text{uq}}_{\frak s}$. 
\ermn
2) If {\rm NF}$_{\frak s}(M_0,M_1,M_2,M_3)$ and $(M_0,M_\ell,a_\ell) \in
K^{3,\text{uq}}_{\frak s}$ 
for $\ell = 1,2$ \ub{then} the following are equivalent: 
\mr
\item "{$(\alpha)$}"  $M_3$ is $\le_{\frak s}$-minimal over $M_1 \cup
M_2$
\sn
\item "{$(\beta)$}"  $(M_1,M_3,a_2) \in K^{3,\text{uq}}_{\frak s}$
\sn
\item "{$(\gamma)$}"  $(M_2,M_3,a_1) \in K^{3,\text{uq}}_{\frak s}$. 
\endroster
\endproclaim
\bigskip

\demo{Proof}  1) If this fails, then by \scite{705-5.19} we can find 
$M'_3 <_{\frak s} M_3,(M_2,M'_3,a_1) \in K^{3,\text{uq}}_\lambda$ 
and $b \in M_3 \backslash M'_3$ such that tp$(b_1,M'_3,M_3)$
is not orthogonal to $M_0$.  Hence by the Definition of ``a type is
orthogonal to a model" there is 
$c \in {\frak C}$ such that tp$(c,M'_3,{\frak C}) \in 
{\Cal S}^{\text{bs}}(M'_3)$ does not fork over $M_0$ and
$c \dsize \uplus_{M'_3} b$.  
By the choice of $c,c \nonfork{}{}_{M_0} M'_3$
but $\nonforkin{a_1}{a_2}_{M_0}^{M'_3}$ hence $\{a_1,a_2,c\}$ is 
independent over $M_0$ hence $\{c,a_2\}$ is independent over 
$(M_0,M_1)$.  Also $c$ is independent over $(M_1,M_2)$ hence $c$ is
independent over 
$(M_0,M_3)$, i.e., tp$(c,M_3,{\frak C})$ does not fork over
$M_0$ hence by monotonicity tp$(c,M_3,{\frak C})$ does not fork over $M'_3$.
By the choice of $b$ this contradicts $c \dsize \uplus_{M'_3} b$. \nl
2)  Similar. \nl
\sn
[\ub{Question}: can move to \S5, ignoring regularity?]
\hfill$\square_{\scite{705-8.1A}}$\margincite{705-8.1A}
\enddemo
\bigskip

\definition{\stag{705-8.2A} Definition}  Assume $M_1,M_2 \in K_{\frak s},
{\Cal P} \subseteq {\Cal P}[M_1] =: \cup \{{\Cal S}(N):
N \le_{{\frak K}[{\frak s}]} M_1,N \in K_{\frak s}\}$. \ub{Then}
$M_1 \le_{{\frak s},{\Cal P}}
M_2$ means that $M_1 \le_{\frak s} M_2$ and if $p \in {\Cal P},p_\ell$ the
nonforking extension of $p$ in ${\Cal S}_{\frak s}(M_\ell)$ for 
$\ell =1,2$, \ub{then} $p_2$
is the unique extension of $p_1$ in ${\Cal S}_{\frak s}(M_2)$.
\enddefinition
\bigskip

\proclaim{\stag{705-8.2B} Claim}  Let $M \in K_{\frak s},{\Cal P}$
as in \scite{705-8.2A}. \nl
1) $\le_{{\frak s},{\Cal P}}$ is a partial order on $\{M':M \le_{\frak K} M'
\in K_{\frak s}\}$. \nl
2) If $\langle M_i:i < \delta \rangle$ is 
$\le_{{\frak K},{\Cal P}}$-increasing continuous, $\delta < \lambda^+$ and
$M_\delta = \dbcu_{i < \delta} M_i$ then $i < \delta \Rightarrow M_i
\le_{{\frak s},{\Cal P}} M_\delta$. \nl
3) If $r \in {\Cal S}^{\text{bs}}_{\frak s}(M)$ is 
orthogonal to every $p \in {\Cal P}$ and $(M,N,a) \in 
K^{3,\text{uq}}_{{\frak s}(+)}$, \nl
{\rm tp}$_{\frak s}(a,M,N) = r$ \ub{then} $M \le_{{\frak s},{\Cal P}}
N$.
\nl
4) $M \le_{{\frak s},{\Cal P}} N$ iff there is a 
{\rm pr}-decomposition $\langle M_i,a_i:i < \alpha \rangle$ of $N$
over $M$ (so letting $M_\alpha =: N,M_i$ is $\le_{\frak s}$-increasing
continuous, $M_0 = M_1 (M_i,M_{i+1},a_i) \in K^{3,\text{pr}}_{\frak
s}$ and) {\rm tp}$_{\frak s}(a_i,M_i,N) \pm {\Cal P}$ for every $i < \alpha$.
\endproclaim
\bigskip

\demo{Proof}  Straight.
\enddemo
\bigskip

\proclaim{\stag{705-8.2C} Claim}  1) Assume $p \in {\Cal S}_{\frak s}(M)$ and
\mr
\item "{$(*)$}"  ${\Cal P}$ is a type base for $M$ which means:
{\roster
\itemitem{ $(a)$ }  ${\Cal P} \subseteq {\Cal P}[M] = \cup\{{\Cal S}
^{\text{bs}}_{\frak s}(N):N \le_{\frak s} M$, (so $N \in K_{\frak s})\}$
\sn
\itemitem{ $(b)$ }  for every $q \in {\Cal S}^{\text{bs}}_{\frak s}(M)$ 
there is $r \in {\Cal P}$ not orthogonal to it.
\endroster}
\ermn
\ub{Then} we can find a decomposition $\langle M_i:\ell \le n \rangle,
\langle a_\ell:\ell < n \rangle$ such that
\mr
\widestnumber\item{$(iii)$}
\item "{$(i)$}"  $M_0 = M$,
\sn
\item "{$(ii)$}"  $p$ is realized in $M_n$, and
\sn
\item "{$(iii)$}"  for each $\ell < n$, either
{\rm tp}$_{\frak s}
(a_\ell,M_\ell,M_{\ell +1})$ is a nonforking extension of 
some $q \in {\Cal P}$ \ub{or} {\rm tp}
$(a_\ell,M_\ell,M_{\ell +1})$ is orthogonal 
to $M$ (can be waived if $K^{3,\text{uq}}_{\frak s} 
= K^{3,\text{pr}}_{\frak s}$).
\ermn
2) Assume that ${\frak s}$ has {\rm NDOP} and {\rm NF}$_{\frak s}
(M_0,M_1,M_2,M_3)$ and
$(M_0,M_\ell,a_\ell) \in K^{3,\text{pr}}_{\frak s}$
for $\ell =1,2$ and $(M_1,M_3,a_2) \in K^{3,\text{pr}}_{\frak s}$.  
\ub{Then}
${\Cal S}^{\text{bs}}_{\frak s}(M_1) 
\cup {\Cal S}^{\text{bs}}_{\frak s}(M_2)$
is a type base for $M_3$. \nl
3) [??] If 
$(\langle M_i:i \le \alpha \rangle,\langle a_i:i < \alpha \rangle)$ is
a decomposition of $M_\alpha$ over $M_0$, so 
$(M_i,M_i,a_i) \in K^{3,\text{pr}}_{\frak s},N_i \le_{\frak K}
M_{i+1}$ \ub{then} 
$\{p:\text{ for some } i < \alpha,p \in
{\Cal S}^{\text{bs}}_{\frak s}(M_{i+1})$ and $p \bot M_i\} \cup 
{\Cal S}^{\text{bs}}_{\frak s}(M_0)$ is a type base for $M_\alpha$.
\endproclaim
\bigskip

\demo{Proof}  Easy.
\enddemo
\bigskip

\proclaim{\stag{705-8.3} Claim}  Assume
\mr
\item "{$(a)$}"  $\langle M^*_\ell:\ell < 4 \rangle,
\langle a_\ell:\ell = 1,2 \rangle,q$ are as in \scite{705-8.1}
\sn
\item "{$(b)$}"  $a^k_\ell \in {\frak C}$ realizes
{\rm tp}$(a_\ell,M_0,M_\ell)$ and 
$\langle a^k_\ell:\ell = 1,2$ and $k=1,2 \rangle$ is independent over $M_0$
\sn
\item "{$(c)$}"   $M^k_\ell \le_{\frak K} {\frak C},f^k_\ell$ 
is an isomorphism from $M_\ell$ onto $M^k_\ell$ over $M_0$ for 
$\ell = 1,2$, \nl
$k=1,2$ and $f^k_\ell(a_\ell) = a^k_\ell$ 
\sn
\item "{$(d)$}"  $(M^{k_1}_1,
M^{k_1,k_2},a^{k_2}_2) \in K^{3,\text{pr}}_{\frak s}$
\sn
\item "{$(e)$}"  $f^{k_1,k_2}$ is an isomorphism from $M_3$ onto
$M^{k_1,k_2} <_{\frak K} {\frak C}$ extending $f^{k_1}_1 \cup f^{k_2}_2$
\sn
\item "{$(f)$}"  $q^{k_1,k_2} = f^{k_1,k_2}(q)$.
\ermn
\ub{Then} the types $q^{1,1},q^{1,2},q^{2,1},
q^{2,2}$ are pairwise orthogonal
and each of them orthogonal in $M^k_\ell$ for $\ell = 1,2,k=1,2$.
\endproclaim
\bigskip

\demo{Proof}  Straightforward.
\enddemo
\bigskip

\proclaim{\stag{705-8.3A} Claim}  ${\frak s}$ has {\rm DOP} iff 
${\frak s}^+$ has {\rm DOP}.
\endproclaim
\bigskip

\demo{Proof}  FILL!  Straight.  Decide!
\enddemo
\bn
Our aim is to get strong nonstructure in $\lambda^{++}$ when ${\frak s}$
has DOP.  [Why in $\lambda^{++}$?  We have quite strong independence
but it speaks on $\lambda$-tuples, hence it is hard to get many models
in $\lambda^+$, and if we deal with $K^{\frak s}_{\lambda^{++}}$, why
not ask $\lambda^+$-saturation.  Using \cite{Sh:F569} we hope to deal
with $K^{\frak s}_{\lambda^+}$, too.]
\definition{\stag{705-8.4} Definition}  1) We call ${\frak a}$ an approximation 
or an ${\frak s}$-approximation
(in symbols ${\frak a} \in {\frak A}$) \ub{if} ${\frak a}$ consists of the
following objects, satisfying the following demands
\mr
\item "{$(a)$}"  $I^{\frak a}_1,I^{\frak a}_2$ disjoint index sets of 
cardinality $\le \lambda^+$
\sn
\item "{$(b)$}"  $R_{\frak a} \subseteq I^{\frak a}_1 \times I^{\frak a}_2$,
we write $sR_{\frak a}t$ for $(s,t) \in R_{\frak a},\neg sR_{\frak a} t$ 
for $s \in I^{\frak a}_1$, \nl
$t \in I^{\frak a}_2,(s,t) \notin R_{\frak a}$
\sn
\item "{$(c)$}"  $M^{\frak a}_\ell$ for $\ell < 4,a^{\frak a}_\ell$ for
$\ell = 1,2$ and $q_{\frak a}$ exemplifying DOP
\sn
\item "{$(d)$}"  $M^{\frak a} \in K_{\lambda^+}$ saturated (so $\in
K_{\lambda^+}[{\frak s}^+]$) such that $M^{\frak a}_0 \le_{\frak K} 
M^{\frak a}$
\sn
\item "{$(e)$}"  $f^{\frak a}_{\ell,t}$ an $\le_{\frak K}$-embeddiing of
$M^{\frak a}_\ell$ into $M^{\frak a}$ for $\ell =1,2,t \in I^{\frak a}_\ell$
and we let $M^{\frak a}_{\ell,t} = f^{\frak a}_{\ell,t}(M^{\frak a}_\ell),
a^{\frak a}_{\ell,t} = f^{\frak a}_{\ell,t}(a^{\frak a}_\ell)$
\sn
\item "{$(f)$}"  $\{a^{\frak a}_{\ell,t}:\ell=1,2$ and $t \in
I^{\frak a}_\ell\} \subseteq \bold I_{M^{\frak a}_0},M^{\frak a}$ is 
independent over $M^{\frak a}_0$; hence \nl
$\langle M^{\frak a}_{\ell,t}:\ell = 1,2,t \in I^{\frak a}_\ell)$ is 
independent
\sn
\item "{$(g)$}"  if $s R_{\frak a} t$ then $f^{\frak a}_{s,t}$ is a
$\le_{\frak K}$-embedding of $M^{\frak a}_3$ into $M^{\frak a}$ extending
$f^{\frak a}_{1,s} \cup f^{\frak a}_{2,t}$; we let $M^{\frak a}_{s,t} =
f^{\frak a}_{s,t}(M_3),q^{\frak a}_{s,t} = f^{\frak a}_{s,t}(q^{\frak
a})$.
\ermn
2) For an approximation ${\frak a}$ let
${\Cal P}^+_{\frak a} = \{q^{\frak a}_{s,t}:sR_{\frak a}t$
hence $s \in I^{\frak a}_1$ and $t \in I^{\frak a}_2\}$ \nl
${\Cal P}^-_{\frak a} = \{f(q_{\frak a}):
\text{for some } f$ and $(s,t) \in I^{\frak a}_1
\times I^{\frak a}_2$ such that $\neg sR_{\frak a} t,f$ is a
$\le_{\frak K}$-embedding of $M^{\frak a}_3$ into $M$ extending
$f^{\frak a}_{1,s} \cup f^{\frak a}_{2,t}\}$.
\nl
3) Let ${\frak A} = {\frak A}_s$ be the class of 
${\frak s}$-approximations. \nl
4) We call ${\frak a}^*$ a DOP witness if it consists of $M^{{\frak
a}^-}_\ell(\ell < 4),a^{{\frak a}^-}_\ell(\ell =1,2),q^{{\frak a}^-}$
which are as above.  If ${\frak b}$ is an approximation let ${\frak
b}^-$ be defined naturally.
\enddefinition
\bigskip

\definition{\stag{705-8.5} Definition}  1) If ${\frak a},{\frak b}$ are
approximations let ${\frak a} \le {\frak b}$ means:
\mr
\item "{$(\alpha)$}"  $M^{\frak a}_\ell = M^{\frak b}_\ell$ for $\ell < 4,
a^{\frak a}_\ell = a^{\frak b}_\ell$ for $\ell = 1,2,q^{\frak a} = 
q^{\frak b}$
\sn
\item "{$(\beta)$}"  $I^{\frak a}_\ell \subseteq I^{\frak a}_b$ for
$\ell=1,2$ and $R_{\frak a} = R_{\frak b} \cap (I^{\frak a}_1 \times 
I^{\frak a}_2)$
\sn
\item "{$(\gamma)$}"  for $\ell = 1,2,t \in I^{\frak a}_1$ we have
$f^{\frak a}_{\ell,t} = f^{\frak b}_{\ell,t}$
\sn 
\item "{$(\delta)$}"  for $(s,t) \in R^{\frak a}$ 
we have $f^{\frak a}_{s,t} =
f^{\frak b}_{s,t}$
\sn
\item "{$(\varepsilon)$}"  $M^{\frak a} \le_{\frak K} M^{\frak b}$,
moreover $M^{\frak a} \le_{{\frak K},{\Cal P}^-_{\frak a}} M^{\frak b}$.
\ermn
2) If $\langle {\frak a}_\zeta:\zeta < \delta \rangle$ is $\le$-increasing
in ${\frak A}$ and $\delta < \lambda^{++}$ let their union ${\frak a} = 
\dbcu_{\zeta < \delta} {\frak a}_\zeta$ be defined by $I^{\frak a}_\ell =
\dbcu_{\zeta < \delta} I^{{\frak a}_\ell},R_{\frak a} = \dbcu_{\zeta <
\delta} R_{{\frak a}_\zeta},f^{\frak a}_{\ell,t} = f^{{\frak a}_\zeta}
_{\ell,t}$ for $\zeta < \delta$ large enough, $f^{\frak a}_{s,t} =
f^{{\frak a}_\zeta}_{s,t}$ for $\zeta < \delta$ large enough when
$sR_{\frak a}t$ and $M^{\frak a} = \cup\{M^{{\frak a}_\zeta}:
\zeta < \delta\}$.
\enddefinition
\bn
Below we restrict ourselves to $K_{{\frak s}(+)}$ for the
application we have in mind but $K_{\lambda^+}$ would be also O.K.

\proclaim{\stag{705-8.6} Claim}  1)  $({\frak s},\le)$ is a partial order. \nl
2) If $\langle {\frak a}_\zeta:\zeta < \delta \rangle$ is increasing in
${\frak A},\delta < \lambda^{++}$ \ub{then} 
${\frak a} = \dbcu_{\zeta < \delta}{\frak a}_\zeta$ belong 
to ${\frak A}_{\frak s}$ is the {\rm lub} of the sequence.
\endproclaim
\bigskip

\demo{Proof}  Straight.
\enddemo
\bigskip

\proclaim{\stag{705-8.7} Claim}  1) If ${\frak a} \in {\frak A}_{\frak s},
p \in {\Cal S}^{\text{bs}}_{\frak s}(M^{\frak a})$ is orthogonal 
to every $q \in {\Cal P}^-_{\frak a}$ and
$(M^{\frak a},N,a) \in K^{3,\text{uq}}_{{\frak s}(+)}$ \ub{then} for some
${\frak b} \in {\frak A}_{\frak s}$ we have 
${\frak a} \le {\frak b}$ and $M^{\frak b} = N$. \nl
2) Assume ${\frak a} \in {\frak A},\ell(*) \in \{1,2\},Y \subseteq
I^{\frak a}_{3 - \ell(*)},t^* \notin I^{\frak a}_{\ell(*)}$ and $\langle
M_\alpha:\alpha < \lambda^+ \rangle$ is a representation of $M^{\frak a}$
such that $M_0 = M^{\frak a}_0$ (and of course $M_{\alpha + 1}$ is
$(\lambda,*$)-brimmed 
over $M_\alpha$ in ${\frak K}_\lambda$).  \ub{Then} we can find
${\frak b},a,\langle N_\alpha:\alpha < \lambda^+ \rangle$ such that:
\mr
\item "{$(A)(a)$}"  $N_\alpha$ is $\le_{\frak K}$-increasing continuous
in $K_\lambda,N_{\alpha +1}$ is $(\lambda,*)$-brimmed over $N_\alpha$
\sn
\item "{$(b)$}"  {\rm NF}$_{\frak s}
(M_\alpha,N_\alpha,M_{\alpha +1},N_{\alpha +1})$
\sn
\item "{$(c)$}"  $(M_\alpha,N_\alpha,a) \in K^{3,\text{uq}}_{\frak s}$
\sn
\item "{$(d)$}"  $N_0$ is isomorphic to $M^{\frak a}_{\ell(*)}$ over
$M^{\frak a}_0$
\sn
\item "{$(e)$}"  $(\dbcu_\alpha M_\alpha,\dbcu_\alpha N_\alpha,a) \in
K^{3,\text{pr}}_{{\frak s}(+)}$
\mn
\item "{$(B)(a)$}"  ${\frak b} \in {\frak A},{\frak a} \le {\frak b}$
\sn
\item "{$(b)$}"  $I^{\frak b}_{\ell(*)} = I^{\frak a}_{\ell(*)} \cup
\{t^*\},I^{\frak b}_{3 - \ell(*)} = I^{\frak a}_{3 - \ell(*)}$
\sn
\item "{$(c)$}"  $R^{\frak b}$ is $R^{\frak a} \cup \{\langle t^*,s
\rangle:s \in Y\}$ if $\ell(*)=1$ and is $R^{\frak a} \cup \{\langle
s,t^* \rangle:s \in y\}$ if $\ell(*)=2$
\sn
\item "{$(d)$}"  $f^{\frak a}_{\ell(*),t^*}$ is an isomorphism from
$M^{\frak a}_{\ell(*)}$ onto $N_0$ mapping $a^{\frak a}_{\ell(*)}$ to
$a$.
\endroster
\endproclaim
\bigskip

\demo{Proof}  1) Easy. \nl
2) First choose $a,N_\alpha$ to satisfy (A).  Then the choice of
${\frak b}$ is actually described in (B); the orthogonality hold by
\scite{705-8.2B}.  \hfill$\square_{\scite{705-8.7}}$\margincite{705-8.7}
\enddemo
\bigskip

\proclaim{\stag{705-8.8} Claim}  Let ${\frak a}$ be a {\rm DOP} witness
and $R \subseteq \lambda^{++} \times \lambda^{++}$
be given.  For $\alpha < \lambda^{++}$ let $I^\alpha_1 = \{i:3i+1 \le
\alpha\},I^\alpha_2 = \{i:3i+2 \le \alpha\},R^*_\alpha =  R \cap (I^\alpha_1
\times I^\alpha_2)$.  We can find $\langle {\frak a}^\alpha:\alpha <
\lambda^{++} \rangle$ such that
\mr
\item "{$(a)$}"  ${\frak a}_\alpha \in {\frak A}_{\frak s}$ 
is increasing continuous and ${\frak a}^-_\alpha = b$
\sn
\item "{$(b)$}"  $(I^{{\frak a}^\alpha}_1,I^{{\frak a}^\alpha}_2,
R_{{\frak a}^\alpha}) = (I^\alpha_1,I^\alpha_2,R_\alpha)$,
\sn
\item "{$(c)$}"  for $(s,t) \in R_\alpha$ for arbitrarily large $\beta \in
(\alpha,\lambda^{++})$ (by some bookkeeping), some $b \in 
M^{{\frak a}^{3 \beta +3}} \backslash M^{{\frak a}^{3\beta +2}}$
the type {\rm tp}$_{\frak s}(b,M^{{\frak a}^{3 \beta +2}},
M^{{\frak a}^{3 \beta +3}})$ is a nonforking
extension of $q^{{\frak a}_\alpha}_{s,t}$
\sn
\item "{$(d)$}"  $a_\alpha$ depends just on $(I^\alpha_1,I^\alpha_2,
R_\alpha)$
\sn
\item "{$(e)$}"  the 
universe of $M^{{\frak a}^\alpha}$ is $\gamma_\alpha <
\lambda^{++}$ (really $\gamma_\alpha = \lambda^* 
\times(1 + \alpha)$ is O.K.
for nontrivial cases.
\endroster
\endproclaim
\bigskip

\demo{Proof}  We choose ${\frak a}^\alpha$ by induction on $\alpha$.  For
$\alpha = 0$ this is trivial, for $\alpha$ limit by \scite{705-8.6}(2), for
$\alpha = 3 \beta +1$ by \scite{705-8.7}(2) for $\ell(*) = 1$, for $\alpha =
3 \beta +2$ by \scite{705-8.7}(2) for $\ell(*) = 2$ for $\alpha = 3 \beta +3$
bookkeeping gives as a pair $(s_\alpha,t_\alpha)$ and we use \scite{705-8.7}(1).
\nl
${{}}$  \hfill$\square_{\scite{705-8.8}}$\margincite{705-8.8}
\enddemo
\bigskip

\proclaim{\stag{705-8.9} Claim}  In \scite{705-8.8} we can add: letting $M^* = 
\dbcu_\alpha\{M^{{\frak a}^\alpha}:\alpha < \lambda^{++}\}$
\mr
\item "{$(*)$}"  for $(s,t) \in \lambda^{++} \times \lambda^{++}$, the 
following are equivalent
{\roster
\itemitem{ $(\alpha)$ }  $(s,t) \in R$
\sn
\itemitem{ $(\beta)$ }  {\rm dim}$(q^{{\frak a}^\alpha}_{s,t},M^*) 
= \lambda^{++}$
when $\alpha = { \text{\rm Max\/}}\{3s+1,3t+2\}$ 
that is, there is a sequence 
$\langle b_\gamma:\gamma < \lambda^{++} \rangle$ independent in
$(M^{{\frak a}^\alpha},M^*)$ 
of elements realizing $q^{{\frak a}^\alpha}_{s,t}$
for any $\alpha > 3s+1,3t+2$ \nl
(check the existence of formal definition, \S5 defines dimension)
\sn
\itemitem{ $(\gamma)$ }  there is a $\le_{\frak K}$-embedding $f$ of
$M^{{\frak a}^\alpha}_3$ into $M^*$, extending $f^{{\frak a}^\alpha}_{1,s}
\cup f^{{\frak a}^\alpha}_{2,t}$ such that {\rm dim}
$(f(q^{{\frak a}^\alpha}),M^*) = 
\lambda^{++}$ for $\alpha = { \text{\rm Max\/}}\{3s+1,3t+2\}$
\sn
\itemitem{ $(\delta)$ }  for no 
$\alpha < \lambda^{++}$ and $f$ as in $(\gamma)$,
we have: $q^* \in {\Cal S}^{\text{bs}}_{\frak s}
(M^{{\frak a}^\alpha})$, the nonforking
extension of $f(q_{\frak b})$ in ${\Cal S}^{\text{bs}}(M^{{\frak a}^\alpha})$
satisfies: for every $\beta \in (\alpha,\lambda^{++}),q^*$ has a unique
extension in ${\Cal S}^{\text{bs}}_{\frak s}(M^{{\frak a}_\beta})$ \nl
(Saharon - check in the reflections].
\endroster}
\endroster
\endproclaim
\bigskip

\demo{Proof}  Easy.
\enddemo
\bigskip

\proclaim{\stag{705-8.10} Claim}  $[2^{\lambda^+} < 2^{\lambda^{++}}]$.  If
${\frak s}$ has {\rm DOP} 
then $I(\lambda^{++},K^{{\frak s}(+)}) = 2^{\lambda^{++}}$.
\endproclaim
\bigskip

\remark{\stag{705-8.11} Remark}  1) It is a strong nonstructure 
(i.e., neither like for deepness, no even like unsuperstable. \nl
2) We can in \scite{705-8.8}, \scite{705-8.9} restrict more the types
realized. \nl
3) We may use here \cite[III]{Sh:e}.  FILL!
\endremark
\bigskip

\demo{Proof} We use the construction above in the framework of
\cite[\S3]{Sh:576}.
\enddemo
\newpage

\head {\S12 Brimmed Systems} \endhead  \resetall \sectno=12
 \spuriousreset
\bn
Saharon: ${\frak s}[\mu]$ the canonical lifting!!!
\sn
This section generalizes \cite{Sh:87b}, \cite[XII,\S4,\S5]{Sh:c}.  Here
every system is in the context of some good frame ${\frak s}$ and
usually we look at models of cardinality $\lambda = \lambda_{\frak s}$
(in this section), but we vary ${\frak s}$.  

The adoption of ``${\frak K}_{\frak s}$ categorical in $\lambda$" (in
\scite{705-12.1}) is very helpful here but there is a price: when we shall
work on ``all $\lambda^{+ \omega}$-models in $K^{\frak s}$" we cannot
just quote the results.  Note that this restriction fits well the
thesis that the main road is first to understand the quite saturated
models.
\bn
\margintag{705-12.0}\ub{\stag{705-12.0} Convention}:  1) Without loss of generality 
always $I \cap {\Cal P}(\text{Dom}(I)) = \emptyset$!! \nl
2) In the cases we assume categoricity of $K^{\frak s}$
in $\lambda_{\frak s}$ we add $*$ (e.g., \scite{705-12.1R}(2)).
\bigskip

\demo{\stag{705-12.1} Hypothesis}  
\mr
\item "{$(a)$}"  ${\frak s}$ is a good $\lambda$-frame
\sn
\item "{$(b)$}"  $\bot = \underset{\text{wk}} {}\to \bot$ is well defined 
\sn
\item "{$(c)$}"  ${\frak s}$ has primes (in the sense of 
$K^{3,\text{qr}}_{\frak s})$ and they are unique
\sn
\item "{$(d)$}"  ${\frak s}$ is full or at least satisfies the 
conclusion  of it ``has enough regulars" (see Definition and Claim \sciteu{705-xxX})
\sn
\item "{$(e)^*$}"  ${\frak K}_{\frak s}$ is categorical in $\lambda =
\lambda_{\frak s}$
\endroster
\enddemo
\bn
Recall (by \sciteu{705-xxX})
\proclaim{\stag{705-12.1R} Claim}  1) The following \footnote{we may be
interested in the case we replace $K^{3,\text{bu}}_{\frak s}$ by
$K^{3,\text{qr}}_{\frak s}$, but we have troubles enough} condition on
$(M,N,\bold J)$ are equivalent 
\mr
\item "{$\circledast_1$}"  $(M,N,\bold J) \in
K^{3,\text{bu}}_{\frak s}$ which means that
$(M,N,\bold J) \in K^{3,\text{bs}}_{\frak s}$, and
$\bold J$ is maximal (used in \sciteu{705-xxX})
\sn
\item "{$\circledast_2$}"  $(M,N,\bold J) \in K^{3,\text{bs}}_{\frak s}$ 
and if $M + \bold J \subseteq N' \le_{\frak s}
N,b \in N \backslash N'$, {\rm tp}$(b,N',N) \in {\Cal S}^{\text{bs}}(N')$ 
\ub{then} {\rm tp}$(b,N',N) \bot M$, 
\sn
\item "{$\circledast_3$}"  $(M,N,\bold J) \in K^{3,\text{vq}}_{\frak s}$. 
\ermn
2) We have $K^{3,\text{bu}}_{\frak s} = K^{3,\text{qr}}_{\frak s}$. \nl
3) $K^{3,\text{qr}}_{\frak s} \subseteq K^{3,\text{bg}}_{\frak s}$.
\nl
4) For every $M \in K_{\frak s},i^* < \lambda^+_{\frak s}$ and $p_i
\in {\Cal S}^{\text{bs}}_{\frak s}(M)$ for $i < i^*$ we can find $N$
and $a_i \, (i < i^*)$ such that: $(M,N,\{a_i:i < i^*\}) \in
K^{3,\text{qr}}_{\frak s} \subseteq
K^{3,\text{bu}}_{\frak s}$ and $p_i = { \text{\rm tp\/}}_{\frak
s}(a_i,M,N)$ and, of course, $\langle a_i:i < i^* \rangle$ is without
repetition. 
\endproclaim
\bigskip

\definition{\stag{705-12.3} Definition}  Let $u_*$ be a set (usually
finite) and $I$ a family of finite subsets of $u_*$ (so for $u_*$
finite this is automatic) satisfying $I \subseteq {\Cal P}(u_*)$
is downward closed; $I,J$ will denote such sets in this 
section; let Dom$(I) = \cup\{u:u \in I\}$. \nl
1) We say $\bold s$ is an $I$-system or $(I,{\frak s})$-system or
$I$-system for ${\frak s}$ \ub{if}:
\mr
\item "{$(a)$}"  $\bold s$ consists of $M_u$ (for $u \in I$), and 
$h_{v,u}$ (for $u \subseteq v \in I$) (mappings, with
$h_{u,u} = \text{ id}_{M_u}$ so we may ignore $\langle h_{u,u}:u \in I
\rangle$ when defining $\bold s$)
\sn
\item "{$(b)$}"  $M_u \in K_{\frak s}$ for $u \in I$
\sn
\item "{$(c) $}"  if $u \subseteq v \in I$ then $h_{v,u}$ is a
$\le_{\frak s}$-embedding of $M_u$ into $M_v$,  
and the diagram of the $h_{u,v}$'s commute and we let
$M_{v,u} =: h_{v,u}(M_u)$ and recall $h_{u,u} = \text{ id}_{M_u}$; so if $u
\subseteq v \subseteq w \in I$ we have $M_{w,u} \le_{\frak s} M_{w,v}$
and $M_{u,u} = M_u$.
\ermn
2) We say $\bold s$ is a $(\mu,I)$-system or $(\mu,I,{\frak s})$-system \ub{if} we replace (b) by
\mr
\item "{$(b)^+$}"  $M_u \in K^{\frak s}_\mu$.
\ermn
Similarly for $(\ge \mu,I,{\frak s})$, etc. \nl
3) We shall write $u^{\bold s}_*$ for Dom$(I)$ and $h^{\bold s}_{v,u},
I^{\bold s},M^{\bold s}_u,$ for $h_{u,v},I,M_u$, respectively.  
If $h_{v,u} = \text{ id}_{M_u}$ 
for $u \subseteq v \in I^{\bold s}$ and $M^s_u \cap M^s_v = M^{\bold
s}_{v \cap u}$ for $u,v \in I$ we call ${\bold s}$ normal.
\nl
4) We say $\bar g$ is an isomorphism from the $I$-system ${\bold s}^1$ onto
the $I$-system ${\bold s}^2$ \ub{if} 
$\bar g = \langle g_u:u \in I \rangle,g_u$
is an isomorphism from $M^{{\bold s}^1}_u$ onto $M^{{\bold s}^2}_4$ such that
$u \subseteq v \in I \Rightarrow h^{{\bold s}^2}_{v,u} \circ g_u = g_v \circ
h^{{\bold s}_1}_{v,u}$. If ${\bold s}^1,{\bold s}^2$ are normal we may say
$g = \dbcu_u g_u$ is an isomorphism from $\bold s^1$ onto $\bold s^2$.  
Similarly $\bar g$ is a $\le_{\frak s}$-embedding of $\bold s^1$ into
$\bold s^2$ if $g_u(M^{{\bold s}^1}_u) \le_{\frak s}
M^{{\bold s}^2}_u$ or \footnote{the difference is meaningful only if
$M^{\bold s^2}_u \in K^{\frak s} \backslash K_{\frak s}$}
$g_u(M^{{\bold s}^1}_u) \le_{{\frak K}[{\frak s}]} 
M^{{\bold s}^2}_u$ for $u \in I$ and $u \subset v
\in I \Rightarrow h^{\bold s^2}_{v,u} \circ g_u =  g_v \circ 
h^{\bold s^1}_{v,u}$.  Let $\bold s_1 \le_{{\frak K}_{\frak s}} \bold
s_2$ if $\langle \text{id}_{M^{\bold s_1}_u}:u \in I \rangle$ is a
$\le_{\frak s}$-embedding of $M^{\bold s^1}_u$ into $M^{\bold s^2_u}$
for $u \in I$, similarly $\le_{{\frak K}[{\frak s}]}$.
\nl
5) We say $\bar f$ is an embedding of an $I$-system $\bold s$ to a
model $M$ \ub{if} $\bar f = \langle f_u:u \in I \rangle,f_u$ is a 
$\le_{\frak s}$-embedding of $M_u$ into $M$ and $u \subseteq v \in I
\Rightarrow f_u = f_v \circ g_{v,u}$. \nl
6) Similarly for $(\ge \mu,I)$-systems, $(\ge \mu,I,{\frak s})$-systems, so
we use $\le_{{\frak K}[{\frak s}]}$-embeddings.
\enddefinition
\bn
The important case here is:
\definition{\stag{705-12.4} Definition}  1) We say that $\bold d$ is an
expanded stable $(I,{\frak s})$-system or $I$-system or
$(\lambda,I)$-system or $(\lambda,I,{\frak s})$-system \ub{if} 
it consists of $\bold s$ and $\bold J_{v,u}$ for 
$u \subseteq v \in I$ such that:
\mr
\item "{$(a)$}"  $\bold s = \langle M_u,h_{v,u}:u \subseteq v \in I
\rangle$ is an $(\lambda,I,{\frak s})$-system
\sn
\item "{$(b)$}"  $\bold J_{v,u} \subseteq \bold I_{M_{v,u},M_v}
\backslash \cup\{M_{v,w}:w \subset v\}$ for $u \subseteq v \in I$ 
so $\bold J_{u,u} = \emptyset$.
\ermn
such that
\mr
\item "{$(c)$}"  if $u_0 \subset u_1 \subset u_2 \in I$ and $c \in \bold
J_{u_2,u_1}$ \ub{then} 
tp$(c,M_{u_2,u_1},M_{u_2})$ is orthogonal to $M_{u_2,u_0}$, recalling
that $M_{u_2,u_1} = f_{u_2,u_1}(M_{u_1})$
\sn
\item "{$(d)$}"  $\bold J_{v,u}$ is a maximal subset of $\{c \in M_v:
c \notin \cup \{\bold J_{v,w}:w \subset u\}$ and tp$_{\frak
s}(c,M_{v,u},M_u)$ belongs to ${\Cal S}^{\text{bs}}_{\frak s}(M)$ and
is orthogonal to $M_{v,w}$ for every $w \subset u\}$ such that $\bold
J_{v,u} \cup \bigcup\{\bold J_{v,w}:w \subset u\}$ is independent in
$(M_{v,u},M_v)$.
\ermn
1A) For $u_0 \subseteq u_1 \subseteq u_2 \in I$ we define
(if $u_2 = u_1$
we may omit it, this catch the ``main action", so $\bold J^0_{u_1,u_1,u_0} = 
\bold J^0_{u_1,u_0} = \bold J_{u_1,u_0}$)
\mr
\item "{$(\alpha)$}"  $\bold J^0_{u_2,u_1,u_0} = \{h_{u_2,u_1}(c):
c \in \bold J_{u_1,u_0}\}$
\sn
\item "{$(\beta)$}"  $\bold J^1_{u_2,u_1,u_0} = \cup\{\bold
J^0_{u_2,u_1,u}:u \subseteq u_0\}$
\sn
\item "{$(\gamma)$}"  $\bold J^2_{u_2,u_1,u_0} = \cup\{\bold
J^0_{u_2,w_1,w_0}:w_0 \subseteq w_1 
\subseteq u_1,w_0 \subseteq u_0$ and $w_1 \subset u_0\}$
\ermn
If we omit clause (d) we say ``an expanded $(I,{\frak s})$-system".
\nl
2) We say $\bold d$ is normal if each $h_{v,u}$ is the
identity (on $M_u$) and $M_u \cap M_v = M_{u \cap v}$, that is $\bold
s$ is normal.  
\nl
3) For an expanded stable $(\lambda,I)$-system $\bold d$ and $M \in
K_{\frak s}$ we say that $\bar f$ is an embedding 
(or $\le_{\frak s}$-embedding) of $\bold d$ into $M$ when:
\mr
\item "{$(A)$}"  $\bar f$ embeds $\bold s^{\bold d}$ into $M$, i.e.
{\roster
\itemitem{ $(a)$ }  $\bar f = \langle f_u:u \in I \rangle$
\sn
\itemitem{ $(b)$ }  $f_u$ is a $\le_{\frak s}$-embedding of $M^{\bold d}_u$
into $M$
\sn
\itemitem{ $(c)$ }  if $u \subset v \in I$ then $f_u =f_v \circ
h^{\bold d}_{v,u}$
\sn
\itemitem{ $(d)$ }  if $u \in I$ then $\cup\{f_u(\bold J^{2,\bold d}_{v,u}):v$
satisfies $ u \subset v \in I\}$, is an independent set in 
$(f_u(M^{\bold d}_u),M)$ and, of course, $u \subseteq v_1 \in I \cap u
\subseteq _2 \in I \cap v_1 \ne v_2 \Rightarrow f_u(J^{2,\bold
d}_{v_1,u}) \cap f_u(\bold J^{2,\bold d}_{v_2,u}) = \emptyset$
\endroster}
\ermn
4) We say that a normal expanded stable $(\lambda,I)$-system $\bold d$ is
embedded into $M \in K_{\frak s}$ \ub{if} $(\bold d$ is normal and)
$\bar f = \langle f_u:u \in I \rangle$ is an embedding of $\bold d$
into $M$ when we choose $f_u = \text{ id}_{M^{\bold d}_u}$.
\nl
5) We say $\bold d$ an expanded stable $I$-system is explicit or
 strongly regular \ub{if}:
\mr
\item "{$(a)$}"  if $u \subset v \in I$ and $c \in \bold J^{\bold
d}_{v,u}$ then tp$_{\frak s}(c,M^{\bold d}_{v,u},M^{\bold d}_v)$ is
regular (see \scite{705-11b.2}(2))
\sn
\item "{$(b)$}"  if $c_1 \ne c_2 \in \bold J^{\bold d}_{v,u}$ and
$u \subset v \in
I$ \ub{then} tp$_{\frak s}(c_1,M^{\bold d}_{v,u},M^{\bold d}_v)$,
tp$_{\frak s}(c_2,M^{\bold d}_{v,u},M^{\bold d}_v)$ are equal or
orthogonal. 
\endroster
\enddefinition
\bn
\margintag{705-12.5}\ub{\stag{705-12.5} Notation}: Above we let $\bold s^{\bold d}
= \bold s[\bold d],M^{\bold d}_u = M^{\bold s}_u,h^{\bold d}_{v,u} =
h^{\bold s}_{v,v},\bold J^{\bold d}_{v,u} = \bold J_{v,u},
\bold J^{\ell,\bold d}_{u_2,u_1,u_0} = 
\bold J^\ell_{u_2,u_1,u_0}$ but we do not write the 
superscript $\bold d$ when clear from the context. 
\bigskip

\definition{\stag{705-11b.2} Definition}
1) We say $\bar f$ is an isomorphism from the expanded stable 
$I$-system $\bold d^1$ onto the expanded stable 
system $\bold d^2$ \ub{if} $\bar f$ is an
isomorphism from $\bold s^{{\bold d}^1}$ onto $\bold s^{{\bold d}^2}$
and $f_v$ maps $\bold J^{{\bold d}^1}_{v,u}$ onto $\bold J^{{\bold
d}^2}_{v,u}$ for $u \subseteq v \in I$. \nl
1A) For expanded stable $(\lambda,I)$-systems $\bold d_1,\bold d_2$, we 
say $\bar f$ is an $\le_{\frak s}$-embedding of $\bold d_1$ into
$\bold d_2$ \ub{if} 
\mr
\item "{$(\alpha)$}"  $\bar f$ is an $\le_{\frak s}$-embedding 
of $\bold s_1$ into $\bold s_2$
\sn
\item "{$(\beta)$}"  for $u \subset v \in I,f_v$ maps $\bold J^{{\bold
d}_1}_{v,u}$ into $\bold J^{{\bold d}_2}_{v,u}$
\sn
\item "{$(\beta)^+$}"  moreover, if 
$u \subset v \in I$ and $c \in \bold J^{{\bold d}_1}_{v,u}$ then 
tp$_{\frak s}(f_v(c),M^{{\bold d}_2}_{v,u},M^{{\bold d}_2}_v)$ does not
fork over $f_v(M^{{\bold d}_1}_{v,u})$.
\ermn
2) An expanded stable $(\lambda,I)$-system $\bold d$ is called regular
\ub{if}
\mr
\item "{$(f)$}"  if $u \subseteq v \in I,c \in \bold J^{\bold
d}_{v,u}$ \ub{then} tp$(c,M_{v,u},M_v)$ is regular.
\ermn
3) A system $\bold s$ is called stable \ub{if} for some expanded stable
system $\bold d$ we have $\bold s = \bold s^{\bold d}$ (called an
expansion of $\bold s$), we call the $\bold J^{\bold d}_{v,u}$'s
witnesses.  Similarly for other properties.
\nl
4) A $(\lambda,I)$-system $\bold s$ (or an expanded stable system
$\bold d$ with $\bold s^{\bold d} = \bold s$) is called very brimmed \ub{if}:
\mr
\item "{$(g)$}"  for every $v \in I,\bold s$ is very brimmed at $v$ which means
\sn
\item "{$(g)^0_v$}"  $M^{\bold s}_v$ is brimmed
\sn 
\item "{$(g)^1_v$}"  $M^{\bold s}_v$ is
$(\lambda_{\frak s},*)$-brimmed over $\cup \{M^{\bold s}_{v,u}:u
\subset v\}$.
\ermn
5) A stable $I$-system $\bold s$ (or an expanded stable system
$\bold d$ with $\bold s^{\bold d} = \bold s$) is weakly brimmed
\ub{if} for every $v \in I$ it is weakly brimmed at $v$ which means
\mr
\item "{$(h)^-_v$}"   $M^{\bold s}_v$ is brimmed
\footnote{note that if $K_{\frak s}$ is categorical (in $\lambda_{\frak
s}$) then this follows from the other demands}. 
\ermn
6) An expanded stable system $\bold d$ is called brimmed \ub{if} for every
$v \in I_s$  it is brimmed at $v$ which means
\mr
\item "{$(h)^0_v$}"  $M^{\bold d}_v$ is brimmed
\sn 
\item "{$(h)^1_v$}"  if $u \subset v \in I$ and $p \in 
{\Cal S}^{\text{bs}}(M_{v,u})$ is orthogonal to $M_{v,w}$ for every 
$w \subset u$, \ub{then} the set 
$\{c \in \bold J_{v,u}:\text{tp}_{\bar{\frak s}}(c,M_{v,u},M_v) \pm
p\}$ has cardinality $\|M_v\|$.
\ermn
7) An $I$-system is [very][weakly] brimmed \ub{if} there is a [very][weakly]
brimmed expanded stable system $\bold d$ expanding it (so the system
is stable).
\enddefinition
\bigskip

\definition{\stag{705-12.bf.3} Definition}  1) We say $\bold s$ is a
$(I,{\frak s})^\ell$-system or a brimmed$^\ell \,I$-system
when it is a $I$-system and: if $\ell=0$ no additional demands; 
if $\ell =1$, it is stable; if $\ell =2$ it is a stable system which
is weakly brimmed, that is each $M^{\bold s}_u$ is brimmed; 
if $\ell=3$, it is stable and brimmed; if $\ell=4$, it is a
very brimmed stable $(I,{\frak s})$-system.
\nl
2) Similarly for $\bold d$, an expanded stable $(I,{\frak s})$-system
(so $\ell=0,\ell=1$ become equivalent). \nl
3) We say $M$ is brimmed$^\ell$ \ub{if} letting $I =
\{\emptyset\},M_\emptyset =M$, we get that $\langle M_t:t \in I
\rangle$ is a brimmed$^\ell$ system. \nl
4) For $\ell=1,2,3,4$ we say that an expanded stable $I$-system $\bold d$ 
is brimmed$^\ell$ in $u \in I$ if the demand in Definition
\scite{705-11b.2} holds for $u$ (so for $\ell=1$: no demand).
\enddefinition
\bn
\ub{Notation}:  Let ${\Cal P}^-(u_*) = \{v \subseteq u_*:v \ne u_*\}$. 
\bigskip

\proclaim{\stag{705-12b.3} Claim} Let $\bold d$ be an
expanded stable $I$-system. \nl
0) If $u \subset v \in I$ \ub{then} $(M^{\bold d}_{v,u},M^{\bold d}_v,
\bold J^{2,\bold d}_{v,u})$ belongs to $K^{3,\text{bu}}_{\frak s}$;
equivalently (see \scite{705-12.1R}) it belongs to
$K^{3,\text{vq}}_{\frak s}$, recalling 
$\bold J^{\bold d}_{v,v,u} = \bold J^{2,\bold d}_{v,u}$.
\nl
1) If $u_1 \subseteq u \in I,u_2 \subseteq u,u_0 = u_1 \cap u_2$,
\ub{then} {\rm NF}$_{\frak s}(M_{u,u_0},M_{u,u_1},M_{u,u_2},M_u)$. \nl
1A) If $u_0,u_1,u_2 \subseteq v \in I,u_0 \subseteq u_1,c \in \bold J^{\bold
d}_{v,u_1,u_0}$ and $u_0 \nsubseteq u_2$ \ub{then}
{\rm tp}$(c,M_{v,u_0},M_{v,u_1})$ is orthogonal to $M_{v,u_2}$. \nl
2)  If $p \in {\Cal S}^{\text{bs}}
(M_{v,u})$ is regular, \ub{then} there is a
unique set $w \subseteq u$ such that $p$ is not orthogonal to
$M_{v,w}$ but $p$ is orthogonal to $M_{v,w'}$ whenever $w' \subseteq u
\and w \nsubseteq w'$ hence even when $w' \subseteq v \and w
\nsubseteq w'$.
\nl
3) $\bold d$ is isomorphic to some normal $\bold d'$ (see Definition
\scite{705-12.3}(3)).
\nl
4) If $\ell \in \{1,2,3\}$ and the expanded stable 
$I$-system $\bold d$ is brimmed$^{\ell +1}$ \ub{then} it is 
brimmed$^\ell$.  If the $I$-system $\bold s$ is brimmed$^{\ell +1}$ and
$\ell = 0,1,2,3$ \ub{then} ${\bold s}$ is brimmed$^\ell$. \nl
5) If $u_0 \subset u_2 \in I$ \ub{then} $\bold J^{\bold d}_{u_2,u_0}$ is
a maximal set $\bold J$ such that
\mr
\item "{$(\alpha)$}"  $\bold J \subseteq M_{u_2}$ is 
disjoint to $\cup\{M_{u_2,w}:w \subset u_2\}$
\sn
\item "{$(\beta)$}"  for each $c \in \bold J$ we have 
{\rm tp}$(c,M_{u_2,u_0},M_{u_2}) 
\in {\Cal S}^{\text{bs}}_{\frak s}(M_{u_2,u_0})$ and
it is orthogonal to $M_{u_2,w}$ whenever $w \subset u_0$
\sn
\item "{$(\gamma)$}"  $\bold J \cup \bigcup 
\{\bold J^0_{u_2,w_1,w_0}:w_1 \subseteq u_2,\neg(w_1=u_2 \wedge
w_0=u_0)$ and $w_1 \nsubseteq u_0,w_0 \subseteq u_0 \cap w_1,w_0 \ne w_1\}$ 
is independent in $(M_{u_2,u_0},M_{u_2})$; note that necessarily there
are no repetitions in the union.
\ermn
6) Assume that in addition to $\bold d$ being an 
expanded stable $I$-system we have
\mr
\item "{$(a)$}"  for $u \subset v \in I,\bold J'_{v,u}$ is a maximal
subset of $\{c \in \bold I_{M^{\bold d}_{v,u},M^{\bold
d}_v}:{\text{\rm tp\/}}_{\frak s}
(c,M^{\bold d}_{v,u},M^{\bold d}_v) \bot M^{\bold d}_w$
for $w \subset u$ and $c \notin M^{\bold d}_{v,w}$ for $w \subset v\}$
such that the set $\bold J'_{v,u} \cup \bigcup\{\bold J^0_{v,w_1,w_0}:w_1
\subseteq v,\neg(w_1 = v \and w_0=u_0),
w_1 \nsubseteq u_0,w_0 \subseteq w_1 \cap u_0,w_0 \ne w_1\}$ is 
independent in $(M^{\bold d}_{v,u},M^{\bold d}_v)$; for the last phrase
alternatively: $\bold J'_{v,u} \cup \{h_{v,w_1}(\bold J'_{w_1,w_0}):
w_1 \subseteq v,\neg(w_1 =v \cap w_0 = u_0),w_1 \nsubseteq u_0,w_0
\subseteq u_0 \cap w_1,w_0 \ne w_1\}$ is independent in 
$(M^{\bold d}_{v,u},M^{\bold d}_v)$
\sn
\item "{$(b)$}"  $\bold d'$ is defined by $\bold s^{{\bold d}'}
= \bold s^{\bold d},\bold J^{{\bold d}'}_{v,u} = \bold J'_{v,u}$ for
$u \subset v \in I$.
\ermn
\ub{Then} ${\bold d}'$ is an expanded stable $I$-system. 
\nl
7) If $\bold d$ is a brimmed$^3$ $I$-system and $u \in I$ \ub{then}
$M^{\bold d}_u \backslash \cup\{M^{\bold d}_{u,w}:w \subset u\}$ has
cardinality $\lambda$; moreover, for $w_1 \subset u,\bold J_{u,w_1}$ is
a subset of $M^{\bold d}_u \backslash \cup\{M^{\bold d}_{u,w_2}:w_2
\subset u\}$ of cardinality $\lambda$.
\endproclaim
\bigskip

\demo{Proof}  0) By \scite{705-12.1R}(1), clause (d) of Definition
\scite{705-12.4}(1) (and the choice of $\bold J^{2,\bold d}_{v,u}$ in
\scite{705-12.4}(1A)). 
 \nl
1) We have
\mr
\widestnumber\item{$(iii)$}
\item "{$(i)$}"  $M_{u,u_1 \cap u_2} \le_{\frak s} M_{u,u_\ell}
\le_{\frak s} M_u$ for $\ell = 1,2$ when $u_1 \cup u_2 \subseteq u \in I$
\nl
[Why?  Each case by a different instance of clause (c) of
\scite{705-12.3}(1)]
\sn
\item "{$(ii)_\ell$}"  $(M_{u,u_1 \cap u_2},M_{u,u_\ell},\bold
J^2_{u,u_\ell,u_0}) \in K^{3,\text{bu}}_{\frak s}$ hence $\in
K^{3,\text{vq}}_{\frak s}$
\nl
[Why?  Without loss of generality $u_1 \ne u_2$ so $u_1 \cap u_2
\subset u_1,u_1 \cap u_2 \subset u_2$.  By part (0), we know 
that $(M_{u_\ell,u_1 \cap u_2},M_{u_\ell},\bold
J^2_{u_\ell,u}) \in K^{3,\text{bu}}_{\frak s}$ and by the definition of
$\bold J^2_{u,u_\ell,u}$ we know that $h_{u,u_\ell}$ maps this triple to
the one mentioned in clause $(ii)_\ell$]
\sn
\item "{$(iii)$}"  $\bold J^2_{u,u_1,u_1 \cap u_2},\bold
J^2_{u,u_2,u_1 \cap u_2}$ are disjoint \nl
[Why?  E.g., as on the one hand
$\bold J^2_{u,u_1,u_1 \cap u_2} \subseteq \bold
J^2_{u,u_2}$ and $\bold J^2_{u,u_2}$ is disjoint to
$M_{u,u_2}$ and on the other hand
$(M_{u,u_2},M_{u,u_1 \cap u_2},\bold J^2_{u,u_1,u_1 \cap u_2}) 
\in K^{3,\text{bu}}_{\frak s}$ hence 
$\bold J^2_{u,u_2,u_1 \cap u_2} \subseteq M_{u,u_2}$]
\sn
\item "{$(iv)$}"  $\bold J^2_{u,u_1,u_1 \cap u_2} \cup
\bold J^2_{u,u_2,u_1 \cap u_2} \subseteq \bold J_{u,u_1 \cap u_2}$
\nl
[Why?  By their definitions]
\sn
\item "{$(v)$}"  $\bold J^2_{u,u_1,u_1 \cap u_2} \cup
\bold J^2_{u_1,u_2,u_1 \cap u_2}$ is independent in $(M_{u,u_1 \cap
u_2},M_u)$ hence in $(M_{u,u_1 \cap u_2},M_{u,u_1 \cup u_2})$ \nl
[Why?  By monotonicity properties.]
\ermn
Now by \scite{705-5.30}(3) we are done. \nl
1A) By Part (1) and \scite{705-5.8}.
\nl
2) As $I$ is finite, there is $w \subset u$ such that $p$ is not
orthogonal to $M_w$ but is orthogonal to $M_{v,w'}$ if $w' \subset
w$.  The ``hence" follows by part (1) and \scite{705-7.6R} (or use part (1A)). \nl
3) By renaming; possible as $M_{v,u_1} \cap M_{v,u_2} = M_{v,u_1 \cap
u_2}$ whenever $u_1,u_2 \subseteq v$ by part (1) and properties of 
NF$_{\bar{\frak s}}$. \nl
4) The least easy part is $\ell=3$.  So we 
have to check clauses $(h)^0_v,(h)^1_v$ of \scite{705-11b.2}(6).  For
the first, clearly every $M^{\bold d}_u (u \in I)$ is brimmed, see clause
$(g)^0_v$ in \scite{705-11b.2}(4) so only the second clause $(h)^1_v$ there
may fail.  \nl
Assume toward contradiction that $u \subset v \in I,p \in 
{\Cal S}^{\text{bs}}_{\frak s}(M^{\bold d}_{v,u})$ and the set 
$\bold J =: \{c \in \bold J^{\bold d}_{v,u}:
\text{tp}(c,M^{\bold d}_{v,u},M^{\bold d}_v) \pm p\}$ has
cardinality $< \lambda_{\frak s}$.  By Definition \scite{705-12.4}(1), $c
\in \bold J^{2,\bold d}_{v,u} \backslash \bold J^{\bold d}_{v,u}
\Rightarrow \text{ tp}(c,M^{\bold d}_{v,u},M_v) \bot p$ hence we have
$\bold J^1 =: \{c \in J^{2,\bold d}_{v,u}:\text{tp}(c,M^{\bold
d}_{v,u},M^{\bold d}_v) \pm p\}$ is equal to $\bold J$ so it still has
cardinality $< \lambda_{\frak s}$.

By the assumption there is $N$
such that $\cup \{M^{\bold d}_{v,w}:w \subset v\} \subseteq N
\le_{\frak s} M^{\bold d}_v$ and $M^{\bold d}_v$ 
is $(\lambda,*)$-brimmed over $N$.  Hence
there is $\bold I$ of cardinality $\lambda_{\frak s}$ independent in
$(M^{\bold d}_{v,u},N,M^{\bold d}_v)$ such that tp$(c,N,M^{\bold
d}_v)$ is a nonforking extension of $p \in {\Cal
S}^{\text{bs}}(M^{\bold d}_{v,u})$ for every $c \in I$.  By dimension
calculus (see Claim \scite{705-10p.16} or \scite{705-5a.3Y}, we 
get contradiction to $(M^{\bold
d}_{v,u},M^{\bold d}_v,\bold J^{2,\bold d}_{v,u}) 
\in K^{3,\text{bu}}_{\frak s}$. \nl
\nl
5),6)  Left to the reader using the dimension calculus (recalling (7)).  \nl
7) The first conclusion follows from the second which holds by parts
(0) + (1) (and clause (d) of Definition \scite{705-12b.3}(1)).
  \hfill$\square_{\scite{705-12b.3}}$\margincite{705-12b.3}
\enddemo
\bigskip

\demo{\stag{705-12b.4} Conclusion}  [Density of explicit expanded 
stable systems, see Definition \scite{705-12.4}(5)]  Assume:
\mr
\item "{$(a)$}"  $\bold s$ is a stable $I$-system 
\sn
\item "{$(b)$}"  for each $u \in I^{\bold s},{\bold P}_u \subseteq \{p
\in {\Cal S}^{\text{bs}}_{\frak s}(M_u):p$ is regular
orthogonal to $M^{\bold d}_{u,w}$ whenever $w \subset u\}$ is 
a maximal subset of pairwise orthogonal regular types.
\ermn
\ub{Then} there is an expanded stable $I$-system $\bold d^*$ such that:
\mr
\item "{$(\alpha)$}"  $\bold s^{{\bold d}^*} = \bold s$
\sn
\item "{$(\beta)$}"  $\bold d^*$ is regular, moreover, $\bold d$
obeys $\bar{\bold P} = \langle {\bold P}_u:
u \in I^{\bold d} \rangle$ which means:
{\roster
\itemitem{ $\boxtimes$ }  if $u \subset v \in I^{\bold d}$ and
$c \in \bold J^{\bold d}_{v,u}$ then for some $q \in {\bold P}_u$ 
we have tp$_{\frak s}(c,M^{\bold d}_{v,u},M^{\bold d}_v)$ is $q$
\endroster}
\sn
\item "{$(\gamma)$}"  $\bold d^*$ is explicit.
\endroster
\enddemo
\bigskip

\demo{Proof}  Note that $(\gamma)$ follows from $(\beta)$.  This holds
by \scite{705-12b.3}(6) and \scite{705-10b.11}(5).  \hfill$\square_{\scite{705-12b.4}}$\margincite{705-12b.4}
\enddemo
\bigskip

\proclaim{\stag{705-12p.9} Claim}  Assume $\ell \in \{1,2,3,4\}$ and
\mr
\item "{$(a)$}"  $\bold d_k$ is an expanded stable
$(\lambda,I)$-system for $k=1,2$
\sn
\item "{$(b)$}"  $\bold s^{\bold d_1} = \bold s^{\bold d_2}$.
\ermn
\ub{Then} $\bold d_1$ is brimmed$^\ell$ iff $\bold d^2$ is brimmed$^\ell$.
\endproclaim
\bigskip

\demo{Proof}  The least easy case is $\ell=3$, see Definition
\scite{705-11b.2}(6), which is similar to \scite{705-12b.4}.
\hfill$\square_{\scite{705-12p.9}}$\margincite{705-12p.9} 
\enddemo
\bigskip

\definition{\stag{705-12f.4A} Definition}  1) For expanded stable
$I$-systems $\bold d_0,\bold d_1$ let $\bold d_0 \le_{\frak s} \bold
d_1$ or $\bold d_0 \le^I_{\frak s} \bold d_1$ mean \footnote{in
Definition \scite{705-11b.2}'s terms this means that
$\langle\text{id}_{M^{\bold d_0}_v}:v \in I \rangle$ embeds $\bold
d_0$ into $\bold d_1$} that:
\mr
\item "{$(a)$}"  $M^{{\bold d}_0}_u \le_{\frak s} M^{{\bold d}_1}_u$
for $u \in I$ and $h^{{\bold d}_0}_{v,u} \subseteq h^{{\bold
d}_1}_{v,u}$ for $u \subset u \in I$
\sn
\item "{$(b)$}"  $\bold J^{{\bold d}_0}_{v,u} \subseteq \bold
J^{{\bold d}_1}_{v,u}$ for $u \subseteq v \in I$
\sn
\item "{$(c)$}"  if $c \in \bold J^{{\bold d}_0}_{v,u}$ then tp$_{\frak
s}(c,M^{{\bold d}_1}_{v,u},M^{{\bold d}_1}_v)$ does not fork over
$M^{{\bold d}_0}_{v,u}$. 
\ermn
2) We say that $J$ is a successor of $I$ \ub{if} for some $t^* \notin
\text{ Dom}(I)$ we have $J= I \cup\{u \cup \{t^*\}:u \in I\}$ and we
call $t^*$ the witness for $J$ being a successor of $I$;
so Dom$(J) = \text{ Dom}(I) \cup \{t^*\}$. \nl
3) For stable $I$-system ${\bold s}_0,{\bold s}_1$ let
\footnote{this relation, $\le_{\frak s}$, is a two-place relation, we
shall prove that it is a partial order.}
${\bold s}_0
\le_{\frak s} {\bold s}_1$ or $\bold s_0 \le^I_{\frak s} \bold s_1$
mean that for some expansions $\bold d_0,\bold d_1$ of $\bold
s_0,\bold s_1$ respectively we have $\bold d_0 \le^I_{\frak s} \bold d_1$. \nl
4) Assume that
$\bold d_0 \le^I_{\frak s} \bold d_1$ and $J$ is a successor of
$I$ with the witness $t^*$ and
\mr
\item "{$(*)$}"  $[u \subset v \in I] \and c \in \bold J^{{\bold
d}_1}_{v,u} \backslash 
\bold J^{{\bold d}_0}_{v,u} \Rightarrow \text{ tp}_{\frak
s}(c,M^{{\bold d}_1}_{v,u},M^{{\bold d}_1}_v)$ either does not fork
over $M^{{\bold d}_0}_{v,u}$ or is orthogonal to $M^{{\bold
d}_0}_{v,u}$.
\ermn
\ub{Then} we let $\bold d \approx \bold d_0 *_J \bold d_1$ mean that
$\bold d = \langle M^{\bold d}_u,h^{\bold d}_{v,u},\bold J^{\bold d}_{v,u}:
u \subseteq v \in I \rangle$ (but note that $\bold d$ is not determined
uniquely by $\bold d_0,\bold d_1,J$ as we have freedom concerning
$\bold J^{\bold d}_{u \cup \{t^*\},u}$ for $u \in I$, still we may use
$\bold d_0 *_J \bold d_1$ for any such $\bold d$) where
\mr
\item "{$(a)$}"  $M^{\bold d}_u \text{ is } 
M^{{\bold d}_0}_u \text{ \ub{if} } u \in I$
\sn
\item "{$(b)$}"  $\bold J^{\bold d}_{v,u} = J^{{\bold d}_0}_{v,u}$ \ub{if}
$u \subseteq v \in I$
\sn
\item "{$(c)$}"  $M^{\bold d}_u$ is 
$M^{{\bold d}_1}_{u \backslash \{t^*\}} \text{ \ub{if} } u \in J \backslash
I$
\sn
\item "{$(d)$}"  $\bold J^{\bold d}_{v,u} = \{c \in \bold J^{{\bold
d}_1}_{v \backslash \{t^*\},u}:c \notin \bold J^{{\bold d}_0}_{v
\backslash \{t^*\},u}$ and tp$_{\frak s}(c,M^{{\bold
d}_1}_{v,u},M^{{\bold d}_1}_v)$ does not fork over $M^{{\bold d}_0}_{v
\backslash \{t^*\},u}\}$ \ub{if} $u \in I,v \in J \backslash I,u
\subset v$
\sn 
\item "{$(e)$}"  $\bold J^{\bold d}_{v,u} = \{c \in \bold J^{{\bold d}_1}
_{v \backslash \{t^*\},u \backslash \{t^*\}}:\text{tp}_{\frak s}
(c,M^{{\bold d}_1}_{v \backslash \{t^*\},u \backslash \{t^*\}},
M^{{\bold d}_1}_v)$ is orthogonal to 
$M^{{\bold d}_0}_{v \backslash \{t^*\},u \backslash \{t^*\}}\}$
\ub{if} $u \subset v$ are both from $J \backslash I$
\sn
\item "{$(f)$}"  $h^{\bold d}_{v,u}$ is: $h^{{\bold d}_0}_{v,u}$ if
$u \subseteq v \in I;h^{\bold d}_{v,v}$ is 
$h^{{\bold d}_1}_{v,u}$ if $t \in u \subseteq
v \in I$ and $h^{\bold d}_{v,v}$
is $h^{{\bold d}_0}_{v_1,u}$ if $v = v_1 \cup \{t(*)\},u
\subseteq v_1 \in I$
\sn
\item "{$(g)$}"  $\bold J^{\bold d}_{u \cup \{t(*)\},u}$ is a maximal
subset of $\{c \in M^{{\bold d}_1}_u:\text{tp}_{\frak s}(c,
M^{{\bold d}_0}_u,M^{{\bold d}_1}_u) \in {\Cal S}^{\text{bs}}_{\frak s}
(M^{{\bold d}_0}_u)$ is orthogonal to $M^{{\bold d}_0}_w$ for every $w
\subset u\}$ which is independent in $(M^{{\bold d}_0}_u,M^{{\bold
d}_1}_u)$ for every $u \in I$.
\ermn
4A) Similarly for $\bold s = \bold s_0 *_J \bold s_1$ (but now $\bold
s$ is uniquely determined). \nl
5) If $\delta \le \lambda^+_{\frak s}$ and  
$\langle \bold d_\alpha:\alpha < \delta \rangle$ is a
$\le^I_{\frak s}$-increasing sequence of expanded stable $I$-systems
(see \scite{705-12f.4B} below) \ub{then} we
let $\bold d = \dbcu_\alpha \bold d_\alpha$
be $\langle M^{\bold d}_u,h^{\bold d}_{v,u},
\bold J^{\bold d}_{v,u}:u \subseteq v \in I
\rangle$ where $M^{\bold d}_u = \cup\{M^{{\bold d}_\alpha}_u:\alpha <
\delta\}$ and $h^{\bold d}_{v,u} = \cup\{M^{{\bold
d}_\alpha}_{v,u}:\alpha < \delta\}$ and
$\bold J^{\bold d}_{v,u} = \cup\{\bold J^{{\bold d}_\alpha}_{v,u}:\alpha <
\delta\}$.  Similarly for $\langle \bold s_\alpha:\alpha < \delta
\rangle$.
\enddefinition
\bigskip

\proclaim{\stag{705-12f.4B} Claim}  1) $\le^I_{\frak s}$ is a partial order
on the family of expanded stable $I$-systems. \nl
2) If $\delta < \lambda^+_{\frak s}$ and $\langle \bold d_\alpha:\alpha <
\delta \rangle$ is $\le^I_{\frak s}$-increasing sequence (of expanded
stable $I$-systems) \ub{then} $\bold d = \dbcu_{\alpha < \delta} \bold
d_\alpha$ is an expanded stable $I$-system and $\alpha < \delta
\Rightarrow \bold d_\alpha \le^I_{\frak s} \bold d$. \nl
3) If $\bold d_0 \le^I_{\frak s} \bold d_1$ and $u \subset v \in I$
\ub{then} {\rm NF}$_{\frak s}(M^{{\bold d}_0}_{v,u},M^{{\bold
d}_0}_v,M^{{\bold d}_1}_{v,u},M^{{\bold d}_1}_v)$. \nl
4) If $\bold d_0 \le^I_{\frak s} \bold d_1$ and $\bold d'_0$ is an
expanded stable $I$-system satisfying $\bold s^{{\bold d}'_0} =
\bold s^{{\bold d}_0}$ \ub{then} we can find an expanded stable
$I$-system $\bold d'_1$ such that $\bold d'_0 \le^I_{\frak s} \bold
d'_1$ and $\bold s^{{\bold d}'_1} = \bold s^{{\bold d}_1}$. \nl
5) The relation $\le^I_{\frak s}$ on the family of stable $I$-systems
is a partial order.
\nl
6) In Definition \scite{705-12f.4A}(4), always there is $\bold d$ such that
$\bold d \approx  \bold d_0 *_J \bold d_1$
is an expanded stable $(\lambda,J)$-system; similarly without
``expanded". \nl
7) For $I$-system $\bold s_1,\bold s_2$ we have
$\bold s_1 \le_{\frak s} \bold s_2$ \ub{iff} $\bold s_1 \le \bold
s_2$ and $u \subset v \in I \Rightarrow { \text{\rm NF\/}}_{\frak s}(M^{\bold
s_1}_u,M^{\bold s_1}_v,M^{\bold s_2}_u,M^{\bold s_2}_v)$. 
\endproclaim
\bigskip

\demo{Proof}  1) Obvious.  Check the definition. \nl
2) First we prove this for stable expanded systems.
The main point is why, for $u \subset v \in I$, the triple
$(M^{\bold d}_{v,u},M^{\bold d}_v,\bold J^{2,\bold d}_{v,u})$ belong to
$K^{3,\text{bu}}_{\frak s}$.  This holds by \scite{705-12.1R}. 

Second, for stable systems, given $\langle {\bold s}_\alpha:\alpha <
\delta \rangle,{\frak s}_\delta =: \bold s$ we choose by induction on
$\alpha \le \delta$ an expanded stable $I$-system $\bold
d_\alpha$ such that $\bold s^{{\bold d}_\alpha} = {\bold s}_\alpha$
and $\langle \bold d_\beta:\beta \le \alpha \rangle$ is increasing
continuous.  For $\alpha = 0$ this is by the definition, for $\alpha =
\beta +1$ by part (4) below (which does not rely on parts (2), (3).
Lastly, for $\alpha$ limit use what we have proved for stable expanded
systems.  For $\alpha = \delta,\bold d_\delta$ witnessed which we
need.
\nl
3) As in the proof of \scite{705-12b.3}(1). \nl
4) As in the proof of \scite{705-12b.3}(6). \nl
5) Being partial order follows by (4) and (1).  The union of
increasing continuous sequences follows by (4) and (1).  [Note that
the cases of non-continuous chains are problematic, those are 
parallel of AxIII$_2$ of a.e.c.; compares with \scite{705-12f.4D}.]  
As for brimmed$^\ell$ for $\bold d$'s version: if $\ell=1$ this is
trivial; if $\ell=2$ this is as ``the union of an increasing chain of
brimmed model" (for ${\frak s}$) of length $< \lambda^+$ is brimmed.
Lastly, for $\ell=3$, by the finite character of
nonforking (= Ax(E)(x)) and the definition.
\nl
6) Easy by \S10. \nl
7) Think (that is let $\bold d_\ell$ be an expanded $I$-system with
$\bold s[\bold d_\ell] = \bold s_\ell$ for $\ell=1,2$.  Now we shall
choose $\bold J'_{v,u}$ for $u \subset v \in \bold I$ such that $\bold
d_1 \le^I_{\frak s} \bold d'_2$ where $\bold s[\bold d'_v] = \bold
s_2,J^{\bold d'_2}_{v,u} = \bold J'_{v,u}$.  We do this by induction
on $(v)$, as in previous cases).  \hfill$\square_{\scite{705-12f.4B}}$\margincite{705-12f.4B}
\enddemo
\bigskip

\proclaim{\stag{705-12f.4C} Claim}  Assume
\mr
\item "{$(a)$}"  $\bold s_1$ is a stable $I$-system
\sn
\item "{$(b)$}"  $\bold s_0$ is an $I$-system
\sn
\item "{$(c)$}"  $M^{{\bold s}_0}_u 
\le_{\frak s} M^{{\bold s}_1}_u$ for $u \in I$ and 
$h^{{\bold s}_0}_{v,u} \subseteq h^{{\bold s}_1}_{v,u}$ for $u \subset
v \in I$
\sn
\item "{$(d)$}"  if $u \subset v \in I$ then {\rm NF}$_{\frak s}(M^{{\bold
s}_0}_{v,u},M^{{\bold s}_0}_v,M^{{\bold s}_1}_{v,u},M^{{\bold s}_1}_v)$.
\ermn
\ub{Then} $\bold s_0$ is a stable $I$-system and $\bold s_0
\le^I_{\frak s} \bold s_1$.
\endproclaim
\bigskip

\demo{Proof}  Let $\bold d_1$ be an expanded stable $I$-system such
that $\bold s^{{\bold d}_1} = \bold s_1$.
For each $u \subset v \in I$ we choose $\bold I_{v,u}$
as a maximal set such that
\mr
\widestnumber\item{$(iii)$}
\item "{$(i)$}"  $\bold I_{v,u} \subseteq M^{{\bold s}_0}_v \backslash
\cup \{M^{{\bold s}_0}_{v,w}:w \subset v\}$ and for any $c \in \bold
I_{v,u}$ we have tp$_{\frak s}(c,M^{\bold s_0}_{v,u},M^{\bold s_0}_v)
\in {\Cal S}^{\text{bs}}(M^{\bold s_0}_v)$ is orthogonal to $M^{\bold
s_0}_{v,w}$ whenever $w \subset u$
\sn
\item "{$(ii)$}"  $\bold I_{v,u} \cup \{\bold J^{2,{\bold
s}_1}_{v,w_1,w_0}:w_0 \subseteq w_1 \subset v,w_0 \subseteq u\}$ is
independent in $(M^{{\bold s}_1}_{v,u},M^{{\bold s}_1}_v)$.
\ermn
Let $\bold J'_{v,u}$ be maximal set such that 
\mr
\widestnumber\item{$(iii)$}
\item "{$(i)$}"  $\bold J'_{v,u} \subseteq M^{{\bold s}_1}_v
\cup\{M^{{\bold s}_1}_{v,w}:w \subset v\}$
\sn
\item "{$(ii)$}"  $\bold J'_{v,u} 
\cup \{\bold J^{2,{\bold s}_1}_{v,w_1,w_0}:
w_0 \subseteq w_1 \subset v,w_2 \subseteq u\}$ is
independent in $(M^{{\bold s}_1}_{v,u},M^{{\bold s}_1}_v)$
\sn
\item "{$(iii)$}"  $\bold I_{v,u} \subseteq \bold J'_{v,u}$.
\ermn 
So by \scite{705-12b.3}(6), $\bold d'_1 =: \langle M^{{\bold s}_1}_u,h^{{\bold
s}_1}_{v,u},\bold J'_{v,u}:u \subseteq v \in I \rangle$ is an expanded
stable $I$-system.  Also $\bold d_0 =: \langle M^{{\bold
s}_0}_u,h^{{\bold s}_0}_{v,u};\bold I_{v,u}:u \subseteq v \in I
\rangle$ is an expanded stable $I$-system.  Now easily $\bold s_\ell =
\bold s^{{\bold d}_\ell}$ and $\bold d_0 \le^I_{\frak s} \bold d'_1$.
\nl
So we are done.  \hfill$\square_{\scite{705-12f.4C}}$\margincite{705-12f.4C} 
\enddemo
\bigskip

\proclaim{\stag{705-12f.4D} Claim}  Assume $J$ is a successor of $I$ with
witness $t^*$. \nl
1) Assume $\bold s_0,\bold s_1$ are stable $I$-systems
satisfying $\bold s_0 \le^I_{\frak s} \bold s_1$.  \ub{Then}
\mr
\item "{$(a)$}"  for one and only one 
$\bold s,\bold s = \bold s_0 *_J \bold s_1$
\sn
\item "{$(b)$}"   $\bold s$ is a stable $J$-system. 
\ermn
2) Assume that $\delta < \lambda^+_{\frak s}$ 
and $\langle \bold s_\alpha:\alpha < \delta \rangle$ is a
$<^I_{\frak s}$-increasing continuous sequence of stable $J$-systems
and $\bold s_\alpha$ is brimmed$^\ell$ for $\alpha <
\delta$ and $\ell \le 3$ \ub{then} $\bold s_\delta$ is brimmed$^\ell$ and 
$\bold s_\alpha \le^I_{\frak s} \bold s_\delta$ for $\alpha < \delta$. \nl
3) Assume that ${\frak s}$ is successful hence
${\frak s}^+$ is as in \scite{705-12.1} and $\ell = 3$.  
If $\langle \bold s_\alpha:\alpha < \lambda^+_{\frak s}
\rangle$ is a $<^I_{\frak s}$-increasing sequence of stable
$I$-systems and $\bold s_\alpha *_J \bold s_{\alpha +1}$ is
brimmed$^\ell$ for each 
$\alpha < \lambda^+$ \ub{then} $\bold s = \cup\{\bold
s_\alpha:\alpha < \lambda^+_{\frak s}\}$ is a brimmed$^\ell \,
(\lambda^+,I,{\frak s})$-system. \nl
4) In part (3), if $\delta = \lambda^+_{\frak s}$ and 
$\ell \in \{0,1,2\}$ and $a \in I \Rightarrow
M^{\bold s}_u \in K_{{\frak s}(+)}$ \ub{then} $\bold s$ is brimmed$^\ell
\, (I,{\frak s}^+)$-system. 
\endproclaim
\bigskip

\demo{Proof}  1) Clearly $\bold s$ is well defined. \nl

Recall that $\bold d_0 \le^I_{\frak s} \bold d_1$ above does not
imply that for some $\bold d,\bold d = \bold d_0 * \bold d_1$. \nl
By Definition \scite{705-12f.4A}(3) there are stable
expanded $I$-system $\bold d_1,\bold d_2$ such that $\bold s^{{\bold
d}_\ell} = \bold s_\ell$ and $\bold d_0 \le_{\frak s} \bold d_1$.  For
each $u \in I$ let ${\bold P}^0_u$ be a maximal set of pairwise
orthogonal regular types from ${\Cal S}^{\text{bs}}_{\frak
s}(M^{{\bold s}_0}_u)$ orthogonal to $M^{{\bold s}_0}_{u,w}$ for every
$w \subset u$.  For each $u \in I$ let ${\bold P}^1_u$ be a maximal set
of pairwise orthogonal regular types $p \in 
{\bold S}^{\text{bs}}_{\frak s}(M^{{\bold s}_1}_u)$ orthogonal to $M^{{\bold
s}_v}_w$ for every $w \subset u$ such that either $p$ is a nonforking
extension of some $q \in {\bold P}^0_u$ or $p \perp M^{{\bold s}_0}_u$.  
By Claims \scite{705-12b.3}(6), \scite{705-12b.4} (and see \scite{705-12f.4B}(4)),
\wilog \, $\{0,1\},w \subset u \subset v \in I,c \in \bold J^{\bold d_i}_{v,u}
\Rightarrow \text{ tp}_{\frak s}(c,M^{{\bold d_i}_\ell}_{v,u},
M^{{\bold d}_i}_{v,u}) \in \bold P^i_u$ and $\bold J^{\bold d_0}_{v,u}
\subseteq \bold J^{\bold d_1}_{v,u}$.

So $\bold d_0,\bold d_1$ are as in Definition \scite{705-12f.4A}(4) (as
$J,t(*)$ are given) hence there is a stable expanded $I$-system $\bold d
\approx \bold d_0 *_J \bold d_1$.  So $\bold s = \bold s^{\bold d}$ is
as required. \nl
2) For $\ell=0$ this is just by properties of a.e.c..  For $\ell=1$ we
use Claim \scite{705-12f.4B}(1) (we can expand $\bold s_\alpha$ to 
$\bold d_\alpha,\bold d_\alpha$ increasing continuous by Claim
\scite{705-12f.4B}(4) and \scite{705-12f.4B}(5).  
For $\ell=2$ note that the union of an 
$\le_{\frak s}$-increasing continuous sequence of brimmed models in brimmed.

For $\ell=3$ use the basic properties of orthogonality. \nl
3) Easy, too. \nl
4) Easy, too.  \hfill$\square_{\scite{705-12f.4D}}$\margincite{705-12f.4D}
\enddemo
\bigskip

\proclaim{\stag{705-12b.10} Lemma}   1) Let $\ell \in \{0,2,3\}$.
Assume that ${\frak s}$ is successful (hence ${\frak s}^+$ 
has the properties required in \scite{705-12.1}) and
\mr
\item "{$(a)$}"  $J$ a successor of $I$, in details
$I \subseteq {\Cal P}(u_*)$ is downward closed,
$u_*$ finite, $u_{**} = u_* \cup 
\{t^*\},t^* \notin u_*,J = I \cup \{u \cup \{t^*\}:u \in I\}$  
\sn
\item "{$(b)$}"  $\bold s$ is a brimmed$^\ell$ 
$(\lambda^+,I,{\frak s}^+)$-system
\sn 
\item "{$(c)$}"  $\langle M^{\bold s}_{u,\alpha}:\alpha < \lambda^+
\rangle$ is a $<_{{\frak K}[{\frak s}]}$-representation of $M^{\bold s}_u$
\sn
\item "{$(d)$}"  for $\alpha < \lambda^+$ we try to define
$(\lambda,I,{\frak s})$-system $\bold s_\alpha$ by $M^{{\bold
s}_\alpha}_u = M^{\bold s}_{u,\alpha},h^{{\bold s}_\alpha}_{v,u} =
h^{\bold s}_{v,u} \restriction M^{\bold s}_{u,\alpha}$
\sn
\item "{$(e)$}"  for $\alpha < \beta < \lambda^+$ we try to define a
$(\lambda,J,{\frak s})$-system ${\bold s}_{\alpha,\beta}$ by

$$
M^{{\bold s}_{\alpha,\beta}}_u = M^{\bold s}_{u,\alpha} \text{ for } u
\in I
$$

$$
M^{{\bold s}_{\alpha,\beta}}_{u \cup \{t^*\}} = M^{\bold s}_\beta
$$

$$
h^{{\bold s}_{\alpha,\beta}}_{v,u} = h^{\bold s}_{v,u} \restriction
M^{{\bold s}_{\alpha,\beta}}_{u \cup \{t^*\}} \text{ if } u \in I,
u \subseteq v \in J
$$

$$
h^{{\bold s}_{\alpha,\beta}}_{v,u \cup \{t^*\}} = h^{\bold s}_{v,u}
\restriction M^{{\bold s}_{\alpha,\beta}}_u \text{ if } u \in I,u \cup
\{t^*\} \subseteq v \in I.
$$
\ermn
\ub{Then} for some club $E$ of $\lambda^+$ we have
\mr
\item "{$(\alpha)$}"  for every $\alpha \in E,\bold s^\alpha$ is an
$(I,{\frak s})$-system which is brimmed$^\ell$
\sn
\item "{$(\beta)$}"   for every $\alpha < \beta$
from $E,\bold s_{\alpha,\beta}$ is a brimmed$^\ell$ 
$(J,{\frak s})$-system. 
\ermn
2) Assume that $\ell \in \{2,3\}$ and (a),(c),(d),(e) of part (1)
holds for $\bold s = \bold s^{\bold d}$ and
\mr
\item "{$(b)'$}"  $\bold d$ is an explicit brimmed$^\ell$ expanded
stable $(\lambda^+,I,{\frak s}^+)$-system.
\ermn
\ub{Then} for some club $E$ of $\lambda^+$,  we can define for $\alpha
\in E,\bold d^\alpha$ as below and we can find 
$\bold d^{\alpha,\beta}$ for $\alpha < \beta$ from $E$ as below such that
\mr
\item "{$(\alpha)$}"  $\bold d_\alpha$ is an explicit 
expanded stable $(\lambda,I,{\frak s})$-system
\sn
\item "{$(\beta)$}"  $\bold d^{\alpha,\beta}$ is an explicit
brimmed$^\ell$ expanded stable $(\lambda^+,J,{\frak s}^+)$-system
\sn
\item "{$(\gamma)$}"  $\bold s[\bold d^{\alpha,\beta}] = \bold
s^{\alpha,\beta}$ from part (1) and $\bold s[\bold d^\alpha] = \bold
s^\alpha$ from part (1) 
\sn
\item "{$(\delta)$}"  for $u \subset v \in I,\bold J^{{\bold
d}^\alpha}_{v,u} = \bold J^{\bold d}_{v,u} \cap M^{\bold d^\alpha}_v$ 
\endroster
\endproclaim
\bigskip

\demo{Proof}  1) We leave the case $\ell=0$ to the reader, so it is
enough to prove part (2). \nl
2) Choose for each $u \in I$ a maximal subset of ${\bold P}_u$ of
${\bold P}^*_u = \{p \in {\Cal S}^{\text{bs}}_{{\frak s}(+)}
(M^{\bold s}_u):p$ regular orthogonal to 
$M^{\bold s}_{u,w}$ for every $w \subset u\}$ of pairwise
orthogonal types.

Let $\bold d$ be the stable expanded $(\lambda,I)$-system and let
$\bold s^{\bold d} = \bold s$.  Now
\mr
\item "{$\circledast_1$}"  if $p \in {\bold P}^*_u$ then for some $\alpha_0 =
\alpha_0(p) \le \alpha_1 = \alpha_1(p) < \lambda^+_{\frak s}$ we have:
{\roster
\itemitem{ $(a)$ }  $M^{\bold s}_{v,\alpha_1}$ is a witness for $p$
\sn
\itemitem{ $(b)$ }  $p \restriction M^{\bold s}_{u,\alpha_1}$ is
orthogonal to $M^{\bold s}_{u,\gamma}$ iff $\gamma < \alpha_0$.
\endroster}
\ermn
For each $p \in {\bold P}_u$ let $\bold J_p$ be a maximal subset of
$\bold I_{M^{\bold s}_{u,\alpha_1(p)},M^{\bold s}_u} = \cup\{\bold
I_{M^{\bold s}_{u,\alpha(p)},M^{\bold s}_{u,\beta}}:\beta \in
[\alpha_1(p),\lambda^+_{\frak s})\}$ of elements
realizing $p \restriction M^{\bold s}_{u,\alpha_1(p)}$ in $M^{\bold
s}_u$ which is independent in $(M^{\bold s}_{u,\alpha_1(p)},M^{\bold
s}_u)$. 

We can find a club $E$ of $\lambda^+_{\frak s}$ such that
\mr
\item "{$\circledast_2$}"  if $\delta \in E$ then
{\roster
\itemitem{ $(a)$ }  for $u \subset v \in I$ we have 
$h^{\bold d}_{v,u}(M^{\bold s}_{u,\delta}) 
= M^{\bold s}_{v,u} \cap M^{\bold s}_{v,\delta}$
\sn
\itemitem{ $(b)$ }  for $u \subset v$ and $c \in 
\bold J^{\bold d}_{v,u}$, if $c \in M^{\bold s}_{v,\delta}$ then
$\delta \ge \alpha_1({\text{\rm tp\/}}(c,M^{\bold s}_{v,u},M^{\bold s}_v))$
\sn
\itemitem{ $(c)$ }  $(M^{\bold s}_{v,u} \cap M^{\bold
s}_{v,\delta},M^{\bold s}_{v,\delta},
\cup\{\bold J^{2,\bold d}_{v,u_1,u_0} 
\cap M^{\bold s}_{v,\delta}:u_0 \subseteq u,u_1
\subseteq v,u_1 \nsubseteq u\})$ belongs to $K^{3,\text{bu}}_{\frak s}$
\sn
\itemitem{ $(d)$ }  if $u \in I,p \in {\bold P}_u,\alpha_1(p) < \delta$
then $\bold J_p \cap M^{\bold s}_{u,\delta}$ is a maximal subset of
$\bold I_{M^{\bold s}_{u,\alpha_1(p)},M^{\bold s}_{u,\delta}}$
of elements realizing $p \restriction M^{\bold s}_{u,\alpha(p)}$ which
is independent in $(M^{\bold s}_{u,\alpha_1(p)},M^{\bold
s}_{u,\delta})$
\sn
\itemitem{ $(e)$ }  if $u \subset v \in I$ and 
$\delta_1 < \delta_2$ are
from $E$ then NF$_{\frak s}(M^{\bold s}_{u,\delta_1},M^{\bold
s}_{v,\delta_1},M^{\bold s}_{u,\delta_2},M^{\bold s}_{v,\delta_2})$
\sn
\itemitem{ $(f)$ }  if $u \subset v \in I,c \in \bold J^{\bold s}_{v,u}$ and
$p = \text{ tp}_{{\frak s}(+)}(c,M^{\bold s}_{v,u},M^{\bold s}_v)$ and
$\alpha_0(p) \le \delta$ \ub{then} $\alpha_1(p) \le \delta$.
\endroster}
\ermn
Let $\alpha < \beta$ be from $E$ and let $\bold s^{\alpha,\beta}$ be
as in part (1) of the claim.  
We define an expanded $(\lambda,J)$-system $\bold d^{\alpha,\beta}$ by
(recalling $\bold d$ is explicit)
\mr
\item "{$(a)$}"  $\bold s^{\bold d^{\alpha,\beta}} = \bold
s^{\alpha,\beta}$ (defined in part (1))
\sn
\item "{$(b)$}"  if $u \subset v \in I$ then 
$\bold J^{\bold d^{\alpha,\beta}}_{u,v} = \bold J^{\bold d}_{u,v} \cap
M^{\bold s}_{v,\alpha}$
\sn
\item "{$(c)$}"  if $u \subset v \in J,u \in I,v = v_1 \cup \{t^*\}$
hence $u \subset v_1 \in I$ then \nl
$\bold J^{\bold d^{\alpha,\beta}}_{u,v} = \{c \in \bold J^{\bold d}_{v,u}:c \in
M^{\bold s^{\alpha,\beta}}_v$ and $\alpha \ge \alpha_1({\text{\rm
tp\/}}(c,M^{\bold d}_{v,u},M^{\bold d}_v))\}$
\sn
\item "{$(d)$}"  if $u \subseteq v \in J,u = u_1 \cup \{t^*\},v = v_1
\cup \{t^*\},u \subset v \in I$ then \nl
$\bold J^{\bold d^{\alpha,\beta}}_{v,u} = \{c \in 
\bold J^{\bold d}_{v,u}:c \in \bold J^{\bold d}_{v,u}$ 
and $\alpha < \alpha_1({\text{\rm tp\/}}
(c,M^{\bold d}_{v,u},M^{\bold d}_v))\}$
\sn
\item "{$(e)$}"  if $u \in I,v = u \cup \{t^*\}$ then
\nl
$\bold J^{\bold d^{\alpha,\beta}}_{v,u} = \{c$: for some $p \in
{\bold P}_u,\alpha_1(p) \le \alpha$ and $c \in 
\bold J_p \cap M^{\bold s}_{u,\beta} \backslash M^{\bold s}_{u,\alpha}\}$.
\ermn
Now check.  \hfill$\square_{\scite{705-12b.10}}$\margincite{705-12b.10}
\enddemo
\bigskip

\definition{\stag{705-12b.5} Definition}  Let $\ell \in \{1,2,3,4\}$.  
\nl
1) We say that ${\frak s}$ satisfies
(or has) the brimmed$^\ell$ weak $(\lambda,n)$-existence property \ub{if}: 
\mn
\ub{Case 1}:  $n=0$.

There is a brimmed$^\ell$ model in $K_{\frak s}$ (so always holds).
\mn
\ub{Case 2}:  \ub{$n = 1$}.

If $M_\emptyset$ is brimmed$^\ell,p_i \in {\Cal
S}^{\text{bs}}(M_\emptyset)$ for $i < i^* < \lambda^+_{\frak s}$ then we can
find $c_i (i < i^*)$ and $M_{\{\emptyset\}}$ such that $p_i = \text{
tp}(c_i,M_\emptyset,M_{\{\emptyset\}})$ and 
$(M_\emptyset,M_{\{\emptyset\}},\{c_i:i
< i^*\}) \in K^{3,\text{bs}}_{\frak s},i < j \Rightarrow c_i \ne c_j$.
\mn
\ub{Case 3}:  $n \ge 2$.
\mn
Every brimmed$^\ell$ expanded stable 
${\Cal P}^-(n)$-system ${\bold d}$ can be completed to an
expanded stable $({\frak s},{\Cal P}(n))$-system $\bold d^+$,
i.e. there is an expanded stable 
${\Cal P}(n)$-system $\bold d^+$ such that ${\bold d}^+ 
\restriction {\Cal P}^-(n) = {\bold d}$; recall 
${\Cal P}^-(n) = \{u:u \in \{0,\dotsc,n-1\}\} = 
{\Cal P}(n) \backslash \{n\}$. 
\nl
2) We say ${\frak s}$ satisfies (or has) brimmed$^\ell$ weak 
$(\lambda,n)$-uniqueness property when:
\mn
\ub{Case 1}:  If $n=0,K_{\frak s}$ has amalgamation.
\mn
\ub{Case 2}:  If $n=1$.  If $(M,N_\ell,\{a^k_i:i < i^*\}) \in
K^{3,\text{bu}}_{\frak s}$ for $k=1,2$ and tp$(a^1_i,M,N_\ell) = 
\text{ tp}(a^2_i,M,N_\ell)$ there is a $\le_{\frak s}$-embedding $f$ of $N_1$
into some $N'_2$ such that $N_2 \le_{\frak s} N'_2$ and $f(a^1_1) = a^2_i$.
\mn
\ub{Case 3}:  $\bold d_1,\bold d_2$ are brimmed$^\ell$ expanded stable 
${\Cal P}(n)$-systems and $\bold d^m = \bold d_m \restriction {\Cal
P}^-(n)$ and $\bar f = \langle f_u:u \in {\Cal P}^-(n) \rangle$ 
is an isomorphism from $\bold s[\bold d^1]$ onto
$\bold s[\bold d^2]$, \ub{then} we can find $(f,N)$ such that:
\mr
\item "{$(a)$}"  $M^{{\bold d}_2}_n \le_{\frak s} N$,
\sn
\item "{$(b)$}"  $f$ is a $\le_{\frak s}$-embedding of $M^{{\bold
d}_1}_n$ into $N$
\sn
\item "{$(c)$}"  for $u \subset n$ we have $f \circ h^{{\bold
d}_1}_{n,u} = h^{{\bold d}_2}_{n,u} \circ f_u$.
\ermn
3) We say that ${\frak s}$ has the brimmed$^\ell$ strong
$(\lambda,n)$-existence \ub{when}:
\mn
\ub{Case 1}:  $n=0$, just $K_{\frak s} \ne \emptyset$.
\mn
\ub{Case 2}:  $n=1$, as in part (1) but
$(M_\emptyset,M_{\{\emptyset\}},\{c_i:i < i^*\}) \in
K^{3,\text{bu}}_{\frak s}$. 
\mn
\ub{Case 3}:  $n \ge 2$.  

For every brimmed$^\ell$ expanded stable 
$(\lambda,{\Cal P}^-(n))$-system $\bold
d$ we can find an expanded stable $(\lambda,{\Cal P}(n))$-system
$\bold d^+$ such that $\bold d^+ \restriction {\Cal P}^-(n) = \bold d$
and $\bold d^+$ is reduced in $n$ which means $u \subset
n \Rightarrow \bold J^{{\bold d}^+}_{n,u} = \emptyset$. 
\nl
4) We say that ${\frak s}$ has the brimmed$^\ell$ 
strong $(\lambda,n)$-uniqueness property when:
\mn
\ub{Case 1}:  If $n=0,K_{\frak s}$ categorical in $\lambda_{\frak s}$.
\mn
\ub{Case 2}:  If $n=1$, uniqueness for $K^{3,\text{bu}}_{\frak s}$ see
\scite{705-5.25}. 
\mn
\ub{Case 3}:  $n \ge 2$ and $\ell \in \{3,4\}$.

In part (2) we add  $N = M^{{\bold d}_2}_n$ and $f$ is onto $N$.
\mn
\ub{Case 4}: $n \ge 2$ and $\ell \in \{1,2\}$.  

The conclusion of Case 1 holds \footnote{in the csae central for this
section the $M^{{\bold d}_m}_n$ are brimmed so the only freedom left
are about dimensions of types from ${\Cal S}^{\text{bs}}_{\frak
s}(M^{{\bold d}_m}_{n,u})$} if we assume (what is said in part (2)
and) that $\bold d_1,\bold d_2$ are regular (or just $c \in \bold
J_{n,u} \Rightarrow \text{ tp}_{\frak s}(c,h_{n,u}(M^{\bold
d_k}_u),M^{\bold d_k}_n)$ is regular):
\mr
\item "{$\boxdot$}"  if $u \subset v \in {\Cal P}^-(n)$ and 
$p \in {\Cal S}^{\text{bs}}(M^{\bold d_1}_{n,u})$ then the 
cardinality of $\{c \in \bold J^{{\bold d}_1}_{n,u}:
\text{tp}_{\frak s}(c,M^{{\bold
d}_1}_{n,u},M^{{\bold d}_1}_{u_n}) \pm p\}$ is equal to the
cardinality of $\{c \in \bold J^{{\bold d}_2}_{n,u}:
\text{tp}_{\frak s}(c,M^{{\bold d}_2}_{n,u},M^{{\bold d}_2}_n) \pm
f(p)\}$ assuming for simplicity $\bold d_1$ is regular.
\ermn
5)  We say that ${\frak s}$ has the brimmed$^\ell$ weak/strong 
$(\lambda,n)$-primeness property when $n=0,1$ or for any 
brimmed$^\ell \,{\Cal P}^-(n)$-system $\bold d_0$ and expanded stable
${\Cal P}(n)$-system $\bold d_1$ which is reduced in $n$ which means $u
= n \Rightarrow \bold J^{{\bold d}_1}_{n,u} = \emptyset$ and $\bold
d_1$ satisfy $\bold d_1 \restriction {\Cal P}^-(n) = \bold d_0$ we
have: $\bold d_1$ is weakly/strongly prime$^\ell$ over $\bold d_0$,
see below. \nl
5A)  We say $\bold d_1$ is weakly/strongly prime$^\ell$ over 
$\bold d_0$ (and also say that $\bold s^{{\bold d}_1}$ is
weakly/strongly
prime$^\ell$ over $\bold s^{{\bold d}_0}$) when for some $n,\bold d_1$ is a
brimmed$^\ell \, {\Cal P}(n)$-system, $\bold d_0 = \bold d_1
\restriction {\Cal P}^-(n)$ and:
\mr
\item "{$(*)$}"  if $\bold d_2$ is a brimmed$^\ell \,
({\Cal P}(n),{\frak s})$-system satisfying $\bold
d_2 \restriction {\Cal P}^-(n) = \bold d_0 \restriction {\Cal
P}^-(n)$, \ub{then} there is an $\le_{\frak s}$-embedding 
$f$ of $M^{{\bold d}_1}_n$ into
$M^{{\bold d}_2}_n$ such that $u \in {\Cal P}^-(n) \Rightarrow f \circ
h^{{\bold d}_1}_{n,u} = h^{{\bold d}_2}_{n,u}$, and in the strong
case, for every regular $p \in {\Cal S}^{\text{bs}}(f(M^{{\bold
d}_1}_n)$ orthogonal to $M^{{\bold d}_2}_{n,u}$ for every $u \subset n$
we have dim$(p,M^{{\bold d}_2}_n) = \lambda$.
\ermn
6) In (1)-(5) we may restrict ourselves to one expanded stable
${\Cal P}^-(n)$-system or ${\Cal P}(n)$-system $\bold d$, i.e.,
consider the property as a property of $\bold d$; so in this case
the brimmed$^\ell$ may refer to only $u = n$! \nl
7) We say that ${\frak s}$ has the brimmed$^\ell$ weak/strong
$(\lambda,n)$-prime existence \ub{if} for every brimmed$^\ell$ expanded
stable $(\lambda,{\Cal P}^-(n))$-system $\bold d_1$ there is an
expanded stable $(\lambda,{\Cal P}(n))$-system $\bold d_2$ which is
weakly/strongly prime$^\ell$ 
over $\bold d_1$, (note: $\bold d_2$ is only brimmed$^1$ and is
reduced in $n$). \nl
8) Writing ``$\ldots (<n) \ldots$ property" we mean ``$\ldots m
\ldots$ property" for every $m<n$.
\enddefinition
\bn
Our main aim is to show that if ${\frak s}$ is excellent and $\langle
2^{\lambda^{+n}_{\frak s}}:n < \omega \rangle$ increasing then every
one of those properties is satisfied by ${\frak s}^{+m}$ for $m <
\omega$ large enough. 
\proclaim{\stag{705-12p.15B} Claim}  1) For $n=0,1,2$, the frame ${\frak s}$ has the
brimmed$^\ell$ weak $(\lambda,n)$-existence property for $\ell \le 4$. \nl
2) For $n=0,1,2$, the frame ${\frak s}$ has the brimmed$^\ell$ weak
$(\lambda,n)$-uniqueness property for $\ell=1,2,3,4$. \nl
3) For $n=0,1$ the frame ${\frak s}$ has the brimmed$^\ell$ strong
$(\lambda,n)$-existence property for $\ell=1,2,3,4$. 
\red
also for $\pm 1$

\endred

4) For $n=0,1$, the frame ${\frak s}$ has the brimmed$^\ell
(\lambda,n)$-primeness property for $\ell=1,2,3,4$. \nl
5) For $n=0,1$ the frame ${\frak s}$ has the brimmed$^\ell
(\lambda,n)$-primeness existence property for $\ell=1,2,3,4$. \nl
6) If ${\frak s}$ has the brimmed$^\ell$ strong 
$(\lambda,n)$-existence property and the brimmed$^\ell$  
weak/strong $(\lambda,n)$-primeness property, 
\ub{then} it has the brimmed$^\ell$ weak/strong 
$(\lambda,n)$-primeness existence property \scite{705-12f.4H}. 
\endproclaim
\bigskip

\demo{Proof}  1) If $n=0$, this is clear, for $n=1$ this is the
existence theorem for $K^{3,\text{bu}}_{\frak s}$, see \sciteu{705-xxX}.
\nl
Lastly, if $n=2$ by \sciteu{705-yyY}; ``the existence of stable amalgamation" we can
find a stable $(\lambda,{\Cal P}(n)$-system $\bold s^+$ such that
$\bold s^+ \restriction {\Cal P}^-(n) = \bold s^{\bold d}$ and we can
extend $M^{{\bold s}^+}_n$ and choose $\bold J_{n,u}$ for $u \subset n$.
\nl
2) For $n=2$, by the uniqueness 
of {\rm NF}$_{\frak s}$-amalgamation, see \sciteu{705-yyY}. \nl
3), 4), 5)  Easy, too. \nl
6) Let $\bold d$ be a brimmed$^\ell \, (\lambda,{\Cal
P}^-(n))$-system.  By the strong $(\lambda,n)^\ell$-existence property
there is an expanded stable $(\lambda,{\Cal P}(n))$-system $\bold
d^*$, reduced in $n$ such that $\bold d^* \restriction {\Cal P}^-(n) =
\bold d$.  By the $(\lambda,n)$-primeness$^\ell$ property $\bold d^*$
is prime$^\ell$ so we are done.   \hfill$\square_{\scite{705-12p.15B}}$\margincite{705-12p.15B}
\enddemo
\bigskip

\remark{Remark}  Assume ${\frak K}_{\frak s}$ is the class of $(A,E),|A| =
\lambda_{\frak s},E$ an equivalence relation.  If $(M,N,a) \in
K^{3,\text{bu}}_{\frak s},a/E^N$ disjoint to $M$, if $a/E^N$ is too
large then $(M,N,a) \notin K^{3,\text{pr}}_{\frak s}$.
\endremark
\bigskip

\proclaim{\stag{705-12p.16A} Claim}  The following properties of $\bold d$
are actually properties of $\bold s$, that is, their satisfaction
depends just on $\bold s^{\bold d}$
\mr
\item "{$(A)$}"  ``$\bar f$ being an embedding of $\bold d$ into $M$" 
\sn
\item "{$(B)$}"  $\bold d$ is brimmed$^\ell$ at $u$ for $\ell=1,2,3,4$
\sn
\item "{$(C)$}"  $\bold d$ has weak/strong uniqueness/existence
property
\sn
\item "{$(D)$}"  $\bold d$ has the weak/strong primeness existence
property
\sn
\item "{$(E)$}"  $\bold d$ has the weak/strong primeness property.
\endroster
\endproclaim
\bigskip

\demo{Proof}  The least easy case is that replacing $\bold J^{\bold
d}_{v,u}$ by similar $\bold J'_{v,u}$ does not make a difference which
holds by \sciteu{705-xxX}.  \hfill$\square_{\scite{705-12p.16A}}$\margincite{705-12p.16A}
\enddemo
\bn
\margintag{705-12.16K}\ub{\stag{705-12.16K} Discussion}:  Why do we define the ``strong
primeness", ``strong prime$^\ell$" existence?

The problem arises in \scite{705-12f.4H}.  Assume $\bold d_0$ is a
brimmed$^3$ expanded stable 
$(\lambda,{\Cal P}(n))$-system, $\bold d_1$ is an expanded
stable $(\lambda,{\Cal P}(n))$-system, $\bold d_1 \restriction {\Cal
P}^-(n) = \bold d_0 \restriction {\Cal P}^-(n),\bold d_1$ is reduced
in $n,M^{{\bold d}_1}_n \le_{\frak s} M^{{\bold d}_0}_n$ and $h^{{\bold
d}_1}_{n,u} = h^{{\bold d}_0}_{n,u}$ for $u \subset n$.  Is
$M^{{\bold d}_0} \, (\lambda,*)$-brimmed over $M^{{\bold d}_1}_n$?  If
${\frak s}$ has the NDOP (and $n \ge 2$) yes, but in general for $p
\in {\Cal S}^{\text{bs}}(M^{{\bold d}_1}_n)$ which is $\perp h^{{\bold
d}_1}_{n,u}(M^{{\bold d}_1}_u)$ for every $u \subset n$, we do not
know that dim$(p,M^{{\bold d}_0}_n) = \lambda_{\frak s}$.  This
motivates the definition of strong primeness.
\bigskip

\proclaim{\stag{705-11b.7} Claim}  1) Assume $I_1 \subseteq I_2$ and
${\frak s}$ has the brimmed$^\ell$ weak $(\lambda,|u|)$-existence property
whenever $u \in I_2 \backslash I_1$.  \ub{Then} for any brimmed$^\ell$
expanded stable $(\lambda,I_1,{\frak s})$-system 
$\bold d_1$ there is a brimmed$^\ell$ expanded stable
$(\lambda,I_2,{\frak s})$-system $\bold d_2$ 
satisfying $\bold d_2 \restriction I_1 = \bold d_1$. 
\nl
2) Let $\ell \in \{3,4\}$.  
Assume that for any $m<n,{\frak s}$ has the brimmed$^\ell$ strong
$(\lambda,m)$-uniqueness property.  \ub{Then} for
any two brimmed$^\ell$ expanded stable $(\lambda,{\Cal P}^-(n))$-systems 
$\bold d_1,\bold d_2$, the systems $\bold s[\bold
d_1],\bold s[\bold d_2]$ are isomorphic.  Similarly for $(\lambda,I)$
if $u \in I \Rightarrow |u| < n$.  \nl
3) In (2) if $\bold d_k$ is a brimmed$^\ell$ expanded stable
$(I_2,{\frak s})$-system, for $k =1,2$ and $I_1 \subseteq I_2$ and
${\frak s}$ has the brimmed$^\ell$ strong 
$(\lambda,|u|)$-uniqueness whenever $u \in I_2 \backslash J_1$ 
and $\bar f = \langle f_u:u \in I_1 \rangle$ is an
isomorphism from $\bold s[\bold d_1 \restriction I_1]$ onto $\bold
s[\bold d_2 \restriction I_1]$, \ub{then} we can find $\bar f'$, an
isomorphism from $\bold s[\bold d_1]$ onto $\bold s[\bold d_2]$ such
that $\bar f' \restriction I_1 = \bar f$.
\endproclaim
\bigskip

\demo{Proof}  Natural.  \hfill$\square_{\scite{705-11b.7}}$\margincite{705-11b.7}
\enddemo
\bigskip

\demo{\stag{705-12b.9} Conclusion}  1) Assume $\ell \in \{3,4\}$ and 
\mr
\item "{$(a)$}"  ${\frak s}$ has the brimmed$^\ell$ strong
$(\lambda,<n)$-uniqueness property
\sn
\item "{$(b)$}"  there is a brimmed$^\ell$ stable 
$(\lambda,{\Cal P}(n))$-system.
\ermn
\ub{Then} ${\frak s}$ has the brimmed$^\ell$ weak
$(\lambda,n)$-existence property. 
\enddemo
\bigskip

\demo{Proof}  1) We prove this by induction on $n$, so we can assume that
${\frak s}$ has the brimmed$^\ell$ weak 
$(\lambda,m)$-existence property for $m < n$.  
The cases $n=0,1$ are trivial, see \scite{705-12p.15B}.  
Now by clause (b) there is a brimmed$^\ell$ expanded 
stable $(\lambda,{\Cal P}(n))$-system $\bold d^*$.  
To prove the brimmed$^\ell$ weak
$(\lambda,n)$-existence property, let $\bold d$ be a 
brimmed$^\ell$ expanded stable $(\lambda,{\Cal P}^-(n))$-system.  
By assumption (a) and \scite{705-11b.7}(2)
and the induction hypothesis, $\bold d$ and 
$\bold d^* \restriction {\Cal P}^-(n)$
are isomorphic hence, by renaming, without loss of generality they equal, so
$\bold d^*$ prove the existence.  \hfill$\square_{\scite{705-12b.9}}$\margincite{705-12b.9}
\enddemo
\bigskip

\proclaim{\stag{705-12b.6} Claim}  1) Let $\ell = 1,2,3,4$.  
The brimmed$^\ell$ weak 
$(\lambda,n)$-existence property is equivalent to: for 
every brimmed$^\ell$ expanded 
stable ${\Cal P}^-(n)$-system $\bold d$ there are 
$M \in {\frak K}_{\frak s}$ and
$\bar f = \langle f_u:u \in {\Cal P}^-(n) \rangle$, such that
\mr
\item "{$\circledast$}"  $\bar f$ is an embedding of $\bold d$ into
$M$ which means
{\roster
\itemitem{ $(a)$ }   $f_u$ is a $\le_{\frak s}$-embedding of $M^{\bold
d}_u$ into $M$
\sn
\itemitem{ $(b)$ }  if $u \subset v \in {\Cal P}^-(n)$ then $f_u = f_v
\circ f^{\bold d}_{v,u}$
\sn
\itemitem{ $(c)$ }  for any $u \in {\Cal P}^-(n)$, the set
$\cup\{f_v(\bold J^{\bold d}_{v,u}):v$ satisfies 
$u \subseteq v \in {\Cal P}^-(n)\}$
is independent in $(f_u(M^{\bold d}_u),M)$, as an indexed set
\footnote{this just means that $\langle \{f_v(c):c \in \bold
J^{\bold d}_{v,u}\}:u \subset v \in {\Cal P}^-(n) \rangle$ is a
sequence of pairwise disjoint sets}.
\endroster}
\ermn
2)
\mr
\item "{$(a)$}"  If ${\frak s}$ has the brimmed$^\ell$ strong
$(\lambda,n)$-existence property \ub{then} ${\frak s}$ has the brimmed$^\ell$
weak $(\lambda,n)$-existence property
\sn
\item "{$(b)$}"  if ${\frak s}$ has the brimmed$^\ell$ strong
$(\lambda,n)$-uniqueness property \ub{then} ${\frak s}$ has the
brimmed$^\ell$ weak $(\lambda,n)$-uniqueness property
\sn
\item "{$(c)$}"  if ${\frak s}$ has the brimmed$^\ell$ strong
$(\lambda,n)$-primeness existence property \ub{then} ${\frak s}$ has the
brimmed$^\ell$ weak $(\lambda,n)$-primeness existence property
\sn
\item "{$(d)$}"  similarly for primeness.
\ermn
3) If $\{(\ell(1),\ell(2))\} \in \{(1,2),(2,3),(3,4)\}$ and
${\frak s}$ has the brimmed$^{\ell(1)}$ weak/strong 
$(\lambda,n)$-existence property
\ub{then} ${\frak s}$ has the brimmed$^{\ell(2)}$ weak/strong
$(\lambda,n)$-existence property. \nl
4) Similarly to part (3) for weak $(\lambda,n)$-uniqueness. \nl
4A)  If $\{(\ell(1),\ell(2))\} \in \{(1,2),(2,3),(3,4))\}$ and ${\frak s}$ 
has the
brimmed$^{\ell(1)}$ strong uniqueness, \ub{then} ${\frak s}$ has the
brimmed$^{\ell(2)}$ strong uniqueness. \nl
5)$^*$ 
\mr
\item "{$(a)$}"  every expanded stable 
$(\lambda_{\frak s},I)$-system for ${\frak s}$ is brimmed$^2$
\sn
\item "{$(b)$}"  for each of the properties defined in Definition
\scite{705-12b.5}, the weak brimmed$^2$ version and the brimmed$^1$ one are
equivalent.
\endroster
\endproclaim
\bigskip

\demo{Proof of \scite{705-12b.6}}  1) First assume that ${\frak s}$ has the
brimmed$^\ell$ weak $(\lambda,n)$-existence property and we shall prove the
condition $\circledast$ in \scite{705-12b.6}(1).  
By the present assumption, we can find
$\bold d'$ as in the definition \scite{705-12b.5}(1), and 
we let $M =  M^{\bold d}_n,f_u =
f^{{\bold d}'}_{n,u}$, clearly they are as required in $\circledast$.

Second assume that ${\frak s}$ satisfies the condition $\circledast$ from
\scite{705-12b.6}(1).  Let $\bold d = \bold d_0$, a brimmed$^\ell$
expanded stable $(\lambda,{\Cal P}^-(n))$-system be given, so 
by our assumption there
are $M,\langle f_u:u \in {\Cal P}^-(n))$ as there.  Now we define an expanded
stable $(\lambda,{\Cal P}(n))$-system $\bold d_1$ as follows:  $\bold
d_1 \restriction {\Cal P}^-(n) = \bold d_0,M^{{\bold d}_1}_n =
M,h^{{\bold d}_1}_{n,u} = f_u$ for $u \subset n$ let $\bold J^{{\bold
d}_1}_{n,u}$ be a maximal subset $\bold J$ of $M \backslash
\cup\{f_v(M_v):v \subset n\}$ such that
$\bold J^{{\bold d}_1}_{v,u} = 
\cup\{f_v(J^{{\bold d}_0}_{v,u}):v \text{ satisfies } u \subset v \subset
n\}$ is independent.  Clearly $\bold d_1$ is as required except that in
the brimmed$^2$ case we are missing ``$M$ is brimmed", and in the
brimmed$^3$ case,
we are missing the $\bold J_{n,u}$ and in the brimmed$^4$ case the relevant
condition.  Choose $M^*$ such that $M \le_{\frak s}
M^*$ and $M^*$ is $(\lambda_{\frak s},*)$-brimmed over $M$ letting
$M^{{\bold d}_1}_n = M^*$ and lastly, for each $u \in {\Cal P}^-(n)$ let 
$\bold J^{{\bold d}_1}_{n,u}$ be a maximal $\bold J$ such that:
\mr
\item "{$(\alpha)$}"   $\bold J \subseteq M^* \backslash
\cup\{f_u(\bold J^{\bold d}_{v,u}):v \text{ satisfies } u 
\subset v \in {\Cal P}^-(n)\}$
\sn
\item "{$(\beta)$}"  $c \in \bold J \Rightarrow \text{ tp}_{\frak
s}(c,f_u(M^{{\bold d}_0}_u),M^*) \in {\Cal S}^{bs}(f_u(M^{{\bold
d}_0}_u)$ is orthogonal to $f_w(M^{{\bold d}_0}_w)$ for every $w
\subset u$
\sn
\item "{$(\gamma)$}"  $\bold J \cup \bigcup\{f_{w_1}(\bold J^{\bold
d}_{w_1,u}):u \subset w_1 \in {\Cal P}^-(n)\}$ is independent in
$(f_u(M^{{\bold d}_\ell}_u),M^*)$
\sn
\item "{$(\delta)$}"  $\bold J$ is maximal under $(\alpha) + (\beta)$.
\ermn
By \scite{705-12b.3}(7) we are done. \nl
2), 3), 4), 5)  Left to the reader.   \hfill$\square_{\scite{705-12b.6}}$\margincite{705-12b.6}
\enddemo
\bn
The following claim will give a crucial ``saving".
\proclaim{\stag{705-12.u.9} Claim}  Assume that ($\ell =3$ or just $\ell \in
\{1,2,3,4\}$ and $n \ge 2$ and)
\mr
\item "{$(a)$}"  ${\frak s}$ has the brimmed$^\ell$ weak
$(\lambda,n)$-uniqueness
\sn
\item "{$(b)$}"  $\bold d$ is an expanded stable 
$(\lambda,{\Cal P}^-(n))$-system
\sn
\item "{$(c)$}"  $\bold d \restriction [n]^{< n-1}$ is a
brimmed$^\ell$ system (but $\bold d$ not necessarily).
\ermn
\ub{Then} $\bold d$ has weak uniqueness, i.e.
\mr
\item "{$\circledast$}"  if $\bold d_1,\bold d_2$ are expanded stable
$(\lambda,{\Cal P}(n))$-system satisfying $\bold d_1 \restriction
{\Cal P}^-(n) = \bold d,\bold d_2 \restriction {\Cal P}^-(n)$ \ub{then} we can
find $(N,f)$ such that $M^{{\bold d}_2}_n \le_{\frak s} N,f$ is a
$\le_{\frak s}$-embedding of $M^{{\bold d}_1}_n$ into $N$ and $u
\subset n \Rightarrow h^{{\bold d}_2}_{n,u} = f \circ h^{{\bold d}_1}_{n,u}$. 
\endroster
\endproclaim
\bigskip

\demo{Proof}  We prove this by induction on $k_{\bold d} = |{\Cal
P}_{\bold d}|$ where ${\Cal P}_{\bold d} = \{v \in {\Cal P}^-(n):\bold
d$ is not brimmed$^\ell$ in $u$ (so $u \in [n]^{n-1})\}$. 

Clearly $k_{\bold d} \le \binom{n}{n-1} = n (\le|{\Cal P}^-(n)| <
2^n)$.  If $k_{\bold d} = 0$ the conclusion follows from assumption
(a).

So assume that $k_{\bold d} >0$ and choose $v_* \in {\Cal P}_{\bold
d}$.  We can find $M$ which is $(\lambda,*)$-brimmed$^\ell$ over $M^{\bold
d}_{v_*}$.  Next we define an expanded $(\lambda,{\Cal
P}^-(n))$-system $\bold d^+$:
\mr
\item "{$\circledast_1$}"  $(a) \quad \bold d^+ \restriction ({\Cal
P}^-(n) \backslash \{v_*\}) = \bold d \restriction ({\Cal P}^-(n)
\backslash \{v^*\})$
\sn
\item "{${{}}$}"  $(b) \quad h^{{\bold d}^+}_{v_*,u} = 
h^{\bold d}_{v_*,u}$ if $u \subset v_*$
\sn
\item "{${{}}$}"  $(c) \quad M^{{\bold d}^+}_{v_*} = M$
\sn
\item "{${{}}$}"  $(d) \quad \bold J^{{\bold d}^+}_{v,u} = \bold
J^{\bold d}_{v,u}$ if $u \subset v \in {\Cal P}^-(n) \backslash
\{v_*\}$
\sn
\item "{${{}}$}"  $(e) \quad \bold J^{{\bold d}^+}_{v_*,u}$, for $u \subset
v_*$, is a maximal subset $\bold J$ of \nl

\hskip35pt $\{c \in M \backslash \cup\{M^{\bold d}_{v_*,u}:u \subset
v_*\}:\text{tp}(c,M^{\bold d}_{v_*,u},M) \in 
{\Cal S}^{\text{bs}}(M_{v,u})$ \nl

\hskip35pt is orthogonal to $M_{v_*,w}$ when 
$w \subset u\}$ such that $\bold J$ is \nl

\hskip35pt independent in $(M_{v_*,u},M_{v_*})$ and 
$\bold J \supseteq \bold J^{\bold d}_{v_*,u}$.
\ermn
Now let $\bold d_1,\bold d_2$ be as in the assumption of
$\circledast$.  Next for $k=1,2$ by the existence of stable amalgamation
there are $(N_k,f_k)$ such that $M^{{\bold d}_k}_n
\le_{\frak s} N_k,f_k$ is a $\le_{\frak s}$-embedding of $M =
M^{{\bold d}^+}_{v_*}$ into $N_k$ extending $f^{{\bold d}_k}_{n,v_*}$
and NF$_{\frak s}(M^{{\bold d}_k}_{n,v_*},M^{{\bold d}_k}_n,f_k(M),N_k)$
holds.

By renaming \wilog \, $f_1 = f_2 = \text{ id}_M$.
We now define an expanded stable $(\lambda,{\Cal P}^-(n))$-system
$\bold d^+_k$
\mr 
\item "{$\circledast^2_k$}"  $(a) \quad \bold d^+_k \restriction
{\Cal P}^-(n) = \bold d^+$
\sn
\item "{${{}}$}"  $(b) \quad M^{{\bold d}^+_k}_n = N_k$.
\sn 
\item "{${{}}$}"  $(c) \quad h^{{\bold d}^+_k}_{n,u} = h^{{\bold
d}_k}_{n,u}$ for $u \in {\Cal P}^-(n) \backslash \{v_*\}$
\sn
\item "{${{}}$}"  $(d) \quad h^{{\bold d}^+}_{n,v_*} = f_k$
\sn
\item "{${{}}$}"  $(e) \quad \bold J^{{\bold d}^+_k}_{n,v}$ for $v
\subset n$ are defined as in previous cases.
\ermn
Now we use the induction hypothesis (on $k_{\bold d}$)
\hfill$\square_{\scite{705-12.u.9}}$\margincite{705-12.u.9} 
\enddemo
\bigskip

\proclaim{\stag{705-12.u10} Claim}  1) Assume that (a) + (b) and: (c)$_1$ 
or (c)$_2$ where;
\mr
\item "{$(a)$}"  $\bold s^k$ is a stable $(\lambda,I)$-system for
$k = 1,2$
\sn
\item "{$(b)$}"  $\bold s^1 \restriction J = \bold s^2 \restriction J$
where $J =: \{u \in I:(\exists v \in I)(u \subset v)\}$
\sn
\item "{$(c)_1$}"  if $u \in I \backslash J$ \ub{then} $\bold s^\ell
\restriction {\Cal P}^-(u)$ has the weak uniqueness property
\sn
\item "{$(c)_2$}"  if $v \in I \backslash J$ \ub{then} $u \subset v
\Rightarrow h^{{\bold s}^1}_{v,u} = h^{{\bold s}^2}_{v,u}$ and $M^{{\bold
s}^\ell}_v \le_{\frak s} M^{{\bold s}^{2-\ell}}_v$ for some $\ell \in
\{1,2\}$.
\ermn
\ub{Then} 
\mr
\item "{$(\alpha)$}"  $\bold s^1$ has the weak existence property
\ub{iff} $\bold s^2$ has the weak existence property
\sn
\item "{$(\beta)$}"  $\bold s^1$ has the weak uniqueness property \ub{iff}
$\bold s^2$ has the weak uniqueness property.
\endroster
\endproclaim
\bigskip

\demo{Proof}  Similar to \scite{705-12.u.9}.
\hfill$\square_{\scite{705-12.u10}}$\margincite{705-12.u10}
\enddemo
\bn
In \scite{705-12b.10} we have proved actually some things on ${\frak s}^+$
concerning Definition \scite{705-12b.5}.  
[Why do we ignore $\ell=1$?  As for ${\frak s}^+$,
brimmed$^1 \Rightarrow \text{ brimmed}^2$.]
\demo{\stag{705-12f.10D} Conclusion}  Let $\ell \in \{0,2,3\}$ and assume
that ${\frak s}$ is successful hence
${\frak s}^+$ satisfies the demands in \scite{705-12.1}.  If for ${\frak
s}^+$ there is a brimmed$^\ell$ stable ${\Cal P}(n)$-system,
\ub{then} for ${\frak s}$ there is a brimmed$^\ell$ stable 
${\Cal P}(n+1)$-system. 
\enddemo
\bigskip

\remark{Remark}  For $\ell=0,2$ we can find trivial examples.
\endremark
\bigskip

\demo{Proof}  By \scite{705-12b.10}.
\enddemo
\bigskip

\proclaim{\stag{705-12p.21A} Claim}   Assume that $\bold d_k$ is an
expanded stable $(\lambda,I)$-system for $k=1,2$ and $\bar f$ is an
embedding of $\bold d_1$ into $\bold d_2$, see Definition
\scite{705-11b.2}(1A). \nl
1) If $I = {\Cal P}^-(n)$ and $\bold d_2$ has the weak existence property
\ub{then} $\bold d_1$ has the weak existence property. \nl
2) If $I = {\Cal P}^-(n)$ and $\bar f \restriction [n]^{< n-1}$ is an
isomorphism from $\bold d_1 \restriction [n]^{< n-1}$ onto $\bold d_2
\restriction [n]^{< n-1}$ and $\bold d_2$ has the weak uniqueness property
\ub{then} $\bold d_1$ has weak uniqueness property.
\endproclaim
\bigskip

\demo{Proof}  1) Easy (and was used inside the proof of \scite{705-12.u.9}). \nl
2) Very similar to the proof of \scite{705-12.u.9} (and we can use part (1)).
\hfill$\square_{\scite{705-12p.21A}}$\margincite{705-12p.21A} 
\enddemo
\bigskip

\proclaim{\stag{705-12b.11} Lemma}  1) Let $\ell =3,n \ge 2$.  Assume
$2^\lambda < 2^{\lambda^+}$ and:
\mr
\item "{$(a)$}"  ${\frak s}$ is successful hence ${\frak s}^+$ has the 
properties required in \scite{705-12.1}
\sn
\item "{$(b)$}"  ${\frak s}$ has the brimmed$^\ell$ weak
$(\lambda,\le n+1)$-existence property
\sn
\item "{$(c)$}"  ${\frak s}$ has the brimmed$^\ell$ strong $(\lambda,\le
n)$-uniqueness property
\sn
\item "{$(d)$}"  ${\frak s}$ does not have the brimmed$^\ell$ weak
$(\lambda,n+1)$-uniqueness property.
\ermn
\ub{Then} ${\frak s}^+$ does not have the brimmed$^\ell$ strong
$(\lambda^+,n)$-uniqueness property. 
\endproclaim
\bigskip

\remark{Remark}  Of course, it would be better to have ``strong" in
clause (d) of the assumption and it
would be better to have weak in the conclusion.  Still we can prove a
slightly stronger claim.
\endremark
\bn
A variant of \sciteu{705-12.b11} is
\proclaim{\stag{705-12v.11} Claim}  1) In \scite{705-12b.11} we can replace clause
$(c)$ by $(c)^-_1 + (c)^-_2$ (which obviously follows from it) where
\mr
\item "{$(c)^-_1$}"  ${\frak s}$ has the brimmed$^\ell$ strong
$(\lambda,<n)$-uniqueness property
\sn
\item "{$(c)^-_2$}"  ${\frak s}$ has the brimmed$^\ell$ weak
$(\lambda,n)$-uniqueness property.
\ermn
2) We can strengthen the conclusion to: ${\frak s}^+$ fail the
brimmed$^3$ weak $(\lambda^+,n)$-uniqueness property.
\endproclaim
\bigskip

\demo{Proof of \scite{705-12b.11}}  1) We choose 
by induction on $\alpha < \lambda^+,\bold
s_\eta$ for every $\eta \in {}^\alpha 2$ such that:
\mr
\item "{$\circledast$}"  $(\alpha) \quad 
\bold s_\eta$ is a normal brimmed$^\ell$ stable $(\lambda,{\Cal P}(n))$-system
\sn
\item "{${{}}$}"   $(\beta) \quad$ the 
universe of $M^{{\bold s}_\eta}_n$ is an ordinal $\gamma_{\ell g(\eta)}
< \lambda^+$
\sn
\item "{${{}}$}"   $(\gamma) \quad$ the sequence $\langle M^{{\bold s}_{\eta
\restriction \gamma}}_u:\gamma \le \alpha \rangle$ is $\le_{\frak
s}$-increasing continuous
\sn
\item "{${{}}$}"  $(\delta) \quad$ if $\alpha = \beta +1$ then 
$\bold s^*_\eta =: \bold s_{\eta \restriction \beta} *_{{\Cal P}(n+1)}
\bold s_\eta$ is a brimmed$^\ell$ stable \nl

\hskip25pt  $(\lambda,{\Cal P}(n+1))$-system
\sn
\item "{${{}}$}"  $(\varepsilon) \quad$ if 
$\alpha = \beta +1,\nu \in {}^\beta 2$, \ub{then}:
\nl

\hskip25pt  $\bold s^*_{\nu \char 94 <0>} \restriction
{\Cal P}^-(n+1) = \bold s^*_{\nu \char 94 <1>} \restriction 
{\Cal P}^-(n+1)$  
\sn
\item "{${{}}$}"  $(\zeta) \quad$ if $\alpha = \beta +1$ and $\nu \in
{}^\beta 2$ \ub{then} for no $f,N$ do we have: \nl

\hskip25pt $f$ is an $\le_{\frak s}$-embedding 
of $M^{{\bold s}^*_{\nu \char 94 <0>}}_{n+1}$ \nl

\hskip25pt into $N,M^{{\bold s}^*_{\nu 
\char 94 <1>}}_{n+1} \le_{\frak s} N$ \nl

\hskip25pt and $f$ is the identity on $M^{{\bold s}^*_{\nu \char 94 <0>}}_u = 
M^{{\bold s}^*_{\nu \char 94 <1>}}_u$ for $u \in {\Cal P}^-(n+1)$ 
\sn
\item "{${{}}$}"  $(\eta) \quad$ if $\nu_1,\nu_2 \in {}^\alpha 2$ then
$\bold s_{\nu_1} \restriction {\Cal P}^-(n) = \bold s_{\nu_2}
\restriction {\Cal P}^-(n)$ \nl

\hskip25pt (this strengthens clause $(\varepsilon)$).
\ermn
Now
\mr
\item "{$(*)_1$}"  we can carry the induction \nl
[Why?  For $\alpha=0$ trivial.  For $\alpha = \beta +1$ and $\eta \in
{}^\beta 2$ by clause (d) of the assumption there are 
brimmed$^\ell \, {\Cal P}(n+1)$-system $\bold s',
\bold s''$ with $\bold s' \restriction {\Cal P}^-(n+1) = 
\bold s'' \restriction {\Cal P}^-(n+1)$ as in clause $(\zeta)$, i.e.,
for no $(N,f)$ do we have $M^{{\bold s}''}_n \le_{\frak s} N$ and $f$
is a $\le_{\frak s}$-embedding of $M^{{\bold s}'_n}$ into $N$ which is
the identity on $M^{{\bold s}'}_u$ for every $u \subset n+1$.  Now for
proving \scite{705-12b.11}, by assumption (c) and claim \scite{705-11b.7}(2) 
renaming we have $\bold
s_{\eta \char 94 <0>},\bold s_{\eta \char 94 <1>}$ as required in
clause $(\varepsilon)$ and $(\zeta)$ of $\circledast$ and by renaming
we have $(\beta)$; by the proof we can get clause $(\eta)$, too.  
(For proving \scite{705-12v.11} this is done in its proof).
For $\alpha$ limit $\bold s_\eta =: \dbcu_{\beta < \alpha}
\bold s_{\eta \restriction \beta}$ is a stable $(\lambda,
{\Cal P}(n))$-system  by Claim \scite{705-12f.4D}(2); moreover is
brimmed$^\ell$.]
\ermn
For $\eta \in {}^{\lambda^+}2$ we define a normal $(\lambda^+,{\Cal
P}(n))$-system for ${\frak s}^+$ called $\bold s_\eta$ by
$M^{{\bold s}_\eta}_u = \cup\{M^{{\bold s}_{\eta \restriction
\alpha}}_u:\alpha < \lambda^+\}$, clearly 
\mr
\item "{$(*)_2$}"   $\bold s_\eta$ is really a $(\lambda^+,
{\Cal P}(n))$-system for the frame ${\frak s}^+$ 
(noting that $M^{{\bold s}_\eta}_u \in K_{{\frak s}(+)}$ as it belongs
to $K^{\frak s}$ and is saturated over $\lambda_{\frak s}$ because by
clause $(\delta)$, $M^{{\bold s}_{\eta \restriction(\alpha +1)}}_u$ is
$(\lambda,*$)-brimmed over $M^{{\bold s}_{\eta \restriction
\alpha}}_u$ for every $\alpha < \lambda^+_{\frak s}$)
\sn
\item "{$(*)_3$}"  $\bold s_\eta$ is a brimmed stable
$(\lambda^+,{\Cal P}(n))$-system for ${\frak s}(+)$ 
\nl
[why?  by \scite{705-12f.4D}(3).]
\sn
\item "{$(*)_4$}"  $\bold s_\eta \restriction {\Cal P}^-(n)$ is the
same for all $\eta \in {}^{\lambda^+} 2$ call it $\bold t^*$. \nl
[Why?  By Clause $(\eta)$ of $\circledast$.]
\ermn
Let $\rho \in {}^{\lambda^+} 2$.  To 
finish the proof it is enough to find $\eta
\in {}^{\lambda^+}2$ such that $h_* = \cup\{\text{id}_{M^{\bold
t}_u}:u \subset n\}$ cannot be extended to an
isomorphism from $M^{{\bold s}_\rho}_n$ onto $M^{{\bold s}_\eta}_n$;
toward contradiction assume that $f_\eta$ is such an isomorphism for
every $\eta \in {}^{\lambda^+}2$.  By the weak diamond, (see
\sectioncite[\S0]{88}) for some
$\eta_0,\eta_1 \in {}^{\lambda^+}2$ and $\delta < \lambda^+$ we have
$\nu = \eta_\ell \restriction \delta,\nu \char 94 \langle \ell \rangle
\triangleleft \eta_\ell$ and $f_{\eta_1} \restriction M^{{\bold
s}_\nu}_n = f_{\eta_2} \restriction M^{{\bold s}_\nu}_n$.  Clearly
we get contradiction to clause $(\zeta)$ in the construction. 
\hfill$\square_{\scite{705-12b.11}}$\margincite{705-12b.11}
\enddemo
\bigskip

\demo{Proof of \scite{705-12v.11}}  In the proof of \scite{705-12b.11} there
one point in which the proofs differ.  We are 
given $\bold s_\eta$ for $\eta \in {}^\beta 2$ and we
know that there are normal brimmed$^\ell {\Cal P}(n+1)$-system $\bold
s',\bold s''$ such that $\bold s' \restriction {\Cal P}^-(n+1) = \bold
s'' \restriction {\Cal P}^-(n+1)$ but there is no $\le_{\frak
s}$-embedding of $M^{{\bold s}'}_{n+2}$ into any $N,M^{{\bold
s}''}_{n+1} \le_{\frak s} N$ over $\dbcu_{u \subset n^+} M^{{\bold
s}'}_u$.  By the amount of uniqueness we have, i.e. by assumption
$(c)^-_1$ \wilog \, $\bold s'
\restriction {\Cal P}^-(n) = \bold s_\eta \restriction {\Cal P}^-(n)$
(for every $\eta \in {}^\beta 2$, hence also $\bold s'' \restriction
{\Cal P}^-(n) = \bold s_\eta$).
Without loss of generality the universe of $M^{\bold s'}_{n+1}$ and of
$M^{{\bold s}''}_{n+1}$ is $\gamma_\beta + \lambda$ (recall
\scite{705-12b.3}(8)) and, of course, the universe of $M^{\bold s'}_n =
M^{{\bold s}''}_n$ is $\gamma_\eta$. 

Now we define $\bold s^*_\eta$, a  
$(\lambda,{\Cal P}^-(n+1))$-system by 
$\bold s^*_\eta \restriction ({\Cal P}(n+1)
\backslash \{n,n+1\}) = \bold s' \restriction ({\Cal P}^*(n+1)
\backslash \{n,n+1\}) = \bold s'' \restriction ({\Cal P}^-(n+1)
\backslash \{n,n+1\})$ and $M^{{\bold s}^*_\eta}_n= 
M^{{\bold s}_n}_\eta$.  Clearly $\bold s^*_\eta$ is stable and brimmed$^\ell$.
\nl
Now by \scite{705-12.u10}, the version with $(c)_1$ as
$\bold s' \restriction {\Cal P}^-(n+1)$ fails the brimmed$^3$ weak
uniqueness, also $\bold s^*_\eta$ fails it.  Hence we can find a
stable $(\lambda,{\Cal P}^*(n+1))$-systems $\bold s'_\eta,\bold
s''_\eta$ witnessing it, to $\bold s'_\eta \restriction {\Cal
P}^-(n+1) = \bold s^*_\eta = \bold s''_n \restriction {\Cal P}^-(n+1)$. By
renaming we take care of clause $(\beta)$ of $\circledast$ (we use
freely \scite{705-12b.3}). \nl
2) We choose in addition to $\bold s_\eta$ also $N_\eta$ such that
\mr
\item "{$\circledast$}"  $(\alpha),(\gamma)-(\eta) \quad$ as in the proof
\sn
\item "{${{}}$}"  $(\beta) \quad N_\eta$ is brimmed over $M^{{\bold
s}_n}_\eta$ the universe of $N_\eta$ is an ordinal \nl

$\qquad$ $\gamma_{\ell g(\eta)}$ (instead of $(\beta))$
\sn 
\item "{${{}}$}"  $(\theta) \quad$ if $\nu \triangleleft \eta$ then
NF$_{\frak s}(M^{{\bold s}_\nu}_n,N_\nu,M^{{\bold s}_\eta}_n,N_\eta)$
and $N_\eta$ is brimmed over $M^{{\bold s}_\eta}_n \cup N_\nu$.
\ermn
In the end for $\eta \in {}^{\lambda^+}2$ we define also $N_\eta =
\cup\{N_{\eta \restriction \alpha}:\alpha < \lambda^+\}$ hence
$M^{\bold s_\eta} \le_{{\frak s}(+)} N_\eta$ and $N_\eta$ is
$(\lambda^+,*)$-brimmed over $M^{\bold s_\eta}_n$.
\nl
So if $\rho,\eta \in {}^{\lambda^+} 2,N_\eta \le_{{\frak s}(+)} N$ and
$f$ is a $\le_{{\frak s}(+)}$-embedding of $M^{{\bold s}_\rho}_n$ into
$N$ over $\cup\{M^{{\bold s}_\rho}_u:u \subset n\}$, \ub{then} $N$ can
be $\le_{{\frak s}(+)}$-embedded into $N_\eta$ over $M^{{\bold
s}_\eta}_n$ so \wilog \, $N \le_{{\frak s}(+)} M^{{\bold
s}_\eta}_n,N_\eta$ is $(\lambda^+,*)$-brimmed over  $M^{\bold
s_\eta}_n$.  The rest should be clear.  \hfill$\square_{\scite{705-12v.11}}$\margincite{705-12v.11}
\enddemo
\bigskip

\proclaim{\stag{705-12f.10E} Claim}   1) Assume $\ell = 2,3,n \ge 1$ and
\mr
\item "{$(a)$}"   ${\frak s}$ has the brimmed$^\ell$ strong 
$(\lambda,\le n)$-existence property
\sn
\item "{$(b)$}"  ${\frak s}$ has the brimmed$^\ell$ weak 
$(\lambda,\le n)$-primeness property
\ermn
\ub{Then} there is an expanded stable ${\Cal P}(n+1)$-system $\bold d$
reduced at $n+1$ such that $\bold d \restriction {\Cal P}^-(n)$ is
brimmed$^3$.
\nl
2) If in addition clause (c) below holds then
${\frak s}$ has the brimmed$^\ell$ strong $(\lambda,n+1)$-existence
property where
\mr
\item "{$(c)$}"  ${\frak s}$ has the brimmed$^\ell$ strong 
$(\lambda,\le n)$-uniqueness property.
\endroster
\endproclaim
\bigskip

\demo{Proof}  1) Let $M_\emptyset \in K_{\frak s}$ be brimmed$^\ell$.  Let $M
\in K_{\frak s}$ be $(\lambda_{\frak s},*)$-brimmed over
$M_\emptyset$.  Let ${\Cal P}_\emptyset \subseteq \{p \in {\Cal
S}^{\text{bs}}(M_\emptyset):p$ regular$\}$ be a maximal family of pairwise
orthogonal types.  Let $\bold J^2_\emptyset$ be a maximal subset of
$\{c \in M:c$ realizes some $p \in {\Cal P}^2_\emptyset$ in $M$
over $M_\emptyset\}$ which is independent in $(M_\emptyset,M)$ 
so $|\bold J_\emptyset| = \lambda_{\frak s}$ and let
$\{\bold J^2_{u,\emptyset}:\emptyset \subset u \subset n+1\}$ be a partition of
$\bold J^2_\emptyset$ to sets each of cardinality $\lambda_{\frak s}$
such that for every $p \in {\Cal P}^2_\emptyset$ and $u \subset n+1$ the
set $\{c \in \bold J_u$: tp$_{\frak s}(c,M_\emptyset,M) = p\}$
has cardinality $\lambda_{\frak s}$.  For $u \in I_1 \backslash I_0$
let $M_u \le_{\frak s} M$ be such that $(M_\emptyset,M_u,\bold
j_{u,\emptyset}) \in K^{3,\text{bu}}_{\frak s}$. 

Let $I_k = \{u \subseteq n+1:|u| \le k\}$ for $k \le n+1$.  We now
choose by induction on $k \in \{1,2,\dotsc,n\}$, the objects $\bold
d_k$ and $\langle {\Cal P}_u,\bold J_{v,u}:u \in I_k,u \subseteq v
\subset n+1 \rangle,\langle N_u,{\Cal P}^1_u:u \in I_k,|u|>1 \rangle,\langle 
\bold J^i_{u,v}:u \subseteq v \subset n+1$ and $(u,v) \ne 
(\emptyset,\emptyset)$ and $i=1,2 \rangle$ such that
\mr
\item "{$\circledast(a)$}"  $\bold d_k$ is a normal brimmed$^\ell$ expanded
stable $(\lambda,I_k)$-system embedded in $M$ and $M^{{\bold
d}_k}_\emptyset = M_\emptyset,M^{{\bold d}_k}_u = M_u$ for $u \in I_1
\backslash I_0,J^{{\bold d}_u}_u = \bold J_{u,\emptyset}$
\sn
\item "{$(b)$}"  $m < k \Rightarrow \bold d_m = \bold d_k \restriction
I_m$
\sn
\item "{$(c)$}"  for $u \in I_k \backslash I_1,N_u$ 
is such that $\cup\{M^{{\bold
d}_k}_w:w \subset u\} \subseteq N_u \le_{\frak s} M$ and $\bold
d^*_u$ is reduced in $u$ where $\bold d^*_u$ is the normal brimmed$^1
{\Cal P}(u)$-system $\bold d^*_u \restriction {\Cal P}^-(u) = \bold d_k
\restriction {\Cal P}^-(u),M^{{\bold d}^*_u}_u = N_u,\bold J^*_{u,w} =
\emptyset$ for $w \subset u$
\sn
\item "{$(d)$}"   for $u \in I_k \backslash I_1,
{\Cal P}^1_u$ is a maximal set of
pairwise orthogonal types from $\subseteq \{p \in {\Cal
S}^{\text{bs}}(N_u):p$ regular orthogonal to $M^{{\bold d}_k}_w$ for $w
\subset u$ (hence to $M^{{\bold d}_k}_w$ for $w \in I_k$ such that $u
\nsubseteq w$ by \scite{705-10p.17})
\sn
\item "{$(e)$}"  for $u \in I_k \backslash I_1$, the set
$\bold J^1_u$ is a maximal subset of $\{c \in M:\text{tp}_{\frak
s}(c,N_u,M) \in {\Cal P}^1_u\}$ independent in $(N_u,M)$
\sn
\item "{$(f)$}"  if $u \in I_k \backslash I_1$, then $\langle \bold J^1_{v,u}:u
\subseteq v \subset n+1$ but $(v,u) \ne (\emptyset,\emptyset) \rangle$ is a
partition of $\bold J^1_u$ to sets each of cardinality $\lambda_{\frak
s}$ such that moreover, for every $p \in {\Cal P}^1_u$ the set $\bold
J^1_{v,u,p} =: \{c \in \bold J^1_{v,u}$: tp$_{\frak s}(c,N_u,M)=p\}$ has
cardinality $\lambda_{\frak s}$
\sn
\item "{$(g)$}"  $(N_u,M^{{\bold d}_k}_u,\dbcu_{w \subseteq u} 
\bold J^1_{u,w} \cup \dbcu_{w \subset u} \bold J^2_{u,w}) 
\in K^{3,\text{bu}}_{\bold s}$ for $u \in I_k
\backslash \{\emptyset\}$ where we stipulate $\bold J^1_{u,w} =
\emptyset$ when $w \subseteq u,|w| \le 1$
\sn
\item "{$(h)$}"  ${\Cal P}^2_u$ is a maximal set of pairwise
orthogonal types from $\{p \in {\Cal S}^{\text{bs}}(M^{{\bold
d}_k}_u):p$ orthogonal to $N_u$ when $u \notin I_1$ and to
$M_\emptyset$ if $u \in I_1 \backslash I_0\}$ when $u \in I_k$ (if $u
= \emptyset$ then ${\Cal P}^2_u$ has already been chosen)
\sn
\item "{$(i)$}"  $\bold J^2_u$ is a maximal set of 
$\{c \in M:\text{tp}_{\frak s}(c,M^{{\bold d}_k}_u,M) \in {\Cal
P}^2_u\}$ independent in $(M^{{\bold d}_k}_u,M)$ when $u \in I_k$,
(if $u = \emptyset,J^2_u$ has already been chosen), note that the
$\bold J^2_{u,w}$ used is clause (g) has already been chosen)
\sn
\item "{$(j)$}"  $\langle \bold J^2_{v,u}:u \subset v \in {\Cal P}^-(n+1)
\rangle$ is a partition of $\bold J^2_u$ such that for every $p \in
{\Cal P}^2_u$ and $v$ such that $u \subset v \in {\Cal P}^-(n+1)$ the
set $\{c \in \bold J^2_{v,u}:\text{tp}_{\frak s}(c,M^{{\bold
d}_k}_u,M)=p\}$ has cardinality $\lambda$ when $u \in I_k$
\sn
\item "{$(k)$}"  $\bold J^{{\bold d}_k}_{v,u} = \bold J^1_{v,u} \cup
\bold J^2_{v,u}$ where $u \subset v \in I_k$
\ermn
For $k=1,\bold d_k$ is defined by clause (a)
and choose ${\Cal P}^2_u,\bold J^2_u,\langle \bold J^2_{v,u}:v$ satisfies $u
\subseteq v \subset n+1 \rangle$ as above for $u \in
I_k \backslash I_0$ (i.e., $u = \{m\},m < n+1$).

For $k=m+1>1$ for each $u \in I_k \backslash I_m$ clearly $\bold d_m
\restriction {\Cal P}^-(u)$ is a normal brimmed$^\ell$ expanded stable ${\Cal
P}^-(u)$-system hence by assumption (a) we can find a normal
brimmed$^1 {\Cal P}(u)$-system $\bold d^*_u$ such that $\bold d^*_u
\restriction {\Cal P}^-(u) = \bold d_u \restriction {\Cal
P}^-(u),M^{{\bold d}^*_u} = N_u$ which is reduced in $u$ 
hence $\bold d^*_u$ is prime over $\bold d_m 
\restriction {\Cal P}^-(u)$, use (b).

Now $\bold d^*_u \restriction {\Cal P}^-(u)$ is embedded in $M$ 
(by clause (a) of $\circledast$) so by the definition of the
primeness \wilog \, $N_u =: M^{{\bold d}^*_u}_u \le_{\frak s} M$.  Now
as $|u| > 1$ choose ${\Cal P}^1_u,\bold J^1_u$ and $\langle \bold J^1_{v,u}:v$
satisfies $u \subseteq v \subset n+1 \rangle$ as required. 
Now let $M^{{\bold d}_k}_u \le_{\frak s} M$ be 
\footnote{exists as for every $w \subset u$, the set $\bold J_{u,w}$
is independent in $(N_u,M)$ because $(M_w,N_u,\cup\{\bold
J_{v_1,u_1}:v_1 \subset u,u_1 \subseteq w,u_1 \subset v_1\}) \in
K^{3,\text{bu}}_{\frak s}$}
such that $(N_u,M^{{\bold d}_k}_u,\cup\{\bold J^1_{u,w} \cup \bold
J^2_{u,w}:w \subset u\}) \in 
K^{3,\text{bu}}_{\frak s}$.  We then choose ${\Cal P}^2_u,\bold
J^2_u,\bold J^2_{v,u} \, (u \subset v \subset n+1)$ as required.

Having carried the induction we define a normal ${\Cal P}(n+1)$-system
$\bold d_{n+1}$ by

$$
\bold d_{n+1} \restriction {\Cal P}^-(n+1) = \bold d_n
$$

$$
M^{{\bold d}_{n+1}}_{n+1} = M
$$

$$
\bold J^{{\bold d}_{n+1}}_{n+1,u} = \emptyset \text{ for } u \subset n+1.
$$
\mn
It is easy to check that $\bold d_{n+1}$ is as required. \nl
2) Easy, by \scite{705-11b.7}(2).  \hfill$\square_{\scite{705-12f.10E}}$\margincite{705-12f.10E}
\enddemo
\bigskip

\proclaim{\stag{705-12p.10EA} Claim}  Assume
\mr
\item "{$(a)$}"  $\bold d$ is an expanded stable $(\lambda,{\Cal
P}(n))$-system
\sn
\item "{$(b)$}"  $\bold d$ is reduced at $n$.
\ermn
\ub{Then} we can find $\bold d'$ such that
\mr
\item "{$(\alpha)$}"  $\bold d'$ is an expanded stable $(\lambda,{\Cal
P}(n))$-system reduced at $n$
\sn
\item "{$(\beta)$}"  $\bold d' \restriction {\Cal P}^-(n) = \bold d
\restriction {\Cal P}^-(n)$
\sn
\item "{$(\gamma)$}"  $h^{{\bold d}'}_{n,u} = h^{\bold d}_{n,u}$ for
$u \subset n$
\sn
\item "{$(\delta)$}"  $M^{\bold d}_n \le_{\frak s} M^{{\bold d}'}_n$
\sn
\item "{$(\varepsilon)$}"  if $p \in {\Cal S}^{\text{bs}}(M^{\bold d}_n)$ is
regular orthogonal to $M^{\bold d}_{n,u}$ for every $u \subset n$
\ub{then} {\rm dim}$(p,M^{{\bold d}'}_n) = \lambda$.
\endroster
\endproclaim
\bigskip

\demo{Proof}  Let $M^+ \in K_{\frak s}$ be $(\lambda,*)$-brimmed over
$M^{\bold d}_n$, let $\bold P = \{p \in {\Cal S}^{\text{bs}}(M^{\bold
d}_n):p$ regular $\perp M^{\bold d}_{n,u}$ for every $u \subset n\}$.
Let $\bold J = \{c_{p,\alpha}:p \in \bold P,\alpha < \lambda\}$ be
such that
\mr
\widestnumber\item{$(iii)$}
\item "{$(i)$}"  $c_{p,\alpha} \in M^+$ realizes $p$
\sn
\item "{$(ii)$}"  $p \in \bold P \wedge \alpha \ne \beta \Rightarrow
c_{p,\alpha} \ne c_{p,\beta}$
\sn
\item "{$(iii)$}"  $\bold J$ is independent in $(M^{\bold d}_n,M^+)$.
\ermn
Now let $M \le_{\frak s} M^+$ be such that $(M^{\bold d}_n,M,\bold J) 
\in K^{3,\text{bu}}_{\frak s}$ and define $\bold d'$ such that
$((\alpha),(\beta),(\gamma)$ above holds and) $M^{{\bold d}'}_n =
M$.  It is easy to check that $\bold d'$ is as required. \nl
${{}}$  \hfill$\square_{\scite{705-12p.10EA}}$\margincite{705-12p.10EA}
\enddemo
\bigskip

\proclaim{\stag{705-12f.10F} Claim}  1) Assume $\ell = 3$ and
\mr
\item "{$(a)$}"  ${\frak s}$ is successful hence
${\frak s}^+$ satisfies the hypothesis \scite{705-12.1}
\sn
\item "{$(b)$}"  ${\frak s}$ has the brimmed$^\ell$ weak 
$(\lambda,\le n+1)$-uniqueness property [actually only $n+1$ and $n$ are used]
\sn
\item "{$(c)$}"  $\bold s$ is a brimmed$^\ell$ stable $({\Cal
P}^-(n),{\frak s}^+)$-system
\sn
\item "{$(d)$}"  $\bold s^*$ is a stable $({\Cal P}(n),
{\frak s}^+)$-system reduced at $n$ such that $\bold s^* 
\restriction {\Cal P}^-(n) = \bold s$.
\ermn
\ub{Then} $\bold s^*$ is prime over $\bold s$ for ${\frak
s}^+$. \nl
2) Moreover, $\bold s^*$ is strongly prime$^\ell$ over $\bold s$ for
${\frak s}^+$.
\endproclaim
\bigskip

\remark{Remark}  This is similar to the proof of the existence of
primes in ${\frak s}^+$.
\endremark
\bigskip

\demo{Proof}  1) Without loss of generality $\bold s$ is normal.  Let
$\langle M^\alpha_u:\alpha < \lambda^+_{\frak s} \rangle$ be $\le_{\frak
s}$-increasing continuous with union $M^{{\bold s}^*}_u$ and let $E$ be a 
thin enough club of $\lambda^+_{\frak s}$.  By \scite{705-12b.10} 
for each $\alpha \in E,\bold s^*_\alpha 
= \langle M^\alpha_u:u \in {\Cal P}(n) \rangle$
is a normal stable ${\Cal P}(n)$-system reduced at $n$ and 
letting ${\Cal P} = \{u \subseteq n+1:n \nsubseteq u\}$ for $\alpha < \beta$
from $E,\bold s^*_{\alpha,\beta} =: \bold s_\alpha *_{{\Cal P}(n+1)}
\bold s_\beta$ is a stable $({\Cal P}(n+1),{\frak s})$-system and
$\bold s_{\alpha,\beta} \restriction {\Cal P}$ is a normal
brimmed$^3 ({\Cal P},{\frak s})$-system, see \scite{705-12b.10}.  
Suppose that $M \in
K_{{\frak s}(+)}$ and $\langle f_u:u \in {\Cal P}^-(n) \rangle$ an
embedding of $\bold s \restriction {\Cal P}^-(n)$ into $M$ (see
Definition \scite{705-12.4}(3)).

As $\bold s$ is normal, we have $u \subseteq v \in {\Cal P}^-(n)
\Rightarrow f_u \subseteq f_v$.  Let $\langle M_\alpha:\alpha <
\lambda^+_{\frak s} \rangle$ be $\le_{\frak s}$-increasing continuous with
union $M$ and \wilog \, $E$ is a thin enough club for this too; by
renaming $E = \lambda^+_{\frak s}$.  Let $f^\alpha_u = f_u
\restriction M^\alpha_u$, so $\bar f^\alpha = \langle f^\alpha_u:u \in
{\Cal P}^-(n) \rangle$ is an embedding of $\bold s_\alpha$ into 
$M_\alpha \le_{{\frak K}[{\frak s}]} M$.
Now we choose $f^\alpha_n$ by induction on $\alpha$ such that
\mr
\item "{$\circledast(i)$}"  $f^\alpha_n$ is a $\le_{{\frak K}[{\frak
s}]}$-embedding of $M^\alpha_n$ into $M$ (hence into
$M_{\beta(\alpha)}$ for some $\beta(\alpha) < \lambda^+_{\frak s})$
\sn
\item "{$(ii)$}"  $f^\alpha_n$ extend $f^\beta_n$ for $\beta < \alpha$
and $f^\alpha_u$ for $u \in {\Cal P}^-(n)$.
\ermn
For $\alpha=0,f^0_n$ exists as ${\frak s}$ has the brimmed$^3$ weak
$(\lambda,n)$-uniqueness property and $M$ is $\lambda^+$-saturated. \nl
For $\alpha$ limit let $f^\alpha_n = \dbcu_{\beta < \alpha}
f^\beta_n$. \nl
For $\alpha = \beta +1$ let $\gamma < \lambda^+_{\frak s}$ be such
that $M_\gamma$ is $(\lambda,*)$-brimmed over Rang$(f^\beta_n) \cup
\bigcup \{\text{Rang}(f^\alpha_u:u \in {\Cal P}^-(n)\}$.  We shall
show that there
is a $\le_{\frak s}$-embedding of $M^\alpha_n$ into $M$ and even into
$M_\gamma$ extending $f_\beta \cup \bigcup\{f^\alpha_u:
u \in {\Cal P}^-(n)\}$.  For this we shall use 
``${\frak s}$ has the brimmed$^3$ weak $(\lambda,n+1)$-uniqueness 
property" defined in \scite{705-12b.5}(2) for the $(\lambda_{\frak s},{\Cal
P}^-(n+1))$-system $\bold s^*_{\alpha,\beta}$.  The embedding are
O.K. as $M^\alpha_n$ does not contribute (as $\bold s$ is reduced in
$n$).  But the assumption is not fully satisfied because $\bold
s^*_{\alpha,\beta}$ the brimmed$^\ell$ demand does not (necessarily)
holds for $u=n$, i.e., for $f^\alpha_n(M^\alpha_n)$; however, by Claim
\scite{705-12.u.9} this is overcomed.  We get that there is a pair $(f,N)$
such that $M_{\beta(\alpha)} \le_{\frak s} N$ and $f$ a 
$\le_{\frak s}$-embedding $M^\alpha_n$ into $N$, but \wilog \, 
$N \le M_\gamma$ so we are done.

Now $f_n = \dbcu_{\alpha < \lambda^+_{\frak s}} f^\alpha_n$ is the
required embedding. \nl
2) By part (1) and \scite{705-12p.10EA}.    \hfill$\square_{\scite{705-12f.10F}}$\margincite{705-12f.10F}
\enddemo
\bigskip

\demo{\stag{705-12f.10G} Conclusion}  Assume $\ell=3$ and
\mr
\item "{$(a)$}"  ${\frak s}$ is successful hence 
${\frak s}^+$ satisfies \scite{705-12.1}
\sn
\item "{$(b)$}"  ${\frak s}$ has the brimmed$^\ell$ weak 
$(\lambda,\le n+1)$-uniqueness property.
\ermn
\ub{Then} ${\frak s}^+$ has the brimmed$^\ell$ strong 
$(\lambda^+,n)$-primeness$^\ell$ property.
\enddemo
\bigskip

\demo{Proof}  By \scite{705-12f.10F}.
\enddemo
\bigskip

\proclaim{\stag{705-12f.4H} Claim}  1) Let $\ell=3$ and
\mr
\item "{$(a)$}"   ${\frak s}$ has the brimmed$^\ell$ 
strong $(\lambda,n)$-existence property
\sn
\item "{$(b)$}"  ${\frak s}$ has the brimmed$^\ell$ strong
$(\lambda,n)$-primeness property.
\ermn
\ub{Then} ${\frak s}$ has the brimmed$^\ell$ strong
$(\lambda,n)$-uniqueness property. \nl
2) Let $\ell=4$ and ${\frak s}$ has the brimmed$^\ell$ weak
$(\lambda,n)$-uniqueness property \ub{then} ${\frak s}$ has the
brimmed$^\ell$ strong $(\lambda,n)$-uniqueness property.
\endproclaim
\bigskip

\demo{Proof}  1) Assume $\bold d_1,\bold d_2$ are brimmed$^\ell$ stable 
$(\lambda,{\Cal P}(n))$-system and 
$\bar f = \langle f_u:u \in {\Cal P}^-(n) \rangle$ be an isomorphism 
from $\bold d_1 \restriction {\Cal P}^-(n)$ onto $\bold d_2
\restriction {\Cal P}^-(n)$.  As ${\frak s}$ has the
brimmed$^\ell$ strong $(\lambda,n)$-existence property, (i.e.,
assumption (a)), clearly for $k=1,2$ there is an expanded stable
$(\lambda,{\Cal P}(n))$-system $\bold d^k,\bold d^k \restriction
{\Cal P}^-(n) = \bold d_k \restriction {\Cal P}^-(n),\bold d^k$
reduced in $n$, i.e.,  such that $u
\subset n \Rightarrow \bold J^{{\bold d}^k}_{n,u} = \emptyset$). \nl
Let $f'_u = f_u$ for $u \subset n$.  Without loss of generality there
is an isomorphism $f'_n$ such that $\langle f'_u:u \in {\Cal P}(n)
\rangle$ is an isomorphism from $\bold d^1$ onto $\bold d^2$.

By clause (b) of the assumption + $(*)$ 
of Definition \scite{705-12b.5} \wilog \, $M^{{\bold d}^k}_n
\le_{\frak s} M^{{\bold d}_k}_n$; recall $u \subset n \Rightarrow
h^{{\bold d}^k}_{n,u} = h^{{\bold d}_k}_{n,u}$.  

Now 
\mr
\item "{$(*)$}"  $M^{\bold d^k}_n \le_{\frak s} M^{{\bold d}_k}_n$
and $p \in {\Cal S}^{\text{bs}}(M^{\bold d^k}_n) \Rightarrow \text{
dim}(p,M^{\bold d_k}_n) = \lambda$ (or just for a dense set of regular
$p \in {\Cal S}^{\text{bs}}(M^{\bold d^k}_n)$.
\ermn
[Why?  Let $p \in {\Cal S}^{\text{bs}}
(M^{{\bold d}^k}_n)$ be regular; \wilog \, such that for some $u(p) 
\subseteq n,p$ does not fork over
$M^{{\bold d}_k}_{n,u(p)}$ and $u \subset u(p)
\Rightarrow p$ orthogonal to $M^{{\bold d}_k}_{n,u}$.  Now 
dim$(p,M^{{\bold d}_k}_n) = \lambda$.  Why?  If $u(p) =n$ by $(**)$
from Definition \scite{705-12b.5}(5),(5A),(6) and if $u(p) \subset n$ as $\bold
d_k$ is brimmed$^\ell$ (i.e., the definition and basic properties of
dimension and regular types).]

We know that if $M$ is $(\lambda_{\frak s},*)$-brimmed$^\ell,N
\le_{\frak s} M$ and $p \in {\Cal S}^{\text{bs}}(N) \Rightarrow \text{
dim}(p,N) = \lambda_{\frak s}$ then $M$ is $(\lambda,*)$-brimmed over
$N$ (see \scite{705-10p.16}).  But this demand by $(*)$ holds with
$(M^{\bold d_k}_n,M^{\bold d^k}_n)$ here standing for $(M,N)$, hence
$M^{{\bold d}_k}_n$ is $(\lambda,*)$-brimmed over
$M^{{\bold d}^k}_n$.  As $f'_n$ is an isomorphism from $M^{{\bold d}^2}_n$ onto
$M^{{\bold d}^2}_n$ there is an isomorphism $f_u$ from $M^{{\bold d}_1}_n$ onto
$M^{{\bold d}_2}$ which extends $f'_n$.  So clearly
$\langle f_u:u \subseteq n \rangle$ is an isomorphism from $\bold d_1$
onto $\bold d_2$.  \hfill$\square_{\scite{705-12f.4H}}$\margincite{705-12f.4H}
\enddemo
\bigskip

\proclaim{\stag{705-12b.14} Theorem}  $[2^{\lambda^{+n}} <
2^{\lambda^{+n+1}}$ for $n < \omega]$.  Let $\ell = 3$.
Assume ${\frak s}$ is $\omega$-successful.  \ub{Then} for every $n \ge
2$ and $m \ge n-2,{\frak s}^{+ m}$ is $n$-excellent$^\ell$; 
(see definition below).
\endproclaim
\bigskip

\definition{\stag{705-12r.15} Definition}  1) We say that ${\frak s}$ is
$n$-beautiful$^\ell$ if:
\mr
\item "{$(a)$}"   ${\frak s}$ has the brimmed$^\ell$ strong
$(\lambda,\le n)$-existence property
\sn
\item "{$(b)$}"   ${\frak s}$ has the brimmed$^\ell$ weak
$(\lambda,\le n)$-uniqueness property
\sn
\item "{$(c)$}"   ${\frak s}$ has the brimmed$^\ell$ strong
$(\lambda,<n)$-uniqueness property
\sn
\item "{$(d)$}"   ${\frak s}$ has the brimmed$^\ell$ strong
$(\lambda,<n)$-primeness property.
\ermn
2) We say that ${\frak s}$ is $\omega$-beautiful$^\ell$ \ub{if}
${\frak s}$ is $n$-beautiful$^\ell$ for every $n$.
\enddefinition
\bigskip

\remark{Remark}  In the Theorem we could restrict our demand to $n \le
n_* (< \omega)$ and get the appropriate conclusion, essentially
${\frak s}^{+m}$ is $n$-excellent$^3$ \ub{if} ${\frak s}$ is
$2n$-successful (and $\langle 2^{\lambda^{+\ell}}:\ell \le 2n \rangle$
is increasing).
\endremark
\bigskip

\demo{Proof of \scite{705-12b.14}}  We know that
\mr
\item "{$(*)_1$}"   ${\frak s}^{+m}$ satisfies the demand in
\scite{705-12.1} and is successful.
\ermn
We now prove by induction on $n \ge 2$ that
\mr
\item "{$\boxtimes_n$}"  ${\frak s}^{+m}$ is $n$-beautiful$^\ell$ if
$m \ge n-2$.
\ermn
First we prove $\boxtimes_2$. 

By \scite{705-12p.15B} for $n=0,1$ the demands in \scite{705-12r.15} holds.
Now ${\frak s}^{+m}$ has the brimmed$^\ell$ weak $(\lambda,2)$-uniqueness as
${\frak s}$ is a good frame (i.e., the uniqueness of NF$_{{\frak
s}(+m)}$-amalgamation.  Lastly, the brimmed$^\ell$ strong
$(\lambda,2)$-existence holds by \scite{705-12f.10E}.

So let $n \ge 2$ and we assume $\boxtimes_n$ and we shall prove
$\boxtimes_{n+1}$, this suffices
\mr
\item "{$(*)_2$}"  there is a brimmed$^\ell ({\Cal P}(n+1),{\frak
s}^{+m})$-system for $m \ge n-2$ \nl
[Why?  By \scite{705-12f.10D}.]
\sn
\item "{$(*)_3$}"   ${\frak s}^{+m}$ has the brimmed$^\ell$ weak
$(\lambda,n+1)$-existence property if $m \ge n-2$.
\ermn
[Why?  By $(*)_2$ there is brimmed$^\ell ({\Cal P}(n+1),{\frak
s}^{+m})$-system call it $\bold d^*$.  Let a $({\Cal P}^-(n+1),
{\frak s}^{+m})$-system $\bold d$ be given.  By $\boxtimes_n$ we know
that ${\frak s}$ has the strong $(\lambda,<n)$-uniqueness property
hence by \scite{705-11b.7}(2) 
\wilog \, $\bold d^*
\restriction [n+1]^{<n},\bold d \restriction [n+1]^{<n}$ are
isomorphic so \wilog \, they are equal.  Now we shall apply clause
$(\alpha)$ of the conclusion of \scite{705-12.u10} to $\bold d$ and $\bold
d^* \restriction {\Cal P}^-(n+1)$, as the latter has the weak
existence property (as $\bold d^*$ exemplify) it suffices to check the
assumptions of \scite{705-12.u10}.  So here $I = {\Cal P}^-(n+1)$ and $J =
[n+1]^{<n}$, clause (a) of \scite{705-12.u10} is obvious, clause $(b)$ was
assumed above and clause $(c)_1$ follows from ``${\frak s}^{+m}$ has
the weak $n$-uniqueness property, which holds as we assume $\boxtimes_n$.]
\mr
\item "{$(*)_4$}"   ${\frak s}^{+m}$ has the brimmed$^\ell$ weak
$(\lambda,n+1)$-uniqueness property if $m \ge n-1$. \nl
[Why?  We try to apply \scite{705-12v.11}(2) hence implicitly
\scite{705-12b.11} + \scite{705-12v.11}(1) to 
${\frak s}^{+m}$ and $n$.
Its conclusion fails by clause (b) of Definition \scite{705-12r.15},
aplied to $({\frak s}^{+m})^+ = {\frak s}^{+m+1}$ for $n$ which holds
as we are assuming $\boxtimes_n$.
Clause (a) from its assumptions holds by $(*)_1$, 
clause (b) holds by $(*)_3$ for $n+1$ and by clause $(a)$ of Definition
\scite{705-12r.15} and for $m \le n$.  \nl
Now clause $(c)^-_2$ holds by clause (c) of \scite{705-12r.15} by
$\boxtimes_n$ applied to ${\frak s}^{+m}$ and clause $(c)^-_1$ holds
by clause (b) of Definition \scite{705-12r.15} by $\boxtimes_n$ applied to
${\frak s}^{+m}$.  So in \scite{705-12v.11} only the fourth assumption
(d), may fail, so as the conclusion fails, (d) there fails.
\nl
Hence clause (d) from \scite{705-12b.11} has to fails which is
the desired conclusion.]
\sn
\item "{$(*)_5$}"   ${\frak s}^{+m}$ has the brimmed$^\ell$ strong
$(\lambda,n)$-primeness property if $m \ge n-1$. \nl
[By \scite{705-12f.10G} applied to $n' = n$ and ${\frak s}' 
= {\frak s}^{+(m-1)}$.  It gives the desired conclusion.  As for its assumption
clause $(a)$ there holds by $(*)_1$ and clause (b) there for $n+1$ 
by $(*)_4$ above and for $0,\dotsc,n$ by $\boxtimes_n$.]
\sn
\item "{$(*)_6$}"   ${\frak s}^{+m}$ has the brimmed$^\ell$ strong
$(\lambda,n)$-uniqueness property if $m \ge n-1$. \nl
[Why?  By \scite{705-12f.4H}, assumption (a) there holds by clause (a) of
the definition \scite{705-12r.15} of excellent 
and clause (b) there holds by $(*)_5$ above.]
\sn
\item "{$(*)_7$}"   ${\frak s}^{+m}$ has the brimmed$^\ell$ strong
$(\lambda^{+m},n+1)$-existence property for $m \ge n-1$. \nl
[Why?  By \scite{705-12f.10E}, its conclusion is what we need, assumption
(a) there holds by $\boxtimes_n$.  Assumption (b) there holds by
$(*)_5$ and \scite{705-12b.6}(2d) (which says that strong prime
$\Rightarrow$ weakly prime.  Lastly, assumption (c) there holds by
$(*)_6$ above so we are done.]
\ermn
So $\boxtimes_{n+1}$ holds.   \hfill$\square_{\scite{705-12b.14}}$\margincite{705-12b.14}
\enddemo
\bn
Recall that \chaptercite{600} has tried generalizing \cite{Sh:87a},
\cite{Sh:87b} but through it give the parallel conclusions about each
$\lambda^{+n}$, it does not say anything on $\mu \ge \lambda^{+
\omega}$.  In the claims (\scite{705-12f.16R}), \scite{705-12f.16A},
\scite{705-12f.17K}, \scite{705-12f.18} below we derived the parallel of several of
the further conclusions of \cite{Sh:87a}, \cite{Sh:87b}. \nl
The aim of the following claim is to help proving for the case
$\ell=3$ that for
non-unidimensional ${\frak s}$, we can prove non-categoricity in
higher cardinals (of course, we shall get better results when we prove
excellency for $\ell=1$).
\sn
\proclaim{\stag{705-12f.16R} Claim}  Assume that
\mr
\item "{$(a)$}"  ${\frak s}$ is a good $\lambda$-frame
\sn
\item "{$(b)$}"  $M^* \in K_{\frak s},\langle c_i:i < \lambda_{\frak
s} \rangle$ list the elements of $M^*$, \nl
$\bold P \subseteq {\Cal S}^{\text{bs}}(M^*)$ is a
non-empty set of regular types 
such that there is $q \in {\Cal S}^{\text{bs}}(M^*)$ orthogonal to
$\bold P$ (so ${\frak s}$ is not weakly unidimensional)
\sn
\item "{$(c)$}"  $\tau^* = \tau \cup \{c_i:c \in M_*\}$
\sn
\item "{$(d)$}"  $K^* = \{M:M$ is a $\tau^*$-model and $M \restriction
\tau \in K^{\frak s},c_i \mapsto c^M_i$ is a $\le_{{\frak K}[{\frak
s}]}$-embedding of $M^*$ into $M \restriction \tau\}$ and if $q \in
{\Cal S}^{\text{bs}}(M^*)$ is regular realized in $M \restriction
\tau$ then it is orthogonal to ${\bold P}$
\sn
\item "{$(e)$}"  $M_1 \le_{{\frak K}^*} M_2$ \ub{iff} $(M_1,M_2 \in
K^*$ and) $M_1 \restriction \tau \le_{{\frak K}[{\frak s}]} M_2
\restriction \tau$
\sn
\item "{$(f)$}"  ${\frak s}^* = ({\frak K}^*,
{\Cal S}^{\text{bs}}_*,\nonfork{}{}_{*})$ where 
{\roster
\itemitem{ $(i)$ }  ${\Cal S}^{\text{bs}}_*(M_1)$ is
essentially 
$\{{\text{\rm tp\/}}(a,M_1,M_2):M_1 \le_{{\frak K}^*} M_2$ is of
cardinality $\lambda_{\frak s}$, {\rm tp}$(a,M_1 \restriction \tau,M_2
\restriction \tau \in {\Cal S}^{\text{bs}}_{\frak s}(M_1 \restriction
\tau)\}$
\sn
\itemitem{ $(ii)$ }   $\nonfork{}{}_{*}$ similarly, i.e.,
$\nonfork{}{}_{*}(M_0,M_1,a,M_3)$ \ub{iff} 
$M_0 \le_{{\frak K}^*} M_1 \le_{{\frak K}^*} M_3,
a \in M_3$ and $\nonfork{}{}_{*}(M_0,M_1,a,M_3)$.
\endroster}
\ermn
\ub{Then}
\mr
\item "{$(\alpha)$}"  ${\frak s}^*$ is a good $\lambda$-frame
\sn
\item "{$(\beta)$}"  $K^{{\frak s}^*} \subseteq K^{\frak s},
K^{{\frak s}^*}_\lambda \ne \emptyset,K^{{\frak s}^*}_{\lambda^+} 
\ne \emptyset$
\sn
\item "{$(\gamma)$}"  if $\bold I(\lambda^{+n+1},K^{{\frak s}^*})
< \mu_{\text{wd}}(\lambda^{+n+1},2^{\lambda^{+n}})$ for $n < \omega$,
\ub{then} ${\frak s}^*$ is $\omega$-successful and is well 
defined for every $\mu > \lambda$.
\endroster
\endproclaim
\bigskip

\demo{Proof} \ub{Clause $(\alpha)$}.

Check
\mn
\ub{Clause $(\beta)$}.

Trivial.
\mn
\ub{Clause $(\gamma)$}.

By \scite{705-12b.14} applied to ${\frak s}^*$.
\enddemo
\bigskip

\demo{\stag{705-12f.16A} Major Conclusion}   Assume that ${\frak s}$
satisfies the conclusion of \scite{705-12b.14} and let ${\frak t} 
= {\frak s}^{+ \omega}$ (see Definition \scite{705-0.X}(4)) 
\ub{then}
\mr
\item "{$(a)$}"  ${\frak s}^{+ \omega} = {\frak s}(+ \omega)$ is a good
$\lambda^{+ \omega}_{\frak s}$-frame (recall that 
${\frak K}_{{\frak s}(+\omega)}$ 
is $\cap \{{\frak K}^{{\frak s}(+n)}_{\lambda^{+ \omega}}:n < \omega\}$)
\sn
\item "{$(b)$}"  ${\frak t} = {\frak s}^{+ \omega}$ is 
$\omega$-beautiful$^3$ 
\sn
\item "{$(c)$}"  for $\mu \ge \lambda^{+ \omega},{\frak t}[\mu] 
= {\frak s}[\mu]$, see Definition \sciteu{705-xxX} is a good
$\mu$-frame, which is $\omega$-beautiful$^3$, categorical in $\mu$ for
every $\mu \ge \lambda^{+ \omega}$, so in particular $\mu \ge \lambda
\Rightarrow K^{\frak s}_\mu \ne \emptyset$ (on ${\frak t}[\mu]$ 
see Definition \scite{705-0.X}
\sn
\item "{$(d)$}"  if ${\frak s}$ is weakly unidimensional and $\mu \ge
\lambda$  \ub{then} ${\frak s}[\mu] = {\frak s}[\mu] = {\frak s}[\mu]$; hence
$K^{{\frak s}(+)}$ is categorical in $\mu$ for every $\mu \ge
\lambda^+_{\frak s}$ (if $K_{\frak s}$ is categorical in
$\lambda_{\frak s}$ then $K^{\frak s}$ is categorical in $\mu$ for
every $\mu \ge \lambda_{\frak s}$
\sn
\item "{$(e)$}"  if ${\frak s}$ is not weakly unidimensional \ub{then}
${\frak s}(\mu)$ is not weakly unidimensional
\sn
\item "{$(f)$}"  ${\frak s}$ has NDOP iff ${\frak t} = {\frak s}^{+ \omega}$
has NDOP iff ${\frak s}(\mu)$ has NDOP (for any $\mu \ge \lambda_{\frak s}$)
\sn
\item "{$(g)$}"  $K^{\frak s}_\mu \ne \emptyset$ for $\mu \ge
\lambda,{\frak s}[\mu]$ well defined for $\mu \ge \lambda$ (on ${\frak
s}[\mu]$ see Definition \scite{705-0.X}(4)).
\endroster
\enddemo
\bigskip

\demo{\stag{705-12f.17K} Conclusion}  Assume $2^{\lambda^{+n}} <
2^{\lambda^{n+1}}$ for $n < \omega$ and
\mr
\item "{$(a)$}"  ${\frak s}$ is a good $\lambda$-frame not weakly
unidimensional
\sn
\item "{$(b)$}"  $I(\lambda^{+n+1},K^{\frak s}) <
\mu_{\text{wd}}(\lambda^{+n},2^{\lambda^{+n}})$ for $n < \omega$.
\ermn
\ub{Then} $K^{\frak s}$ is not categorical in $\mu$, for every $\mu > \lambda$.
\enddemo
\bigskip

\demo{Proof}  For each $\mu > \lambda$, there is $M^1_\mu \in
K_{{\frak s}[\mu]}$ see \scite{705-12f.16A}.  But we can 
define ${\frak t}$ as ${\frak s}^*[+ \omega]$, where ${\frak s}^*$ is 
as in \scite{705-12f.16R} and apply it to \scite{705-12b.14},
\scite{705-12f.16A}, and get $M^2_\mu \in K_{{\frak t}[\mu]}$.  Looking at
$\le_{K[{\frak s}]}$-submodels of $M^1_\mu,M^2_\mu$ it is clear that
$M^1_\mu \approx M^2_\mu$ so we are done.
\enddemo
\bn
We can sum up
\demo{\stag{705-12f.18} Conclusion}  Assume $2^{\lambda^{+n}} <
2^{\lambda^{+n+1}}$ for $n < \omega$.  If an a.e.c. ${\frak K}$ with
LS$({\frak K}) \le \lambda$, is categorical in
$\lambda,\lambda^+,1 \le I(\lambda^{++},K)$ and $I(\lambda^{+n+2},K) <
\mu_{\text{wd}}(\lambda^{n+2},2^{\lambda^{+n+1}})$ for $n < \omega$ \ub{then}
${\frak K}$ is categorical in every $\mu \ge \lambda$.
\enddemo
\bigskip

\remark{\stag{705-12f.19K} Remark}  1) This through light on \cite{MaSh:285},
\cite{KlSh:362}, \cite{Sh:472}, \cite{Sh:394}. \nl
In those works we start with an appropriate a.e.c. ${\frak K}$ and
assume that it is categorical in $\lambda$ large enough then
LS$({\frak K})$ and prove that for some $\alpha_* (2^{\text{LS}({\frak
K})})^+$ that the class is categorical in every $\lambda' \in
[\beth_{\alpha_*},\lambda)$, but nothing is said about $\lambda' >
\lambda$.  However, if for some $\mu,\mu^{+ \omega} \in
(\beth_{\alpha_*},\lambda]$ then by \scite{705-12f.18} we are done.  This
weak set theoretic assumption will be eliminated in a sequel. \nl
2) Moreover, we can eliminate the ``$\lambda$ successor" assumption.
\nl
3) We can say much more: $\omega$-successful frames are very much like
superstable first order classes and more.  We delay this.
\endremark
\newpage

     \shlhetal 

\nocite{ignore-this-bibtex-warning} 
\newpage
    
REFERENCES.  
\bibliographystyle{lit-plain}
\bibliography{lista,listb,listx,listf,liste}

\enddocument